\documentclass[a4paper, reqno]{amsbook}

\usepackage[utf8]{inputenc}
\usepackage[english]{babel}
\usepackage{amsmath}
\usepackage{amssymb}
\usepackage{graphicx}
\usepackage{xcolor}
\usepackage{caption}
\usepackage{tikz}
\usetikzlibrary{arrows.meta,calc}
\usepackage[all]{xy}
\usepackage{fontawesome5}
\usepackage{etoolbox}
\usepackage{hyperref}

\makeatletter
\newcommand*\bigcdot{\mathpalette\bigcdot@{.5}}
\newcommand*\bigcdot@[2]{\mathbin{\vcenter{\hbox{\scalebox{#2}{$\m@th#1\bullet$}}}}}
\makeatother

\newcommand{\pr}{\mathrm{pr}}
\DeclareMathOperator{\Hom}{Hom}

\DeclareMathOperator{\G}{G}

\DeclareMathOperator{\bconcat}{\curlywedge}
\newcommand{\concat}{\ensuremath{\text{\Large $\curlywedge$}}}
\newcommand{\bb}[1]{\ensuremath{\lvert #1 \rvert}}
\newcommand{\vs}[1]{\mathbin{\downarrow}#1}

\newcommand{\sig}[1]{\mathbf{#1}}             
\newcommand{\calg}[1]{\boldsymbol{\mathsf{#1}}} 
\newcommand{\spec}[1]{\boldsymbol{\mathcal{#1}}} 
\newcommand{\eeq}{\overset{\mathrm e}{=}}      
\newcommand{\QE}{\mathrm{QE}}
\newcommand{\Dom}{\mathrm{Dom}}
\newcommand{\gd}{\mathfrak{gd}}              
\newcommand{\htg}{\mathrm{htg}}                

\newcommand{\id}{\mathrm{id}}
\newcommand{\Id}{\mathrm{Id}}
\newcommand{\src}{\mathrm{sc}}
\newcommand{\tgt}{\mathrm{tg}}

\newcommand{\Gd}{\mathrm{Gd}}
\newcommand{\Uf}{\mathrm{Uf}}
\newcommand{\Ds}{\mathrm{Ds}}
\newcommand{\Tr}{\mathrm{Tr}}

\newcommand{\bGd}{\mathbf{Gd}}
\newcommand{\bUf}{\mathbf{Uf}}
\newcommand{\bDs}{\mathbf{Ds}}
\newcommand{\bTr}{\mathbf{Tr}}

\NoCompileMatrices

\xymatrixrowsep = {8ex}
\xymatrixcolsep = {10ex}

\newdir{+>}{@{|}*@{-}!/-6.5pt/@{>}*!/-13pt/{}}
\newdir{ +}{{}*!/-5pt/@{+}}

\newsavebox{\xymor}
\newsavebox{\xydmor}

\newcommand{\xynode}{\makebox[0ex]{}}
\savebox{\xymor}{\ensuremath{%
		\xymatrix@1@C=19pt{\xynode \ar@{>}[r] & \xynode }}}
\savebox{\xydmor}{\ensuremath{%
		\xymatrix@1@C=19pt{\xynode \ar@{{}{-}_{>}}[r] & \xynode }}}

\newcommand{\mor}{\usebox{\xymor}}    
\newcommand{\dmor}{\usebox{\xydmor}}  

\newcommand{\function}[4]{
	\begin{array}{@{\:}c@{\:}c@{\:}l}
		#1 &\mor& #2 \\
		#3 &\longmapsto& #4
\end{array} }
\newcommand{\nfunction}[4]
{\left\{
	\function{#1}{#2}{#3}{#4}
	\right. }

\allowdisplaybreaks
\theoremstyle{plain}
\newtheorem{theorem}{Theorem}[chapter]
\newtheorem{proposition}[theorem]{Proposition}
\newtheorem{corollary}[theorem]{Corollary}
\newtheorem{lemma}[theorem]{Lemma}
\newtheorem{claim}[theorem]{Claim}
\newtheorem{remark}[theorem]{Remark}

\theoremstyle{definition}
\newtheorem{definition}[theorem]{Definition}
\newtheorem{example}[theorem]{Example}
\newtheorem*{convention}{Convention}

\newcommand*{\leandocs}{https://mariovagomarzal.github.io/categorial-algebras/docs/find/?pattern=}
\definecolor{lean@equivalent}{HTML}{009E73}
\definecolor{lean@weaker}{HTML}{D55E00}
\definecolor{lean@stronger}{HTML}{0072B2}
\makeatletter
\newcommand*{\CA@leanpending}{}
\newcommand*{\CA@leanrelation}{}
\newcommand*{\lean}[2]{%
  \ifcsname\string\color@lean@#1\endcsname\else
    \PackageError{main}{Unknown relation `#1' in \string\lean}{Use equivalent, weaker or stronger.}%
  \fi
  \gdef\CA@leanrelation{#1}%
  \gdef\CA@leanpending{#2}}
\newcommand*{\CA@leanlogo}{\textsf{L$\exists\forall$N}}
\newcommand*{\leansample}[1]{{\normalfont\footnotesize\color{lean@#1}\CA@leanlogo}}
\newcommand*{\CA@leanlink}[2]{%
  \in@{CategoryTheory.}{#1}%
  \ifin@\def\CA@leanfull{#1}\else\def\CA@leanfull{CategorialAlgebras.#1}\fi
  {\hypersetup{pdfborder={0 0 0}}\href{\leandocs\CA@leanfull\#doc}{#2}}}
\newcount\CA@leancount
\newcommand*{\CA@leanlinks}{%
  \global\CA@leancount=0
  \@for\CA@leanname:=\CA@leanpending\do{%
    \global\advance\CA@leancount by 1
    \ifnum\CA@leancount=1
      \expandafter\CA@leanlink\expandafter{\CA@leanname}{\CA@leanlogo}%
    \else
      \hspace{0.25em}\expandafter\CA@leanlink\expandafter{\CA@leanname}{%
        {\scriptsize\the\CA@leancount}}%
    \fi}}
\newcommand*{\CA@leanmark}{%
  \ifx\CA@leanpending\@empty\else
    \ifhmode\unskip\else\noindent\fi
    {\nobreak\hfil\penalty50\hskip1em\null\nobreak\hfil
      \mbox{\normalfont\footnotesize\color{lean@\CA@leanrelation}\CA@leanlinks}%
      \parfillskip=\z@\finalhyphendemerits=\z@\par}%
    \global\let\CA@leanpending\@empty
  \fi}
\@for\CA@env:=theorem,proposition,lemma,corollary,definition,remark,example\do{%
  \expandafter\AtEndEnvironment\expandafter{\CA@env}{\CA@leanmark}}
\makeatother

\begin{document}

\hypersetup{pdftitle={Higher-order many-sorted categorial algebras via generalized hep-derivors},
            pdfauthor={J. Climent Vidal, E. Cosme Llópez, R. Ruiz Mora, M. Vago Marzal}}

\frontmatter

\thispagestyle{empty}
\begin{center}
\vspace*{2.5cm}
{\Huge\bfseries Higher-order many-sorted\\[10pt] categorial algebras\\[10pt] via generalized hep-derivors\par}

\vspace{3cm}

{\bfseries J. Climent Vidal\par}
\vspace{4pt}
{\small\itshape Universitat de Val\`{e}ncia\\
Departament de L\`{o}gica i Filosofia de la Ci\`{e}ncia\par}
{\small Juan.B.Climent@uv.es\par}

\vspace{0.6cm}

{\bfseries E. Cosme Ll\'{o}pez\par}
\vspace{4pt}
{\small\itshape Universitat de Val\`{e}ncia\\
Departament de Matem\`{a}tiques\par}
{\small Enric.Cosme@uv.es (\faEnvelope[regular])\par}
\vspace{3pt}
{\small\itshape Nantong University\\
School of Mathematics and Statistics\par}

\vspace{0.6cm}

{\bfseries R. Ruiz Mora\par}
\vspace{4pt}
{\small\itshape Universitat de Val\`{e}ncia\\
Departament de Matem\`{a}tiques\par}
{\small Raul.Ruiz-Mora@uv.es\par}

\vspace{0.6cm}

{\bfseries M. Vago Marzal\par}
\vspace{4pt}
{\small\itshape Universitat de Val\`{e}ncia\\
M\`{a}ster en Investigaci\'{o} Matem\`{a}tica\par}
{\small mavamar9@alumni.uv.es\par}

\vfill
{\small November 2026\par}
\end{center}
\cleardoublepage

\thispagestyle{empty}
\begin{center}
\vspace*{2cm}
{\bfseries Abstract\par}
\end{center}
\vspace{0.5cm}
Higher-order many-sorted categorial algebras are many-sorted algebras whose carriers are strict higher-order categories and whose operations are higher-order functors.
They admit a homogeneous and a heterogeneous presentation, both by partial algebras whose operation domains are given by existence equations.
To relate these presentations and their dimensions, we introduce generalized hep-derivors between hierarchical equationally partial signatures.
A generalized hep-derivor selects carriers by existence equations in one variable and interprets operation symbols by terms, with operation domains evaluated recursively in the transformed algebra.
Every generalized hep-derivor induces a functor between categories of many-sorted partial algebras. The values of this functor satisfy the domain conditions wherever their operations are defined, and we characterize when the operations are defined exactly where these conditions hold.
The functor induced by a composite receives a natural inclusion from the composite of the induced functors; this inclusion, which can be strict, makes the induced functors a normal lax indexed category.
In finite dimensions, generalized hep-derivors realize graduation and unification as an equivalence between the homogeneous and heterogeneous presentations, and discretization as a left adjoint to truncation.
In dimension $\omega$, we reconstruct the categories of $\omega$-categorial algebras as projective limits of their finite-dimensional truncations and obtain unification by a quotient construction; we also establish algebraic squares and compatibility with change of signature.
Finally, for nonempty sets of sorts, we show that the finite-dimension condition in the homogeneous $\omega$-presentation is not axiomatizable by quasi-existence equations on its fixed signature.

\bigskip
\noindent\textit{2020 Mathematics Subject Classification.} Primary 18N30; Secondary 08A55, 08A68, 18C10, 18D30.

\smallskip
\noindent\textit{Key words and phrases.} Many-sorted partial algebra, equationally partial signature, derivor, strict higher-order category, categorial algebra, graduation, truncation, Grothendieck construction.
\cleardoublepage

\setcounter{tocdepth}{1}
\tableofcontents

\mainmatter

\chapter{Introduction}\label{SIntroduction}

Higher-order many-sorted categorial algebras are many-sorted algebras whose carriers are strict higher-order categories, one for each sort, and whose operations are higher-order functors.
For a set of sorts $S$, an $S$-sorted signature $\Sigma$ and $n\in\mathbb{N}$, the forgetful functor from the category $\mathsf{nCatAlg}(\sig{\Sigma})$ of $n$-categorial $\sig{\Sigma}$-algebras to $\mathsf{Set}^{S}$ factors through the category of $S$-sorted strict $n$-categories; for $S=1$ and the empty signature, $\mathsf{nCatAlg}(\sig{\Sigma})$ is the category $\mathsf{nCat}$ of strict $n$-categories.
These algebras arise naturally in the study of higher-order many-sorted rewriting systems, since the paths of a rewriting system of order $n$, up to suitable congruences, form an $n$-categorial algebra.
For $n=1$ they are the structures that interpret the rewrite theories of rewriting logic, in which the cells are the rewrites, the composition is their sequential composition, and the operations act on rewrites in parallel (see~\cite{Meseguer1992}).

As with strict higher-order categories, there are two natural ways of presenting higher-order categorial algebras.
In the \emph{homogeneous} presentation, which goes back to~\cite{Street1987}, all the cells of a given sort, whatever their dimension, form a single set, and the dimension of a cell is read off from its sources and targets; in the \emph{heterogeneous} presentation (see, e.g.,~\cite{AraEtAl2025}), the cells of every dimension form a separate sort, and there are sources, targets and identities between dimensions.
In both cases, the composition of cells is only partially defined, so higher-order categorial algebras are partial algebras, and both presentations can be described by quasi-existence equations, or $\QE$-equations, in the sense of Burmeister~\cite{Burmeister1986, Burmeister2002}.
Moreover, the domains of the compositions are given by existence equations, so the corresponding signatures are equationally partial signatures in the sense of Reichel~\cite{Reichel1987}.

The main aim of this work is to relate, by means of functors induced by morphisms between signatures, the two presentations and the different dimensions of higher-order categorial algebras.
Graduation and unification change the presentation, and discretization and truncation change the dimension.
Two features of these functors determine the kind of morphism between signatures that we need.
First, some of them do not keep all the elements of the given algebra: for instance, the cells of sort $(k,s)$ of the graduation of a homogeneous $n$-categorial algebra are only the $k$-cells of sort $s$, i.e., the cells $x$ with $\src^{k}_{s}(x)=x$, and the truncation to a dimension $m$ keeps only the $m$-cells.
These subsets are described by existence equations in one variable.
Second, some of them restrict the domains of the operations in a way that no term can express: the discretization of an $m$-categorial algebra into an $n$-categorial algebra, for $m<n$, interprets the $k$-composition, for $m\leq k<n$, by means of a variable, which is defined everywhere, whereas, in the discretized algebra, two cells are $k$-composable only if they are equal.
This restriction is imposed by the domain conditions of the composition, evaluated in the transformed algebra itself.

To deal with both features we introduce generalized hep-derivors, which combine two ideas.
On the one hand, the derivors of~\cite{ClimentSoliveres2010b}, which send every operation symbol of a signature to a term over another, induce functors between the categories of algebras in the opposite direction.
On the other hand, the general derivors of~\cite{LucasGutierrez2018} additionally send every sort to a formula that delimits, inside the interpretation of its image, the elements that interpret it.
A generalized hep-derivor between hierarchical equationally partial signatures, or hep-signatures, sends every sort to a set of existence equations in one variable, its \emph{domain}, and every operation symbol to a term.
The domains of the transformed operations are then determined by the domain conditions of the operation symbols, and the hierarchical character of the signatures allows us to define the transformed algebras by well-founded recursion.

\section{Main results}

The general results of this work are those of Chapter~\ref{SGenHepDerivors}.
For a generalized hep-derivor $\mathbf{c}$ from a hep-signature $\sig{\Sigma}$ to a hep-signature $\sig{\Lambda}$, we construct, by well-founded recursion, the transformed algebra $\mathbf{c}^{\ast}_{\gd}(\mathbf{B})$ of every partial $\Lambda$-algebra $\mathbf{B}$; this algebra is always a weak $\sig{\Sigma}$-algebra, and we prove that the construction is functorial.
We characterize when the transformed algebra is a $\sig{\Sigma}$-algebra, and we prove that hep-signatures and generalized hep-derivors form a category $\mathsf{HSig}_{\gd}$.
There is a natural inclusion of the composite of the induced functors into the functor induced by the composite, and we give a criterion for this inclusion to be an identity.
The inclusion can be strict. For instance, on $\{0,1\}$, suppose that we select only $0$ and interpret a unary operation by the constant $1$.
The transformed operation is nowhere defined because its value leaves the selected carrier.
If we then select the points at which this operation is defined, we obtain the empty set, whereas if we first translate the two selections, we obtain $\{0\}$.
These inclusions make the induced functors a normal lax indexed category with values in the categories of weak ep-algebras; its Grothendieck construction contains the category $\mathsf{Alg}_{\mathfrak{d}}$ of~\cite{ClimentSoliveres2010b}.
For many-sorted signatures and derivors, regarded as hep-signatures and generalized hep-derivors with empty domain conditions and empty domains, the construction is the usual one of~\cite{ClimentSoliveres2010b}.

The rest of this work applies these results to higher-order categorial algebras.
In finite dimensions, graduation, unification, discretization and truncation are all induced by generalized hep-derivors: graduation and unification are mutually quasi-inverse equivalences, discretization is left adjoint to truncation, and the four adjunctions form an algebraic square.
In dimension $\omega$, unification cannot be induced by a generalized hep-derivor when the set of sorts is nonempty, and it is obtained instead from the compatible finite-dimensional unifications by a quotient construction; the categories of $\omega$-categorial algebras are the projective limits of their finite-dimensional truncations, and the same equivalences, adjunctions and algebraic squares hold.
Finally, all these functors commute with the change of signature along a derivor, both in finite dimensions and in dimension $\omega$. Hence the families of these functors are natural transformations between split indexed categories over the category of many-sorted signatures and derivors; these natural transformations extend to the corresponding Grothendieck constructions, and the equivalences, adjunctions and algebraic squares are compatible with the change of signature.

\section{Relation to other work}

For $S=1$ and the empty signature, several of the facts that we obtain are well known: the equivalence between homogeneous and heterogeneous strict $n$-categories, the adjunction between discretization and truncation and, in the heterogeneous presentation, the description of the category of strict $\omega$-categories as the projective limit of the categories of strict $n$-categories (see, e.g., \cite{AraEtAl2025}, Section~14.2).
Our contribution is the explicit realization of these constructions through generalized hep-derivors, with recursively evaluated operation domains, weak transformed algebras and a lax comparison under composition, which may be strict.
The extension to arbitrary many-sorted algebraic operations is compatible with the finite-product structure: such operations are maps from finite products of carriers, and equivalences preserve these products up to canonical isomorphism.
The many-sorted extension is therefore the setting of our applications; what we add is the explicit syntactic realization and the verification of exact domains and comparison maps.

An $n$-categorial $\sig{\Sigma}$-algebra is a $\Sigma$-algebra whose carriers are strict $n$-categories and whose operations are strict $n$-functors from the cartesian products of the carriers, in the spirit of the structured categories of Ehresmann~\cite{Ehresmann1964}.
Equivalently, for $n=1$, a $1$-categorial $\sig{\Sigma}$-algebra is a category internal to the category $\mathsf{Alg}(\Sigma)$ of $\Sigma$-algebras, and internal categories and groupoids in varieties of algebras have been studied extensively (see, e.g., \cite{Bourn2017}).
Higher-order categorical structures endowed with additional algebraic structure, and the comparison between different presentations of them, have a long history (see, e.g., \cite{BrownHiggins1981} for the equivalence between $\infty$-groupoids and crossed complexes and \cite{Leinster2004} for a comparison of definitions of higher-order categories).

\textbf{General derivors.}
A general derivor of~\cite{LucasGutierrez2018} assigns formulas to sorts and operation symbols.
Its soundness conditions ensure that the operation formulas determine unique outputs in the selected carriers; preservation of a source theory additionally requires satisfaction of its translation in the target structure (\cite{LucasGutierrez2018}, Definition~2 and Theorem~1).
A generalized hep-derivor, instead, sends operation symbols to terms, and the transformed operations are partial: an operation symbol $\sigma$ is defined at a family of cells if, and only if, the family satisfies the domain conditions of $\sigma$ in the transformed algebra, the term $c(\sigma)$ is defined at it and its value is a cell.
This construction requires no prior soundness assumption and always produces a weak ep-algebra.
Exactness is characterized separately: the interpreting terms must be defined and take values in the selected carriers whenever the source domain conditions hold.
The axioms of a source specification must then also be verified to obtain a model.
This separation of construction, exactness and satisfaction of axioms is a matter of organization, not a claim of greater expressive power than that of sound general derivors.
When the domain conditions are empty and, for every operation symbol $\sigma$, $c(\sigma)$ is defined at every family of cells and takes cells as values, $\sigma$ is interpreted by the function defined by $c(\sigma)$, as for a sound general derivor whose formula for $\sigma$ is the graph of the term operation of $c(\sigma)$.

\textbf{Essentially algebraic theories.}
Our ep-specifications are closely related to essentially algebraic theories (see~\cite{Reichel1987} and~\cite{AdamekRosicky1994}, Section~3.D): in both cases the operations are partial and their domains are described by equations.
There are two differences.
In~\cite{AdamekRosicky1994}, the domains are described by equations between total operations, whereas the domain conditions of a hep-signature may involve partial operation symbols lower in the hierarchy, as in~\cite{Reichel1987}.
Moreover, we allow an arbitrary set of domain conditions for every operation symbol. Hence, even if every axiom has finitely many variables, an ep-specification need not be finitary in the usual sense, which requires the sets of domain conditions to be finite. For instance, a unary operation symbol $p$ may have as domain conditions the equations $u_{i}(v_{0})\eeq v_{0}$, for infinitely many total unary operation symbols $u_{i}$.
We neither need nor claim that the categories $\mathsf{EPAlg}(\spec{E})$ are locally finitely presentable or that they arise from finite-limit theories.

An interpretation of an essentially algebraic theory in another induces a functor between their categories of models.
A generalized hep-derivor differs from such an interpretation in three respects, all of which are used in this work.
First, it induces a functor on all partial algebras, and a model need not be sent to a model: a single construction gives the four functors of this work in finite dimensions, and whether they send categorial algebras to categorial algebras is then checked by means of the exactness criterion.
Second, the domain conditions of an operation symbol are evaluated in the transformed algebra, so the domain of a transformed operation can be smaller than that of the term that interprets it; in the discretization this restriction cannot be expressed by any term.
Third, composition is preserved only up to natural inclusions, which we record by means of a normal lax indexed category.
We do not claim that the results of Chapters~\ref{SGdUf}--\ref{SChangeSig} could not be obtained by other means; the point is that generalized hep-derivors provide a single mechanism for all of them. There is one exception: in dimension $\omega$ with a nonempty set of sorts, unification cannot be induced by a generalized hep-derivor, and it is obtained instead from the finite-dimensional unifications by a quotient of the disjoint union of the sets of cells of all dimensions.

\section{Structure of this work}

In Chapter~\ref{SPreliminaries} we collect the notions concerning many-sorted sets, free monoids, many-sorted algebras, terms and variables, derivors, many-sorted partial algebras, existence and quasi-existence equations, adjunctions and the Grothendieck construction that we use in the rest of this work.

In Chapter~\ref{SGenHepDerivors} we introduce ep-signatures, ep-algebras and their varieties, hep-signatures and the generalized hep-derivors between them, and we prove the general results described above.

In Chapters~\ref{SHomogeneous} and~\ref{SHeterogeneous} we recall homogeneous and heterogeneous higher-order categorial algebras, and we present them as the categories canonically associated with the ep-varieties of suitable hep-specifications.
In dimension $\omega$, heterogeneous $\omega$-categorial algebras also form such a category, whereas, for a nonempty set of sorts in the homogeneous presentation, the condition that every cell has a finite dimension cannot be expressed by $\QE$-equations, since the class of homogeneous $\omega$-categorial algebras is not closed under products of partial $\Sigma^{(\omega)}$-algebras; this is a property of the homogeneous presentation, not of the $\omega$-categorial algebras themselves.

In Chapter~\ref{SGdUf} we define, by means of generalized hep-derivors, the graduation functor, which distributes the cells of a homogeneous $n$-categorial algebra into sorts according to their dimension, and the unification functor, which keeps only the cells of top dimension of a heterogeneous $n$-categorial algebra.
We prove that they are mutually quasi-inverse equivalences of categories and, more precisely, that graduation is left adjoint to unification in an adjoint equivalence whose unit is the identity.
In Chapter~\ref{SDsTr} we define, in both presentations, the discretization functors, which regard an $m$-categorial algebra as an $n$-categorial algebra whose cells of dimension greater than $m$ are identities, and the truncation functors, which forget the cells of dimension greater than $m$.
We prove that discretization is left adjoint to truncation, that these functors compose as expected and that they are compatible with graduation and unification, so graduation, unification, discretization and truncation fit together in an algebraic square between two finite dimensions.
We end this chapter by defining the truncation and discretization functors in dimension $\omega$.

In Chapter~\ref{SProjLimits} we prove that the categories of homogeneous and of heterogeneous $\omega$-categorial algebras are the projective limits of the projective systems formed by the categories of $n$-categorial algebras and the truncation functors.
In Chapter~\ref{SFunctorsOmega} we obtain the graduation functor in dimension $\omega$ by the universal property of the projective limit and the unification functor in dimension $\omega$ by means of an inductive limit of sets of cells, and we prove that, also in dimension $\omega$, graduation and unification are mutually quasi-inverse equivalences, discretization is left adjoint to truncation, and these adjunctions form an algebraic square between a finite dimension and $\omega$.
Finally, in Chapter~\ref{SChangeSig} we show that every derivor $\mathbf{c}$ between many-sorted signatures determines generalized hep-derivors $\mathbf{c}^{(n)}$ and $\mathbf{c}^{\htg(n)}$ between the corresponding homogeneous and heterogeneous $n$-categorial signatures, that the induced functors send categorial algebras to categorial algebras and form split indexed categories, and that they commute with graduation, unification, discretization and truncation, both in finite dimensions and in dimension $\omega$. Hence the families of graduation, unification, discretization and truncation functors are natural transformations between split indexed categories, and they extend to the Grothendieck constructions.

The following diagram summarizes, for $m\leq n$ in $\mathbb{N}$, the functors of this work.
The horizontal arrows change the presentation and the vertical ones change the dimension.
Graduation and unification are mutually quasi-inverse equivalences (see Chapter~\ref{SGdUf}); discretization is left adjoint to truncation (see Chapter~\ref{SDsTr}); and the four adjunctions form an algebraic square (see Section~\ref{SAlgSquare}).
The same holds with $\omega$ in place of $n$ (see Chapter~\ref{SFunctorsOmega}).
Figure~\ref{FFunctorsSummary}, in the concluding remarks, gathers both squares, with the dimensions $m$, $n$ and $\omega$, together with the change of signature of Chapter~\ref{SChangeSig}.
$$
\xymatrix@C=70pt@R=45pt{
\mathsf{mCatAlg}(\sig{\Sigma})
    \ar@<1.1ex>[r]^-{\Gd^{(m)}_{\sig{\Sigma}}}
    \ar@{}[r]|-{\perp}
    \ar@<1.1ex>[d]^-{\Ds^{(n,m)}_{\sig{\Sigma}}}
    \ar@{}[d]|-{\vdash}
& \mathsf{mCat}^{\htg}\mathsf{Alg}(\sig{\Sigma})
    \ar@<1.1ex>[l]^-{\Uf^{(m)}_{\sig{\Sigma}}}
    \ar@<1.1ex>[d]^-{\Ds^{\htg(n,m)}_{\sig{\Sigma}}}
    \ar@{}[d]|-{\vdash}
\\
\mathsf{nCatAlg}(\sig{\Sigma})
    \ar@<1.1ex>[r]^-{\Gd^{(n)}_{\sig{\Sigma}}}
    \ar@{}[r]|-{\perp}
    \ar@<1.1ex>[u]^-{\Tr^{(m,n)}_{\sig{\Sigma}}}
& \mathsf{nCat}^{\htg}\mathsf{Alg}(\sig{\Sigma})
    \ar@<1.1ex>[l]^-{\Uf^{(n)}_{\sig{\Sigma}}}
    \ar@<1.1ex>[u]^-{\Tr^{\htg(m,n)}_{\sig{\Sigma}}}
}
$$

\textbf{How to read this work.}
Chapter~\ref{SPreliminaries} can be consulted when needed.
Chapter~\ref{SGenHepDerivors} contains the general theory on which all the constructions rest.
A reader mainly interested in higher-order categorial algebras may read Chapters~\ref{SHomogeneous} and~\ref{SHeterogeneous} first, and return to Chapter~\ref{SGenHepDerivors} for the transformed algebras, the criterion for them to be ep-algebras, the induced functors and the comparison between the functor induced by a composite and the composite of the induced functors.
Chapter~\ref{SDsTr} uses from Chapter~\ref{SGdUf} only two lemmas on cells, except for the section that relates the functors of the two chapters.
Chapters~\ref{SProjLimits} and~\ref{SFunctorsOmega} rely on the last section of Chapter~\ref{SDsTr}; the algebraic squares of Sections~\ref{SAlgSquare} and~\ref{SAlgSquareOmega} combine the results of Chapters~\ref{SGdUf} and~\ref{SDsTr}, and of Chapters~\ref{SGdUf}--\ref{SFunctorsOmega}, respectively; and Chapter~\ref{SChangeSig} uses all the previous ones.

\section{Prerequisites, notation and terminology}

The only prerequisite for reading this work is familiarity with category theory, universal algebra (in particular, many-sorted partial algebras), ordered sets and set theory.
Suitable references are, e.g., \cite{HerrlichStrecker1973, MacLane1998} for category theory, \cite{Bergman2015, Burmeister1982, Burmeister1986, Burmeister2002, BurrisSankappanavar1981, GoguenEtAl1978, Gratzer2008, Kurosh1967, Schmidt1966, Schmidt1970, Wechler1992} for universal algebra, \cite{Blyth2005, DaveyPriestley2002, Roman2008} for ordered sets, and \cite{Bourbaki1970, Enderton1977} for set theory.
The many-sorted algebraic conventions follow~\cite{ClimentCosme2019, ClimentCosme2020, ClimentCosme2023}.

\textbf{Terminology.} We use \emph{many-sorted} when the set of sorts is left implicit and \emph{$S$-sorted} when it is specified.
The adjective \emph{many-sorted} may be omitted when the signature or context fixes the sorts, in both the homogeneous and heterogeneous presentations.
The symbols of the heterogeneous presentation carry the superscript $\htg$, e.g., $\mathsf{nCat}^{\htg}$ and $\Sigma^{\htg(n)}$, whereas those of the homogeneous presentation carry no marker, e.g., $\mathsf{nCat}$ and $\Sigma^{(n)}$; accordingly, an $n$-categorial algebra without further qualification is a homogeneous one.

\textbf{Set theory.} We work in $\mathbf{ZFC}$ with a fixed Grothendieck universe $\boldsymbol{\mathcal{U}}$, see~\cite{MacLane1998}, pp.~21--24.
The elements of $\boldsymbol{\mathcal{U}}$ are called $\boldsymbol{\mathcal{U}}$-small sets and the subsets of $\boldsymbol{\mathcal{U}}$ are called $\boldsymbol{\mathcal{U}}$-large sets or classes.

\textbf{Natural numbers.} We write $\mathbb{N}=\omega$ for the set of natural numbers and identify every $n\in\mathbb{N}$ with $\{0,\ldots,n-1\}$.

\textbf{Functions and mappings.} For sets $A$ and $B$, a \emph{function} from $A$ to $B$ is a subset $F$ of $A\times B$ such that, for every $x\in A$, there exists a unique $y\in B$ with $(x,y)\in F$.
We denote by $\mathrm{Fnc}(A,B)$ the set of all such functions.
A \emph{mapping} from $A$ to $B$ is an ordered triple $f=(A,F,B)$ with $F\in\mathrm{Fnc}(A,B)$, denoted by $f\colon A\mor B$.
A \emph{partial mapping} from $A$ to $B$, denoted by $f\colon A\dmor B$, is a triple $f=(A,F,B)$ with $F\subseteq A\times B$ such that, for every $x\in A$, there is at most one $y\in B$ with $(x,y)\in F$.
We denote by $\mathrm{Hom}(A,B)$, or by $B^{A}$, and by $\mathrm{Hom}_{\mathrm{p}}(A,B)$ the sets of all mappings and of all partial mappings from $A$ to $B$, respectively.
For a partial mapping $f=(A,F,B)$, its \emph{domain of definition} is $\Dom(f)=\{x\in A\mid\exists\,y\in B\,((x,y)\in F)\}$; it equals $A$ when $f$ is a mapping.
When a mapping whose codomain is $\boldsymbol{\mathcal{U}}$, such as an $S$-sorted set, see Definition~\ref{DSortedSet}, has to be regarded as an element of $\boldsymbol{\mathcal{U}}$, we code it by its graph $F$, since the triple $(A,F,\boldsymbol{\mathcal{U}})$ is not an element of $\boldsymbol{\mathcal{U}}$.

\textbf{Categories.} We let $\mathsf{Set}$ stand for the category whose set of objects is $\boldsymbol{\mathcal{U}}$ and whose morphisms are the mappings between $\boldsymbol{\mathcal{U}}$-small sets, and $\mathsf{Cat}$ for the category of all $\boldsymbol{\mathcal{U}}$-categories, i.e., the category whose objects are the categories $\mathsf{C}$ such that $\mathrm{Ob}(\mathsf{C})\subseteq\boldsymbol{\mathcal{U}}$ and, for every $x,y\in\mathrm{Ob}(\mathsf{C})$, $\mathrm{Hom}(x,y)\in\boldsymbol{\mathcal{U}}$, and whose morphisms are the functors between them.
For a category $\mathsf{C}$, we denote by $\mathrm{Ob}(\mathsf{C})$ its class of objects and, for $x,y\in\mathrm{Ob}(\mathsf{C})$, by $\mathrm{Hom}_{\mathsf{C}}(x,y)$, or by $\mathrm{Hom}(x,y)$, the set of morphisms from $x$ to $y$.

\section{Formalization}\label{SFormalization}

The definitions and results of this work have been formalized in the proof assistant Lean~4, on top of its mathematical library Mathlib.
Its source code and its documentation are available at
\begin{center}
\url{https://github.com/mariovagomarzal/categorial-algebras}\quad and\quad \url{https://mariovagomarzal.github.io/categorial-algebras/docs/},
\end{center}
respectively.
It contains no unproved statements, and its theorems use no axioms beyond the three standard ones of Lean: propositional extensionality, the axiom of choice and the soundness of quotients.
The formalization was produced automatically by AI agents, see Section~\ref{SAIDisclosure}.

Every definition and result with a formal counterpart ends, flush right, with a mark \leansample{equivalent}, which links to the documentation of the corresponding Lean declaration; when there are several declarations, the mark is followed by small numerals that link to the others.
The colour of the mark states how the declarations relate to the statement in the text.
\begin{itemize}
\item \leansample{equivalent}\enspace\emph{Formally equivalent}: the declarations state the same as the text, up to the usual encodings of set-theoretic notions in type theory; for instance, many-sorted sets are families of types, words are lists and partial operations take partial values.
Many constructions are stated for an arbitrary type of dimensions, of which $n$, $n+1$ and $\omega$ are instances, and we also regard these as equivalent.
\item \leansample{weaker}\enspace\emph{Formally weaker}: the declarations follow from the statement in the text, but not conversely.
In most cases, an equality between algebras, functors or natural transformations becomes an isomorphism that is the identity on the elements.
In Lean, the carrier of a transformed algebra is the type of its cells, a subtype of the carrier of the given algebra, so two algebras with the same cells are isomorphic rather than equal, and so are the functors and the categories built from them.
A few statements also need the set of sorts to be nonempty.
\item \leansample{stronger}\enspace\emph{Formally stronger}: the declarations imply the statement in the text, for instance because they hold in greater generality.
\end{itemize}
When the relation is not an equivalence, the documentation of the declarations explains the difference briefly.

The parts of the text without a formal counterpart are some standard preliminary notions, which the formalization takes from Mathlib or does not need (products and coproducts of many-sorted sets, birestrictions, free monoids, the algebra of rows and split indexed categories), and a few remarks: on the relation between hep-algebras and the algebras and derivors of~\cite{ClimentSoliveres2010b}, on the size of the categories of categorial algebras, which universes of types take care of in Lean, on the extension to the Grothendieck constructions of the natural transformations of Chapter~\ref{SChangeSig} in finite dimensions, and on the behaviour under change of signature of the natural isomorphisms and counits of Chapters~\ref{SGdUf}--\ref{SFunctorsOmega}.

\section{Acknowledgements}
The second and third authors were supported by the PID2024-159495NB-I00 grant from the
Ministerio de Ciencia e Innovaci\'{o}n, Spain, and the CIAICO/2023/007 grant from the
Conselleria d'Educaci\'{o}, Universitats i Ocupaci\'{o}, Generalitat Valenciana.
The second author has held a Specially Appointed Professor position at Nantong University during the completion of this work.
The third author was supported by the CIACIF/2022/489 grant from the Conselleria d'Educaci\'{o}, Universitats i Ocupaci\'{o}, Generalitat Valenciana co-funded by the European Social Fund.
Special thanks to Asunci\'{o}n de Montesa for her categorical knowledge.

\section{AI disclosure}\label{SAIDisclosure}
The authors used Claude (Anthropic) to check proofs and to improve the presentation and the \LaTeX{} source.
The authors reviewed all its suggestions and take full responsibility for the content of this work.

The formalization described in Section~\ref{SFormalization} was produced by AI agents.
Its initial stages were carried out with a preliminary version of Stemma, a software infrastructure for mathematical research with AI agents that the second and fourth authors are developing, see~\cite{VagoCosme2026}.
It was then completed with Claude Opus~5.5 (Anthropic) and ChatGPT Astra~6 (OpenAI).

\chapter{Preliminaries}\label{SPreliminaries}

This chapter fixes the algebraic and categorical conventions used in the rest of this work.
Most of these notions are standard, and the results that we state without proof are well known or follow directly from the definitions.

The chapter proceeds in three stages.
Sections~\ref{SManySortedSets}--\ref{SDerivors} deal with total structures: many-sorted sets, words on a set of sorts, many-sorted algebras, terms, variables and substitutions, and the derivors between many-sorted signatures.
Sections~\ref{SManySortedPartialAlg} and~\ref{STermsEq} deal with partial structures: many-sorted partial algebras, their partial term operations, and existence and quasi-existence equations, for which we follow~\cite{Burmeister1986, Burmeister2002}.
Finally, Section~\ref{SAdjGrothendieck} recalls the categorical notions that we use, namely adjunctions, algebraic squares, together with the relation between their two natural isomorphisms, corestrictions of functors and the Grothendieck construction, both for split and for normal lax indexed categories.

\section{Many-sorted sets}\label{SManySortedSets}

\begin{definition}
\label{DSortedSet}\lean{equivalent}{SortedSet}
Let $S$ be a set of sorts in $\boldsymbol{\mathcal{U}}$.
An $S$-\emph{sorted set} is a mapping $A = (A_{s})_{s\in S}$ from $S$ to $\boldsymbol{\mathcal{U}}$.
If $A$ and $X$ are $S$-sorted sets, then we say that $X$ is an $S$-\emph{sorted subset} or, to abbreviate, a \emph{subset} of $A$, denoted by $X\subseteq A$, if, for every $s\in S$, $X_{s}\subseteq A_{s}$.
We denote by $\mathrm{Sub}(A)$ the set of all $S$-sorted sets $X$ such that $X\subseteq A$.
\end{definition}

\begin{definition}
\label{DSortedMapping}\lean{equivalent}{SortedSet.Map}
Let $S$ be a set of sorts and $A$ and $B$ $S$-sorted sets.
The \emph{cartesian product of} $A$ \emph{and} $B$, denoted by $A\times B$, is the $S$-sorted set $\left(A_{s}\times B_{s}\right)_{s\in S}$.
An $S$-\emph{sorted function from} $A$ \emph{to} $B$ is an $S$-sorted set $F\subseteq A\times B$ such that, for every $s\in S$, $F_{s}$ is a function from $A_{s}$ to $B_{s}$, and an $S$-\emph{sorted mapping from} $A$ \emph{to} $B$ is a triple $f = (A,F,B)$ in which $F$ is an $S$-sorted function from $A$ to $B$.
Thus an $S$-sorted mapping from $A$ to $B$ amounts to a family $f=(f_{s})_{s\in S}$ of mappings $f_{s}\colon A_{s}\mor B_{s}$.
We denote by $\Hom(A,B)$ or by $B^{A}$ the set of all $S$-sorted mappings from $A$ to $B$, and we write $f\colon A\mor B$ for $f\in\Hom(A,B)$.
Composition and identities are defined componentwise, and we denote by $\mathsf{Set}^{S}$ the category of $S$-sorted sets and $S$-sorted mappings.
For $f\colon A\mor B$, $X\subseteq A$ and $Y\subseteq B$, the \emph{direct image} $f[X]$, the \emph{inverse image} $f^{-1}[Y]$ and the \emph{image} $\mathrm{Im}(f)=f[A]$ are defined componentwise.
If $X\subseteq A$, then the \emph{restriction of} $f$ \emph{to} $X$, denoted by $f|_{X}$, is $f\circ \mathrm{in}^{X,A}$, where $\mathrm{in}^{X,A} = (\mathrm{in}^{X_{s},A_{s}})_{s\in S}$ is the canonical embedding of $X$ into $A$.
\end{definition}

Every mapping between sets of sorts induces a functor between the corresponding categories of many-sorted sets.

\begin{proposition}[\cite{Soliveres2000}, Proposition~1.3.5]
\label{PDeltaFunctor}\lean{equivalent}{SortedSet.delta}
Let $S$ and $T$ be sets and $\varphi \colon S \mor T$ a mapping.
We let $\Delta_{\varphi}$ stand for the assignment from $\mathsf{Set}^{T}$ to $\mathsf{Set}^{S}$ that sends a $T$-sorted set $A = (A_{t})_{t \in T}$ to the $S$-sorted set $A_{\varphi} = (A_{\varphi(s)})_{s \in S}$, and a $T$-sorted mapping $f \colon A \mor B$ to the $S$-sorted mapping $f_{\varphi} = (f_{\varphi(s)})_{s \in S}$.
Then $\Delta_{\varphi}$ is a functor from $\mathsf{Set}^{T}$ to $\mathsf{Set}^{S}$.
\end{proposition}

\begin{definition}\label{DSortedProd}
Let $S$ be a set of sorts, $I$ a set in $\boldsymbol{\mathcal{U}}$ and $(A^{i})_{i\in I}$ an $I$-indexed family of $S$-sorted sets.
Then the \emph{product} $\prod_{i\in I}A^{i}$ and the \emph{coproduct} $\coprod_{i\in I}A^{i}$ of $(A^{i})_{i\in I}$ are the $S$-sorted sets defined, for every $s\in S$, as
\begin{align*}
\textstyle
\left(\prod_{i\in I}A^{i}\right)_{s} &= \textstyle\left\{(a_{i})_{i\in I}\in\mathrm{Fnc}\left(I,\bigcup_{i\in I}A^{i}_{s}\right)\mid \forall\,i\in I\,\left(a_{i}\in A^{i}_{s}\right)\right\},
\\
\textstyle\left(\coprod_{i\in I}A^{i}\right)_{s} &= \textstyle\bigcup_{i\in I}(A^{i}_{s}\times\{i\}),
\end{align*}
with the canonical projections $\mathrm{pr}^{i}\colon\prod_{i\in I}A^{i}\mor A^{i}$ and the canonical injections $\mathrm{in}^{i}\colon A^{i}\mor\coprod_{i\in I}A^{i}$, which have the usual universal properties.
The remaining set-theoretic operations on $S$-sorted sets, such as $\amalg$ (binary coproduct), $\bigcup$, $\cup$, $\bigcap$, $\cap$, $-$ and $\complement_{A}$, are defined componentwise.
\end{definition}

\begin{proposition}\label{PImFactorization}\lean{equivalent}{SortedSet.range_subset_iff_exists_unique}
Let $S$ be a set of sorts, $A$ and $B$ $S$-sorted sets, $Y$ a subset of $B$, and $f$ an $S$-sorted mapping from $A$ to $B$.
Then $\mathrm{Im}(f)\subseteq Y$ if, and only if, there exists an $S$-sorted mapping $h$ from $A$ to $Y$ such that $\mathrm{in}^{Y,B}\circ h=f$.
In this case $h$ is unique, and we call it the \emph{corestriction of} $f$ \emph{to} $Y$ and denote it by $f|^{Y}$.
\end{proposition}

\begin{definition}\label{DBirestriction}
Let $S$ be a set of sorts, $f\colon A\mor B$ an $S$-sorted mapping, $X\subseteq A$ and $Y\subseteq B$ such that $f[X]\subseteq Y$.
Then we denote by $f|_{X}^{Y}$ the $S$-sorted mapping $(f|_{X})|^{Y}$ from $X$ to $Y$, which is well defined since $f[X]\subseteq Y$.
We call this $S$-sorted mapping the \emph{birestriction of} $f$ to $X$ and $Y$.
\end{definition}

\section{Free monoids}\label{SFreeMonoids}

We next recall the concept of free monoid on a set and several notions associated with it, which are used to construct the free algebra on an $S$-sorted set and to define the variables of the terms.

\begin{definition}\label{DFreeMonoid}
For a set $A$, the \emph{free monoid} $\mathbf{A}^{\star}=(A^{\star},\curlywedge,\lambda)$ consists of the finite words $\mathbf{a}\colon n\mor A$, with concatenation and the empty word $\lambda$.
We write $\mathbf{a}=(a_i)_{i\in\bb{\mathbf{a}}}$, where $\bb{\mathbf{a}}=n$ is its length.
The insertion $\eta^{A}\colon A\mor A^{\star}$ sends $a$ to the one-letter word $(a)$, and the pair $(\eta^{A},\mathbf{A}^{\star})$ has the following universal property: for every monoid $\mathbf{M}=(M,\cdot,e)$ and every mapping $f\colon A\mor M$, there exists a unique monoid homomorphism $f^{\sharp}\colon\mathbf{A}^{\star}\mor\mathbf{M}$ such that $f^{\sharp}\circ\eta^{A}=f$, namely, the one that sends $\lambda$ to $e$ and a word $\mathbf{a}$ with $\bb{\mathbf{a}}\geq 1$ to $f(a_{0})\cdot\ldots\cdot f(a_{\bb{\mathbf{a}}-1})$.
In particular, the length $\bb{\,\cdot\,}\colon\mathbf{A}^{\star}\mor(\mathbb{N},+,0)$ is the unique monoid homomorphism taking every letter to $1$, i.e., the monoid homomorphism for the constant mapping with value $1$.
\end{definition}

\begin{definition}\label{DSubword}\lean{equivalent}{Word.position}
For a word $\mathbf{a}\in A^{\star}$ and $a\in A$, put
\[
\bb{\mathbf{a}}_{a}=\bigl|\{i\in\bb{\mathbf{a}}\mid a_i=a\}\bigr|.
\]
The positions at which $a$ occurs are enumerated in increasing order as $(i_j)_{j\in\bb{\mathbf{a}}_{a}}$.
This enumeration is empty when $a$ does not occur in $\mathbf{a}$.
Only finitely many of the numbers $\bb{\mathbf{a}}_{a}$ are nonzero, and $\bb{\mathbf{a}}=\sum_{a\in A}\bb{\mathbf{a}}_{a}$.
\end{definition}

\begin{convention}
In what follows, for a set of sorts $S$, we denote by $\mathbf{s}$, i.e., a lower case bold type $s$, an arbitrary word on $S$, and by $s$ an arbitrary letter in $S$.
\end{convention}

For a mapping $\varphi\colon S\mor T$, we let $\varphi^{\star}$ stand for the unique homomorphism of monoids from $\mathbf{S}^{\star}$ to $\mathbf{T}^{\star}$ such that $\varphi^{\star}\circ\eta^{S}=\eta^{T}\circ\varphi$.
Thus, for every $\mathbf{s}\in S^{\star}$, $\varphi^{\star}(\mathbf{s})=(\varphi(s_{j}))_{j\in\bb{\mathbf{s}}}$ and, in particular, $\varphi^{\star}(\lambda)=\lambda$.

\section{Many-sorted algebras}\label{SManySortedAlg}

We next recall the notions of many-sorted signature, many-sorted algebra and homomorphism between many-sorted algebras.

\begin{definition}\label{DSig}\lean{equivalent}{Signature}
Let $S$ be a set of sorts.
An $S$-\emph{sorted signature} is a mapping $\Sigma$ from $S^{\star}\times S$ to $\boldsymbol{\mathcal{U}}$ which sends a pair $(\mathbf{s},s)\in S^{\star}\times S$ to the set $\Sigma_{\mathbf{s},s}$ of the \emph{formal operations} of \emph{arity} $\mathbf{s}$, \emph{sort} (or \emph{coarity}) $s$, and \emph{rank} (or \emph{biarity}) $(\mathbf{s},s)$.
A \emph{many-sorted signature} is an ordered pair $\sig{\Sigma}=(S,\Sigma)$, where $S$ is a set of sorts and $\Sigma$ an $S$-sorted signature.
\end{definition}

\begin{definition}\label{DAlg}\lean{equivalent}{Algebra}
Let $S$ be a set of sorts and $\Sigma$ an $S$-sorted signature.
The $S^{\star}\times S$-sorted set of the \emph{finitary operations on} an $S$-sorted set $A$ is $(\mathrm{Hom}(A_{\mathbf{s}},A_{s}))_{(\mathbf{s},s)\in S^{\star}\times S}$, where, for every $\mathbf{s}\in S^{\star}$, $A_{\mathbf{s}} = \prod_{j\in \bb{\mathbf{s}}}A_{s_{j}}$, with $\bb{\mathbf{s}}$ denoting the length of the word $\mathbf{s}$ (if $\mathbf{s} = \lambda$, then $A_{\lambda}$ is a final set).
A \emph{structure of} $\Sigma$-\emph{algebra on} an $S$-sorted set $A$ is a family $(F_{\mathbf{s},s})_{(\mathbf{s},s)\in S^{\star}\times S}$, denoted by $F$, where, for $(\mathbf{s},s)\in S^{\star}\times S$, $F_{\mathbf{s},s}$ is a mapping from $\Sigma_{\mathbf{s},s}$ to $\mathrm{Hom}(A_{\mathbf{s}},A_{s})$ (if $(\mathbf{s},s) = (\lambda,s)$ and $\sigma\in \Sigma_{\lambda,s}$, then $F_{\lambda,s}(\sigma)$ picks out an element of $A_{s}$).
A \emph{many-sorted} $\Sigma$-\emph{algebra} (or, to abbreviate, $\Sigma$-\emph{algebra}) is a pair $(A,F)$, denoted by $\mathbf{A}$, where $A$ is an $S$-sorted set and $F$ a structure of $\Sigma$-algebra on $A$.
For a pair $(\mathbf{s},s)\in S^{\star}\times S$ and a formal operation $\sigma\in \Sigma_{\mathbf{s},s}$, in order to simplify the notation, the operation $F_{\mathbf{s},s}(\sigma)$ from $A_{\mathbf{s}}$ to $A_{s}$ is written as $F_{\sigma}$.
In some cases, to avoid mistakes, we denote by $F^{\mathbf{A}}$ the structure of $\Sigma$-algebra on $A$, and, for $(\mathbf{s},s)\in S^{\star}\times S$ and $\sigma\in \Sigma_{\mathbf{s},s}$, by $F^{\mathbf{A}}_{\sigma}$, or simply by $\sigma^{\mathbf{A}}$, the corresponding operation.
Moreover, for $s\in S$ and $\sigma\in\Sigma_{\lambda,s}$, we usually denote by $\sigma^{\mathbf{A}}$ the value of the mapping $F^{\mathbf{A}}_{\sigma}\colon A_{\lambda}\mor A_{s}$ at the unique element in $A_{\lambda}$.

A $\Sigma$-\emph{homomorphism} (or, to abbreviate, \emph{homomorphism}) from $\mathbf{A}$ to $\mathbf{B}$, where $\mathbf{B} = (B,G)$, is a triple $(\mathbf{A},f,\mathbf{B})$, denoted by $f\colon \mathbf{A}\mor \mathbf{B}$, where $f$ is an $S$-sorted mapping from $A$ to $B$ such that, for every $(\mathbf{s},s)\in S^{\star}\times S$, every $\sigma\in \Sigma_{\mathbf{s},s}$, and every $(a_{j})_{j\in \bb{\mathbf{s}}}\in A_{\mathbf{s}}$, we have
$
f_{s}(\sigma^{\mathbf{A}}((a_{j})_{j\in \bb{\mathbf{s}}})) = \sigma^{\mathbf{B}}(f_{\mathbf{s}}((a_{j})_{j\in \bb{\mathbf{s}}})),
$
where $f_{\mathbf{s}}$ is the mapping $\prod_{j\in \bb{\mathbf{s}}}f_{s_{j}}$ from $A_{\mathbf{s}}$ to $B_{\mathbf{s}}$ that sends $(a_{j})_{j\in \bb{\mathbf{s}}}$ in $A_{\mathbf{s}}$ to $(f_{s_{j}}(a_{j}))_{j\in \bb{\mathbf{s}}}$ in $B_{\mathbf{s}}$.
For a $\Sigma$-algebra $\mathbf{A}$ the triple $(\mathbf{A},\id^{A},\mathbf{A})$, where $\id^{A}$ is the identity at $A$, which is a homomorphism from $\mathbf{A}$ to $\mathbf{A}$, is denoted by $\id^{\mathbf{A}}$ and called the \emph{identity} homomorphism at $\mathbf{A}$.
We denote by $\mathsf{Alg}(\Sigma)$ the category of $\Sigma$-algebras and homomorphisms and by $\mathrm{Alg}(\Sigma)$ the set of objects of $\mathsf{Alg}(\Sigma)$.
\end{definition}

\section{The free algebra on a many-sorted set}\label{SFreeAlg}

We next recall that the forgetful functor $\mathrm{G}_{\Sigma}$ from $\mathsf{Alg}(\Sigma)$ to $\mathsf{Set}^{S}$ has a left adjoint $\mathbf{T}_{\Sigma}$, which assigns to an $S$-sorted set $X$ the free $\Sigma$-algebra $\mathbf{T}_{\Sigma}(X)$ on $X$, the algebra of terms.
The construction is standard, see, e.g., \cite{GoguenEtAl1978, Wechler1992} and~\cite{Soliveres2000}, Section~2.6, and we only recall what we use.

\begin{definition}\label{DRowAlg}
Let $S$ be a set of sorts, $\Sigma$ an $S$-sorted signature and $X$ an $S$-sorted set.
The \emph{algebra of $\Sigma$-rows in} $X$, denoted by $\mathbf{W}_{\Sigma}(X)$, is defined as follows:
\begin{enumerate}
\item The underlying $S$-sorted set of $\mathbf{W}_{\Sigma}(X)$, written as $\mathrm{W}_{\Sigma}(X)$, is precisely the $S$-sorted set $((\coprod\Sigma\amalg\coprod X)^{\star})_{s\in S}$, i.e., the mapping from $S$ to $\boldsymbol{\mathcal{U}}$ constantly $(\coprod\Sigma\amalg\coprod X)^{\star}$, where $(\coprod\Sigma\amalg\coprod X)^{\star}$ is the set of all words on the set $\coprod\Sigma\amalg\coprod X$, i.e., on the set
$$
\textstyle
\left[\left(\bigcup_{(\mathbf{s},s)\in S^{\star}\times S}(\Sigma_{\mathbf{s},s}\times\{(\mathbf{s},s)\})\right)\times\{0\}\right]\cup\left[\left(\bigcup_{s\in S}(X_{s}\times\{s\})\right)\times\{1\}\right].
$$
\item For every $(\mathbf{s},s)\in S^{\star}\times S$ and every $\sigma\in\Sigma_{\mathbf{s},s}$, the structural operation $F_{\sigma}$ associated with $\sigma$ is the mapping from $\mathrm{W}_{\Sigma}(X)_{\mathbf{s}}$ to $\mathrm{W}_{\Sigma}(X)_{s}$ which sends $(P_{j})_{j\in\bb{\mathbf{s}}}\in\mathrm{W}_{\Sigma}(X)_{\mathbf{s}}$ to $(\sigma)\bconcat\concat_{j\in\bb{\mathbf{s}}}P_{j}\in\mathrm{W}_{\Sigma}(X)_{s}$, where, for every $(\mathbf{s},s)\in S^{\star}\times S$ and every $\sigma\in\Sigma_{\mathbf{s},s}$, $(\sigma)$ stands for $(((\sigma,(\mathbf{s},s)),0))$, which is the value at $\sigma$ of the canonical mapping from $\Sigma_{\mathbf{s},s}$ to $(\coprod\Sigma\amalg\coprod X)^{\star}$.
\end{enumerate}
\end{definition}

\begin{definition}\label{DTermAlg}\lean{equivalent}{Term.algebra}
The \emph{free} $\Sigma$-\emph{algebra on} $X$, denoted by $\mathbf{T}_{\Sigma}(X)$, is the $\Sigma$-algebra determined by the subalgebra of $\mathbf{W}_{\Sigma}(X)$ generated by $(\{(x)\mid x\in X_{s}\})_{s\in S}$, where, for every $s\in S$ and every $x\in X_{s}$, $(x)$ stands for $(((x,s),1))$, which is the value at $x$ of the canonical mapping from $X_{s}$ to $(\coprod\Sigma\amalg\coprod X)^{\star}$.
We denote by $\mathrm{T}_{\Sigma}(X)$ the underlying $S$-sorted set of $\mathbf{T}_{\Sigma}(X)$ and, for $s\in S$, we call the elements of $\mathrm{T}_{\Sigma}(X)_{s}$ \emph{terms of type} $s$ \emph{with variables in} $X$ or $(X,s)$-\emph{terms}.
\end{definition}

When no confusion can arise, we write $x$ for $(x)$ and $\sigma((P_{j})_{j\in\bb{\mathbf{s}}})$ for $\sigma^{\mathbf{T}_{\Sigma}(X)}((P_{j})_{j\in\bb{\mathbf{s}}})$, which is the word $(\sigma)\bconcat\concat_{j\in\bb{\mathbf{s}}}P_{j}$.
The terms have the unique readability property, which allows definitions by recursion and proofs by induction on terms.

\begin{proposition}\label{PTermChar}\lean{equivalent}{Term}
Let $S$ be a set of sorts, $\Sigma$ an $S$-sorted signature and $X$ an $S$-sorted set.
Then, for every $s\in S$, every term $P\in \mathrm{T}_{\Sigma}(X)_{s}$ is either a variable $x\in X_{s}$, or $\sigma^{\mathbf{T}_{\Sigma}(X)}((P_{j})_{j\in\bb{\mathbf{s}}})$ for a unique $\mathbf{s}\in S^{\star}$, a unique $\sigma\in\Sigma_{\mathbf{s},s}$ and a unique family $(P_{j})_{j\in\bb{\mathbf{s}}}\in\mathrm{T}_{\Sigma}(X)_{\mathbf{s}}$, and these possibilities are mutually exclusive.
\end{proposition}

The universal property of the free $\Sigma$-algebra on an $S$-sorted set $X$, stated in the following proposition, follows immediately from Proposition~\ref{PTermChar}.

\begin{proposition}\label{PTermAlgUniv}\lean{equivalent}{Term.exists_unique_lift}
Let $S$ be a set of sorts, $\Sigma$ an $S$-sorted signature and $X$ an $S$-sorted set.
Then the pair $(\eta^{X},\mathbf{T}_{\Sigma}(X))$, where $\eta^{X}\colon X\mor\mathrm{T}_{\Sigma}(X)$ is the \emph{insertion of the generators}, has the following universal property: for every $\Sigma$-algebra $\mathbf{A}$ and every $S$-sorted mapping $f\colon X\mor A$, there exists a unique homomorphism $f^{\sharp}\colon\mathbf{T}_{\Sigma}(X)\mor\mathbf{A}$ such that $f^{\sharp}\circ \eta^{X} = f$.
\end{proposition}
\begin{proof}
Define $f^{\sharp}_{s}(x)=f_s(x)$ on variables and
\[
f^{\sharp}_{s}\bigl(\sigma^{\mathbf{T}_{\Sigma}(X)}((P_j)_{j\in\bb{\mathbf{s}}})\bigr)
=\sigma^{\mathbf{A}}\bigl((f^{\sharp}_{s_j}(P_j))_{j\in\bb{\mathbf{s}}}\bigr)
\]
on compound terms, including nullary symbols.
Unique readability makes this well defined, and induction on terms proves uniqueness.
\end{proof}

Proposition~\ref{PTermAlgUniv} allows us to carry out definitions by algebraic recursion on terms, as we do throughout this work.

\begin{corollary}\label{CTermAlgAdj}\lean{equivalent}{Term.adjunction}
Let $S$ be a set of sorts and $\Sigma$ an $S$-sorted signature.
Then the functor $\mathbf{T}_{\Sigma}$, which sends an $S$-sorted set $X$ to $\mathbf{T}_{\Sigma}(X)$ and an $S$-sorted mapping $f$ from $X$ to $Y$ to $f^{@} (= (\eta^{Y}\circ f)^{\sharp})$, the unique homomorphism from $\mathbf{T}_{\Sigma}(X)$ to $\mathbf{T}_{\Sigma}(Y)$ such that $f^{@}\circ\eta^{X} = \eta^{Y}\circ f$, is left adjoint for the forgetful functor $\mathrm{G}_{\Sigma}$ from $\mathsf{Alg}(\Sigma)$ to $\mathsf{Set}^{S}$.
\end{corollary}

The following two consequences of Proposition~\ref{PTermAlgUniv} are used in the treatment of substitutions and derivors.

\begin{proposition}\label{PTermAlgUnivId}\lean{equivalent}{Term.lift_var}
Let $S$ be a set of sorts, $\Sigma$ an $S$-sorted signature and $X$ an $S$-sorted set.
Then $(\eta^{X})^{\sharp} = \id^{\mathbf{T}_{\Sigma}(X)}$.
\end{proposition}

\begin{proposition}\label{PTermAlgUnivComp}\lean{equivalent}{Term.lift_comp}
Let $S$ be a set of sorts, $\Sigma$ an $S$-sorted signature, $X$ an $S$-sorted set, $\mathbf{A}$ and $\mathbf{B}$ $\Sigma$-algebras, $f$ a mapping from $X$ to $A$ and $g$ a $\Sigma$-homomorphism from $\mathbf{A}$ to $\mathbf{B}$.
Then $(g \circ f)^{\sharp} = g \circ f^{\sharp}$.
\end{proposition}

\section{Variables and substitutions}\label{SVariables}

We next recall the $S$-sorted sets of variables that we associate to the words on $S$, and the identifications of the families of elements indexed by a word with the valuations of these variables.

\begin{definition}\label{DVariables}\lean{equivalent}{Word.Variables}
Let $S$ be a set of sorts.
For every word $\mathbf{s}\in S^{\star}$, we let $\vs{\mathbf{s}}$ stand for the $S$-sorted set defined, for every $s\in S$, as $(\vs{\mathbf{s}})_{s} = \mathbf{s}^{-1}[\{s\}]$, i.e., the set of the positions at which $s$ occurs in $\mathbf{s}$.

Let $V^{S} = (V^{S}_{s})_{s\in S}$ be an arbitrary, but fixed, $S$-sorted set of variables with countably infinitely many variables in each component, say $V^{S}_{s} = \{\,v^{s}_{n}\mid n\in \mathbb{N}\,\}$, for every $s\in S$.
Then, for every word $\mathbf{s}\in S^{\star}$, we let $V^{S}_{\mathbf{s}}$ stand for the finite subset of $V^{S}$ defined, for every $s\in S$, as
$$
V^{S}_{\mathbf{s},s} = \{\,v^{s}_{i}\in V^{S}_{s}\mid i\in \bb{\mathbf{s}}_{s}\,\}
$$
(observe that $V^{S}_{\mathbf{s},s}$ is empty for those sorts $s$ that do not occur in the word $\mathbf{s}$).
\end{definition}

\begin{proposition}\label{PVariablesIso}\lean{weaker}{Word.valuationEquiv,Word.valuation_equiv_natural}
Let $S$ be a set of sorts and $\mathbf{s}$ a word of $S^{\star}$.
Then the following statements hold:
\begin{enumerate}
\item there exists an isomorphism $\vartheta_{\mathbf{s}}$ from $\vs{\mathbf{s}}$ to $V^{S}_{\mathbf{s}}$, which, for every $s\in S$ and every $j\in\bb{\mathbf{s}}_{s}$, sends $i_{j}$ to $v^{s}_{j}$, where $(i_{j})_{j\in\bb{\mathbf{s}}_{s}}$ is the enumeration in ascending order of the occurrences of $s$ in $\mathbf{s}$, see Definition~\ref{DSubword}; and
\item for every $S$-sorted set $A$, there exists an isomorphism $\theta_{\mathbf{s},A}$ from $A_{\vs{\mathbf{s}}}=\mathrm{Hom}(\vs{\mathbf{s}},A)$ to $A_{\mathbf{s}}$, which sends $a$ to $(a_{s_{j}}(j))_{j\in\bb{\mathbf{s}}}$, and, for every $S$-sorted mapping $f\colon A\mor B$, $f_{\mathbf{s}}\circ\theta_{\mathbf{s},A}=\theta_{\mathbf{s},B}\circ f_{\vs{\mathbf{s}}}$, where $f_{\vs{\mathbf{s}}}$ sends $a\in A_{\vs{\mathbf{s}}}$ to $f\circ a\in B_{\vs{\mathbf{s}}}$.
\end{enumerate}
\end{proposition}

In virtue of the previous proposition, for every word $\mathbf{s}\in S^{\star}$ and every $S$-sorted set $A$, the sets $A_{\vs{\mathbf{s}}}$, $A_{\mathbf{s}}$ and $A_{V^{S}_{\mathbf{s}}}$ are naturally isomorphic.
In what follows, if there is no ambiguity, we do not distinguish notationally between the elements of $A_{\vs{\mathbf{s}}}$ and those of $A_{\mathbf{s}}$ or those of $A_{V^{S}_{\mathbf{s}}}$.
Also, for simplicity of notation, instead of using $V^{S}_{\mathbf{s}}$ we use $\vs{\mathbf{s}}$.
Accordingly, for $j\in\bb{\mathbf{s}}$, we let $v^{\mathbf{s}}_{j}$ stand for the variable $\vartheta_{\mathbf{s},s_{j}}(j)$ of $V^{S}_{\mathbf{s}}$ that corresponds to the $j$-th letter of $\mathbf{s}$, i.e., $v^{\mathbf{s}}_{j}=v^{s_{j}}_{i}$ if the $j$-th letter of $\mathbf{s}$ is the $i$-th occurrence, counting from $0$, of the sort $s_{j}$ in $\mathbf{s}$; and a family $\mathbf{a}=(a_{j})_{j\in\bb{\mathbf{s}}}\in A_{\mathbf{s}}$ is identified with the valuation of $\vs{\mathbf{s}}$ in $A$ that sends $v^{\mathbf{s}}_{j}$ to $a_{j}$, for every $j\in\bb{\mathbf{s}}$.

We next define the substitutions, which are the homomorphisms between free algebras determined by assigning terms to variables, see~\cite{ClimentCosme2020}, Subsection~3.2, for a description of substitutions in terms of the occurrences of the variables.

\begin{definition}[cf.~\cite{ClimentCosme2020}, Definition~3.4]\label{DSubst}\lean{equivalent}{Term.subst}
Let $S$ be a set of sorts, $\Sigma$ an $S$-sorted signature, $X$ and $Y$ $S$-sorted sets and $f\colon Y\mor\mathrm{T}_{\Sigma}(X)$ an $S$-sorted mapping.
The \emph{substitution determined by} $f$ is the unique homomorphism $f^{\sharp}\colon\mathbf{T}_{\Sigma}(Y)\mor\mathbf{T}_{\Sigma}(X)$ such that $f^{\sharp}\circ\eta^{Y}=f$, see Proposition~\ref{PTermAlgUniv}.
That is, for every $s\in S$ and every $P\in\mathrm{T}_{\Sigma}(Y)_{s}$, $f^{\sharp}_{s}(P)$ is the term of sort $s$ obtained from $P$ by substituting, for every $t\in S$ and every $y\in Y_{t}$, the term $f_{t}(y)$ for every occurrence of $y$ in $P$.

When $Y=\vs{\mathbf{t}}$, for a word $\mathbf{t}\in S^{\star}$, the mapping $f$ is identified, by Proposition~\ref{PVariablesIso}, with the family $\mathcal{Q}=(Q_{j})_{j\in\bb{\mathbf{t}}}\in\mathrm{T}_{\Sigma}(X)_{\mathbf{t}}$, where $Q_{j}=f(v^{\mathbf{t}}_{j})$, for every $j\in\bb{\mathbf{t}}$.
In this case we denote $f^{\sharp}$ by $(Q_{j})^{\sharp}_{j\in\bb{\mathbf{t}}}$, or by $\mathcal{Q}^{\sharp}$, and, for $P\in\mathrm{T}_{\Sigma}(\vs{\mathbf{t}})_{s}$, we call $(Q_{j})^{\sharp}_{j\in\bb{\mathbf{t}}}(P)$ the \emph{substitution of $Q_{j}$ for $v^{\mathbf{t}}_{j}$ in $P$, for every $j\in\bb{\mathbf{t}}$}.
\end{definition}

For substitutions determined by families of terms, Proposition~\ref{PTermAlgUnivComp} takes the following form.

\begin{corollary}\label{CSubstComp}\lean{equivalent}{Term.subst_family_comp}
Let $S$ be a set of sorts, $\Sigma$ an $S$-sorted signature, $Z$ an $S$-sorted set, $\mathbf{r},\mathbf{t}\in S^{\star}$, $\mathcal{Q}=(Q_{j})_{j\in\bb{\mathbf{t}}}\in\mathrm{T}_{\Sigma}(\vs{\mathbf{r}})_{\mathbf{t}}$ and $\mathcal{R}\in\mathrm{T}_{\Sigma}(Z)_{\mathbf{r}}$.
Then $\mathcal{R}^{\sharp}\circ\mathcal{Q}^{\sharp}=(\mathcal{R}^{\sharp}(Q_{j}))^{\sharp}_{j\in\bb{\mathbf{t}}}$.
\end{corollary}

\section{Derivors}\label{SDerivors}

In this section we recall the derivors between many-sorted signatures and the category $\mathsf{Sig}_{\mathfrak{d}}$ of many-sorted signatures and derivors, see~\cite{ClimentSoliveres2010b} and~\cite{Soliveres2000}, Section~3.5.
A derivor sends the sorts of a signature to sorts of another one, and its operation symbols to terms of the other one.

\begin{definition}\label{DDerivor}\lean{equivalent}{Derivor}
Let $\sig{\Sigma} = (S,\Sigma)$ and $\sig{\Lambda} = (T,\Lambda)$ be many-sorted signatures.
A \emph{derivor} from $\sig{\Sigma}$ to $\sig{\Lambda}$ is an ordered pair $\mathbf{c} = (\varphi,c)$, where $\varphi\colon S\mor T$ is a mapping and $c=(c_{\mathbf{s},s})_{(\mathbf{s},s)\in S^{\star}\times S}$ is a family of mappings
$$
c_{\mathbf{s},s}\colon \Sigma_{\mathbf{s},s}\mor\mathrm{T}_{\Lambda}(\vs{\varphi^{\star}(\mathbf{s})})_{\varphi(s)}.
$$
Thus, for every $\sigma\in \Sigma_{\mathbf{s},s}$, $c_{\mathbf{s},s}(\sigma)$ is a term of sort $\varphi(s)$ with variables in $\vs{\varphi^{\star}(\mathbf{s})}$.
We write $c(\sigma)$ instead of $c_{\mathbf{s},s}(\sigma)$ when this is unlikely to cause confusion.
\end{definition}

A derivor translates terms into terms, see~\cite{Soliveres2000}, Section~3.5.
To define this translation by the universal property of the free algebra, we first transform the free algebras of the target signature along a derivor.

\begin{definition}\label{DTermAlgTransformed}\lean{equivalent}{Derivor.transformed}
Let $\mathbf{c}=(\varphi,c)$ be a derivor from $\sig{\Sigma}=(S,\Sigma)$ to $\sig{\Lambda}=(T,\Lambda)$ and $X$ a $T$-sorted set.
We denote by $\mathbf{T}_{\Lambda}(X)^{\mathbf{c}}$ the $\Sigma$-algebra whose underlying $S$-sorted set is $\mathrm{T}_{\Lambda}(X)_{\varphi}$ and in which, for every $(\mathbf{s},s)\in S^{\star}\times S$ and every $\sigma\in\Sigma_{\mathbf{s},s}$, the operation $\sigma^{\mathbf{T}_{\Lambda}(X)^{\mathbf{c}}}$ sends $\mathcal{Q}=(Q_{j})_{j\in\bb{\mathbf{s}}}\in(\mathrm{T}_{\Lambda}(X)_{\varphi})_{\mathbf{s}}=\mathrm{T}_{\Lambda}(X)_{\varphi^{\star}(\mathbf{s})}$ to $\mathcal{Q}^{\sharp}(c(\sigma))\in\mathrm{T}_{\Lambda}(X)_{\varphi(s)}$, i.e., to the term $c(\sigma)$ with $Q_{j}$ substituted for $v^{\varphi^{\star}(\mathbf{s})}_{j}$, for every $j\in\bb{\mathbf{s}}$.
\end{definition}

The substitutions of the target signature are homomorphisms between the transformed algebras.

\begin{lemma}\label{LSubstTransformedHom}\lean{stronger}{Derivor.substHom}
Let $\mathbf{d}=(\psi,d)$ be a derivor from $\sig{\Lambda}=(T,\Lambda)$ to $\sig{\Omega}=(U,\Omega)$, $\mathbf{t}\in T^{\star}$, $Z$ a $U$-sorted set and $\mathcal{R}\in\mathrm{T}_{\Omega}(Z)_{\psi^{\star}(\mathbf{t})}$.
Then $(\mathcal{R}^{\sharp})_{\psi}$ is a $\Lambda$-homomorphism from $\mathbf{T}_{\Omega}(\vs{\psi^{\star}(\mathbf{t})})^{\mathbf{d}}$ to $\mathbf{T}_{\Omega}(Z)^{\mathbf{d}}$.
\end{lemma}
\begin{proof}
Let $(\mathbf{q},t)\in T^{\star}\times T$, $\tau\in\Lambda_{\mathbf{q},t}$ and $\mathcal{P}=(P_{i})_{i\in\bb{\mathbf{q}}}\in\mathrm{T}_{\Omega}(\vs{\psi^{\star}(\mathbf{t})})_{\psi^{\star}(\mathbf{q})}$.
By Definition~\ref{DTermAlgTransformed} and Corollary~\ref{CSubstComp},
\begin{align*}
\mathcal{R}^{\sharp}\left(\tau^{\mathbf{T}_{\Omega}(\vs{\psi^{\star}(\mathbf{t})})^{\mathbf{d}}}(\mathcal{P})\right)
&=\mathcal{R}^{\sharp}\left(\mathcal{P}^{\sharp}(d(\tau))\right)
\\&=\left(\mathcal{R}^{\sharp}(P_{i})\right)^{\sharp}_{i\in\bb{\mathbf{q}}}(d(\tau))
\\&=\tau^{\mathbf{T}_{\Omega}(Z)^{\mathbf{d}}}\left(\left(\mathcal{R}^{\sharp}(P_{i})\right)_{i\in\bb{\mathbf{q}}}\right).
\end{align*}
\end{proof}

We next define the translation of terms along a derivor.

\begin{definition}\label{DDerivorExtension}\lean{equivalent}{Derivor.translateWord}
Let $\mathbf{c}=(\varphi,c)$ be a derivor from $\sig{\Sigma}=(S,\Sigma)$ to $\sig{\Lambda}=(T,\Lambda)$ and $\mathbf{s}\in S^{\star}$.
We let $c^{\sharp}$ stand for the unique $\Sigma$-homomorphism from $\mathbf{T}_{\Sigma}(\vs{\mathbf{s}})$ to $\mathbf{T}_{\Lambda}(\vs{\varphi^{\star}(\mathbf{s})})^{\mathbf{c}}$ that sends $v^{\mathbf{s}}_{j}$ to $v^{\varphi^{\star}(\mathbf{s})}_{j}$, for every $j\in\bb{\mathbf{s}}$, see Proposition~\ref{PTermAlgUniv}.
Thus $c^{\sharp}$ is the $S$-sorted mapping from $\mathrm{T}_{\Sigma}(\vs{\mathbf{s}})$ to $\mathrm{T}_{\Lambda}(\vs{\varphi^{\star}(\mathbf{s})})_{\varphi}$ determined by the following conditions:
\begin{enumerate}
\item[(Der1)] $c^{\sharp}(v^{\mathbf{s}}_{j})=v^{\varphi^{\star}(\mathbf{s})}_{j}$, for every $j\in\bb{\mathbf{s}}$; and
\item[(Der2)] $c^{\sharp}(\sigma^{\mathbf{T}_{\Sigma}(\vs{\mathbf{s}})}((R_{j})_{j\in\bb{\mathbf{r}}}))=(c^{\sharp}(R_{j}))^{\sharp}_{j\in\bb{\mathbf{r}}}(c(\sigma))$, for every $(\mathbf{r},s')\in S^{\star}\times S$, every $\sigma\in\Sigma_{\mathbf{r},s'}$ and every $(R_{j})_{j\in\bb{\mathbf{r}}}\in\mathrm{T}_{\Sigma}(\vs{\mathbf{s}})_{\mathbf{r}}$.
\end{enumerate}
That is, $c^{\sharp}(\sigma^{\mathbf{T}_{\Sigma}(\vs{\mathbf{s}})}((R_{j})_{j\in\bb{\mathbf{r}}}))$ is the term $c(\sigma)$ with $c^{\sharp}(R_{j})$ substituted for $v^{\varphi^{\star}(\mathbf{r})}_{j}$.
In particular, by Proposition~\ref{PTermAlgUnivId}, $c^{\sharp}(\sigma^{\mathbf{T}_{\Sigma}(\vs{\mathbf{r}})}((v^{\mathbf{r}}_{j})_{j\in\bb{\mathbf{r}}}))=(v^{\varphi^{\star}(\mathbf{r})}_{j})^{\sharp}_{j\in\bb{\mathbf{r}}}(c(\sigma))=c(\sigma)$.
Equivalently, $c^{\sharp}$ is the underlying mapping of the extension of $c$ to a homomorphism of Hall algebras, see~\cite{Soliveres2000}, Section~3.5.
\end{definition}

We next define the composition and the identities of derivors.

\begin{definition}\label{DDerivorCompId}\lean{equivalent}{Derivor.comp,Derivor.identity}
Let $\mathbf{c}=(\varphi,c)\colon\sig{\Sigma}\mor\sig{\Lambda}$ and $\mathbf{d}=(\psi,d)\colon\sig{\Lambda}\mor\sig{\Omega}$ be derivors, where $\sig{\Omega}=(U,\Omega)$ is a many-sorted signature.
The \emph{composition} of $\mathbf{c}$ and $\mathbf{d}$ is the derivor $\mathbf{d}\circ \mathbf{c}=(\psi\circ\varphi,d^{\sharp}\circ c)$ from $\sig{\Sigma}$ to $\sig{\Omega}$.
The \emph{identity} at $\sig{\Sigma}=(S,\Sigma)$ is the derivor $(\id_{S},\eta^{\Sigma})$, where $\eta^{\Sigma}(\sigma)=\sigma^{\mathbf{T}_{\Sigma}(\vs{\mathbf{s}})}((v^{\mathbf{s}}_{j})_{j\in\bb{\mathbf{s}}})$, for every $(\mathbf{s},s)\in S^{\star}\times S$ and every $\sigma\in\Sigma_{\mathbf{s},s}$.
\end{definition}

In the following lemma we prove that the translation of terms along a derivor commutes with substitutions.

\begin{lemma}[cf.~\cite{Soliveres2000}, Section~3.5]\label{LDerivorSubst}\lean{equivalent}{Derivor.translate_word_subst_family}
Let $\mathbf{d}=(\psi,d)$ be a derivor from $\sig{\Lambda}=(T,\Lambda)$ to $\sig{\Omega}=(U,\Omega)$, $\mathbf{r},\mathbf{t}\in T^{\star}$, $t\in T$, $P\in\mathrm{T}_{\Lambda}(\vs{\mathbf{t}})_{t}$ and $(Q_{j})_{j\in\bb{\mathbf{t}}}\in\mathrm{T}_{\Lambda}(\vs{\mathbf{r}})_{\mathbf{t}}$.
Then the following equality holds
$$
d^{\sharp}\left((Q_{j})^{\sharp}_{j\in\bb{\mathbf{t}}}(P)\right)=\left(d^{\sharp}(Q_{j})\right)^{\sharp}_{j\in\bb{\mathbf{t}}}\left(d^{\sharp}(P)\right).
$$
\end{lemma}
\begin{proof}
Let us write $\mathcal{Q}$ for $(Q_{j})_{j\in\bb{\mathbf{t}}}$ and $d^{\sharp}(\mathcal{Q})$ for $(d^{\sharp}(Q_{j}))_{j\in\bb{\mathbf{t}}}$.
The mappings $d^{\sharp}\circ\mathcal{Q}^{\sharp}$ and $(d^{\sharp}(\mathcal{Q})^{\sharp})_{\psi}\circ d^{\sharp}$ are $\Lambda$-homomorphisms from $\mathbf{T}_{\Lambda}(\vs{\mathbf{t}})$ to $\mathbf{T}_{\Omega}(\vs{\psi^{\star}(\mathbf{r})})^{\mathbf{d}}$, by Definitions~\ref{DSubst} and~\ref{DDerivorExtension} and Lemma~\ref{LSubstTransformedHom}.
By~(Der1), both send $v^{\mathbf{t}}_{j}$ to $d^{\sharp}(Q_{j})$, for every $j\in\bb{\mathbf{t}}$.
Hence they are equal, by Proposition~\ref{PTermAlgUniv}.
\end{proof}

As a consequence, the translation of terms is compatible with the composition and the identities of derivors.

\begin{corollary}\label{CDerivorCatLaws}\lean{equivalent}{Derivor.translate_comp,Derivor.translate_word_identity,Derivor.translate_symbol}
Let $\mathbf{c}=(\varphi,c)\colon\sig{\Sigma}\mor\sig{\Lambda}$, $\mathbf{d}=(\psi,d)\colon\sig{\Lambda}\mor\sig{\Omega}$ and $\mathbf{e}=(\chi,e)\colon\sig{\Omega}\mor\sig{\Xi}$ be derivors.
Then the following statements hold:
\begin{enumerate}
\item $(e^{\sharp}\circ d)^{\sharp}=e^{\sharp}\circ d^{\sharp}$;
\item $(\eta^{\Lambda})^{\sharp}$ is the identity; and
\item $c^{\sharp}\circ\eta^{\Sigma}=c$.
\end{enumerate}
\end{corollary}
\begin{proof}
Let $\mathbf{r}\in T^{\star}$.

(1) By Lemma~\ref{LDerivorSubst} and Definition~\ref{DTermAlgTransformed}, $(e^{\sharp})_{\psi}$ is a $\Lambda$-homomorphism from $\mathbf{T}_{\Omega}(\vs{\psi^{\star}(\mathbf{r})})^{\mathbf{d}}$ to $\mathbf{T}_{\Xi}(\vs{\chi^{\star}(\psi^{\star}(\mathbf{r}))})^{\mathbf{e}\circ\mathbf{d}}$, since, for every $(\mathbf{q},t)\in T^{\star}\times T$, every $\tau\in\Lambda_{\mathbf{q},t}$ and every $\mathcal{P}\in\mathrm{T}_{\Omega}(\vs{\psi^{\star}(\mathbf{r})})_{\psi^{\star}(\mathbf{q})}$, we have that $e^{\sharp}(\mathcal{P}^{\sharp}(d(\tau)))=e^{\sharp}(\mathcal{P})^{\sharp}(e^{\sharp}(d(\tau)))$.
Hence the mappings $(e^{\sharp})_{\psi}\circ d^{\sharp}$ and $(e^{\sharp}\circ d)^{\sharp}$ are $\Lambda$-homomorphisms from $\mathbf{T}_{\Lambda}(\vs{\mathbf{r}})$ to $\mathbf{T}_{\Xi}(\vs{\chi^{\star}(\psi^{\star}(\mathbf{r}))})^{\mathbf{e}\circ\mathbf{d}}$ which, by~(Der1), send $v^{\mathbf{r}}_{j}$ to $v^{\chi^{\star}(\psi^{\star}(\mathbf{r}))}_{j}$, for every $j\in\bb{\mathbf{r}}$.
Therefore they are equal, by Proposition~\ref{PTermAlgUniv}.

(2) The $\Lambda$-algebra $\mathbf{T}_{\Lambda}(\vs{\mathbf{r}})^{(\id_{T},\eta^{\Lambda})}$ is $\mathbf{T}_{\Lambda}(\vs{\mathbf{r}})$, since, for every $(\mathbf{q},t)\in T^{\star}\times T$, every $\tau\in\Lambda_{\mathbf{q},t}$ and every $\mathcal{P}=(P_{i})_{i\in\bb{\mathbf{q}}}\in\mathrm{T}_{\Lambda}(\vs{\mathbf{r}})_{\mathbf{q}}$, the substitution $\mathcal{P}^{\sharp}$ is a homomorphism that sends $v^{\mathbf{q}}_{i}$ to $P_{i}$, so that $\mathcal{P}^{\sharp}(\eta^{\Lambda}(\tau))=\tau^{\mathbf{T}_{\Lambda}(\vs{\mathbf{r}})}(\mathcal{P})$.
Hence $(\eta^{\Lambda})^{\sharp}$ is the endomorphism of $\mathbf{T}_{\Lambda}(\vs{\mathbf{r}})$ that fixes the variables, i.e., the identity, by Proposition~\ref{PTermAlgUnivId}.

(3) It is the last statement of Definition~\ref{DDerivorExtension}, since $\eta^{\Sigma}(\sigma)=\sigma^{\mathbf{T}_{\Sigma}(\vs{\mathbf{s}})}((v^{\mathbf{s}}_{j})_{j\in\bb{\mathbf{s}}})$, for every $\sigma\in\Sigma_{\mathbf{s},s}$.
\end{proof}

From the above corollary we obtain that the derivors form a category.

\begin{corollary}\label{CSigdCat}\lean{equivalent}{instCategoryManySortedSignature}
The many-sorted signatures together with the derivors, the composition and the identities of Definition~\ref{DDerivorCompId} constitute a category, denoted by $\mathsf{Sig}_{\mathfrak{d}}$.
\end{corollary}
\begin{proof}
By Corollary~\ref{CDerivorCatLaws}~(1), the composition is associative, since $(e^{\sharp}\circ d)^{\sharp}\circ c=e^{\sharp}\circ(d^{\sharp}\circ c)$, and, by Corollary~\ref{CDerivorCatLaws}~(2) and~(3), the identities are neutral for it.
\end{proof}

\section{Many-sorted partial algebras}\label{SManySortedPartialAlg}

We next recall the notions of many-sorted partial algebra and homomorphism between many-sorted partial algebras.

\begin{definition}\label{DPartialAlg}\lean{equivalent}{PartialAlgebra}
Let $S$ be a set of sorts and $\Sigma$ an $S$-sorted signature.
A \emph{structure of partial} $\Sigma$-\emph{algebra on} an $S$-sorted set $A$ is a family $(F_{\mathbf{s},s})_{(\mathbf{s},s)\in S^{\star}\times S}$, denoted by $F$, where, for $(\mathbf{s},s)\in S^{\star}\times S$, $F_{\mathbf{s},s}$ is a mapping from $\Sigma_{\mathbf{s},s}$ to $\mathrm{Hom}_{\mathrm{p}}(A_{\mathbf{s}},A_{s})$.
If $\Dom(F_{\mathbf{s},s}(\sigma)) = A_{\mathbf{s}}$, then $F_{\mathbf{s},s}(\sigma)$ is a \emph{total operation} from $A_{\mathbf{s}}$ to $A_{s}$.
If $\Dom(F_{\mathbf{s},s}(\sigma)) = \varnothing$, then $F_{\mathbf{s},s}(\sigma)$ is called a \emph{discrete operation} or an \emph{empty (partial) operation} from $A_{\mathbf{s}}$ to $A_{s}$.
If $(\mathbf{s},s) = (\lambda,s)$, $\sigma\in \Sigma_{\lambda,s}$, and $\Dom(F_{\lambda,s}(\sigma))\neq\varnothing$, then $F_{\lambda,s}(\sigma)$ picks out an element of $A_{s}$.
Therefore a partial operation $F_{\lambda,s}(\sigma)$ from $A_{\lambda}$ to $A_{s}$, also called a \emph{partial constant} $F_{\lambda,s}(\sigma)$ from $A_{\lambda}$ to $A_{s}$, is either discrete or picks out exactly one element of $A_{s}$.
A \emph{many-sorted partial} $\Sigma$-\emph{algebra} (or, to abbreviate, \emph{partial $\Sigma$-algebra}) is a pair $(A,F)$, abbreviated to $\mathbf{A}$, where $A$ is an $S$-sorted set and $F$ a structure of partial $\Sigma$-algebra on $A$.
For partial algebras we use the notational conventions of Definition~\ref{DAlg}; in particular, we write $F_{\sigma}$, $F^{\mathbf{A}}_{\sigma}$ or $\sigma^{\mathbf{A}}$ for the partial operation $F_{\mathbf{s},s}(\sigma)$.

A $\Sigma$-\emph{homomorphism} (or, to abbreviate, \emph{homomorphism}) from $\mathbf{A}$ to $\mathbf{B}$, where $\mathbf{B} = (B,G)$, is a triple $(\mathbf{A},f,\mathbf{B})$, abbreviated to $f\colon \mathbf{A}\mor \mathbf{B}$, where $f$ is an $S$-sorted mapping from $A$ to $B$ such that, for every $(\mathbf{s},s)\in S^{\star}\times S$, every $\sigma\in \Sigma_{\mathbf{s},s}$, and every $(a_{j})_{j\in \bb{\mathbf{s}}}\in A_{\mathbf{s}}$, if $(a_{j})_{j\in \bb{\mathbf{s}}}\in\Dom(\sigma^{\mathbf{A}})$, then $(f_{s_{j}}(a_{j}))_{j\in \bb{\mathbf{s}}}\in\Dom(\sigma^{\mathbf{B}})$ and
$
f_{s}(\sigma^{\mathbf{A}}((a_{j})_{j\in \bb{\mathbf{s}}})) = \sigma^{\mathbf{B}}((f_{s_{j}}(a_{j}))_{j\in \bb{\mathbf{s}}}).
$
The \emph{identity} homomorphism $\id^{\mathbf{A}}$ at $\mathbf{A}$ is defined as for total algebras.
We denote by $\mathsf{PAlg}(\Sigma)$ the category of partial $\Sigma$-algebras and homomorphisms and by $\mathrm{PAlg}(\Sigma)$ the set of objects of $\mathsf{PAlg}(\Sigma)$.
\end{definition}

Let us note that homomorphisms preserve, but need not reflect, definedness.
The category $\mathsf{Alg}(\Sigma)$ is a full subcategory of the category $\mathsf{PAlg}(\Sigma)$.
In the following proposition we recall a characterization of the isomorphisms of partial algebras.

\begin{proposition}\label{PPartialAlgIso}\lean{equivalent}{PartialAlgebra.Hom.isIso_iff_inverse}
Let $S$ be a set of sorts, $\Sigma$ an $S$-sorted signature and $f$ a homomorphism from a partial $\Sigma$-algebra $\mathbf{A}$ to a partial $\Sigma$-algebra $\mathbf{B}$.
Then the following statements are equivalent:
\begin{enumerate}
\item $f$ is an isomorphism from $\mathbf{A}$ to $\mathbf{B}$, i.e., $f$ is a section and a retraction.
\item $f$ is bijective and $(\mathbf{B},f^{-1},\mathbf{A})$, abbreviated to $f^{-1}\colon \mathbf{B}\mor \mathbf{A}$, is also a homomorphism from $\mathbf{B}$ to $\mathbf{A}$.
\end{enumerate}
\end{proposition}

In particular, a bijective homomorphism $f\colon\mathbf{A}\mor\mathbf{B}$ is an isomorphism if, and only if, it \emph{reflects definedness}, i.e., for every operation symbol $\sigma\in\Sigma_{\mathbf{s},s}$ and every $(a_{j})_{j\in\bb{\mathbf{s}}}\in A_{\mathbf{s}}$, if $(f_{s_{j}}(a_{j}))_{j\in\bb{\mathbf{s}}}\in\Dom(\sigma^{\mathbf{B}})$, then $(a_{j})_{j\in\bb{\mathbf{s}}}\in\Dom(\sigma^{\mathbf{A}})$.

We next recall the product of a family of many-sorted partial algebras.

\begin{definition}\label{DPartialAlgProd}\lean{equivalent}{PartialAlgebra.product}
Let $S$ be a set of sorts, $\Sigma$ an $S$-sorted signature and $(\mathbf{A}^{i})_{i\in I}$ a family of partial $\Sigma$-algebras.
Then
\begin{enumerate}
\item The \emph{product} of $(\mathbf{A}^{i})_{i\in I}$, $\prod_{i\in I}\mathbf{A}^{i}$, is the partial $\Sigma$-algebra which has as $S$-sorted underlying set $\prod_{i\in I}A^{i}$, and where, for every $(\mathbf{s},s)\in S^{\star}\times S$ and every $\sigma\in \Sigma_{\mathbf{s},s}$, the structural operation $\sigma^{\prod_{i\in I}\mathbf{A}^{i}}$ is defined as follows:
      \begin{enumerate}
      \item $\Dom(\sigma^{\prod_{i\in I}\mathbf{A}^{i}})$ is equal to 
	  \[
		\textstyle
	  \{(a_{j})_{j\in \bb{\mathbf{s}}}\in
      (\prod_{i\in I}A^{i})_{\mathbf{s}}\mid \forall\,i\in I\,
      ((a_{j}(i))_{j\in \bb{\mathbf{s}}}\in \Dom(\sigma^{\mathbf{A}^{i}}))\}.\]
      \item
      $$
      \sigma^{\prod_{i\in I}\mathbf{A}^{i}}
      \nfunction
      {\Dom(\sigma^{\prod_{i\in I}\mathbf{A}^{i}})}
      {\prod_{i\in I}A^{i}_{s}}
      {(a_{j})_{j\in \bb{\mathbf{s}}}}
      {(\sigma^{\mathbf{A}^{i}}((a_{j}(i))_{j\in \bb{\mathbf{s}}})_{i\in I}}
      $$
      \end{enumerate}
\item For every $i\in I$, the $i$-th \emph{canonical projection} is the homomorphism from $\prod_{i\in I}\mathbf{A}^{i}$ to $\mathbf{A}^{i}$ determined by the $S$-sorted mapping $\pr^{i}$ which, for every $s\in S$, is defined as follows
      $$
      \pr^{i}_{s}
      \nfunction
      {\prod_{i\in I}A^{i}_{s}}
      {A^{i}_{s}}
      {(a_{l})_{l\in I}}
      {a_{i}}
      $$
\end{enumerate}
\end{definition}

\begin{proposition}\label{PPartialAlgProd}\lean{equivalent}{PartialAlgebra.productIsLimit}
Let $S$ be a set of sorts, $\Sigma$ an $S$-sorted signature and $(\mathbf{A}^{i})_{i\in I}$ a family of partial $\Sigma$-algebras.
Then the pair $(\prod_{i\in I}\mathbf{A}^{i},(\pr^{i})_{i\in I})$ is a product in $\mathsf{PAlg}(\Sigma)$.
\end{proposition}

\section{Terms and equations}\label{STermsEq}

We next recall the partial term operations of a partial algebra.
The value of a term at a valuation in a partial algebra can be defined by means of the Schmidt construction, see~\cite{Schmidt1970, ClimentCosme2023}; for the purposes of this work, the following equivalent recursive description suffices.

\begin{definition}\label{DTermOp}\lean{equivalent}{Term.partialEval}
Let $S$ be a set of sorts, $\Sigma$ an $S$-sorted signature and $X$ an $S$-sorted set.
For an $S$-sorted set $A$, we let $A_{X}=\mathrm{Hom}(X,A)$ stand for the set of the valuations of $X$ in $A$.
For a partial $\Sigma$-algebra $\mathbf{A}$, $s\in S$ and $P\in\mathrm{T}_{\Sigma}(X)_{s}$, the \emph{partial term operation} $P^{\mathbf{A}}$ induced by $P$ is the partial mapping from $A_{X}=\mathrm{Hom}(X,A)$ to $A_{s}$ defined by recursion on $P$ as follows:
\begin{enumerate}
\item if $P=x$, for a variable $x\in X_{s}$, then $P^{\mathbf{A}}$ is total and sends $a$ to $a_{s}(x)$; and
\item if $P=\sigma^{\mathbf{T}_{\Sigma}(X)}((P_{j})_{j\in\bb{\mathbf{s}}})$, for $\mathbf{s}\in S^{\star}$, $\sigma\in\Sigma_{\mathbf{s},s}$ and $(P_{j})_{j\in\bb{\mathbf{s}}}\in\mathrm{T}_{\Sigma}(X)_{\mathbf{s}}$, then $a\in\Dom(P^{\mathbf{A}})$ if, and only if, $a\in\Dom(P_{j}^{\mathbf{A}})$ for every $j\in\bb{\mathbf{s}}$ and $(P_{j}^{\mathbf{A}}(a))_{j\in\bb{\mathbf{s}}}\in\Dom(\sigma^{\mathbf{A}})$, in which case
$$
P^{\mathbf{A}}(a)=\sigma^{\mathbf{A}}\left((P_{j}^{\mathbf{A}}(a))_{j\in\bb{\mathbf{s}}}\right).
$$
\end{enumerate}
\end{definition}

We next define the notion of existence equation.

\begin{definition}\label{DExistenceEq}\lean{equivalent}{EEquation}
Let $S$ be a set of sorts and $\Sigma$ an $S$-sorted signature.
\begin{enumerate}
\item If $(X,s)\in\boldsymbol{\mathcal{U}}^{S}\times S$, then an \emph{existence equation} (briefly $\mathrm{E}$-\emph{equation} or $\mathrm{Eeqt}$) of type $(X,s)$ is a pair $(P,Q)$, also denoted by $P\eeq Q$, where $P$ and $Q$ are $\Sigma$-terms of type $(X,s)$, i.e., elements of $\mathrm{T}_{\Sigma}(X)_{s}$.
If an $\mathrm{E}$-equation of type $(X,s)$ is such that $X$ is $S$-finite, i.e., for every $s\in S$, $X_{s}$ is finite, then we call it $S$\emph{-finitary}, and if it is such that $X$ is finite, then we call it \emph{finitary}.
\item We denote by $\mathrm{Eeqt}(\Sigma)_{X}$ the $S$-sorted set $(\mathrm{T}_{\Sigma}(X)_{s}^{2})_{s\in S}$ of all $\mathrm{E}$-equations with variables in an $S$-sorted set $X$.
A \emph{set of $\mathrm{E}$-equations} with variables in $X$ is an $S$-sorted subset $\mathcal{E}\subseteq\mathrm{Eeqt}(\Sigma)_{X}$.
\end{enumerate}
\end{definition}

Substitutions and derivors act on sets of $\mathrm{E}$-equations componentwise.
For $S$-sorted sets $X$ and $Y$, an $S$-sorted mapping $f\colon Y\mor\mathrm{T}_{\Sigma}(X)$ and $\mathcal{E}\subseteq\mathrm{Eeqt}(\Sigma)_{Y}$, we let $f^{\sharp}[\mathcal{E}]\subseteq\mathrm{Eeqt}(\Sigma)_{X}$ stand for the $S$-sorted set $(\{(f^{\sharp}_{s}(P),f^{\sharp}_{s}(Q))\mid(P,Q)\in\mathcal{E}_{s}\})_{s\in S}$.
For a derivor $\mathbf{c}=(\varphi,c)$ from $\sig{\Sigma}=(S,\Sigma)$ to $\sig{\Lambda}=(T,\Lambda)$, $\mathbf{s}\in S^{\star}$ and $\mathcal{H}\subseteq\mathrm{Eeqt}(\Sigma)_{\vs{\mathbf{s}}}$, we let $c^{\sharp}[\mathcal{H}]\subseteq\mathrm{Eeqt}(\Lambda)_{\vs{\varphi^{\star}(\mathbf{s})}}$ stand for the set of $\mathrm{E}$-equations
$$
c^{\sharp}[\mathcal{H}]=\left\lbrace\left(c^{\sharp}(P),c^{\sharp}(Q)\right)\,\middle|\,(P,Q)\in\mathcal{H}\right\rbrace.
$$

For $(X,s)\in\boldsymbol{\mathcal{U}}^{S}\times S$, we next define the ternary relation of satisfiability $\cdot \models^{\Sigma}_{X,s}\!\!\cdot\, [\cdot]$ between a partial $\Sigma$-algebra, an $\mathrm{E}$-equation of type $(X,s)$, and a valuation of $X$ in the underlying $S$-sorted set of the partial $\Sigma$-algebra, and the derived binary relation of validity.

\begin{definition}\label{DExistenceEqSat}\lean{equivalent}{EEquation.Satisfies}
Let $S$ be a set of sorts, $\Sigma$ an $S$-sorted signature, $\mathbf{A}$ a partial $\Sigma$-algebra, $(X,s)\in\boldsymbol{\mathcal{U}}^{S}\times S$, $(P,Q)\in \mathrm{Eeqt}(\Sigma)_{X,s}$, $\mathcal{E}\subseteq \mathrm{Eeqt}(\Sigma)_{X}$, and $a\in A_{X} = \mathrm{Hom}(X,A)$.
\begin{enumerate}
\item We say that $\mathbf{A}$ \emph{satisfies the $\mathrm{E}$-equation} $(P,Q)$ \emph{with respect to the valuation} $a$, denoted by
      $$
      \mathbf{A}\models^{\Sigma}_{X,s}(P,Q) [a],
      $$
      if, and only if, $a\in\Dom(P^{\mathbf{A}})$, $a\in\Dom(Q^{\mathbf{A}})$, and $P^{\mathbf{A}}(a) = Q^{\mathbf{A}}(a)$, see Definition~\ref{DTermOp}.
\item We say that $\mathbf{A}$ \emph{satisfies} $\mathcal{E}$ \emph{with respect to the valuation} $a$, denoted by $\mathbf{A}\models^{\Sigma}_{X}\mathcal{E}[a]$, if, for every $s\in S$ and every $(P,Q)\in\mathcal{E}_{s}$, $\mathbf{A}\models^{\Sigma}_{X,s}(P,Q)[a]$.
\item We say that $(P,Q)$ is \emph{valid} in $\mathbf{A}$, denoted by $\mathbf{A}\models^{\Sigma}_{X,s}(P,Q)$, or by $\mathbf{A}\models(P,Q)$, when this is unlikely to cause confusion, if, for every $a\in A_{X}$, $\mathbf{A}\models^{\Sigma}_{X,s}(P,Q) [a]$.
\end{enumerate}
We omit the sub- and superscripts when this is unlikely to cause confusion.
\end{definition}

Let us note that $\mathbf{A}\models^{\Sigma}_{X,s}(P,P)[a]$ if, and only if, $a\in\Dom(P^{\mathbf{A}})$, i.e., the $\mathrm{E}$-equation $P\eeq P$ states that $P$ is defined.

For $(X,(\mathbf{s},s))$ in $\boldsymbol{\mathcal{U}}^{S}\times (S^{\star}\times S)$, we next define the notion of quasi-existence equation of type $(X,(\mathbf{s},s))$, its satisfaction and the varieties determined by sets of them.

\begin{definition}\label{DQEEq}\lean{equivalent}{QEEquation,QESpecification,IsQEVariety}
Let $S$ be a set of sorts and $\Sigma$ an $S$-sorted signature.
\begin{enumerate}
\item If $(X,(\mathbf{s},s))\in\boldsymbol{\mathcal{U}}^{S}\times (S^{\star}\times S)$, then a \emph{quasi-existence equation} (briefly $\QE$-\emph{equation} or $\mathrm{QEeqt}$) of type $(X,(\mathbf{s},s))$ is a pair
$$
\left(\left(
P_{j},Q_{j}
\right)_{j\in \bb{\mathbf{s}}},(P,Q)\right),
$$
also denoted by
$$
\textstyle
\bigwedge_{j\in\bb{\mathbf{s}}}(P_{j},Q_{j})\to (P,Q)\, \text{ or by }\,
\bigwedge_{j\in\bb{\mathbf{s}}}(P_{j}\eeq Q_{j})\to (P\eeq Q),
$$
where, for every $j\in \bb{\mathbf{s}}$, $(P_{j},Q_{j})\in\mathrm{T}_{\Sigma}(X)^{2}_{s_{j}}$ and $(P,Q)\in \mathrm{T}_{\Sigma}(X)^{2}_{s}$.
The pairs $(P_{j},Q_{j})$, for $j\in \bb{\mathbf{s}}$, are called the \emph{premisses} and $(P,Q)$ the \emph{conclusion} of the $\QE$-equation.
If a $\QE$-equation of type $(X,(\mathbf{s}, s))$ is such that $X$ is $S$-finite, then we call it $S$-\emph{finitary}, and if it is such that $X$ is finite, then we call it \emph{finitary}.
\item We denote by $\mathrm{QEeqt}(\Sigma)_{X}$ the $S^{\star}\times S$-sorted set of all $\QE$-equations with variables in $X$.
\item Let $\mathbf{A}$ be a partial $\Sigma$-algebra and $a\in A_{X}$.
We say that $\mathbf{A}$ \emph{satisfies} the $\QE$-equation $\bigwedge_{j\in\bb{\mathbf{s}}}(P_{j},Q_{j})\to (P,Q)$ \emph{with respect to the valuation} $a$, denoted by
      $$
      \textstyle
      \mathbf{A}\models^{\Sigma}_{X,\mathbf{s},s}\bigwedge_{j\in\bb{\mathbf{s}}}(P_{j},Q_{j})\to (P,Q) [a],
      $$
      if $\mathbf{A}\models^{\Sigma}_{X,s} (P,Q)[a]$ whenever $\mathbf{A}\models^{\Sigma}_{X,s_{j}} (P_{j},Q_{j})[a]$ for every $j\in \bb{\mathbf{s}}$.
We say that the $\QE$-equation \emph{holds} in $\mathbf{A}$, or that it is \emph{valid} in $\mathbf{A}$, if $\mathbf{A}$ satisfies it with respect to every valuation $a\in A_{X}$.
\item A \emph{$\QE$-specification} is a quadruple $(S,\Sigma,X,\mathcal{E})$ with $\mathcal{E}\subseteq\mathrm{QEeqt}(\Sigma)_{X}$.
We denote by $\mathrm{Mod}_{\Sigma}(\mathcal{E})$ the class of all partial $\Sigma$-algebras in which every $\QE$-equation of $\mathcal{E}$ holds, the \emph{models} of the specification.
We say that a class $\boldsymbol{\mathcal{K}}$ of partial $\Sigma$-algebras is a \emph{finitary} $\QE$-\emph{variety} if there exists a set $\mathcal{E}$ of finitary $\QE$-equations such that $\boldsymbol{\mathcal{K}} = \mathrm{Mod}_{\Sigma}(\mathcal{E})$.
For a $\QE$-specification $(S,\Sigma,X,\mathcal{E})$, we also denote $\mathrm{Mod}_{\Sigma}(\mathcal{E})$ by $\mathcal{V}(\mathcal{E})$ and call it the \emph{$\QE$-variety} determined by the specification.
\end{enumerate}
\end{definition}

Every $\QE$-variety is closed under products, see~\cite{Burmeister1986}; cf.~\cite{Soliveres2000}, Propositions~2.11.18~(3) and~2.11.22, for equations and total algebras.
In fact, if $(\mathbf{A}^{i})_{i\in I}$ is a family of partial $\Sigma$-algebras in which a $\QE$-equation holds and $a$ is a valuation in $\prod_{i\in I}\mathbf{A}^{i}$ that satisfies its premisses, then, by Definitions~\ref{DPartialAlgProd} and~\ref{DTermOp}, for every $i\in I$, $\pr^{i}\circ a$ satisfies the premisses in $\mathbf{A}^{i}$, hence also the conclusion, and, again by Definitions~\ref{DPartialAlgProd} and~\ref{DTermOp}, $a$ satisfies the conclusion in $\prod_{i\in I}\mathbf{A}^{i}$.

The following two facts on the partial term operations are proved by induction on terms.

\begin{remark}\label{RTermOpSubst}\lean{equivalent}{Term.partialEval_subst,EEquationSet.satisfies_subst}
Let $S$ be a set of sorts, $\Sigma$ an $S$-sorted signature, $\mathbf{A}$ a partial $\Sigma$-algebra, $X$ and $Y$ $S$-sorted sets, $f\colon Y\mor\mathrm{T}_{\Sigma}(X)$ an $S$-sorted mapping and $a\in A_{X}$ a valuation such that, for every $s\in S$ and every $y\in Y_{s}$, $a\in\Dom(f_{s}(y)^{\mathbf{A}})$.
Let $f^{\mathbf{A}}(a)\in A_{Y}$ be the valuation that sends $y\in Y_{s}$ to $f_{s}(y)^{\mathbf{A}}(a)$.
Then, for every $s\in S$ and every $R\in\mathrm{T}_{\Sigma}(Y)_{s}$, $a\in\Dom(f^{\sharp}_{s}(R)^{\mathbf{A}})$ if, and only if, $f^{\mathbf{A}}(a)\in\Dom(R^{\mathbf{A}})$, and then $f^{\sharp}_{s}(R)^{\mathbf{A}}(a)=R^{\mathbf{A}}(f^{\mathbf{A}}(a))$.
Consequently, for every $\mathcal{E}\subseteq\mathrm{Eeqt}(\Sigma)_{Y}$, $\mathbf{A}\models f^{\sharp}[\mathcal{E}][a]$ if, and only if, $\mathbf{A}\models\mathcal{E}[f^{\mathbf{A}}(a)]$.
\end{remark}

\begin{remark}\label{RTermOpHom}\lean{equivalent}{Term.partialEval_hom,EEquationSet.satisfies_hom}
Let $S$ be a set of sorts, $\Sigma$ an $S$-sorted signature, $h\colon\mathbf{A}\mor\mathbf{B}$ a homomorphism of partial $\Sigma$-algebras, $X$ an $S$-sorted set, $a\in A_{X}$ a valuation, $s\in S$ and $P\in\mathrm{T}_{\Sigma}(X)_{s}$.
If $a\in\Dom(P^{\mathbf{A}})$, then $h\circ a\in\Dom(P^{\mathbf{B}})$ and $P^{\mathbf{B}}(h\circ a)=h_{s}(P^{\mathbf{A}}(a))$.
Consequently, for every $\mathcal{E}\subseteq\mathrm{Eeqt}(\Sigma)_{X}$, if $\mathbf{A}\models\mathcal{E}[a]$, then $\mathbf{B}\models\mathcal{E}[h\circ a]$.
\end{remark}

\section[Adjunctions, indexed categories, Grothendieck construction]{Adjunctions, indexed categories and the Grothendieck construction}\label{SAdjGrothendieck}

We let ``$\ast$'' and ``$\circ$'' stand for the horizontal and vertical composition of natural transformations, respectively; we denote by $\Id_{\mathsf{X}}$ the identity functor at $\mathsf{X}$; and we denote by $\id_{G}$ the identity natural transformation at $G\colon \mathsf{A}\mor \mathsf{X}$.
Moreover, for categories $\mathsf{A}$ and $\mathsf{X}$, functors $F,G\colon\mathsf{A}\mor\mathsf{X}$ and a natural transformation $\theta$ from $F$ to $G$, we write $\theta=(\theta_{a})_{a\in\mathrm{Ob}(\mathsf{A})}$, where, for every object $a$ of $\mathsf{A}$, $\theta_{a}\colon F(a)\mor G(a)$ is the component of $\theta$ at $a$.

\begin{definition}\label{DAdj}\lean{equivalent}{CategoryTheory.Adjunction}
An \emph{adjunction} $F\dashv G$, from $\mathsf{X}$ to $\mathsf{A}$, consists of functors $F\colon\mathsf{X}\mor\mathsf{A}$ and $G\colon\mathsf{A}\mor\mathsf{X}$, a \emph{unit} $\eta\colon\Id_{\mathsf{X}}\Longrightarrow G\circ F$ and a \emph{counit} $\varepsilon\colon F\circ G\Longrightarrow\Id_{\mathsf{A}}$ satisfying
\begin{enumerate}
\item[(TI1)] $(G\ast\varepsilon)\circ(\eta\ast G)=\id_G$;
\item[(TI2)] $(\varepsilon\ast F)\circ(F\ast\eta)=\id_F$.
\end{enumerate}
We also denote it by $(G,F,\eta,\varepsilon)$; $F$ is the left adjoint and $G$ the right adjoint.
\end{definition}

Four adjunctions form an \emph{algebraic square} when there exist natural isomorphisms between the appropriate composites, as described in the following definition.

\begin{definition}[\cite{Soliveres2000}, Definition~4.3.60]\label{DAlgSquare}\lean{equivalent}{AlgebraicSquare}
Let $\mathsf{C}$, $\mathsf{C}'$, $\mathsf{D}$ and $\mathsf{D}'$ be four categories.
Moreover let
\begin{enumerate}
    \item
    $(G, F, \eta, \varepsilon)$ be an adjunction with $F\colon\mathsf{C}\mor\mathsf{D}$ left adjoint to $G\colon\mathsf{D}\mor\mathsf{C}$.

    \item
    $(G', F', \eta', \varepsilon')$ be an adjunction with $F'\colon\mathsf{C}'\mor\mathsf{D}'$ left adjoint to $G'\colon\mathsf{D}'\mor\mathsf{C}'$.

    \item
    $(K, J, \eta^{\mathsf{C}}, \varepsilon^{\mathsf{C}})$ be an adjunction with $J\colon\mathsf{C}\mor\mathsf{C}'$ left adjoint to $K\colon\mathsf{C}'\mor\mathsf{C}$.

    \item
    $(I, H, \eta^{\mathsf{D}}, \varepsilon^{\mathsf{D}})$ be an adjunction with $H\colon\mathsf{D}\mor\mathsf{D}'$ left adjoint to $I\colon\mathsf{D}'\mor\mathsf{D}$.
\end{enumerate}
We say that these adjunctions form an \emph{algebraic square} if there exists a pair of natural isomorphisms $\alpha$ from $F' \circ J$ to $H \circ F$ and $\beta$ from $G \circ I$ to $K \circ G'$.

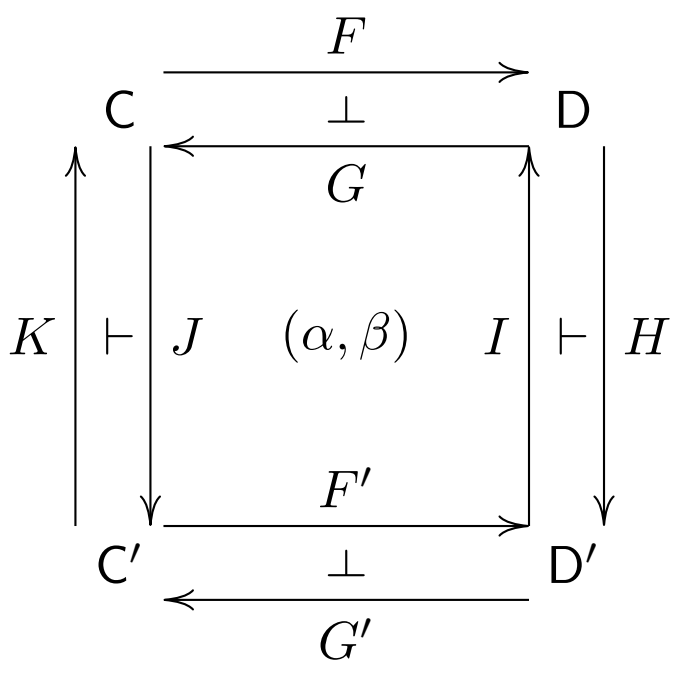
\begin{center}
\begin{tikzpicture}
    [ACliment/.style={-{To [angle'=45, length=5.75pt, width=4pt, round]}},
    RHACliment/.style={right hook-{To [angle'=45, length=5.75pt, width=4pt, round]}, font=\scriptsize},
    scale=1
    ]

    \node[] (C') at (0,0) {\textcolor{white}{$AAAA$}}; 
    \node[] () at (0,0) {$\mathsf{C}'$};
    \node[] (D') at (3,0) {\textcolor{white}{$AAAA$}}; 
    \node[] () at (3,0) {$\mathsf{D}'$};
    \node[] () at (1.5,0) {\rotatebox{90}{$\vdash$}};

    \draw[ACliment]  (C'.40) to node [above] {$F'$} (D'.140);
    \draw[ACliment]  (D'.220) to node [below] {$G'$} (C'.320);

    \node[] (C) at (0,3) {\textcolor{white}{$AAAA$}}; 
    \node[] () at (0,3) {$\mathsf{C}$};
    \node[] (D) at (3,3) {\textcolor{white}{$AAAA$}}; 
    \node[] () at (3,3) {$\mathsf{D}$};
    \node[] () at (1.5,3) {\rotatebox{90}{$\vdash$}};
    \node[] () at (0,1.5) {\rotatebox{0}{$\vdash$}};
    \node[] () at (3,1.5) {\rotatebox{0}{$\vdash$}};

    \node[] () at (1.5,1.5) {$(\alpha, \beta)$};

    \draw[ACliment]  (C.40) to node [above] {$F$} (D.140);
    \draw[ACliment]  (D.220) to node [below] {$G$} (C.320);

    \draw[ACliment]  (C.310) to node [right] {$J$} (C'.50);
    \draw[ACliment]  (C'.140) to node [left] {$K$} (C.220);

    \draw[ACliment]  (D.310) to node [right] {$H$} (D'.50);
    \draw[ACliment]  (D'.140) to node [left] {$I$} (D.220);
\end{tikzpicture}
\end{center}
\end{definition}

Informally, an algebraic square states that the two ways of going around the square agree up to natural isomorphism.
The two natural isomorphisms of an algebraic square are related as follows.

\begin{remark}\label{RAlgSquareMate}\lean{equivalent}{AlgebraicSquare.conjugateIso,AlgebraicSquare.conjugateIso_hom_app}
Given the four adjunctions of Definition~\ref{DAlgSquare}, the existence of $\alpha$ implies that of $\beta$.
In fact, $F'\circ J$ is left adjoint to $K\circ G'$ and $H\circ F$ is left adjoint to $G\circ I$, with the composite units and counits, and, if $F'\circ J\cong H\circ F$, then their right adjoints are isomorphic, by the uniqueness of right adjoints up to natural isomorphism, see~\cite{MacLane1998}, Section~IV.1, Corollary~1, p.~85, applied by duality.
Explicitly, if $\eta^{\circ}$ is the unit of $F'\circ J\dashv K\circ G'$ and $\varepsilon^{\circ}$ is the counit of $H\circ F\dashv G\circ I$, then
$$
\beta^{\alpha}=\left(K\circ G'\ast\varepsilon^{\circ}\right)\circ\left(K\circ G'\ast\alpha\ast G\circ I\right)\circ\left(\eta^{\circ}\ast G\circ I\right)
$$
is a natural isomorphism from $G\circ I$ to $K\circ G'$.
Thus the content of an algebraic square lies in the explicit description of $\alpha$ and $\beta$, and in the compatibility between them, rather than in their existence.
\end{remark}

\begin{definition}\label{DCorestriction}\lean{equivalent}{corestrict,restrict}
Let $F$ be a functor from a category $\mathsf{A}$ to a category $\mathsf{X}$, $\mathsf{A}'$ a full subcategory of $\mathsf{A}$ and $\mathsf{X}'$ a full subcategory of $\mathsf{X}$.
The \emph{restriction} of $F$ to $\mathsf{A}'$ is the functor $F\circ\mathrm{In}_{\mathsf{A}'}$ from $\mathsf{A}'$ to $\mathsf{X}$, where $\mathrm{In}_{\mathsf{A}'}$ is the inclusion functor of $\mathsf{A}'$ into $\mathsf{A}$.
We say that $F$ \emph{corestricts} to $\mathsf{X}'$ if, for every object $a$ of $\mathsf{A}$, $F(a)$ is an object of $\mathsf{X}'$.
In this case, since $\mathsf{X}'$ is full, there exists a unique functor $F'$ from $\mathsf{A}$ to $\mathsf{X}'$ such that $\mathrm{In}_{\mathsf{X}'}\circ F'=F$, the \emph{corestriction} of $F$ to $\mathsf{X}'$.
When this is unlikely to cause confusion, we denote the corestriction of the restriction of $F$ to $\mathsf{A}'$, if it exists, by $F$ itself, i.e., we write $F\colon\mathsf{A}'\mor\mathsf{X}'$.
\end{definition}

We next recall the split indexed categories and the Grothendieck construction, which assigns to a split indexed category a split fibration, see~\cite{Grothendieck1971} and, for split fibrations, \cite{BarrWells2012}, Sections~12.1 and~12.2.
We follow the notation of~\cite{JohnsonYau2021}, Chapter~10, but add the domain of $F$ as a superscript to indicate the contravariance of $F$.

\begin{definition}[cf.~\cite{Soliveres2000}, Propositions~3.2.1 and~3.5.7]\label{DSplitIndexed}
Let $\mathsf{I}$ be a category.
A \emph{split} $\mathsf{I}$-\emph{indexed category} is a contravariant functor $F$ from $\mathsf{I}$ to $\mathsf{Cat}$.
Given an object $i\in \mathrm{Ob}(\mathsf{I})$, we write $\mathsf{F}_{i}$ for the category $F(i)$, and given a morphism $\varphi\in \mathrm{Hom}_{\mathsf{I}}(i,j)$, we write $F_{\varphi}$ for the functor $F(\varphi)\colon \mathsf{F}_{j}\mor \mathsf{F}_{i}$.
A \emph{split indexed category} is an ordered pair $(\mathsf{I},F)$ in which $\mathsf{I}$ is a category and $F$ a split $\mathsf{I}$-indexed category.
\end{definition}

\begin{definition}[cf.~\cite{Soliveres2000}, Definitions~3.1.2, 3.2.2 and~3.5.8]\label{DGrothendieck}\lean{equivalent}{SplitIndexedCategory.Grothendieck}
Let $(\mathsf{I},F)$ be a split indexed category.
Then the \emph{Grothen\-dieck construction} at $(\mathsf{I},F)$ is $(\int^{\mathsf{I}}F,\pi_{F})$, where $\int^{\mathsf{I}}F$ is the category defined as follows:
\begin{enumerate}
\item $\mathrm{Ob}(\int^{\mathsf{I}}F) = \bigcup_{i\in\mathrm{Ob}(\mathsf{I})}(\{i\}\times\mathrm{Ob}(\mathsf{F}_{i}))$.
\item For every $(i,x), (j,y)\in \mathrm{Ob}(\int^{\mathsf{I}}F)$, $\mathrm{Hom}_{\int^{\mathsf{I}}F}((i,x),(j,y))$ is the set of all ordered pairs $(\varphi,f)$, where $\varphi\in \mathrm{Hom}_{\mathsf{I}}(i,j)$ and $f\in \mathrm{Hom}_{\mathsf{F}_{i}}(x,F_{\varphi}(y))$.
\item For every $(i,x)\in \mathrm{Ob}(\int^{\mathsf{I}}F)$, the identity morphism at $(i,x)$ is $(\id_{i},\id_{x})$.
\item For every $(i,x), (j,y), (k,z)\in \mathrm{Ob}(\int^{\mathsf{I}}F)$, every $(\varphi,f)\colon (i,x)\mor (j,y)$, and every $(\psi,g)\colon (j,y)\mor (k,z)$, the composite morphism $(\psi,g)\circ (\varphi,f)$ from $(i,x)$ to $(k,z)$ is
$$
(\psi,g)\circ (\varphi,f) = (\psi\circ \varphi,F_{\varphi}(g)\circ f)\colon (i,x)\mor (k,z).
$$
Here $F_{\varphi}F_{\psi}=F_{\psi\circ\varphi}$, so the composite has the required codomain.
\end{enumerate}
And $\pi_{F}$ is the functor from $\int^{\mathsf{I}}F$ to $\mathsf{I}$ that sends an object $(i,x)$ of $\int^{\mathsf{I}}F$ to $i$ and a morphism $(\varphi,f)\colon (i,x)\mor (j,y)$ of $\int^{\mathsf{I}}F$ to $\varphi\colon i\mor j$.
The functor $\pi_{F}$ is a split fibration, the \emph{canonical split fibration} from $\int^{\mathsf{I}}F$ to $\mathsf{I}$ (determined by $(\mathsf{I},F)$).
\end{definition}

The Grothendieck construction is functorial in the split indexed category, see~\cite{JohnsonYau2021}, Chapter~10.
We only use the following particular case.

\begin{remark}\label{RGrothendieckNat}\lean{equivalent}{SplitIndexedCategory.map}
Let $(\mathsf{I},F)$ and $(\mathsf{I},G)$ be split indexed categories and $\alpha=(\alpha_{i})_{i\in\mathrm{Ob}(\mathsf{I})}\colon F\Longrightarrow G$ a natural transformation, i.e., a family of functors $\alpha_{i}\colon\mathsf{F}_{i}\mor\mathsf{G}_{i}$ such that, for every morphism $\varphi\colon i\mor j$ of $\mathsf{I}$, $\alpha_{i}\circ F_{\varphi}=G_{\varphi}\circ\alpha_{j}$.
Then the assignment $\int^{\mathsf{I}}\alpha$ that sends an object $(i,x)$ of $\int^{\mathsf{I}}F$ to $(i,\alpha_{i}(x))$ and a morphism $(\varphi,f)\colon(i,x)\mor(j,y)$ to $(\varphi,\alpha_{i}(f))$ is a functor from $\int^{\mathsf{I}}F$ to $\int^{\mathsf{I}}G$ such that $\pi_{G}\circ\int^{\mathsf{I}}\alpha=\pi_{F}$.
In fact, $\alpha_{i}(f)$ is a morphism from $\alpha_{i}(x)$ to $\alpha_{i}(F_{\varphi}(y))=G_{\varphi}(\alpha_{j}(y))$; $\int^{\mathsf{I}}\alpha$ preserves identities, since every $\alpha_{i}$ does; and, for $(\psi,g)\colon(j,y)\mor(k,z)$, the following equality holds
$$
\alpha_{i}\left(F_{\varphi}(g)\circ f\right)=G_{\varphi}\left(\alpha_{j}(g)\right)\circ\alpha_{i}(f),
$$
so that $\int^{\mathsf{I}}\alpha$ preserves composites, see Definition~\ref{DGrothendieck}~(4).
\end{remark}

In Chapter~\ref{SGenHepDerivors} we will need the Grothendieck construction for assignments that preserve identities, but only preserve composites up to a coherent family of natural transformations, which need not be isomorphisms.
We next state the corresponding notion and the corresponding construction.

\begin{definition}\label{DLaxIndexed}\lean{equivalent}{NormalLaxIndexedCategory,LaxIndexedCategory}
Let $\mathsf{I}$ be a category.
A \emph{normal lax} $\mathsf{I}$-\emph{indexed category} is an ordered pair $(F,\kappa)$ in which
\begin{enumerate}
\item $F$ assigns to every object $i$ of $\mathsf{I}$ a category $\mathsf{F}_{i}$, and to every morphism $\varphi\colon i\mor j$ of $\mathsf{I}$ a functor $F_{\varphi}\colon\mathsf{F}_{j}\mor\mathsf{F}_{i}$, in such a way that, for every object $i$ of $\mathsf{I}$, $F_{\id_{i}}=\Id_{\mathsf{F}_{i}}$; and
\item $\kappa$ assigns to every pair of morphisms $\varphi\colon i\mor j$ and $\psi\colon j\mor k$ of $\mathsf{I}$ a natural transformation
$$
\kappa^{\varphi,\psi}=\left(\kappa^{\varphi,\psi}_{z}\right)_{z\in\mathrm{Ob}(\mathsf{F}_{k})}\colon F_{\varphi}\circ F_{\psi}\Longrightarrow F_{\psi\circ\varphi},
$$
\end{enumerate}
subject to the following conditions:
\begin{enumerate}
\item[(L1)] for every morphism $\varphi\colon i\mor j$ of $\mathsf{I}$, $\kappa^{\id_{i},\varphi}=\id_{F_{\varphi}}=\kappa^{\varphi,\id_{j}}$; and
\item[(L2)] for every triple of morphisms $\varphi\colon i\mor j$, $\psi\colon j\mor k$ and $\chi\colon k\mor l$ of $\mathsf{I}$, the following equality holds
$$
\kappa^{\varphi,\chi\circ\psi}\circ\left(F_{\varphi}\ast\kappa^{\psi,\chi}\right)
=
\kappa^{\psi\circ\varphi,\chi}\circ\left(\kappa^{\varphi,\psi}\ast F_{\chi}\right),
$$
i.e., for every object $w$ of $\mathsf{F}_{l}$, $\kappa^{\varphi,\chi\circ\psi}_{w}\circ F_{\varphi}(\kappa^{\psi,\chi}_{w})=\kappa^{\psi\circ\varphi,\chi}_{w}\circ\kappa^{\varphi,\psi}_{F_{\chi}(w)}$.
\end{enumerate}
Let us note that, if every $\kappa^{\varphi,\psi}$ is an identity, then $F$ is a split $\mathsf{I}$-indexed category, see Definition~\ref{DSplitIndexed}.
\end{definition}

\begin{definition}[cf.~\cite{Soliveres2000}, Definition~3.5.26, for pseudo-functors]\label{DLaxGrothendieck}\lean{equivalent}{NormalLaxIndexedCategory.Grothendieck}
Let $\mathsf{I}$ be a category and $(F,\kappa)$ a normal lax $\mathsf{I}$-indexed category.
Then $\int^{\mathsf{I}}(F,\kappa)$ stands for the following data:
\begin{enumerate}
\item the objects are the ordered pairs $(i,x)$ with $i\in\mathrm{Ob}(\mathsf{I})$ and $x\in\mathrm{Ob}(\mathsf{F}_{i})$;
\item the morphisms from $(i,x)$ to $(j,y)$ are the ordered pairs $(\varphi,f)$, where $\varphi\in \mathrm{Hom}_{\mathsf{I}}(i,j)$ and $f\in \mathrm{Hom}_{\mathsf{F}_{i}}(x,F_{\varphi}(y))$;
\item the identity morphism at $(i,x)$ is $(\id_{i},\id_{x})$; and
\item for every $(\varphi,f)\colon (i,x)\mor (j,y)$ and every $(\psi,g)\colon (j,y)\mor (k,z)$, the composite morphism $(\psi,g)\circ (\varphi,f)$ from $(i,x)$ to $(k,z)$ is
$$
(\psi,g)\circ (\varphi,f) = \left(\psi\circ \varphi,\kappa^{\varphi,\psi}_{z}\circ F_{\varphi}(g)\circ f\right).
$$
\end{enumerate}
Moreover, $\pi_{(F,\kappa)}$ stands for the assignment that sends an object $(i,x)$ to $i$ and a morphism $(\varphi,f)$ to $\varphi$.
\end{definition}

Let us note that, since $f\colon x\mor F_{\varphi}(y)$ and $g\colon y\mor F_{\psi}(z)$, the composite $\kappa^{\varphi,\psi}_{z}\circ F_{\varphi}(g)\circ f$ is a morphism of $\mathsf{F}_{i}$ from $x$ to $F_{\psi\circ\varphi}(z)$, so that the above definition is correct.
If every $\kappa^{\varphi,\psi}$ is an identity, then $\int^{\mathsf{I}}(F,\kappa)$ is the Grothendieck construction of Definition~\ref{DGrothendieck}.

\begin{proposition}\label{PLaxGrothendieck}\lean{equivalent}{NormalLaxIndexedCategory.Grothendieck.instCategory,NormalLaxIndexedCategory.Grothendieck.projection}
Let $\mathsf{I}$ be a category and $(F,\kappa)$ a normal lax $\mathsf{I}$-indexed category.
Then $\int^{\mathsf{I}}(F,\kappa)$ is a category and $\pi_{(F,\kappa)}$ is a functor from $\int^{\mathsf{I}}(F,\kappa)$ to $\mathsf{I}$.
\end{proposition}
\begin{proof}
The identity laws follow from normality: $F_{\id_i}=\Id_{\mathsf{F}_i}$ and the two unit comparisons in~(L1) are identities.

\textsf{(Associativity)} Let
$$
(\varphi,f)\colon(i,x)\mor(j,y),\qquad(\psi,g)\colon(j,y)\mor(k,z),\qquad(\chi,h)\colon(k,z)\mor(l,w)
$$
be morphisms of $\int^{\mathsf{I}}(F,\kappa)$.
Both $(\chi,h)\circ((\psi,g)\circ(\varphi,f))$ and $((\chi,h)\circ(\psi,g))\circ(\varphi,f)$ have $\chi\circ\psi\circ\varphi$ as first coordinate.
For the second coordinates, the following chain of equalities holds
\begin{flushleft}
$
\kappa^{\psi\circ\varphi,\chi}_{w}\circ F_{\psi\circ\varphi}(h)\circ\kappa^{\varphi,\psi}_{z}\circ F_{\varphi}(g)\circ f
$
\begin{align*}
&=
\kappa^{\psi\circ\varphi,\chi}_{w}\circ\kappa^{\varphi,\psi}_{F_{\chi}(w)}\circ F_{\varphi}(F_{\psi}(h))\circ F_{\varphi}(g)\circ f
\tag{1}
\\&=
\kappa^{\varphi,\chi\circ\psi}_{w}\circ F_{\varphi}(\kappa^{\psi,\chi}_{w})\circ F_{\varphi}(F_{\psi}(h))\circ F_{\varphi}(g)\circ f
\tag{2}
\\&=
\kappa^{\varphi,\chi\circ\psi}_{w}\circ F_{\varphi}\left(\kappa^{\psi,\chi}_{w}\circ F_{\psi}(h)\circ g\right)\circ f.
\tag{3}
\end{align*}
\end{flushleft}

The first equality follows from the naturality of $\kappa^{\varphi,\psi}$, applied to the morphism $h\colon z\mor F_{\chi}(w)$ of $\mathsf{F}_{k}$; the second equality follows from condition~(L2) of Definition~\ref{DLaxIndexed}; finally, the last equality follows from the fact that $F_{\varphi}$ is a functor.
By Definition~\ref{DLaxGrothendieck}, the left-hand side of the above chain is the second coordinate of $(\chi,h)\circ((\psi,g)\circ(\varphi,f))$ and its right-hand side is the second coordinate of $((\chi,h)\circ(\psi,g))\circ(\varphi,f)$.

\textsf{(The projection)} By Definition~\ref{DLaxGrothendieck}, $\pi_{(F,\kappa)}$ preserves identities and composites.

This completes the proof.
\end{proof}

We call $(\int^{\mathsf{I}}(F,\kappa),\pi_{(F,\kappa)})$ the \emph{Grothendieck construction} at $(\mathsf{I},(F,\kappa))$.

\chapter{Generalized hep-derivors}\label{SGenHepDerivors}

In this chapter we extend the derivors in two directions, as required by the functors between higher-order categorial algebras.
First, the signatures are hierarchical equationally partial signatures, or hep-signatures, i.e., equationally partial signatures in the sense of~\cite{Reichel1987}, in which every operation symbol carries a set of existence equations that describes its domain of definition, with the additional requirement that these domain conditions are organized hierarchically.
Second, a generalized hep-derivor sends every sort $s$ not only to a sort $\varphi(s)$, but also to a set of existence equations in one variable, the \emph{domain} of $s$, which selects, among the elements of sort $\varphi(s)$, those that interpret $s$.
As with derivors, every operation symbol $\sigma$ is sent to a term $c(\sigma)$; the domain of the operation that interprets $\sigma$ in the transformed algebras is determined by the domain conditions of $\sigma$, evaluated in the transformed algebra itself.
Sorts interpreted by the elements that satisfy a formula appear also in the general derivors of~\cite{LucasGutierrez2018}, which are morphisms between order-sorted signatures with predicates.

The chapter proceeds as follows.
In Section~\ref{SEpSig} we recall the ep-signatures and the ep-algebras, and we introduce the hierarchical ep-signatures, or hep-signatures, together with their specifications, their ep-varieties and the categories canonically associated with them.
In Section~\ref{SGenHepDerivorsDef} we define the generalized hep-derivors and the algebras transformed along them, and we prove that every generalized hep-derivor induces a functor between the categories of partial algebras.
In Section~\ref{SHepSigCat} we define the composition of generalized hep-derivors and we prove that the hep-signatures and the generalized hep-derivors form a category, $\mathsf{HSig}_{\gd}$.
In Section~\ref{SHepCompFunctors} we compare the functor induced by a composition with the composite of the induced functors.
Finally, in Section~\ref{SHepAlgCat} we assemble the ep-algebras of all the hep-signatures into a category, $\mathsf{HAlg}_{\gd}$.
We use the notions concerning terms, substitutions and derivors recalled in Sections~\ref{SVariables}, \ref{SDerivors} and~\ref{STermsEq}.

\section{ep-signatures, hep-signatures and their specifications}\label{SEpSig}

We begin by defining the notion of ep-signature.

\begin{definition}\label{DEpSig}\lean{equivalent}{EPSignature}
An \emph{equationally partial signature}, or \emph{ep-signature}, is an ordered triple $\sig{\Sigma}=(S,\Sigma,\mathrm{def}^{(S,\Sigma)})$ in which $(S,\Sigma)$ is a many-sorted signature and $\mathrm{def}^{(S,\Sigma)}$ assigns to every $(\mathbf{s},s)\in S^{\star}\times S$ and every $\sigma\in\Sigma_{\mathbf{s},s}$ a set
$$
\mathrm{def}^{(S,\Sigma)}(\sigma)\subseteq\mathrm{Eeqt}(\Sigma)_{\vs{\mathbf{s}}}
$$
of $\mathrm{E}$-equations with variables in $\vs{\mathbf{s}}$.
We call $\mathrm{def}^{(S,\Sigma)}(\sigma)$ the set of \emph{domain conditions} of $\sigma$.
\end{definition}

In the following definition we introduce the partial algebras whose operations are defined exactly on the families that satisfy the corresponding domain conditions.

\begin{definition}\label{DEpAlg}\lean{equivalent}{EPSignature.IsEPAlgebra}
Let $\sig{\Sigma}=(S,\Sigma,\mathrm{def}^{(S,\Sigma)})$ be an ep-signature and $\mathbf{A}$ a partial $\Sigma$-algebra, see Definition~\ref{DPartialAlg}.
For every $(\mathbf{s},s)\in S^{\star}\times S$ and every $\sigma\in\Sigma_{\mathbf{s},s}$, we let $A_{\mathrm{def}^{(S,\Sigma)}(\sigma)}$ stand for the subset of $A_{\mathbf{s}}$ defined as
$$
A_{\mathrm{def}^{(S,\Sigma)}(\sigma)}=\left\lbrace\mathbf{a}\in A_{\mathbf{s}}\,\middle|\,\mathbf{A}\models\mathrm{def}^{(S,\Sigma)}(\sigma)[\mathbf{a}]\right\rbrace.
$$
We say that $\mathbf{A}$ is an \emph{ep-algebra of type $\sig{\Sigma}$}, or a \emph{$\sig{\Sigma}$-algebra}, if, for every $(\mathbf{s},s)\in S^{\star}\times S$ and every $\sigma\in\Sigma_{\mathbf{s},s}$, the following equality holds
$$
\Dom\left(\sigma^{\mathbf{A}}\right)=A_{\mathrm{def}^{(S,\Sigma)}(\sigma)}.
$$
We denote by $\mathsf{EPAlg}(\sig{\Sigma})$ the full subcategory of $\mathsf{PAlg}(\Sigma)$, the category of partial $\Sigma$-algebras and their homomorphisms, see Definition~\ref{DPartialAlg}, whose objects are the ep-algebras of type $\sig{\Sigma}$.
\end{definition}

\begin{remark}\label{REpAlg}\lean{weaker}{EPSignature.ofSignature_isEPAlgebra_iff}
Let us note that the domain conditions of an operation symbol $\sigma$ are evaluated in $\mathbf{A}$ itself, so they may involve the other partial operations of $\mathbf{A}$, and even $\sigma$ itself.
Let us also note that a many-sorted signature $(S,\Sigma)$ is identified with the ep-signature $(S,\Sigma,\varnothing)$, in which every operation symbol has an empty set of domain conditions.
Its ep-algebras are exactly the (total) $\Sigma$-algebras, so that
$$
\mathsf{EPAlg}(S,\Sigma,\varnothing)=\mathsf{Alg}(\Sigma).
$$
\end{remark}

We also need the following weakening of the notion of ep-algebra, in which an operation may be defined only on families that satisfy the corresponding domain conditions, but need not be defined on all of them.

\begin{definition}\label{DWeakEpAlg}\lean{equivalent}{EPSignature.IsWeakEPAlgebra}
Let $\sig{\Sigma}=(S,\Sigma,\mathrm{def}^{(S,\Sigma)})$ be an ep-signature and $\mathbf{A}$ a partial $\Sigma$-algebra.
We say that $\mathbf{A}$ is a \emph{weak ep-algebra of type $\sig{\Sigma}$}, or a \emph{weak $\sig{\Sigma}$-algebra}, if, for every $(\mathbf{s},s)\in S^{\star}\times S$ and every $\sigma\in\Sigma_{\mathbf{s},s}$, the following inclusion holds
$$
\Dom\left(\sigma^{\mathbf{A}}\right)\subseteq A_{\mathrm{def}^{(S,\Sigma)}(\sigma)}.
$$
We denote by $\mathsf{WEPAlg}(\sig{\Sigma})$ the full subcategory of $\mathsf{PAlg}(\Sigma)$ whose objects are the weak ep-algebras of type $\sig{\Sigma}$.
\end{definition}

Let us note that every $\sig{\Sigma}$-algebra is a weak $\sig{\Sigma}$-algebra, so that $\mathsf{EPAlg}(\sig{\Sigma})$ is a full subcategory of $\mathsf{WEPAlg}(\sig{\Sigma})$.

In order to define the algebras transformed along a generalized hep-derivor by recursion, we need the domain conditions to be organized hierarchically.
To this end, we first assign to every term the set of the operation symbols that occur in it.

\begin{definition}\label{DOpSymbols}\lean{equivalent}{Term.symbols}
Let $\Sigma$ be an $S$-sorted signature and $X$ an $S$-sorted set.
We let $\mathrm{Op}$ stand for the mapping that assigns to every $s\in S$ and every term $P\in\mathrm{T}_{\Sigma}(X)_{s}$ the set $\mathrm{Op}(P)$ of the \emph{operation symbols that occur in $P$}, defined by recursion on $P$ as follows:
\begin{enumerate}
\item if $P=x$, for a variable $x\in X_{s}$, then $\mathrm{Op}(P)=\varnothing$; and
\item if $P=\sigma^{\mathbf{T}_{\Sigma}(X)}((P_{j})_{j\in\bb{\mathbf{s}}})$, for $\mathbf{s}\in S^{\star}$, $\sigma\in\Sigma_{\mathbf{s},s}$ and $(P_{j})_{j\in\bb{\mathbf{s}}}\in\mathrm{T}_{\Sigma}(X)_{\mathbf{s}}$, then
$$
\mathrm{Op}(P)=\{\sigma\}\cup\textstyle\bigcup_{j\in\bb{\mathbf{s}}}\mathrm{Op}(P_{j}).
$$
\end{enumerate}
Moreover, for every $s\in S$ and every $\mathrm{E}$-equation $(P,Q)\in\mathrm{Eeqt}(\Sigma)_{X,s}$, we let $\mathrm{Op}(P,Q)$ stand for $\mathrm{Op}(P)\cup\mathrm{Op}(Q)$ and, for every set of $\mathrm{E}$-equations $\mathcal{H}\subseteq\mathrm{Eeqt}(\Sigma)_{X}$, we let $\mathrm{Op}(\mathcal{H})$ stand for $\bigcup_{s\in S}\bigcup_{(P,Q)\in\mathcal{H}_{s}}\mathrm{Op}(P,Q)$.
\end{definition}

We are now in a position to introduce the hierarchical ep-signatures.

\begin{definition}\label{DHepSig}\lean{equivalent}{HepSignature}
Let $\sig{\Sigma}=(S,\Sigma,\mathrm{def}^{(S,\Sigma)})$ be an ep-signature.
We let $\lhd$ stand for the binary relation on the set of the operation symbols of $\Sigma$ defined, for every pair of operation symbols $\tau$ and $\sigma$ of $\Sigma$, as follows:
\begin{align*}
\tau&\lhd\sigma
&\mbox{if, and only if,}&&
&\tau\in\mathrm{Op}\left(\mathrm{def}^{(S,\Sigma)}(\sigma)\right),
\end{align*}
see Definition~\ref{DOpSymbols}, i.e., if, and only if, $\tau$ occurs in some equation of $\mathrm{def}^{(S,\Sigma)}(\sigma)$.
We say that $\sig{\Sigma}$ is \emph{hierarchical}, or that it is a \emph{hierarchical equationally partial signature}, or a \emph{hep-signature}, if $\lhd$ is well-founded.
\end{definition}

Let us note that, in a hep-signature, no operation symbol occurs in its own domain conditions, since, otherwise, $\sigma\lhd\sigma$ would give rise to an infinite descending chain.

Since $\lhd$ is well-founded, the operation symbols of a hep-signature admit proofs by induction and definitions by recursion along $\lhd$.

\begin{proposition}[Principle of Proof by Well-Founded Induction]\label{PWellFoundedInduction}\lean{equivalent}{HepSignature.symbol_induction}
Let $\sig{\Sigma}=(S,\Sigma,\mathrm{def}^{(S,\Sigma)})$ be a hep-signature.
Then to prove that every operation symbol of $\Sigma$ has a property $\mathcal{P}$ it suffices to show that, for every operation symbol $\sigma$ of $\Sigma$, if every $\tau\lhd\sigma$ has $\mathcal{P}$, then $\sigma$ has $\mathcal{P}$ (well-founded induction step).
\end{proposition}

No separate induction basis is required.
If $\sigma$ is $\lhd$-minimal, i.e., if no operation symbol occurs in $\mathrm{def}^{(S,\Sigma)}(\sigma)$, then the hypothesis of the well-founded induction step is vacuously satisfied, and the well-founded induction step proves, without any assumption, that $\sigma$ has $\mathcal{P}$.
Likewise, by well-founded recursion along $\lhd$, an assignment of a datum to every operation symbol of $\Sigma$ is uniquely determined by a rule that determines the datum of every operation symbol $\sigma$ from the data of the operation symbols $\tau\lhd\sigma$; for the $\lhd$-minimal operation symbols, the rule determines the datum from no previous data.

We next introduce the specifications on ep-signatures.

\begin{definition}\label{DEpSpec}\lean{equivalent}{EPSpecification}
An \emph{ep-specification} is an ordered pair $(\sig{\Sigma},\mathcal{E})$, where $\sig{\Sigma}=(S,\Sigma,\mathrm{def}^{(S,\Sigma)})$ is an ep-signature and $\mathcal{E}$ is a set of \emph{contextual $\QE$-equations}, i.e., of ordered pairs $(X,e)$ in which $X$ is a finite $S$-sorted subset of $V^{S}$, the \emph{context} of $e$, and $e\in\mathrm{QEeqt}(\Sigma)_{X}$ is a $\QE$-equation with variables in $X$, see Definition~\ref{DQEEq}.
A contextual $\QE$-equation $(X,e)$ \emph{holds} in a partial $\Sigma$-algebra $\mathbf{A}$ if $e$ holds in $\mathbf{A}$ in the sense of Definition~\ref{DQEEq}, i.e., if $\mathbf{A}$ satisfies $e$ with respect to every valuation $a\in A_{X}$.
When we display the axioms of an ep-specification we only write the $\QE$-equation $e$; unless otherwise stated, its context is the $S$-sorted set of the variables that occur in $e$.
We denote the ep-specification $(\sig{\Sigma},\mathcal{E})$ simply by $\spec{E}$, and we say that $\spec{E}$ is a \emph{hep-specification} if $\sig{\Sigma}$ is a hep-signature, see Definition~\ref{DHepSig}.
An ep-signature $\sig{\Sigma}$ is identified with the ep-specification $(\sig{\Sigma},\varnothing)$, and a hep-signature with the hep-specification $(\sig{\Sigma},\varnothing)$.
\end{definition}

The context of a $\QE$-equation matters when carriers may be empty, see~\cite{GoguenMeseguer1985} and, for equations and total algebras, \cite{Soliveres2000}, Propositions~2.11.10 and~2.11.11 and Example~2.11.13.
If a context $X$ contains a variable of a sort $t$ and $A_{t}=\varnothing$, then there is no valuation of $X$ in $A$, so that every $\QE$-equation with context $X$ holds in $\mathbf{A}$, even if it only constrains the elements of other sorts.
For this reason the axioms of an ep-specification are not evaluated in the whole $S$-sorted set of variables $V^{S}$: otherwise, a single empty carrier would make all of them hold.

Following the case of the $\QE$-specifications, see Definition~\ref{DQEEq}, an ep-specification determines a class of ep-algebras, its ep-variety, and a category canonically associated with it.

\begin{definition}\label{DEpVariety}\lean{equivalent}{EPSpecification.IsModel,EPSpecification.variety}
Let $\spec{E}=(\sig{\Sigma},\mathcal{E})$ be an ep-specification.
We say that a $\sig{\Sigma}$-algebra $\mathbf{A}$ is a \emph{model} of $\spec{E}$ if every contextual $\QE$-equation of $\mathcal{E}$ holds in $\mathbf{A}$, see Definition~\ref{DEpSpec}.
We let $\mathcal{V}(\spec{E})$ stand for the \emph{ep-variety} determined by the ep-specification $\spec{E}$, i.e., the class of all the models of $\spec{E}$, and $\mathsf{EPAlg}(\spec{E})$ for the category canonically associated with $\mathcal{V}(\spec{E})$, i.e., the full subcategory of $\mathsf{EPAlg}(\sig{\Sigma})$ whose objects are the $\sig{\Sigma}$-algebras in $\mathcal{V}(\spec{E})$.
\end{definition}

Let us note that $\mathcal{V}(\sig{\Sigma},\varnothing)$ is the class of all the $\sig{\Sigma}$-algebras, so that $\mathsf{EPAlg}(\sig{\Sigma},\varnothing)=\mathsf{EPAlg}(\sig{\Sigma})$.

\section{Generalized hep-derivors}\label{SGenHepDerivorsDef}

We next define the notion of generalized hep-derivor.

\begin{definition}\label{DGenHepDerivor}\lean{equivalent}{GeneralizedHepDerivor}
Let $\sig{\Sigma}=(S,\Sigma,\mathrm{def}^{(S,\Sigma)})$ and $\sig{\Lambda}=(T,\Lambda,\mathrm{def}^{(T,\Lambda)})$ be hep-signatures, see Definition~\ref{DHepSig}.
A \emph{generalized hep-derivor} from $\sig{\Sigma}$ to $\sig{\Lambda}$ is an ordered triple $\mathbf{c}=(\varphi,c,\Dom^{(\varphi,c)})$ such that
\begin{enumerate}
\item $(\varphi,c)$ is a derivor from $(S,\Sigma)$ to $(T,\Lambda)$, see Section~\ref{SDerivors}, i.e., $\varphi$ is a mapping from $S$ to $T$ and $c$ assigns to every $(\mathbf{s},s)\in S^{\star}\times S$ and every $\sigma\in\Sigma_{\mathbf{s},s}$ a term $c(\sigma)\in\mathrm{T}_{\Lambda}(\vs{\varphi^{\star}(\mathbf{s})})_{\varphi(s)}$; and
\item $\Dom^{(\varphi,c)}=(\Dom^{(\varphi,c)}_{s})_{s\in S}$ assigns to every sort $s\in S$ a set
$$
\Dom^{(\varphi,c)}_{s}\subseteq\mathrm{Eeqt}(\Lambda)_{\vs{\varphi(s)}}
$$
of $\mathrm{E}$-equations in the single variable $v^{\varphi(s)}_{0}$, the \emph{domain} of $s$.
\end{enumerate}
\end{definition}

\begin{convention}
In this chapter, unless otherwise stated, we write the hep-signatures as $\sig{\Sigma}=(S,\Sigma,\mathrm{def}^{(S,\Sigma)})$, $\sig{\Lambda}=(T,\Lambda,\mathrm{def}^{(T,\Lambda)})$, $\sig{\Omega}=(U,\Omega,\mathrm{def}^{(U,\Omega)})$ and $\sig{\Xi}=(W,\Xi,\mathrm{def}^{(W,\Xi)})$, and the generalized hep-derivors as $\mathbf{c}=(\varphi,c,\Dom^{(\varphi,c)})$, $\mathbf{d}=(\psi,d,\Dom^{(\psi,d)})$ and $\mathbf{e}=(\chi,e,\Dom^{(\chi,e)})$.
\end{convention}

Let $\mathbf{c}$ be a generalized hep-derivor from $\sig{\Sigma}$ to $\sig{\Lambda}$ and $\mathbf{B}$ a partial $\Lambda$-algebra.
We next associate to $\mathbf{B}$ the $S$-sorted set of the elements that satisfy the domains of $\mathbf{c}$.
Since it is obtained from the $S$-sorted set $B_{\varphi}=(B_{\varphi(s)})_{s\in S}$ by restricting every $B_{\varphi(s)}$ to the elements that satisfy $\Dom^{(\varphi,c)}_{s}$, we denote it, as the sets $A_{\mathrm{def}^{(S,\Sigma)}(\sigma)}$ of Definition~\ref{DEpAlg}, by subscripting the domains.
Let us recall that, by the identifications of Section~\ref{SVariables}, for every $s\in S$, an element $b\in B_{\varphi(s)}$ is identified with the valuation of $\vs{(\varphi(s))}$ in $B$ that sends $v^{\varphi(s)}_{0}$ to $b$.

\begin{definition}\label{DDomainSet}\lean{equivalent}{GeneralizedHepDerivor.Cell}
Let $\sig{\Sigma}$ and $\sig{\Lambda}$ be hep-signatures, $\mathbf{c}$ a generalized hep-derivor from $\sig{\Sigma}$ to $\sig{\Lambda}$ and $\mathbf{B}$ a partial $\Lambda$-algebra.
We let $B_{\Dom^{(\varphi,c)}}$ stand for the $S$-sorted subset of $B_{\varphi}$ defined, for every $s\in S$, as
$$
\left(B_{\Dom^{(\varphi,c)}}\right)_{s}=\left\lbrace b\in B_{\varphi(s)}\,\middle|\,\mathbf{B}\models\Dom^{(\varphi,c)}_{s}[b]\right\rbrace.
$$
The elements of $(B_{\Dom^{(\varphi,c)}})_{s}$ are called the \emph{cells of sort $s$} of $\mathbf{B}$ relative to $\mathbf{c}$.
\end{definition}

The terminology anticipates the generalized hep-derivors of Chapters~\ref{SGdUf}--\ref{SChangeSig}, for which the cells of a given sort are precisely the cells of a fixed dimension of a higher-order categorial algebra, see, e.g., Definition~\ref{DGdDerivor}.

In the following lemma we prove that the homomorphisms preserve the cells.

\begin{lemma}\label{LDomainSetHom}\lean{equivalent}{GeneralizedHepDerivor.isCell_map}
Let $\sig{\Sigma}$ and $\sig{\Lambda}$ be hep-signatures, $\mathbf{c}$ a generalized hep-derivor from $\sig{\Sigma}$ to $\sig{\Lambda}$ and $f\colon\mathbf{B}\mor\mathbf{B}'$ a homomorphism of partial $\Lambda$-algebras.
Then, for every $s\in S$, $f_{\varphi(s)}$ sends the cells of sort $s$ of $\mathbf{B}$ to cells of sort $s$ of $\mathbf{B}'$.
\end{lemma}
\begin{proof}
Let $b$ be a cell of sort $s$ of $\mathbf{B}$, i.e., $\mathbf{B}\models\Dom^{(\varphi,c)}_{s}[b]$.
Then, by Remark~\ref{RTermOpHom}, $\mathbf{B}'\models\Dom^{(\varphi,c)}_{s}[f_{\varphi(s)}(b)]$, i.e., $f_{\varphi(s)}(b)$ is a cell of sort $s$ of $\mathbf{B}'$.
\end{proof}

By Lemma~\ref{LDomainSetHom} and Proposition~\ref{PImFactorization}, for every homomorphism $f\colon\mathbf{B}\mor\mathbf{B}'$ of partial $\Lambda$-algebras, we can consider the $S$-sorted mapping from $B_{\Dom^{(\varphi,c)}}$ to $B'_{\Dom^{(\varphi,c)}}$
$$
f_{\Dom^{(\varphi,c)}}
=
\left(f_{\varphi(s)}\bigr|^{(B'_{\Dom^{(\varphi,c)}})_{s}}_{(B_{\Dom^{(\varphi,c)}})_{s}}\right)_{s\in S},
$$
i.e., the family of the birestrictions of the mappings $f_{\varphi(s)}$ to the cells.
In the following proposition we prove that this construction is functorial.

\begin{proposition}\label{PDomainSetFunctor}\lean{equivalent}{GeneralizedHepDerivor.cellsFunctor}
Let $\sig{\Sigma}$ and $\sig{\Lambda}$ be hep-signatures and $\mathbf{c}$ a generalized hep-derivor from $\sig{\Sigma}$ to $\sig{\Lambda}$.
Then $\G_{\mathbf{c}}$, the assignment from $\mathsf{PAlg}(\Lambda)$ to $\mathsf{Set}^{S}$ defined as follows:
\begin{enumerate}
\item for every partial $\Lambda$-algebra $\mathbf{B}$, $\G_{\mathbf{c}}(\mathbf{B})=B_{\Dom^{(\varphi,c)}}$; and
\item for every homomorphism $f\colon\mathbf{B}\mor\mathbf{B}'$ of partial $\Lambda$-algebras, $\G_{\mathbf{c}}(f)=f_{\Dom^{(\varphi,c)}}$;
\end{enumerate}
is a functor.
\end{proposition}
\begin{proof}
By Lemma~\ref{LDomainSetHom}, $\G_{\mathbf{c}}(f)$ is an $S$-sorted mapping from $\G_{\mathbf{c}}(\mathbf{B})$ to $\G_{\mathbf{c}}(\mathbf{B}')$.
Since, for every $s\in S$, the birestriction of the identity of $B_{\varphi(s)}$ to $(B_{\Dom^{(\varphi,c)}})_{s}$ is the identity of $(B_{\Dom^{(\varphi,c)}})_{s}$, $\G_{\mathbf{c}}$ preserves identities.
Since, for homomorphisms $f\colon\mathbf{B}\mor\mathbf{B}'$ and $g\colon\mathbf{B}'\mor\mathbf{B}''$, $f_{\varphi(s)}$ sends the cells of sort $s$ of $\mathbf{B}$ to cells of sort $s$ of $\mathbf{B}'$, by Lemma~\ref{LDomainSetHom}, the birestriction of $g_{\varphi(s)}\circ f_{\varphi(s)}$ to the cells is the composite of the birestrictions of $g_{\varphi(s)}$ and of $f_{\varphi(s)}$, so that $\G_{\mathbf{c}}$ preserves composites.
\end{proof}

When all the domains of a generalized hep-derivor are empty, every element is a cell, and the functor $\G_{\mathbf{c}}$ reduces to the usual reindexing along $\varphi$ of the underlying sorted sets.

\begin{remark}\label{RDomainSetEmpty}\lean{weaker}{GeneralizedHepDerivor.cellsFunctorIsoDelta}
Let $\sig{\Sigma}$ and $\sig{\Lambda}$ be hep-signatures and $\mathbf{c}$ a generalized hep-derivor from $\sig{\Sigma}$ to $\sig{\Lambda}$ such that $\Dom^{(\varphi,c)}_{s}=\varnothing$, for every $s\in S$.
Then, for every partial $\Lambda$-algebra $\mathbf{B}$, $B_{\Dom^{(\varphi,c)}}=B_{\varphi}$, and, for every homomorphism $f$, $f_{\Dom^{(\varphi,c)}}=f_{\varphi}$.
Therefore
$$
\G_{\mathbf{c}}=\Delta_{\varphi}\circ\G_{\Lambda},
$$
where $\G_{\Lambda}$ is the forgetful functor from $\mathsf{PAlg}(\Lambda)$ to $\mathsf{Set}^{T}$ and $\Delta_{\varphi}$ is the functor of Proposition~\ref{PDeltaFunctor}.
\end{remark}

We next equip the $S$-sorted set of cells with a structure of partial $\Sigma$-algebra.
Every operation symbol $\sigma$ is interpreted as the partial term operation $c(\sigma)^{\mathbf{B}}$, restricted to the families of cells that satisfy the domain conditions of $\sigma$ and at which $c(\sigma)^{\mathbf{B}}$ is defined and takes a cell as value.
Since the domain conditions of $\sigma$ are evaluated in the algebra that we are defining, the operations are defined by well-founded recursion along $\lhd$.

\begin{proposition}\label{PTransformedAlg}\lean{equivalent}{GeneralizedHepDerivor.transformed,GeneralizedHepDerivor.mem_transformed_op}
Let $\sig{\Sigma}$ and $\sig{\Lambda}$ be hep-signatures, $\mathbf{c}$ a generalized hep-derivor from $\sig{\Sigma}$ to $\sig{\Lambda}$ and $\mathbf{B}$ a partial $\Lambda$-algebra.
Then the construction below, by well-founded recursion along $\lhd$, defines a weak $\sig{\Sigma}$-algebra, see Definition~\ref{DWeakEpAlg}, with underlying $S$-sorted set $B_{\Dom^{(\varphi,c)}}$, which we denote by $\mathbf{c}^{\ast}_{\gd}(\mathbf{B})$.
We call $\mathbf{c}^{\ast}_{\gd}(\mathbf{B})$ the \emph{transformed algebra} of $\mathbf{B}$ along $\mathbf{c}$.
\end{proposition}
\begin{proof}
Let us denote by $\mathbf{c}^{\ast}_{\gd}(\mathbf{B})$ the partial $\Sigma$-algebra defined as follows.

\textsf{(1)} The underlying $S$-sorted set of $\mathbf{c}^{\ast}_{\gd}(\mathbf{B})$ is $B_{\Dom^{(\varphi,c)}}$, see Definition~\ref{DDomainSet}.

\textsf{(2)} For every $(\mathbf{s},s)\in S^{\star}\times S$ and every operation symbol $\sigma\in\Sigma_{\mathbf{s},s}$, the operation $\sigma^{\mathbf{c}^{\ast}_{\gd}(\mathbf{B})}$ is defined by well-founded recursion along $\lhd$, see Proposition~\ref{PWellFoundedInduction} and the note after it, as follows.
Let us assume that, for every $\tau\lhd\sigma$, the operation $\tau^{\mathbf{c}^{\ast}_{\gd}(\mathbf{B})}$ has already been defined.
Then $\sigma^{\mathbf{c}^{\ast}_{\gd}(\mathbf{B})}$ is the partial operation from $(B_{\Dom^{(\varphi,c)}})_{\mathbf{s}}$ to $(B_{\Dom^{(\varphi,c)}})_{s}$ with domain
$$
\Dom\left(\sigma^{\mathbf{c}^{\ast}_{\gd}(\mathbf{B})}\right)
=
\left\lbrace
\mathbf{b}\in\left(B_{\Dom^{(\varphi,c)}}\right)_{\mathbf{s}}
\,\middle|\,
\begin{array}{l}
\mathbf{c}^{\ast}_{\gd}(\mathbf{B})\models\mathrm{def}^{(S,\Sigma)}(\sigma)[\mathbf{b}],\\
\mathbf{b}\in\Dom\left(c(\sigma)^{\mathbf{B}}\right)\text{ and}\\
c(\sigma)^{\mathbf{B}}(\mathbf{b})\in\left(B_{\Dom^{(\varphi,c)}}\right)_{s}
\end{array}
\right\rbrace
$$
which sends every family $\mathbf{b}$ in its domain to $c(\sigma)^{\mathbf{B}}(\mathbf{b})$.

Let us check that this is a well-defined recursion.
By Definition~\ref{DHepSig}, the operation symbols that occur in $\mathrm{def}^{(S,\Sigma)}(\sigma)$ are the $\tau\lhd\sigma$.
Therefore, by Definitions~\ref{DTermOp} and~\ref{DExistenceEq}, whether a family $\mathbf{b}\in(B_{\Dom^{(\varphi,c)}})_{\mathbf{s}}$ satisfies $\mathrm{def}^{(S,\Sigma)}(\sigma)$ in $\mathbf{c}^{\ast}_{\gd}(\mathbf{B})$ only depends on the operations $\tau^{\mathbf{c}^{\ast}_{\gd}(\mathbf{B})}$ with $\tau\lhd\sigma$, which have already been defined.
The other two conditions only depend on $\mathbf{B}$ and on $\mathbf{c}$.
Hence $\Dom(\sigma^{\mathbf{c}^{\ast}_{\gd}(\mathbf{B})})$ is well defined.
Moreover, since $c(\sigma)^{\mathbf{B}}$ sends every family in $\Dom(\sigma^{\mathbf{c}^{\ast}_{\gd}(\mathbf{B})})$ to a cell of sort $s$, $\sigma^{\mathbf{c}^{\ast}_{\gd}(\mathbf{B})}$ is a well-defined partial operation from $(B_{\Dom^{(\varphi,c)}})_{\mathbf{s}}$ to $(B_{\Dom^{(\varphi,c)}})_{s}$.
Since $\lhd$ is well-founded, this defines $\sigma^{\mathbf{c}^{\ast}_{\gd}(\mathbf{B})}$ for every operation symbol $\sigma$ of $\Sigma$.

This completes the definition of the partial $\Sigma$-algebra $\mathbf{c}^{\ast}_{\gd}(\mathbf{B})$.
Finally, by the first condition that defines $\Dom(\sigma^{\mathbf{c}^{\ast}_{\gd}(\mathbf{B})})$, for every operation symbol $\sigma$ of $\Sigma$,
$$
\Dom\left(\sigma^{\mathbf{c}^{\ast}_{\gd}(\mathbf{B})}\right)\subseteq\left(B_{\Dom^{(\varphi,c)}}\right)_{\mathrm{def}^{(S,\Sigma)}(\sigma)}.
$$
Hence $\mathbf{c}^{\ast}_{\gd}(\mathbf{B})$ is a weak $\sig{\Sigma}$-algebra.
\end{proof}

By the proof of Proposition~\ref{PTransformedAlg}, a family of cells $\mathbf{b}$ is in the domain of $\sigma^{\mathbf{c}^{\ast}_{\gd}(\mathbf{B})}$ if, and only if, it satisfies the domain conditions of $\sigma$ in $\mathbf{c}^{\ast}_{\gd}(\mathbf{B})$, $c(\sigma)^{\mathbf{B}}$ is defined at $\mathbf{b}$, and its value is a cell; and then $\sigma^{\mathbf{c}^{\ast}_{\gd}(\mathbf{B})}(\mathbf{b})=c(\sigma)^{\mathbf{B}}(\mathbf{b})$.
In what follows we refer to these three requirements as the three conditions of Proposition~\ref{PTransformedAlg}.

Let us note that, in $\mathbf{c}^{\ast}_{\gd}(\mathbf{B})$, the domain of $\sigma$ is not part of the data, and it may be smaller than that of $c(\sigma)^{\mathbf{B}}$.
For instance, for the discretization functors of Chapter~\ref{SDsTr}, the term assigned to $\circ^{k}_{s}$, for $k$ greater than or equal to the dimension of the algebra, is the variable $v^{s}_{0}$, which is defined everywhere, whereas, in the transformed algebras, the domain condition of $\circ^{k}_{s}$ is satisfied only by the pairs of equal cells.
In~\cite{Reichel1987}, the morphisms of ep-signatures are required, moreover, to translate the domain conditions of $\sigma$ into a set of equations equivalent to the definedness of $c(\sigma)$.

In the following proposition we characterize when the transformed algebra is a $\sig{\Sigma}$-algebra.

\begin{proposition}\label{PTransformedEpAlg}\lean{equivalent}{GeneralizedHepDerivor.transformed_isEPAlgebra_iff}
Let $\sig{\Sigma}$ and $\sig{\Lambda}$ be hep-signatures, $\mathbf{c}$ a generalized hep-derivor from $\sig{\Sigma}$ to $\sig{\Lambda}$ and $\mathbf{B}$ a partial $\Lambda$-algebra.
Then the following statements are equivalent:
\begin{enumerate}
\item $\mathbf{c}^{\ast}_{\gd}(\mathbf{B})$ is a $\sig{\Sigma}$-algebra.
\item For every $(\mathbf{s},s)\in S^{\star}\times S$, every $\sigma\in\Sigma_{\mathbf{s},s}$ and every family of cells $\mathbf{b}\in(B_{\Dom^{(\varphi,c)}})_{\mathbf{s}}$ that satisfies $\mathrm{def}^{(S,\Sigma)}(\sigma)$ in $\mathbf{c}^{\ast}_{\gd}(\mathbf{B})$, $c(\sigma)^{\mathbf{B}}$ is defined at $\mathbf{b}$ and $c(\sigma)^{\mathbf{B}}(\mathbf{b})$ is a cell of sort $s$.
\end{enumerate}
\end{proposition}
\begin{proof}
By Definition~\ref{DEpAlg} and Proposition~\ref{PTransformedAlg}, $\mathbf{c}^{\ast}_{\gd}(\mathbf{B})$ is a $\sig{\Sigma}$-algebra if, and only if, for every operation symbol $\sigma$ of $\Sigma$, every family of cells that satisfies $\mathrm{def}^{(S,\Sigma)}(\sigma)$ in $\mathbf{c}^{\ast}_{\gd}(\mathbf{B})$ is in the domain of $\sigma^{\mathbf{c}^{\ast}_{\gd}(\mathbf{B})}$.
By Proposition~\ref{PTransformedAlg}, such a family is in this domain if, and only if, it satisfies the remaining two conditions of Proposition~\ref{PTransformedAlg}, which are those of statement~(2).
\end{proof}

In the following proposition we prove that the birestrictions of homomorphisms introduced after Lemma~\ref{LDomainSetHom} are homomorphisms between the transformed algebras.

\begin{proposition}\label{PGenHepDerivorHom}\lean{equivalent}{GeneralizedHepDerivor.transformedHom}
Let $\sig{\Sigma}$ and $\sig{\Lambda}$ be hep-signatures, $\mathbf{c}$ a generalized hep-derivor from $\sig{\Sigma}$ to $\sig{\Lambda}$ and $f\colon\mathbf{B}\mor\mathbf{B}'$ a homomorphism of partial $\Lambda$-algebras.
Then $f_{\Dom^{(\varphi,c)}}$ is a homomorphism of partial $\Sigma$-algebras from $\mathbf{c}^{\ast}_{\gd}(\mathbf{B})$ to $\mathbf{c}^{\ast}_{\gd}(\mathbf{B}')$.
\end{proposition}
\begin{proof}
By Lemma~\ref{LDomainSetHom}, $f_{\Dom^{(\varphi,c)}}$ is an $S$-sorted mapping from $B_{\Dom^{(\varphi,c)}}$ to $B'_{\Dom^{(\varphi,c)}}$.
We prove, by well-founded induction along $\lhd$, see Proposition~\ref{PWellFoundedInduction}, that, for every $(\mathbf{s},s)\in S^{\star}\times S$, every $\sigma\in\Sigma_{\mathbf{s},s}$ and every $\mathbf{b}\in\Dom(\sigma^{\mathbf{c}^{\ast}_{\gd}(\mathbf{B})})$, $f_{\Dom^{(\varphi,c)}}\circ\mathbf{b}\in\Dom(\sigma^{\mathbf{c}^{\ast}_{\gd}(\mathbf{B}')})$ and
$$
(f_{\Dom^{(\varphi,c)}})_{s}\left(\sigma^{\mathbf{c}^{\ast}_{\gd}(\mathbf{B})}(\mathbf{b})\right)=\sigma^{\mathbf{c}^{\ast}_{\gd}(\mathbf{B}')}\left(f_{\Dom^{(\varphi,c)}}\circ\mathbf{b}\right).
$$
Let $\sigma\in\Sigma_{\mathbf{s},s}$, let us assume that the above holds for every $\tau\lhd\sigma$, and let $\mathbf{b}\in\Dom(\sigma^{\mathbf{c}^{\ast}_{\gd}(\mathbf{B})})$.

We first prove that $f_{\Dom^{(\varphi,c)}}\circ\mathbf{b}\in\Dom(\sigma^{\mathbf{c}^{\ast}_{\gd}(\mathbf{B}')})$, i.e., that $f_{\Dom^{(\varphi,c)}}\circ\mathbf{b}$ satisfies the three conditions of Proposition~\ref{PTransformedAlg}.
By the induction hypothesis, $f_{\Dom^{(\varphi,c)}}$ commutes with the operations $\tau$ that occur in $\mathrm{def}^{(S,\Sigma)}(\sigma)$, which are the $\tau\lhd\sigma$, i.e., $f_{\Dom^{(\varphi,c)}}$ is a homomorphism between the reducts of $\mathbf{c}^{\ast}_{\gd}(\mathbf{B})$ and $\mathbf{c}^{\ast}_{\gd}(\mathbf{B}')$ to these operation symbols.
Since $\mathbf{c}^{\ast}_{\gd}(\mathbf{B})\models\mathrm{def}^{(S,\Sigma)}(\sigma)[\mathbf{b}]$, by Remark~\ref{RTermOpHom}, applied to these reducts, it follows that $\mathbf{c}^{\ast}_{\gd}(\mathbf{B}')\models\mathrm{def}^{(S,\Sigma)}(\sigma)[f_{\Dom^{(\varphi,c)}}\circ\mathbf{b}]$.
Next, since $\mathbf{b}\in\Dom(c(\sigma)^{\mathbf{B}})$ and $f$ is a homomorphism, by Remark~\ref{RTermOpHom}, $f_{\Dom^{(\varphi,c)}}\circ\mathbf{b}=f_{\varphi}\circ\mathbf{b}$ is in $\Dom(c(\sigma)^{\mathbf{B}'})$ and
$$
c(\sigma)^{\mathbf{B}'}(f_{\Dom^{(\varphi,c)}}\circ\mathbf{b})=f_{\varphi(s)}\left(c(\sigma)^{\mathbf{B}}(\mathbf{b})\right).
$$
Finally, since $c(\sigma)^{\mathbf{B}}(\mathbf{b})$ is a cell of sort $s$ of $\mathbf{B}$, by Lemma~\ref{LDomainSetHom}, $c(\sigma)^{\mathbf{B}'}(f_{\Dom^{(\varphi,c)}}\circ\mathbf{b})$ is a cell of sort $s$ of $\mathbf{B}'$.

Now we prove the equality.
The following chain of equalities holds
\begin{align*}
(f_{\Dom^{(\varphi,c)}})_{s}\left(\sigma^{\mathbf{c}^{\ast}_{\gd}(\mathbf{B})}\left(\mathbf{b}\right)\right)
&=
f_{\varphi(s)}\left(c(\sigma)^{\mathbf{B}}\left(\mathbf{b}\right)\right)
\tag{1}
\\&=
c(\sigma)^{\mathbf{B}'}\left(f_{\Dom^{(\varphi,c)}}\circ\mathbf{b}\right)
\tag{2}
\\&=
\sigma^{\mathbf{c}^{\ast}_{\gd}(\mathbf{B}')}\left(f_{\Dom^{(\varphi,c)}}\circ\mathbf{b}\right).
\tag{3}
\end{align*}

The first equality unravels the definitions of $f_{\Dom^{(\varphi,c)}}$ and of $\sigma^{\mathbf{c}^{\ast}_{\gd}(\mathbf{B})}$, see Proposition~\ref{PTransformedAlg}; the second equality has been proved above; finally, the last equality follows from Proposition~\ref{PTransformedAlg}, since $f_{\Dom^{(\varphi,c)}}\circ\mathbf{b}\in\Dom(\sigma^{\mathbf{c}^{\ast}_{\gd}(\mathbf{B}')})$.

All in all, since $\lhd$ is well-founded, $f_{\Dom^{(\varphi,c)}}$ is a homomorphism from $\mathbf{c}^{\ast}_{\gd}(\mathbf{B})$ to $\mathbf{c}^{\ast}_{\gd}(\mathbf{B}')$. This completes the proof.
\end{proof}

We are now in a position to prove that every generalized hep-derivor induces a functor between the categories of partial algebras.
Let us note that, in general, the transformed algebra of a $\sig{\Lambda}$-algebra need not be a $\sig{\Sigma}$-algebra, see Proposition~\ref{PTransformedEpAlg} and Remark~\ref{RGdNotEpAlg}.
Hence a generalized hep-derivor does not induce, in general, a functor from $\mathsf{EPAlg}(\sig{\Lambda})$ to $\mathsf{EPAlg}(\sig{\Sigma})$, and we define the induced functor on all the partial $\Lambda$-algebras.

\begin{proposition}\label{PGenHepDerivorFunctor}\lean{equivalent}{GeneralizedHepDerivor.transformedFunctor}
Let $\sig{\Sigma}$ and $\sig{\Lambda}$ be hep-signatures and $\mathbf{c}$ a generalized hep-derivor from $\sig{\Sigma}$ to $\sig{\Lambda}$.
Then $\mathbf{c}^{\ast}_{\gd}$, the assignment from $\mathsf{PAlg}(\Lambda)$ to $\mathsf{PAlg}(\Sigma)$ defined as follows:
\begin{enumerate}
\item for every partial $\Lambda$-algebra $\mathbf{B}$, $\mathbf{c}^{\ast}_{\gd}(\mathbf{B})$ is the transformed algebra of Proposition~\ref{PTransformedAlg}; and
\item for every homomorphism $f\colon\mathbf{B}\mor\mathbf{B}'$ of partial $\Lambda$-algebras, $\mathbf{c}^{\ast}_{\gd}(f)=f_{\Dom^{(\varphi,c)}}$, see Proposition~\ref{PGenHepDerivorHom};
\end{enumerate}
is a functor.
\end{proposition}
\begin{proof}
By Propositions~\ref{PTransformedAlg} and~\ref{PGenHepDerivorHom}, it suffices to prove that $\mathbf{c}^{\ast}_{\gd}$ preserves identities and composites.
Since $\mathbf{c}^{\ast}_{\gd}$ and $\G_{\mathbf{c}}$ act in the same way on homomorphisms, this follows from Proposition~\ref{PDomainSetFunctor}.
\end{proof}

From Propositions~\ref{PDomainSetFunctor} and~\ref{PGenHepDerivorFunctor} we obtain the following corollary.

\begin{corollary}\label{CGenHepDerivorForgetful}\lean{equivalent}{GeneralizedHepDerivor.transformedFunctor_comp_forget}
Let $\sig{\Sigma}$ and $\sig{\Lambda}$ be hep-signatures and $\mathbf{c}$ a generalized hep-derivor from $\sig{\Sigma}$ to $\sig{\Lambda}$.
Then the following diagram commutes
$$
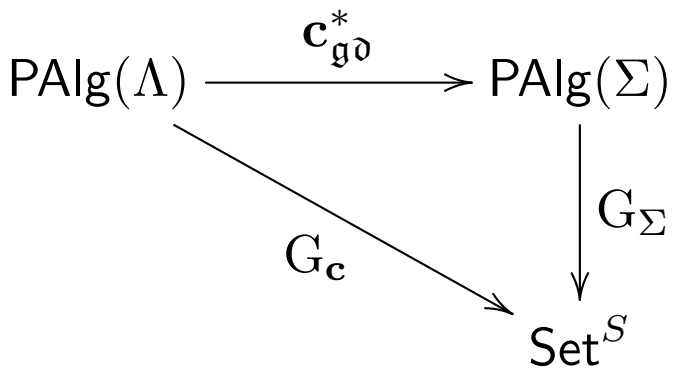
\xymatrix@C=50pt@R=35pt{
\mathsf{PAlg}(\Lambda)\ar[r]^-{\mathbf{c}^{\ast}_{\gd}}\ar[rd]_-{\G_{\mathbf{c}}} &
\mathsf{PAlg}(\Sigma)\ar[d]^-{\G_{\Sigma}}\\
& \mathsf{Set}^{S}
}
$$
where $\G_{\Sigma}$ is the forgetful functor.
\end{corollary}
\begin{proof}
It follows from Propositions~\ref{PDomainSetFunctor}, \ref{PTransformedAlg} and~\ref{PGenHepDerivorFunctor}, since the underlying $S$-sorted set of $\mathbf{c}^{\ast}_{\gd}(\mathbf{B})$ is $B_{\Dom^{(\varphi,c)}}$ and the underlying $S$-sorted mapping of $\mathbf{c}^{\ast}_{\gd}(f)$ is $f_{\Dom^{(\varphi,c)}}$.
\end{proof}

\section{The category \texorpdfstring{$\mathsf{HSig}_{\gd}$}{HSiggd}}\label{SHepSigCat}

In this section we define the composition of generalized hep-derivors and we prove that the hep-signatures and the generalized hep-derivors form a category.
We use the notation $d^{\sharp}[\mathcal{H}]$ introduced after Definition~\ref{DExistenceEq} for the translation of a set of $\mathrm{E}$-equations $\mathcal{H}$ along a derivor $(\psi,d)$.

\begin{definition}\label{DGenHepDerivorCompId}\lean{equivalent}{GeneralizedHepDerivor.comp,GeneralizedHepDerivor.identity}
Let $\sig{\Sigma}$, $\sig{\Lambda}$ and $\sig{\Omega}$ be hep-signatures, $\mathbf{c}$ a generalized hep-derivor from $\sig{\Sigma}$ to $\sig{\Lambda}$ and $\mathbf{d}$ a generalized hep-derivor from $\sig{\Lambda}$ to $\sig{\Omega}$.
The \emph{composition} of $\mathbf{c}$ and $\mathbf{d}$ is the generalized hep-derivor
$$
\mathbf{d}\circ\mathbf{c}=\left(\psi\circ\varphi,d^{\sharp}\circ c,\Dom^{\mathbf{d}\circ\mathbf{c}}\right)
$$
from $\sig{\Sigma}$ to $\sig{\Omega}$, where $(\psi\circ\varphi,d^{\sharp}\circ c)$ is the composition of the derivors $(\varphi,c)$ and $(\psi,d)$, see Definition~\ref{DDerivorCompId}, and, for every $s\in S$,
$$
\Dom^{\mathbf{d}\circ\mathbf{c}}_{s}=\Dom^{(\psi,d)}_{\varphi(s)}\cup d^{\sharp}\left[\Dom^{(\varphi,c)}_{s}\right].
$$
The \emph{identity} at $\sig{\Sigma}$ is the generalized hep-derivor $\id_{\sig{\Sigma}}=(\id_{S},\eta^{\Sigma},\varnothing)$ from $\sig{\Sigma}$ to $\sig{\Sigma}$, where $(\id_{S},\eta^{\Sigma})$ is the identity derivor at $(S,\Sigma)$, see Definition~\ref{DDerivorCompId}, and every domain is empty.
\end{definition}

Let us note that $\Dom^{(\psi,d)}_{\varphi(s)}$ and $d^{\sharp}[\Dom^{(\varphi,c)}_{s}]$ are sets of $\mathrm{E}$-equations with variables in $\vs{\psi(\varphi(s))}$, since $\psi^{\star}((\varphi(s)))=(\psi(\varphi(s)))$, so that the above definition is correct.

We are now in a position to prove that the hep-signatures and the generalized hep-derivors form a category.

\begin{proposition}\label{PHepSigCat}\lean{equivalent}{instCategoryBundledHepSignature}
The hep-signatures, together with the generalized hep-derivors, the composition and the identities of Definition~\ref{DGenHepDerivorCompId}, form a category, which we denote by $\mathsf{HSig}_{\gd}$.
\end{proposition}
\begin{proof}
Let $\sig{\Sigma}$, $\sig{\Lambda}$, $\sig{\Omega}$ and $\sig{\Xi}$ be hep-signatures, and $\mathbf{c}$, $\mathbf{d}$ and $\mathbf{e}$ generalized hep-derivors from $\sig{\Sigma}$ to $\sig{\Lambda}$, from $\sig{\Lambda}$ to $\sig{\Omega}$ and from $\sig{\Omega}$ to $\sig{\Xi}$, respectively.
By Corollary~\ref{CSigdCat}, the underlying derivors of $\mathbf{e}\circ(\mathbf{d}\circ\mathbf{c})$ and of $(\mathbf{e}\circ\mathbf{d})\circ\mathbf{c}$ are equal, and the underlying derivors of $\id_{\sig{\Lambda}}\circ\mathbf{c}$ and of $\mathbf{c}\circ\id_{\sig{\Sigma}}$ are $(\varphi,c)$.
It remains to compare the domains.
Let $s\in S$.

We first prove the associativity.
The following chain of equalities holds
\begin{align*}
\Dom^{\mathbf{e}\circ(\mathbf{d}\circ\mathbf{c})}_{s}
&=
\Dom^{(\chi,e)}_{\psi(\varphi(s))}\cup e^{\sharp}\left[\Dom^{(\psi,d)}_{\varphi(s)}\cup d^{\sharp}\left[\Dom^{(\varphi,c)}_{s}\right]\right]
\tag{1}
\\&=
\Dom^{(\chi,e)}_{\psi(\varphi(s))}\cup e^{\sharp}\left[\Dom^{(\psi,d)}_{\varphi(s)}\right]\cup e^{\sharp}\left[d^{\sharp}\left[\Dom^{(\varphi,c)}_{s}\right]\right]
\tag{2}
\\&=
\Dom^{\mathbf{e}\circ\mathbf{d}}_{\varphi(s)}\cup\left(e^{\sharp}\circ d\right)^{\sharp}\left[\Dom^{(\varphi,c)}_{s}\right]
\tag{3}
\\&=
\Dom^{(\mathbf{e}\circ\mathbf{d})\circ\mathbf{c}}_{s}.
\tag{4}
\end{align*}

The first equality unravels Definition~\ref{DGenHepDerivorCompId}; the second equality follows from the fact that $e^{\sharp}[\,\cdot\,]$ preserves unions; the third equality follows from Definition~\ref{DGenHepDerivorCompId} and the equality $(e^{\sharp}\circ d)^{\sharp}=e^{\sharp}\circ d^{\sharp}$, see Corollary~\ref{CDerivorCatLaws}; finally, the last equality unravels Definition~\ref{DGenHepDerivorCompId}.

Now we prove that the identities are neutral.
By Definitions~\ref{DGenHepDerivorCompId} and~\ref{DDerivorCompId}, the following equalities hold
\begin{align*}
\Dom^{\id_{\sig{\Lambda}}\circ\mathbf{c}}_{s}
&=\varnothing\cup(\eta^{\Lambda})^{\sharp}\left[\Dom^{(\varphi,c)}_{s}\right]=\Dom^{(\varphi,c)}_{s},
\\
\Dom^{\mathbf{c}\circ\id_{\sig{\Sigma}}}_{s}
&=\Dom^{(\varphi,c)}_{s}\cup c^{\sharp}[\varnothing]=\Dom^{(\varphi,c)}_{s},
\end{align*}
since $(\eta^{\Lambda})^{\sharp}$ is the identity, by Corollary~\ref{CDerivorCatLaws}.

This completes the proof.
\end{proof}

\section{Composition and induced functors}\label{SHepCompFunctors}

We now compare the functor induced by the composition of two generalized hep-derivors with the composite of the induced functors.
The composition of Definition~\ref{DGenHepDerivorCompId} translates the domains of $\mathbf{c}$ by means of the terms of $\mathbf{d}$, but these terms may be defined in a partial $\Omega$-algebra $\mathbf{C}$ on more families than the corresponding operations of $\mathbf{d}^{\ast}_{\gd}(\mathbf{C})$.
Therefore the two functors need not coincide; we prove that the second one is always included in the first one, and that they coincide under a condition that we will check in Chapters~\ref{SGdUf} and~\ref{SDsTr}.
We begin by relating the partial term operations of $\mathbf{d}^{\ast}_{\gd}(\mathbf{C})$ with those of $\mathbf{C}$.

\begin{lemma}\label{LGenHepDerivorTermOp}\lean{stronger}{GeneralizedHepDerivor.partialEval_translate}
Let $\sig{\Lambda}$ and $\sig{\Omega}$ be hep-signatures, $\mathbf{d}$ a generalized hep-derivor from $\sig{\Lambda}$ to $\sig{\Omega}$ and $\mathbf{C}$ a partial $\Omega$-algebra.
Let $\mathbf{t}\in T^{\star}$, $t\in T$, $P\in\mathrm{T}_{\Lambda}(\vs{\mathbf{t}})_{t}$ and $\mathbf{x}\in(C_{\Dom^{(\psi,d)}})_{\mathbf{t}}$.
If $\mathbf{x}\in\Dom(P^{\mathbf{d}^{\ast}_{\gd}(\mathbf{C})})$, then $\mathbf{x}\in\Dom(d^{\sharp}(P)^{\mathbf{C}})$ and
$$
d^{\sharp}(P)^{\mathbf{C}}(\mathbf{x})=P^{\mathbf{d}^{\ast}_{\gd}(\mathbf{C})}(\mathbf{x}).
$$
\end{lemma}
\begin{proof}
We proceed by induction on the term $P$.

If $P=v^{\mathbf{t}}_{j}$, for some $j\in\bb{\mathbf{t}}$, then, by~(Der1), $d^{\sharp}(P)=v^{\psi^{\star}(\mathbf{t})}_{j}$, so that 
\[d^{\sharp}(P)^{\mathbf{C}}(\mathbf{x})=x_{j}=P^{\mathbf{d}^{\ast}_{\gd}(\mathbf{C})}(\mathbf{x}).\]

Let us assume that $P=\tau^{\mathbf{T}_{\Lambda}(\vs{\mathbf{t}})}((P_{j})_{j\in\bb{\mathbf{r}}})$, with $\mathbf{r}\in T^{\star}$, $\tau\in\Lambda_{\mathbf{r},t}$ and $(P_{j})_{j\in\bb{\mathbf{r}}}\in\mathrm{T}_{\Lambda}(\vs{\mathbf{t}})_{\mathbf{r}}$, and that the statement holds for every $P_{j}$.
Let $\mathbf{x}\in\Dom(P^{\mathbf{d}^{\ast}_{\gd}(\mathbf{C})})$ and let $\mathbf{y}$ stand for the family $(P_{j}^{\mathbf{d}^{\ast}_{\gd}(\mathbf{C})}(\mathbf{x}))_{j\in\bb{\mathbf{r}}}$, which, by Definition~\ref{DTermOp}, is in $\Dom(\tau^{\mathbf{d}^{\ast}_{\gd}(\mathbf{C})})$.
By Proposition~\ref{PTransformedAlg}, $\mathbf{y}\in\Dom(d(\tau)^{\mathbf{C}})$ and $\tau^{\mathbf{d}^{\ast}_{\gd}(\mathbf{C})}(\mathbf{y})=d(\tau)^{\mathbf{C}}(\mathbf{y})$.
By~(Der2), $d^{\sharp}(P)=(d^{\sharp}(P_{j}))^{\sharp}_{j\in\bb{\mathbf{r}}}(d(\tau))$.
Therefore, by the induction hypothesis and Remark~\ref{RTermOpSubst}, $\mathbf{x}\in\Dom(d^{\sharp}(P)^{\mathbf{C}})$, and the following chain of equalities holds
\begin{align*}
d^{\sharp}(P)^{\mathbf{C}}\left(\mathbf{x}\right)
&=
\left(\left(d^{\sharp}(P_{j})\right)^{\sharp}_{j\in\bb{\mathbf{r}}}\left(d(\tau)\right)\right)^{\mathbf{C}}\left(\mathbf{x}\right)
\tag{1}
\\&=
d(\tau)^{\mathbf{C}}\left(\left(d^{\sharp}(P_{j})^{\mathbf{C}}(\mathbf{x})\right)_{j\in\bb{\mathbf{r}}}\right)
\tag{2}
\\&=
d(\tau)^{\mathbf{C}}\left(\mathbf{y}\right)
\tag{3}
\\&=
\tau^{\mathbf{d}^{\ast}_{\gd}(\mathbf{C})}\left(\mathbf{y}\right)
\tag{4}
\\&=
P^{\mathbf{d}^{\ast}_{\gd}(\mathbf{C})}\left(\mathbf{x}\right).
\tag{5}
\end{align*}

The first equality follows from~(Der2); the second equality follows from Remark~\ref{RTermOpSubst}; the third equality follows from the induction hypothesis; the fourth equality follows from Proposition~\ref{PTransformedAlg}; finally, the last equality unravels the definition of the partial term operation $P^{\mathbf{d}^{\ast}_{\gd}(\mathbf{C})}$, see Definition~\ref{DTermOp}.

This completes the proof.
\end{proof}

From the above lemma we obtain the corresponding fact on the $\mathrm{E}$-equations.

\begin{corollary}\label{CGenHepDerivorEq}\lean{equivalent}{GeneralizedHepDerivor.satisfies_translateWord}
Let $\sig{\Lambda}$ and $\sig{\Omega}$ be hep-signatures, $\mathbf{d}$ a generalized hep-derivor from $\sig{\Lambda}$ to $\sig{\Omega}$, $\mathbf{C}$ a partial $\Omega$-algebra, $\mathbf{t}\in T^{\star}$, $\mathbf{x}\in(C_{\Dom^{(\psi,d)}})_{\mathbf{t}}$ and $\mathcal{H}\subseteq\mathrm{Eeqt}(\Lambda)_{\vs{\mathbf{t}}}$.
If $\mathbf{d}^{\ast}_{\gd}(\mathbf{C})\models\mathcal{H}[\mathbf{x}]$, then $\mathbf{C}\models d^{\sharp}[\mathcal{H}][\mathbf{x}]$.
\end{corollary}
\begin{proof}
Let $(P,Q)\in\mathcal{H}$ be such that $\mathbf{d}^{\ast}_{\gd}(\mathbf{C})\models(P,Q)[\mathbf{x}]$.
Then $P^{\mathbf{d}^{\ast}_{\gd}(\mathbf{C})}(\mathbf{x})$ and $Q^{\mathbf{d}^{\ast}_{\gd}(\mathbf{C})}(\mathbf{x})$ are defined and equal, so that, by Lemma~\ref{LGenHepDerivorTermOp}, $d^{\sharp}(P)^{\mathbf{C}}(\mathbf{x})$ and $d^{\sharp}(Q)^{\mathbf{C}}(\mathbf{x})$ are defined and equal, i.e., $\mathbf{C}\models(d^{\sharp}(P),d^{\sharp}(Q))[\mathbf{x}]$.
\end{proof}

In the following lemma we prove that the transformed algebra only depends on the domains through the cells they determine.

\begin{lemma}\label{LGenHepDerivorSameCells}\lean{weaker}{GeneralizedHepDerivor.sameDerivorIso,GeneralizedHepDerivor.sameCellsIso_naturality}
Let $\sig{\Sigma}$ and $\sig{\Lambda}$ be hep-signatures, $\mathbf{c}=(\varphi,c,\Dom)$ and $\mathbf{c}'=(\varphi,c,\Dom')$ generalized hep-derivors from $\sig{\Sigma}$ to $\sig{\Lambda}$ with the same underlying derivor, and $\mathbf{B}$ a partial $\Lambda$-algebra such that $B_{\Dom}=B_{\Dom'}$.
Then $\mathbf{c}^{\ast}_{\gd}(\mathbf{B})=\mathbf{c}'^{\ast}_{\gd}(\mathbf{B})$.
Moreover, for every homomorphism $f\colon\mathbf{B}\mor\mathbf{B}'$ of partial $\Lambda$-algebras such that $B'_{\Dom}=B'_{\Dom'}$, $\mathbf{c}^{\ast}_{\gd}(f)=\mathbf{c}'^{\ast}_{\gd}(f)$.
\end{lemma}
\begin{proof}
Both partial $\Sigma$-algebras have the same underlying $S$-sorted set.
We prove, by well-founded induction along $\lhd$, see Proposition~\ref{PWellFoundedInduction}, that every $\sigma\in\Sigma_{\mathbf{s},s}$ has the same interpretation in both.
Let us assume that every $\tau\lhd\sigma$ does.
Then, by Definition~\ref{DHepSig}, a family of cells satisfies $\mathrm{def}^{(S,\Sigma)}(\sigma)$ in $\mathbf{c}^{\ast}_{\gd}(\mathbf{B})$ if, and only if, it satisfies it in $\mathbf{c}'^{\ast}_{\gd}(\mathbf{B})$.
The remaining two conditions of Proposition~\ref{PTransformedAlg} only depend on $c(\sigma)$, on $\mathbf{B}$ and on the common cells.
Hence, by Proposition~\ref{PTransformedAlg}, $\sigma$ has the same domain in both algebras and, on it, both interpretations send $\mathbf{b}$ to $c(\sigma)^{\mathbf{B}}(\mathbf{b})$.
Finally, by Proposition~\ref{PGenHepDerivorFunctor}, $\mathbf{c}^{\ast}_{\gd}(f)$ and $\mathbf{c}'^{\ast}_{\gd}(f)$ are the birestrictions of $f_{\varphi}$ to the same cells.

This completes the proof.
\end{proof}

In the following proposition we prove that the identity generalized hep-derivor leaves the weak $\sig{\Sigma}$-algebras and their homomorphisms unchanged.

\begin{proposition}\label{PGenHepDerivorId}\lean{weaker}{GeneralizedHepDerivor.identityIso,GeneralizedHepDerivor.weakIdentityIso}
Let $\sig{\Sigma}$ be a hep-signature.
Then, for every weak $\sig{\Sigma}$-algebra $\mathbf{A}$ and every homomorphism $f$ between weak $\sig{\Sigma}$-algebras, see Definition~\ref{DWeakEpAlg}, $(\id_{\sig{\Sigma}})^{\ast}_{\gd}(\mathbf{A})=\mathbf{A}$ and $(\id_{\sig{\Sigma}})^{\ast}_{\gd}(f)=f$.
\end{proposition}
\begin{proof}
Let $\mathbf{A}$ be a weak $\sig{\Sigma}$-algebra.
Since the domains of $\id_{\sig{\Sigma}}$ are empty, by Remark~\ref{RDomainSetEmpty}, the underlying $S$-sorted set of $(\id_{\sig{\Sigma}})^{\ast}_{\gd}(\mathbf{A})$ is $A$, and $(\id_{\sig{\Sigma}})^{\ast}_{\gd}(f)=f$.
We prove, by well-founded induction along $\lhd$, see Proposition~\ref{PWellFoundedInduction}, that every $\sigma\in\Sigma_{\mathbf{s},s}$ has the same interpretation in $\mathbf{A}$ and in $(\id_{\sig{\Sigma}})^{\ast}_{\gd}(\mathbf{A})$.

Let $\sigma\in\Sigma_{\mathbf{s},s}$ and assume that every $\tau\lhd\sigma$ has the same interpretation in both algebras.
Then, by Definition~\ref{DHepSig}, a family $\mathbf{a}\in A_{\mathbf{s}}$ satisfies $\mathrm{def}^{(S,\Sigma)}(\sigma)$ in $\mathbf{A}$ if, and only if, it satisfies it in $(\id_{\sig{\Sigma}})^{\ast}_{\gd}(\mathbf{A})$.
Moreover, by Definitions~\ref{DDerivorCompId} and~\ref{DTermOp}, $\eta^{\Sigma}(\sigma)^{\mathbf{A}}=\sigma^{\mathbf{A}}$, since $\eta^{\Sigma}(\sigma)=\sigma^{\mathbf{T}_{\Sigma}(\vs{\mathbf{s}})}((v^{\mathbf{s}}_{j})_{j\in\bb{\mathbf{s}}})$.
Therefore, the following chain of equalities holds
\begin{align*}
\Dom\left(\sigma^{(\id_{\sig{\Sigma}})^{\ast}_{\gd}(\mathbf{A})}\right)
&=
A_{\mathrm{def}^{(S,\Sigma)}(\sigma)}\cap\Dom\left(\eta^{\Sigma}(\sigma)^{\mathbf{A}}\right)
\tag{1}
\\&=
A_{\mathrm{def}^{(S,\Sigma)}(\sigma)}\cap\Dom\left(\sigma^{\mathbf{A}}\right)
\tag{2}
\\&=
\Dom\left(\sigma^{\mathbf{A}}\right),
\tag{3}
\end{align*}
and, on this domain, $\sigma^{(\id_{\sig{\Sigma}})^{\ast}_{\gd}(\mathbf{A})}$ acts as $\eta^{\Sigma}(\sigma)^{\mathbf{A}}=\sigma^{\mathbf{A}}$.

The first equality follows from Proposition~\ref{PTransformedAlg}, taking into account the induction hypothesis and the fact that all the elements of $A$ are cells, so that the last condition of Proposition~\ref{PTransformedAlg} always holds; the second equality follows from the fact that $\eta^{\Sigma}(\sigma)^{\mathbf{A}}=\sigma^{\mathbf{A}}$; finally, the last equality follows from the fact that $\mathbf{A}$ is a weak $\sig{\Sigma}$-algebra, so that $\Dom(\sigma^{\mathbf{A}})\subseteq A_{\mathrm{def}^{(S,\Sigma)}(\sigma)}$, see Definition~\ref{DWeakEpAlg}.

This completes the proof.
\end{proof}

From the above proposition we obtain that the identity generalized hep-derivor induces, on the weak $\sig{\Sigma}$-algebras, and, in particular, on the $\sig{\Sigma}$-algebras, the identity functor.

\begin{corollary}\label{CGenHepDerivorIdFunctor}\lean{weaker}{GeneralizedHepDerivor.epIdentityIso}
Let $\sig{\Sigma}$ be a hep-signature.
Then the restriction of $(\id_{\sig{\Sigma}})^{\ast}_{\gd}$ to $\mathsf{WEPAlg}(\sig{\Sigma})$ is $\Id_{\mathsf{WEPAlg}(\sig{\Sigma})}$ and its restriction to $\mathsf{EPAlg}(\sig{\Sigma})$ is $\Id_{\mathsf{EPAlg}(\sig{\Sigma})}$.
\end{corollary}
\begin{proof}
It follows from Proposition~\ref{PGenHepDerivorId}, since $\mathsf{WEPAlg}(\sig{\Sigma})$ and $\mathsf{EPAlg}(\sig{\Sigma})$ are full subcategories of $\mathsf{PAlg}(\Sigma)$ and every $\sig{\Sigma}$-algebra is a weak $\sig{\Sigma}$-algebra.
\end{proof}

We are now in a position to compare the functor induced by a composition with the composite of the induced functors.

\begin{proposition}\label{PGenHepDerivorComp}\lean{weaker}{GeneralizedHepDerivor.comparison,GeneralizedHepDerivor.comparison_naturality}
Let $\sig{\Sigma}$, $\sig{\Lambda}$ and $\sig{\Omega}$ be hep-signatures, $\mathbf{c}$ a generalized hep-derivor from $\sig{\Sigma}$ to $\sig{\Lambda}$, $\mathbf{d}$ a generalized hep-derivor from $\sig{\Lambda}$ to $\sig{\Omega}$ and $\mathbf{C}$ a partial $\Omega$-algebra.
Then the underlying $S$-sorted set of $\mathbf{c}^{\ast}_{\gd}(\mathbf{d}^{\ast}_{\gd}(\mathbf{C}))$ is an $S$-sorted subset of that of $(\mathbf{d}\circ\mathbf{c})^{\ast}_{\gd}(\mathbf{C})$, and the inclusion $\kappa^{\mathbf{c},\mathbf{d}}_{\mathbf{C}}$ of the former into the latter is a homomorphism of partial $\Sigma$-algebras
$$
\kappa^{\mathbf{c},\mathbf{d}}_{\mathbf{C}}\colon\mathbf{c}^{\ast}_{\gd}\left(\mathbf{d}^{\ast}_{\gd}(\mathbf{C})\right)\mor(\mathbf{d}\circ\mathbf{c})^{\ast}_{\gd}(\mathbf{C}).
$$
Moreover, for every homomorphism $f\colon\mathbf{C}\mor\mathbf{C}'$ of partial $\Omega$-algebras,
$$
(\mathbf{d}\circ\mathbf{c})^{\ast}_{\gd}(f)\circ\kappa^{\mathbf{c},\mathbf{d}}_{\mathbf{C}}=\kappa^{\mathbf{c},\mathbf{d}}_{\mathbf{C}'}\circ\mathbf{c}^{\ast}_{\gd}\left(\mathbf{d}^{\ast}_{\gd}(f)\right).
$$
\end{proposition}
\begin{proof}
In this proof we let $\mathbf{A}$ and $\mathbf{A}'$ stand for $\mathbf{c}^{\ast}_{\gd}(\mathbf{d}^{\ast}_{\gd}(\mathbf{C}))$ and $(\mathbf{d}\circ\mathbf{c})^{\ast}_{\gd}(\mathbf{C})$, respectively.

We first prove that the cells of $\mathbf{A}$ are cells of $\mathbf{A}'$.
Let $s\in S$ and $x$ a cell of sort $s$ of $\mathbf{A}$, i.e., $x\in(C_{\Dom^{(\psi,d)}})_{\varphi(s)}$ and $\mathbf{d}^{\ast}_{\gd}(\mathbf{C})\models\Dom^{(\varphi,c)}_{s}[x]$.
By Corollary~\ref{CGenHepDerivorEq}, $\mathbf{C}\models d^{\sharp}[\Dom^{(\varphi,c)}_{s}][x]$.
Moreover, since $x$ is a cell of sort $\varphi(s)$ of $\mathbf{C}$ relative to $\mathbf{d}$, $\mathbf{C}\models\Dom^{(\psi,d)}_{\varphi(s)}[x]$.
Hence, by Definition~\ref{DGenHepDerivorCompId}, $\mathbf{C}\models\Dom^{\mathbf{d}\circ\mathbf{c}}_{s}[x]$, i.e., $x$ is a cell of sort $s$ of $\mathbf{A}'$.

Now we prove that the inclusion is a homomorphism.
We prove, by well-founded induction along $\lhd$, see Proposition~\ref{PWellFoundedInduction}, that, for every $\sigma\in\Sigma_{\mathbf{s},s}$, every $\mathbf{a}\in\Dom(\sigma^{\mathbf{A}})$ is in $\Dom(\sigma^{\mathbf{A}'})$ and $\sigma^{\mathbf{A}'}(\mathbf{a})=\sigma^{\mathbf{A}}(\mathbf{a})$.
Let us assume that this holds for every $\tau\lhd\sigma$ and let $\mathbf{a}\in\Dom(\sigma^{\mathbf{A}})$.
We check the three conditions of Proposition~\ref{PTransformedAlg} for $\mathbf{d}\circ\mathbf{c}$.
First, by the induction hypothesis, the inclusion commutes with the operations that occur in $\mathrm{def}^{(S,\Sigma)}(\sigma)$, so that, since $\mathbf{A}\models\mathrm{def}^{(S,\Sigma)}(\sigma)[\mathbf{a}]$, by Remark~\ref{RTermOpHom}, applied to the reducts of $\mathbf{A}$ and $\mathbf{A}'$ to these operation symbols, $\mathbf{A}'\models\mathrm{def}^{(S,\Sigma)}(\sigma)[\mathbf{a}]$.
Second, by Proposition~\ref{PTransformedAlg}, $\mathbf{a}\in\Dom(c(\sigma)^{\mathbf{d}^{\ast}_{\gd}(\mathbf{C})})$, so that, by Lemma~\ref{LGenHepDerivorTermOp}, $\mathbf{a}\in\Dom(d^{\sharp}(c(\sigma))^{\mathbf{C}})$ and $d^{\sharp}(c(\sigma))^{\mathbf{C}}(\mathbf{a})=c(\sigma)^{\mathbf{d}^{\ast}_{\gd}(\mathbf{C})}(\mathbf{a})$.
Third, by Proposition~\ref{PTransformedAlg}, $c(\sigma)^{\mathbf{d}^{\ast}_{\gd}(\mathbf{C})}(\mathbf{a})=\sigma^{\mathbf{A}}(\mathbf{a})$ is a cell of $\mathbf{A}$, hence, by the first part of the proof, a cell of $\mathbf{A}'$.
Therefore, by Proposition~\ref{PTransformedAlg}, $\mathbf{a}\in\Dom(\sigma^{\mathbf{A}'})$ and the following chain of equalities holds
\begin{align*}
\sigma^{\mathbf{A}'}\left(\mathbf{a}\right)
&=
d^{\sharp}(c(\sigma))^{\mathbf{C}}\left(\mathbf{a}\right)
\tag{1}
\\&=
c(\sigma)^{\mathbf{d}^{\ast}_{\gd}(\mathbf{C})}\left(\mathbf{a}\right)
\tag{2}
\\&=
\sigma^{\mathbf{A}}\left(\mathbf{a}\right).
\tag{3}
\end{align*}

The first equality follows from Proposition~\ref{PTransformedAlg} applied to $\mathbf{d}\circ\mathbf{c}$; the second equality follows from Lemma~\ref{LGenHepDerivorTermOp}; finally, the last equality follows from Proposition~\ref{PTransformedAlg} applied to $\mathbf{c}$.

Finally, we prove the naturality.
By Proposition~\ref{PGenHepDerivorFunctor}, for every $s\in S$, both $((\mathbf{d}\circ\mathbf{c})^{\ast}_{\gd}(f)\circ\kappa^{\mathbf{c},\mathbf{d}}_{\mathbf{C}})_{s}$ and $(\kappa^{\mathbf{c},\mathbf{d}}_{\mathbf{C}'}\circ\mathbf{c}^{\ast}_{\gd}(\mathbf{d}^{\ast}_{\gd}(f)))_{s}$ send a cell $x$ of sort $s$ of $\mathbf{A}$ to $f_{\psi(\varphi(s))}(x)$. This completes the proof.
\end{proof}

From the above proposition we obtain a sufficient condition for the functor induced by a composition to agree with the composite of the induced functors.

\begin{corollary}\label{CGenHepDerivorCompEqual}\lean{weaker}{GeneralizedHepDerivor.comparisonIso}
Let $\sig{\Sigma}$, $\sig{\Lambda}$ and $\sig{\Omega}$ be hep-signatures, $\mathbf{c}$ a generalized hep-derivor from $\sig{\Sigma}$ to $\sig{\Lambda}$, $\mathbf{d}$ a generalized hep-derivor from $\sig{\Lambda}$ to $\sig{\Omega}$ and $\mathbf{C}$ a partial $\Omega$-algebra such that $\mathbf{c}^{\ast}_{\gd}(\mathbf{d}^{\ast}_{\gd}(\mathbf{C}))$ and $(\mathbf{d}\circ\mathbf{c})^{\ast}_{\gd}(\mathbf{C})$ have the same underlying $S$-sorted set and $\mathbf{c}^{\ast}_{\gd}(\mathbf{d}^{\ast}_{\gd}(\mathbf{C}))$ is a $\sig{\Sigma}$-algebra.
Then
$$
(\mathbf{d}\circ\mathbf{c})^{\ast}_{\gd}(\mathbf{C})=\mathbf{c}^{\ast}_{\gd}\left(\mathbf{d}^{\ast}_{\gd}(\mathbf{C})\right).
$$
Moreover, if $f\colon\mathbf{C}\mor\mathbf{C}'$ is a homomorphism of partial $\Omega$-algebras and $\mathbf{C}'$ satisfies the same hypotheses, then $(\mathbf{d}\circ\mathbf{c})^{\ast}_{\gd}(f)=\mathbf{c}^{\ast}_{\gd}(\mathbf{d}^{\ast}_{\gd}(f))$.
\end{corollary}
\begin{proof}
As in the proof of Proposition~\ref{PGenHepDerivorComp}, we let $\mathbf{A}$ and $\mathbf{A}'$ stand for $\mathbf{c}^{\ast}_{\gd}(\mathbf{d}^{\ast}_{\gd}(\mathbf{C}))$ and $(\mathbf{d}\circ\mathbf{c})^{\ast}_{\gd}(\mathbf{C})$, respectively.
By Proposition~\ref{PGenHepDerivorComp}, for every $\sigma\in\Sigma_{\mathbf{s},s}$, $\Dom(\sigma^{\mathbf{A}})\subseteq\Dom(\sigma^{\mathbf{A}'})$, and both operations agree on $\Dom(\sigma^{\mathbf{A}})$.
We prove, by well-founded induction along $\lhd$, see Proposition~\ref{PWellFoundedInduction}, the reverse inclusion.
Let us assume that every $\tau\lhd\sigma$ has the same interpretation in $\mathbf{A}$ and in $\mathbf{A}'$, and let $\mathbf{a}\in\Dom(\sigma^{\mathbf{A}'})$.
Then, by Proposition~\ref{PTransformedAlg}, $\mathbf{A}'\models\mathrm{def}^{(S,\Sigma)}(\sigma)[\mathbf{a}]$, hence, by the induction hypothesis and Definition~\ref{DHepSig}, $\mathbf{A}\models\mathrm{def}^{(S,\Sigma)}(\sigma)[\mathbf{a}]$.
Since $\mathbf{A}$ is a $\sig{\Sigma}$-algebra, by Definition~\ref{DEpAlg}, $\mathbf{a}\in\Dom(\sigma^{\mathbf{A}})$.
Therefore $\mathbf{A}=\mathbf{A}'$, i.e., $\kappa^{\mathbf{c},\mathbf{d}}_{\mathbf{C}}$ is the identity.
If $\mathbf{C}'$ satisfies the same hypotheses, then $\kappa^{\mathbf{c},\mathbf{d}}_{\mathbf{C}'}$ is also the identity and, by Proposition~\ref{PGenHepDerivorComp}, $(\mathbf{d}\circ\mathbf{c})^{\ast}_{\gd}(f)=\mathbf{c}^{\ast}_{\gd}(\mathbf{d}^{\ast}_{\gd}(f))$. 
This completes the proof.
\end{proof}

By Proposition~\ref{PGenHepDerivorComp}, the hypothesis on the cells in Corollary~\ref{CGenHepDerivorCompEqual} amounts to the following condition: for every $s\in S$ and every cell $x$ of sort $\varphi(s)$ of $\mathbf{C}$ relative to $\mathbf{d}$, if $\mathbf{C}\models d^{\sharp}[\Dom^{(\varphi,c)}_{s}][x]$, then $\mathbf{d}^{\ast}_{\gd}(\mathbf{C})\models\Dom^{(\varphi,c)}_{s}[x]$.
It trivially holds when the domains of $\mathbf{c}$ are empty.

The following example shows that the inclusion $\kappa^{\mathbf{c},\mathbf{d}}_{\mathbf{C}}$ of Proposition~\ref{PGenHepDerivorComp} can be strict, even for very small hep-signatures.

\begin{example}\label{EStrictComparison}\lean{equivalent}{StrictComparison.comparison_not_isIso,exists_noninvertible_comparison}
Let $\sig{\Sigma}$ be the hep-signature with a single sort $s$ and no operation symbols, $\sig{\Lambda}$ the hep-signature with a single sort $t$ and a single unary operation symbol $u\in\Lambda_{(t),t}$ with $\mathrm{def}^{(T,\Lambda)}(u)=\varnothing$, and $\sig{\Omega}$ the hep-signature with a single sort $w$ and two constants $z,o\in\Omega_{\lambda,w}$ with empty sets of domain conditions.
Let $\mathbf{c}=(\varphi,c,\Dom^{(\varphi,c)})$ be the generalized hep-derivor from $\sig{\Sigma}$ to $\sig{\Lambda}$ with $\varphi(s)=t$ and $\Dom^{(\varphi,c)}_{s}=\{(u(v^{t}_{0}),u(v^{t}_{0}))\}$, which selects the elements at which $u$ is defined, and let $\mathbf{d}=(\psi,d,\Dom^{(\psi,d)})$ be the generalized hep-derivor from $\sig{\Lambda}$ to $\sig{\Omega}$ with $\psi(t)=w$, $d(u)=o$ and $\Dom^{(\psi,d)}_{t}=\{(v^{w}_{0},z)\}$.
Let $\mathbf{C}$ be the $\Omega$-algebra with $C_{w}=\{0,1\}$, $z^{\mathbf{C}}=0$ and $o^{\mathbf{C}}=1$.
Then the only cell of $\mathbf{C}$ relative to $\mathbf{d}$ is $0$, and, by Proposition~\ref{PTransformedAlg}, $u$ is nowhere defined in $\mathbf{d}^{\ast}_{\gd}(\mathbf{C})$, since the value $1$ of $d(u)=o$ is not a cell; thus $\mathbf{d}^{\ast}_{\gd}(\mathbf{C})$ is a weak $\sig{\Lambda}$-algebra but not a $\sig{\Lambda}$-algebra.
Consequently, $\mathbf{c}^{\ast}_{\gd}(\mathbf{d}^{\ast}_{\gd}(\mathbf{C}))$ has no cells.
On the other hand, by Definition~\ref{DGenHepDerivorCompId}, $\Dom^{\mathbf{d}\circ\mathbf{c}}_{s}=\{(v^{w}_{0},z),(o,o)\}$, and $0$ satisfies it in $\mathbf{C}$, so that $(\mathbf{d}\circ\mathbf{c})^{\ast}_{\gd}(\mathbf{C})$ has the single cell $0$.
Hence $\kappa^{\mathbf{c},\mathbf{d}}_{\mathbf{C}}$ is the strict inclusion of $\varnothing$ into $\{0\}$.
This is why, in Section~\ref{SHepAlgCat}, the induced functors form a lax, and not a split, indexed category, and why the hypotheses of Corollary~\ref{CGenHepDerivorCompEqual} are needed.
\end{example}

In the following lemma we combine Corollary~\ref{CGenHepDerivorCompEqual} and Lemma~\ref{LGenHepDerivorSameCells} in the form in which they are used to compare the two composites of a square of induced functors.

\begin{lemma}\label{LGenHepDerivorCommSquare}\lean{weaker}{GeneralizedHepDerivor.squareNatIso}
Let $\sig{\Sigma}$, $\sig{\Lambda}$, $\sig{\Lambda}'$ and $\sig{\Omega}$ be hep-signatures, $\mathbf{c}$ a generalized hep-derivor from $\sig{\Sigma}$ to $\sig{\Lambda}$, $\mathbf{d}$ a generalized hep-derivor from $\sig{\Lambda}$ to $\sig{\Omega}$, $\mathbf{c}'$ a generalized hep-derivor from $\sig{\Sigma}$ to $\sig{\Lambda}'$ and $\mathbf{d}'$ a generalized hep-derivor from $\sig{\Lambda}'$ to $\sig{\Omega}$ such that $\mathbf{d}\circ\mathbf{c}$ and $\mathbf{d}'\circ\mathbf{c}'$ have the same underlying derivor.
Let $\mathbf{C}$ be a partial $\Omega$-algebra such that $\mathbf{c}^{\ast}_{\gd}(\mathbf{d}^{\ast}_{\gd}(\mathbf{C}))$, $\mathbf{c}'^{\ast}_{\gd}(\mathbf{d}'^{\ast}_{\gd}(\mathbf{C}))$, $(\mathbf{d}\circ\mathbf{c})^{\ast}_{\gd}(\mathbf{C})$ and $(\mathbf{d}'\circ\mathbf{c}')^{\ast}_{\gd}(\mathbf{C})$ have the same underlying $S$-sorted set and the first two are $\sig{\Sigma}$-algebras.
Then
$$
\mathbf{c}^{\ast}_{\gd}\left(\mathbf{d}^{\ast}_{\gd}(\mathbf{C})\right)=\mathbf{c}'^{\ast}_{\gd}\left(\mathbf{d}'^{\ast}_{\gd}(\mathbf{C})\right).
$$
Moreover, if $f\colon\mathbf{C}\mor\mathbf{C}'$ is a homomorphism of partial $\Omega$-algebras and $\mathbf{C}'$ satisfies the same hypotheses, then $\mathbf{c}^{\ast}_{\gd}(\mathbf{d}^{\ast}_{\gd}(f))=\mathbf{c}'^{\ast}_{\gd}(\mathbf{d}'^{\ast}_{\gd}(f))$.
\end{lemma}
\begin{proof}
By Corollary~\ref{CGenHepDerivorCompEqual}, applied to $\mathbf{c}$ and $\mathbf{d}$, and to $\mathbf{c}'$ and $\mathbf{d}'$, $\mathbf{c}^{\ast}_{\gd}(\mathbf{d}^{\ast}_{\gd}(\mathbf{C}))=(\mathbf{d}\circ\mathbf{c})^{\ast}_{\gd}(\mathbf{C})$ and $\mathbf{c}'^{\ast}_{\gd}(\mathbf{d}'^{\ast}_{\gd}(\mathbf{C}))=(\mathbf{d}'\circ\mathbf{c}')^{\ast}_{\gd}(\mathbf{C})$.
By Proposition~\ref{PTransformedAlg}, the underlying $S$-sorted sets of the last two partial $\Sigma$-algebras are the sets of cells of $\mathbf{C}$ relative to $\mathbf{d}\circ\mathbf{c}$ and to $\mathbf{d}'\circ\mathbf{c}'$, which, by hypothesis, are equal.
Hence, by Lemma~\ref{LGenHepDerivorSameCells}, $(\mathbf{d}\circ\mathbf{c})^{\ast}_{\gd}(\mathbf{C})=(\mathbf{d}'\circ\mathbf{c}')^{\ast}_{\gd}(\mathbf{C})$.
The statement on homomorphisms follows in the same way from the corresponding statements of Corollary~\ref{CGenHepDerivorCompEqual} and Lemma~\ref{LGenHepDerivorSameCells}.
\end{proof}

In particular, if $\mathbf{d}\circ\mathbf{c}=\mathbf{d}'\circ\mathbf{c}'$ in $\mathsf{HSig}_{\gd}$, the hypothesis on $\mathbf{C}$ amounts to the following: $\mathbf{c}^{\ast}_{\gd}(\mathbf{d}^{\ast}_{\gd}(\mathbf{C}))$ and $\mathbf{c}'^{\ast}_{\gd}(\mathbf{d}'^{\ast}_{\gd}(\mathbf{C}))$ are $\sig{\Sigma}$-algebras with the same underlying $S$-sorted set as $(\mathbf{d}\circ\mathbf{c})^{\ast}_{\gd}(\mathbf{C})$.
The general form, in which the two composites need only have the same underlying derivor and select the same cells of $\mathbf{C}$, is used in Section~\ref{SDsTrCompatGdUf}, where the domains of the two composites differ.

The situation of the general form of Lemma~\ref{LGenHepDerivorCommSquare} recurs in the applications.

\begin{remark}\label{RSameCellsPattern}\lean{weaker}{GeneralizedHepDerivor.sameCellsNatIso}
In Chapters~\ref{SDsTr} and~\ref{SFunctorsOmega}, it often happens that two generalized hep-derivors between the same hep-signatures, a composite and a single generalized hep-derivor, or two composites, have the same underlying derivor but different domains, so that they are different morphisms of $\mathsf{HSig}_{\gd}$.
On the higher-order categorial algebras, however, both domains select the same cells, so that, by Lemma~\ref{LGenHepDerivorSameCells}, together with Corollary~\ref{CGenHepDerivorCompEqual} or Lemma~\ref{LGenHepDerivorCommSquare}, both induce the same functor.
This is the reason why several equalities of functors hold without a corresponding equality in $\mathsf{HSig}_{\gd}$, see Remarks~\ref{RTrCompHomogeneous}, \ref{RTrDsIdHomogeneous}, \ref{RTrDsCompHomogeneous}, \ref{RCompatGdTrDerivors}, \ref{RCompatGdDsDerivors} and~\ref{RCompatGdTrOmegaDerivors}.
\end{remark}

\section{The category \texorpdfstring{$\mathsf{HAlg}_{\gd}$}{HAlggd}}\label{SHepAlgCat}

In~\cite{ClimentSoliveres2010b, Soliveres2000}, the category $\mathsf{Alg}_{\mathfrak{d}}$, whose objects are the pairs formed by a many-sorted signature and an algebra for it, is obtained by means of the Grothendieck construction, see Definition~\ref{DGrothendieck}, for the split indexed category that sends a many-sorted signature $\sig{\Sigma}$ to $\mathsf{Alg}(\Sigma)$ and a derivor to the functor that it induces.
Following this pattern, we next assemble the ep-algebras of all the hep-signatures into a single category.
However, by Proposition~\ref{PGenHepDerivorComp}, the composite $\mathbf{c}^{\ast}_{\gd}\circ\mathbf{d}^{\ast}_{\gd}$ is, in general, only included in $(\mathbf{d}\circ\mathbf{c})^{\ast}_{\gd}$, and the inclusions $\kappa^{\mathbf{c},\mathbf{d}}_{\mathbf{C}}$ need not be isomorphisms.
Therefore the induced functors do not constitute a split, nor even a pseudo, indexed category, but a normal lax indexed category, see Definition~\ref{DLaxIndexed}, provided that we take as fibres the categories of weak ep-algebras, see Definition~\ref{DWeakEpAlg}, on which, by Proposition~\ref{PGenHepDerivorId}, the identity generalized hep-derivors induce identity functors.
We begin by proving that the transformed algebras are weak ep-algebras.

\begin{lemma}\label{LTransformedWeakEpAlg}\lean{equivalent}{GeneralizedHepDerivor.transformed_isWeak}
Let $\sig{\Sigma}$ and $\sig{\Lambda}$ be hep-signatures and $\mathbf{c}$ a generalized hep-derivor from $\sig{\Sigma}$ to $\sig{\Lambda}$.
Then, for every partial $\Lambda$-algebra $\mathbf{B}$, $\mathbf{c}^{\ast}_{\gd}(\mathbf{B})$ is a weak $\sig{\Sigma}$-algebra.
Consequently, the restriction of $\mathbf{c}^{\ast}_{\gd}$ to $\mathsf{WEPAlg}(\sig{\Lambda})$ corestricts to $\mathsf{WEPAlg}(\sig{\Sigma})$, see Definition~\ref{DCorestriction}, and we denote the corestriction by
$$
\mathbf{c}^{\ast}_{\gd}\colon\mathsf{WEPAlg}(\sig{\Lambda})\mor\mathsf{WEPAlg}(\sig{\Sigma}).
$$
\end{lemma}
\begin{proof}
The first statement is part of Proposition~\ref{PTransformedAlg}.
The last statement follows from Definition~\ref{DCorestriction}, since $\mathsf{WEPAlg}(\sig{\Lambda})$ and $\mathsf{WEPAlg}(\sig{\Sigma})$ are full subcategories of $\mathsf{PAlg}(\Lambda)$ and of $\mathsf{PAlg}(\Sigma)$, respectively.
\end{proof}

In the following proposition we prove that the categories of weak ep-algebras, together with the functors induced by the generalized hep-derivors and the inclusions of Proposition~\ref{PGenHepDerivorComp}, constitute a normal lax indexed category.

\begin{proposition}\label{PHepLaxIndexed}\lean{weaker}{BundledHepSignature.weakAlgebras}
Let $\mathrm{WEPAlg}$ stand for the assignment that sends every hep-signature $\sig{\Sigma}$ to the category $\mathsf{WEPAlg}(\sig{\Sigma})$ and every generalized hep-derivor $\mathbf{c}$ from $\sig{\Sigma}$ to $\sig{\Lambda}$ to the functor $\mathbf{c}^{\ast}_{\gd}\colon\mathsf{WEPAlg}(\sig{\Lambda})\mor\mathsf{WEPAlg}(\sig{\Sigma})$ of Lemma~\ref{LTransformedWeakEpAlg}, and let $\kappa$ stand for the assignment that sends every pair of generalized hep-derivors $\mathbf{c}\colon\sig{\Sigma}\mor\sig{\Lambda}$ and $\mathbf{d}\colon\sig{\Lambda}\mor\sig{\Omega}$ to the family
$$
\kappa^{\mathbf{c},\mathbf{d}}=\left(\kappa^{\mathbf{c},\mathbf{d}}_{\mathbf{C}}\right)_{\mathbf{C}\in\mathrm{Ob}(\mathsf{WEPAlg}(\sig{\Omega}))}
$$
of the inclusions of Proposition~\ref{PGenHepDerivorComp}.
Then $(\mathrm{WEPAlg},\kappa)$ is a normal lax $\mathsf{HSig}_{\gd}$-indexed category, see Definition~\ref{DLaxIndexed}.
\end{proposition}
\begin{proof}
\textsf{(Identities)} By Proposition~\ref{PGenHepDerivorId}, for every hep-signature $\sig{\Sigma}$, \[(\id_{\sig{\Sigma}})^{\ast}_{\gd}=\Id_{\mathsf{WEPAlg}(\sig{\Sigma})}.\]

\textsf{(Naturality)} Let $\mathbf{c}\colon\sig{\Sigma}\mor\sig{\Lambda}$ and $\mathbf{d}\colon\sig{\Lambda}\mor\sig{\Omega}$ be generalized hep-derivors.
By Lemma~\ref{LTransformedWeakEpAlg}, for every weak $\sig{\Omega}$-algebra $\mathbf{C}$, both $\mathbf{c}^{\ast}_{\gd}(\mathbf{d}^{\ast}_{\gd}(\mathbf{C}))$ and $(\mathbf{d}\circ\mathbf{c})^{\ast}_{\gd}(\mathbf{C})$ are weak $\sig{\Sigma}$-algebras, so that, by Proposition~\ref{PGenHepDerivorComp}, $\kappa^{\mathbf{c},\mathbf{d}}_{\mathbf{C}}$ is a morphism of the full subcategory $\mathsf{WEPAlg}(\sig{\Sigma})$ of $\mathsf{PAlg}(\Sigma)$.
Moreover, again by Proposition~\ref{PGenHepDerivorComp}, $\kappa^{\mathbf{c},\mathbf{d}}$ is a natural transformation from $\mathbf{c}^{\ast}_{\gd}\circ\mathbf{d}^{\ast}_{\gd}$ to $(\mathbf{d}\circ\mathbf{c})^{\ast}_{\gd}$.

Before checking the conditions (L1) and (L2) of Definition~\ref{DLaxIndexed}, let us note that a homomorphism of partial algebras is determined by its domain, its codomain and its underlying many-sorted mapping, and that, by Proposition~\ref{PGenHepDerivorComp} and Proposition~\ref{PGenHepDerivorFunctor}, the underlying many-sorted mapping of every $\kappa^{\mathbf{c},\mathbf{d}}_{\mathbf{C}}$, and of its image under every functor $\mathbf{e}^{\ast}_{\gd}$, which is a birestriction of it, is an inclusion.

\textsf{(L1)} Let $\mathbf{c}\colon\sig{\Sigma}\mor\sig{\Lambda}$ be a generalized hep-derivor and $\mathbf{B}$ a weak $\sig{\Lambda}$-algebra.
By Proposition~\ref{PHepSigCat}, $\mathbf{c}\circ\id_{\sig{\Sigma}}=\mathbf{c}$ and, by Lemma~\ref{LTransformedWeakEpAlg} and Proposition~\ref{PGenHepDerivorId}, $(\id_{\sig{\Sigma}})^{\ast}_{\gd}(\mathbf{c}^{\ast}_{\gd}(\mathbf{B}))=\mathbf{c}^{\ast}_{\gd}(\mathbf{B})$.
Hence $\kappa^{\id_{\sig{\Sigma}},\mathbf{c}}_{\mathbf{B}}$ is the inclusion of the underlying $S$-sorted set of $\mathbf{c}^{\ast}_{\gd}(\mathbf{B})$ into itself, i.e., $\kappa^{\id_{\sig{\Sigma}},\mathbf{c}}_{\mathbf{B}}=\id_{\mathbf{c}^{\ast}_{\gd}(\mathbf{B})}$.
Similarly, $\id_{\sig{\Lambda}}\circ\mathbf{c}=\mathbf{c}$ and, by Proposition~\ref{PGenHepDerivorId}, $\mathbf{c}^{\ast}_{\gd}((\id_{\sig{\Lambda}})^{\ast}_{\gd}(\mathbf{B}))=\mathbf{c}^{\ast}_{\gd}(\mathbf{B})$, so that $\kappa^{\mathbf{c},\id_{\sig{\Lambda}}}_{\mathbf{B}}=\id_{\mathbf{c}^{\ast}_{\gd}(\mathbf{B})}$.
Therefore $\kappa^{\id_{\sig{\Sigma}},\mathbf{c}}=\id_{\mathbf{c}^{\ast}_{\gd}}=\kappa^{\mathbf{c},\id_{\sig{\Lambda}}}$.

\textsf{(L2)} Let $\mathbf{c}\colon\sig{\Sigma}\mor\sig{\Lambda}$, $\mathbf{d}\colon\sig{\Lambda}\mor\sig{\Omega}$ and $\mathbf{e}\colon\sig{\Omega}\mor\sig{\Xi}$ be generalized hep-derivors and $\mathbf{D}$ a weak $\sig{\Xi}$-algebra.
By Proposition~\ref{PHepSigCat}, $(\mathbf{e}\circ\mathbf{d})\circ\mathbf{c}=\mathbf{e}\circ(\mathbf{d}\circ\mathbf{c})$.
Therefore both
$$
\kappa^{\mathbf{c},\mathbf{e}\circ\mathbf{d}}_{\mathbf{D}}\circ\mathbf{c}^{\ast}_{\gd}\left(\kappa^{\mathbf{d},\mathbf{e}}_{\mathbf{D}}\right)
\qquad\text{and}\qquad
\kappa^{\mathbf{d}\circ\mathbf{c},\mathbf{e}}_{\mathbf{D}}\circ\kappa^{\mathbf{c},\mathbf{d}}_{\mathbf{e}^{\ast}_{\gd}(\mathbf{D})}
$$
are homomorphisms from $\mathbf{c}^{\ast}_{\gd}(\mathbf{d}^{\ast}_{\gd}(\mathbf{e}^{\ast}_{\gd}(\mathbf{D})))$ to $(\mathbf{e}\circ\mathbf{d}\circ\mathbf{c})^{\ast}_{\gd}(\mathbf{D})$ whose underlying $S$-sorted mappings are composites of inclusions, i.e., the inclusion of the underlying $S$-sorted set of the former into that of the latter.
Hence both homomorphisms are equal. This completes the proof.
\end{proof}

We are now in a position to define the category of the ep-algebras of all the hep-signatures.

\begin{definition}\label{DHepAlgCat}\lean{equivalent}{HepAlgebra}
We let $\mathsf{HAlg}_{\gd}$ stand for the full subcategory of
$$
\int^{\mathsf{HSig}_{\gd}}(\mathrm{WEPAlg},\kappa),
$$
see Definition~\ref{DLaxGrothendieck} and Proposition~\ref{PHepLaxIndexed}, whose objects are the ordered pairs $(\sig{\Sigma},\mathbf{A})$, where $\sig{\Sigma}$ is a hep-signature and $\mathbf{A}$ a $\sig{\Sigma}$-algebra.
Explicitly, $\mathsf{HAlg}_{\gd}$ is defined as follows:
\begin{enumerate}
\item its objects are the ordered pairs $(\sig{\Sigma},\mathbf{A})$, where $\sig{\Sigma}$ is a hep-signature and $\mathbf{A}$ a $\sig{\Sigma}$-algebra;
\item its morphisms from $(\sig{\Sigma},\mathbf{A})$ to $(\sig{\Lambda},\mathbf{B})$ are the ordered pairs $(\mathbf{c},f)$, where $\mathbf{c}$ is a generalized hep-derivor from $\sig{\Sigma}$ to $\sig{\Lambda}$ and $f\colon\mathbf{A}\mor\mathbf{c}^{\ast}_{\gd}(\mathbf{B})$ a homomorphism of partial $\Sigma$-algebras;
\item the composite of $(\mathbf{c},f)\colon(\sig{\Sigma},\mathbf{A})\mor(\sig{\Lambda},\mathbf{B})$ and $(\mathbf{d},g)\colon(\sig{\Lambda},\mathbf{B})\mor(\sig{\Omega},\mathbf{C})$ is given by
\[
(\mathbf{d},g)\circ(\mathbf{c},f)=
\left(\mathbf{d}\circ\mathbf{c},\kappa^{\mathbf{c},\mathbf{d}}_{\mathbf{C}}\circ\mathbf{c}^{\ast}_{\gd}(g)\circ f\right),
\]
see Definition~\ref{DGenHepDerivorCompId} and Proposition~\ref{PGenHepDerivorComp}; and
\item the identity at $(\sig{\Sigma},\mathbf{A})$ is $(\id_{\sig{\Sigma}},\id_{\mathbf{A}})$.
\end{enumerate}
\end{definition}

Let us note that, since every $\sig{\Sigma}$-algebra is a weak $\sig{\Sigma}$-algebra and $\mathsf{WEPAlg}(\sig{\Sigma})$ is a full subcategory of $\mathsf{PAlg}(\Sigma)$, the morphisms of $\int^{\mathsf{HSig}_{\gd}}(\mathrm{WEPAlg},\kappa)$ between objects of $\mathsf{HAlg}_{\gd}$ are exactly the pairs described in the second item of the above definition, so that both descriptions of $\mathsf{HAlg}_{\gd}$ agree.

\begin{proposition}\label{PHepAlgCat}\lean{equivalent}{HepAlgebra.projection}
$\mathsf{HAlg}_{\gd}$ is a category.
Moreover, the assignment that sends every object $(\sig{\Sigma},\mathbf{A})$ to $\sig{\Sigma}$ and every morphism $(\mathbf{c},f)$ to $\mathbf{c}$ is a functor from $\mathsf{HAlg}_{\gd}$ to $\mathsf{HSig}_{\gd}$.
\end{proposition}
\begin{proof}
By Propositions~\ref{PLaxGrothendieck} and~\ref{PHepLaxIndexed}, $\int^{\mathsf{HSig}_{\gd}}(\mathrm{WEPAlg},\kappa)$ is a category and $\pi_{(\mathrm{WEPAlg},\kappa)}$ is a functor from it to $\mathsf{HSig}_{\gd}$.
Therefore its full subcategory $\mathsf{HAlg}_{\gd}$ is a category, and the assignment of the statement, which is the restriction of $\pi_{(\mathrm{WEPAlg},\kappa)}$ to $\mathsf{HAlg}_{\gd}$, is a functor.

This completes the proof.
\end{proof}

\begin{remark}\label{RHepAlgCat}
Let us consider the many-sorted signatures as hep-signatures with empty sets of domain conditions, see Remark~\ref{REpAlg}, and the derivors $(\varphi,c)$ between them as the generalized hep-derivors $(\varphi,c,\varnothing)$.
For a many-sorted signature $(S,\Sigma)$, the weak $(S,\Sigma,\varnothing)$-algebras are all the partial $\Sigma$-algebras and the $(S,\Sigma,\varnothing)$-algebras are the $\Sigma$-algebras.
Moreover, since the partial term operations of a total algebra are total, for every derivor $(\varphi,c)$ from $(S,\Sigma)$ to $(T,\Lambda)$, the restriction of $(\varphi,c,\varnothing)^{\ast}_{\gd}$ to $\mathsf{Alg}(\Lambda)$ corestricts to $\mathsf{Alg}(\Sigma)$, and this corestriction is the usual functor that sends a $\Lambda$-algebra $\mathbf{B}$ to the $\Sigma$-algebra on $B_{\varphi}$ in which every $\sigma$ is interpreted as the term operation $c(\sigma)^{\mathbf{B}}$, see~\cite{ClimentSoliveres2010b} and~\cite{Soliveres2000}, Proposition~3.5.6; and, since the domains of these generalized hep-derivors are empty, by Corollary~\ref{CGenHepDerivorCompEqual} and the note after it, the inclusions $\kappa^{\mathbf{c},\mathbf{d}}_{\mathbf{C}}$ are identities when $\mathbf{C}$ is a total algebra.
Therefore the subcategory of $\mathsf{HAlg}_{\gd}$ whose objects are the pairs $((S,\Sigma,\varnothing),\mathbf{A})$, with $\mathbf{A}$ a $\Sigma$-algebra, and whose morphisms are the pairs $((\varphi,c,\varnothing),f)$, is the category $\mathsf{Alg}_{\mathfrak{d}}$ of~\cite{ClimentSoliveres2010b}.
\end{remark}

\chapter[Homogeneous presentation]{Homogeneous presentation of higher-order categorial algebras}\label{SHomogeneous}

We now turn to the higher-order categorial algebras, to which the general theory of Chapter~\ref{SGenHepDerivors} is applied in the rest of this work.
In this chapter, we recall the definitions of homogeneous category, homogeneous $2$-category, homogeneous $n$-category, and many-sorted $n$-categorial algebra, together with their $\omega$-dimensional counterparts.
These are the homogeneous presentations of the corresponding categorical structures, in which all the cells, whatever their dimension, form a single set and the dimension of a cell is read off from its sources and targets, as in~\cite{Street1987}.
Next we introduce the homogeneous $n$-categorial signature determined by a many-sorted signature and a hep-specification on it, and we show that the many-sorted $n$-categorial algebras are exactly the models of this hep-specification.
We end this chapter by showing that, in dimension $\omega$, the condition that every cell has a finite dimension cannot be imposed by $\QE$-equations on the homogeneous $\omega$-categorial signature.

We begin by recalling the notion of homogeneous category.

\begin{definition}\label{DHomogeneousCat}\lean{equivalent}{Homogeneous.CategoryStructure}
A \emph{homogeneous category} is an ordered pair $\mathsf{C} = (C,\xi)$ where $C$ is a set and $\xi = (\circ,\src,\tgt)$ a structure of category on $C$, i.e., $\circ\colon C\times C\dmor C$ is a partial operation of \emph{composition}, and $\src$, $\tgt$ are unary operations on $C$ which assign, to every morphism in $C$, its \emph{source} and its \emph{target}, respectively, subject to the following conditions:
\begin{enumerate}
\item[(C1)] For every $f\in C$,
\begin{align*}
\src(\src(f)) &= \src(f),&
 \src(\tgt(f)) &= \tgt(f),\\
\tgt(\src(f)) &= \src(f),&
 \tgt(\tgt(f)) &= \tgt(f).
\end{align*}
\item[(C2)] For every $f,g\in C$,
\begin{align*}
g\circ f &\mbox{ is defined}&
 \mbox{if, and only if,}&&
  \src(g)& = \tgt(f).
\end{align*}
\item[(C3)] For every $f,g\in C$, if $\src(g) = \tgt(f)$, then
\begin{align*}
\src(g\circ f) &= \src(f),
&
\tgt(g\circ f) &= \tgt(g).
\end{align*}
\item[(C4)] For every $f\in C$,
\begin{align*}
f\circ \src(f) &= f,
&
\tgt(f)\circ f &= f.
\end{align*}
\item[(C5)] For every $f,g,h\in C$, if $\src(h)=\tgt(g)$ and $\src(g)=\tgt(f)$, then
\begin{align*}
h\circ (g\circ f) = (h\circ g)\circ f.
\end{align*}
\end{enumerate}

We call the $f\in C$ such that $\src(f) = f$ or, equivalently, by~(C1), $\tgt(f) = f$, the \emph{identities} of $\mathsf{C}$.

Let $\mathsf{C}=(C,\xi)$ and $\mathsf{C}'=(C',\xi')$ be two homogeneous categories, and let us write $\xi'=(\circ',\src',\tgt')$.
A \emph{homogeneous functor} from $\mathsf{C}$ to $\mathsf{C}'$ is an ordered triple $(\mathsf{C},F,\mathsf{C}')$, denoted by
$F\colon \mathsf{C}\mor\mathsf{C}'$, where $F$ is a mapping from $C$ to $C'$ such that
\begin{itemize}
\item[(Ci)] If $f\in C$, then
\begin{align*}
F(\src(f))&=\src'(F(f)),
&
F(\tgt(f))&=\tgt'(F(f)).
\end{align*}
\item[(Cii)] If $f,g\in C$ are such that $\src(g)=\tgt(f)$, then
\begin{align*}
F(g\circ f)=F(g)\circ' F(f).
\end{align*}
\end{itemize}

We say that a homogeneous category $\mathsf{C}$ is \emph{discrete} if every $f\in C$ is an identity.

Homogeneous functors can be composed.
In fact, if $F$ is a homogeneous functor from $\mathsf{C}$ to $\mathsf{D}$ and $G$ is a homogeneous functor from $\mathsf{D}$ to $\mathsf{E}$, then $G\circ F$ is a homogeneous functor from $\mathsf{C}$ to $\mathsf{E}$, which we call the \emph{composite} of $F$ and $G$.
This composition is associative.
Moreover, for every homogeneous category $\mathsf{C}$, there is an \emph{identity} homogeneous functor $(\mathsf{C},\Id_{C},\mathsf{C})$, denoted by $\Id_{\mathsf{C}}\colon \mathsf{C}\mor\mathsf{C}$, where $\Id_{C}$ is the identity mapping at $C$, which is the neutral element for the composition.
We denote by $\mathsf{1Cat}$ the category consisting of homogeneous categories and homogeneous functors.
\end{definition}

Let us note that we write $\mathsf{1Cat}$, and not $\mathsf{Cat}$, for the category of homogeneous categories.
On the one hand, $\mathsf{1Cat}$ is the case $n=1$ of the categories $\mathsf{nCat}$ of homogeneous $n$-categories, see Definition~\ref{DHomogeneousnCat}, and it agrees with the notation $\mathsf{1CatAlg}(\sig{\Sigma})$ for the $1$-categorial $\sig{\Sigma}$-algebras, see Definition~\ref{DCatAlg}.
On the other hand, the symbol $\mathsf{Cat}$ is reserved for the category of all $\boldsymbol{\mathcal{U}}$-categories, see Chapter~\ref{SIntroduction}, which is the codomain of the indexed categories, see Definition~\ref{DSplitIndexed}, and whose objects, such as $\mathsf{nCatAlg}(\sig{\Sigma})$, are categories in the usual sense, given by a class of objects and sets of morphisms, and not homogeneous categories.
Thus $\mathsf{1Cat}$ and $\mathsf{Cat}$ are different categories.

We next recall the notion of homogeneous $2$-category, in which two structures of homogeneous category on the same set interact as in a $2$-category.

\begin{definition}\label{DHomogeneous2Cat}\lean{equivalent}{Homogeneous.HigherCategory}
A \emph{homogeneous $2$-category} is an ordered pair $\mathsf{C} = (C,(\xi_{k})_{k\in 2})$
where $C$ is a set and, for every $k\in 2$, $\xi_{k}= (\circ^{k},\src^{k},\tgt^{k})$ a structure of homogeneous category on $C$, such that the following conditions are satisfied:
\begin{enumerate}
\item[(2C1)] For every $f\in C$,
\begin{align*}
\src^{1}(\src^{0}(f))&=\src^{0}(f)
=\src^{0}(\src^{1}(f))=\src^{0}(\tgt^{1}(f)),\\
\tgt^{1}(\tgt^{0}(f))&=\tgt^{0}(f)
=\tgt^{0}(\tgt^{1}(f))=\tgt^{0}(\src^{1}(f)).
\end{align*}
\item[(2C2)] For every $f,g\in C$, if $\src^{0}(g)=\tgt^{0}(f)$, then
\begin{align*}
\src^{1}(g{\circ^{0}}f)&=\src^{1}(g){\circ^{0}}\src^{1}(f),
&
\tgt^{1}(g{\circ^{0}}f)&=\tgt^{1}(g){\circ^{0}}\tgt^{1}(f).
\end{align*}
\item[(2C3)] For every $f,f',g,g'\in C$, if
\begin{align*}
\src^{1}(g')&=\tgt^{1}(g);
&
\src^{1}(f')&=\tgt^{1}(f);\\
\src^{0}(g')&=\tgt^{0}(f');
&
\src^{0}(g)&=\tgt^{0}(f),
\end{align*}
then
\begin{align*}
(g'\circ^{0}f')\circ^{1}(g\circ^{0} f)
=
(g'\circ^{1}g)\circ^{0}(f'\circ^{1}f)
.
\end{align*}
\end{enumerate}

We identify the family $(\xi_{k})_{k\in 2}$ with the ordered pair $(\xi_{0},\xi_{1})$.

Let $\mathsf{C}=(C,(\xi_{k})_{k\in 2})$ and $\mathsf{C}'=(C',(\xi'_{k})_{k\in 2})$ be two homogeneous $2$-categories.
A \emph{homogeneous $2$-functor} from $\mathsf{C}$ to $\mathsf{C}'$ is an ordered triple $(\mathsf{C},F,\mathsf{C}')$, denoted by $F\colon \mathsf{C}\mor\mathsf{C}'$, where $F$ is a mapping from $C$ to $C'$ such that, for every $k\in 2$, $F$ is a homogeneous functor from $(C,\xi_{k})$ to $(C',\xi'_{k})$.

Homogeneous $2$-functors can be composed.
In fact, if $F$ is a homogeneous $2$-functor from $\mathsf{C}$ to $\mathsf{D}$ and $G$ is a homogeneous $2$-functor from $\mathsf{D}$ to $\mathsf{E}$, then $G\circ F$ is a homogeneous $2$-functor from $\mathsf{C}$ to $\mathsf{E}$, which we call the \emph{composite} of $F$ and $G$.
This composition is associative.
Moreover, for every homogeneous $2$-category $\mathsf{C}$, there is an \emph{identity} homogeneous $2$-functor $(\mathsf{C},\Id_{C},\mathsf{C})$, denoted by $\Id_{\mathsf{C}}\colon \mathsf{C}\mor\mathsf{C}$, where $\Id_{C}$ is the identity mapping at $C$, which is the neutral element for the composition.
We denote by $\mathsf{2Cat}$ the category consisting of homogeneous
$2$-categories and homogeneous $2$-functors.
\end{definition}

We next recall the notion of homogeneous $n$-category, for $n\in\mathbb{N}$, in which every pair of the $n$ structures of homogeneous category forms a homogeneous $2$-category.

\begin{definition}\label{DHomogeneousnCat}\lean{equivalent}{Homogeneous.NCategory,Homogeneous.HigherCategory}
A \emph{homogeneous $n$-category}, for $n\in\mathbb{N}$, is an ordered pair $\mathsf{C} = (C,(\xi_{k})_{k\in n})$ where $C$ is a set and, for every $k\in n$, $\xi_{k}$ a structure of homogeneous category on $C$ such that, for every $j,k\in n$ with $j<k$, the ordered pair $(C,(\xi_{i})_{i\in\{j,k\}})$ is a homogeneous $2$-category.
For $j<k$, we identify $(\xi_{i})_{i\in \{j,k\}}$ with $(\xi_{j},\xi_{k})$.
Thus a homogeneous $0$-category is just a set, a homogeneous $1$-category is a homogeneous category in the sense of Definition~\ref{DHomogeneousCat}, and a homogeneous $2$-category is one in the sense of Definition~\ref{DHomogeneous2Cat}.

For $k\in n$, we write $\xi_{k}=(\circ^{k},\src^{k},\tgt^{k})$, and we say that $f\in C$ is a \emph{$k$-cell} if $f=\src^{k}(f)$ or, equivalently, by~(C1) of Definition~\ref{DHomogeneousCat}, if $f=\tgt^{k}(f)$.
We denote by $C_{k}$ the set of $k$-cells of $\mathsf{C}$, i.e.,
$$
C_{k}=\{f\in C\mid f=\src^{k}(f)\}.
$$
Moreover, we let $C_{n}$ stand for $C$ and we call the elements of $C_{n}$ the \emph{$n$-cells} of $\mathsf{C}$.

Let $\mathsf{C}=(C,(\xi_{k})_{k\in n})$ and $\mathsf{C}'=(C',(\xi'_{k})_{k\in n})$ be two homogeneous $n$-categories.
A \emph{homogeneous $n$-functor} from $\mathsf{C}$ to $\mathsf{C}'$, denoted by
$F\colon \mathsf{C}\mor\mathsf{C}'$, is a mapping $F$ from $C$ to $C'$ such that, for every $k\in n$, $F$ is a homogeneous functor from $(C, \xi_{k})$ to $(C', \xi'_{k})$.

Homogeneous $n$-functors can be composed.
In fact, if $F$ is a homogeneous $n$-functor from $\mathsf{C}$ to $\mathsf{D}$ and $G$ is a homogeneous $n$-functor from $\mathsf{D}$ to $\mathsf{E}$, then $G\circ F$ is a homogeneous $n$-functor from $\mathsf{C}$ to $\mathsf{E}$, which we call the \emph{composite} of $F$ and $G$.
This composition is associative.
Moreover, for every homogeneous $n$-category $\mathsf{C}$, there is an \emph{identity} homogeneous $n$-functor $(\mathsf{C},\Id_{C},\mathsf{C})$, denoted by $\Id_{\mathsf{C}}\colon \mathsf{C}\mor\mathsf{C}$, where $\Id_{C}$ is the identity mapping at $C$, which is the neutral element for the composition.
We denote by $\mathsf{nCat}$ the category consisting of homogeneous $n$-categories and homogeneous $n$-functors.
\end{definition}

In the following definition we introduce the many-sorted $n$-categorial algebras, i.e., the many-sorted algebras whose carriers are homogeneous $n$-categories and whose operations are homogeneous $n$-functors.

\begin{definition}\label{DCatAlg}\lean{equivalent}{Homogeneous.NCategorialAlgebra,Homogeneous.CategorialAlgebra}
Let $S$ be a set of sorts and $n\in\mathbb{N}$.
We call the objects of the functor category $\mathsf{nCat}^{S}$ \emph{$S$-sorted homogeneous $n$-categories}.
Thus, an $S$-sorted homogeneous $n$-category is an $S$-indexed family $\mathsf{A}=(\mathsf{A}_{s})_{s\in S}$ such that, for every sort $s\in S$, $\mathsf{A}_{s}$ is a homogeneous $n$-category.
Given two $S$-sorted homogeneous $n$-categories
$\mathsf{A}$ and $\mathsf{B}$, we call the morphisms in $\mathsf{nCat}^{S}$ from $\mathsf{A}$ to
$\mathsf{B}$ \emph{$S$-sorted homogeneous $n$-functors}.
Thus, an $S$-sorted homogeneous $n$-functor from $\mathsf{A}$ to $\mathsf{B}$ is an $S$-indexed family $F=(F_{s})_{s\in S}$ such that, for every sort $s\in S$, $F_{s}$ is a homogeneous $n$-functor from $\mathsf{A}_{s}$ to $\mathsf{B}_{s}$.
We denote by $\mathrm{nFunc}(\mathsf{A},\mathsf{B})$ the set of all $S$-sorted homogeneous $n$-functors from $\mathsf{A}$ to $\mathsf{B}$.
We consider the expressions $F\in\mathrm{nFunc}(\mathsf{A},\mathsf{B})$ and $F\colon\mathsf{A}\mor\mathsf{B}$ as synonyms.

Let $\mathsf{A}$ be an $S$-sorted homogeneous $n$-category.
The $S^{\star}\times S$-sorted set of finitary operations on $\mathsf{A}$ is $(\mathrm{nFunc}(\mathsf{A}_{\mathbf{s}},\mathsf{A}_{s}))_{(\mathbf{s},s)\in S^{\star}\times S}$, where, for every $\mathbf{s}\in S^{\star}$, $\mathsf{A}_{\mathbf{s}}=\prod_{j\in\bb{\mathbf{s}}}\mathsf{A}_{s_{j}}$, with $\bb{\mathbf{s}}$ denoting the length of the word $\mathbf{s}$.
The finite products $\mathsf{A}_{\mathbf{s}}$ are taken componentwise: sources and targets are computed componentwise, and $(\mathsf{g}_{j})_{j}\circ^{k}(\mathsf{f}_{j})_{j}$ is defined if, and only if, every $\mathsf{g}_{j}\circ^{k}\mathsf{f}_{j}$ is defined, and then it is $(\mathsf{g}_{j}\circ^{k}\mathsf{f}_{j})_{j}$; for $\mathbf{s}=\lambda$, $\mathsf{A}_{\lambda}$ is a final homogeneous $n$-category.

Let $\sig{\Sigma}=(S,\Sigma)$ be a many-sorted signature.
A \emph{structure of many-sorted $n$-categorial $\Sigma$-algebra on} an $S$-sorted homogeneous $n$-category $\mathsf{A}$ is a family $F = (F_{\mathbf{s},s})_{(\mathbf{s},s)\in S^{\star}\times S}$ where, for every $(\mathbf{s},s)\in S^{\star}\times S$, $F_{\mathbf{s},s}$ is a mapping from $\Sigma_{\mathbf{s},s}$ to $\mathrm{nFunc}(\mathsf{A}_{\mathbf{s}},\mathsf{A}_{s})$.
If $\sigma\in \Sigma_{\lambda,s}$, then $F_{\lambda,s}(\sigma)$ picks out a cell of $\mathsf{A}_{s}$ which is a $k$-cell for every $k\in n$, since $\mathsf{A}_{\lambda}$ is a final homogeneous $n$-category.
An \emph{$S$-sorted homogeneous $n$-categorial $\Sigma$-algebra} or, simply, a \emph{many-sorted $n$-categorial $\sig{\Sigma}$-algebra}, is a pair $(\mathsf{A},F)$, abbreviated to $\calg{A}$, where $\mathsf{A}=(\mathsf{A}_{s})_{s\in S}$ is an $S$-sorted homogeneous $n$-category and $F$ a structure of many-sorted $n$-categorial $\Sigma$-algebra on $\mathsf{A}$.
For a pair $(\mathbf{s},s)\in S^{\star}\times S$ and an operation symbol $\sigma\in \Sigma_{\mathbf{s},s}$, in order to simplify the notation, the homogeneous $n$-functor $F_{\mathbf{s},s}(\sigma)$ from $\mathsf{A}_{\mathbf{s}}$ to $\mathsf{A}_{s}$ will be written as $F_{\sigma}$, or simply as $\sigma^{\calg{A}}$.

An \emph{$S$-sorted homogeneous $n$-categorial $\Sigma$-homomorphism} or, simply, a \emph{many-sorted $n$-categorial $\sig{\Sigma}$-ho\-mo\-mor\-phism}, from $\calg{A}$ to $\calg{B}$, where $\calg{B}=(\mathsf{B},G)$, is a triple $(\calg{A},H,\calg{B})$, abbreviated to $H\colon\calg{A}\mor\calg{B}$, where $H$ is an $S$-sorted homogeneous $n$-functor from $\mathsf{A}$ to $\mathsf{B}$ such that, for every $(\mathbf{s},s)\in S^{\star}\times S$, every $\sigma\in\Sigma_{\mathbf{s},s}$ and every family of cells $(\mathsf{a}_{j})_{j\in\bb{\mathbf{s}}}\in\mathsf{A}_{\mathbf{s}}$, of any dimension,
$$
H_{s}\left(\sigma^{\calg{A}}
\left(\left(\mathsf{a}_{j}
\right)_{j\in\bb{\mathbf{s}}}\right)\right)
=
\sigma^{\calg{B}}
\left(\left(
H_{s_{j}}\left(
\mathsf{a}_{j}
\right)
\right)_{j\in\bb{\mathbf{s}}}
\right).
$$
We denote by $\mathsf{nCatAlg}(\sig{\Sigma})$ the category of $n$-categorial $\sig{\Sigma}$-algebras and $n$-categorial $\sig{\Sigma}$-homomorphisms.
\end{definition}

We next state the $\omega$-dimensional counterparts of these notions, which are used in Chapters~\ref{SDsTr} and~\ref{SProjLimits}.

\begin{definition}\label{DHomogeneousOmegaCat}\lean{equivalent}{Homogeneous.OmegaCategory}
A \emph{homogeneous $\omega$-category} is an ordered pair $\mathsf{C} = (C,(\xi_{k})_{k\in\omega})$ where $C$ is a set and, for every $k\in\omega$, $\xi_{k}$ a structure of homogeneous category on $C$ such that, for every $j,k\in\omega$ with $j<k$, the ordered pair $(C,(\xi_{j},\xi_{k}))$ is a homogeneous $2$-category, and such that every $f\in C$ is a $k$-cell, i.e., $\src^{k}(f)=f$, for some $k\in\omega$.
By~(2C1), a $j$-cell is a $k$-cell for every $k\geq j$.

Let $\mathsf{C}=(C,(\xi_{k})_{k\in\omega})$ and $\mathsf{C}'=(C',(\xi'_{k})_{k\in\omega})$ be two homogeneous $\omega$-categories.
A \emph{homogeneous $\omega$-functor} from $\mathsf{C}$ to $\mathsf{C}'$ is a mapping $F\colon C\mor C'$ which is a homogeneous functor from $(C,\xi_{k})$ to $(C',\xi'_{k})$ for every $k\in\omega$.
We denote by $\omega\mathsf{Cat}$ the category of homogeneous $\omega$-categories and homogeneous $\omega$-functors.
\end{definition}

We next introduce the $\omega$-dimensional counterpart of Definition~\ref{DCatAlg}.

\begin{definition}\label{DCatAlgOmega}\lean{equivalent}{Homogeneous.OmegaCategorialAlgebra}
Let $\sig{\Sigma}=(S,\Sigma)$ be a many-sorted signature.
An \emph{$\omega$-categorial $\sig{\Sigma}$-algebra} is a pair $\calg{A}=(\mathsf{A},F)$ where $\mathsf{A}=(\mathsf{A}_{s})_{s\in S}$ is an $S$-sorted homogeneous $\omega$-category and $F$ assigns to every $(\mathbf{s},s)\in S^{\star}\times S$ and every $\sigma\in\Sigma_{\mathbf{s},s}$ a homogeneous $\omega$-functor $\sigma^{\calg{A}}\colon\mathsf{A}_{\mathbf{s}}\mor\mathsf{A}_{s}$.
Here $\mathsf{A}_{\mathbf{s}}$ is a homogeneous $\omega$-category: if each $\mathsf{a}_{j}$ is a $k_{j}$-cell, then $(\mathsf{a}_{j})_{j\in\bb{\mathbf{s}}}$ is a $k$-cell for $k=\max_{j}k_{j}$ (and, for $\mathbf{s}=\lambda$, the unique cell of $\mathsf{A}_{\lambda}$ is a $0$-cell).
An \emph{$\omega$-categorial $\sig{\Sigma}$-homomorphism} is an $S$-sorted homogeneous $\omega$-functor commuting with the operations, as in Definition~\ref{DCatAlg}.
We denote by $\omega\mathsf{CatAlg}(\sig{\Sigma})$ the category of $\omega$-categorial $\sig{\Sigma}$-algebras and $\omega$-categorial $\sig{\Sigma}$-homomorphisms.
\end{definition}

In the homogeneous presentation, an $n$-categorial $\sig{\Sigma}$-algebra is regarded as a partial algebra on a single $S$-sorted set of cells, for a signature that enlarges $\Sigma$ with sources, targets and compositions of every dimension below $n$.

\begin{definition}\label{DCatSig}\lean{equivalent}{Homogeneous.epSignature,Homogeneous.hepSignature}
Let $\Sigma$ be an $S$-sorted signature and $n\in\mathbb{N}$.
The \emph{homogeneous $n$-categorial signature determined by $\Sigma$}, denoted by $\Sigma^{(n)}$, is the $S$-sorted signature defined, for every $(\mathbf{s},s)\in S^{\star}\times S$, as follows:
$$
{\Sigma}^{(n)}_{\mathbf{s},s}=
\begin{cases}
{\Sigma}_{\mathbf{s},s},
& \text{if $\mathbf{s}\not\in \{(s),(s,s)\}$;}\\
{\Sigma}_{\mathbf{s},s}\amalg\{\src^{k}_{s}, \tgt^{k}_{s}\mid k\in n\},
& \text{if $\mathbf{s}= (s)$;}\\
{\Sigma}_{\mathbf{s},s}\amalg\{\circ^{k}_{s}\mid k\in n\},
& \text{if $\mathbf{s}= (s,s)$.}
\end{cases}
$$
That is, ${\Sigma}^{(n)}$ is the expansion of $\Sigma$ obtained by adding, for every $k\in n$ and every sort $s\in S$, two unary operation symbols $\src^{k}_{s}$ and $\tgt^{k}_{s}$, the \emph{$k$-source} and the \emph{$k$-target}, of arity $(s)$ and coarity $s$, and a binary operation symbol $\circ^{k}_{s}$, the \emph{$k$-composition}, of arity $(s,s)$ and coarity $s$.
We write $x\circ^{k}_{s}y$ for $\circ^{k}_{s}(x,y)$.

Moreover, we regard $\Sigma^{(n)}$ as the hep-signature, see Definition~\ref{DHepSig}, $\sig{\Sigma}^{(n)}=(S,\Sigma^{(n)},\mathrm{def}^{(S,\Sigma^{(n)})})$ in which the $k$-compositions are the only operation symbols with a nonempty set of domain conditions: for every $k\in n$ and every $s\in S$,
$$
\mathrm{def}^{(S,\Sigma^{(n)})}(\circ^{k}_{s})=\{(\src^{k}_{s}(v^{s}_{0}),\tgt^{k}_{s}(v^{s}_{1}))\},
$$
where $v^{s}_{0}$ and $v^{s}_{1}$ are the variables of $\vs{(s,s)}$ that correspond to the first and to the second argument, respectively, so that $v^{s}_{0}\circ^{k}_{s}v^{s}_{1}$ is defined if, and only if, the $k$-source of $v^{s}_{0}$ is equal to the $k$-target of $v^{s}_{1}$.
For every $\sigma\in\Sigma$, every $k\in n$ and every $s\in S$, $\mathrm{def}^{(S,\Sigma^{(n)})}(\sigma)=\mathrm{def}^{(S,\Sigma^{(n)})}(\src^{k}_{s})=\mathrm{def}^{(S,\Sigma^{(n)})}(\tgt^{k}_{s})=\varnothing$.
The relation $\lhd$ of $\sig{\Sigma}^{(n)}$, see Definition~\ref{DHepSig}, is well-founded, since the domain conditions of the compositions only involve sources and targets, which have empty sets of domain conditions; hence $\sig{\Sigma}^{(n)}$ is indeed a hep-signature.
\end{definition}

Let us note that $\Sigma^{(0)}=\Sigma$ and that, for every $n,m\in\mathbb{N}$ with $n\leq m$, $\Sigma^{(n)}\subseteq\Sigma^{(m)}$ and $\mathrm{def}^{(S,\Sigma^{(m)})}$ extends $\mathrm{def}^{(S,\Sigma^{(n)})}$.

In the following definition we introduce the $\omega$-dimensional counterpart of the homogeneous $n$-categorial signature.

\begin{definition}\label{DCatSigOmega}\lean{equivalent}{Homogeneous.epSignature}
Let $\Sigma$ be an $S$-sorted signature.
The \emph{homogeneous $\omega$-categorial signature determined by $\Sigma$} is the $S$-sorted signature $\Sigma^{(\omega)}=\bigcup_{n\in\mathbb{N}}\Sigma^{(n)}$, and $\sig{\Sigma}^{(\omega)}=(S,\Sigma^{(\omega)},\mathrm{def}^{(S,\Sigma^{(\omega)})})$ is the hep-signature in which $\mathrm{def}^{(S,\Sigma^{(\omega)})}$ is the union of the $\mathrm{def}^{(S,\Sigma^{(n)})}$, for $n\in\mathbb{N}$.
Thus $\Sigma^{(\omega)}$ has the operation symbols $\src^{k}_{s}$, $\tgt^{k}_{s}$ and $\circ^{k}_{s}$ for every $k\in\omega$, with the same domain conditions as in Definition~\ref{DCatSig}.
\end{definition}

We next present the $n$-categorial $\sig{\Sigma}$-algebras as a $\QE$-variety of partial $\Sigma^{(n)}$-algebras.
The axiom (A0) states that the operations of $\Sigma$ are total, the axioms (A$^{k}$1)--(A$^{k}$8), stated for every dimension $k\in n$, describe the $k$-dimensional structure, and the axioms (AB$^{l,k}$1)--(AB$^{l,k}$3), stated for every pair of dimensions $l<k$ in $n$, relate the $l$-dimensional and the $k$-dimensional structures.
Informally, these equations state that, for every sort $s\in S$, the sources, targets and compositions of sort $s$ of every dimension behave as in an $n$-category, and that the operations of $\Sigma$ are compatible with them.

\begin{definition}\label{DCatSpec}\lean{equivalent}{Homogeneous.specification}
Let $\sig{\Sigma}=(S,\Sigma)$ be a many-sorted signature and $n\in\mathbb{N}$.
We denote by $\spec{C}^{\sig{\Sigma}^{(n)}}=(\sig{\Sigma}^{(n)},\mathcal{C}^{\sig{\Sigma}^{(n)}})$ the ep-specification in which $\sig{\Sigma}^{(n)}$ is the homogeneous $n$-categorial hep-signature introduced in Definition~\ref{DCatSig} and $\mathcal{C}^{\sig{\Sigma}^{(n)}}$ is the set of contextual $\QE$-equations, see Definition~\ref{DEpSpec}, listed below, each of them with the variables that occur in it as context.
The variables are taken from the $S$-sorted set $V^{S}$ of Section~\ref{SVariables}.
In these axioms, a family of variables indexed by a word is a tuple of pairwise distinct variables of the indicated sorts, and not the finite $S$-sorted context $V^{S}_{\mathbf{s}}$ of Definition~\ref{DVariables}; in the axioms that involve two families, such as (A$^{k}$8), the two families are independent, i.e., all their variables are pairwise distinct.
It suffices to take one such choice of variables for every axiom, since every assignment of values to them is a valuation of its context.

For every $(\mathbf{s},s)\in S^{\star}\times S$, every operation symbol $\sigma\in \Sigma_{\mathbf{s},s}$, and every family $(x_{j})_{j\in \lvert\mathbf{s}\rvert}$ of pairwise distinct variables with $x_{j}\in V^{S}_{s_{j}}$, the operation $\sigma$ applied to the family $(x_{j})_{j\in \lvert\mathbf{s}\rvert}$ is always defined.
Formally,
\begin{align*}
\sigma((x_{j})_{j\in \lvert\mathbf{s}\rvert})\eeq\sigma((x_{j})_{j\in \lvert\mathbf{s}\rvert}).\label{ACatA0}\tag{A0}
\end{align*}

For every $k\in n$, the following equations describe the $k$-dimensional structure.

For every sort $s\in S$ and every variable $x\in V^{S}_{s}$, the $k$-source and $k$-target of $x$ are always defined.
Formally,
\begin{align*}
\src^{k}_{s}(x)&\eeq\src^{k}_{s}(x);
&
\tgt^{k}_{s}(x)&\eeq\tgt^{k}_{s}(x).\label{ACatAk1}\tag{A$^{k}$1}
\end{align*}

For every sort $s\in S$ and every variable $x\in V^{S}_{s}$, we have the following equations:
\begin{align*}
\src^{k}_{s}(\src^{k}_{s}(x))&\eeq\src^{k}_{s}(x);
&
\src^{k}_{s}(\tgt^{k}_{s}(x))&\eeq\tgt^{k}_{s}(x);\\
\tgt^{k}_{s}(\src^{k}_{s}(x))&\eeq\src^{k}_{s}(x);
&
\tgt^{k}_{s}(\tgt^{k}_{s}(x))&\eeq\tgt^{k}_{s}(x).\label{ACatAk2}\tag{A$^{k}$2}
\end{align*}
In other words, $\src^{k}_{s}$ and $\tgt^{k}_{s}$ are right zeros.
In particular, $\src^{k}_{s}$ and $\tgt^{k}_{s}$ are idempotent.

For every sort $s\in S$ and every pair of variables $x,y\in V^{S}_{s}$, $x\circ^{k}_{s} y$ is defined if, and only if, the $k$-target of $y$ is equal to the $k$-source of $x$.
Formally,
\begin{align*}
x\circ^{k}_{s} y\eeq x\circ^{k}_{s} y \;\;
&\longrightarrow \;\;
\src^{k}_{s}(x)\eeq\tgt^{k}_{s}(y);\\
\src^{k}_{s}(x)\eeq\tgt^{k}_{s}(y) \;\;
&\longrightarrow \;\;
x\circ^{k}_{s} y\eeq x\circ^{k}_{s}y.\label{ACatAk3}\tag{A$^{k}$3}
\end{align*}

For every sort $s\in S$ and every pair of variables $x,y\in V^{S}_{s}$, if $x\circ^{k}_{s} y$ is defined, then the $k$-source of $x\circ^{k}_{s}y$ is that of $y$ and the $k$-target of $x\circ^{k}_{s}y$ is that of $x$.
Formally,
\begin{align*}
x\circ^{k}_{s}y\eeq x\circ^{k}_{s}y \;\;
&\longrightarrow \;\;
\src^{k}_{s}(x\circ^{k}_{s}y)\eeq\src^{k}_{s}(y);\\
x\circ^{k}_{s}y\eeq x\circ^{k}_{s}y \;\;
&\longrightarrow \;\;
\tgt^{k}_{s}(x\circ^{k}_{s}y)\eeq\tgt^{k}_{s}(x).\label{ACatAk4}\tag{A$^{k}$4}
\end{align*}

For every sort $s\in S$ and every variable $x\in V^{S}_{s}$, the compositions $x\circ^{k}_{s}\src^{k}_{s}(x)$ and $\tgt^{k}_{s}(x)\circ^{k}_{s}x$ are always defined and are equal to $x$, i.e., $\src^{k}_{s}(x)$ is a right unit element for the $k$-composition with $x$ and $\tgt^{k}_{s}(x)$ is a left unit element for the $k$-composition with $x$.
Formally,
\begin{align*}
x\circ^{k}_{s}\src^{k}_{s}(x)&\eeq x;
&
\tgt^{k}_{s}(x)\circ^{k}_{s}x&\eeq x.\label{ACatAk5}\tag{A$^{k}$5}
\end{align*}

For every sort $s\in S$ and every triple of variables $x, y, z\in V^{S}_{s}$, if the $k$-compositions $x\circ^{k}_{s}y$ and $y\circ^{k}_{s}z$ are defined, then the $k$-compositions $x\circ^{k}_{s}(y\circ^{k}_{s}z)$ and $(x\circ^{k}_{s}y)\circ^{k}_{s}z$ are defined and they are equal, i.e., the $k$-composition, when defined, is associative.
Formally,
\begin{align*}
(x\circ^{k}_{s}y\eeq x\circ^{k}_{s}y)
\wedge
(y\circ^{k}_{s}z\eeq y\circ^{k}_{s}z)\;\;
&\longrightarrow \;\;
(x\circ^{k}_{s}y)\circ^{k}_{s}z\eeq x\circ^{k}_{s}(y\circ^{k}_{s}z).\label{ACatAk6}\tag{A$^{k}$6}
\end{align*}

For every $(\mathbf{s},s)\in S^{\star}\times S$, every operation symbol $\sigma\in \Sigma_{\mathbf{s},s}$, and every family $(x_{j})_{j\in \lvert\mathbf{s}\rvert}$ of pairwise distinct variables with $x_{j}\in V^{S}_{s_{j}}$, the $k$-source of $\sigma((x_{j})_{j\in \lvert\mathbf{s}\rvert})$ is equal to $\sigma$ applied to the family $((\src^{k}_{s_{j}}(x_{j}))_{j\in \lvert\mathbf{s}\rvert})$, and the $k$-target of $\sigma((x_{j})_{j\in \lvert\mathbf{s}\rvert})$ is equal to $\sigma$ applied to the family $((\tgt^{k}_{s_{j}}(x_{j}))_{j\in \lvert\mathbf{s}\rvert})$.
Formally,
\begin{align*}
\src^{k}_{s}(\sigma((x_{j})_{j\in \lvert\mathbf{s}\rvert}))&\eeq\sigma((\src^{k}_{s_{j}}(x_{j}))_{j\in \lvert\mathbf{s}\rvert});\\
\tgt^{k}_{s}(\sigma((x_{j})_{j\in \lvert\mathbf{s}\rvert}))&\eeq\sigma((\tgt^{k}_{s_{j}}(x_{j}))_{j\in \lvert\mathbf{s}\rvert}).\label{ACatAk7}\tag{A$^{k}$7}
\end{align*}

For every $(\mathbf{s},s)\in S^{\star}\times S$, every operation symbol $\sigma\in \Sigma_{\mathbf{s},s}$, and every pair of families $(x_{j})_{j\in \lvert\mathbf{s}\rvert}$ and $(y_{j})_{j\in \lvert\mathbf{s}\rvert}$ of variables with $x_{j},y_{j}\in V^{S}_{s_{j}}$, all of them pairwise distinct, if, for every $j\in \lvert\mathbf{s}\rvert$, the $k$-compositions $x_{j}\circ^{k}_{s_{j}}y_{j}$ are defined, then the $k$-composition $\sigma((x_{j})_{j\in \lvert\mathbf{s}\rvert})\circ^{k}_{s}\sigma((y_{j})_{j\in \lvert\mathbf{s}\rvert})$ is defined and it is equal to $\sigma$ applied to the family $(x_{j}\circ^{k}_{s_{j}}y_{j})_{j\in \lvert\mathbf{s}\rvert}$.
Formally,
\begin{multline*}
\textstyle
\bigwedge_{j\in \lvert\mathbf{s}\rvert}
(x_{j}\circ^{k}_{s_{j}} y_{j}\eeq x_{j}\circ^{k}_{s_{j}} y_{j})
\;\;\longrightarrow\;\;\\
\sigma((x_{j}\circ^{k}_{s_{j}}y_{j})_{j\in \lvert\mathbf{s}\rvert})
\eeq
\sigma((x_{j})_{j\in \lvert\mathbf{s}\rvert})\circ^{k}_{s}\sigma((y_{j})_{j\in \lvert\mathbf{s}\rvert}).
\label{ACatAk8}\tag{A$^{k}$8}
\end{multline*}

For every $l,k\in n$ with $l<k$, the following equations relate the $l$-dimensional and the $k$-dimensional structures.

For every sort $s\in S$ and every variable $x\in V^{S}_{s}$, the elements $\src^{k}_{s}(\src^{l}_{s}(x))$, $\src^{l}_{s}(\src^{k}_{s}(x))$ and $\src^{l}_{s}(\tgt^{k}_{s}(x))$ are always defined and are equal to $\src^{l}_{s}(x)$.
Similarly, the elements $\tgt^{k}_{s}(\tgt^{l}_{s}(x))$, $\tgt^{l}_{s}(\tgt^{k}_{s}(x))$, and $\tgt^{l}_{s}(\src^{k}_{s}(x))$ are always defined and are equal to $\tgt^{l}_{s}(x)$.
Formally,
\begin{align*}
\src^{k}_{s}(\src^{l}_{s}(x))&\eeq\src^{l}_{s}(x);
&
\tgt^{k}_{s}(\tgt^{l}_{s}(x))&\eeq\tgt^{l}_{s}(x);\\
\src^{l}_{s}(\src^{k}_{s}(x))&\eeq\src^{l}_{s}(x);
&
\tgt^{l}_{s}(\tgt^{k}_{s}(x))&\eeq\tgt^{l}_{s}(x);\\
\src^{l}_{s}(\tgt^{k}_{s}(x))&\eeq\src^{l}_{s}(x);
&
\tgt^{l}_{s}(\src^{k}_{s}(x))&\eeq\tgt^{l}_{s}(x).\label{ACatABlk1}\tag{AB$^{l,k}$1}
\end{align*}

For every sort $s\in S$ and every pair of variables $x,y\in V^{S}_{s}$, if $x\circ^{l}_{s}y$ is defined, then the $k$-source of $x\circ^{l}_{s}y$ is the $l$-composition of the $k$-source of $x$ with the $k$-source of $y$ and the $k$-target of $x\circ^{l}_{s}y$ is the $l$-composition of the $k$-target of $x$ with the $k$-target of $y$.
Formally,
\begin{align*}
x\circ^{l}_{s}y&\eeq x\circ^{l}_{s}y \;\;\longrightarrow\;\;
\src^{k}_{s}(x\circ^{l}_{s}y)\eeq\src^{k}_{s}(x)\circ^{l}_{s}\src^{k}_{s}(y);\\
x\circ^{l}_{s}y&\eeq x\circ^{l}_{s}y \;\;\longrightarrow\;\;
\tgt^{k}_{s}(x\circ^{l}_{s}y)\eeq\tgt^{k}_{s}(x)\circ^{l}_{s}\tgt^{k}_{s}(y).
\label{ACatABlk2}\tag{AB$^{l,k}$2}
\end{align*}

For every sort $s\in S$ and every quadruple of variables $x,y,z,t\in V^{S}_{s}$, if the $l$-compositions $x\circ^{l}_{s}y$ and $z\circ^{l}_{s}t$ are defined and the $k$-compositions $x\circ^{k}_{s}z$ and $y\circ^{k}_{s}t$ are defined, then the $k$-composition of $x\circ^{l}_{s}y$ with $z\circ^{l}_{s}t$ is equal to the $l$-composition of $x\circ^{k}_{s}z$ with $y\circ^{k}_{s}t$.
Formally,
\begin{multline*}
(x\circ^{l}_{s}y\eeq x\circ^{l}_{s}y)
\wedge
(z\circ^{l}_{s}t\eeq z\circ^{l}_{s}t)
\wedge
(x\circ^{k}_{s}z\eeq x\circ^{k}_{s}z)
\wedge
(y\circ^{k}_{s}t\eeq y\circ^{k}_{s}t)
\\
\;\;\longrightarrow\;\;
(x\circ^{l}_{s}y)\circ^{k}_{s}(z\circ^{l}_{s}t)\eeq
(x\circ^{k}_{s}z)\circ^{l}_{s}(y\circ^{k}_{s}t).
\label{ACatABlk3}\tag{AB$^{l,k}$3}
\end{multline*}

A model of the axioms (A$^{k}$1)--(A$^{k}$6), for $k\in n$, and (AB$^{l,k}$1)--(AB$^{l,k}$3), for $l<k$ in $n$, is an $S$-sorted homogeneous $n$-category, see Definition~\ref{DCatAlg}, presented with all the cells of a given sort, whatever their dimension, in a single set.

Since $\sig{\Sigma}^{(n)}$ is a hep-signature, see Definition~\ref{DCatSig}, $\spec{C}^{\sig{\Sigma}^{(n)}}$ is a hep-specification, see Definition~\ref{DEpSpec}, and we work with the ep-variety $\mathcal{V}(\spec{C}^{\sig{\Sigma}^{(n)}})$ and with its canonically associated category $\mathsf{EPAlg}(\spec{C}^{\sig{\Sigma}^{(n)}})$, see Definition~\ref{DEpVariety}.
\end{definition}

For $n=\omega$, we let $\spec{C}^{\sig{\Sigma}^{(\omega)}}=(\sig{\Sigma}^{(\omega)},\mathcal{C}^{\sig{\Sigma}^{(\omega)}})$ stand for the hep-specification, see Definitions~\ref{DCatSigOmega} and~\ref{DEpSpec}, in which $\mathcal{C}^{\sig{\Sigma}^{(\omega)}}$ consists of the contextual $\QE$-equations (A0), (A$^{k}$1)--(A$^{k}$8), for every $k\in\omega$, and (AB$^{l,k}$1)--(AB$^{l,k}$3), for every $l<k$ in $\omega$.

We next relate the many-sorted $n$-categorial $\sig{\Sigma}$-algebras with the $\sig{\Sigma}^{(n)}$-algebras.
Let $\sig{\Sigma}=(S,\Sigma)$ be a many-sorted signature and $n\in\mathbb{N}$.
In what follows, a pair $(\mathsf{A},F)$, where $\mathsf{A}=(\mathsf{A}_{s})_{s\in S}$ is an $S$-sorted homogeneous $n$-category, with $\mathsf{A}_{s}=(A_{s},(\circ^{k}_{s},\src^{k}_{s},\tgt^{k}_{s})_{k\in n})$, and $F$ assigns to every $(\mathbf{s},s)\in S^{\star}\times S$ and every $\sigma\in\Sigma_{\mathbf{s},s}$ a mapping $F_{\sigma}\colon A_{\mathbf{s}}\mor A_{s}$, is identified with the partial $\Sigma^{(n)}$-algebra with underlying $S$-sorted set $A=(A_{s})_{s\in S}$ in which, for every $k\in n$ and every $s\in S$, the operation symbols $\src^{k}_{s}$, $\tgt^{k}_{s}$ and $\circ^{k}_{s}$ are interpreted as $\src^{k}_{s}$, $\tgt^{k}_{s}$ and $\circ^{k}_{s}$, respectively, and every $\sigma\in\Sigma$ is interpreted as $F_{\sigma}$.

\begin{remark}\label{RCatAlgAsEpAlg}\lean{equivalent}{Homogeneous.modelIsomorphism}
Let $\sig{\Sigma}=(S,\Sigma)$ be a many-sorted signature and $n\in\mathbb{N}$.
Then, under the above identification,
$$
\mathsf{nCatAlg}(\sig{\Sigma})=\mathsf{EPAlg}(\spec{C}^{\sig{\Sigma}^{(n)}}).
$$
In fact, for every $k\in n$, the conditions (C1)--(C5) of Definition~\ref{DHomogeneousCat} for the $k$-th structure of category are the axioms (A$^{k}$2)--(A$^{k}$6), the totality of $\src^{k}_{s}$ and $\tgt^{k}_{s}$ being the axiom (A$^{k}$1); for every $l<k$ in $n$, the conditions (2C1)--(2C3) of Definition~\ref{DHomogeneous2Cat} are the axioms (AB$^{l,k}$1)--(AB$^{l,k}$3); and, since the products $\mathsf{A}_{\mathbf{s}}$ are taken componentwise, the operations of $\Sigma$ are total homogeneous $n$-functors if, and only if, the axioms (A0), (A$^{k}$7) and (A$^{k}$8) hold, for every $k\in n$.
Moreover, the $S$-sorted homogeneous $n$-functors that commute with the operations of $\Sigma$, i.e., the $n$-categorial $\sig{\Sigma}$-homomorphisms, are exactly the homomorphisms of partial $\Sigma^{(n)}$-algebras, see Definition~\ref{DPartialAlg}, since the latter preserve the defined $k$-compositions.
Finally, every model of $\spec{C}^{\sig{\Sigma}^{(n)}}$, i.e., every partial $\Sigma^{(n)}$-algebra in which the axioms of $\mathcal{C}^{\sig{\Sigma}^{(n)}}$ hold, is a $\sig{\Sigma}^{(n)}$-algebra, see Definitions~\ref{DEpAlg} and~\ref{DCatSig}: the operation symbols of $\Sigma$ and the $k$-sources and $k$-targets have empty sets of domain conditions, and they are total, by (A0) and (A$^{k}$1); and, by (A$^{k}$3), $x\circ^{k}_{s}y$ is defined if, and only if, $(x,y)$ satisfies the domain condition $\src^{k}_{s}(v^{s}_{0})\eeq\tgt^{k}_{s}(v^{s}_{1})$ of $\circ^{k}_{s}$.
Hence the objects of $\mathsf{EPAlg}(\spec{C}^{\sig{\Sigma}^{(n)}})$ are exactly the models of $\spec{C}^{\sig{\Sigma}^{(n)}}$.
Let us note that this identification relies on the contexts of Definition~\ref{DEpSpec}: every axiom of $\mathcal{C}^{\sig{\Sigma}^{(n)}}$ only involves variables of the sorts that it constrains, so that it holds in a partial $\Sigma^{(n)}$-algebra $\mathbf{A}$ if, and only if, the corresponding condition holds for the elements of these sorts, also when some carriers of $\mathbf{A}$ are empty.
For instance, for $n=1$ and $S=\{s,t\}$, a partial $\Sigma^{(1)}$-algebra with $A_{t}=\varnothing$ in which $\src^{0}_{s}$ is not idempotent is not a model of $\spec{C}^{\sig{\Sigma}^{(1)}}$, since (A$^{0}$2), whose context consists of a single variable $x$ of sort $s$, fails at some element of $A_{s}$.
\end{remark}

In dimension $\omega$, a pair $(\mathsf{A},F)$, where $\mathsf{A}$ is an $S$-sorted homogeneous $\omega$-category, is identified in the same way, with $\omega$ in place of $n$, with a partial $\Sigma^{(\omega)}$-algebra, and the identification works as in Remark~\ref{RCatAlgAsEpAlg}, except for the condition that every cell has a finite dimension.

\begin{remark}\label{RCatAlgAsEpAlgOmega}\lean{equivalent}{Homogeneous.omegaModelIsomorphism}
Let $\sig{\Sigma}=(S,\Sigma)$ be a many-sorted signature.
Then, under the above identification, and by the same arguments as in Remark~\ref{RCatAlgAsEpAlg}, applied for every $k\in\omega$ and every $l<k$ in $\omega$, $\omega\mathsf{CatAlg}(\sig{\Sigma})$ is the full subcategory of $\mathsf{EPAlg}(\spec{C}^{\sig{\Sigma}^{(\omega)}})$ whose objects are the $\sig{\Sigma}^{(\omega)}$-algebras in which, for every $s\in S$, every cell $x$ of sort $s$ is a $k$-cell, i.e., satisfies $\src^{k}_{s}(x)=x$, for some $k\in\omega$, see Definitions~\ref{DHomogeneousOmegaCat} and~\ref{DCatAlgOmega}.
\end{remark}

One could try to add to $\mathcal{C}^{\sig{\Sigma}^{(\omega)}}$ some axioms expressing that every cell has a finite dimension, in order to present $\omega\mathsf{CatAlg}(\sig{\Sigma})$ as the category canonically associated with an ep-variety, as in the finite-dimensional case.
In the following proposition we prove that this is impossible.

\begin{proposition}\label{PFiniteDimNotQE}\lean{equivalent}{Homogeneous.not_qe_axiomatizable}
Let $\sig{\Sigma}=(S,\Sigma)$ be a many-sorted signature with $S\neq\varnothing$.
Then there is no set $\mathcal{F}$ of $\QE$-equations for $\Sigma^{(\omega)}$ such that the objects of $\omega\mathsf{CatAlg}(\sig{\Sigma})$, under the above identification, are exactly the partial $\Sigma^{(\omega)}$-algebras in which every $\QE$-equation of $\mathcal{F}$ holds.
\end{proposition}
\begin{proof}
Let us assume, towards a contradiction, that such a set $\mathcal{F}$ exists.
Then, by the note after Definition~\ref{DQEEq}, the class of the objects of $\omega\mathsf{CatAlg}(\sig{\Sigma})$ is closed under products of partial $\Sigma^{(\omega)}$-algebras, see Definition~\ref{DPartialAlgProd}.
We next construct a family $(\calg{A}_{k})_{k\in\omega}$ of $\omega$-categorial $\sig{\Sigma}$-algebras whose product is not an $\omega$-categorial $\sig{\Sigma}$-algebra.

Let $s_{0}\in S$.
For every $k\in\omega$, let $\mathsf{G}_{k}$ be a homogeneous $\omega$-category with a cell $g_{k}$ which is a $(k+1)$-cell but not a $k$-cell and such that every cell of $\mathsf{G}_{k}$ is a $(k+1)$-cell; for instance, the homogeneous $\omega$-category freely generated by a $(k+1)$-dimensional globe, whose cells are $g_{k}$ and its iterated $j$-sources and $j$-targets, for $j\leq k$.
Let $\calg{A}_{k}$ be the $\omega$-categorial $\sig{\Sigma}$-algebra in which, for every $s\in S$, the homogeneous $\omega$-category of sort $s$ is $\mathsf{G}_{k}$, and in which every operation symbol $\sigma\in\Sigma_{\mathbf{s},s}$ is interpreted as the constant mapping with value the $0$-cell $u_{k}=\src^{0}(g_{k})$.
Every such constant mapping is a homogeneous $\omega$-functor, since $u_{k}$ is a $j$-cell for every $j\in\omega$, so that $\src^{j}(u_{k})=u_{k}=\tgt^{j}(u_{k})$ and, by~(C4), $u_{k}\circ^{j}u_{k}=u_{k}$.
Moreover, every cell of $\calg{A}_{k}$ is a $(k+1)$-cell, so that $\calg{A}_{k}$ is an object of $\omega\mathsf{CatAlg}(\sig{\Sigma})$, see Remark~\ref{RCatAlgAsEpAlgOmega}.

Let $\mathbf{A}=\prod_{k\in\omega}\calg{A}_{k}$ be the product of the corresponding partial $\Sigma^{(\omega)}$-algebras and $g=(g_{k})_{k\in\omega}$, which is an element of $A_{s_{0}}$.
By Definition~\ref{DPartialAlgProd}, for every $j\in\omega$, $\src^{j\,\mathbf{A}}_{s_{0}}(g)=(\src^{j}(g_{k}))_{k\in\omega}$, which is different from $g$, since $\src^{j}(g_{j})\neq g_{j}$.
Hence $g$ is not a $j$-cell of $\mathbf{A}$ for any $j\in\omega$, so that $\mathbf{A}$ is not an object of $\omega\mathsf{CatAlg}(\sig{\Sigma})$, a contradiction. This completes the proof.
\end{proof}

From the proof of the above proposition we obtain the following corollary.

\begin{corollary}\label{CProperSubCatAlgOmega}\lean{equivalent}{Homogeneous.omegaModels_proper}
Let $\sig{\Sigma}=(S,\Sigma)$ be a many-sorted signature with $S\neq\varnothing$.
Then $\omega\mathsf{CatAlg}(\sig{\Sigma})$ is a proper full subcategory of $\mathsf{EPAlg}(\spec{C}^{\sig{\Sigma}^{(\omega)}})$.
\end{corollary}
\begin{proof}
By Remark~\ref{RCatAlgAsEpAlgOmega}, $\omega\mathsf{CatAlg}(\sig{\Sigma})$ is a full subcategory of $\mathsf{EPAlg}(\spec{C}^{\sig{\Sigma}^{(\omega)}})$.
Let $\mathbf{A}$ be the product constructed in the proof of Proposition~\ref{PFiniteDimNotQE}.
Since $\mathcal{C}^{\sig{\Sigma}^{(\omega)}}$ is a set of $\QE$-equations, by the note after Definition~\ref{DQEEq}, $\mathbf{A}$ is a model of $\spec{C}^{\sig{\Sigma}^{(\omega)}}$, hence, by the argument at the end of Remark~\ref{RCatAlgAsEpAlg}, a $\sig{\Sigma}^{(\omega)}$-algebra, i.e., an object of $\mathsf{EPAlg}(\spec{C}^{\sig{\Sigma}^{(\omega)}})$ which, by the proof of Proposition~\ref{PFiniteDimNotQE}, is not an object of $\omega\mathsf{CatAlg}(\sig{\Sigma})$.
\end{proof}

Let us note that the same argument shows that the condition that every cell has a finite dimension cannot be expressed by any set of axioms whose class of models is closed under products.
The products involved are the products of partial $\Sigma^{(\omega)}$-algebras, see Definition~\ref{DPartialAlgProd}, i.e., the ambient products of the fixed homogeneous presentation.
The obstruction is therefore a property of this presentation: it does not mean that $\omega\mathsf{CatAlg}(\sig{\Sigma})$ lacks products, and, by Remark~\ref{RCatAlgAsEpAlgHtg}, the heterogeneous presentation of the same structures is the category canonically associated with an ep-variety.
Thus, when $S\neq\varnothing$, the fixed homogeneous signature does not admit additional $\QE$-axioms imposing finite-dimensional cells.
In Chapter~\ref{SProjLimits} we describe both presentations by projective limits of their finite-dimensional categories.

\chapter[Heterogeneous presentation]{Heterogeneous presentation of higher-order categorial algebras}\label{SHeterogeneous}

In the heterogeneous presentation, cells of different dimensions occupy separate sorts.
We follow the same pattern as in Chapter~\ref{SHomogeneous}: after recalling the heterogeneous $\omega$-categories and the heterogeneous $n$-categories, we introduce the heterogeneous $n$-categorial algebras, the heterogeneous $n$-categorial signature determined by a many-sorted signature, and a hep-specification on this signature.
We end this chapter by showing that the heterogeneous $n$-categorial algebras are exactly the models of this hep-specification and that, in contrast with the homogeneous presentation, the same holds in dimension $\omega$.

We begin by defining heterogeneous $\omega$-categories, which present $\omega$-categories with separate sorts for each dimension, together with sources, targets and identities between dimensions, as in~\cite{AraEtAl2025}.
A heterogeneous $n$-category is then a heterogeneous $\omega$-category truncated at level $n$.

\begin{definition}\label{DHeterogeneousOmegaCat}\lean{equivalent}{Heterogeneous.OmegaCategory,Heterogeneous.HigherCategory}
A \emph{heterogeneous $\omega$-category} is an ordered pair $\mathsf{C} = ((C_{k})_{k\in\omega},\zeta)$ consisting of an $\omega$-sorted set $(C_{k})_{k\in\omega}$ and a $\coprod_{k\in \omega}k$-indexed family
\[
\zeta = (\zeta_{j,k})_{(j,k)\in \coprod_{k\in \omega}k},
\]
where, for every $k\in \omega$ and every $j\in k$, $\zeta_{j,k} = (\circ^{(j,k)},\src^{(j,k)},\tgt^{(j,k)},\mathrm{i}^{(k,j)})$ with
\begin{enumerate}
\item $\src^{(j,k)}$ and $\tgt^{(j,k)}$ mappings from $C_{k}$ to $C_{j}$
which send $f\in C_{k}$ to its \emph{$(j,k)$-source} and \emph{$(j,k)$-target}, respectively;
\item $\mathrm{i}^{(k,j)}$ a mapping from $C_{j}$ to $C_{k}$
which sends $f\in C_{j}$ to its \emph{$(k,j)$-identity}; and
\item $\circ^{(j,k)}$ a partial binary operation
$$
\circ^{(j,k)}\colon C_{k}\times C_{k} \dmor C_{k},
$$
whose domain of definition is
$$
\Dom\left(\circ^{(j,k)}\right) =
\left\lbrace
(g,f)\in C_{k}\times C_{k}
\,\middle|\, \src^{(j,k)}(g)=\tgt^{(j,k)}(f) \right\rbrace,
$$
i.e., the fibered product of $C_{k}$ and $C_{k}$ over $\src^{(j,k)}$ and $\tgt^{(j,k)}$, and
which sends a pair $(g,f)\in \Dom(\circ^{(j,k)})$ to its \emph{$j$-composition of $k$-cells} $g\circ^{(j,k)}f$.
\end{enumerate}
These data are required to satisfy the following conditions for all $i,j,k\in \omega$ with $i<j<k$; the conditions (MS2), (MS6) and (MS7) only involve the two dimensions $j<k$, and the condition (MS5) only the two dimensions $i<j$, and these four conditions are required for every such pair of dimensions:
\begin{itemize}
\item[(MS1)] For every $f\in C_{k}$,
\begin{align*}
\src^{(i,j)}(\src^{(j,k)}(f))&=\src^{(i,k)}(f)=\src^{(i,j)}(\tgt^{(j,k)}(f)),\\
\tgt^{(i,j)}(\tgt^{(j,k)}(f))&=\tgt^{(i,k)}(f)=\tgt^{(i,j)}(\src^{(j,k)}(f)).
\end{align*}
\item[(MS2)] If $g,f\in C_{k}$ are such that $\src^{(j,k)}(g)=\tgt^{(j,k)}(f)$, then
\begin{align*}
\src^{(j,k)}(g\circ^{(j,k)}f)&=\src^{(j,k)}(f), &
\tgt^{(j,k)}(g\circ^{(j,k)}f)&=\tgt^{(j,k)}(g).
\end{align*}
\item[(MS3)] If $g,f\in C_{k}$ are such that $\src^{(i,k)}(g)=\tgt^{(i,k)}(f)$, then
\begin{align*}
\src^{(j,k)}(g\circ^{(i,k)}f)&=\src^{(j,k)}(g)\circ^{(i,j)}\src^{(j,k)}(f),\\
\tgt^{(j,k)}(g\circ^{(i,k)}f)&=\tgt^{(j,k)}(g)\circ^{(i,j)}\tgt^{(j,k)}(f).
\end{align*}
\item[(MS4)] If $f\in C_{i}$, then
$$
\mathrm{i}^{(k,j)}(\mathrm{i}^{(j,i)}(f))=\mathrm{i}^{(k,i)}(f).
$$
\item[(MS5)] If $f\in C_{i}$, then
$$
\src^{(i,j)}(\mathrm{i}^{(j,i)}(f))=f=
\tgt^{(i,j)}(\mathrm{i}^{(j,i)}(f)).
$$
\item[(MS6)] If $f\in C_{k}$, then
$$
f\circ^{(j,k)}
(\mathrm{i}^{(k,j)}(\src^{(j,k)}(f)))
=
f
=
(\mathrm{i}^{(k,j)}(\tgt^{(j,k)}(f)))
\circ^{(j,k)}
f.
$$
\item[(MS7)] If $h,g,f\in C_{k}$ are such that
\begin{align*}
\src^{(j,k)}(h)&=\tgt^{(j,k)}(g); &
\src^{(j,k)}(g)&=\tgt^{(j,k)}(f),
\end{align*}
then
$$
h\circ^{(j,k)}(g\circ^{(j,k)}f)=
(h\circ^{(j,k)}g)\circ^{(j,k)}f.
$$
\item[(MS8)] If $g,f\in C_{j}$ are such that $\src^{(i,j)}(g)=\tgt^{(i,j)}(f)$, then
$$
\mathrm{i}^{(k,j)}(g\circ^{(i,j)}f)=\mathrm{i}^{(k,j)}(g)\circ^{(i,k)}\mathrm{i}^{(k,j)}(f).
$$
\item[(MS9)] If $g',g,f',f\in C_{k}$ are such that
\begin{align*}
\src^{(j,k)}(g')&=\tgt^{(j,k)}(g); &
\src^{(j,k)}(f')&=\tgt^{(j,k)}(f);\\
\src^{(i,k)}(g')&=\tgt^{(i,k)}(f'); &
\src^{(i,k)}(g)&=\tgt^{(i,k)}(f),
\end{align*}
then
$$
(g'\circ^{(i,k)}f')\circ^{(j,k)}(g\circ^{(i,k)} f)
=
(g'\circ^{(j,k)}g)\circ^{(i,k)}(f'\circ^{(j,k)}f)
.
$$
\end{itemize}

Let $\mathsf{C}=((C_{k})_{k\in\omega},\zeta)$ and $\mathsf{C}'=((C'_{k})_{k\in\omega},\zeta')$ be two heterogeneous $\omega$-categories, and let us write $\zeta'_{j,k}=(\circ'^{(j,k)},\src'^{(j,k)},\tgt'^{(j,k)},\mathrm{i}'^{(k,j)})$.
A \emph{heterogeneous $\omega$-functor} from $\mathsf{C}$ to $\mathsf{C}'$ is a family $F=(F_{k})_{k\in\omega}$ of mappings $F_{k}\colon C_{k}\mor C'_{k}$ such that, for all $j,k\in\omega$ with $j<k$,
\begin{align*}
F_{j}\circ \src^{(j,k)} &= \src'^{(j,k)}\circ F_{k}, &
F_{j}\circ \tgt^{(j,k)} &= \tgt'^{(j,k)}\circ F_{k},\\
F_{k}\circ \mathrm{i}^{(k,j)} &= \mathrm{i}'^{(k,j)}\circ F_{j}, &
F_{k}(g\circ^{(j,k)}f) &= F_{k}(g)\circ'^{(j,k)} F_{k}(f),
\end{align*}
the last one for every $(g,f)\in\Dom(\circ^{(j,k)})$.
We denote by $\omega\mathsf{Cat}^{\htg}$ the category of heterogeneous $\omega$-categories and heterogeneous $\omega$-functors.

For $n\in\mathbb{N}$, a \emph{heterogeneous $n$-category} is an ordered pair $\mathsf{C} = ((C_{k})_{k\in n+1},\zeta)$ consisting of an $(n+1)$-sorted set $(C_{k})_{k\in n+1}$ and a structure of heterogeneous $n$-category
\[
\zeta = (\zeta_{j,k})_{(j,k)\in \coprod_{k\in n+1}k}
\]
on $(C_{k})_{k\in n+1}$, where every $\zeta_{j,k}=(\circ^{(j,k)},\src^{(j,k)},\tgt^{(j,k)},\mathrm{i}^{(k,j)})$, for $k\in n+1$ and $j\in k$, consists of the same data as above, required to satisfy the conditions (MS1), (MS3), (MS4), (MS8) and (MS9) for all $i,j,k\in n+1$ with $i<j<k$, and the conditions (MS2), (MS5), (MS6) and (MS7) for every pair of dimensions in $n+1$ to which they apply, i.e., (MS2), (MS6) and (MS7) for all $j<k$ in $n+1$ and (MS5) for all $i<j$ in $n+1$.
In particular, for $n=1$, a heterogeneous $1$-category is a category presented by its set of objects $C_{0}$ and its set of morphisms $C_{1}$, and (MS5) states that the identity of an object has this object as source and as target.

Heterogeneous $n$-functors are defined as heterogeneous $\omega$-functors, truncated at level $n$: for two heterogeneous $n$-categories $\mathsf{C}$ and $\mathsf{C}'$, a \emph{heterogeneous $n$-functor} from $\mathsf{C}$ to $\mathsf{C}'$ is a family $F=(F_{k})_{k\in n+1}$ of mappings $F_{k}\colon C_{k}\mor C'_{k}$ that satisfies the above equations for all $j,k\in n+1$ with $j<k$.
We denote by $\mathsf{nCat}^{\htg}$ the category of heterogeneous $n$-categories and heterogeneous $n$-functors.
\end{definition}

We next define the $S$-sorted heterogeneous $n$-categories and the $S$-sorted heterogeneous $n$-functors.

\begin{definition}\label{DSortedHeterogeneousnCat}\lean{equivalent}{Heterogeneous.SortedCategory}
Let $S$ be a set of sorts and $n\in\mathbb{N}$.
We call the objects of the functor category $(\mathsf{nCat}^{\htg})^{S}$ \emph{$S$-sorted heterogeneous $n$-categories}.
Thus, an $S$-sorted heterogeneous $n$-category is an $S$-indexed family $\mathsf{A}=(\mathsf{A}_{s})_{s\in S}$ such that, for every sort $s\in S$, $\mathsf{A}_{s}$ is a heterogeneous $n$-category.
Given two $S$-sorted heterogeneous $n$-categories
$\mathsf{A}$ and $\mathsf{B}$, we call the morphisms in $(\mathsf{nCat}^{\htg})^{S}$ from $\mathsf{A}$ to
$\mathsf{B}$ \emph{$S$-sorted heterogeneous $n$-functors}.
Thus, an $S$-sorted heterogeneous $n$-functor from $\mathsf{A}$ to $\mathsf{B}$ is an $S$-indexed family $F=(F_{s})_{s\in S}$ such that, for every sort $s\in S$, $F_{s}$ is a heterogeneous $n$-functor from $\mathsf{A}_{s}$ to $\mathsf{B}_{s}$.
We denote by $\mathrm{nFunc}^{\htg}(\mathsf{A},\mathsf{B})$ the set of all $S$-sorted heterogeneous $n$-functors from $\mathsf{A}$ to $\mathsf{B}$.
We consider the expressions $F\in\mathrm{nFunc}^{\htg}(\mathsf{A},\mathsf{B})$ and $F\colon\mathsf{A}\mor\mathsf{B}$ as synonyms.
Let us recall that, for $F\colon \mathsf{A}\mor\mathsf{B}$ and $G\colon \mathsf{B}\mor\mathsf{C}$, the \emph{composite} $G\circ F$ is $(G_{s}\circ F_{s})_{s\in S}$, and that, for an $S$-sorted heterogeneous $n$-category $\mathsf{A}$, the \emph{identity} $S$-sorted heterogeneous $n$-functor at $\mathsf{A}$ is $\Id_{\mathsf{A}}=(\Id_{\mathsf{A}_{s}})_{s\in S}$.
\end{definition}

Next we define the heterogeneous $n$-categorial algebras as in Definition~\ref{DCatAlg}.

\begin{definition}\label{DCatAlgHtg}\lean{equivalent}{Heterogeneous.NCategorialAlgebra,Heterogeneous.CategorialAlgebra}
Let $S$ be a set of sorts and $n\in\mathbb{N}$.
For an $S$-sorted heterogeneous $n$-category $\mathsf{A}$, the $S^{\star}\times S$-sorted set of finitary operations on $\mathsf{A}$ is $(\mathrm{nFunc}^{\htg}(\mathsf{A}_{\mathbf{s}},\mathsf{A}_{s}))_{(\mathbf{s},s)\in S^{\star}\times S}$, where, for every $\mathbf{s}\in S^{\star}$, $\mathsf{A}_{\mathbf{s}}=\prod_{j\in\bb{\mathbf{s}}}\mathsf{A}_{s_{j}}$, with $\bb{\mathbf{s}}$ denoting the length of the word $\mathbf{s}$.
The finite products $\mathsf{A}_{\mathbf{s}}$ are taken dimensionwise: the $k$-cells of $\mathsf{A}_{\mathbf{s}}$ are the families of $k$-cells, and all the operations are computed componentwise; for $\mathbf{s}=\lambda$, $\mathsf{A}_{\lambda}$ is a final heterogeneous $n$-category, with a single cell in each dimension.

Let $\sig{\Sigma}=(S,\Sigma)$ be a many-sorted signature.
A \emph{structure of heterogeneous $n$-categorial $\Sigma$-algebra on} an $S$-sorted heterogeneous $n$-category $\mathsf{A}$ is a family $F = (F_{\mathbf{s},s})_{(\mathbf{s},s)\in S^{\star}\times S}$ where, for every $(\mathbf{s},s)\in S^{\star}\times S$, $F_{\mathbf{s},s}$ is a mapping from $\Sigma_{\mathbf{s},s}$ to $\mathrm{nFunc}^{\htg}(\mathsf{A}_{\mathbf{s}},\mathsf{A}_{s})$.
If $\sigma\in \Sigma_{\lambda,s}$, then $F_{\lambda,s}(\sigma)$ picks out a $0$-cell $u$ of $\mathsf{A}_{s}$ together with its identities $\mathrm{i}^{(k,0)}(u)$, for $0<k\leq n$, since $\mathsf{A}_{\lambda}$ is a final heterogeneous $n$-category.
An \emph{$S$-sorted heterogeneous $n$-categorial $\Sigma$-algebra} or, simply, a \emph{heterogeneous $n$-categorial $\sig{\Sigma}$-algebra}, is a pair $(\mathsf{A},F)$, abbreviated to $\calg{A}$, where $\mathsf{A}=(\mathsf{A}_{s})_{s\in S}$ is an $S$-sorted heterogeneous $n$-category and $F$ a structure of heterogeneous $n$-categorial $\Sigma$-algebra on $\mathsf{A}$.
For a pair $(\mathbf{s},s)\in S^{\star}\times S$ and an operation symbol $\sigma\in \Sigma_{\mathbf{s},s}$, in order to simplify the notation, the heterogeneous $n$-functor $F_{\mathbf{s},s}(\sigma)$ from $\mathsf{A}_{\mathbf{s}}$ to $\mathsf{A}_{s}$ will be written as $F_{\sigma}$, or simply as $\sigma^{\calg{A}}$.

An \emph{$S$-sorted heterogeneous $n$-categorial $\Sigma$-homomorphism} or, simply, a \emph{heterogeneous $n$-categorial $\sig{\Sigma}$-homomorphism}, from $\calg{A}$ to $\calg{B}$, where $\calg{B}=(\mathsf{B},G)$, is a triple $(\calg{A},H,\calg{B})$, abbreviated to $H\colon\calg{A}\mor\calg{B}$, where $H$ is an $S$-sorted heterogeneous $n$-functor from $\mathsf{A}$ to $\mathsf{B}$ such that, for every $k\in n+1$, every $(\mathbf{s},s)\in S^{\star}\times S$, every $\sigma\in\Sigma_{\mathbf{s},s}$ and every family of $k$-cells $(\mathsf{a}_{j})_{j\in\bb{\mathbf{s}}}\in\mathsf{A}_{\mathbf{s}}$,
$$
H_{s,k}\left(\sigma^{\calg{A}}_{k}
\left(\left(\mathsf{a}_{j}
\right)_{j\in\bb{\mathbf{s}}}\right)\right)
=
\sigma^{\calg{B}}_{k}
\left(\left(
H_{s_{j},k}\left(
\mathsf{a}_{j}
\right)
\right)_{j\in\bb{\mathbf{s}}}
\right),
$$
where $H_{s,k}$ is the component at $k$-cells of the heterogeneous $n$-functor $H_{s}$ and $\sigma^{\calg{A}}_{k}$ that of $\sigma^{\calg{A}}$.
We denote by $\mathsf{nCat}^{\htg}\mathsf{Alg}(\sig{\Sigma})$ the category of heterogeneous $n$-categorial $\sig{\Sigma}$-algebras and heterogeneous $n$-categorial $\sig{\Sigma}$-homomorphisms.
\end{definition}

We next introduce the $\omega$-dimensional counterpart of Definition~\ref{DCatAlgHtg}.

\begin{definition}\label{DCatAlgOmegaHtg}\lean{equivalent}{Heterogeneous.OmegaCategorialAlgebra}
Let $\sig{\Sigma}=(S,\Sigma)$ be a many-sorted signature.
A \emph{heterogeneous $\omega$-categorial $\sig{\Sigma}$-algebra} is a pair $\calg{A}=(\mathsf{A},F)$ where $\mathsf{A}=(\mathsf{A}_{s})_{s\in S}$ is an $S$-sorted heterogeneous $\omega$-category and $F$ assigns to every $(\mathbf{s},s)\in S^{\star}\times S$ and every $\sigma\in\Sigma_{\mathbf{s},s}$ a heterogeneous $\omega$-functor $\sigma^{\calg{A}}\colon\mathsf{A}_{\mathbf{s}}\mor\mathsf{A}_{s}$.
A \emph{heterogeneous $\omega$-categorial $\sig{\Sigma}$-homomorphism} is an $S$-sorted heterogeneous $\omega$-functor satisfying the equation of Definition~\ref{DCatAlgHtg} for every $k\in\omega$.
We denote by $\omega\mathsf{Cat}^{\htg}\mathsf{Alg}(\sig{\Sigma})$ the category of heterogeneous $\omega$-categorial $\sig{\Sigma}$-algebras and heterogeneous $\omega$-categorial $\sig{\Sigma}$-homomorphisms.
\end{definition}

In the heterogeneous presentation, a heterogeneous $n$-categorial $\sig{\Sigma}$-algebra is regarded as a partial algebra on an $(n+1)\times S$-sorted set of cells, for a signature that enlarges $\Sigma$ with a copy of every operation symbol for every dimension, and with sources, targets, identities and compositions between dimensions.

\begin{definition}\label{DCatSigHtg}\lean{equivalent}{Heterogeneous.epSignature,Heterogeneous.hepSignature}
Let $\Sigma$ be an $S$-sorted signature and $n\in\mathbb{N}$.
For every $k\in n+1$ and every $\mathbf{s}\in S^{\star}$, we let $(k,\mathbf{s})$ stand for the word $((k,s_{j}))_{j\in\lvert\mathbf{s}\rvert}$ in $((n+1)\times S)^{\star}$; in particular, $(k,\lambda)=\lambda$.
The \emph{heterogeneous $n$-categorial signature determined by $\Sigma$}, denoted by $\Sigma^{\htg(n)}$, is the $(n+1)\times S$-sorted signature defined, for every $(\mathbf{w},(k,s))\in((n+1)\times S)^{\star}\times((n+1)\times S)$, as follows:
$$
{\Sigma}^{\htg(n)}_{\mathbf{w},(k,s)}=
\begin{cases}
\{\sigma_{k}\mid\sigma\in{\Sigma}_{\mathbf{s},s}\},
& \text{if $\mathbf{w}=(k,\mathbf{s})$ with $\mathbf{s}\neq(s,s)$;}\\
\{\sigma_{k}\mid\sigma\in{\Sigma}_{(s,s),s}\}\amalg\{\circ^{(j,k)}_{s}\mid j\in k\},
& \text{if $\mathbf{w}=((k,s),(k,s))$;}\\
\{\src^{(k,m)}_{s}, \tgt^{(k,m)}_{s}\},
& \text{if $\mathbf{w}=((m,s))$ with $k<m$;}\\
\{\mathrm{i}^{(k,j)}_{s}\},
& \text{if $\mathbf{w}=((j,s))$ with $j<k$;}\\
\varnothing,
& \text{otherwise.}
\end{cases}
$$
That is, the sort $(k,s)$ of ${\Sigma}^{\htg(n)}$ stands for the $k$-cells of sort $s$, and ${\Sigma}^{\htg(n)}$ is obtained from $\Sigma$ by taking, for every $k\in n+1$ and every operation symbol $\sigma\in\Sigma_{\mathbf{s},s}$, a copy $\sigma_{k}$ of $\sigma$ of arity $(k,\mathbf{s})$ and coarity $(k,s)$, which acts on the $k$-cells, and by adding, for every $j,k\in n+1$ with $j<k$ and every sort $s\in S$, two unary operation symbols $\src^{(j,k)}_{s}$ and $\tgt^{(j,k)}_{s}$, the \emph{$(j,k)$-source} and the \emph{$(j,k)$-target}, of arity $((k,s))$ and coarity $(j,s)$, a unary operation symbol $\mathrm{i}^{(k,j)}_{s}$, the \emph{$(k,j)$-identity}, of arity $((j,s))$ and coarity $(k,s)$, and a binary operation symbol $\circ^{(j,k)}_{s}$, the \emph{$j$-composition of $k$-cells}, of arity $((k,s),(k,s))$ and coarity $(k,s)$.
As in Definition~\ref{DCatSig}, we write $x\circ^{(j,k)}_{s}y$ for $\circ^{(j,k)}_{s}(x,y)$.

Moreover, we regard $\Sigma^{\htg(n)}$ as the hep-signature, see Definition~\ref{DHepSig}, $\sig{\Sigma}^{\htg(n)}=((n+1)\times S,\Sigma^{\htg(n)},\mathrm{def}^{((n+1)\times S,\Sigma^{\htg(n)})})$ in which the compositions are the only operation symbols with a nonempty set of domain conditions: for every $j,k\in n+1$ with $j<k$ and every $s\in S$,
$$
\mathrm{def}^{((n+1)\times S,\Sigma^{\htg(n)})}(\circ^{(j,k)}_{s})=\{(\src^{(j,k)}_{s}(v^{(k,s)}_{0}),\tgt^{(j,k)}_{s}(v^{(k,s)}_{1}))\},
$$
where $v^{(k,s)}_{0}$ and $v^{(k,s)}_{1}$ are the variables of $\vs{((k,s),(k,s))}$, so that $v^{(k,s)}_{0}\circ^{(j,k)}_{s}v^{(k,s)}_{1}$ is defined if, and only if, the $(j,k)$-source of $v^{(k,s)}_{0}$ is equal to the $(j,k)$-target of $v^{(k,s)}_{1}$; and $\mathrm{def}^{((n+1)\times S,\Sigma^{\htg(n)})}(\tau)=\varnothing$ for every other operation symbol $\tau$ of $\Sigma^{\htg(n)}$.
The relation $\lhd$ of $\sig{\Sigma}^{\htg(n)}$, see Definition~\ref{DHepSig}, is well-founded, since the domain conditions of the compositions only involve sources and targets, which have empty sets of domain conditions; hence $\sig{\Sigma}^{\htg(n)}$ is indeed a hep-signature.

Finally, the \emph{heterogeneous $\omega$-categorial signature determined by $\Sigma$}, $\Sigma^{\htg(\omega)}$, and the hep-signature $\sig{\Sigma}^{\htg(\omega)}$ are the $\omega\times S$-sorted signature and the hep-signature defined by the same clauses, with $\omega$ in place of $n+1$.
\end{definition}

In contrast with the homogeneous case, for $n\leq m$ in $\mathbb{N}$, $\Sigma^{\htg(n)}$ is not a subset of $\Sigma^{\htg(m)}$, since the sets of sorts are different; however, the inclusion of $(n+1)\times S$ into $(m+1)\times S$ induces a canonical embedding of $\Sigma^{\htg(n)}$ into $\Sigma^{\htg(m)}$.
Likewise, for every $n\in\mathbb{N}$, the inclusion of $(n+1)\times S$ into $\omega\times S$ induces a canonical embedding of $\Sigma^{\htg(n)}$ into $\Sigma^{\htg(\omega)}$, and $\Sigma^{\htg(\omega)}$ is the union of the images of these embeddings.
Let us also note that, for $n=0$, $\Sigma^{\htg(0)}$ is, up to the identification of $1\times S$ with $S$, the signature $\Sigma$.

The specification of the heterogeneous presentation follows the pattern of that of the homogeneous presentation, see Definition~\ref{DCatSpec}.
Now, however, the role played in the homogeneous case by a dimension $k\in n$ is played by a pair of dimensions $j<k$ in $n+1$, since the $j$-composition acts on the $k$-cells and the $(k,j)$-identities take the place of the $k$-sources and $k$-targets as unit elements; and the role played by a pair of dimensions $l<k$ in $n$ is played by a triple of dimensions $h<j<k$ in $n+1$.
Informally, the equations state that the operations of $\Sigma$ are always defined, that, for every sort $s\in S$ and every pair $j<k$, the $j$-cells and the $k$-cells of sort $s$ are the objects and the morphisms of a category, that these categories fit together as in an $n$-category, and that the operations of $\Sigma$ are compatible with them.

\begin{definition}\label{DCatSpecHtg}\lean{equivalent}{Heterogeneous.specification}
Let $\sig{\Sigma}=(S,\Sigma)$ be a many-sorted signature and $n\in\mathbb{N}$.
We denote by $\spec{C}^{\sig{\Sigma}^{\htg(n)}}=(\sig{\Sigma}^{\htg(n)},\mathcal{C}^{\sig{\Sigma}^{\htg(n)}})$ the ep-specification in which $\sig{\Sigma}^{\htg(n)}$ is the heterogeneous $n$-categorial hep-signature introduced in Definition~\ref{DCatSigHtg} and $\mathcal{C}^{\sig{\Sigma}^{\htg(n)}}$ is the set of contextual $\QE$-equations, see Definition~\ref{DEpSpec}, listed below, each of them with the variables that occur in it as context.
The variables are taken from the $(n+1)\times S$-sorted set $V^{(n+1)\times S}$ of Section~\ref{SVariables}, which we abbreviate to $V$.
As in Definition~\ref{DCatSpec}, a family of variables indexed by a word is a tuple of pairwise distinct variables of the indicated sorts, and the families that occur in the same axiom are independent.

For every $k\in n+1$, every $(\mathbf{s},s)\in S^{\star}\times S$, every operation symbol $\sigma\in \Sigma_{\mathbf{s},s}$, and every family $(x_{p})_{p\in \lvert\mathbf{s}\rvert}$ of pairwise distinct variables with $x_{p}\in V_{(k,s_{p})}$, the operation $\sigma_{k}$ applied to the family $(x_{p})_{p\in \lvert\mathbf{s}\rvert}$ is always defined.
Formally,
\begin{align*}
\sigma_{k}((x_{p})_{p\in \lvert\mathbf{s}\rvert})\eeq\sigma_{k}((x_{p})_{p\in \lvert\mathbf{s}\rvert}).\label{ACatHtgEA0}\tag{EA0}
\end{align*}

For every $j,k\in n+1$ with $j<k$, the following equations describe the $j$-composition of $k$-cells.

For every sort $s\in S$, every variable $x\in V_{(k,s)}$ and every variable $u\in V_{(j,s)}$, the $(j,k)$-source and $(j,k)$-target of $x$ and the $(k,j)$-identity of $u$ are always defined.
Formally,
\begin{align*}
\src^{(j,k)}_{s}(x)&\eeq\src^{(j,k)}_{s}(x);
&
\tgt^{(j,k)}_{s}(x)&\eeq\tgt^{(j,k)}_{s}(x);
&
\mathrm{i}^{(k,j)}_{s}(u)&\eeq\mathrm{i}^{(k,j)}_{s}(u).\label{ACatHtgEAjk1}\tag{EA$^{j,k}$1}
\end{align*}

For every sort $s\in S$ and every variable $u\in V_{(j,s)}$, the $(j,k)$-source and the $(j,k)$-target of the $(k,j)$-identity of $u$ are equal to $u$.
Formally,
\begin{align*}
\src^{(j,k)}_{s}(\mathrm{i}^{(k,j)}_{s}(u))&\eeq u;
&
\tgt^{(j,k)}_{s}(\mathrm{i}^{(k,j)}_{s}(u))&\eeq u.\label{ACatHtgEAjk2}\tag{EA$^{j,k}$2}
\end{align*}
In other words, $\mathrm{i}^{(k,j)}_{s}$ is a common section of $\src^{(j,k)}_{s}$ and $\tgt^{(j,k)}_{s}$.
This is the heterogeneous counterpart of the fact that $\src^{k}_{s}$ and $\tgt^{k}_{s}$ are right zeros.

For every sort $s\in S$ and every pair of variables $x,y\in V_{(k,s)}$, $x\circ^{(j,k)}_{s} y$ is defined if, and only if, the $(j,k)$-target of $y$ is equal to the $(j,k)$-source of $x$.
Formally,
\begin{align*}
x\circ^{(j,k)}_{s} y\eeq x\circ^{(j,k)}_{s} y \;\;
&\longrightarrow \;\;
\src^{(j,k)}_{s}(x)\eeq\tgt^{(j,k)}_{s}(y);\\
\src^{(j,k)}_{s}(x)\eeq\tgt^{(j,k)}_{s}(y) \;\;
&\longrightarrow \;\;
x\circ^{(j,k)}_{s} y\eeq x\circ^{(j,k)}_{s}y.\label{ACatHtgEAjk3}\tag{EA$^{j,k}$3}
\end{align*}

For every sort $s\in S$ and every pair of variables $x,y\in V_{(k,s)}$, if $x\circ^{(j,k)}_{s} y$ is defined, then the $(j,k)$-source of $x\circ^{(j,k)}_{s}y$ is that of $y$ and the $(j,k)$-target of $x\circ^{(j,k)}_{s}y$ is that of $x$.
Formally,
\begin{align*}
x\circ^{(j,k)}_{s}y\eeq x\circ^{(j,k)}_{s}y \;\;
&\longrightarrow \;\;
\src^{(j,k)}_{s}(x\circ^{(j,k)}_{s}y)\eeq\src^{(j,k)}_{s}(y);\\
x\circ^{(j,k)}_{s}y\eeq x\circ^{(j,k)}_{s}y \;\;
&\longrightarrow \;\;
\tgt^{(j,k)}_{s}(x\circ^{(j,k)}_{s}y)\eeq\tgt^{(j,k)}_{s}(x).\label{ACatHtgEAjk4}\tag{EA$^{j,k}$4}
\end{align*}

For every sort $s\in S$ and every variable $x\in V_{(k,s)}$, the compositions $x\circ^{(j,k)}_{s}\mathrm{i}^{(k,j)}_{s}(\src^{(j,k)}_{s}(x))$ and $\mathrm{i}^{(k,j)}_{s}(\tgt^{(j,k)}_{s}(x))\circ^{(j,k)}_{s}x$ are always defined and are equal to $x$, i.e., $\mathrm{i}^{(k,j)}_{s}(\src^{(j,k)}_{s}(x))$ is a right unit element for the $j$-composition with $x$ and $\mathrm{i}^{(k,j)}_{s}(\tgt^{(j,k)}_{s}(x))$ is a left unit element for the $j$-composition with $x$.
Formally,
\begin{align*}
x\circ^{(j,k)}_{s}\mathrm{i}^{(k,j)}_{s}(\src^{(j,k)}_{s}(x))&\eeq x;
&
\mathrm{i}^{(k,j)}_{s}(\tgt^{(j,k)}_{s}(x))\circ^{(j,k)}_{s}x&\eeq x.\label{ACatHtgEAjk5}\tag{EA$^{j,k}$5}
\end{align*}

For every sort $s\in S$ and every triple of variables $x, y, z\in V_{(k,s)}$, if the compositions $x\circ^{(j,k)}_{s}y$ and $y\circ^{(j,k)}_{s}z$ are defined, then the compositions $x\circ^{(j,k)}_{s}(y\circ^{(j,k)}_{s}z)$ and $(x\circ^{(j,k)}_{s}y)\circ^{(j,k)}_{s}z$ are defined and they are equal, i.e., the $j$-composition of $k$-cells, when defined, is associative.
Formally,
\begin{multline*}
(x\circ^{(j,k)}_{s}y\eeq x\circ^{(j,k)}_{s}y)
\wedge
(y\circ^{(j,k)}_{s}z\eeq y\circ^{(j,k)}_{s}z)\;\;
\longrightarrow \;\;\\
(x\circ^{(j,k)}_{s}y)\circ^{(j,k)}_{s}z\eeq x\circ^{(j,k)}_{s}(y\circ^{(j,k)}_{s}z).\label{ACatHtgEAjk6}\tag{EA$^{j,k}$6}
\end{multline*}

For every $(\mathbf{s},s)\in S^{\star}\times S$, every operation symbol $\sigma\in \Sigma_{\mathbf{s},s}$, every family $(x_{p})_{p\in \lvert\mathbf{s}\rvert}$ of pairwise distinct variables with $x_{p}\in V_{(k,s_{p})}$ and every family $(u_{p})_{p\in \lvert\mathbf{s}\rvert}$ of pairwise distinct variables with $u_{p}\in V_{(j,s_{p})}$, the $(j,k)$-source of $\sigma_{k}((x_{p})_{p\in \lvert\mathbf{s}\rvert})$ is equal to $\sigma_{j}$ applied to the family $((\src^{(j,k)}_{s_{p}}(x_{p}))_{p\in \lvert\mathbf{s}\rvert})$, the $(j,k)$-target of $\sigma_{k}((x_{p})_{p\in \lvert\mathbf{s}\rvert})$ is equal to $\sigma_{j}$ applied to the family $((\tgt^{(j,k)}_{s_{p}}(x_{p}))_{p\in \lvert\mathbf{s}\rvert})$, and the $(k,j)$-identity of $\sigma_{j}((u_{p})_{p\in \lvert\mathbf{s}\rvert})$ is equal to $\sigma_{k}$ applied to the family $((\mathrm{i}^{(k,j)}_{s_{p}}(u_{p}))_{p\in \lvert\mathbf{s}\rvert})$.
Formally,
\begin{align*}
\src^{(j,k)}_{s}(\sigma_{k}((x_{p})_{p\in \lvert\mathbf{s}\rvert}))&\eeq\sigma_{j}((\src^{(j,k)}_{s_{p}}(x_{p}))_{p\in \lvert\mathbf{s}\rvert});\\
\tgt^{(j,k)}_{s}(\sigma_{k}((x_{p})_{p\in \lvert\mathbf{s}\rvert}))&\eeq\sigma_{j}((\tgt^{(j,k)}_{s_{p}}(x_{p}))_{p\in \lvert\mathbf{s}\rvert});\\
\mathrm{i}^{(k,j)}_{s}(\sigma_{j}((u_{p})_{p\in \lvert\mathbf{s}\rvert}))&\eeq\sigma_{k}((\mathrm{i}^{(k,j)}_{s_{p}}(u_{p}))_{p\in \lvert\mathbf{s}\rvert}).\label{ACatHtgEAjk7}\tag{EA$^{j,k}$7}
\end{align*}

For every $(\mathbf{s},s)\in S^{\star}\times S$, every operation symbol $\sigma\in \Sigma_{\mathbf{s},s}$, and every pair of families $(x_{p})_{p\in \lvert\mathbf{s}\rvert}$ and $(y_{p})_{p\in \lvert\mathbf{s}\rvert}$ of variables with $x_{p},y_{p}\in V_{(k,s_{p})}$, all of them pairwise distinct, if, for every $p\in \lvert\mathbf{s}\rvert$, the compositions $x_{p}\circ^{(j,k)}_{s_{p}}y_{p}$ are defined, then the composition $\sigma_{k}((x_{p})_{p\in \lvert\mathbf{s}\rvert})\circ^{(j,k)}_{s}\sigma_{k}((y_{p})_{p\in \lvert\mathbf{s}\rvert})$ is defined and it is equal to $\sigma_{k}$ applied to the family $(x_{p}\circ^{(j,k)}_{s_{p}}y_{p})_{p\in \lvert\mathbf{s}\rvert}$.
Formally,
\begin{multline*}
\textstyle
\bigwedge_{p\in \lvert\mathbf{s}\rvert}
(x_{p}\circ^{(j,k)}_{s_{p}} y_{p}\eeq x_{p}\circ^{(j,k)}_{s_{p}} y_{p})
\;\;\longrightarrow\;\;\\
\sigma_{k}((x_{p}\circ^{(j,k)}_{s_{p}}y_{p})_{p\in \lvert\mathbf{s}\rvert})
\eeq
\sigma_{k}((x_{p})_{p\in \lvert\mathbf{s}\rvert})\circ^{(j,k)}_{s}\sigma_{k}((y_{p})_{p\in \lvert\mathbf{s}\rvert}).
\label{ACatHtgEAjk8}\tag{EA$^{j,k}$8}
\end{multline*}

For every $h,j,k\in n+1$ with $h<j<k$, the following equations relate the $h$-, $j$- and $k$-cells.

For every sort $s\in S$ and every variable $x\in V_{(k,s)}$, the $(h,j)$-source of both the $(j,k)$-source and the $(j,k)$-target of $x$ is equal to the $(h,k)$-source of $x$, and the $(h,j)$-target of both the $(j,k)$-source and the $(j,k)$-target of $x$ is equal to the $(h,k)$-target of $x$.
Formally,
\begin{align*}
\src^{(h,j)}_{s}(\src^{(j,k)}_{s}(x))&\eeq\src^{(h,k)}_{s}(x);
&
\tgt^{(h,j)}_{s}(\tgt^{(j,k)}_{s}(x))&\eeq\tgt^{(h,k)}_{s}(x);\\
\src^{(h,j)}_{s}(\tgt^{(j,k)}_{s}(x))&\eeq\src^{(h,k)}_{s}(x);
&
\tgt^{(h,j)}_{s}(\src^{(j,k)}_{s}(x))&\eeq\tgt^{(h,k)}_{s}(x).\label{ACatHtgEAhjk1}\tag{EA$^{h,j,k}$1}
\end{align*}

For every sort $s\in S$ and every variable $u\in V_{(h,s)}$, the $(k,j)$-identity of the $(j,h)$-identity of $u$ is the $(k,h)$-identity of $u$.
Formally,
\begin{align*}
\mathrm{i}^{(k,j)}_{s}(\mathrm{i}^{(j,h)}_{s}(u))&\eeq\mathrm{i}^{(k,h)}_{s}(u).\label{ACatHtgEAhjk2}\tag{EA$^{h,j,k}$2}
\end{align*}

For every sort $s\in S$ and every pair of variables $x,y\in V_{(k,s)}$, if $x\circ^{(h,k)}_{s}y$ is defined, then the $(j,k)$-source of $x\circ^{(h,k)}_{s}y$ is the $h$-composition of the $(j,k)$-source of $x$ with the $(j,k)$-source of $y$, and the $(j,k)$-target of $x\circ^{(h,k)}_{s}y$ is the $h$-composition of the $(j,k)$-target of $x$ with the $(j,k)$-target of $y$.
Formally,
\begin{align*}
&x\circ^{(h,k)}_{s}y\eeq x\circ^{(h,k)}_{s}y \;\;\longrightarrow\\
&\qquad\src^{(j,k)}_{s}(x\circ^{(h,k)}_{s}y)\eeq\src^{(j,k)}_{s}(x)\circ^{(h,j)}_{s}\src^{(j,k)}_{s}(y);\\
&x\circ^{(h,k)}_{s}y\eeq x\circ^{(h,k)}_{s}y \;\;\longrightarrow\\
&\qquad\tgt^{(j,k)}_{s}(x\circ^{(h,k)}_{s}y)\eeq\tgt^{(j,k)}_{s}(x)\circ^{(h,j)}_{s}\tgt^{(j,k)}_{s}(y).
\label{ACatHtgEAhjk3}\tag{EA$^{h,j,k}$3}
\end{align*}

For every sort $s\in S$ and every quadruple of variables $x,y,z,t\in V_{(k,s)}$, if the compositions $x\circ^{(j,k)}_{s}y$, $z\circ^{(j,k)}_{s}t$, $x\circ^{(h,k)}_{s}z$ and $y\circ^{(h,k)}_{s}t$ are defined, then the $h$-composition of $x\circ^{(j,k)}_{s}y$ with $z\circ^{(j,k)}_{s}t$ is equal to the $j$-composition of $x\circ^{(h,k)}_{s}z$ with $y\circ^{(h,k)}_{s}t$.
Formally,
\begin{multline*}
(x\circ^{(j,k)}_{s}y\eeq x\circ^{(j,k)}_{s}y)
\wedge
(z\circ^{(j,k)}_{s}t\eeq z\circ^{(j,k)}_{s}t)
\\
\wedge
(x\circ^{(h,k)}_{s}z\eeq x\circ^{(h,k)}_{s}z)
\wedge
(y\circ^{(h,k)}_{s}t\eeq y\circ^{(h,k)}_{s}t)
\\
\;\;\longrightarrow\;\;
(x\circ^{(j,k)}_{s}y)\circ^{(h,k)}_{s}(z\circ^{(j,k)}_{s}t)\eeq
(x\circ^{(h,k)}_{s}z)\circ^{(j,k)}_{s}(y\circ^{(h,k)}_{s}t).
\label{ACatHtgEAhjk4}\tag{EA$^{h,j,k}$4}
\end{multline*}

For every sort $s\in S$ and every pair of variables $u,v\in V_{(j,s)}$, if $u\circ^{(h,j)}_{s}v$ is defined, then the $(k,j)$-identity of $u\circ^{(h,j)}_{s}v$ is the $h$-composition of the $(k,j)$-identities of $u$ and $v$.
Formally,
\begin{align*}
u\circ^{(h,j)}_{s}v\eeq u\circ^{(h,j)}_{s}v \;\;\longrightarrow\;\;
\mathrm{i}^{(k,j)}_{s}(u\circ^{(h,j)}_{s}v)\eeq\mathrm{i}^{(k,j)}_{s}(u)\circ^{(h,k)}_{s}\mathrm{i}^{(k,j)}_{s}(v).
\label{ACatHtgEAhjk5}\tag{EA$^{h,j,k}$5}
\end{align*}

The axioms (EA$^{j,k}$2), (EA$^{j,k}$4)--(EA$^{j,k}$6) and (EA$^{h,j,k}$1)--(EA$^{h,j,k}$5) are the counterparts of the conditions (MS1)--(MS9) of Definition~\ref{DHeterogeneousOmegaCat}, see the correspondence at the end of Remark~\ref{RCatAlgAsEpAlgHtg}; the axioms (EA$^{j,k}$1) and (EA$^{j,k}$3) state that the sources, targets and identities are total and that the domain of the $j$-composition of $k$-cells is the fibered product of Definition~\ref{DHeterogeneousOmegaCat}.
Hence a model of the axioms (EA$^{j,k}$1)--(EA$^{j,k}$6), for $j<k$ in $n+1$, and (EA$^{h,j,k}$1)--(EA$^{h,j,k}$5), for $h<j<k$ in $n+1$, is an $S$-sorted heterogeneous $n$-category, see Definition~\ref{DSortedHeterogeneousnCat}, presented with the cells of each dimension in a separate sort.

Since $\sig{\Sigma}^{\htg(n)}$ is a hep-signature, see Definition~\ref{DCatSigHtg}, $\spec{C}^{\sig{\Sigma}^{\htg(n)}}$ is a hep-specification, see Definition~\ref{DEpSpec}, and we work with the ep-variety $\mathcal{V}(\spec{C}^{\sig{\Sigma}^{\htg(n)}})$ and with its canonically associated category $\mathsf{EPAlg}(\spec{C}^{\sig{\Sigma}^{\htg(n)}})$, see Definition~\ref{DEpVariety}.
\end{definition}

For $n=\omega$, we let $\spec{C}^{\sig{\Sigma}^{\htg(\omega)}}=(\sig{\Sigma}^{\htg(\omega)},\mathcal{C}^{\sig{\Sigma}^{\htg(\omega)}})$ stand for the hep-specification in which $\mathcal{C}^{\sig{\Sigma}^{\htg(\omega)}}$ consists of the contextual $\QE$-equations (EA0), for every $k\in\omega$, (EA$^{j,k}$1)--(EA$^{j,k}$8), for every $j<k$ in $\omega$, and (EA$^{h,j,k}$1)--(EA$^{h,j,k}$5), for every $h<j<k$ in $\omega$.

We next relate the heterogeneous $n$-categorial $\sig{\Sigma}$-algebras with the $\sig{\Sigma}^{\htg(n)}$-algebras.
Let $\sig{\Sigma}=(S,\Sigma)$ be a many-sorted signature and $n\in\mathbb{N}$.
In what follows, a pair $(\mathsf{A},F)$, where $\mathsf{A}=(\mathsf{A}_{s})_{s\in S}$ is an $S$-sorted heterogeneous $n$-category and $F$ assigns to every $(\mathbf{s},s)\in S^{\star}\times S$ and every $\sigma\in\Sigma_{\mathbf{s},s}$ a family of mappings $(F_{\sigma,k})_{k\in n+1}$, one for each dimension, is identified with the partial $\Sigma^{\htg(n)}$-algebra whose underlying set of sort $(k,s)$ is the set of $k$-cells of $\mathsf{A}_{s}$, in which $\src^{(j,k)}_{s}$, $\tgt^{(j,k)}_{s}$, $\mathrm{i}^{(k,j)}_{s}$ and $\circ^{(j,k)}_{s}$ are interpreted as the corresponding operations of $\mathsf{A}_{s}$, and $\sigma_{k}$ is interpreted as $F_{\sigma,k}$.
In dimension $\omega$, the heterogeneous $\omega$-categorial $\sig{\Sigma}$-algebras are identified in the same way, with $\omega$ in place of $n+1$, with partial $\Sigma^{\htg(\omega)}$-algebras.

\begin{remark}\label{RCatAlgAsEpAlgHtg}\lean{equivalent}{Heterogeneous.nModelIsomorphism,Heterogeneous.omegaModelIsomorphism}
Let $\sig{\Sigma}=(S,\Sigma)$ be a many-sorted signature.
Then, under the above identification, the following equalities hold:
\begin{enumerate}
	\item
	for every $n \in \mathbb{N}$, $\mathsf{nCat}^{\htg}\mathsf{Alg}(\sig{\Sigma})=\mathsf{EPAlg}(\spec{C}^{\sig{\Sigma}^{\htg(n)}})$; and

	\item
	$\omega\mathsf{Cat}^{\htg}\mathsf{Alg}(\sig{\Sigma})=\mathsf{EPAlg}(\spec{C}^{\sig{\Sigma}^{\htg(\omega)}})$.
\end{enumerate}
In fact, in a heterogeneous $n$-category, the $(j,k)$-sources, $(j,k)$-targets and $(k,j)$-identities are total mappings, i.e., (EA$^{j,k}$1) holds, and the domain of $\circ^{(j,k)}$ is the fibered product of Definition~\ref{DHeterogeneousOmegaCat}, i.e., (EA$^{j,k}$3) holds; the conditions (MS1)--(MS9) of Definition~\ref{DHeterogeneousOmegaCat} correspond to the remaining axioms on the heterogeneous structure, as listed at the end of this remark; since the products $\mathsf{A}_{\mathbf{s}}$ are taken dimensionwise, a family $(F_{\sigma,k})_{k}$ of partial mappings is a total heterogeneous $n$-functor if, and only if, the axioms (EA0), (EA$^{j,k}$7) and (EA$^{j,k}$8) hold for $\sigma$; and the heterogeneous $n$-categorial $\sig{\Sigma}$-homomorphisms are exactly the homomorphisms of partial $\Sigma^{\htg(n)}$-algebras, see Definitions~\ref{DPartialAlg} and~\ref{DCatAlgHtg}.
Finally, every model of $\spec{C}^{\sig{\Sigma}^{\htg(n)}}$ is a $\sig{\Sigma}^{\htg(n)}$-algebra, see Definitions~\ref{DEpAlg} and~\ref{DCatSigHtg}: the operation symbols $\sigma_{k}$, $\src^{(j,k)}_{s}$, $\tgt^{(j,k)}_{s}$ and $\mathrm{i}^{(k,j)}_{s}$ have empty sets of domain conditions, and they are total, by (EA0) and (EA$^{j,k}$1); and, by (EA$^{j,k}$3), $x\circ^{(j,k)}_{s}y$ is defined if, and only if, $(x,y)$ satisfies the domain condition of $\circ^{(j,k)}_{s}$.
Hence the objects of $\mathsf{EPAlg}(\spec{C}^{\sig{\Sigma}^{\htg(n)}})$ are exactly the models of $\spec{C}^{\sig{\Sigma}^{\htg(n)}}$.
As in Remark~\ref{RCatAlgAsEpAlg}, every axiom only involves variables of the sorts that it constrains, see Definition~\ref{DEpSpec}, so that the identification also holds when some carriers are empty.
The same arguments apply in dimension $\omega$.
In contrast with the homogeneous presentation, see Proposition~\ref{PFiniteDimNotQE}, in the heterogeneous presentation every cell has a dimension, namely the first coordinate of its sort, so that no further condition is needed in dimension $\omega$.
The correspondences between the conditions of Definition~\ref{DHeterogeneousOmegaCat} and the axioms are the following:
\begin{align*}
&\text{(MS1)}\leftrightarrow\text{(EA$^{h,j,k}$1)};
&&\text{(MS2)}\leftrightarrow\text{(EA$^{j,k}$4)};
&&\text{(MS3)}\leftrightarrow\text{(EA$^{h,j,k}$3)};
\\
&\text{(MS4)}\leftrightarrow\text{(EA$^{h,j,k}$2)};
&&\text{(MS5)}\leftrightarrow\text{(EA$^{j,k}$2)};
&&\text{(MS6)}\leftrightarrow\text{(EA$^{j,k}$5)};
\\
&\text{(MS7)}\leftrightarrow\text{(EA$^{j,k}$6)};
&&\text{(MS8)}\leftrightarrow\text{(EA$^{h,j,k}$5)};
&&\text{(MS9)}\leftrightarrow\text{(EA$^{h,j,k}$4)}.
\end{align*}
\end{remark}

\chapter[Graduation and unification]{Graduation and unification of higher-order categorial algebras}\label{SGdUf}

Throughout this chapter, $\sig{\Sigma}=(S,\Sigma)$ is a many-sorted signature and $n\in\mathbb{N}$; they are fixed in all the definitions and results of this chapter.
In this chapter, we relate the homogeneous presentation of Chapter~\ref{SHomogeneous} and the heterogeneous presentation of Chapter~\ref{SHeterogeneous} of the $n$-categorial $\sig{\Sigma}$-algebras.
Graduation associates to a homogeneous $n$-categorial $\sig{\Sigma}$-algebra a heterogeneous $n$-categorial $\sig{\Sigma}$-algebra and, conversely, unification associates to a heterogeneous $n$-categorial $\sig{\Sigma}$-algebra a homogeneous $n$-categorial $\sig{\Sigma}$-algebra.
For instance, for $n=1$, a homogeneous $1$-categorial $\sig{\Sigma}$-algebra is an $S$-sorted family of homogeneous categories together with operations that are functors.
Its graduation has, as cells of sort $(0,s)$, the identities of the category of sort $s$, i.e., its objects, and, as cells of sort $(1,s)$, all its morphisms; this is the usual presentation of a category by a set of objects and a set of morphisms.
Its unification goes back: it keeps only the morphisms, and every object is replaced by its identity.
In Sections~\ref{SGd} and~\ref{SUf} we define the graduation and the unification functors, and in Section~\ref{SGdUfRel} we show that they are mutually quasi-inverse equivalences of categories.
In this chapter and in Chapter~\ref{SDsTr} we only consider finite dimensions; the $\omega$-dimensional functors are considered in Section~\ref{SDsTrOmega} and in Chapter~\ref{SFunctorsOmega}.

Both functors are obtained from the general theory of Chapter~\ref{SGenHepDerivors}.
By Remarks~\ref{RCatAlgAsEpAlg} and~\ref{RCatAlgAsEpAlgHtg}, the following equalities hold:
\begin{enumerate}
	\item
	$\mathsf{nCatAlg}(\sig{\Sigma})=\mathsf{EPAlg}(\spec{C}^{\sig{\Sigma}^{(n)}})$; and

	\item
	$\mathsf{nCat}^{\htg}\mathsf{Alg}(\sig{\Sigma})=\mathsf{EPAlg}(\spec{C}^{\sig{\Sigma}^{\htg(n)}})$.
\end{enumerate}
Thus, the categories of higher-order categorial algebras are the categories canonically associated with the ep-varieties $\mathcal{V}(\spec{C}^{\sig{\Sigma}^{(n)}})$ and $\mathcal{V}(\spec{C}^{\sig{\Sigma}^{\htg(n)}})$ of hep-specifications.
Graduation and unification are induced by generalized hep-derivors between the corresponding hep-signatures, i.e., by morphisms of $\mathsf{HSig}_{\gd}$, see Proposition~\ref{PHepSigCat}.
They are defined on all the partial algebras, and we prove, by means of Proposition~\ref{PTransformedEpAlg}, that their restrictions to the categories of higher-order categorial algebras corestrict to categories of higher-order categorial algebras, see Definition~\ref{DCorestriction}.
By abuse of notation, we denote these corestrictions as the functors themselves.

For the generalized hep-derivors of this work we use the following notation, which we explain for the graduation.
The graduation functor $\Gd^{(n)}_{\sig{\Sigma}}$ is induced by the generalized hep-derivor
$$
\bGd^{(n)}_{\sig{\Sigma}}=\left(\varphi^{\Gd^{(n)}_{\sig{\Sigma}}},\Gd^{(n)}_{\sig{\Sigma}},\Dom^{\Gd^{(n)}_{\sig{\Sigma}}}\right),
$$
whose mapping on sorts is $\varphi^{\Gd^{(n)}_{\sig{\Sigma}}}$, whose assignment of terms to operation symbols is denoted, as the functor, by $\Gd^{(n)}_{\sig{\Sigma}}$, and whose domains are $\Dom^{\Gd^{(n)}_{\sig{\Sigma}}}$.
Since the assignment of terms is applied to operation symbols and the functor to algebras and homomorphisms, this causes no confusion.
The same convention applies to the unification, discretization and truncation derivors $\bUf$, $\bDs$ and $\bTr$.
For the compositions of these generalized hep-derivors, see Definition~\ref{DGenHepDerivorCompId}, we keep the notation $\Dom^{\mathbf{d}\circ\mathbf{c}}$, e.g., $\Dom^{\bGd^{(n)}_{\sig{\Sigma}}\circ\bUf^{(n)}_{\sig{\Sigma}}}$.

\section{Graduation of higher-order categorial algebras}\label{SGd}

Let us recall that the dimensions of $\Sigma^{\htg(n)}$ are the elements of $n+1$.
Graduating a homogeneous $n$-categorial $\sig{\Sigma}$-algebra $\calg{A}$ consists in distributing its cells into sorts according to their dimension: for $k\in n$, the $k$-cells of sort $s$, i.e., the cells $x$ of sort $s$ with $\src^{k}_{s}(x)=x$, become the cells of sort $(k,s)$, and all the cells of sort $s$ become the cells of sort $(n,s)$.
The $(j,k)$-source and the $(j,k)$-target of a $k$-cell are its $j$-source and its $j$-target, the $(k,j)$-identity of a $j$-cell is the cell itself, the $j$-composition of two $k$-cells is their $j$-composition, and $\sigma_{k}$ acts on $k$-cells as $\sigma$.
All these data are given by a single generalized hep-derivor, which we define next.

\begin{definition}\label{DGdDerivor}\lean{equivalent}{Presentations.graduation}
The \emph{graduation derivor} is the generalized hep-derivor
$$
\bGd^{(n)}_{\sig{\Sigma}}=\left(\varphi^{\Gd^{(n)}_{\sig{\Sigma}}},\Gd^{(n)}_{\sig{\Sigma}},\Dom^{\Gd^{(n)}_{\sig{\Sigma}}}\right)
$$
from $\sig{\Sigma}^{\htg(n)}$ to $\sig{\Sigma}^{(n)}$ defined as follows:
\begin{enumerate}
\item for every $k\in n+1$ and every $s\in S$, $\varphi^{\Gd^{(n)}_{\sig{\Sigma}}}(k,s)=s$;
\item for every $(\mathbf{s},s)\in S^{\star}\times S$, every $\sigma\in\Sigma_{\mathbf{s},s}$ and every $k\in n+1$,
$$
\Gd^{(n)}_{\sig{\Sigma}}(\sigma_{k})=\sigma((v^{\mathbf{s}}_{p})_{p\in\bb{\mathbf{s}}});
$$
\item for every pair $j<k$ in $n+1$ and every $s\in S$,
\begin{align*}
\Gd^{(n)}_{\sig{\Sigma}}(\src^{(j,k)}_{s})&=\src^{j}_{s}(v^{s}_{0}), &
\Gd^{(n)}_{\sig{\Sigma}}(\tgt^{(j,k)}_{s})&=\tgt^{j}_{s}(v^{s}_{0}),\\
\Gd^{(n)}_{\sig{\Sigma}}(\mathrm{i}^{(k,j)}_{s})&=v^{s}_{0}, &
\Gd^{(n)}_{\sig{\Sigma}}(\circ^{(j,k)}_{s})&=v^{s}_{0}\circ^{j}_{s}v^{s}_{1};
\end{align*}
\item for every $s\in S$, $\Dom^{\Gd^{(n)}_{\sig{\Sigma}}}_{(n,s)}=\varnothing$; and
\item for every $k\in n$ and every $s\in S$,
$$
\Dom^{\Gd^{(n)}_{\sig{\Sigma}}}_{(k,s)}=\left\lbrace\left(v^{s}_{0},\src^{k}_{s}(v^{s}_{0})\right)\right\rbrace.
$$
\end{enumerate}
\end{definition}

Let us note that $(\varphi^{\Gd^{(n)}_{\sig{\Sigma}}})^{\star}((k,\mathbf{s}))=\mathbf{s}$, so that $\Gd^{(n)}_{\sig{\Sigma}}(\sigma_{k})$ is a term with variables in $\vs{\mathbf{s}}$, and that $j\in n$ whenever $j<k$, so that $\src^{j}_{s}$, $\tgt^{j}_{s}$ and $\circ^{j}_{s}$ are operation symbols of $\Sigma^{(n)}$.
Moreover, by Definition~\ref{DExistenceEq}, for every partial $\Sigma^{(n)}$-algebra $\mathbf{B}$, every $k\in n$, every $s\in S$ and every $x\in B_{s}$, $\mathbf{B}\models\Dom^{\Gd^{(n)}_{\sig{\Sigma}}}_{(k,s)}[x]$ if, and only if, $x=\src^{k\,\mathbf{B}}_{s}(x)$, i.e., $x$ is a $k$-cell; and every element of $B_{s}$ satisfies $\Dom^{\Gd^{(n)}_{\sig{\Sigma}}}_{(n,s)}$.

We next define the graduation functor.

\begin{definition}\label{DGdFunctor}\lean{equivalent}{Presentations.graduationPartialFunctor}
The \emph{graduation functor} $\Gd^{(n)}_{\sig{\Sigma}}$ is the functor $(\bGd^{(n)}_{\sig{\Sigma}})^{\ast}_{\gd}$ induced by the generalized hep-derivor $\bGd^{(n)}_{\sig{\Sigma}}$ from $\mathsf{PAlg}(\Sigma^{(n)})$ to $\mathsf{PAlg}(\Sigma^{\htg(n)})$, see Proposition~\ref{PGenHepDerivorFunctor}.
\end{definition}

In the following lemma we collect the facts about cells of a homogeneous $n$-categorial $\sig{\Sigma}$-algebra that we need to graduate it.

\begin{lemma}\label{LGdCells}\lean{equivalent}{Presentations.operation_isCell,Presentations.source_isCell,Presentations.target_isCell,Presentations.isFiniteCell_mono,Presentations.comp_isCell}
Let $\mathbf{B}$ be an object in $\mathsf{EPAlg}(\spec{C}^{\sig{\Sigma}^{(n)}})$.
Then the following statements hold:
\begin{enumerate}
\item for every $k\in n$, every $(\mathbf{s},s)\in S^{\star}\times S$, every $\sigma\in\Sigma_{\mathbf{s},s}$ and every family $(x_{p})_{p\in\bb{\mathbf{s}}}$ of $k$-cells, with $x_{p}$ of sort $s_{p}$, $\sigma^{\mathbf{B}}((x_{p})_{p\in\bb{\mathbf{s}}})$ is a $k$-cell;
\item for every $j\in n$, every $s\in S$ and every cell $x$ of sort $s$, $\src^{j\,\mathbf{B}}_{s}(x)$ and $\tgt^{j\,\mathbf{B}}_{s}(x)$ are $j$-cells;
\item for every $j<k$ in $n$ and every $s\in S$, every $j$-cell of sort $s$ is a $k$-cell;
\item for every $j<k$ in $n$ and every $s\in S$, if $x$, $y$ are $k$-cells of sort $s$ such that $\src^{j\,\mathbf{B}}_{s}(x)=\tgt^{j\,\mathbf{B}}_{s}(y)$, then $x\circ^{j\,\mathbf{B}}_{s}y$ is defined and it is a $k$-cell.
\end{enumerate}
\end{lemma}
\begin{proof}
(1) The following chain of equalities holds
\begin{align*}
\src^{k\,\mathbf{B}}_{s}\left(\sigma^{\mathbf{B}}\left((x_{p})_{p\in\bb{\mathbf{s}}}\right)\right)
&=
\sigma^{\mathbf{B}}\left(\left(\src^{k\,\mathbf{B}}_{s_{p}}(x_{p})\right)_{p\in\bb{\mathbf{s}}}\right)
\tag{1}
\\&=
\sigma^{\mathbf{B}}\left((x_{p})_{p\in\bb{\mathbf{s}}}\right).
\tag{2}
\end{align*}

The first equality follows from~(A$^{k}$7); the second equality follows from the fact that every $x_{p}$ is a $k$-cell.

(2) By~(A$^{j}$2), the following equalities hold
\begin{align*}
\src^{j\,\mathbf{B}}_{s}\left(\src^{j\,\mathbf{B}}_{s}(x)\right)
&=
\src^{j\,\mathbf{B}}_{s}(x),
&
\src^{j\,\mathbf{B}}_{s}\left(\tgt^{j\,\mathbf{B}}_{s}(x)\right)
&=
\tgt^{j\,\mathbf{B}}_{s}(x).
\end{align*}

(3) Let $u$ be a $j$-cell of sort $s$.
Then the following chain of equalities holds
\begin{align*}
\src^{k\,\mathbf{B}}_{s}(u)
&=
\src^{k\,\mathbf{B}}_{s}\left(\src^{j\,\mathbf{B}}_{s}(u)\right)
\tag{1}
\\&=
\src^{j\,\mathbf{B}}_{s}(u)
\tag{2}
\\&=
u.
\tag{3}
\end{align*}

The first equality follows from the fact that $u$ is a $j$-cell; the second equality follows from~(AB$^{j,k}$1); finally, the last equality follows, again, from the fact that $u$ is a $j$-cell.

(4) By~(A$^{j}$3), $x\circ^{j\,\mathbf{B}}_{s}y$ is defined.
Moreover, the following chain of equalities holds
\begin{align*}
\src^{k\,\mathbf{B}}_{s}\left(x\circ^{j\,\mathbf{B}}_{s}y\right)
&=
\src^{k\,\mathbf{B}}_{s}(x)\circ^{j\,\mathbf{B}}_{s}\src^{k\,\mathbf{B}}_{s}(y)
\tag{1}
\\&=
x\circ^{j\,\mathbf{B}}_{s}y.
\tag{2}
\end{align*}

The first equality follows from~(AB$^{j,k}$2); the second equality follows from the fact that $x$ and $y$ are $k$-cells. This completes the proof.
\end{proof}

In the following lemma we describe the transformed algebra along the graduation derivor.

\begin{lemma}\label{LGdDescription}\lean{equivalent}{Presentations.graduation_transformed,Presentations.graduatedAlgebra}
Let $\mathbf{B}$ be an object in $\mathsf{EPAlg}(\spec{C}^{\sig{\Sigma}^{(n)}})$.
Then $\Gd^{(n)}_{\sig{\Sigma}}(\mathbf{B})$ is a $\sig{\Sigma}^{\htg(n)}$-algebra whose cells of sort $(k,s)$ are the $k$-cells of sort $s$ of $\mathbf{B}$, for $k\in n$, and all the cells of sort $s$ of $\mathbf{B}$, for $k=n$.
For every $(\mathbf{s},s)\in S^{\star}\times S$, every $\sigma\in\Sigma_{\mathbf{s},s}$ and every $k\in n+1$, $\sigma_k$ is total in $\Gd^{(n)}_{\sig{\Sigma}}(\mathbf{B})$ and acts as the restriction of $\sigma^{\mathbf{B}}$ to the $k$-cells.
Moreover, for every pair $j<k$ in $n+1$ and every $s\in S$, the following statements hold:
\begin{enumerate}
\item $\src^{(j,k)}_{s}$, $\tgt^{(j,k)}_{s}$ and $\mathrm{i}^{(k,j)}_{s}$ are total in $\Gd^{(n)}_{\sig{\Sigma}}(\mathbf{B})$ and act as $\src^{j\,\mathbf{B}}_{s}$, $\tgt^{j\,\mathbf{B}}_{s}$ and the inclusion, respectively; and
\item for every pair of $k$-cells $x$, $y$ of sort $s$, $x\circ^{(j,k)}_{s}y$ is defined in $\Gd^{(n)}_{\sig{\Sigma}}(\mathbf{B})$ if, and only if, $\src^{j\,\mathbf{B}}_{s}(x)=\tgt^{j\,\mathbf{B}}_{s}(y)$, and, in this case, the $j$-composition of $x$ and $y$ in $\Gd^{(n)}_{\sig{\Sigma}}(\mathbf{B})$ is $x\circ^{j\,\mathbf{B}}_{s}y$.
\end{enumerate}
\end{lemma}
\begin{proof}
The description of the cells follows from Definitions~\ref{DDomainSet} and~\ref{DGdDerivor} and the note after Definition~\ref{DGdDerivor}.
By Proposition~\ref{PTransformedEpAlg}, to prove that $\Gd^{(n)}_{\sig{\Sigma}}(\mathbf{B})$ is a $\sig{\Sigma}^{\htg(n)}$-algebra, we have to show that, for every operation symbol $\tau$ of $\Sigma^{\htg(n)}$ and every family $\mathbf{x}$ of cells of $\Gd^{(n)}_{\sig{\Sigma}}(\mathbf{B})$ that satisfies the domain conditions of $\tau$ in $\Gd^{(n)}_{\sig{\Sigma}}(\mathbf{B})$, the partial term operation $\Gd^{(n)}_{\sig{\Sigma}}(\tau)^{\mathbf{B}}$ is defined at $\mathbf{x}$ and its value is a cell of the coarity of $\tau$.
We consider the different kinds of operation symbols of $\Sigma^{\htg(n)}$.
We first consider $\sigma_k$, for every $(\mathbf{s},s)\in S^{\star}\times S$, every $\sigma\in\Sigma_{\mathbf{s},s}$ and every $k\in n+1$.
The set of domain conditions of $\sigma_{k}$ is empty and the partial term operation $\sigma((v^{\mathbf{s}}_{p})_{p\in\bb{\mathbf{s}}})^{\mathbf{B}}$, which is $\sigma^{\mathbf{B}}$, is total, by~(A0).
If $k=n$, there is nothing more to prove.
If $k\in n$, then, by Lemma~\ref{LGdCells}~(1), the value at a family of $k$-cells is a $k$-cell.

For the structural operations, fix $j<k$ in $n+1$ and $s\in S$.
We begin with the sources and the targets.
The sets of domain conditions of $\src^{(j,k)}_{s}$ and $\tgt^{(j,k)}_{s}$ are empty and the partial term operations $\src^{j}_{s}(v^{s}_{0})^{\mathbf{B}}$ and $\tgt^{j}_{s}(v^{s}_{0})^{\mathbf{B}}$, which are $\src^{j\,\mathbf{B}}_{s}$ and $\tgt^{j\,\mathbf{B}}_{s}$, are total, by~(A$^{j}$1).
By Lemma~\ref{LGdCells}~(2), their values are $j$-cells.

Next, we consider the identities.
The set of domain conditions of $\mathrm{i}^{(k,j)}_{s}$ is empty and the partial term operation $(v^{s}_{0})^{\mathbf{B}}$, which is the identity of $B_{s}$, is total.
If $k=n$, there is nothing more to prove.
If $k\in n$, then, by Lemma~\ref{LGdCells}~(3), every $j$-cell is a $k$-cell.

Therefore, in $\Gd^{(n)}_{\sig{\Sigma}}(\mathbf{B})$, $\src^{(j,k)}_{s}$, $\tgt^{(j,k)}_{s}$ and $\mathrm{i}^{(k,j)}_{s}$ are total and act as $\src^{j\,\mathbf{B}}_{s}$, $\tgt^{j\,\mathbf{B}}_{s}$ and the inclusion, respectively.

Finally, we consider the compositions.
By the above, two $k$-cells $x$, $y$ of sort $s$ satisfy the domain condition of $\circ^{(j,k)}_{s}$ in $\Gd^{(n)}_{\sig{\Sigma}}(\mathbf{B})$ if, and only if, $\src^{j\,\mathbf{B}}_{s}(x)=\tgt^{j\,\mathbf{B}}_{s}(y)$.
In this case, by~(A$^{j}$3), $x\circ^{j\,\mathbf{B}}_{s}y$ is defined.
If $k=n$, there is nothing more to prove.
If $k\in n$, then, by Lemma~\ref{LGdCells}~(4), $x\circ^{j\,\mathbf{B}}_{s}y$ is a $k$-cell.

All in all, by Proposition~\ref{PTransformedEpAlg}, $\Gd^{(n)}_{\sig{\Sigma}}(\mathbf{B})$ is a $\sig{\Sigma}^{\htg(n)}$-algebra and, by Proposition~\ref{PTransformedAlg}, its operations are as stated. This completes the proof.
\end{proof}

We next prove that the graduation of a homogeneous $n$-categorial $\sig{\Sigma}$-algebra is a heterogeneous $n$-categorial $\sig{\Sigma}$-algebra.

\begin{proposition}\label{PGdEpVariety}\lean{equivalent}{Presentations.graduation_isModel,Presentations.graduationModels}
The restriction of $\Gd^{(n)}_{\sig{\Sigma}}$ to $\mathsf{EPAlg}(\spec{C}^{\sig{\Sigma}^{(n)}})$ corestricts to $\mathsf{EPAlg}(\spec{C}^{\sig{\Sigma}^{\htg(n)}})$.
\end{proposition}
\begin{proof}
Let $\mathbf{B}$ be an object in $\mathsf{EPAlg}(\spec{C}^{\sig{\Sigma}^{(n)}})$.
By Lemma~\ref{LGdDescription}, $\Gd^{(n)}_{\sig{\Sigma}}(\mathbf{B})$ is a $\sig{\Sigma}^{\htg(n)}$-algebra.
It remains to check the axioms of $\spec{C}^{\sig{\Sigma}^{\htg(n)}}$.

By Lemma~\ref{LGdDescription}, the operations of $\Gd^{(n)}_{\sig{\Sigma}}(\mathbf{B})$ are restrictions of those of $\mathbf{B}$, its identities are inclusions, and $x\circ^{(j,k)\,\Gd^{(n)}_{\sig{\Sigma}}(\mathbf{B})}_{s}y$ is defined if, and only if, $x\circ^{j\,\mathbf{B}}_{s}y$ is defined.
Hence the axioms become statements on $\mathbf{B}$.
By Lemma~\ref{LGdDescription}, the axioms (EA0) hold for every operation symbol of $\Sigma$ and every $k\in n+1$.
Now fix $j<k$ in $n+1$.
The axioms (EA$^{j,k}$1) and (EA$^{j,k}$3) hold by the same lemma.
The axiom (EA$^{j,k}$2) holds since every $j$-cell $u$ satisfies $\src^{j\,\mathbf{B}}_{s}(u)=u=\tgt^{j\,\mathbf{B}}_{s}(u)$.
The axioms (EA$^{j,k}$4)--(EA$^{j,k}$8) become the axioms (A$^{j}$4)--(A$^{j}$8) restricted to $k$-cells; the last equation of (EA$^{j,k}$7) becomes trivial.
For the remaining axioms, fix $h<j<k$ in $n+1$.
The axioms (EA$^{h,j,k}$2) and (EA$^{h,j,k}$5) become trivial, since the identities are inclusions.

Finally, the axioms (EA$^{h,j,k}$1), (EA$^{h,j,k}$3) and (EA$^{h,j,k}$4) become the axioms (AB$^{h,j}$1), (AB$^{h,j}$2) and (AB$^{h,j}$3) restricted to $k$-cells, the last one after exchanging $y$ and $z$.
Here we use that, by Lemma~\ref{LGdCells}~(2), the $j$-sources and the $j$-targets are $j$-cells.
For instance, for $k$-cells $x$, $y$ of sort $s$ such that $x\circ^{(h,k)\,\Gd^{(n)}_{\sig{\Sigma}}(\mathbf{B})}_{s}y$ is defined, the following chain of equalities holds
\begin{flushleft}
	$\src^{(j,k)\,\Gd^{(n)}_{\sig{\Sigma}}(\mathbf{B})}_{s}\left(x\circ^{(h,k)\,\Gd^{(n)}_{\sig{\Sigma}}(\mathbf{B})}_{s}y\right)$
\begin{align*}
\quad
&=
\src^{j\,\mathbf{B}}_{s}\left(x\circ^{h\,\mathbf{B}}_{s}y\right)
\tag{1}
\\&=
\src^{j\,\mathbf{B}}_{s}(x)\circ^{h\,\mathbf{B}}_{s}\src^{j\,\mathbf{B}}_{s}(y)
\tag{2}
\\&=
\src^{(j,k)\,\Gd^{(n)}_{\sig{\Sigma}}(\mathbf{B})}_{s}(x)
\circ^{(h,j)\,\Gd^{(n)}_{\sig{\Sigma}}(\mathbf{B})}_{s}\src^{(j,k)\,\Gd^{(n)}_{\sig{\Sigma}}(\mathbf{B})}_{s}(y).
\tag{3}
\end{align*}
\end{flushleft}

The first and the last equalities follow from Lemma~\ref{LGdDescription}; the second equality follows from~(AB$^{h,j}$2). This completes the proof.
\end{proof}

\begin{remark}\label{RGdNotEpAlg}\lean{equivalent}{Presentations.GraduationCounterexample.graduation_not_isEPAlgebra}
Assume that $S\neq\varnothing$ and $n\geq 1$.
The restriction of $\Gd^{(n)}_{\sig{\Sigma}}$ to $\mathsf{EPAlg}(\sig{\Sigma}^{(n)})$ does not corestrict to $\mathsf{EPAlg}(\sig{\Sigma}^{\htg(n)})$.
Indeed, let $\mathbf{B}$ be a $\sig{\Sigma}^{(n)}$-algebra, $j<k$ in $n+1$ and $s\in S$.
The operation $\src^{(j,k)}_{s}$ has no domain conditions, so that it must be total in $\Gd^{(n)}_{\sig{\Sigma}}(\mathbf{B})$, and its values must be cells of sort $(j,s)$, i.e., $j$-cells of $\mathbf{B}$.
If $\mathbf{B}$ does not satisfy~(A$^{j}$2), i.e., if $\src^{j\,\mathbf{B}}_{s}$ is not idempotent, the $j$-source of a cell need not be a $j$-cell; the axiom~(A$^{j}$2) is precisely what is used in the proof of Lemma~\ref{LGdCells}~(2).
For instance, let $\mathbf{B}$ be the partial $\Sigma^{(n)}$-algebra in which, for every $t\in S$, $B_{t}=\{0,1\}$, $\src^{0\,\mathbf{B}}_{t}$ and $\tgt^{0\,\mathbf{B}}_{t}$ are both the transposition $q$ of $\{0,1\}$, $\src^{k\,\mathbf{B}}_{t}$ and $\tgt^{k\,\mathbf{B}}_{t}$ are the identity of $\{0,1\}$, for $1\leq k<n$, and, for every $k\in n$, $x\circ^{k\,\mathbf{B}}_{t}y$ is defined if, and only if, $x=y$, and then it is equal to $x$; and in which, for every operation symbol $\sigma$ of $\Sigma$, $\sigma^{\mathbf{B}}$ is the constant mapping with value $0$.
Since, for every $k\in n$, $\src^{k\,\mathbf{B}}_{t}(x)=\tgt^{k\,\mathbf{B}}_{t}(y)$ if, and only if, $x=y$, $\mathbf{B}$ is a $\sig{\Sigma}^{(n)}$-algebra.
It does not satisfy~(A$^{0}$2), since $q\circ q$ is the identity of $\{0,1\}$ and not $q$.
Moreover, no element of $B_{s}$, for $s\in S$, is a $0$-cell, since $q$ has no fixed points, so that $\Gd^{(n)}_{\sig{\Sigma}}(\mathbf{B})$ has no cells of sort $(0,s)$ and $\src^{(0,n)}_{s}$ cannot be total in $\Gd^{(n)}_{\sig{\Sigma}}(\mathbf{B})$ on the nonempty set $B_{s}$ of its cells of sort $(n,s)$.
For $n=0$ the graduation only relabels the sorts, and it trivially corestricts.
This is why we restrict $\Gd^{(n)}_{\sig{\Sigma}}$ to the category canonically associated with the ep-variety $\mathcal{V}(\spec{C}^{\sig{\Sigma}^{(n)}})$, whose objects satisfy~(A$^{j}$2), for every $j\in n$.
\end{remark}

From Proposition~\ref{PGdEpVariety} we obtain the graduation functor between the categories of higher-order categorial algebras.

\begin{corollary}\label{CGdFunctor}\lean{equivalent}{Presentations.graduationFunctor}
The restriction of $\Gd^{(n)}_{\sig{\Sigma}}$ to $\mathsf{nCatAlg}(\sig{\Sigma})$ corestricts to $\mathsf{nCat}^{\htg}\mathsf{Alg}(\sig{\Sigma})$.
\end{corollary}
\begin{proof}
It follows from Proposition~\ref{PGdEpVariety} and Remarks~\ref{RCatAlgAsEpAlg} and~\ref{RCatAlgAsEpAlgHtg}.
\end{proof}

By abuse of notation, we denote the corestriction of Corollary~\ref{CGdFunctor} by
$$
\Gd^{(n)}_{\sig{\Sigma}}\colon\mathsf{nCatAlg}(\sig{\Sigma})\mor\mathsf{nCat}^{\htg}\mathsf{Alg}(\sig{\Sigma}).
$$
Explicitly, by Lemma~\ref{LGdDescription}, $\Gd^{(n)}_{\sig{\Sigma}}(\calg{A})$ is the graduated algebra described at the beginning of this section, and, by Proposition~\ref{PGenHepDerivorFunctor}, $\Gd^{(n)}_{\sig{\Sigma}}$ sends an $n$-categorial $\sig{\Sigma}$-homomorphism $F$ to the family of its restrictions to the $k$-cells.

\section{Unification of higher-order categorial algebras}\label{SUf}

Unifying a heterogeneous $n$-categorial $\sig{\Sigma}$-algebra $\calg{A}$ consists in keeping only its $n$-cells, every $k$-cell being represented by its $(n,k)$-identity: the $k$-source of an $n$-cell is the $(n,k)$-identity of its $(k,n)$-source, the $k$-composition is the $k$-composition of $n$-cells, and $\sigma$ acts as $\sigma_{n}$.
These data are given by the following generalized hep-derivor.

\begin{definition}\label{DUfDerivor}\lean{equivalent}{Presentations.unification}
The \emph{unification derivor} is the generalized hep-derivor
$$
\bUf^{(n)}_{\sig{\Sigma}}=\left(\varphi^{\Uf^{(n)}_{\sig{\Sigma}}},\Uf^{(n)}_{\sig{\Sigma}},\Dom^{\Uf^{(n)}_{\sig{\Sigma}}}\right)
$$
from $\sig{\Sigma}^{(n)}$ to $\sig{\Sigma}^{\htg(n)}$ defined as follows.
For every $s\in S$, put $\varphi^{\Uf^{(n)}_{\sig{\Sigma}}}(s)=(n,s)$ and $\Dom^{\Uf^{(n)}_{\sig{\Sigma}}}_{s}=\varnothing$.
For every $(\mathbf{s},s)\in S^{\star}\times S$ and every $\sigma\in\Sigma_{\mathbf{s},s}$, put
$$
\Uf^{(n)}_{\sig{\Sigma}}(\sigma)=\sigma_{n}((v^{(n,\mathbf{s})}_{p})_{p\in\bb{\mathbf{s}}}).
$$
For every $k\in n$ and every $s\in S$, put
\begin{align*}
\Uf^{(n)}_{\sig{\Sigma}}(\src^{k}_{s})&=\mathrm{i}^{(n,k)}_{s}(\src^{(k,n)}_{s}(v^{(n,s)}_{0})),\\
\Uf^{(n)}_{\sig{\Sigma}}(\tgt^{k}_{s})&=\mathrm{i}^{(n,k)}_{s}(\tgt^{(k,n)}_{s}(v^{(n,s)}_{0})),\\
\Uf^{(n)}_{\sig{\Sigma}}(\circ^{k}_{s})&=v^{(n,s)}_{0}\circ^{(k,n)}_{s}v^{(n,s)}_{1}.
\end{align*}
\end{definition}

When $n=0$, graduation and unification only relabel the sorts $s\leftrightarrow(0,s)$: their domains are empty, and they send $\sigma_0$ to $\sigma$ applied to its variables and $\sigma$ to $\sigma_0$ applied to its variables, respectively, also when $\sigma$ is nullary.
When $\Sigma$ is empty, only the clauses for the sources, targets, identities and compositions remain in both definitions.

We next define the unification functor.

\begin{definition}\label{DUfFunctor}\lean{equivalent}{Presentations.unificationPartialFunctor}
The \emph{unification functor} $\Uf^{(n)}_{\sig{\Sigma}}$ is the functor $(\bUf^{(n)}_{\sig{\Sigma}})^{\ast}_{\gd}$ induced by the generalized hep-derivor $\bUf^{(n)}_{\sig{\Sigma}}$ from $\mathsf{PAlg}(\Sigma^{\htg(n)})$ to $\mathsf{PAlg}(\Sigma^{(n)})$, see Proposition~\ref{PGenHepDerivorFunctor}.
\end{definition}

In the following lemma we prove that the identities of a heterogeneous $n$-categorial $\sig{\Sigma}$-algebra are injective.

\begin{lemma}\label{LUfIdInj}\lean{equivalent}{Presentations.identity_injective}
Let $\mathbf{A}$ be an object in $\mathsf{EPAlg}(\spec{C}^{\sig{\Sigma}^{\htg(n)}})$, $j<k$ in $n+1$ and $s\in S$.
Then, for every $j$-cell $u$ of sort $s$, the following equalities hold:
$$
\src^{(j,k)\,\mathbf{A}}_{s}\left(\mathrm{i}^{(k,j)\,\mathbf{A}}_{s}(u)\right)=u=\tgt^{(j,k)\,\mathbf{A}}_{s}\left(\mathrm{i}^{(k,j)\,\mathbf{A}}_{s}(u)\right).
$$
In particular, $\mathrm{i}^{(k,j)\,\mathbf{A}}_{s}$ is injective.
\end{lemma}
\begin{proof}
The equalities are the axiom (EA$^{j,k}$2).
Hence $\mathrm{i}^{(k,j)\,\mathbf{A}}_{s}$ has a left inverse, so that it is injective.
\end{proof}

In the following lemma we describe the transformed algebra along the unification derivor.

\begin{lemma}\label{LUfDescription}\lean{weaker}{Presentations.unificationRealization,Presentations.unifiedAlgebra}
Let $\mathbf{A}$ be an object in $\mathsf{EPAlg}(\spec{C}^{\sig{\Sigma}^{\htg(n)}})$.
Then $\Uf^{(n)}_{\sig{\Sigma}}(\mathbf{A})$ is a $\sig{\Sigma}^{(n)}$-algebra whose cells of sort $s$ are the $n$-cells of sort $s$ of $\mathbf{A}$.
For every $(\mathbf{s},s)\in S^{\star}\times S$ and every $\sigma\in\Sigma_{\mathbf{s},s}$, the operation $\sigma$ is total in $\Uf^{(n)}_{\sig{\Sigma}}(\mathbf{A})$ and acts as $\sigma_n^{\mathbf{A}}$.
Moreover, for every $k\in n$ and every $s\in S$, the following statements hold:
\begin{enumerate}
\item $\src^{k}_{s}$ and $\tgt^{k}_{s}$ are total in $\Uf^{(n)}_{\sig{\Sigma}}(\mathbf{A})$ and act as $\mathrm{i}^{(n,k)\,\mathbf{A}}_{s}\circ\src^{(k,n)\,\mathbf{A}}_{s}$ and $\mathrm{i}^{(n,k)\,\mathbf{A}}_{s}\circ\tgt^{(k,n)\,\mathbf{A}}_{s}$, respectively; and
\item for every pair of $n$-cells $x$, $y$ of sort $s$, $x\circ^{k}_{s}y$ is defined in $\Uf^{(n)}_{\sig{\Sigma}}(\mathbf{A})$ if, and only if, $\src^{(k,n)\,\mathbf{A}}_{s}(x)=\tgt^{(k,n)\,\mathbf{A}}_{s}(y)$, and, in this case, the $k$-composition of $x$ and $y$ in $\Uf^{(n)}_{\sig{\Sigma}}(\mathbf{A})$ is $x\circ^{(k,n)\,\mathbf{A}}_{s}y$.
\end{enumerate}
\end{lemma}
\begin{proof}
Since the domains of $\bUf^{(n)}_{\sig{\Sigma}}$ are empty, the cells of sort $s$ of $\Uf^{(n)}_{\sig{\Sigma}}(\mathbf{A})$ are the $n$-cells of sort $s$ of $\mathbf{A}$.
By Proposition~\ref{PTransformedEpAlg}, it suffices to check that, for every operation symbol $\tau$ of $\Sigma^{(n)}$, the partial term operation $\Uf^{(n)}_{\sig{\Sigma}}(\tau)^{\mathbf{A}}$ is defined at the families of cells that satisfy the domain conditions of $\tau$ in $\Uf^{(n)}_{\sig{\Sigma}}(\mathbf{A})$; since every $n$-cell is a cell of $\Uf^{(n)}_{\sig{\Sigma}}(\mathbf{A})$, its values are then cells.

For every $(\mathbf{s},s)\in S^{\star}\times S$ and every $\sigma\in\Sigma_{\mathbf{s},s}$, the domain conditions of $\sigma$ are empty and the partial term operation $\sigma_n((v^{(n,\mathbf{s})}_p)_{p\in\bb{\mathbf{s}}})^{\mathbf{A}}$, which is $\sigma_n^{\mathbf{A}}$, is total by~(EA0), including when $n=0$.
For the structural operations, fix $k\in n$ and $s\in S$.
The sets of domain conditions of $\src^{k}_{s}$ and $\tgt^{k}_{s}$ are empty, and the partial term operations $\mathrm{i}^{(n,k)}_s(\src^{(k,n)}_s(v^{(n,s)}_0))^{\mathbf{A}}$ and $\mathrm{i}^{(n,k)}_s(\tgt^{(k,n)}_s(v^{(n,s)}_0))^{\mathbf{A}}$, which are $\mathrm{i}^{(n,k)\,\mathbf{A}}_{s}\circ\src^{(k,n)\,\mathbf{A}}_{s}$ and $\mathrm{i}^{(n,k)\,\mathbf{A}}_{s}\circ\tgt^{(k,n)\,\mathbf{A}}_{s}$, are total by~(EA$^{k,n}$1).

Now we consider the compositions.
By the above, two $n$-cells $x$, $y$ of sort $s$ satisfy the domain condition of $\circ^{k}_{s}$ in $\Uf^{(n)}_{\sig{\Sigma}}(\mathbf{A})$ if, and only if,
$$
\mathrm{i}^{(n,k)\,\mathbf{A}}_{s}\left(\src^{(k,n)\,\mathbf{A}}_{s}(x)\right)=\mathrm{i}^{(n,k)\,\mathbf{A}}_{s}\left(\tgt^{(k,n)\,\mathbf{A}}_{s}(y)\right),
$$
i.e., by Lemma~\ref{LUfIdInj}, if, and only if, $\src^{(k,n)\,\mathbf{A}}_{s}(x)=\tgt^{(k,n)\,\mathbf{A}}_{s}(y)$.
In this case, by~(EA$^{k,n}$3), $x\circ^{(k,n)\,\mathbf{A}}_{s}y$ is defined.

All in all, by Proposition~\ref{PTransformedEpAlg}, $\Uf^{(n)}_{\sig{\Sigma}}(\mathbf{A})$ is a $\sig{\Sigma}^{(n)}$-algebra and, by Proposition~\ref{PTransformedAlg}, its operations are as stated. This completes the proof.
\end{proof}

We next prove that the unification of a heterogeneous $n$-categorial $\sig{\Sigma}$-algebra is a homogeneous $n$-categorial $\sig{\Sigma}$-algebra.

\begin{proposition}\label{PUfEpVariety}\lean{equivalent}{Presentations.unification_isModel,Presentations.unificationModels}
The restriction of $\Uf^{(n)}_{\sig{\Sigma}}$ to $\mathsf{EPAlg}(\spec{C}^{\sig{\Sigma}^{\htg(n)}})$ corestricts to $\mathsf{EPAlg}(\spec{C}^{\sig{\Sigma}^{(n)}})$.
\end{proposition}
\begin{proof}
Let $\mathbf{A}$ be an object in $\mathsf{EPAlg}(\spec{C}^{\sig{\Sigma}^{\htg(n)}})$.
By Lemma~\ref{LUfDescription}, $\Uf^{(n)}_{\sig{\Sigma}}(\mathbf{A})$ is a $\sig{\Sigma}^{(n)}$-algebra.
It remains to check the axioms of $\spec{C}^{\sig{\Sigma}^{(n)}}$.

By Lemma~\ref{LUfDescription}, these axioms are equivalent to their translations, i.e., to the statements on the $n$-cells of $\mathbf{A}$ obtained from them by replacing $\src^{k\,\Uf^{(n)}_{\sig{\Sigma}}(\mathbf{A})}_{s}$ by $\mathrm{i}^{(n,k)\,\mathbf{A}}_{s}\circ\src^{(k,n)\,\mathbf{A}}_{s}$, $\tgt^{k\,\Uf^{(n)}_{\sig{\Sigma}}(\mathbf{A})}_{s}$ by $\mathrm{i}^{(n,k)\,\mathbf{A}}_{s}\circ\tgt^{(k,n)\,\mathbf{A}}_{s}$, $\circ^{k\,\Uf^{(n)}_{\sig{\Sigma}}(\mathbf{A})}_{s}$ by $\circ^{(k,n)\,\mathbf{A}}_{s}$ and $\sigma^{\Uf^{(n)}_{\sig{\Sigma}}(\mathbf{A})}$ by $\sigma_{n}^{\mathbf{A}}$.
The translation of (A0) follows from (EA0), for every operation symbol of $\Sigma$, including when $n=0$.
Fix $k\in n$ and $s\in S$.
The translations of (A$^{k}$1) and (A$^{k}$3) follow from (EA$^{k,n}$1) and (EA$^{k,n}$3), together with Lemma~\ref{LUfIdInj}.
The translation of (A$^{k}$2) follows from Lemma~\ref{LUfIdInj}.
The translations of (A$^{k}$4), (A$^{k}$5), (A$^{k}$6) and (A$^{k}$8) follow from (EA$^{k,n}$4), (EA$^{k,n}$5), (EA$^{k,n}$6) and (EA$^{k,n}$8).
The translation of (A$^{k}$7) follows from the first and the third equations of~(EA$^{k,n}$7).

We next consider the axioms that relate two dimensions, fixing $l<k$ in $n$ and $s\in S$.
By~(EA$^{l,k,n}$1) and Lemma~\ref{LUfIdInj}, for every $k$-cell $w$ of sort $s$, $\src^{(l,n)\,\mathbf{A}}_{s}(\mathrm{i}^{(n,k)\,\mathbf{A}}_{s}(w))=\src^{(l,k)\,\mathbf{A}}_{s}(w)$, and similarly for the targets.
With this, the translation of (AB$^{l,k}$1) follows from~(EA$^{l,k,n}$1) and~(EA$^{l,k,n}$2).
The translation of (AB$^{l,k}$3) is the axiom (EA$^{l,k,n}$4), after exchanging $y$ and $z$.
Finally, for $n$-cells $x$, $y$ of sort $s$ such that $x\circ^{(l,n)\,\mathbf{A}}_{s}y$ is defined, the following chain of equalities holds
\begin{flushleft}
$
\mathrm{i}^{(n,k)\,\mathbf{A}}_{s}\left(\src^{(k,n)\,\mathbf{A}}_{s}\left(x\circ^{(l,n)\,\mathbf{A}}_{s}y\right)\right)
$
\begin{align*}
&=
\mathrm{i}^{(n,k)\,\mathbf{A}}_{s}\left(\src^{(k,n)\,\mathbf{A}}_{s}(x)\circ^{(l,k)\,\mathbf{A}}_{s}\src^{(k,n)\,\mathbf{A}}_{s}(y)\right)
\tag{1}
\\&=
\mathrm{i}^{(n,k)\,\mathbf{A}}_{s}\left(\src^{(k,n)\,\mathbf{A}}_{s}(x)\right)\circ^{(l,n)\,\mathbf{A}}_{s}\mathrm{i}^{(n,k)\,\mathbf{A}}_{s}\left(\src^{(k,n)\,\mathbf{A}}_{s}(y)\right).
\tag{2}
\end{align*}
\end{flushleft}

The first equality follows from~(EA$^{l,k,n}$3) and the second from~(EA$^{l,k,n}$5).
This is the translation of (AB$^{l,k}$2) for the sources, and the case of the targets is similar.
This completes the proof.
\end{proof}

From Proposition~\ref{PUfEpVariety} we obtain the unification functor between the categories of higher-order categorial algebras.

\begin{corollary}\label{CUfFunctor}\lean{equivalent}{Presentations.unificationFunctor}
The restriction of $\Uf^{(n)}_{\sig{\Sigma}}$ to $\mathsf{nCat}^{\htg}\mathsf{Alg}(\sig{\Sigma})$ corestricts to $\mathsf{nCatAlg}(\sig{\Sigma})$.
\end{corollary}
\begin{proof}
It follows from Proposition~\ref{PUfEpVariety} and Remarks~\ref{RCatAlgAsEpAlg} and~\ref{RCatAlgAsEpAlgHtg}.
\end{proof}

By abuse of notation, we denote the corestriction of Corollary~\ref{CUfFunctor} by
$$
\Uf^{(n)}_{\sig{\Sigma}}\colon\mathsf{nCat}^{\htg}\mathsf{Alg}(\sig{\Sigma})\mor\mathsf{nCatAlg}(\sig{\Sigma}).
$$
Explicitly, by Lemma~\ref{LUfDescription}, the cells of sort $s$ of $\Uf^{(n)}_{\sig{\Sigma}}(\calg{A})$ are the $n$-cells of sort $s$ of $\calg{A}$, $\src^{k}_{s}=\mathrm{i}^{(n,k)}_{s}\circ\src^{(k,n)}_{s}$, $\tgt^{k}_{s}=\mathrm{i}^{(n,k)}_{s}\circ\tgt^{(k,n)}_{s}$, $\circ^{k}_{s}=\circ^{(k,n)}_{s}$ and $\sigma$ acts as $\sigma_{n}$; and, by Proposition~\ref{PGenHepDerivorFunctor}, $\Uf^{(n)}_{\sig{\Sigma}}(F)$ is the family of the components of $F$ at the $n$-cells.

\section{On the relationship between graduation and unification}\label{SGdUfRel}

In this section we investigate the relation between the functors $\Gd^{(n)}_{\sig{\Sigma}}$ and $\Uf^{(n)}_{\sig{\Sigma}}$.
\begin{center}
	\begin{tikzpicture}
		[ACliment/.style={-{To [angle'=45, length=5.75pt, width=4pt, round]}},
		scale=1]

		\node[] (0R) at (0.5,0) [color=white] {};
		\node[] (1L) at (2.3,0)  [color=white] {};

		\node[] (C0) at (-0.5,0) [] {$\mathsf{nCatAlg}(\sig{\Sigma})$};
		\node[] (C1) at (3.5,0) [] {$\mathsf{nCat}^{\htg}\mathsf{Alg}(\sig{\Sigma})$};

		\draw[ACliment] (0R.north east) to node [above] {$\Gd^{(n)}_{\sig{\Sigma}}$} (1L.north west);
		\draw[ACliment] (1L.south west) to node [below] {$\Uf^{(n)}_{\sig{\Sigma}}$} (0R.south east);
	\end{tikzpicture}
\end{center}

Both functors are induced by generalized hep-derivors.
Therefore we compute their composites by means of the compositions $\bGd^{(n)}_{\sig{\Sigma}}\circ\bUf^{(n)}_{\sig{\Sigma}}$ and $\bUf^{(n)}_{\sig{\Sigma}}\circ\bGd^{(n)}_{\sig{\Sigma}}$, see Definition~\ref{DGenHepDerivorCompId}, Proposition~\ref{PGenHepDerivorComp} and Corollary~\ref{CGenHepDerivorCompEqual}.
We start by proving that the composite of the unification derivor and the graduation derivor is the identity.

\begin{proposition}\label{PGdUfDerivorId}\lean{equivalent}{Presentations.unification_comp_graduation}
The following equality holds in $\mathsf{HSig}_{\gd}$
$$
\bGd^{(n)}_{\sig{\Sigma}}\circ\bUf^{(n)}_{\sig{\Sigma}}=\id_{\sig{\Sigma}^{(n)}}.
$$
\end{proposition}
\begin{proof}
By Definition~\ref{DGenHepDerivorCompId}, we compare the underlying derivor and the domains of $\bGd^{(n)}_{\sig{\Sigma}}\circ\bUf^{(n)}_{\sig{\Sigma}}$ with those of $\id_{\sig{\Sigma}^{(n)}}$.
\begin{claim}\label{ClGdUfDerivorIdA}
The underlying derivor of $\bGd^{(n)}_{\sig{\Sigma}}\circ\bUf^{(n)}_{\sig{\Sigma}}$ is the identity derivor at $(S,\Sigma^{(n)})$.
\end{claim}
The underlying derivor is $(\varphi^{\Gd^{(n)}_{\sig{\Sigma}}}\circ\varphi^{\Uf^{(n)}_{\sig{\Sigma}}},(\Gd^{(n)}_{\sig{\Sigma}})^{\sharp}\circ\Uf^{(n)}_{\sig{\Sigma}})$.
For every $s\in S$, $\varphi^{\Gd^{(n)}_{\sig{\Sigma}}}(\varphi^{\Uf^{(n)}_{\sig{\Sigma}}}(s))=\varphi^{\Gd^{(n)}_{\sig{\Sigma}}}(n,s)=s$.
By~(Der1)--(Der2) and Definitions~\ref{DGdDerivor} and~\ref{DUfDerivor}, for every $(\mathbf{s},s)\in S^{\star}\times S$, every $\sigma\in\Sigma_{\mathbf{s},s}$, every $k\in n$ and every $s\in S$,
\begin{align*}
(\Gd^{(n)}_{\sig{\Sigma}})^{\sharp}\left(\sigma_{n}\left((v^{(n,\mathbf{s})}_{p})_{p\in\bb{\mathbf{s}}}\right)\right)&=\sigma\left((v^{\mathbf{s}}_{p})_{p\in\bb{\mathbf{s}}}\right),\\
(\Gd^{(n)}_{\sig{\Sigma}})^{\sharp}\left(\mathrm{i}^{(n,k)}_{s}\left(\src^{(k,n)}_{s}(v^{(n,s)}_{0})\right)\right)&=\src^{k}_{s}(v^{s}_{0}),\\
(\Gd^{(n)}_{\sig{\Sigma}})^{\sharp}\left(\mathrm{i}^{(n,k)}_{s}\left(\tgt^{(k,n)}_{s}(v^{(n,s)}_{0})\right)\right)&=\tgt^{k}_{s}(v^{s}_{0}),\\
(\Gd^{(n)}_{\sig{\Sigma}})^{\sharp}\left(v^{(n,s)}_{0}\circ^{(k,n)}_{s}v^{(n,s)}_{1}\right)&=v^{s}_{0}\circ^{k}_{s}v^{s}_{1},
\end{align*}
since $\Gd^{(n)}_{\sig{\Sigma}}$ sends $\mathrm{i}^{(n,k)}_{s}$ to $v^{s}_{0}$.
These are the terms that the identity derivor assigns to $\sigma$, $\src^{k}_{s}$, $\tgt^{k}_{s}$ and $\circ^{k}_{s}$, see Definition~\ref{DDerivorCompId}.

\begin{claim}\label{ClGdUfDerivorIdB}
The domains of $\bGd^{(n)}_{\sig{\Sigma}}\circ\bUf^{(n)}_{\sig{\Sigma}}$ are empty.
\end{claim}
For every $s\in S$, by Definitions~\ref{DGenHepDerivorCompId}, \ref{DGdDerivor} and~\ref{DUfDerivor},
\[
\Dom^{\bGd^{(n)}_{\sig{\Sigma}}\circ\bUf^{(n)}_{\sig{\Sigma}}}_s
=\Dom^{\Gd^{(n)}_{\sig{\Sigma}}}_{(n,s)}
 \cup(\Gd^{(n)}_{\sig{\Sigma}})^{\sharp}\left[\Dom^{\Uf^{(n)}_{\sig{\Sigma}}}_s\right]
=\varnothing.
\]
Since the domains of $\id_{\sig{\Sigma}^{(n)}}$ are also empty, the two claims give the stated equality, also for $n=0$.
\end{proof}

From the above proposition we obtain that graduation followed by unification is the identity functor.

\begin{corollary}\label{CUfGdId}\lean{weaker}{Presentations.homogeneousPresentationRoundTrip}
The following equality holds
$$
\Uf^{(n)}_{\sig{\Sigma}}\circ\Gd^{(n)}_{\sig{\Sigma}}=\Id_{\mathsf{nCatAlg}(\sig{\Sigma})}.
$$
\end{corollary}
\begin{proof}
Let $\calg{A}$ be an $n$-categorial $\sig{\Sigma}$-algebra.
By Lemmas~\ref{LGdDescription} and~\ref{LUfDescription}, the cells of sort $s$ of $\Uf^{(n)}_{\sig{\Sigma}}(\Gd^{(n)}_{\sig{\Sigma}}(\calg{A}))$ are the cells of sort $(n,s)$ of $\Gd^{(n)}_{\sig{\Sigma}}(\calg{A})$, i.e., all the cells of sort $s$ of $\calg{A}$.
By Proposition~\ref{PGdUfDerivorId} and Remark~\ref{RDomainSetEmpty}, these are also the cells of sort $s$ of $(\bGd^{(n)}_{\sig{\Sigma}}\circ\bUf^{(n)}_{\sig{\Sigma}})^{\ast}_{\gd}(\calg{A})$.
Moreover, by Propositions~\ref{PGdEpVariety} and~\ref{PUfEpVariety}, $\Uf^{(n)}_{\sig{\Sigma}}(\Gd^{(n)}_{\sig{\Sigma}}(\calg{A}))$ is a $\sig{\Sigma}^{(n)}$-algebra.
These are the hypotheses of Corollary~\ref{CGenHepDerivorCompEqual}; the equality of derivors of Proposition~\ref{PGdUfDerivorId} alone would not suffice, since, in general, the composite of the transformations is only included in the transformation along the composite, see Proposition~\ref{PGenHepDerivorComp}.
Hence, by Corollary~\ref{CGenHepDerivorCompEqual} and Propositions~\ref{PGdUfDerivorId} and~\ref{PGenHepDerivorId},
$$
\Uf^{(n)}_{\sig{\Sigma}}\left(\Gd^{(n)}_{\sig{\Sigma}}(\calg{A})\right)
=\left(\id_{\sig{\Sigma}^{(n)}}\right)^{\ast}_{\gd}(\calg{A})
=\calg{A},
$$
and, by the same results, $\Uf^{(n)}_{\sig{\Sigma}}(\Gd^{(n)}_{\sig{\Sigma}}(F))=F$ for every $n$-categorial $\sig{\Sigma}$-homomorphism $F$.
\end{proof}

In contrast to the previous result, the process of unification followed by that of graduation does not give rise to the identity functor on the category of heterogeneous $n$-categorial algebras, since the $k$-cells of $\Gd^{(n)}_{\sig{\Sigma}}(\Uf^{(n)}_{\sig{\Sigma}}(\calg{A}))$ are the $(n,k)$-identities of the $k$-cells of $\calg{A}$.
However, $\Gd^{(n)}_{\sig{\Sigma}}\circ\Uf^{(n)}_{\sig{\Sigma}}$ is naturally isomorphic to the identity functor.

\begin{proposition}\label{PGdUfIso}\lean{equivalent}{Presentations.theta,Presentations.heterogeneousPresentationTheta}
For every heterogeneous $n$-categorial $\sig{\Sigma}$-algebra $\calg{A}$, let $\theta_{\calg{A}}$ be the family whose component at the sort $(k,s)\in(n+1)\times S$ sends $u\in A_{(k,s)}$ to $\mathrm{i}^{(n,k)\,\calg{A}}_{s}(u)$, where $\mathrm{i}^{(n,n)\,\calg{A}}_{s}$ is the identity of $A_{(n,s)}$.
Then $\theta=(\theta_{\calg{A}})_{\calg{A}\in\mathrm{Ob}(\mathsf{nCat}^{\htg}\mathsf{Alg}(\sig{\Sigma}))}$ is a natural isomorphism from $\Id_{\mathsf{nCat}^{\htg}\mathsf{Alg}(\sig{\Sigma})}$ to $\Gd^{(n)}_{\sig{\Sigma}}\circ\Uf^{(n)}_{\sig{\Sigma}}$.
Thus
\[\Gd^{(n)}_{\sig{\Sigma}}\circ\Uf^{(n)}_{\sig{\Sigma}}\cong \Id_{\mathsf{nCat}^{\htg}\mathsf{Alg}(\sig{\Sigma})}.\]
\end{proposition}
\begin{proof}
Let $\calg{A}$ be a heterogeneous $n$-categorial $\sig{\Sigma}$-algebra.
In this proof, for every $s\in S$, we let $\mathrm{i}^{(n,n)\,\calg{A}}_{s}$ and $\src^{(n,n)\,\calg{A}}_{s}$ stand for the identity mapping of $A_{(n,s)}$, so that the formulas below also cover the dimension $k=n$.
We first describe the composite generalized hep-derivor $\bUf^{(n)}_{\sig{\Sigma}}\circ\bGd^{(n)}_{\sig{\Sigma}}$ from $\sig{\Sigma}^{\htg(n)}$ to $\sig{\Sigma}^{\htg(n)}$, see Definition~\ref{DGenHepDerivorCompId}.

\begin{claim}\label{ClGdUfIsoA}
The mapping on sorts of $\bUf^{(n)}_{\sig{\Sigma}}\circ\bGd^{(n)}_{\sig{\Sigma}}$ sends $(k,s)$ to $(n,s)$, for every $k\in n+1$ and every $s\in S$.
Moreover, for every $(\mathbf{s},s)\in S^{\star}\times S$, every $\sigma\in\Sigma_{\mathbf{s},s}$, every $j<k$ in $n+1$ and every $s\in S$, its terms for $\sigma_k$, $\src^{(j,k)}_s$, $\tgt^{(j,k)}_s$, $\mathrm{i}^{(k,j)}_s$ and $\circ^{(j,k)}_s$ are, respectively,
\begin{gather*}
\sigma_n\left((v^{(n,\mathbf{s})}_p)_{p\in\bb{\mathbf{s}}}\right),\qquad
\mathrm{i}^{(n,j)}_s\left(\src^{(j,n)}_s(v^{(n,s)}_{0})\right),\qquad
\mathrm{i}^{(n,j)}_s\left(\tgt^{(j,n)}_s(v^{(n,s)}_{0})\right),\\
v^{(n,s)}_{0},\qquad v^{(n,s)}_{0}\circ^{(j,n)}_s v^{(n,s)}_{1}.
\end{gather*}
\end{claim}
It follows from Definitions~\ref{DGdDerivor}, \ref{DUfDerivor} and~\ref{DGenHepDerivorCompId}, by~(Der1)--(Der2).
For instance, $(\Uf^{(n)}_{\sig{\Sigma}})^{\sharp}(\Gd^{(n)}_{\sig{\Sigma}}(\src^{(j,k)}_s))=(\Uf^{(n)}_{\sig{\Sigma}})^{\sharp}(\src^{j}_{s}(v^{s}_{0}))=\mathrm{i}^{(n,j)}_s(\src^{(j,n)}_s(v^{(n,s)}_{0}))$.

\begin{claim}\label{ClGdUfIsoB}
For every $k\in n+1$ and every $s\in S$, the cells of sort $(k,s)$ of $\Gd^{(n)}_{\sig{\Sigma}}(\Uf^{(n)}_{\sig{\Sigma}}(\calg{A}))$ are the elements of the image of $\mathrm{i}^{(n,k)\,\calg{A}}_s$, and the operations of $\Gd^{(n)}_{\sig{\Sigma}}(\Uf^{(n)}_{\sig{\Sigma}}(\calg{A}))$ are the restrictions of the partial term operations in $\calg{A}$ of the terms of Claim~\ref{ClGdUfIsoA}.
\end{claim}
By Lemmas~\ref{LGdDescription} and~\ref{LUfDescription}, for $k\in n$, the cells of sort $(k,s)$ of $\Gd^{(n)}_{\sig{\Sigma}}(\Uf^{(n)}_{\sig{\Sigma}}(\calg{A}))$ are the $n$-cells $x$ of sort $s$ of $\calg{A}$ such that $x=\mathrm{i}^{(n,k)\,\calg{A}}_s(\src^{(k,n)\,\calg{A}}_s(x))$, and, for $k=n$, they are all the $n$-cells of sort $s$.
By Lemma~\ref{LUfIdInj}, the $n$-cells $x$ such that $x=\mathrm{i}^{(n,k)\,\calg{A}}_s(\src^{(k,n)\,\calg{A}}_s(x))$ are exactly the elements of the image of $\mathrm{i}^{(n,k)\,\calg{A}}_s$.
By Definition~\ref{DGenHepDerivorCompId}, the same condition defines the domain of sort $(k,s)$ of $\bUf^{(n)}_{\sig{\Sigma}}\circ\bGd^{(n)}_{\sig{\Sigma}}$, so that the iterated and the composite transformations of $\calg{A}$ have the same cells.
Moreover, by Propositions~\ref{PUfEpVariety} and~\ref{PGdEpVariety}, $\Gd^{(n)}_{\sig{\Sigma}}(\Uf^{(n)}_{\sig{\Sigma}}(\calg{A}))$ is a $\sig{\Sigma}^{\htg(n)}$-algebra.
Hence, by Corollary~\ref{CGenHepDerivorCompEqual}, it is the transformation of $\calg{A}$ along the composite, and, by Proposition~\ref{PTransformedAlg}, its operations, including their domains, are as stated.

For every $k\in n+1$ and every $s\in S$, we let $\theta_{\calg{A},(k,s)}$ stand for the component of $\theta_{\calg{A}}$ at $(k,s)$, regarded as a mapping from $A_{(k,s)}$ to the cells of sort $(k,s)$ of $\Gd^{(n)}_{\sig{\Sigma}}(\Uf^{(n)}_{\sig{\Sigma}}(\calg{A}))$.

\begin{claim}\label{ClGdUfIsoC}
$\theta_{\calg{A}}$ is an $(n+1)\times S$-sorted bijection whose inverse sends, at every sort $(k,s)$, a cell $x$ to $\src^{(k,n)\,\calg{A}}_s(x)$.
\end{claim}
For every $k$-cell $u$ of sort $s$, $\src^{(k,n)\,\calg{A}}_s(\mathrm{i}^{(n,k)\,\calg{A}}_s(u))=u$, by Lemma~\ref{LUfIdInj}; and, for every cell $x$ of sort $(k,s)$ of $\Gd^{(n)}_{\sig{\Sigma}}(\Uf^{(n)}_{\sig{\Sigma}}(\calg{A}))$, $\mathrm{i}^{(n,k)\,\calg{A}}_s(\src^{(k,n)\,\calg{A}}_s(x))=x$, by Claim~\ref{ClGdUfIsoB}.

\begin{claim}\label{ClGdUfIsoD}
$\theta_{\calg{A}}$ preserves all the operations and reflects the definedness of the compositions; hence it is an isomorphism of partial $\Sigma^{\htg(n)}$-algebras from $\calg{A}$ to $\Gd^{(n)}_{\sig{\Sigma}}(\Uf^{(n)}_{\sig{\Sigma}}(\calg{A}))$.
\end{claim}
By Claims~\ref{ClGdUfIsoA} and~\ref{ClGdUfIsoB}, we have to compare the operations of $\calg{A}$ with the term operations of Claim~\ref{ClGdUfIsoA}.
For $k\in n$, the identities $\mathrm{i}^{(n,k)\,\calg{A}}_s$ commute with the operations of $\Sigma$, by~(EA$^{k,n}$7).
For $j<k$ in $n$, by Lemma~\ref{LUfIdInj} and~(EA$^{j,k,n}$1)--(EA$^{j,k,n}$2), the following equalities hold:
\begin{align*}
\src^{(j,n)\,\calg{A}}_s\circ\mathrm{i}^{(n,k)\,\calg{A}}_s&=\src^{(j,k)\,\calg{A}}_s, &
\tgt^{(j,n)\,\calg{A}}_s\circ\mathrm{i}^{(n,k)\,\calg{A}}_s&=\tgt^{(j,k)\,\calg{A}}_s,
\\
\mathrm{i}^{(n,k)\,\calg{A}}_s\circ\mathrm{i}^{(k,j)\,\calg{A}}_s&=\mathrm{i}^{(n,j)\,\calg{A}}_s.
\end{align*}
For $k=n$, these equalities hold trivially, since $\mathrm{i}^{(n,n)\,\calg{A}}_s$ is the identity mapping.
Together with Claim~\ref{ClGdUfIsoA}, they show that $\theta_{\calg{A}}$ preserves the operations $\sigma_k$, the sources, the targets and the identities.

For $j<k$ in $n+1$ and $k$-cells $u$, $u'$ of sort $s$, the pair $(\theta_{\calg{A},(k,s)}(u),\theta_{\calg{A},(k,s)}(u'))$ is in the domain of the $j$-composition of $\Gd^{(n)}_{\sig{\Sigma}}(\Uf^{(n)}_{\sig{\Sigma}}(\calg{A}))$ if, and only if, by the above equalities,
\[
\mathrm{i}^{(n,j)\,\calg{A}}_s\left(\src^{(j,k)\,\calg{A}}_s(u)\right)
=\mathrm{i}^{(n,j)\,\calg{A}}_s\left(\tgt^{(j,k)\,\calg{A}}_s(u')\right),
\]
i.e., by Lemma~\ref{LUfIdInj}, if, and only if, $u\circ^{(j,k)\,\calg{A}}_s u'$ is defined.
For such a pair, by~(EA$^{j,k,n}$5) if $k\in n$, and trivially if $k=n$,
\[
\mathrm{i}^{(n,k)\,\calg{A}}_s\left(u\circ^{(j,k)\,\calg{A}}_s u'\right)
=\mathrm{i}^{(n,k)\,\calg{A}}_s(u)\circ^{(j,n)\,\calg{A}}_s\mathrm{i}^{(n,k)\,\calg{A}}_s(u').
\]
All the other operations are total, so that $\theta_{\calg{A}}$ is a homomorphism whose inverse, by Claim~\ref{ClGdUfIsoC}, is also a homomorphism.
Hence, by Proposition~\ref{PPartialAlgIso}, $\theta_{\calg{A}}$ is an isomorphism.

\begin{claim}\label{ClGdUfIsoE}
The family $\theta=(\theta_{\calg{A}})_{\calg{A}}$ is a natural transformation from $\Id_{\mathsf{nCat}^{\htg}\mathsf{Alg}(\sig{\Sigma})}$ to $\Gd^{(n)}_{\sig{\Sigma}}\circ\Uf^{(n)}_{\sig{\Sigma}}$.
\end{claim}
Let $F\colon\calg{A}\mor\calg{A}'$ be a heterogeneous $n$-categorial $\sig{\Sigma}$-homomorphism.
By Proposition~\ref{PGenHepDerivorFunctor}, the component at the sort $(k,s)$ of $\Gd^{(n)}_{\sig{\Sigma}}(\Uf^{(n)}_{\sig{\Sigma}}(F))$ is the restriction of $F_{(n,s)}$.
Hence naturality amounts to the equality $F_{(n,s)}\circ\mathrm{i}^{(n,k)\,\calg{A}}_s=\mathrm{i}^{(n,k)\,\calg{A}'}_s\circ F_{(k,s)}$, which holds because $F$ preserves the identities.

By Claims~\ref{ClGdUfIsoD} and~\ref{ClGdUfIsoE}, $\theta$ is a natural isomorphism. 
\end{proof}

We prove that graduation is left adjoint to unification, with the identity as unit.

\begin{proposition}\label{PGdUfAdj}\lean{weaker}{Presentations.presentationAdjunction}
Let $\theta$ be the natural isomorphism from $\Id_{\mathsf{nCat}^{\htg}\mathsf{Alg}(\sig{\Sigma})}$ to $\Gd^{(n)}_{\sig{\Sigma}}\circ\Uf^{(n)}_{\sig{\Sigma}}$ of Proposition~\ref{PGdUfIso}.
Then $(\Uf^{(n)}_{\sig{\Sigma}},\Gd^{(n)}_{\sig{\Sigma}},\id,\theta^{-1})$ is an adjunction from $\mathsf{nCatAlg}(\sig{\Sigma})$ to $\mathsf{nCat}^{\htg}\mathsf{Alg}(\sig{\Sigma})$, see Definition~\ref{DAdj}, i.e., $\Gd^{(n)}_{\sig{\Sigma}}\dashv\Uf^{(n)}_{\sig{\Sigma}}$, with unit the identity and counit $\theta^{-1}=(\theta^{-1}_{\calg{A}})_{\calg{A}\in\mathrm{Ob}(\mathsf{nCat}^{\htg}\mathsf{Alg}(\sig{\Sigma}))}$.
\end{proposition}
\begin{proof}
By Corollary~\ref{CUfGdId}, $\Uf^{(n)}_{\sig{\Sigma}}\circ\Gd^{(n)}_{\sig{\Sigma}}=\Id_{\mathsf{nCatAlg}(\sig{\Sigma})}$, so that the identity $\id=(\id_{\calg{B}})_{\calg{B}\in\mathrm{Ob}(\mathsf{nCatAlg}(\sig{\Sigma}))}$ is a natural transformation from $\Id_{\mathsf{nCatAlg}(\sig{\Sigma})}$ to $\Uf^{(n)}_{\sig{\Sigma}}\circ\Gd^{(n)}_{\sig{\Sigma}}$; and, by Proposition~\ref{PGdUfIso}, $\theta^{-1}$ is a natural transformation from $\Gd^{(n)}_{\sig{\Sigma}}\circ\Uf^{(n)}_{\sig{\Sigma}}$ to $\Id_{\mathsf{nCat}^{\htg}\mathsf{Alg}(\sig{\Sigma})}$.
By Definition~\ref{DAdj}, we have to prove the triangular identities
\begin{enumerate}
	\item[(TI1)]
	$\left(\Uf^{(n)}_{\sig{\Sigma}}\ast\theta^{-1}\right)\circ\left(\id\ast\Uf^{(n)}_{\sig{\Sigma}}\right)=\id_{\Uf^{(n)}_{\sig{\Sigma}}}$; and

	\item[(TI2)]
	$\left(\theta^{-1}\ast\Gd^{(n)}_{\sig{\Sigma}}\right)\circ\left(\Gd^{(n)}_{\sig{\Sigma}}\ast\id\right)=\id_{\Gd^{(n)}_{\sig{\Sigma}}}$.
\end{enumerate}
Since the unit is the identity, the triangular identities state that, for every heterogeneous $n$-categorial $\sig{\Sigma}$-algebra $\calg{A}$ and every $n$-categorial $\sig{\Sigma}$-algebra $\calg{B}$,
\begin{enumerate}
	\item[(TI1)]
	$\Uf^{(n)}_{\sig{\Sigma}}(\theta^{-1}_{\calg{A}})$ is the identity of $\Uf^{(n)}_{\sig{\Sigma}}(\calg{A})$; and

	\item[(TI2)]
	$\theta^{-1}_{\Gd^{(n)}_{\sig{\Sigma}}(\calg{B})}$ is the identity of $\Gd^{(n)}_{\sig{\Sigma}}(\calg{B})$.
\end{enumerate}
We recall that a homomorphism is determined by its domain, its codomain and its underlying many-sorted mapping.
We also recall, from Proposition~\ref{PGdUfIso}, that the component of $\theta_{\calg{A}}$ at the sort $(k,s)$ is $\mathrm{i}^{(n,k)\,\calg{A}}_{s}$.

(TI1) By Proposition~\ref{PGenHepDerivorFunctor}, $\Uf^{(n)}_{\sig{\Sigma}}(\theta^{-1}_{\calg{A}})$ is the family of the components of $\theta^{-1}_{\calg{A}}$ at the sorts $(n,s)$.
By Proposition~\ref{PGdUfIso}, the components of $\theta_{\calg{A}}$ at these sorts are identity mappings, hence so are those of $\theta^{-1}_{\calg{A}}$.

(TI2) By Lemma~\ref{LGdDescription}~(1), the $(n,k)$-identities of $\Gd^{(n)}_{\sig{\Sigma}}(\calg{B})$ are inclusions.
Since the cells of sort $(k,s)$ and $(n,s)$ of $\Gd^{(n)}_{\sig{\Sigma}}(\Uf^{(n)}_{\sig{\Sigma}}(\Gd^{(n)}_{\sig{\Sigma}}(\calg{B})))=\Gd^{(n)}_{\sig{\Sigma}}(\calg{B})$ are cells of $\calg{B}$, $\theta_{\Gd^{(n)}_{\sig{\Sigma}}(\calg{B})}$ is the identity.
Hence so is its inverse. This completes the proof.
\end{proof}

From Corollary~\ref{CUfGdId} and Propositions~\ref{PGdUfIso} and~\ref{PGdUfAdj} we obtain the following corollary.

\begin{corollary}\label{CGdUfEquiv}\lean{weaker}{Presentations.presentationEquivalence}
The functors $\Gd^{(n)}_{\sig{\Sigma}}$ and $\Uf^{(n)}_{\sig{\Sigma}}$ are mutually quasi-inverse equivalences of categories between $\mathsf{nCatAlg}(\sig{\Sigma})$ and $\mathsf{nCat}^{\htg}\mathsf{Alg}(\sig{\Sigma})$, and the adjunction of Proposition~\ref{PGdUfAdj} is an adjoint equivalence.
\end{corollary}
\begin{proof}
By Corollary~\ref{CUfGdId} and Proposition~\ref{PGdUfIso}, $\Uf^{(n)}_{\sig{\Sigma}}\circ\Gd^{(n)}_{\sig{\Sigma}}=\Id_{\mathsf{nCatAlg}(\sig{\Sigma})}$ and $\Gd^{(n)}_{\sig{\Sigma}}\circ\Uf^{(n)}_{\sig{\Sigma}}\cong\Id_{\mathsf{nCat}^{\htg}\mathsf{Alg}(\sig{\Sigma})}$.
Moreover, the unit and the counit of the adjunction of Proposition~\ref{PGdUfAdj} are natural isomorphisms, so that it is an adjoint equivalence.
\end{proof}

\chapter[Discretization and truncation]{Discretization and truncation of higher-order categorial algebras}\label{SDsTr}

Throughout this chapter, $\sig{\Sigma}=(S,\Sigma)$ is a many-sorted signature, which is fixed in all the definitions and results of this chapter; in Sections~\ref{SDs}--\ref{SAlgSquare}, $m,n\in\mathbb{N}$ with $m\leq n$ are also fixed.
In this chapter we move between dimensions.
The discretization functor regards an $m$-categorial $\sig{\Sigma}$-algebra as an $n$-categorial $\sig{\Sigma}$-algebra whose cells of dimension greater than $m$ are identities, and the truncation functor forgets the cells of dimension greater than $m$ of an $n$-categorial $\sig{\Sigma}$-algebra.
In Sections~\ref{SDs} and~\ref{STr} we define the discretization and the truncation functors and prove that they compose as expected, in Section~\ref{SDsTrRel} we prove that discretization is left adjoint to truncation and compute the composite of a discretization followed by a truncation, and in Section~\ref{SDsTrCompatGdUf} we show that they are compatible with graduation and unification.
In Section~\ref{SAlgSquare} we show that graduation, unification, discretization and truncation form an algebraic square.
Finally, in Section~\ref{SDsTrOmega} we define the truncation and the discretization functors between dimension $\omega$ and the finite dimensions.

As in Chapter~\ref{SGdUf}, all these functors are induced by generalized hep-derivors, in the homogeneous and in the heterogeneous presentations, and we use the notation for generalized hep-derivors introduced at the beginning of Chapter~\ref{SGdUf}.

\section{Discretization of higher-order categorial algebras}\label{SDs}

We begin by defining the discretization derivors, in the homogeneous and in the heterogeneous presentations.
In the heterogeneous presentation, a cell of dimension $k$ of the discretization is a cell of dimension $k$ of the given algebra if $k\leq m$, and a cell of dimension $m$ if $k>m$.
Accordingly, every dimension $k\in n+1$ of $\Sigma^{\htg(n)}$ is sent to the dimension $\min\{k,m\}\in m+1$ of $\Sigma^{\htg(m)}$.

\begin{definition}\label{DDsDerivor}\lean{equivalent}{Dimensions.homogeneousDiscretization,Dimensions.heterogeneousDiscretization}
If $m<n$, we define the following generalized hep-derivors.
\begin{enumerate}
\item The \emph{homogeneous discretization derivor} is the generalized hep-derivor
$$
\bDs^{(n,m)}_{\sig{\Sigma}}=\left(\varphi^{\Ds^{(n,m)}_{\sig{\Sigma}}},\Ds^{(n,m)}_{\sig{\Sigma}},\Dom^{\Ds^{(n,m)}_{\sig{\Sigma}}}\right)
$$
from $\sig{\Sigma}^{(n)}$ to $\sig{\Sigma}^{(m)}$ defined as follows:
\begin{enumerate}
\item the mapping on sorts is $\varphi^{\Ds^{(n,m)}_{\sig{\Sigma}}}=\id_{S}$;
\item for every $(\mathbf{s},s)\in S^{\star}\times S$ and every $\sigma\in\Sigma_{\mathbf{s},s}$,
$$
\Ds^{(n,m)}_{\sig{\Sigma}}(\sigma)=\sigma((v^{\mathbf{s}}_{p})_{p\in\bb{\mathbf{s}}});
$$
\item for every $k\in m$ and every $s\in S$,
\begin{align*}
\Ds^{(n,m)}_{\sig{\Sigma}}(\src^{k}_{s})&=\src^{k}_{s}(v^{s}_{0}),\\
\Ds^{(n,m)}_{\sig{\Sigma}}(\tgt^{k}_{s})&=\tgt^{k}_{s}(v^{s}_{0}),\\
\Ds^{(n,m)}_{\sig{\Sigma}}(\circ^{k}_{s})&=v^{s}_{0}\circ^{k}_{s}v^{s}_{1};
\end{align*}
\item for every $k\in n$ with $k\geq m$ and every $s\in S$,
$$
\Ds^{(n,m)}_{\sig{\Sigma}}(\src^{k}_{s})=\Ds^{(n,m)}_{\sig{\Sigma}}(\tgt^{k}_{s})=\Ds^{(n,m)}_{\sig{\Sigma}}(\circ^{k}_{s})=v^{s}_{0};
$$
and
\item for every $s\in S$, the domain is $\Dom^{\Ds^{(n,m)}_{\sig{\Sigma}}}_{s}=\varnothing$.
\end{enumerate}
\item The \emph{heterogeneous discretization derivor} is the generalized hep-derivor
$$
\bDs^{\htg(n,m)}_{\sig{\Sigma}}=\left(\varphi^{\Ds^{\htg(n,m)}_{\sig{\Sigma}}},\Ds^{\htg(n,m)}_{\sig{\Sigma}},\Dom^{\Ds^{\htg(n,m)}_{\sig{\Sigma}}}\right)
$$
from $\sig{\Sigma}^{\htg(n)}$ to $\sig{\Sigma}^{\htg(m)}$ defined as follows:
\begin{enumerate}
\item for every dimension $k\in n+1$ and every $s\in S$, the mapping on sorts sends $(k,s)$ to
$$
\varphi^{\Ds^{\htg(n,m)}_{\sig{\Sigma}}}(k,s)=(\min\{k,m\},s);
$$
\item for every $(\mathbf{s},s)\in S^{\star}\times S$, every $\sigma\in\Sigma_{\mathbf{s},s}$ and every $k\in n+1$,
$$
\Ds^{\htg(n,m)}_{\sig{\Sigma}}(\sigma_{k})=\sigma_{\min\{k,m\}}\left(\left(v^{(\min\{k,m\},\mathbf{s})}_{p}\right)_{p\in\bb{\mathbf{s}}}\right);
$$
\item for every pair of dimensions $j<k$ in $n+1$ with $j<m$, so that $j<\min\{k,m\}$, and every $s\in S$,
\begin{align*}
\Ds^{\htg(n,m)}_{\sig{\Sigma}}(\src^{(j,k)}_{s})&=\src^{(j,\min\{k,m\})}_{s}\left(v^{(\min\{k,m\},s)}_{0}\right),\\
\Ds^{\htg(n,m)}_{\sig{\Sigma}}(\tgt^{(j,k)}_{s})&=\tgt^{(j,\min\{k,m\})}_{s}\left(v^{(\min\{k,m\},s)}_{0}\right),\\
\Ds^{\htg(n,m)}_{\sig{\Sigma}}(\mathrm{i}^{(k,j)}_{s})&=\mathrm{i}^{(\min\{k,m\},j)}_{s}\left(v^{(j,s)}_{0}\right),\\
\Ds^{\htg(n,m)}_{\sig{\Sigma}}(\circ^{(j,k)}_{s})&=v^{(\min\{k,m\},s)}_{0}\circ^{(j,\min\{k,m\})}_{s}v^{(\min\{k,m\},s)}_{1};
\end{align*}
\item for every pair of dimensions $j<k$ in $n+1$ with $j\geq m$, so that $\min\{j,m\}=\min\{k,m\}=m$, and every $s\in S$,
\begin{multline*}
\Ds^{\htg(n,m)}_{\sig{\Sigma}}(\src^{(j,k)}_{s})=
\Ds^{\htg(n,m)}_{\sig{\Sigma}}(\tgt^{(j,k)}_{s})=
\Ds^{\htg(n,m)}_{\sig{\Sigma}}(\mathrm{i}^{(k,j)}_{s})
\\=
\Ds^{\htg(n,m)}_{\sig{\Sigma}}(\circ^{(j,k)}_{s})=v^{(m,s)}_{0};
\end{multline*}
and
\item for every dimension $k\in n+1$ and every $s\in S$, the domain is
$$
\Dom^{\Ds^{\htg(n,m)}_{\sig{\Sigma}}}_{(k,s)}=\varnothing.
$$
\end{enumerate}
\end{enumerate}
For $m=n$, we let $\bDs^{(n,n)}_{\sig{\Sigma}}$ and $\bDs^{\htg(n,n)}_{\sig{\Sigma}}$ be the identities $\id_{\sig{\Sigma}^{(n)}}$ and $\id_{\sig{\Sigma}^{\htg(n)}}$, see Definition~\ref{DGenHepDerivorCompId}.
\end{definition}

Let us note that, in the homogeneous case, for every $k\in n$ with $k\geq m$ and every $s\in S$, $\Ds^{(n,m)}_{\sig{\Sigma}}(\circ^{k}_{s})$ is total but, in the transformed algebra, $\src^{k}_{s}$ and $\tgt^{k}_{s}$ act as identities, so that two cells satisfy the domain condition $\src^{k}_{s}(v^{s}_{0})=\tgt^{k}_{s}(v^{s}_{1})$ of $\circ^{k}_{s}$ if, and only if, they are equal, see Proposition~\ref{PTransformedAlg}.
Thus the $k$-composition of two cells is defined if, and only if, they are equal.
The next example shows that this restriction cannot in general be represented by a term operation alone.

\begin{example}\label{EDiscretizationDiagonal}\lean{equivalent}{Dimensions.Diagonal.composition_dom,Dimensions.Diagonal.no_term_domain}
Suppose that $S=\{s\}$, $\Sigma=\varnothing$, $m=0$ and $n=1$, so that a partial $\Sigma^{(0)}$-algebra is just a set $B$.
By Definition~\ref{DDsDerivor}, $\bDs^{(1,0)}_{\sig{\Sigma}}$ sends $\src^{0}_{s}$, $\tgt^{0}_{s}$ and $\circ^{0}_{s}$ to the variable $v^{s}_{0}$, and its domains are empty, so that every element of $B$ is a cell.
By Proposition~\ref{PTransformedAlg}, in $\Ds^{(1,0)}_{\sig{\Sigma}}(B)$ the source and the target are the identity of $B$, and the composition is defined at a pair $(x,y)$ if, and only if, $(x,y)$ satisfies the domain condition $\src^{0}_{s}(v^{(s,s)}_{0})\eeq\tgt^{0}_{s}(v^{(s,s)}_{1})$ of $\circ^{0}_{s}$ in the transformed algebra, i.e., if, and only if, $x=y$, and then $x\circ^{0}_{s}x=x$.
Thus $\Ds^{(1,0)}_{\sig{\Sigma}}(B)$ is the discrete category on $B$.
The term $v^{s}_{0}$ that interprets $\circ^{0}_{s}$ is defined at every pair, and, since $\Sigma^{(0)}$ has no operation symbols, every term operation of $B$ is total: the restriction of the composition to the diagonal is imposed by the domain condition alone.
\end{example}

In the heterogeneous case, on the dimensions $k\leq m$, $\bDs^{\htg(n,m)}_{\sig{\Sigma}}$ is the identity, and, for every $s\in S$, every source, target, identity and composition of sort $s$ whose lower dimension $j$ satisfies $j\geq m$ is sent to a variable of sort $(m,s)$.

\begin{definition}\label{DDsFunctor}\lean{equivalent}{Dimensions.homogeneousDiscretizationPartialFunctor,Dimensions.heterogeneousDiscretizationPartialFunctor}
The \emph{homogeneous discretization functor} $\Ds^{(n,m)}_{\sig{\Sigma}}$ is the functor $(\bDs^{(n,m)}_{\sig{\Sigma}})^{\ast}_{\gd}$ from $\mathsf{PAlg}(\Sigma^{(m)})$ to $\mathsf{PAlg}(\Sigma^{(n)})$, and the \emph{heterogeneous discretization functor} $\Ds^{\htg(n,m)}_{\sig{\Sigma}}$ is the functor $(\bDs^{\htg(n,m)}_{\sig{\Sigma}})^{\ast}_{\gd}$ from $\mathsf{PAlg}(\Sigma^{\htg(m)})$ to $\mathsf{PAlg}(\Sigma^{\htg(n)})$, see Definition~\ref{DDsDerivor} and Proposition~\ref{PGenHepDerivorFunctor}.
\end{definition}

In the following lemma we describe the transformed algebras along the discretization derivors.

\begin{lemma}\label{LDsDescription}\lean{weaker}{Dimensions.homogeneousDiscretizationRealization,Dimensions.heterogeneousDimensionRealization}
Assume that $m<n$.
Then the following statements hold:
\begin{enumerate}
\item for every object $\mathbf{B}$ in $\mathsf{EPAlg}(\spec{C}^{\sig{\Sigma}^{(m)}})$, $\Ds^{(n,m)}_{\sig{\Sigma}}(\mathbf{B})$ is a $\sig{\Sigma}^{(n)}$-algebra whose cells of sort $s$, for every $s\in S$, are all the elements of $B_{s}$, in which, for every $(\mathbf{s},s)\in S^{\star}\times S$ and every $\sigma\in\Sigma_{\mathbf{s},s}$, $\sigma$ acts as $\sigma^{\mathbf{B}}$, and, for every $k\in n$ and every $s\in S$,
\begin{enumerate}
\item if $k<m$, $\src^{k}_{s}$, $\tgt^{k}_{s}$ and $\circ^{k}_{s}$ act as $\src^{k\,\mathbf{B}}_{s}$, $\tgt^{k\,\mathbf{B}}_{s}$ and $\circ^{k\,\mathbf{B}}_{s}$, respectively; and
\item if $k\geq m$, $\src^{k}_{s}$ and $\tgt^{k}_{s}$ act as the identity of $B_{s}$ and, for every $x,y\in B_{s}$, $x\circ^{k}_{s}y$ is defined if, and only if, $x=y$, in which case it is $x$; and
\end{enumerate}
\item for every object $\mathbf{B}$ in $\mathsf{EPAlg}(\spec{C}^{\sig{\Sigma}^{\htg(m)}})$, $\Ds^{\htg(n,m)}_{\sig{\Sigma}}(\mathbf{B})$ is a $\sig{\Sigma}^{\htg(n)}$-algebra whose cells of sort $(k,s)$, for every $k\in n+1$ and every $s\in S$, are all the elements of $B_{(\min\{k,m\},s)}$, in which, for every $(\mathbf{s},s)\in S^{\star}\times S$, every $\sigma\in\Sigma_{\mathbf{s},s}$ and every $k\in n+1$, $\sigma_{k}$ acts as $\sigma^{\mathbf{B}}_{\min\{k,m\}}$, and, for every $j<k$ in $n+1$ and every $s\in S$,
\begin{enumerate}
\item if $j<m$, $\src^{(j,k)}_{s}$, $\tgt^{(j,k)}_{s}$, $\mathrm{i}^{(k,j)}_{s}$ and $\circ^{(j,k)}_{s}$ act as $\src^{(j,\min\{k,m\})\,\mathbf{B}}_{s}$, $\tgt^{(j,\min\{k,m\})\,\mathbf{B}}_{s}$, $\mathrm{i}^{(\min\{k,m\},j)\,\mathbf{B}}_{s}$ and $\circ^{(j,\min\{k,m\})\,\mathbf{B}}_{s}$, respectively; and
\item if $j\geq m$, $\src^{(j,k)}_{s}$, $\tgt^{(j,k)}_{s}$ and $\mathrm{i}^{(k,j)}_{s}$ act as the identity of $B_{(m,s)}$ and, for every $x,y\in B_{(m,s)}$, $x\circ^{(j,k)}_{s}y$ is defined if, and only if, $x=y$, in which case it is $x$.
\end{enumerate}
\end{enumerate}
\end{lemma}
\begin{proof}

(1) Let $\mathbf{B}$ be an object in $\mathsf{EPAlg}(\spec{C}^{\sig{\Sigma}^{(m)}})$.
Since the domains of $\bDs^{(n,m)}_{\sig{\Sigma}}$ are empty and its mapping on sorts is $\id_{S}$, for every $s\in S$, the cells of sort $s$ of $\mathbf{B}$ relative to $\bDs^{(n,m)}_{\sig{\Sigma}}$ are all the elements of $B_{s}$, see Definition~\ref{DDomainSet}.
In particular, every value of a partial term operation of $\mathbf{B}$ is a cell.
Hence, by Proposition~\ref{PTransformedEpAlg}, to prove that $\Ds^{(n,m)}_{\sig{\Sigma}}(\mathbf{B})$ is a $\sig{\Sigma}^{(n)}$-algebra, we have to show that, for every operation symbol $\tau$ of $\Sigma^{(n)}$ and every family $\mathbf{x}$ of elements of $\mathbf{B}$ that satisfies the domain conditions of $\tau$ in $\Ds^{(n,m)}_{\sig{\Sigma}}(\mathbf{B})$, the partial term operation $\Ds^{(n,m)}_{\sig{\Sigma}}(\tau)^{\mathbf{B}}$ is defined at $\mathbf{x}$.

For every $(\mathbf{s},s)\in S^{\star}\times S$ and every $\sigma\in\Sigma_{\mathbf{s},s}$, the set of domain conditions of $\sigma$ is empty and the partial term operation $\sigma((v^{\mathbf{s}}_{p})_{p\in\bb{\mathbf{s}}})^{\mathbf{B}}$, which is $\sigma^{\mathbf{B}}$, is total, by~(A0).

Let $k\in n$ and $s\in S$.
The sets of domain conditions of $\src^{k}_{s}$ and $\tgt^{k}_{s}$ are empty.
If $k<m$, the partial term operations $\src^{k}_{s}(v^{s}_{0})^{\mathbf{B}}$ and $\tgt^{k}_{s}(v^{s}_{0})^{\mathbf{B}}$, which are $\src^{k\,\mathbf{B}}_{s}$ and $\tgt^{k\,\mathbf{B}}_{s}$, are total, by~(A$^{k}$1).
If $k\geq m$, the partial term operation $(v^{s}_{0})^{\mathbf{B}}$ is the identity of $B_{s}$, which is total.
Hence, in $\Ds^{(n,m)}_{\sig{\Sigma}}(\mathbf{B})$, $\src^{k}_{s}$ and $\tgt^{k}_{s}$ act as $\src^{k\,\mathbf{B}}_{s}$ and $\tgt^{k\,\mathbf{B}}_{s}$, if $k<m$, and as the identity of $B_{s}$, if $k\geq m$.

The domain condition of $\circ^{k}_{s}$ is $\src^{k}_{s}(v^{s}_{0})\eeq\tgt^{k}_{s}(v^{s}_{1})$.
If $k<m$, by the above, a pair $(x,y)$ of elements of $B_{s}$ satisfies it in $\Ds^{(n,m)}_{\sig{\Sigma}}(\mathbf{B})$ if, and only if, $\src^{k\,\mathbf{B}}_{s}(x)=\tgt^{k\,\mathbf{B}}_{s}(y)$, and, in this case, by~(A$^{k}$3), $(v^{s}_{0}\circ^{k}_{s}v^{s}_{1})^{\mathbf{B}}$ is defined at $(x,y)$, with value $x\circ^{k\,\mathbf{B}}_{s}y$.
If $k\geq m$, by the above, $(x,y)$ satisfies it if, and only if, $x=y$, and the partial term operation $(v^{s}_{0})^{\mathbf{B}}$, as a term operation in the variables $v^{s}_{0}$ and $v^{s}_{1}$, is total and sends $(x,x)$ to $x$.

All in all, by Proposition~\ref{PTransformedEpAlg}, $\Ds^{(n,m)}_{\sig{\Sigma}}(\mathbf{B})$ is a $\sig{\Sigma}^{(n)}$-algebra and, by Proposition~\ref{PTransformedAlg}, its operations are as stated.

(2) Let $\mathbf{B}$ be an object in $\mathsf{EPAlg}(\spec{C}^{\sig{\Sigma}^{\htg(m)}})$.
Since the domains of $\bDs^{\htg(n,m)}_{\sig{\Sigma}}$ are empty, for every $k\in n+1$ and every $s\in S$, the cells of sort $(k,s)$ of $\mathbf{B}$ relative to $\bDs^{\htg(n,m)}_{\sig{\Sigma}}$ are all the elements of $B_{(\min\{k,m\},s)}$, see Definition~\ref{DDomainSet}.
As in~(1), by Proposition~\ref{PTransformedEpAlg}, we have to show that the interpreting partial term operations are defined at the families of cells that satisfy the domain conditions in $\Ds^{\htg(n,m)}_{\sig{\Sigma}}(\mathbf{B})$.

For every $(\mathbf{s},s)\in S^{\star}\times S$, every $\sigma\in\Sigma_{\mathbf{s},s}$ and every $k\in n+1$, the set of domain conditions of $\sigma_{k}$ is empty and the partial term operation $\sigma_{\min\{k,m\}}((v^{(\min\{k,m\},\mathbf{s})}_{p})_{p\in\bb{\mathbf{s}}})^{\mathbf{B}}$, which is $\sigma^{\mathbf{B}}_{\min\{k,m\}}$, is total, by~(EA0).

Let $j<k$ in $n+1$ and $s\in S$.
The sets of domain conditions of $\src^{(j,k)}_{s}$, $\tgt^{(j,k)}_{s}$ and $\mathrm{i}^{(k,j)}_{s}$ are empty.
If $j<m$, then $j<\min\{k,m\}\leq m$, and the partial term operations $\src^{(j,\min\{k,m\})}_{s}(v^{(\min\{k,m\},s)}_{0})^{\mathbf{B}}$, $\tgt^{(j,\min\{k,m\})}_{s}(v^{(\min\{k,m\},s)}_{0})^{\mathbf{B}}$ and $\mathrm{i}^{(\min\{k,m\},j)}_{s}(v^{(j,s)}_{0})^{\mathbf{B}}$, which are $\src^{(j,\min\{k,m\})\,\mathbf{B}}_{s}$, $\tgt^{(j,\min\{k,m\})\,\mathbf{B}}_{s}$ and $\mathrm{i}^{(\min\{k,m\},j)\,\mathbf{B}}_{s}$, are total, by~(EA$^{j,\min\{k,m\}}$1).
If $j\geq m$, then $\min\{j,m\}=\min\{k,m\}=m$, and the partial term operation $(v^{(m,s)}_{0})^{\mathbf{B}}$ is the identity of $B_{(m,s)}$, which is total.

The domain condition of $\circ^{(j,k)}_{s}$ is $\src^{(j,k)}_{s}(v^{(k,s)}_{0})\eeq\tgt^{(j,k)}_{s}(v^{(k,s)}_{1})$.
If $j<m$, by the above, a pair $(x,y)$ of elements of $B_{(\min\{k,m\},s)}$ satisfies it in $\Ds^{\htg(n,m)}_{\sig{\Sigma}}(\mathbf{B})$ if, and only if, $\src^{(j,\min\{k,m\})\,\mathbf{B}}_{s}(x)=\tgt^{(j,\min\{k,m\})\,\mathbf{B}}_{s}(y)$, and, in this case, by~(EA$^{j,\min\{k,m\}}$3), $x\circ^{(j,\min\{k,m\})\,\mathbf{B}}_{s}y$ is defined.
If $j\geq m$, by the above, $(x,y)$ satisfies it if, and only if, $x=y$, and the partial term operation $(v^{(m,s)}_{0})^{\mathbf{B}}$, as a term operation in the variables $v^{(m,s)}_{0}$ and $v^{(m,s)}_{1}$, is total and sends $(x,x)$ to $x$.

All in all, by Proposition~\ref{PTransformedEpAlg}, $\Ds^{\htg(n,m)}_{\sig{\Sigma}}(\mathbf{B})$ is a $\sig{\Sigma}^{\htg(n)}$-algebra and, by Proposition~\ref{PTransformedAlg}, its operations are as stated. This completes the proof.
\end{proof}

We next prove that discretization sends categorial algebras to categorial algebras.

\begin{proposition}\label{PDsEpVariety}\lean{equivalent}{Dimensions.homogeneousDiscretization_isModel,Dimensions.heterogeneousDiscretization_isModel}
The following statements hold:
\begin{enumerate}
\item the restriction of $\Ds^{(n,m)}_{\sig{\Sigma}}$ to $\mathsf{EPAlg}(\spec{C}^{\sig{\Sigma}^{(m)}})$ corestricts to $\mathsf{EPAlg}(\spec{C}^{\sig{\Sigma}^{(n)}})$; and
\item the restriction of $\Ds^{\htg(n,m)}_{\sig{\Sigma}}$\\ to $\mathsf{EPAlg}(\spec{C}^{\sig{\Sigma}^{\htg(m)}})$ corestricts to $\mathsf{EPAlg}(\spec{C}^{\sig{\Sigma}^{\htg(n)}})$.
\end{enumerate}
\end{proposition}
\begin{proof}
If $m=n$, both statements follow from Corollary~\ref{CGenHepDerivorIdFunctor}.
Assume that $m<n$.

(1) Let $\mathbf{B}$ be an object in $\mathsf{EPAlg}(\spec{C}^{\sig{\Sigma}^{(m)}})$.
By Lemma~\ref{LDsDescription}~(1), $\Ds^{(n,m)}_{\sig{\Sigma}}(\mathbf{B})$ is a $\sig{\Sigma}^{(n)}$-algebra.
It remains to check the axioms of $\spec{C}^{\sig{\Sigma}^{(n)}}$.
We do it according to the dimensions involved, by means of the description of the operations of $\Ds^{(n,m)}_{\sig{\Sigma}}(\mathbf{B})$ in Lemma~\ref{LDsDescription}~(1).

The axiom (A0) and the axioms in which all the dimensions are smaller than $m$ are axioms of $\spec{C}^{\sig{\Sigma}^{(m)}}$, which hold in $\mathbf{B}$.

The axioms (A$^{k}$1)--(A$^{k}$8), for $m\leq k<n$, become trivial, since, in $\Ds^{(n,m)}_{\sig{\Sigma}}(\mathbf{B})$, $\src^{k}_{s}$ and $\tgt^{k}_{s}$ are identities and $x\circ^{k}_{s}y$ is only defined for $x=y$, with value $x$.
For instance, (A$^{k}$5) becomes $x=x$.

The axioms (AB$^{l,k}$1)--(AB$^{l,k}$3), for $l<k$ in $n$ with $k\geq m$, also become trivial.
For instance, if $x\circ^{l}_{s}y$ is defined, then both sides of (AB$^{l,k}$2) are $x\circ^{l}_{s}y$, since $\src^{k}_{s}$ is the identity.
And the premises of (AB$^{l,k}$3) force $x=z$ and $y=t$.

(2) Let $\mathbf{B}$ be an object in $\mathsf{EPAlg}(\spec{C}^{\sig{\Sigma}^{\htg(m)}})$.
By Lemma~\ref{LDsDescription}~(2), $\Ds^{\htg(n,m)}_{\sig{\Sigma}}(\mathbf{B})$ is a $\sig{\Sigma}^{\htg(n)}$-algebra, and we check the axioms of $\spec{C}^{\sig{\Sigma}^{\htg(n)}}$ as in~(1), by means of the description of its operations in Lemma~\ref{LDsDescription}~(2).
The axioms in which all the dimensions are at most $m$ are axioms of $\spec{C}^{\sig{\Sigma}^{\htg(m)}}$, which hold in $\mathbf{B}$.
For $m<k$, the axiom (EA0) for $\sigma_{k}$ becomes (EA0) for $\sigma_{m}$.
For $j<m<k$, the axioms (EA$^{j,k}$1)--(EA$^{j,k}$8) become (EA$^{j,m}$1)--(EA$^{j,m}$8), and, for $h<j<m<k$, the axioms (EA$^{h,j,k}$1)--(EA$^{h,j,k}$5) become (EA$^{h,j,m}$1)--(EA$^{h,j,m}$5).
The remaining axioms are (EA$^{j,k}$1)--(EA$^{j,k}$8) and (EA$^{h,j,k}$1)--(EA$^{h,j,k}$5) with $j\geq m$.
In them the $(j,k)$-sources, $(j,k)$-targets and $(k,j)$-identities are identities and $x\circ^{(j,k)}_{s}y$ is only defined for $x=y$, with value $x$, so that they become trivial, as in~(1).
For instance, if $h<m\leq j$, then both sides of (EA$^{h,j,k}$4) are $x\circ^{(h,m)}_{s}z$. This completes the proof.
\end{proof}

From Proposition~\ref{PDsEpVariety} we obtain the discretization functors between the categories of higher-order categorial algebras.

\begin{corollary}\label{CDsFunctor}\lean{equivalent}{Dimensions.homogeneousDiscretizationFunctor,Dimensions.heterogeneousDiscretizationFunctor}
The restriction of $\Ds^{(n,m)}_{\sig{\Sigma}}$ to $\mathsf{mCatAlg}(\sig{\Sigma})$ corestricts to $\mathsf{nCatAlg}(\sig{\Sigma})$, and the restriction of $\Ds^{\htg(n,m)}_{\sig{\Sigma}}$ to $\mathsf{mCat}^{\htg}\mathsf{Alg}(\sig{\Sigma})$ corestricts to $\mathsf{nCat}^{\htg}\mathsf{Alg}(\sig{\Sigma})$.
\end{corollary}
\begin{proof}
It follows from Proposition~\ref{PDsEpVariety} and Remarks~\ref{RCatAlgAsEpAlg} and~\ref{RCatAlgAsEpAlgHtg}.
\end{proof}

By abuse of notation, we denote the corestrictions of Corollary~\ref{CDsFunctor} by
\begin{enumerate}
	\item
	$\Ds^{(n,m)}_{\sig{\Sigma}}\colon\mathsf{mCatAlg}(\sig{\Sigma})\mor\mathsf{nCatAlg}(\sig{\Sigma})$; and

	\item
	$\Ds^{\htg(n,m)}_{\sig{\Sigma}}\colon\mathsf{mCat}^{\htg}\mathsf{Alg}(\sig{\Sigma})\mor\mathsf{nCat}^{\htg}\mathsf{Alg}(\sig{\Sigma})$.
\end{enumerate}
By Lemma~\ref{LDsDescription}, $\Ds^{(n,m)}_{\sig{\Sigma}}$ keeps the cells and the structure in the dimensions $k<m$, and makes all the cells $k$-cells, for $k\geq m$.
Similarly, $\Ds^{\htg(n,m)}_{\sig{\Sigma}}$ takes as $k$-cells, for $k>m$, the $m$-cells.
Both functors send a homomorphism to itself, in the heterogeneous case repeating its component at the $m$-cells in the dimensions $k>m$.
By Corollary~\ref{CGenHepDerivorIdFunctor}, $\Ds^{(n,n)}_{\sig{\Sigma}}$ and $\Ds^{\htg(n,n)}_{\sig{\Sigma}}$ are identities.

We end this section by proving that the discretization derivors and functors compose as expected.

\begin{proposition}\label{PDsDerivorComp}\lean{weaker}{Dimensions.homogeneousDiscretization_comp,Dimensions.heterogeneousDiscretization_comp,Dimensions.homogeneousDiscretizationFunctorComp,Dimensions.heterogeneousDiscretizationFunctorComp}
Let $p\in\mathbb{N}$ with $n\leq p$.
Then the following equalities hold in $\mathsf{HSig}_{\gd}$
\begin{align*}
\bDs^{(n,m)}_{\sig{\Sigma}}\circ\bDs^{(p,n)}_{\sig{\Sigma}}&=\bDs^{(p,m)}_{\sig{\Sigma}},\\
\bDs^{\htg(n,m)}_{\sig{\Sigma}}\circ\bDs^{\htg(p,n)}_{\sig{\Sigma}}&=\bDs^{\htg(p,m)}_{\sig{\Sigma}},
\end{align*}
and the following equalities of functors hold
\begin{enumerate}
	\item
	$\Ds^{(p,n)}_{\sig{\Sigma}}\circ\Ds^{(n,m)}_{\sig{\Sigma}}=\Ds^{(p,m)}_{\sig{\Sigma}}$; and

	\item
	$\Ds^{\htg(p,n)}_{\sig{\Sigma}}\circ\Ds^{\htg(n,m)}_{\sig{\Sigma}}=\Ds^{\htg(p,m)}_{\sig{\Sigma}}$.
\end{enumerate}
\end{proposition}
\begin{proof}
If $m=n$ or $n=p$, one of the generalized hep-derivors and one of the functors involved in each equality is an identity, by Definition~\ref{DDsDerivor} and Corollary~\ref{CGenHepDerivorIdFunctor}, and the statement follows from Proposition~\ref{PHepSigCat}.
Assume that $m<n<p$.
We first prove the equalities in $\mathsf{HSig}_{\gd}$.
All the domains of the discretization derivors are empty, so that, by Definition~\ref{DGenHepDerivorCompId}, the domains of the composites are empty.
It remains to compare the underlying derivors, which we do by means of~(Der1)--(Der2) and Definition~\ref{DDsDerivor}.

In the homogeneous case, the three mappings on sorts are $\id_{S}$.
For every $(\mathbf{s},s)\in S^{\star}\times S$ and every $\sigma\in\Sigma_{\mathbf{s},s}$, both $\bDs^{(p,n)}_{\sig{\Sigma}}$ and $\bDs^{(n,m)}_{\sig{\Sigma}}$ send $\sigma$ to $\sigma((v^{\mathbf{s}}_{p'})_{p'\in\bb{\mathbf{s}}})$, as $\bDs^{(p,m)}_{\sig{\Sigma}}$ does.
Let $k\in p$ and $s\in S$.
If $k<m$, then $k<n$, and both derivors send $\src^{k}_{s}$, $\tgt^{k}_{s}$ and $\circ^{k}_{s}$ to $\src^{k}_{s}(v^{s}_{0})$, $\tgt^{k}_{s}(v^{s}_{0})$ and $v^{s}_{0}\circ^{k}_{s}v^{s}_{1}$, respectively, so that the composite does the same.
If $m\leq k<n$, $\bDs^{(p,n)}_{\sig{\Sigma}}$ sends them to the same terms, and $(\Ds^{(n,m)}_{\sig{\Sigma}})^{\sharp}$ sends these terms to $v^{s}_{0}$.
If $n\leq k$, $\bDs^{(p,n)}_{\sig{\Sigma}}$ sends them to $v^{s}_{0}$, which $(\Ds^{(n,m)}_{\sig{\Sigma}})^{\sharp}$ sends to $v^{s}_{0}$.
In all the cases, the composite sends them to the terms assigned by $\bDs^{(p,m)}_{\sig{\Sigma}}$.

In the heterogeneous case, for every $k\in p+1$ and every $s\in S$, the composite mapping on sorts sends $(k,s)$ to $(\min\{\min\{k,n\},m\},s)=(\min\{k,m\},s)$.
For every $(\mathbf{s},s)\in S^{\star}\times S$, every $\sigma\in\Sigma_{\mathbf{s},s}$ and every $k\in p+1$, the composite sends $\sigma_{k}$ to $\sigma_{\min\{\min\{k,n\},m\}}=\sigma_{\min\{k,m\}}$ applied to variables.
Let $j<k$ in $p+1$ and $s\in S$.
If $j<m$, then $j<n$, and $\bDs^{\htg(p,n)}_{\sig{\Sigma}}$ sends $\src^{(j,k)}_{s}$, $\tgt^{(j,k)}_{s}$, $\mathrm{i}^{(k,j)}_{s}$ and $\circ^{(j,k)}_{s}$ to the corresponding operation symbols with $\min\{k,n\}$ in place of $k$, applied to variables, which $(\Ds^{\htg(n,m)}_{\sig{\Sigma}})^{\sharp}$ sends to the corresponding operation symbols with $\min\{\min\{k,n\},m\}=\min\{k,m\}$ in place of $k$, applied to variables.
If $m\leq j<n$, $\bDs^{\htg(p,n)}_{\sig{\Sigma}}$ sends them to the corresponding operation symbols with $\min\{k,n\}$ in place of $k$, applied to variables, whose lower dimension $j$ satisfies $j\geq m$, so that $(\Ds^{\htg(n,m)}_{\sig{\Sigma}})^{\sharp}$ sends them to $v^{(m,s)}_{0}$.
If $n\leq j$, $\bDs^{\htg(p,n)}_{\sig{\Sigma}}$ sends them to $v^{(n,s)}_{0}$, which $(\Ds^{\htg(n,m)}_{\sig{\Sigma}})^{\sharp}$ sends to $v^{(m,s)}_{0}$, since $\min\{n,m\}=m$.
In all the cases, the composite sends them to the terms assigned by $\bDs^{\htg(p,m)}_{\sig{\Sigma}}$.

Since the domains of the discretization derivors are empty, for every $m$-categorial $\sig{\Sigma}$-algebra $\calg{A}$, $\Ds^{(p,n)}_{\sig{\Sigma}}(\Ds^{(n,m)}_{\sig{\Sigma}}(\calg{A}))$ and $(\bDs^{(n,m)}_{\sig{\Sigma}}\circ\bDs^{(p,n)}_{\sig{\Sigma}})^{\ast}_{\gd}(\calg{A})$ have the same carriers as $\calg{A}$, and, by Proposition~\ref{PDsEpVariety}, the first one is a $\sig{\Sigma}^{(p)}$-algebra; and similarly in the heterogeneous case.
Hence, by Corollary~\ref{CGenHepDerivorCompEqual}, the functors induced by the above composites are, on categorial algebras, the composites of the induced functors, and the stated equalities of functors follow from the above equalities in $\mathsf{HSig}_{\gd}$. This completes the proof.
\end{proof}

\section{Truncation of higher-order categorial algebras}\label{STr}

We next define the truncation derivors, in the homogeneous and in the heterogeneous presentations.

\begin{definition}\label{DTrDerivor}\lean{equivalent}{Dimensions.homogeneousTruncation,Dimensions.heterogeneousTruncation}
If $m<n$, we define the following generalized hep-derivors.
\begin{enumerate}
\item The \emph{homogeneous truncation derivor} is the generalized hep-derivor
$$
\bTr^{(m,n)}_{\sig{\Sigma}}=\left(\varphi^{\Tr^{(m,n)}_{\sig{\Sigma}}},\Tr^{(m,n)}_{\sig{\Sigma}},\Dom^{\Tr^{(m,n)}_{\sig{\Sigma}}}\right)
$$
from $\sig{\Sigma}^{(m)}$ to $\sig{\Sigma}^{(n)}$ defined as follows:
\begin{enumerate}
\item the mapping on sorts is $\varphi^{\Tr^{(m,n)}_{\sig{\Sigma}}}=\id_{S}$;
\item the assignment of terms is the inclusion, i.e., for every operation symbol $\tau$ of $\Sigma^{(m)}$ of arity $\mathbf{s}$,
$$
\Tr^{(m,n)}_{\sig{\Sigma}}(\tau)=\tau((v^{\mathbf{s}}_{p})_{p\in\bb{\mathbf{s}}});
$$
and
\item for every $s\in S$, the domain is
$$
\Dom^{\Tr^{(m,n)}_{\sig{\Sigma}}}_{s}=\left\lbrace\left(v^{s}_{0},\src^{m}_{s}(v^{s}_{0})\right)\right\rbrace.
$$
\end{enumerate}
\item The \emph{heterogeneous truncation derivor} is the generalized hep-derivor
$$
\bTr^{\htg(m,n)}_{\sig{\Sigma}}=\left(\varphi^{\Tr^{\htg(m,n)}_{\sig{\Sigma}}},\Tr^{\htg(m,n)}_{\sig{\Sigma}},\Dom^{\Tr^{\htg(m,n)}_{\sig{\Sigma}}}\right)
$$
from $\sig{\Sigma}^{\htg(m)}$ to $\sig{\Sigma}^{\htg(n)}$ determined by the canonical embedding of $\Sigma^{\htg(m)}$ into $\Sigma^{\htg(n)}$, see the note after Definition~\ref{DCatSigHtg}, i.e., defined as follows:
\begin{enumerate}
\item for every $(k,s)\in(m+1)\times S$, the mapping on sorts sends $(k,s)$ to $\varphi^{\Tr^{\htg(m,n)}_{\sig{\Sigma}}}(k,s)=(k,s)$;
\item for every operation symbol $\tau$ of $\Sigma^{\htg(m)}$ of arity $\mathbf{s}$,
$$
\Tr^{\htg(m,n)}_{\sig{\Sigma}}(\tau)=\tau((v^{\mathbf{s}}_{p})_{p\in\bb{\mathbf{s}}});
$$
and
\item for every $(k,s)\in(m+1)\times S$, the domain is $\Dom^{\Tr^{\htg(m,n)}_{\sig{\Sigma}}}_{(k,s)}=\varnothing$.
\end{enumerate}
\end{enumerate}
For $m=n$, we let $\bTr^{(n,n)}_{\sig{\Sigma}}$ and $\bTr^{\htg(n,n)}_{\sig{\Sigma}}$ be the identities $\id_{\sig{\Sigma}^{(n)}}$ and $\id_{\sig{\Sigma}^{\htg(n)}}$, see Definition~\ref{DGenHepDerivorCompId}.
\end{definition}

Let us note that, by Definition~\ref{DDomainSet}, for every partial $\Sigma^{(n)}$-algebra $\mathbf{B}$ and every $s\in S$, the cells of sort $s$ of $\mathbf{B}$ relative to $\bTr^{(m,n)}_{\sig{\Sigma}}$ are its $m$-cells of sort $s$, i.e., the elements $x$ of $B_{s}$ such that $x=\src^{m\,\mathbf{B}}_{s}(x)$, and that, for every partial $\Sigma^{\htg(n)}$-algebra $\mathbf{B}$, the cells of sort $(k,s)$ of $\mathbf{B}$ relative to $\bTr^{\htg(m,n)}_{\sig{\Sigma}}$ are all the elements of $B_{(k,s)}$.

\begin{definition}\label{DTrFunctor}\lean{equivalent}{Dimensions.homogeneousTruncationPartialFunctor,Dimensions.heterogeneousTruncationPartialFunctor}
The \emph{homogeneous truncation functor} $\Tr^{(m,n)}_{\sig{\Sigma}}$ is the functor $(\bTr^{(m,n)}_{\sig{\Sigma}})^{\ast}_{\gd}$ from $\mathsf{PAlg}(\Sigma^{(n)})$ to $\mathsf{PAlg}(\Sigma^{(m)})$, and the \emph{heterogeneous truncation functor} $\Tr^{\htg(m,n)}_{\sig{\Sigma}}$ is the functor $(\bTr^{\htg(m,n)}_{\sig{\Sigma}})^{\ast}_{\gd}$ from $\mathsf{PAlg}(\Sigma^{\htg(n)})$ to $\mathsf{PAlg}(\Sigma^{\htg(m)})$, see Definition~\ref{DTrDerivor} and Proposition~\ref{PGenHepDerivorFunctor}.
\end{definition}

In the following lemma we prove that the $m$-cells of a homogeneous $n$-categorial $\sig{\Sigma}$-algebra are closed under the operations of $\Sigma^{(m)}$.

\begin{lemma}\label{LTrCells}\lean{equivalent}{Dimensions.homogeneousTruncation_transformed}
Assume that $m<n$ and let $\mathbf{B}$ be an object in $\mathsf{EPAlg}(\spec{C}^{\sig{\Sigma}^{(n)}})$.
Then the following statements hold:
\begin{enumerate}
\item for every $(\mathbf{s},s)\in S^{\star}\times S$, every $\sigma\in\Sigma_{\mathbf{s},s}$ and every family $(x_{p})_{p\in\bb{\mathbf{s}}}$ of $m$-cells, with $x_{p}$ of sort $s_{p}$, $\sigma^{\mathbf{B}}((x_{p})_{p\in\bb{\mathbf{s}}})$ is an $m$-cell;
\item for every $k<m$, every $s\in S$ and every cell $x$ of sort $s$, $\src^{k\,\mathbf{B}}_{s}(x)$ and $\tgt^{k\,\mathbf{B}}_{s}(x)$ are $m$-cells; and
\item for every $k<m$ and every $s\in S$, if $x$, $y$ are $m$-cells of sort $s$ such that $\src^{k\,\mathbf{B}}_{s}(x)=\tgt^{k\,\mathbf{B}}_{s}(y)$, then $x\circ^{k\,\mathbf{B}}_{s}y$ is defined and it is an $m$-cell.
\end{enumerate}
\end{lemma}
\begin{proof}
(1) It follows from Lemma~\ref{LGdCells}~(1), with $m$ in place of $k$.

(2) By Lemma~\ref{LGdCells}~(2), $\src^{k\,\mathbf{B}}_{s}(x)$ and $\tgt^{k\,\mathbf{B}}_{s}(x)$ are $k$-cells and, by Lemma~\ref{LGdCells}~(3), they are $m$-cells.

(3) It follows from Lemma~\ref{LGdCells}~(4), with $k$ and $m$ in place of $j$ and $k$.
\end{proof}

We next prove that truncation sends categorial algebras to categorial algebras.

\begin{proposition}\label{PTrEpVariety}\lean{equivalent}{Dimensions.homogeneousTruncation_isModel,Dimensions.heterogeneousTruncation_isModel}
The following statements hold:
\begin{enumerate}
\item the restriction of $\Tr^{(m,n)}_{\sig{\Sigma}}$ to $\mathsf{EPAlg}(\spec{C}^{\sig{\Sigma}^{(n)}})$ corestricts to $\mathsf{EPAlg}(\spec{C}^{\sig{\Sigma}^{(m)}})$; and
\item the restriction of $\Tr^{\htg(m,n)}_{\sig{\Sigma}}$\\ to $\mathsf{EPAlg}(\spec{C}^{\sig{\Sigma}^{\htg(n)}})$ corestricts to $\mathsf{EPAlg}(\spec{C}^{\sig{\Sigma}^{\htg(m)}})$.
\end{enumerate}
\end{proposition}
\begin{proof}
If $m=n$, both statements follow from Corollary~\ref{CGenHepDerivorIdFunctor}.
Assume that $m<n$.

(1) Let $\mathbf{B}$ be an object in $\mathsf{EPAlg}(\spec{C}^{\sig{\Sigma}^{(n)}})$.
By the note after Definition~\ref{DTrDerivor}, the cells of $\Tr^{(m,n)}_{\sig{\Sigma}}(\mathbf{B})$ are the $m$-cells of $\mathbf{B}$, and every operation symbol $\tau$ of $\Sigma^{(m)}$ is interpreted by means of $\tau^{\mathbf{B}}$, since $\Tr^{(m,n)}_{\sig{\Sigma}}(\tau)$ is $\tau$ applied to variables.
By Lemma~\ref{LTrCells}, the $m$-cells are closed under the operations of $\mathbf{B}$ of the symbols of $\Sigma^{(m)}$.
Moreover, by~(A0), (A$^{k}$1) and~(A$^{k}$3), for $k<m$, these operations of $\mathbf{B}$ are defined on the families of $m$-cells that satisfy their domain conditions.
Hence, by Proposition~\ref{PTransformedEpAlg}, $\Tr^{(m,n)}_{\sig{\Sigma}}(\mathbf{B})$ is a $\sig{\Sigma}^{(m)}$-algebra.
In fact, it is the relative subalgebra of the reduct of $\mathbf{B}$ to $\Sigma^{(m)}$ on the $m$-cells, and this relative subalgebra is closed, i.e., the domain of every operation is the set of the families of $m$-cells in the domain of the corresponding operation of $\mathbf{B}$, whose values are $m$-cells.
Hence, by induction on terms, for every valuation in the $m$-cells, a term of $\Sigma^{(m)}$ is defined in $\Tr^{(m,n)}_{\sig{\Sigma}}(\mathbf{B})$ if, and only if, it is defined in $\mathbf{B}$, and then both values coincide.
Therefore every axiom of $\spec{C}^{\sig{\Sigma}^{(m)}}$, which is an axiom of $\spec{C}^{\sig{\Sigma}^{(n)}}$, holds in it.

(2) Let $\mathbf{B}$ be an object in $\mathsf{EPAlg}(\spec{C}^{\sig{\Sigma}^{\htg(n)}})$.
By the note after Definition~\ref{DTrDerivor}, $\Tr^{\htg(m,n)}_{\sig{\Sigma}}(\mathbf{B})$ is the reduct of $\mathbf{B}$ to the sorts of dimension at most $m$ and to the operation symbols of $\Sigma^{\htg(m)}$.
These operation symbols have the same domain conditions in $\sig{\Sigma}^{\htg(m)}$ and in $\sig{\Sigma}^{\htg(n)}$.
Hence it is a $\sig{\Sigma}^{\htg(m)}$-algebra.
Moreover, the axioms of $\spec{C}^{\sig{\Sigma}^{\htg(m)}}$ are axioms of $\spec{C}^{\sig{\Sigma}^{\htg(n)}}$, so that they hold in it. This completes the proof.
\end{proof}

From Proposition~\ref{PTrEpVariety} we obtain the truncation functors between the categories of higher-order categorial algebras.

\begin{corollary}\label{CTrFunctor}\lean{equivalent}{Dimensions.homogeneousTruncationFunctor,Dimensions.heterogeneousTruncationFunctor}
The restriction of $\Tr^{(m,n)}_{\sig{\Sigma}}$ to $\mathsf{nCatAlg}(\sig{\Sigma})$ corestricts to $\mathsf{mCatAlg}(\sig{\Sigma})$, and the restriction of $\Tr^{\htg(m,n)}_{\sig{\Sigma}}$ to $\mathsf{nCat}^{\htg}\mathsf{Alg}(\sig{\Sigma})$ corestricts to $\mathsf{mCat}^{\htg}\mathsf{Alg}(\sig{\Sigma})$.
\end{corollary}
\begin{proof}
It follows from Proposition~\ref{PTrEpVariety} and Remarks~\ref{RCatAlgAsEpAlg} and~\ref{RCatAlgAsEpAlgHtg}.
\end{proof}

By abuse of notation, we denote the corestrictions of Corollary~\ref{CTrFunctor} by
\begin{enumerate}
	\item
	$\Tr^{(m,n)}_{\sig{\Sigma}}\colon\mathsf{nCatAlg}(\sig{\Sigma})\mor\mathsf{mCatAlg}(\sig{\Sigma})$; and

	\item
	$\Tr^{\htg(m,n)}_{\sig{\Sigma}}\colon\mathsf{nCat}^{\htg}\mathsf{Alg}(\sig{\Sigma})\mor\mathsf{mCat}^{\htg}\mathsf{Alg}(\sig{\Sigma})$.
\end{enumerate}
The first one sends an $n$-categorial $\sig{\Sigma}$-algebra to the $m$-categorial $\sig{\Sigma}$-algebra of its $m$-cells, and a homomorphism to its restriction.
The second one forgets the cells of dimension greater than $m$, and the corresponding components of the homomorphisms.
By Corollary~\ref{CGenHepDerivorIdFunctor}, $\Tr^{(n,n)}_{\sig{\Sigma}}$ and $\Tr^{\htg(n,n)}_{\sig{\Sigma}}$ are identities.

We end this section by proving that the truncation derivors and functors compose as expected; this is the fact used in Chapter~\ref{SProjLimits} to form the projective systems of categorial algebras.

\begin{proposition}\label{PTrComp}\lean{weaker}{Dimensions.heterogeneousTruncation_comp,Dimensions.homogeneousTruncationFunctorComp,Dimensions.heterogeneousTruncationFunctorComp}
Let $p\in\mathbb{N}$ with $n\leq p$.
Then the following equality holds in $\mathsf{HSig}_{\gd}$
$$
\bTr^{\htg(n,p)}_{\sig{\Sigma}}\circ\bTr^{\htg(m,n)}_{\sig{\Sigma}}=\bTr^{\htg(m,p)}_{\sig{\Sigma}},
$$
and the following equalities of functors hold
\begin{enumerate}
	\item
	$\Tr^{(m,n)}_{\sig{\Sigma}}\circ \Tr^{(n,p)}_{\sig{\Sigma}}=\Tr^{(m,p)}_{\sig{\Sigma}}$; and

	\item
	$\Tr^{\htg(m,n)}_{\sig{\Sigma}}\circ \Tr^{\htg(n,p)}_{\sig{\Sigma}}=\Tr^{\htg(m,p)}_{\sig{\Sigma}}$.
\end{enumerate}
\end{proposition}
\begin{proof}
If $m=n$ or $n=p$, one of the generalized hep-derivors and one of the functors involved in each equality is an identity, by Definition~\ref{DTrDerivor} and Corollary~\ref{CGenHepDerivorIdFunctor}, and the statement is trivial.
Assume that $m<n<p$.

In both presentations, the mappings on sorts of the truncation derivors are inclusions and every operation symbol $\tau$ is sent to $\tau$ applied to variables, see Definition~\ref{DTrDerivor}.
Hence, by~(Der1)--(Der2), the underlying derivor of $\bTr^{(n,p)}_{\sig{\Sigma}}\circ\bTr^{(m,n)}_{\sig{\Sigma}}$ is that of $\bTr^{(m,p)}_{\sig{\Sigma}}$, and the underlying derivor of $\bTr^{\htg(n,p)}_{\sig{\Sigma}}\circ\bTr^{\htg(m,n)}_{\sig{\Sigma}}$ is that of $\bTr^{\htg(m,p)}_{\sig{\Sigma}}$.
In the heterogeneous case all the domains are empty, so that, by Definition~\ref{DGenHepDerivorCompId}, the domains of the composite are empty, which proves the stated equality in $\mathsf{HSig}_{\gd}$.

(1) By Definition~\ref{DGenHepDerivorCompId}, for every $s\in S$, the domain of sort $s$ of $\bTr^{(n,p)}_{\sig{\Sigma}}\circ\bTr^{(m,n)}_{\sig{\Sigma}}$ is
\[
\Dom^{\Tr^{(n,p)}_{\sig{\Sigma}}}_{s}\cup\left(\Tr^{(n,p)}_{\sig{\Sigma}}\right)^{\sharp}\left[\Dom^{\Tr^{(m,n)}_{\sig{\Sigma}}}_{s}\right]
=
\left\lbrace\left(v^s_0,\src^n_s(v^s_0)\right),\left(v^s_0,\src^m_s(v^s_0)\right)\right\rbrace,
\]
whereas that of $\bTr^{(m,p)}_{\sig{\Sigma}}$ is $\{(v^s_0,\src^m_s(v^s_0))\}$; thus the two generalized hep-derivors need not be equal, see Remark~\ref{RTrCompHomogeneous}.
Let $\calg{A}$ be a $p$-categorial $\sig{\Sigma}$-algebra.
By Lemma~\ref{LGdCells}~(3), applied to $\calg{A}$, every $m$-cell of $\calg{A}$ is an $n$-cell, so that both domains select precisely the $m$-cells of $\calg{A}$.
By the note after Definition~\ref{DTrDerivor}, these are also the cells of $\Tr^{(m,n)}_{\sig{\Sigma}}(\Tr^{(n,p)}_{\sig{\Sigma}}(\calg{A}))$, which, by Proposition~\ref{PTrEpVariety}, is a $\sig{\Sigma}^{(m)}$-algebra.
Therefore, by Corollary~\ref{CGenHepDerivorCompEqual} and Lemma~\ref{LGenHepDerivorSameCells}, the following equalities hold
$$
\Tr^{(m,n)}_{\sig{\Sigma}}\left(\Tr^{(n,p)}_{\sig{\Sigma}}(\calg{A})\right)
=
\left(\bTr^{(n,p)}_{\sig{\Sigma}}\circ\bTr^{(m,n)}_{\sig{\Sigma}}\right)^{\ast}_{\gd}(\calg{A})
=
\left(\bTr^{(m,p)}_{\sig{\Sigma}}\right)^{\ast}_{\gd}(\calg{A})
=
\Tr^{(m,p)}_{\sig{\Sigma}}(\calg{A}).
$$
By the same results, the equality also holds for the homomorphisms.

(2) Let $\calg{A}$ be a heterogeneous $p$-categorial $\sig{\Sigma}$-algebra.
By the note after Definition~\ref{DTrDerivor}, $\Tr^{\htg(m,n)}_{\sig{\Sigma}}(\Tr^{\htg(n,p)}_{\sig{\Sigma}}(\calg{A}))$ and $(\bTr^{\htg(m,p)}_{\sig{\Sigma}})^{\ast}_{\gd}(\calg{A})$ have the same carriers, those of $\calg{A}$ of dimension at most $m$, and, by Proposition~\ref{PTrEpVariety}, the first one is a $\sig{\Sigma}^{\htg(m)}$-algebra.
Therefore, by Corollary~\ref{CGenHepDerivorCompEqual} and the above equality in $\mathsf{HSig}_{\gd}$, the stated equality holds on objects and on homomorphisms. This completes the proof.
\end{proof}

\begin{remark}\label{RTrCompHomogeneous}\lean{weaker}{Dimensions.homogeneousTruncation_comp_derivor,Dimensions.homogeneousTruncation_comp_ne,Dimensions.homogeneousTruncation_comp_sameCells}
In Proposition~\ref{PTrComp}, the equality in $\mathsf{HSig}_{\gd}$ is stated only for the heterogeneous truncation derivors, since its homogeneous analogue fails in general.
In fact, if $m<n<p$, then $\bTr^{(n,p)}_{\sig{\Sigma}}\circ\bTr^{(m,n)}_{\sig{\Sigma}}$ and $\bTr^{(m,p)}_{\sig{\Sigma}}$ have the same underlying derivor but different domains: for every $s\in S$, the domain of sort $s$ of the first one is $\{(v^s_0,\src^n_s(v^s_0)),(v^s_0,\src^m_s(v^s_0))\}$, whereas that of the second one is $\{(v^s_0,\src^m_s(v^s_0))\}$.
On every $p$-categorial $\sig{\Sigma}$-algebra both domains select the $m$-cells, see Remark~\ref{RSameCellsPattern}.
\end{remark}

\section{On the relationship between discretization and truncation}\label{SDsTrRel}

In this section we prove that discretization is left adjoint to truncation and we compute the composite of a discretization followed by a truncation.
As in Section~\ref{SGdUfRel}, we compute the composites of these functors by composing the corresponding generalized hep-derivors, by means of Corollary~\ref{CGenHepDerivorCompEqual} and Lemma~\ref{LGenHepDerivorSameCells}.

\begin{proposition}\label{PTrDsId}\lean{weaker}{Dimensions.heterogeneousTruncation_discretization,Dimensions.homogeneousRoundTripFunctor,Dimensions.heterogeneousRoundTripFunctor}
The following equality holds in $\mathsf{HSig}_{\gd}$
$$
\bDs^{\htg(n,m)}_{\sig{\Sigma}}\circ\bTr^{\htg(m,n)}_{\sig{\Sigma}}=\id_{\sig{\Sigma}^{\htg(m)}},
$$
and the following equalities of functors hold
\begin{enumerate}
	\item
	$\Tr^{(m,n)}_{\sig{\Sigma}}\circ\Ds^{(n,m)}_{\sig{\Sigma}}=\Id_{\mathsf{mCatAlg}(\sig{\Sigma})}$; and

	\item
	$\Tr^{\htg(m,n)}_{\sig{\Sigma}}\circ\Ds^{\htg(n,m)}_{\sig{\Sigma}}=\Id_{\mathsf{mCat}^{\htg}\mathsf{Alg}(\sig{\Sigma})}$.
\end{enumerate}
\end{proposition}
\begin{proof}
If $m=n$, all the generalized hep-derivors and functors involved are identities, by Definitions~\ref{DDsDerivor} and~\ref{DTrDerivor} and Corollary~\ref{CGenHepDerivorIdFunctor}.
Assume that $m<n$.

(1) We first compare the underlying derivors.
Both mappings on sorts are $\id_{S}$.
Moreover, $\Tr^{(m,n)}_{\sig{\Sigma}}$ sends every operation symbol $\tau$ of $\Sigma^{(m)}$ to itself, applied to variables, and so does $\Ds^{(n,m)}_{\sig{\Sigma}}$.
Hence, by~(Der1) and~(Der2), $(\Ds^{(n,m)}_{\sig{\Sigma}})^{\sharp}\circ\Tr^{(m,n)}_{\sig{\Sigma}}=\eta^{\Sigma^{(m)}}$.
That is, the underlying derivor of $\bDs^{(n,m)}_{\sig{\Sigma}}\circ\bTr^{(m,n)}_{\sig{\Sigma}}$ is the identity derivor, see Definition~\ref{DDerivorCompId}.

We next compare the cells.
By Definition~\ref{DGenHepDerivorCompId}, for every $s\in S$, the following equality holds
$$
\Dom^{\bDs^{(n,m)}_{\sig{\Sigma}}\circ\bTr^{(m,n)}_{\sig{\Sigma}}}_{s}=\left(\Ds^{(n,m)}_{\sig{\Sigma}}\right)^{\sharp}\left[\left\lbrace\left(v^{s}_{0},\src^{m}_{s}(v^{s}_{0})\right)\right\rbrace\right]=\left\lbrace\left(v^{s}_{0},v^{s}_{0}\right)\right\rbrace.
$$
Thus $\bDs^{(n,m)}_{\sig{\Sigma}}\circ\bTr^{(m,n)}_{\sig{\Sigma}}$ and $\id_{\sig{\Sigma}^{(m)}}$ need not be equal, see Remark~\ref{RTrDsIdHomogeneous}.
Let $\calg{A}$ be an $m$-categorial $\sig{\Sigma}$-algebra.
Every element of $A_{s}$ satisfies this domain, as it does the empty domain of $\id_{\sig{\Sigma}^{(m)}}$.
Moreover, by Lemma~\ref{LDsDescription}~(1), $\src^{m}_{s}$ acts as the identity of $A_{s}$ in $\Ds^{(n,m)}_{\sig{\Sigma}}(\calg{A})$.
Hence every element of $A_{s}$ is a cell of $\Tr^{(m,n)}_{\sig{\Sigma}}(\Ds^{(n,m)}_{\sig{\Sigma}}(\calg{A}))$.

Now, by Propositions~\ref{PDsEpVariety} and~\ref{PTrEpVariety}, $\Tr^{(m,n)}_{\sig{\Sigma}}(\Ds^{(n,m)}_{\sig{\Sigma}}(\calg{A}))$ is a $\sig{\Sigma}^{(m)}$-algebra.
Therefore, by Corollary~\ref{CGenHepDerivorCompEqual}, Lemma~\ref{LGenHepDerivorSameCells} and Proposition~\ref{PGenHepDerivorId}, the following equalities hold
$$
\Tr^{(m,n)}_{\sig{\Sigma}}\left(\Ds^{(n,m)}_{\sig{\Sigma}}(\calg{A})\right)
=
\left(\bDs^{(n,m)}_{\sig{\Sigma}}\circ\bTr^{(m,n)}_{\sig{\Sigma}}\right)^{\ast}_{\gd}(\calg{A})
=
\left(\id_{\sig{\Sigma}^{(m)}}\right)^{\ast}_{\gd}(\calg{A})
=
\calg{A}.
$$
By the same results, the equality also holds for the homomorphisms.

(2) We first prove the stated equality in $\mathsf{HSig}_{\gd}$.
For every $(k,s)\in(m+1)\times S$, the composite mapping on sorts sends $(k,s)$ to $(\min\{k,m\},s)=(k,s)$.
Both derivors send the operation symbols of $\Sigma^{\htg(m)}$ to themselves, applied to variables, since, for every such symbol, its lower dimension $j$ satisfies $j<m$ and its upper dimension $k$ satisfies $\min\{k,m\}=k$, see Definition~\ref{DDsDerivor}~(2).
Moreover, all the domains are empty, so that the domains of the composite are empty.
Hence $\bDs^{\htg(n,m)}_{\sig{\Sigma}}\circ\bTr^{\htg(m,n)}_{\sig{\Sigma}}=\id_{\sig{\Sigma}^{\htg(m)}}$.
Let $\calg{A}$ be a heterogeneous $m$-categorial $\sig{\Sigma}$-algebra.
By Propositions~\ref{PDsEpVariety} and~\ref{PTrEpVariety}, $\Tr^{\htg(m,n)}_{\sig{\Sigma}}(\Ds^{\htg(n,m)}_{\sig{\Sigma}}(\calg{A}))$ is a $\sig{\Sigma}^{\htg(m)}$-algebra, and, by Lemma~\ref{LDsDescription}~(2) and the note after Definition~\ref{DTrDerivor}, its carriers are those of $\calg{A}$.
Hence the equality follows, as in~(1), from Corollary~\ref{CGenHepDerivorCompEqual} and Proposition~\ref{PGenHepDerivorId}. This completes the proof.
\end{proof}

\begin{remark}\label{RTrDsIdHomogeneous}\lean{weaker}{Dimensions.homogeneousTruncation_discretization_derivor,Dimensions.homogeneousTruncation_discretization_ne}
In Proposition~\ref{PTrDsId}, the equality in $\mathsf{HSig}_{\gd}$ is stated only for the heterogeneous derivors, since its homogeneous analogue fails.
In fact, if $m<n$, then $\bDs^{(n,m)}_{\sig{\Sigma}}\circ\bTr^{(m,n)}_{\sig{\Sigma}}$ and $\id_{\sig{\Sigma}^{(m)}}$ have the same underlying derivor, the identity derivor, but different domains: for every $s\in S$, the domain of sort $s$ of the first one is $\{(v^s_0,v^s_0)\}$, whereas that of the second one is empty.
Both domains are satisfied by every element of every partial $\Sigma^{(m)}$-algebra, see Remark~\ref{RSameCellsPattern}.
\end{remark}

In the following lemma we introduce the counits of the adjunctions between discretization and truncation.

\begin{lemma}\label{LDsTrCounit}\lean{equivalent}{Dimensions.homogeneousCounit,Dimensions.heterogeneousCounit}
The following statements hold:
\begin{enumerate}
\item for every $n$-categorial $\sig{\Sigma}$-algebra $\calg{B}$, the inclusion $\varepsilon_{\calg{B}}$ of the $m$-cells of $\calg{B}$ into $\calg{B}$ is a homomorphism from $\Ds^{(n,m)}_{\sig{\Sigma}}(\Tr^{(m,n)}_{\sig{\Sigma}}(\calg{B}))$ to $\calg{B}$, natural in $\calg{B}$; and
\item for every heterogeneous $n$-categorial $\sig{\Sigma}$-algebra $\calg{B}$, the family $\varepsilon^{\htg}_{\calg{B}}$ whose component at the sort $(k,s)$ is the identity of $B_{(k,s)}$, if $k\leq m$, and $\mathrm{i}^{(k,m)\,\calg{B}}_{s}$, if $k>m$, is a homomorphism from $\Ds^{\htg(n,m)}_{\sig{\Sigma}}(\Tr^{\htg(m,n)}_{\sig{\Sigma}}(\calg{B}))$ to $\calg{B}$, natural in $\calg{B}$.
\end{enumerate}
\end{lemma}
\begin{proof}
If $m=n$, then $\varepsilon_{\calg{B}}$ and $\varepsilon^{\htg}_{\calg{B}}$ are identities, and there is nothing to prove.
Assume that $m<n$.

(1) Let $\calg{B}$ be an $n$-categorial $\sig{\Sigma}$-algebra.
By Remark~\ref{RCatAlgAsEpAlg}, it suffices to prove that $\varepsilon_{\calg{B}}$ is a homomorphism of partial $\Sigma^{(n)}$-algebras.
By the proof of Proposition~\ref{PTrEpVariety}, the operations of $\Tr^{(m,n)}_{\sig{\Sigma}}(\calg{B})$ are the restrictions of those of $\calg{B}$ to the $m$-cells.
Hence, by Lemma~\ref{LDsDescription}~(1), the operations of $\Sigma$ and, for every $k<m$ and every $s\in S$, the operations $\src^{k}_{s}$, $\tgt^{k}_{s}$ and $\circ^{k}_{s}$ of $\Ds^{(n,m)}_{\sig{\Sigma}}(\Tr^{(m,n)}_{\sig{\Sigma}}(\calg{B}))$ are restrictions of those of $\calg{B}$, so that $\varepsilon_{\calg{B}}$ commutes with them.

Now let $k\in n$ with $k\geq m$, $s\in S$ and $x$ an $m$-cell of sort $s$ of $\calg{B}$.
By Lemma~\ref{LDsDescription}~(1), in $\Ds^{(n,m)}_{\sig{\Sigma}}(\Tr^{(m,n)}_{\sig{\Sigma}}(\calg{B}))$, the $k$-source and the $k$-target of $x$ are $x$, and the $k$-composition is defined only at the pairs $(x,x)$, with $x\circ^{k}_{s}x=x$.
In $\calg{B}$, $x$ is a $k$-cell, trivially if $k=m$ and by Lemma~\ref{LGdCells}~(3) if $m<k$, so that $\src^{k\,\calg{B}}_{s}(x)=x=\tgt^{k\,\calg{B}}_{s}(x)$ and, by~(A$^{k}$5), $x\circ^{k\,\calg{B}}_{s}x=x\circ^{k\,\calg{B}}_{s}\src^{k\,\calg{B}}_{s}(x)=x$.
Hence $\varepsilon_{\calg{B}}$ is a homomorphism.

Finally, let $H$ be an $n$-categorial $\sig{\Sigma}$-homomorphism from $\calg{B}$ to $\calg{B}'$.
By Proposition~\ref{PGenHepDerivorFunctor}, $\Ds^{(n,m)}_{\sig{\Sigma}}(\Tr^{(m,n)}_{\sig{\Sigma}}(H))$ is the birestriction of $H$ to the $m$-cells, so that $H\circ\varepsilon_{\calg{B}}=\varepsilon_{\calg{B}'}\circ\Ds^{(n,m)}_{\sig{\Sigma}}(\Tr^{(m,n)}_{\sig{\Sigma}}(H))$, i.e., $\varepsilon$ is natural.

(2) Let $\calg{B}$ be a heterogeneous $n$-categorial $\sig{\Sigma}$-algebra.
By Remark~\ref{RCatAlgAsEpAlgHtg}, it suffices to prove that $\varepsilon^{\htg}_{\calg{B}}$ is a homomorphism of partial $\Sigma^{\htg(n)}$-algebras.
By the note after Definition~\ref{DTrDerivor} and Lemma~\ref{LDsDescription}~(2), the cells of sort $(k,s)$ of $\Ds^{\htg(n,m)}_{\sig{\Sigma}}(\Tr^{\htg(m,n)}_{\sig{\Sigma}}(\calg{B}))$ are the elements of $B_{(\min\{k,m\},s)}$, and its operations of dimensions at most $m$ are those of $\calg{B}$.
Since the components of $\varepsilon^{\htg}_{\calg{B}}$ at the dimensions $k\leq m$ are identities, $\varepsilon^{\htg}_{\calg{B}}$ commutes with these operations.
It remains to consider the operation symbols of $\Sigma^{\htg(n)}$ whose upper dimension $k$ satisfies $k>m$.
Let $k\in n+1$ with $k>m$ and $s\in S$.

By the third equation of~(EA$^{m,k}$7), for every $(\mathbf{s},s)\in S^{\star}\times S$, every $\sigma\in\Sigma_{\mathbf{s},s}$ and every family $(u_{p})_{p\in\bb{\mathbf{s}}}$ of $m$-cells, with $u_{p}$ of sort $s_{p}$,
$$
\mathrm{i}^{(k,m)\,\calg{B}}_{s}\left(\sigma^{\calg{B}}_{m}((u_{p})_{p\in\bb{\mathbf{s}}})\right)=\sigma^{\calg{B}}_{k}\left(\left(\mathrm{i}^{(k,m)\,\calg{B}}_{s_{p}}(u_{p})\right)_{p\in\bb{\mathbf{s}}}\right),
$$
so that, by Lemma~\ref{LDsDescription}~(2), $\varepsilon^{\htg}_{\calg{B}}$ commutes with $\sigma_{k}$.

Let $j<m$.
By Lemma~\ref{LDsDescription}~(2), in $\Ds^{\htg(n,m)}_{\sig{\Sigma}}(\Tr^{\htg(m,n)}_{\sig{\Sigma}}(\calg{B}))$, $\src^{(j,k)}_{s}$, $\tgt^{(j,k)}_{s}$, $\mathrm{i}^{(k,j)}_{s}$ and $\circ^{(j,k)}_{s}$ act as $\src^{(j,m)\,\calg{B}}_{s}$, $\tgt^{(j,m)\,\calg{B}}_{s}$, $\mathrm{i}^{(m,j)\,\calg{B}}_{s}$ and $\circ^{(j,m)\,\calg{B}}_{s}$, respectively.
For every $m$-cell $u$ of sort $s$, by~(EA$^{j,m,k}$1) and~(EA$^{m,k}$2), the following equalities hold
\begin{align*}
\src^{(j,k)\,\calg{B}}_{s}\left(\mathrm{i}^{(k,m)\,\calg{B}}_{s}(u)\right)
&=\src^{(j,m)\,\calg{B}}_{s}\left(\src^{(m,k)\,\calg{B}}_{s}\left(\mathrm{i}^{(k,m)\,\calg{B}}_{s}(u)\right)\right)
=\src^{(j,m)\,\calg{B}}_{s}(u),\\
\tgt^{(j,k)\,\calg{B}}_{s}\left(\mathrm{i}^{(k,m)\,\calg{B}}_{s}(u)\right)
&=\tgt^{(j,m)\,\calg{B}}_{s}\left(\tgt^{(m,k)\,\calg{B}}_{s}\left(\mathrm{i}^{(k,m)\,\calg{B}}_{s}(u)\right)\right)
=\tgt^{(j,m)\,\calg{B}}_{s}(u).
\end{align*}
By~(EA$^{j,m,k}$2), for every $j$-cell $w$ of sort $s$, $\mathrm{i}^{(k,m)\,\calg{B}}_{s}(\mathrm{i}^{(m,j)\,\calg{B}}_{s}(w))=\mathrm{i}^{(k,j)\,\calg{B}}_{s}(w)$.
By~(EA$^{j,m,k}$5), for every pair of $m$-cells $u$, $u'$ of sort $s$ such that $u\circ^{(j,m)\,\calg{B}}_{s}u'$ is defined, $\mathrm{i}^{(k,m)\,\calg{B}}_{s}(u\circ^{(j,m)\,\calg{B}}_{s}u')=\mathrm{i}^{(k,m)\,\calg{B}}_{s}(u)\circ^{(j,k)\,\calg{B}}_{s}\mathrm{i}^{(k,m)\,\calg{B}}_{s}(u')$.
Hence $\varepsilon^{\htg}_{\calg{B}}$ commutes with $\src^{(j,k)}_{s}$, $\tgt^{(j,k)}_{s}$, $\mathrm{i}^{(k,j)}_{s}$ and $\circ^{(j,k)}_{s}$.

Let $j\geq m$, so that $m\leq j<k$.
By Lemma~\ref{LDsDescription}~(2), in $\Ds^{\htg(n,m)}_{\sig{\Sigma}}(\Tr^{\htg(m,n)}_{\sig{\Sigma}}(\calg{B}))$, $\src^{(j,k)}_{s}$, $\tgt^{(j,k)}_{s}$ and $\mathrm{i}^{(k,j)}_{s}$ act as the identity of $B_{(m,s)}$, and $\circ^{(j,k)}_{s}$ is defined only at the pairs $(u,u)$, with value $u$.
The component of $\varepsilon^{\htg}_{\calg{B}}$ at $(j,s)$ is $\mathrm{i}^{(j,m)\,\calg{B}}_{s}$, where we let $\mathrm{i}^{(m,m)\,\calg{B}}_{s}$ stand for the identity of $B_{(m,s)}$.
For every $m$-cell $u$ of sort $s$, by~(EA$^{m,j,k}$2), if $j>m$, and trivially, if $j=m$,
$$
\mathrm{i}^{(k,m)\,\calg{B}}_{s}(u)=\mathrm{i}^{(k,j)\,\calg{B}}_{s}\left(\mathrm{i}^{(j,m)\,\calg{B}}_{s}(u)\right),
$$
so that $\varepsilon^{\htg}_{\calg{B}}$ commutes with $\mathrm{i}^{(k,j)}_{s}$.
Hence, by~(EA$^{j,k}$2),
$$
\src^{(j,k)\,\calg{B}}_{s}\left(\mathrm{i}^{(k,m)\,\calg{B}}_{s}(u)\right)=\mathrm{i}^{(j,m)\,\calg{B}}_{s}(u)=\tgt^{(j,k)\,\calg{B}}_{s}\left(\mathrm{i}^{(k,m)\,\calg{B}}_{s}(u)\right),
$$
so that $\varepsilon^{\htg}_{\calg{B}}$ commutes with $\src^{(j,k)}_{s}$ and $\tgt^{(j,k)}_{s}$.
Finally, let $x=\mathrm{i}^{(k,m)\,\calg{B}}_{s}(u)$.
By the above, $\mathrm{i}^{(k,j)\,\calg{B}}_{s}(\src^{(j,k)\,\calg{B}}_{s}(x))=x$, so that, by~(EA$^{j,k}$5), $x\circ^{(j,k)\,\calg{B}}_{s}x=x\circ^{(j,k)\,\calg{B}}_{s}\mathrm{i}^{(k,j)\,\calg{B}}_{s}(\src^{(j,k)\,\calg{B}}_{s}(x))=x$.
Hence $\varepsilon^{\htg}_{\calg{B}}$ commutes with $\circ^{(j,k)}_{s}$, and it is a homomorphism.

Finally, let $H$ be a heterogeneous $n$-categorial $\sig{\Sigma}$-homomorphism from $\calg{B}$ to $\calg{B}'$.
By Proposition~\ref{PGenHepDerivorFunctor}, the component of $\Ds^{\htg(n,m)}_{\sig{\Sigma}}(\Tr^{\htg(m,n)}_{\sig{\Sigma}}(H))$ at $(k,s)$ is $H_{(\min\{k,m\},s)}$.
Since $H$ commutes with the identities, for every $k>m$ and every $s\in S$, $H_{(k,s)}\circ\mathrm{i}^{(k,m)\,\calg{B}}_{s}=\mathrm{i}^{(k,m)\,\calg{B}'}_{s}\circ H_{(m,s)}$, so that $\varepsilon^{\htg}$ is natural. This completes the proof.
\end{proof}

We are now in a position to prove that discretization is left adjoint to truncation.

\begin{proposition}\label{PDsTrAdj}\lean{weaker}{Dimensions.homogeneousAdjunction,Dimensions.heterogeneousAdjunction}
The following statements hold:
\begin{enumerate}
\item $(\Tr^{(m,n)}_{\sig{\Sigma}},\Ds^{(n,m)}_{\sig{\Sigma}},\id,\varepsilon)$ is an adjunction from $\mathsf{mCatAlg}(\sig{\Sigma})$ to $\mathsf{nCatAlg}(\sig{\Sigma})$; and
\item $(\Tr^{\htg(m,n)}_{\sig{\Sigma}},\Ds^{\htg(n,m)}_{\sig{\Sigma}},\id,\varepsilon^{\htg})$ is an adjunction from $\mathsf{mCat}^{\htg}\mathsf{Alg}(\sig{\Sigma})$ to $\mathsf{nCat}^{\htg}\mathsf{Alg}(\sig{\Sigma})$.
\end{enumerate}
Thus discretization is left adjoint to truncation, with the identity as unit and the counits of Lemma~\ref{LDsTrCounit}.
\end{proposition}
\begin{proof}
By Proposition~\ref{PTrDsId}, the identity is a natural transformation from the identity functor to $\Tr^{(m,n)}_{\sig{\Sigma}}\circ\Ds^{(n,m)}_{\sig{\Sigma}}$.
By Lemma~\ref{LDsTrCounit}, $\varepsilon$ is a natural transformation from $\Ds^{(n,m)}_{\sig{\Sigma}}\circ\Tr^{(m,n)}_{\sig{\Sigma}}$ to the identity functor.
The same holds in the heterogeneous case.
By Definition~\ref{DAdj}, it remains to prove the triangular identities.
If $m=n$, all the functors and natural transformations involved are identities.
Assume that $m<n$.

Since the unit is the identity, the triangular identities state that, for every $n$-categorial $\sig{\Sigma}$-algebra $\calg{B}$ and every $m$-categorial $\sig{\Sigma}$-algebra $\calg{A}$,
\begin{enumerate}
	\item[(TI1)]
	$\Tr^{(m,n)}_{\sig{\Sigma}}(\varepsilon_{\calg{B}})$ is the identity of $\Tr^{(m,n)}_{\sig{\Sigma}}(\calg{B})$; and

	\item[(TI2)]
	$\varepsilon_{\Ds^{(n,m)}_{\sig{\Sigma}}(\calg{A})}$ is the identity of $\Ds^{(n,m)}_{\sig{\Sigma}}(\calg{A})$,
\end{enumerate}
and similarly in the heterogeneous case.
We recall that a homomorphism is determined by its domain, its codomain and its underlying many-sorted mapping.

(TI1) By Proposition~\ref{PGenHepDerivorFunctor}, $\Tr^{(m,n)}_{\sig{\Sigma}}(\varepsilon_{\calg{B}})$ is the birestriction of $\varepsilon_{\calg{B}}$ to the $m$-cells.
Since $\varepsilon_{\calg{B}}$ is the inclusion of the $m$-cells, it is the identity.
In the heterogeneous case, $\Tr^{\htg(m,n)}_{\sig{\Sigma}}(\varepsilon^{\htg}_{\calg{B}})$ is the family of the components of $\varepsilon^{\htg}_{\calg{B}}$ at the dimensions $k\leq m$, which are identities.

(TI2) By Lemma~\ref{LDsDescription}~(1), every cell of $\Ds^{(n,m)}_{\sig{\Sigma}}(\calg{A})$ is an $m$-cell.
Hence $\varepsilon_{\Ds^{(n,m)}_{\sig{\Sigma}}(\calg{A})}$, the inclusion of the $m$-cells, is the identity.
In the heterogeneous case, the components of $\varepsilon^{\htg}_{\Ds^{\htg(n,m)}_{\sig{\Sigma}}(\calg{A})}$ at the dimensions $k>m$ are the $(k,m)$-identities of $\Ds^{\htg(n,m)}_{\sig{\Sigma}}(\calg{A})$, which, by Lemma~\ref{LDsDescription}~(2), are identities. This completes the proof.
\end{proof}

The composite of a discretization followed by a truncation is a discretization or a truncation, depending on the relative position of the dimensions.

\begin{proposition}\label{PTrDsComp}\lean{weaker}{Dimensions.heterogeneousMixed_upper,Dimensions.heterogeneousMixed_lower,Dimensions.homogeneousMixed_lower,Dimensions.homogeneousMixedUpperFunctor,Dimensions.heterogeneousMixedUpperFunctor,Dimensions.homogeneousMixedLowerFunctor,Dimensions.heterogeneousMixedLowerFunctor}
Let $p\in\mathbb{N}$ with $n\leq p$.
Then the following equalities hold in $\mathsf{HSig}_{\gd}$
\begin{align*}
\bDs^{\htg(p,m)}_{\sig{\Sigma}}\circ\bTr^{\htg(n,p)}_{\sig{\Sigma}}&=\bDs^{\htg(n,m)}_{\sig{\Sigma}},\\
\bDs^{\htg(p,n)}_{\sig{\Sigma}}\circ\bTr^{\htg(m,p)}_{\sig{\Sigma}}&=\bTr^{\htg(m,n)}_{\sig{\Sigma}},
\end{align*}
and, if $m<n$, the following equality also holds in $\mathsf{HSig}_{\gd}$
$$
\bDs^{(p,n)}_{\sig{\Sigma}}\circ\bTr^{(m,p)}_{\sig{\Sigma}}=\bTr^{(m,n)}_{\sig{\Sigma}}.
$$
Moreover, the following equalities of functors hold
\begin{enumerate}
	\item
	$\Tr^{(n,p)}_{\sig{\Sigma}}\circ\Ds^{(p,m)}_{\sig{\Sigma}}=\Ds^{(n,m)}_{\sig{\Sigma}}$ and $\Tr^{\htg(n,p)}_{\sig{\Sigma}}\circ\Ds^{\htg(p,m)}_{\sig{\Sigma}}=\Ds^{\htg(n,m)}_{\sig{\Sigma}}$; and

	\item
	$\Tr^{(m,p)}_{\sig{\Sigma}}\circ\Ds^{(p,n)}_{\sig{\Sigma}}=\Tr^{(m,n)}_{\sig{\Sigma}}$ and $\Tr^{\htg(m,p)}_{\sig{\Sigma}}\circ\Ds^{\htg(p,n)}_{\sig{\Sigma}}=\Tr^{\htg(m,n)}_{\sig{\Sigma}}$.
\end{enumerate}
\end{proposition}
\begin{proof}
If $n=p$, then $\bTr^{(n,p)}_{\sig{\Sigma}}$, $\bTr^{\htg(n,p)}_{\sig{\Sigma}}$, $\bDs^{(p,n)}_{\sig{\Sigma}}$ and $\bDs^{\htg(p,n)}_{\sig{\Sigma}}$ are identities, by Definitions~\ref{DDsDerivor} and~\ref{DTrDerivor}, and so are the functors that they induce, by Corollary~\ref{CGenHepDerivorIdFunctor}; hence the statement follows from Proposition~\ref{PHepSigCat}.
If $m=n$, both heterogeneous equalities in $\mathsf{HSig}_{\gd}$ state that $\bDs^{\htg(p,n)}_{\sig{\Sigma}}\circ\bTr^{\htg(n,p)}_{\sig{\Sigma}}=\id_{\sig{\Sigma}^{\htg(n)}}$, and both items state that $\Tr^{(n,p)}_{\sig{\Sigma}}\circ\Ds^{(p,n)}_{\sig{\Sigma}}$ and $\Tr^{\htg(n,p)}_{\sig{\Sigma}}\circ\Ds^{\htg(p,n)}_{\sig{\Sigma}}$ are identities, which is Proposition~\ref{PTrDsId}, with $n$ and $p$ in place of $m$ and $n$, respectively.
Assume that $m<n<p$.

(1) We first consider the homogeneous case.
The generalized hep-derivor $\bDs^{(p,m)}_{\sig{\Sigma}}\circ\bTr^{(n,p)}_{\sig{\Sigma}}$ goes from $\sig{\Sigma}^{(n)}$ to $\sig{\Sigma}^{(m)}$, and its mapping on sorts is $\id_{S}$.
Every operation symbol $\tau$ of $\Sigma^{(n)}$ is sent by $\bTr^{(n,p)}_{\sig{\Sigma}}$ to $\tau$ applied to variables and then, by~(Der1)--(Der2), to the term that $\bDs^{(p,m)}_{\sig{\Sigma}}$ assigns to $\tau$, which, since the dimension of $\tau$, if any, belongs to $n$, is the term that $\bDs^{(n,m)}_{\sig{\Sigma}}$ assigns to $\tau$, see Definition~\ref{DDsDerivor}.
Hence both generalized hep-derivors have the same underlying derivor.
By Definition~\ref{DGenHepDerivorCompId}, for every $s\in S$, the domain of sort $s$ of the composite is
\[
\left(\Ds^{(p,m)}_{\sig{\Sigma}}\right)^{\sharp}\left[\left\lbrace\left(v^{s}_{0},\src^{n}_{s}(v^{s}_{0})\right)\right\rbrace\right]=\left\lbrace\left(v^{s}_{0},v^{s}_{0}\right)\right\rbrace,
\]
since $n\geq m$, and every element satisfies it, as it does the empty domain of $\bDs^{(n,m)}_{\sig{\Sigma}}$; thus these two generalized hep-derivors need not be equal, see Remark~\ref{RTrDsCompHomogeneous}.
Let $\calg{A}$ be an $m$-categorial $\sig{\Sigma}$-algebra.
By Lemma~\ref{LDsDescription}~(1), $\src^{n}_{s}$ acts as the identity of $A_{s}$ in $\Ds^{(p,m)}_{\sig{\Sigma}}(\calg{A})$, so that the carriers of $\Tr^{(n,p)}_{\sig{\Sigma}}(\Ds^{(p,m)}_{\sig{\Sigma}}(\calg{A}))$ are those of $\calg{A}$.
Moreover, by Propositions~\ref{PDsEpVariety} and~\ref{PTrEpVariety}, it is a $\sig{\Sigma}^{(n)}$-algebra.
Therefore, by Corollary~\ref{CGenHepDerivorCompEqual} and Lemma~\ref{LGenHepDerivorSameCells}, $\Tr^{(n,p)}_{\sig{\Sigma}}(\Ds^{(p,m)}_{\sig{\Sigma}}(\calg{A}))=\Ds^{(n,m)}_{\sig{\Sigma}}(\calg{A})$, and, by the same results, the equality also holds for the homomorphisms.

In the heterogeneous case, the mapping on sorts of $\bDs^{\htg(p,m)}_{\sig{\Sigma}}\circ\bTr^{\htg(n,p)}_{\sig{\Sigma}}$ sends $(k,s)$ to $(\min\{k,m\},s)$, for every $k\in n+1$ and every $s\in S$, and, as above, the term assigned to every operation symbol of $\Sigma^{\htg(n)}$ is the one assigned by $\bDs^{\htg(n,m)}_{\sig{\Sigma}}$, since the clauses of Definition~\ref{DDsDerivor}~(2) depend only on the dimensions of the operation symbol and on $m$.
All the domains are empty, so that $\bDs^{\htg(p,m)}_{\sig{\Sigma}}\circ\bTr^{\htg(n,p)}_{\sig{\Sigma}}=\bDs^{\htg(n,m)}_{\sig{\Sigma}}$.
For every heterogeneous $m$-categorial $\sig{\Sigma}$-algebra $\calg{A}$, by Lemma~\ref{LDsDescription}~(2) and the note after Definition~\ref{DTrDerivor}, $\Tr^{\htg(n,p)}_{\sig{\Sigma}}(\Ds^{\htg(p,m)}_{\sig{\Sigma}}(\calg{A}))$ has as carrier of sort $(k,s)$ the carrier $A_{(\min\{k,m\},s)}$, as $\Ds^{\htg(n,m)}_{\sig{\Sigma}}(\calg{A})$ does, and, by Propositions~\ref{PDsEpVariety} and~\ref{PTrEpVariety}, it is a $\sig{\Sigma}^{\htg(n)}$-algebra.
Hence the equality follows from Corollary~\ref{CGenHepDerivorCompEqual}.

(2) We first consider the homogeneous case.
The generalized hep-derivor $\bDs^{(p,n)}_{\sig{\Sigma}}\circ\bTr^{(m,p)}_{\sig{\Sigma}}$ goes from $\sig{\Sigma}^{(m)}$ to $\sig{\Sigma}^{(n)}$, and its mapping on sorts is $\id_{S}$.
Every operation symbol $\tau$ of $\Sigma^{(m)}$ has dimension, if any, smaller than $m<n$, so that $\bDs^{(p,n)}_{\sig{\Sigma}}$ sends $\tau$ applied to variables to itself, and the composite sends $\tau$ to $\tau$ applied to variables, as $\bTr^{(m,n)}_{\sig{\Sigma}}$ does.
By Definition~\ref{DGenHepDerivorCompId}, for every $s\in S$, the domain of sort $s$ of the composite is
\[
\left(\Ds^{(p,n)}_{\sig{\Sigma}}\right)^{\sharp}\left[\left\lbrace\left(v^{s}_{0},\src^{m}_{s}(v^{s}_{0})\right)\right\rbrace\right]=\left\lbrace\left(v^{s}_{0},\src^{m}_{s}(v^{s}_{0})\right)\right\rbrace,
\]
since $m<n$, which is the domain of sort $s$ of $\bTr^{(m,n)}_{\sig{\Sigma}}$.
Hence $\bDs^{(p,n)}_{\sig{\Sigma}}\circ\bTr^{(m,p)}_{\sig{\Sigma}}=\bTr^{(m,n)}_{\sig{\Sigma}}$.
Let $\calg{A}$ be an $n$-categorial $\sig{\Sigma}$-algebra.
By Lemma~\ref{LDsDescription}~(1), $\src^{m}_{s}$ acts in $\Ds^{(p,n)}_{\sig{\Sigma}}(\calg{A})$ as $\src^{m\,\calg{A}}_{s}$, so that the cells of $\Tr^{(m,p)}_{\sig{\Sigma}}(\Ds^{(p,n)}_{\sig{\Sigma}}(\calg{A}))$ are the $m$-cells of $\calg{A}$, which are the cells of $\Tr^{(m,n)}_{\sig{\Sigma}}(\calg{A})$.
Moreover, by Propositions~\ref{PDsEpVariety} and~\ref{PTrEpVariety}, it is a $\sig{\Sigma}^{(m)}$-algebra.
Therefore, by Corollary~\ref{CGenHepDerivorCompEqual}, $\Tr^{(m,p)}_{\sig{\Sigma}}(\Ds^{(p,n)}_{\sig{\Sigma}}(\calg{A}))=\Tr^{(m,n)}_{\sig{\Sigma}}(\calg{A})$, and, by the same result, the equality also holds for the homomorphisms.

In the heterogeneous case, every operation symbol of $\Sigma^{\htg(m)}$ has dimensions $j<k\leq m<n$, so that, as above, $\bDs^{\htg(p,n)}_{\sig{\Sigma}}\circ\bTr^{\htg(m,p)}_{\sig{\Sigma}}$ sends $(k,s)$ to $(k,s)$ and every operation symbol $\tau$ to $\tau$ applied to variables, and, since all the domains are empty, it is $\bTr^{\htg(m,n)}_{\sig{\Sigma}}$.
For every heterogeneous $n$-categorial $\sig{\Sigma}$-algebra $\calg{A}$, by Lemma~\ref{LDsDescription}~(2) and the note after Definition~\ref{DTrDerivor}, $\Tr^{\htg(m,p)}_{\sig{\Sigma}}(\Ds^{\htg(p,n)}_{\sig{\Sigma}}(\calg{A}))$ has the carriers of $\calg{A}$ of dimension at most $m$, as $\Tr^{\htg(m,n)}_{\sig{\Sigma}}(\calg{A})$ does, and, by Propositions~\ref{PDsEpVariety} and~\ref{PTrEpVariety}, it is a $\sig{\Sigma}^{\htg(m)}$-algebra.
Hence the equality follows from Corollary~\ref{CGenHepDerivorCompEqual}. This completes the proof.
\end{proof}

\begin{remark}\label{RTrDsCompHomogeneous}\lean{weaker}{Dimensions.homogeneousMixed_upper_derivor,Dimensions.homogeneousMixed_upper_ne,Dimensions.homogeneousTruncation_discretization_ne}
Let $p\in\mathbb{N}$ with $n\leq p$.
In Proposition~\ref{PTrDsComp}, no equality in $\mathsf{HSig}_{\gd}$ is stated for the homogeneous derivors of item~(1), and the homogeneous equality of item~(2) is stated only for $m<n$, since otherwise they fail.
In fact, if $n<p$, then $\bDs^{(p,m)}_{\sig{\Sigma}}\circ\bTr^{(n,p)}_{\sig{\Sigma}}$ and $\bDs^{(n,m)}_{\sig{\Sigma}}$ have the same underlying derivor but different domains: for every $s\in S$, the domain of sort $s$ of the first one is $\{(v^s_0,v^s_0)\}$, whereas that of the second one is empty.
For $m=n<p$, the homogeneous equality of item~(2) would state that $\bDs^{(p,n)}_{\sig{\Sigma}}\circ\bTr^{(n,p)}_{\sig{\Sigma}}=\id_{\sig{\Sigma}^{(n)}}$, which fails for the same reason, see Remark~\ref{RTrDsIdHomogeneous}.
In both cases the domains are satisfied by every element of every partial algebra, see Remark~\ref{RSameCellsPattern}.
\end{remark}

\begin{remark}\label{RTrDsCompOverlap}\lean{weaker}{Dimensions.homogeneousTruncation_self,Dimensions.homogeneousDiscretization_self,Dimensions.heterogeneousTruncation_self,Dimensions.heterogeneousDiscretization_self,GeneralizedHepDerivor.epIdentityIso}
Let $p\in\mathbb{N}$ with $n\leq p$.
In Proposition~\ref{PTrDsComp}, item~(1) describes a discretization followed by a truncation to a dimension at least as large as the initial one, and item~(2) a discretization followed by a truncation to a dimension at most as large as the initial one.
For $m=n$, both items coincide and recover Proposition~\ref{PTrDsId}, with $n$ and $p$ in place of $m$ and $n$, respectively.
Moreover, when $n=p$, all the functors $\Tr^{(n,p)}_{\sig{\Sigma}}$, $\Tr^{\htg(n,p)}_{\sig{\Sigma}}$, $\Ds^{(p,n)}_{\sig{\Sigma}}$ and $\Ds^{\htg(p,n)}_{\sig{\Sigma}}$ are identities.
\end{remark}

\section{Compatibility with graduation and unification}\label{SDsTrCompatGdUf}

In this section we prove that the functors of this chapter commute with graduation and unification, strictly or up to a natural isomorphism.
These facts are used in Section~\ref{SAlgSquare} and in Chapter~\ref{SFunctorsOmega}.
All the functors involved act on homomorphisms by restriction to the corresponding cells, so that, in order to prove the equalities between them, it suffices to compare their values on objects.
Moreover, for $m=n$ all the statements are trivial, since all the truncation and discretization derivors and functors involved are identities.

We begin with the compatibility of truncation with graduation.

\begin{proposition}\label{PCompatGdTr}\lean{weaker}{Dimensions.truncationGraduationFunctor}
The following equality holds
$$
\Tr^{\htg(m,n)}_{\sig{\Sigma}}\circ\Gd^{(n)}_{\sig{\Sigma}}=\Gd^{(m)}_{\sig{\Sigma}}\circ \Tr^{(m,n)}_{\sig{\Sigma}}.
$$
\end{proposition}
\begin{proof}
We may assume that $m<n$.
The two sides are induced by the generalized hep-derivors $\bGd^{(n)}_{\sig{\Sigma}}\circ\bTr^{\htg(m,n)}_{\sig{\Sigma}}$ and $\bTr^{(m,n)}_{\sig{\Sigma}}\circ\bGd^{(m)}_{\sig{\Sigma}}$ from $\sig{\Sigma}^{\htg(m)}$ to $\sig{\Sigma}^{(n)}$, which have the same underlying derivor: both mappings on sorts send $(k,s)$ to $s$, and, since the truncation derivors are inclusions, both send every operation symbol of $\Sigma^{\htg(m)}$ to the term assigned to it in Definition~\ref{DGdDerivor}, which does not depend on the dimension of the graduation derivor.
By Definitions~\ref{DGenHepDerivorCompId}, \ref{DGdDerivor} and~\ref{DTrDerivor}, for every $k<m$ and every $s\in S$, the domain of sort $(k,s)$ of the first one is $\{(v^s_0,\src^k_s(v^s_0))\}$, whereas that of the second one is $\{(v^s_0,\src^m_s(v^s_0)),(v^s_0,\src^k_s(v^s_0))\}$; for $k=m$, both domains are $\{(v^s_0,\src^m_s(v^s_0))\}$.

Let $\calg{A}$ be an $n$-categorial $\sig{\Sigma}$-algebra.
Since, by Lemma~\ref{LGdCells}~(3), every $k$-cell of $\calg{A}$ is an $m$-cell, for every $k\in m+1$ and every $s\in S$, both domains of sort $(k,s)$ select the $k$-cells of sort $s$ of $\calg{A}$.
By Lemma~\ref{LGdDescription} and the note after Definition~\ref{DTrDerivor}, these are also the cells of sort $(k,s)$ of $\Tr^{\htg(m,n)}_{\sig{\Sigma}}(\Gd^{(n)}_{\sig{\Sigma}}(\calg{A}))$ and of $\Gd^{(m)}_{\sig{\Sigma}}(\Tr^{(m,n)}_{\sig{\Sigma}}(\calg{A}))$, and, by Propositions~\ref{PGdEpVariety} and~\ref{PTrEpVariety}, both are $\sig{\Sigma}^{\htg(m)}$-algebras.
Therefore, by Lemma~\ref{LGenHepDerivorCommSquare}, both are equal, and the same holds for the homomorphisms. This completes the proof.
\end{proof}

\begin{remark}\label{RCompatGdTrDerivors}\lean{weaker}{Dimensions.truncation_graduation_derivor,Dimensions.truncation_graduation_eq_iff,Dimensions.truncation_graduation_sameCells}
The equality of Proposition~\ref{PCompatGdTr} does not, in general, come from an equality in $\mathsf{HSig}_{\gd}$.
By the computation of the domains in its proof, the generalized hep-derivors $\bGd^{(n)}_{\sig{\Sigma}}\circ\bTr^{\htg(m,n)}_{\sig{\Sigma}}$ and $\bTr^{(m,n)}_{\sig{\Sigma}}\circ\bGd^{(m)}_{\sig{\Sigma}}$ are equal if $m=0$ or $m=n$, but not if $0<m<n$.
Both domains select the same cells of every $n$-categorial $\sig{\Sigma}$-algebra, see Remark~\ref{RSameCellsPattern}.
\end{remark}

We next consider the compatibility of discretization with graduation.

\begin{proposition}\label{PCompatGdDs}\lean{weaker}{Dimensions.discretizationGraduationFunctor}
The following equality holds
$$
\Gd^{(n)}_{\sig{\Sigma}}\circ\Ds^{(n,m)}_{\sig{\Sigma}}=\Ds^{\htg(n,m)}_{\sig{\Sigma}}\circ\Gd^{(m)}_{\sig{\Sigma}}.
$$
\end{proposition}
\begin{proof}
We may assume that $m<n$.
The two sides are induced by the generalized hep-derivors $\bDs^{(n,m)}_{\sig{\Sigma}}\circ\bGd^{(n)}_{\sig{\Sigma}}$ and $\bGd^{(m)}_{\sig{\Sigma}}\circ\bDs^{\htg(n,m)}_{\sig{\Sigma}}$ from $\sig{\Sigma}^{\htg(n)}$ to $\sig{\Sigma}^{(m)}$, which have the same underlying derivor: both mappings on sorts send $(k,s)$ to $s$, and, by Definitions~\ref{DGdDerivor} and~\ref{DDsDerivor}, both send $\sigma_{k}$ to $\sigma((v^{\mathbf{s}}_{p})_{p\in\bb{\mathbf{s}}})$, $\mathrm{i}^{(k,j)}_{s}$ to $v^{s}_{0}$, and $\src^{(j,k)}_{s}$, $\tgt^{(j,k)}_{s}$ and $\circ^{(j,k)}_{s}$ to $\src^{j}_{s}(v^{s}_{0})$, $\tgt^{j}_{s}(v^{s}_{0})$ and $v^{s}_{0}\circ^{j}_{s}v^{s}_{1}$, respectively, if $j<m$, and to $v^{s}_{0}$, if $j\geq m$.
By Definitions~\ref{DGenHepDerivorCompId}, \ref{DGdDerivor} and~\ref{DDsDerivor}, for every $s\in S$ and every $k\in n+1$, the domain of sort $(k,s)$ of the first one is $\{(v^s_0,\src^k_s(v^s_0))\}$, if $k<m$, $\{(v^s_0,v^s_0)\}$, if $m\leq k<n$, and empty, if $k=n$, whereas that of the second one is $\{(v^s_0,\src^k_s(v^s_0))\}$, if $k<m$, and empty, if $k\geq m$.

Let $\calg{A}$ be an $m$-categorial $\sig{\Sigma}$-algebra.
Since the condition $(v^s_0,v^s_0)$ is satisfied by every element, for every $s\in S$, both domains of sort $(k,s)$ select the $k$-cells of sort $s$ of $\calg{A}$, for $k<m$, and all the cells of sort $s$ of $\calg{A}$, for $m\leq k\leq n$.
By Lemmas~\ref{LGdDescription} and~\ref{LDsDescription}, these are also the cells of sort $(k,s)$ of $\Gd^{(n)}_{\sig{\Sigma}}(\Ds^{(n,m)}_{\sig{\Sigma}}(\calg{A}))$ and of $\Ds^{\htg(n,m)}_{\sig{\Sigma}}(\Gd^{(m)}_{\sig{\Sigma}}(\calg{A}))$, and, by Propositions~\ref{PGdEpVariety} and~\ref{PDsEpVariety}, both are $\sig{\Sigma}^{\htg(n)}$-algebras.
Therefore, by Lemma~\ref{LGenHepDerivorCommSquare}, both are equal, and the same holds for the homomorphisms. This completes the proof.
\end{proof}

\begin{remark}\label{RCompatGdDsDerivors}\lean{weaker}{Dimensions.discretization_graduation_derivor,Dimensions.discretization_graduation_eq_iff,Dimensions.discretization_graduation_sameCells}
The equality of Proposition~\ref{PCompatGdDs} does not, in general, come from an equality in $\mathsf{HSig}_{\gd}$.
By the computation of the domains in its proof, the generalized hep-derivors $\bDs^{(n,m)}_{\sig{\Sigma}}\circ\bGd^{(n)}_{\sig{\Sigma}}$ and $\bGd^{(m)}_{\sig{\Sigma}}\circ\bDs^{\htg(n,m)}_{\sig{\Sigma}}$ are equal if $m=n$, but not if $m<n$.
Both domains select the same cells of every $m$-categorial $\sig{\Sigma}$-algebra, see Remark~\ref{RSameCellsPattern}.
\end{remark}

We next consider the compatibility of discretization with unification.

\begin{proposition}\label{PCompatUfDs}\lean{weaker}{Dimensions.discretization_unification,Dimensions.discretizationUnificationFunctor}
The following equality holds in $\mathsf{HSig}_{\gd}$
$$
\bDs^{\htg(n,m)}_{\sig{\Sigma}}\circ\bUf^{(n)}_{\sig{\Sigma}}=\bUf^{(m)}_{\sig{\Sigma}}\circ\bDs^{(n,m)}_{\sig{\Sigma}},
$$
and the following equality of functors holds
$$
\Uf^{(n)}_{\sig{\Sigma}}\circ\Ds^{\htg(n,m)}_{\sig{\Sigma}}=\Ds^{(n,m)}_{\sig{\Sigma}}\circ\Uf^{(m)}_{\sig{\Sigma}}.
$$
\end{proposition}
\begin{proof}
We may assume that $m<n$.

We first prove the equality in $\mathsf{HSig}_{\gd}$.
Both generalized hep-derivors go from $\sig{\Sigma}^{(n)}$ to $\sig{\Sigma}^{\htg(m)}$.
All the domains of the unification and of the discretization derivors are empty, so that, by Definition~\ref{DGenHepDerivorCompId}, the domains of both composites are empty.
It remains to compare the underlying derivors, which we do by means of~(Der1)--(Der2) and Definitions~\ref{DUfDerivor} and~\ref{DDsDerivor}.

For every $s\in S$, the first mapping on sorts sends $s$ to $(n,s)$ and then to $(\min\{n,m\},s)=(m,s)$, and the second one sends $s$ to $s$ and then to $(m,s)$.

For every $(\mathbf{s},s)\in S^{\star}\times S$ and every $\sigma\in\Sigma_{\mathbf{s},s}$, the first composite sends $\sigma$ to $\sigma_{n}((v^{(n,\mathbf{s})}_{p})_{p\in\bb{\mathbf{s}}})$ and then to $\sigma_{\min\{n,m\}}((v^{(m,\mathbf{s})}_{p})_{p\in\bb{\mathbf{s}}})=\sigma_{m}((v^{(m,\mathbf{s})}_{p})_{p\in\bb{\mathbf{s}}})$, and the second one sends $\sigma$ to $\sigma((v^{\mathbf{s}}_{p})_{p\in\bb{\mathbf{s}}})$ and then to $\sigma_{m}((v^{(m,\mathbf{s})}_{p})_{p\in\bb{\mathbf{s}}})$.

Let $k\in n$ and $s\in S$.
The first composite sends $\src^{k}_{s}$ to $\mathrm{i}^{(n,k)}_{s}(\src^{(k,n)}_{s}(v^{(n,s)}_{0}))$ and then, if $k<m$, to $\mathrm{i}^{(m,k)}_{s}(\src^{(k,m)}_{s}(v^{(m,s)}_{0}))$, and, if $k\geq m$, to $v^{(m,s)}_{0}$, since in this case both $\mathrm{i}^{(n,k)}_{s}$ and $\src^{(k,n)}_{s}$ are sent to variables.
The second composite sends $\src^{k}_{s}$, if $k<m$, to $\src^{k}_{s}(v^{s}_{0})$ and then to $\mathrm{i}^{(m,k)}_{s}(\src^{(k,m)}_{s}(v^{(m,s)}_{0}))$, and, if $k\geq m$, to $v^{s}_{0}$ and then to $v^{(m,s)}_{0}$.
Similarly, the first composite sends $\tgt^{k}_{s}$ to $\mathrm{i}^{(n,k)}_{s}(\tgt^{(k,n)}_{s}(v^{(n,s)}_{0}))$ and then, if $k<m$, to $\mathrm{i}^{(m,k)}_{s}(\tgt^{(k,m)}_{s}(v^{(m,s)}_{0}))$, and, if $k\geq m$, to $v^{(m,s)}_{0}$; and the second one sends $\tgt^{k}_{s}$, if $k<m$, to $\tgt^{k}_{s}(v^{s}_{0})$ and then to $\mathrm{i}^{(m,k)}_{s}(\tgt^{(k,m)}_{s}(v^{(m,s)}_{0}))$, and, if $k\geq m$, to $v^{s}_{0}$ and then to $v^{(m,s)}_{0}$.
Finally, the first composite sends $\circ^{k}_{s}$ to $v^{(n,s)}_{0}\circ^{(k,n)}_{s}v^{(n,s)}_{1}$ and then, if $k<m$, to $v^{(m,s)}_{0}\circ^{(k,m)}_{s}v^{(m,s)}_{1}$, and, if $k\geq m$, to $v^{(m,s)}_{0}$; and the second one sends $\circ^{k}_{s}$, if $k<m$, to $v^{s}_{0}\circ^{k}_{s}v^{s}_{1}$ and then to $v^{(m,s)}_{0}\circ^{(k,m)}_{s}v^{(m,s)}_{1}$, and, if $k\geq m$, to $v^{s}_{0}$ and then to $v^{(m,s)}_{0}$.
Hence both underlying derivors are equal, and the stated equality holds in $\mathsf{HSig}_{\gd}$.

We next prove the equality of functors.
Let $\calg{B}$ be a heterogeneous $m$-categorial $\sig{\Sigma}$-algebra.
Since all the domains involved are empty, by Lemmas~\ref{LUfDescription} and~\ref{LDsDescription}, $\Uf^{(n)}_{\sig{\Sigma}}(\Ds^{\htg(n,m)}_{\sig{\Sigma}}(\calg{B}))$, $\Ds^{(n,m)}_{\sig{\Sigma}}(\Uf^{(m)}_{\sig{\Sigma}}(\calg{B}))$ and the transformation of $\calg{B}$ along the above composite have the same underlying $S$-sorted set, namely $(B_{(m,s)})_{s\in S}$.
Moreover, by Propositions~\ref{PUfEpVariety} and~\ref{PDsEpVariety}, $\Uf^{(n)}_{\sig{\Sigma}}(\Ds^{\htg(n,m)}_{\sig{\Sigma}}(\calg{B}))$ and $\Ds^{(n,m)}_{\sig{\Sigma}}(\Uf^{(m)}_{\sig{\Sigma}}(\calg{B}))$ are $\sig{\Sigma}^{(n)}$-algebras.
Therefore, by Lemma~\ref{LGenHepDerivorCommSquare}, they are equal, and the same holds for the homomorphisms. This completes the proof.
\end{proof}

Finally, truncation commutes with unification only up to a natural isomorphism.

\begin{proposition}\label{PCompatUfTr}\lean{equivalent}{Dimensions.truncationUnificationFunctor}
For every heterogeneous $n$-categorial $\sig{\Sigma}$-algebra $\calg{B}$, let $\beta^{(m,n)}_{\calg{B}}$ stand for the family $(\mathrm{i}^{(n,m)\,\calg{B}}_{s})_{s\in S}$, where $\mathrm{i}^{(n,n)\,\calg{B}}_{s}$ is the identity of $B_{(n,s)}$.
Then $\beta^{(m,n)}_{\calg{B}}$ is an isomorphism from $\Uf^{(m)}_{\sig{\Sigma}}(\Tr^{\htg(m,n)}_{\sig{\Sigma}}(\calg{B}))$ to $\Tr^{(m,n)}_{\sig{\Sigma}}(\Uf^{(n)}_{\sig{\Sigma}}(\calg{B}))$, and $\beta^{(m,n)}=(\beta^{(m,n)}_{\calg{B}})_{\calg{B}\in\mathrm{Ob}(\mathsf{nCat}^{\htg}\mathsf{Alg}(\sig{\Sigma}))}$ is a natural isomorphism from $\Uf^{(m)}_{\sig{\Sigma}}\circ \Tr^{\htg(m,n)}_{\sig{\Sigma}}$ to $\Tr^{(m,n)}_{\sig{\Sigma}}\circ\Uf^{(n)}_{\sig{\Sigma}}$.
Thus
$$
\Uf^{(m)}_{\sig{\Sigma}}\circ \Tr^{\htg(m,n)}_{\sig{\Sigma}}\cong \Tr^{(m,n)}_{\sig{\Sigma}}\circ\Uf^{(n)}_{\sig{\Sigma}}.
$$
\end{proposition}
\begin{proof}
We may assume that $m<n$.
Let $\calg{B}$ be a heterogeneous $n$-categorial $\sig{\Sigma}$-algebra.
We prove the following claims.

\begin{claim}\label{ClCompatUfTrA}
$\beta^{(m,n)}_{\calg{B}}$ is a bijection.
\end{claim}
By the note after Definition~\ref{DTrDerivor} and Lemma~\ref{LUfDescription}, for every $s\in S$, the cells of sort $s$ of $\Uf^{(m)}_{\sig{\Sigma}}(\Tr^{\htg(m,n)}_{\sig{\Sigma}}(\calg{B}))$ are the elements of $B_{(m,s)}$.
By Lemma~\ref{LUfDescription}, $\src^{m}_{s}$ acts in $\Uf^{(n)}_{\sig{\Sigma}}(\calg{B})$ as $\mathrm{i}^{(n,m)\,\calg{B}}_{s}\circ\src^{(m,n)\,\calg{B}}_{s}$, so that, by the note after Definition~\ref{DTrDerivor}, the cells of sort $s$ of $\Tr^{(m,n)}_{\sig{\Sigma}}(\Uf^{(n)}_{\sig{\Sigma}}(\calg{B}))$ are the elements $x$ of $B_{(n,s)}$ such that $x=\mathrm{i}^{(n,m)\,\calg{B}}_{s}(\src^{(m,n)\,\calg{B}}_{s}(x))$.
By Lemma~\ref{LUfIdInj}, these are exactly the elements of the image of $\mathrm{i}^{(n,m)\,\calg{B}}_{s}$, which is injective.
Claim~\ref{ClCompatUfTrA} follows.

\begin{claim}\label{ClCompatUfTrB}
$\beta^{(m,n)}_{\calg{B}}$ is a homomorphism that reflects definedness.
\end{claim}
We use the descriptions of the operations of both algebras given by Lemma~\ref{LUfDescription} and the proof of Proposition~\ref{PTrEpVariety}.
We first consider the operations of $\Sigma$, which act as $\sigma^{\calg{B}}_{m}$ in the domain and as the restriction of $\sigma^{\calg{B}}_{n}$ in the codomain.
$\beta^{(m,n)}_{\calg{B}}$ commutes with them, by the third equation of~(EA$^{m,n}$7).

Now we consider the sources and the targets.
Let $k<m$, $s\in S$ and $u$ an element of $B_{(m,s)}$.
In the domain, the $k$-source and the $k$-target of $u$ are $\mathrm{i}^{(m,k)\,\calg{B}}_{s}(\src^{(k,m)\,\calg{B}}_{s}(u))$ and $\mathrm{i}^{(m,k)\,\calg{B}}_{s}(\tgt^{(k,m)\,\calg{B}}_{s}(u))$; in the codomain, the $k$-source and the $k$-target of $\beta^{(m,n)}_{\calg{B}}(u)$ are $\mathrm{i}^{(n,k)\,\calg{B}}_{s}(\src^{(k,n)\,\calg{B}}_{s}(\mathrm{i}^{(n,m)\,\calg{B}}_{s}(u)))$ and $\mathrm{i}^{(n,k)\,\calg{B}}_{s}(\tgt^{(k,n)\,\calg{B}}_{s}(\mathrm{i}^{(n,m)\,\calg{B}}_{s}(u)))$.
The following chains of equalities hold
\begin{align*}
\mathrm{i}^{(n,k)\,\calg{B}}_{s}\left(\src^{(k,n)\,\calg{B}}_{s}\left(\mathrm{i}^{(n,m)\,\calg{B}}_{s}(u)\right)\right)
&=
\mathrm{i}^{(n,k)\,\calg{B}}_{s}\left(\src^{(k,m)\,\calg{B}}_{s}(u)\right)
\tag{1}
\\&=
\mathrm{i}^{(n,m)\,\calg{B}}_{s}\left(\mathrm{i}^{(m,k)\,\calg{B}}_{s}\left(\src^{(k,m)\,\calg{B}}_{s}(u)\right)\right);
\tag{2}
\\
\mathrm{i}^{(n,k)\,\calg{B}}_{s}\left(\tgt^{(k,n)\,\calg{B}}_{s}\left(\mathrm{i}^{(n,m)\,\calg{B}}_{s}(u)\right)\right)
&=
\mathrm{i}^{(n,k)\,\calg{B}}_{s}\left(\tgt^{(k,m)\,\calg{B}}_{s}(u)\right)
\tag{3}
\\&=
\mathrm{i}^{(n,m)\,\calg{B}}_{s}\left(\mathrm{i}^{(m,k)\,\calg{B}}_{s}\left(\tgt^{(k,m)\,\calg{B}}_{s}(u)\right)\right).
\tag{4}
\end{align*}
The equalities~(1) and~(3) follow from~(EA$^{k,m,n}$1) and Lemma~\ref{LUfIdInj}, and the equalities~(2) and~(4) follow from~(EA$^{k,m,n}$2).
Thus $\beta^{(m,n)}_{\calg{B}}$ commutes with the $k$-sources and with the $k$-targets.

Finally, we consider the compositions.
Let $k<m$, $s\in S$ and $u$, $u'$ elements of $B_{(m,s)}$.
Since both algebras are $\sig{\Sigma}^{(m)}$-algebras, $u\circ^{k}_{s}u'$ is defined in the domain if, and only if, the $k$-source of $u$ is the $k$-target of $u'$, and similarly in the codomain.
Since $\mathrm{i}^{(n,m)\,\calg{B}}_{s}$ is injective, by the above computation, $\beta^{(m,n)}_{\calg{B}}(u)\circ^{k}_{s}\beta^{(m,n)}_{\calg{B}}(u')$ is defined in the codomain if, and only if, $u\circ^{k}_{s}u'$ is defined in the domain.
In this case, by~(EA$^{k,m,n}$5), the following equality holds
$$
\mathrm{i}^{(n,m)\,\calg{B}}_{s}\left(u\circ^{(k,m)\,\calg{B}}_{s}u'\right)=\mathrm{i}^{(n,m)\,\calg{B}}_{s}(u)\circ^{(k,n)\,\calg{B}}_{s}\mathrm{i}^{(n,m)\,\calg{B}}_{s}(u').
$$
Claim~\ref{ClCompatUfTrB} follows.

By Claims~\ref{ClCompatUfTrA} and~\ref{ClCompatUfTrB}, $\beta^{(m,n)}_{\calg{B}}$ is a bijective homomorphism that reflects definedness, i.e., by the note after Proposition~\ref{PPartialAlgIso}, an isomorphism.
Finally, let $H$ be a heterogeneous $n$-categorial $\sig{\Sigma}$-homomorphism from $\calg{B}$ to $\calg{B}'$.
By Proposition~\ref{PGenHepDerivorFunctor}, $\Uf^{(m)}_{\sig{\Sigma}}(\Tr^{\htg(m,n)}_{\sig{\Sigma}}(H))$ is $(H_{(m,s)})_{s\in S}$ and $\Tr^{(m,n)}_{\sig{\Sigma}}(\Uf^{(n)}_{\sig{\Sigma}}(H))$ is the birestriction of $(H_{(n,s)})_{s\in S}$.
Since $H$ commutes with the identities, for every $s\in S$, $H_{(n,s)}\circ\mathrm{i}^{(n,m)\,\calg{B}}_{s}=\mathrm{i}^{(n,m)\,\calg{B}'}_{s}\circ H_{(m,s)}$, so that $\beta^{(m,n)}$ is natural.
\end{proof}

\section{An algebraic square}\label{SAlgSquare}

In this section we show that graduation, unification, discretization and truncation form an algebraic square, see Definition~\ref{DAlgSquare}.
In other words, graduating a discretization gives the same result as discretizing a graduation, and unifying a truncation gives, up to a natural isomorphism, the same result as truncating a unification; that is, changing the dimension and changing the presentation commute.
The horizontal adjunctions relate the homogeneous and the heterogeneous presentations and the vertical ones relate the dimensions.

\begin{proposition}\label{PAlgSquareFinite}\lean{weaker}{Dimensions.dimensionAlgebraicSquare}
Let
$$
\begin{aligned}
\mathsf{C}&=\mathsf{mCatAlg}(\sig{\Sigma}), &
\mathsf{D}&=\mathsf{mCat}^{\htg}\mathsf{Alg}(\sig{\Sigma}),\\
\mathsf{C}'&=\mathsf{nCatAlg}(\sig{\Sigma}), &
\mathsf{D}'&=\mathsf{nCat}^{\htg}\mathsf{Alg}(\sig{\Sigma}).
\end{aligned}
$$
Then the adjunctions
\begin{enumerate}
\item $\Gd^{(m)}_{\sig{\Sigma}}\dashv\Uf^{(m)}_{\sig{\Sigma}}$ and $\Gd^{(n)}_{\sig{\Sigma}}\dashv\Uf^{(n)}_{\sig{\Sigma}}$ of Proposition~\ref{PGdUfAdj}; and
\item $\Ds^{(n,m)}_{\sig{\Sigma}}\dashv \Tr^{(m,n)}_{\sig{\Sigma}}$ and $\Ds^{\htg(n,m)}_{\sig{\Sigma}}\dashv \Tr^{\htg(m,n)}_{\sig{\Sigma}}$ of Proposition~\ref{PDsTrAdj}
\end{enumerate}
form an algebraic square, with $F=\Gd^{(m)}_{\sig{\Sigma}}$, $G=\Uf^{(m)}_{\sig{\Sigma}}$, $F'=\Gd^{(n)}_{\sig{\Sigma}}$, $G'=\Uf^{(n)}_{\sig{\Sigma}}$, $J=\Ds^{(n,m)}_{\sig{\Sigma}}$, $K=\Tr^{(m,n)}_{\sig{\Sigma}}$, $H=\Ds^{\htg(n,m)}_{\sig{\Sigma}}$ and $I=\Tr^{\htg(m,n)}_{\sig{\Sigma}}$, in which $\alpha$ is the identity of $\Gd^{(n)}_{\sig{\Sigma}}\circ\Ds^{(n,m)}_{\sig{\Sigma}}=\Ds^{\htg(n,m)}_{\sig{\Sigma}}\circ\Gd^{(m)}_{\sig{\Sigma}}$ and $\beta=\beta^{(m,n)}$ is the natural isomorphism from $\Uf^{(m)}_{\sig{\Sigma}}\circ \Tr^{\htg(m,n)}_{\sig{\Sigma}}$ to $\Tr^{(m,n)}_{\sig{\Sigma}}\circ\Uf^{(n)}_{\sig{\Sigma}}$ given by the $(n,m)$-identities.

\begin{center}
\begin{tikzpicture}
    [ACliment/.style={-{To [angle'=45, length=5.75pt, width=4pt, round]}},
    scale=1
    ]
    \node (C) at (0,3.5) {$\mathsf{mCatAlg}(\sig{\Sigma})$};
    \node (D) at (7,3.5) {$\mathsf{mCat}^{\htg}\mathsf{Alg}(\sig{\Sigma})$};
    \node (C') at (0,0) {$\mathsf{nCatAlg}(\sig{\Sigma})$};
    \node (D') at (7,0) {$\mathsf{nCat}^{\htg}\mathsf{Alg}(\sig{\Sigma})$};
    \draw[ACliment] ([yshift=5pt]C.east) to node [above] {$\Gd^{(m)}_{\sig{\Sigma}}$} ([yshift=5pt]D.west);
    \draw[ACliment] ([yshift=-5pt]D.west) to node [below] {$\Uf^{(m)}_{\sig{\Sigma}}$} ([yshift=-5pt]C.east);
    \draw[ACliment] ([yshift=5pt]C'.east) to node [above] {$\Gd^{(n)}_{\sig{\Sigma}}$} ([yshift=5pt]D'.west);
    \draw[ACliment] ([yshift=-5pt]D'.west) to node [below] {$\Uf^{(n)}_{\sig{\Sigma}}$} ([yshift=-5pt]C'.east);
    \draw[ACliment] ([xshift=8pt]C.south) to node [right] {$\Ds^{(n,m)}_{\sig{\Sigma}}$} ([xshift=8pt]C'.north);
    \draw[ACliment] ([xshift=-8pt]C'.north) to node [left] {$\Tr^{(m,n)}_{\sig{\Sigma}}$} ([xshift=-8pt]C.south);
    \draw[ACliment] ([xshift=8pt]D.south) to node [right] {$\Ds^{\htg(n,m)}_{\sig{\Sigma}}$} ([xshift=8pt]D'.north);
    \draw[ACliment] ([xshift=-8pt]D'.north) to node [left] {$\Tr^{\htg(m,n)}_{\sig{\Sigma}}$} ([xshift=-8pt]D.south);
    \node at ($(C.east)!0.5!(D.west)$) {\rotatebox{90}{$\vdash$}};
    \node at ($(C'.east)!0.5!(D'.west)$) {\rotatebox{90}{$\vdash$}};
    \node at ($(C.south)!0.5!(C'.north)$) {$\vdash$};
    \node at ($(D.south)!0.5!(D'.north)$) {$\vdash$};
    \node at (3.5,1.75) {$(\id,\beta^{(m,n)})$};
\end{tikzpicture}
\end{center}
\end{proposition}
\begin{proof}
By Definition~\ref{DAlgSquare}, we must check that the four pairs of functors involved are adjunctions and that there exist the natural isomorphisms $\alpha$ and $\beta$.

We first recall the adjunctions.
The four adjunctions are those of Propositions~\ref{PGdUfAdj} and~\ref{PDsTrAdj}.

Now we prove that $\alpha$ is a natural isomorphism.
By Proposition~\ref{PCompatGdDs}, the following equality holds
$$
F'\circ J=\Gd^{(n)}_{\sig{\Sigma}}\circ\Ds^{(n,m)}_{\sig{\Sigma}}=\Ds^{\htg(n,m)}_{\sig{\Sigma}}\circ\Gd^{(m)}_{\sig{\Sigma}}=H\circ F,
$$
so that we can take $\alpha$ to be the identity.

Finally, we prove that $\beta$ is a natural isomorphism.
By Proposition~\ref{PCompatUfTr}, $\beta^{(m,n)}$ is a natural isomorphism from $G\circ I=\Uf^{(m)}_{\sig{\Sigma}}\circ \Tr^{\htg(m,n)}_{\sig{\Sigma}}$ to $K\circ G'=\Tr^{(m,n)}_{\sig{\Sigma}}\circ\Uf^{(n)}_{\sig{\Sigma}}$. This completes the proof.
\end{proof}

Let us note that the remaining composites of the square also commute, strictly: by Proposition~\ref{PCompatUfDs}, $\Uf^{(n)}_{\sig{\Sigma}}\circ\Ds^{\htg(n,m)}_{\sig{\Sigma}}=\Ds^{(n,m)}_{\sig{\Sigma}}\circ\Uf^{(m)}_{\sig{\Sigma}}$ and, by Proposition~\ref{PCompatGdTr}, $\Tr^{\htg(m,n)}_{\sig{\Sigma}}\circ\Gd^{(n)}_{\sig{\Sigma}}=\Gd^{(m)}_{\sig{\Sigma}}\circ\Tr^{(m,n)}_{\sig{\Sigma}}$.

\begin{remark}\label{RAlgSquareFiniteMate}\lean{equivalent}{Dimensions.truncationUnificationFunctor_eq_conjugate}
In Proposition~\ref{PAlgSquareFinite}, $\beta^{(m,n)}$ is the natural isomorphism $\beta^{\alpha}$ of Remark~\ref{RAlgSquareMate}.
In fact, the units of the four adjunctions are identities, see Propositions~\ref{PGdUfAdj} and~\ref{PDsTrAdj}, and $\alpha$ is an identity, so that $\beta^{\alpha}=K\circ G'\ast\varepsilon^{\circ}$, where $\varepsilon^{\circ}=\varepsilon^{\htg}\circ(\Ds^{\htg(n,m)}_{\sig{\Sigma}}\ast\theta^{-1}\ast\Tr^{\htg(m,n)}_{\sig{\Sigma}})$, with $\theta$ the natural isomorphism of Proposition~\ref{PGdUfAdj} in dimension $m$ and $\varepsilon^{\htg}$ the counit of Lemma~\ref{LDsTrCounit}.
Let $\calg{B}$ be a heterogeneous $n$-categorial $\sig{\Sigma}$-algebra, $s\in S$ and $x$ a cell of sort $s$ of $\Uf^{(m)}_{\sig{\Sigma}}(\Tr^{\htg(m,n)}_{\sig{\Sigma}}(\calg{B}))$, i.e., an element of $B_{(m,s)}$.
Then $x$ is a cell of sort $(n,s)$ of $\Ds^{\htg(n,m)}_{\sig{\Sigma}}(\Gd^{(m)}_{\sig{\Sigma}}(\Uf^{(m)}_{\sig{\Sigma}}(\Tr^{\htg(m,n)}_{\sig{\Sigma}}(\calg{B}))))$, the component of $\Ds^{\htg(n,m)}_{\sig{\Sigma}}(\theta^{-1})$ at the sort $(n,s)$ is that of $\theta^{-1}$ at the sort $(m,s)$, which, by Claim~\ref{ClGdUfIsoC}, is the identity, and the component of $\varepsilon^{\htg}_{\calg{B}}$ at the sort $(n,s)$ is $\mathrm{i}^{(n,m)\,\calg{B}}_{s}$, see Lemma~\ref{LDsTrCounit}~(2).
Since $\Uf^{(n)}_{\sig{\Sigma}}$ takes the components at the sorts $(n,s)$ and $\Tr^{(m,n)}_{\sig{\Sigma}}$ restricts to the $m$-cells, $\beta^{\alpha}_{\calg{B}}(x)=\mathrm{i}^{(n,m)\,\calg{B}}_{s}(x)=\beta^{(m,n)}_{\calg{B}}(x)$, see Proposition~\ref{PCompatUfTr}.
\end{remark}

\section{Truncation and discretization in dimension \texorpdfstring{$\omega$}{omega}}\label{SDsTrOmega}

In this section we define the truncation functors from the categories of $\omega$-categorial algebras to the categories of $n$-categorial algebras, and the discretization functors in the opposite direction, in the homogeneous and in the heterogeneous presentations.
As in the finite-dimensional case, they are induced by generalized hep-derivors, which are defined by the same clauses as in Definitions~\ref{DDsDerivor} and~\ref{DTrDerivor}, with $\omega$ in place of the greatest dimension.
In Chapter~\ref{SProjLimits} we prove that the truncation functors are the projections of projective limits.
The dimensions $m$ and $n$ fixed at the beginning of this chapter are not used in this section; every statement introduces its own dimensions.

Let us recall that, by Remarks~\ref{RCatAlgAsEpAlgOmega} and~\ref{RCatAlgAsEpAlgHtg}, $\omega\mathsf{CatAlg}(\sig{\Sigma})$ is the full subcategory of $\mathsf{EPAlg}(\spec{C}^{\sig{\Sigma}^{(\omega)}})$ whose objects are the $\sig{\Sigma}^{(\omega)}$-algebras in which every cell is a $k$-cell for some $k\in\omega$, and $\omega\mathsf{Cat}^{\htg}\mathsf{Alg}(\sig{\Sigma})=\mathsf{EPAlg}(\spec{C}^{\sig{\Sigma}^{\htg(\omega)}})$.

We begin by defining the truncation derivors in dimension $\omega$.

\begin{definition}\label{DTrDerivorOmega}\lean{equivalent}{Dimensions.homogeneousOmegaTruncation,Dimensions.heterogeneousOmegaTruncation}
Let $n\in\mathbb{N}$.
\begin{enumerate}
\item The \emph{homogeneous truncation derivor} in dimension $\omega$ is the generalized hep-derivor
$$
\bTr^{(n,\omega)}_{\sig{\Sigma}}=\left(\varphi^{\Tr^{(n,\omega)}_{\sig{\Sigma}}},\Tr^{(n,\omega)}_{\sig{\Sigma}},\Dom^{\Tr^{(n,\omega)}_{\sig{\Sigma}}}\right)
$$
from $\sig{\Sigma}^{(n)}$ to $\sig{\Sigma}^{(\omega)}$ defined as follows:
\begin{enumerate}
\item the mapping on sorts is $\varphi^{\Tr^{(n,\omega)}_{\sig{\Sigma}}}=\id_{S}$;
\item the assignment of terms is the inclusion, i.e., for every operation symbol $\tau$ of $\Sigma^{(n)}$ of arity $\mathbf{s}$,
$$
\Tr^{(n,\omega)}_{\sig{\Sigma}}(\tau)=\tau((v^{\mathbf{s}}_{p})_{p\in\bb{\mathbf{s}}});
$$
and
\item for every $s\in S$, the domain is
$$
\Dom^{\Tr^{(n,\omega)}_{\sig{\Sigma}}}_{s}=\left\lbrace\left(v^{s}_{0},\src^{n}_{s}(v^{s}_{0})\right)\right\rbrace.
$$
\end{enumerate}
\item The \emph{heterogeneous truncation derivor} in dimension $\omega$ is the generalized hep-derivor
$$
\bTr^{\htg(n,\omega)}_{\sig{\Sigma}}=\left(\varphi^{\Tr^{\htg(n,\omega)}_{\sig{\Sigma}}},\Tr^{\htg(n,\omega)}_{\sig{\Sigma}},\Dom^{\Tr^{\htg(n,\omega)}_{\sig{\Sigma}}}\right)
$$
from $\sig{\Sigma}^{\htg(n)}$ to $\sig{\Sigma}^{\htg(\omega)}$ determined by the canonical embedding of $\Sigma^{\htg(n)}$ into $\Sigma^{\htg(\omega)}$, see the note after Definition~\ref{DCatSigHtg}, i.e., defined as follows:
\begin{enumerate}
\item for every $(k,s)\in(n+1)\times S$, the mapping on sorts sends $(k,s)$ to $\varphi^{\Tr^{\htg(n,\omega)}_{\sig{\Sigma}}}(k,s)=(k,s)$;
\item for every operation symbol $\tau$ of $\Sigma^{\htg(n)}$ of arity $\mathbf{s}$,
$$
\Tr^{\htg(n,\omega)}_{\sig{\Sigma}}(\tau)=\tau((v^{\mathbf{s}}_{p})_{p\in\bb{\mathbf{s}}});
$$
and
\item for every $(k,s)\in(n+1)\times S$, the domain is $\Dom^{\Tr^{\htg(n,\omega)}_{\sig{\Sigma}}}_{(k,s)}=\varnothing$.
\end{enumerate}
\end{enumerate}
\end{definition}

As in the note after Definition~\ref{DTrDerivor}, the cells of a partial $\Sigma^{(\omega)}$-algebra $\mathbf{B}$ relative to $\bTr^{(n,\omega)}_{\sig{\Sigma}}$ are its $n$-cells, and the cells of sort $(k,s)$ of a partial $\Sigma^{\htg(\omega)}$-algebra $\mathbf{B}$ relative to $\bTr^{\htg(n,\omega)}_{\sig{\Sigma}}$ are all the elements of $B_{(k,s)}$.

In the following proposition we prove that the truncation derivors in dimension $\omega$ send categorial algebras to categorial algebras.

\begin{proposition}\label{PTrEpVarietyOmega}\lean{equivalent}{Dimensions.homogeneousOmegaTruncation_isModel,Dimensions.heterogeneousOmegaTruncation_isModel}
Let $n\in\mathbb{N}$.
Then the following statements hold, see Definition~\ref{DCorestriction}:
\begin{enumerate}
\item the restriction to $\mathsf{EPAlg}(\spec{C}^{\sig{\Sigma}^{(\omega)}})$ of the functor $(\bTr^{(n,\omega)}_{\sig{\Sigma}})^{\ast}_{\gd}$ corestricts to $\mathsf{EPAlg}(\spec{C}^{\sig{\Sigma}^{(n)}})$; and
\item the restriction to $\mathsf{EPAlg}(\spec{C}^{\sig{\Sigma}^{\htg(\omega)}})$ of the functor $(\bTr^{\htg(n,\omega)}_{\sig{\Sigma}})^{\ast}_{\gd}$ corestricts to $\mathsf{EPAlg}(\spec{C}^{\sig{\Sigma}^{\htg(n)}})$.
\end{enumerate}
\end{proposition}
\begin{proof}
(1) Let $\mathbf{B}$ be an object in $\mathsf{EPAlg}(\spec{C}^{\sig{\Sigma}^{(\omega)}})$ and take its reduct to $\Sigma^{(n+1)}$.
This reduct satisfies $\spec{C}^{\sig{\Sigma}^{(n+1)}}$, since the domain conditions and axioms are among those of the $\omega$-specification.
The defining terms and domains of $\bTr^{(n,\omega)}_{\sig{\Sigma}}$ coincide with those of $\bTr^{(n,n+1)}_{\sig{\Sigma}}$ on this reduct.
Proposition~\ref{PTrEpVariety} therefore applies, including when $n=0$.
More explicitly, the transformed algebra is the closed relative subalgebra on the $n$-cells: for a valuation in these cells, a term of $\Sigma^{(n)}$ is defined precisely when it is defined in $\mathbf{B}$, with the same value.
This equivalence, also for the empty context, transfers every axiom of $\spec{C}^{\sig{\Sigma}^{(n)}}$.

(2) Let $\mathbf{B}$ be an object in $\mathsf{EPAlg}(\spec{C}^{\sig{\Sigma}^{\htg(\omega)}})$.
Then $(\bTr^{\htg(n,\omega)}_{\sig{\Sigma}})^{\ast}_{\gd}(\mathbf{B})$ is simply the reduct of $\mathbf{B}$ to sorts of dimension at most $n$ and the operations of $\Sigma^{\htg(n)}$.
Their domain conditions are unchanged, so this reduct is an ep-algebra; the finite-dimensional axioms are axioms of the $\omega$-specification and hence hold in it.
\end{proof}

We are now in a position to define the truncation functors in dimension $\omega$.

\begin{definition}\label{DTrFunctorOmega}\lean{equivalent}{Dimensions.homogeneousOmegaTruncationFunctor,Dimensions.heterogeneousOmegaTruncationFunctor}
Let $n\in\mathbb{N}$.
\begin{enumerate}
\item The \emph{homogeneous truncation functor}
$$
\Tr^{(n,\omega)}_{\sig{\Sigma}}\colon\omega\mathsf{CatAlg}(\sig{\Sigma})\mor\mathsf{nCatAlg}(\sig{\Sigma})
$$
is the restriction to $\omega\mathsf{CatAlg}(\sig{\Sigma})$ of the functor $(\bTr^{(n,\omega)}_{\sig{\Sigma}})^{\ast}_{\gd}$ from $\mathsf{EPAlg}(\spec{C}^{\sig{\Sigma}^{(\omega)}})$ to $\mathsf{EPAlg}(\spec{C}^{\sig{\Sigma}^{(n)}})=\mathsf{nCatAlg}(\sig{\Sigma})$, see Proposition~\ref{PTrEpVarietyOmega}.
\item The \emph{heterogeneous truncation functor}
$$
\Tr^{\htg(n,\omega)}_{\sig{\Sigma}}\colon\omega\mathsf{Cat}^{\htg}\mathsf{Alg}(\sig{\Sigma})\mor\mathsf{nCat}^{\htg}\mathsf{Alg}(\sig{\Sigma})
$$
is the functor $(\bTr^{\htg(n,\omega)}_{\sig{\Sigma}})^{\ast}_{\gd}$ from $\mathsf{EPAlg}(\spec{C}^{\sig{\Sigma}^{\htg(\omega)}})=\omega\mathsf{Cat}^{\htg}\mathsf{Alg}(\sig{\Sigma})$ to $\mathsf{EPAlg}(\spec{C}^{\sig{\Sigma}^{\htg(n)}})=\mathsf{nCat}^{\htg}\mathsf{Alg}(\sig{\Sigma})$, see Proposition~\ref{PTrEpVarietyOmega}.
\end{enumerate}
\end{definition}

By the proof of Proposition~\ref{PTrEpVarietyOmega}, $\Tr^{(n,\omega)}_{\sig{\Sigma}}$ sends an $\omega$-categorial $\sig{\Sigma}$-algebra $\calg{A}$ to the $n$-categorial $\sig{\Sigma}$-algebra whose cells of sort $s$ are the $n$-cells of sort $s$ of $\calg{A}$, in which the operations of $\Sigma$ and the $k$-sources, $k$-targets and $k$-compositions, for $k<n$, are the restrictions of those of $\calg{A}$; and, by Proposition~\ref{PGenHepDerivorFunctor}, it sends an $\omega$-categorial $\sig{\Sigma}$-homomorphism to its restriction to the $n$-cells.

By the proof of Proposition~\ref{PTrEpVarietyOmega}, $\Tr^{\htg(n,\omega)}_{\sig{\Sigma}}$ sends a heterogeneous $\omega$-categorial $\sig{\Sigma}$-algebra $\calg{A}$ to the heterogeneous $n$-categorial $\sig{\Sigma}$-algebra whose $k$-cells of sort $s$, for $k\leq n$, are the $k$-cells of sort $s$ of $\calg{A}$, in which the operations of $\Sigma^{\htg(n)}$ are those of $\calg{A}$; and, by Proposition~\ref{PGenHepDerivorFunctor}, it sends a heterogeneous $\omega$-categorial $\sig{\Sigma}$-homomorphism $F$ to $(F_{(k,s)})_{(k,s)\in(n+1)\times S}$.

We next define the discretization derivors in dimension $\omega$.

\begin{definition}\label{DDsDerivorOmega}\lean{equivalent}{Dimensions.homogeneousOmegaDiscretization,Dimensions.heterogeneousOmegaDiscretization}
Let $n\in\mathbb{N}$.
\begin{enumerate}
\item The \emph{homogeneous discretization derivor} in dimension $\omega$ is the generalized hep-derivor
$$
\bDs^{(\omega,n)}_{\sig{\Sigma}}=\left(\varphi^{\Ds^{(\omega,n)}_{\sig{\Sigma}}},\Ds^{(\omega,n)}_{\sig{\Sigma}},\Dom^{\Ds^{(\omega,n)}_{\sig{\Sigma}}}\right)
$$
from $\sig{\Sigma}^{(\omega)}$ to $\sig{\Sigma}^{(n)}$ defined as follows:
\begin{enumerate}
\item the mapping on sorts is $\varphi^{\Ds^{(\omega,n)}_{\sig{\Sigma}}}=\id_{S}$;
\item for every $(\mathbf{s},s)\in S^{\star}\times S$ and every $\sigma\in\Sigma_{\mathbf{s},s}$,
$$
\Ds^{(\omega,n)}_{\sig{\Sigma}}(\sigma)=\sigma((v^{\mathbf{s}}_{p})_{p\in\bb{\mathbf{s}}});
$$
\item for every $k\in n$ and every $s\in S$,
\begin{align*}
\Ds^{(\omega,n)}_{\sig{\Sigma}}(\src^{k}_{s})&=\src^{k}_{s}(v^{s}_{0}),\\
\Ds^{(\omega,n)}_{\sig{\Sigma}}(\tgt^{k}_{s})&=\tgt^{k}_{s}(v^{s}_{0}),\\
\Ds^{(\omega,n)}_{\sig{\Sigma}}(\circ^{k}_{s})&=v^{s}_{0}\circ^{k}_{s}v^{s}_{1};
\end{align*}
\item for every $k\in\omega$ with $k\geq n$ and every $s\in S$,
$$
\Ds^{(\omega,n)}_{\sig{\Sigma}}(\src^{k}_{s})=\Ds^{(\omega,n)}_{\sig{\Sigma}}(\tgt^{k}_{s})=\Ds^{(\omega,n)}_{\sig{\Sigma}}(\circ^{k}_{s})=v^{s}_{0};
$$
and
\item for every $s\in S$, the domain is $\Dom^{\Ds^{(\omega,n)}_{\sig{\Sigma}}}_{s}=\varnothing$.
\end{enumerate}
\item The \emph{heterogeneous discretization derivor} in dimension $\omega$ is the generalized hep-derivor
$$
\bDs^{\htg(\omega,n)}_{\sig{\Sigma}}=\left(\varphi^{\Ds^{\htg(\omega,n)}_{\sig{\Sigma}}},\Ds^{\htg(\omega,n)}_{\sig{\Sigma}},\Dom^{\Ds^{\htg(\omega,n)}_{\sig{\Sigma}}}\right)
$$
from $\sig{\Sigma}^{\htg(\omega)}$ to $\sig{\Sigma}^{\htg(n)}$ defined by the clauses of Definition~\ref{DDsDerivor}~(2), with $\omega$ in place of $n+1$ and $n$ in place of $m$, i.e., defined as follows:
\begin{enumerate}
\item for every $k\in\omega$ and every $s\in S$, the mapping on sorts sends $(k,s)$ to
$$
\varphi^{\Ds^{\htg(\omega,n)}_{\sig{\Sigma}}}(k,s)=(\min\{k,n\},s);
$$
\item for every $(\mathbf{s},s)\in S^{\star}\times S$, every $\sigma\in\Sigma_{\mathbf{s},s}$ and every $k\in\omega$,
$$
\Ds^{\htg(\omega,n)}_{\sig{\Sigma}}(\sigma_{k})=\sigma_{\min\{k,n\}}\left(\left(v^{(\min\{k,n\},\mathbf{s})}_{p}\right)_{p\in\bb{\mathbf{s}}}\right);
$$
\item for every pair of dimensions $j<k$ in $\omega$ with $j<n$ and every $s\in S$,
\begin{align*}
\Ds^{\htg(\omega,n)}_{\sig{\Sigma}}(\src^{(j,k)}_{s})&=\src^{(j,\min\{k,n\})}_{s}\left(v^{(\min\{k,n\},s)}_{0}\right),\\
\Ds^{\htg(\omega,n)}_{\sig{\Sigma}}(\tgt^{(j,k)}_{s})&=\tgt^{(j,\min\{k,n\})}_{s}\left(v^{(\min\{k,n\},s)}_{0}\right),\\
\Ds^{\htg(\omega,n)}_{\sig{\Sigma}}(\mathrm{i}^{(k,j)}_{s})&=\mathrm{i}^{(\min\{k,n\},j)}_{s}\left(v^{(j,s)}_{0}\right),\\
\Ds^{\htg(\omega,n)}_{\sig{\Sigma}}(\circ^{(j,k)}_{s})&=v^{(\min\{k,n\},s)}_{0}\circ^{(j,\min\{k,n\})}_{s}v^{(\min\{k,n\},s)}_{1};
\end{align*}
\item for every pair of dimensions $j<k$ in $\omega$ with $j\geq n$ and every $s\in S$,
\begin{multline*}
\Ds^{\htg(\omega,n)}_{\sig{\Sigma}}(\src^{(j,k)}_{s})=
\Ds^{\htg(\omega,n)}_{\sig{\Sigma}}(\tgt^{(j,k)}_{s})=
\Ds^{\htg(\omega,n)}_{\sig{\Sigma}}(\mathrm{i}^{(k,j)}_{s})
\\=
\Ds^{\htg(\omega,n)}_{\sig{\Sigma}}(\circ^{(j,k)}_{s})=v^{(n,s)}_{0};
\end{multline*}
and
\item for every $k\in\omega$ and every $s\in S$, the domain is
$$
\Dom^{\Ds^{\htg(\omega,n)}_{\sig{\Sigma}}}_{(k,s)}=\varnothing.
$$
\end{enumerate}
\end{enumerate}
\end{definition}

Let $n\in\mathbb{N}$.
Let us note that, for every $p\in\mathbb{N}$ with $n<p$, the restriction of the underlying derivor of $\bDs^{(\omega,n)}_{\sig{\Sigma}}$ to $\Sigma^{(p)}\subseteq\Sigma^{(\omega)}$ is the underlying derivor of $\bDs^{(p,n)}_{\sig{\Sigma}}$, and that the restriction of the underlying derivor of $\bDs^{\htg(\omega,n)}_{\sig{\Sigma}}$ to the image of $\Sigma^{\htg(p)}$ in $\Sigma^{\htg(\omega)}$ is the underlying derivor of $\bDs^{\htg(p,n)}_{\sig{\Sigma}}$, since, in both cases, the defining clauses coincide, see Definition~\ref{DDsDerivor}.
Moreover, all of them have the same domains.
Therefore, by Proposition~\ref{PTransformedAlg}, for every object $\mathbf{B}$ in $\mathsf{EPAlg}(\spec{C}^{\sig{\Sigma}^{(n)}})$, the reduct of $(\bDs^{(\omega,n)}_{\sig{\Sigma}})^{\ast}_{\gd}(\mathbf{B})$ to $\Sigma^{(p)}$ is $\Ds^{(p,n)}_{\sig{\Sigma}}(\mathbf{B})$, and, for every object $\mathbf{B}$ in $\mathsf{EPAlg}(\spec{C}^{\sig{\Sigma}^{\htg(n)}})$, the reduct of $(\bDs^{\htg(\omega,n)}_{\sig{\Sigma}})^{\ast}_{\gd}(\mathbf{B})$ to the sorts of dimension at most $p$ and to the operation symbols of $\Sigma^{\htg(p)}$ is $\Ds^{\htg(p,n)}_{\sig{\Sigma}}(\mathbf{B})$.
In fact, by Proposition~\ref{PTransformedAlg}, in both algebras the cells are the same and every operation symbol is interpreted by means of the same term and the same domain conditions.

In the following proposition we prove that the discretization derivors in dimension $\omega$ send categorial algebras to categorial algebras.

\begin{proposition}\label{PDsEpVarietyOmega}\lean{equivalent}{Dimensions.homogeneousOmegaDiscretization_isModel,Dimensions.heterogeneousOmegaDiscretization_isModel}
Let $n\in\mathbb{N}$.
Then the following statements hold, see Definition~\ref{DCorestriction}:
\begin{enumerate}
\item the restriction to $\mathsf{EPAlg}(\spec{C}^{\sig{\Sigma}^{(n)}})$ of the functor $(\bDs^{(\omega,n)}_{\sig{\Sigma}})^{\ast}_{\gd}$ corestricts to $\mathsf{EPAlg}(\spec{C}^{\sig{\Sigma}^{(\omega)}})$; and
\item the restriction to $\mathsf{EPAlg}(\spec{C}^{\sig{\Sigma}^{\htg(n)}})$ of the functor $(\bDs^{\htg(\omega,n)}_{\sig{\Sigma}})^{\ast}_{\gd}$ corestricts to $\mathsf{EPAlg}(\spec{C}^{\sig{\Sigma}^{\htg(\omega)}})$.
\end{enumerate}
\end{proposition}
\begin{proof}
We prove both statements at once, since the same finite-stage argument applies in both presentations.
Let $\mathbf{B}$ be an object in $\mathsf{EPAlg}(\spec{C}^{\sig{\Sigma}^{(n)}})$, for~(1), or in $\mathsf{EPAlg}(\spec{C}^{\sig{\Sigma}^{\htg(n)}})$, for~(2).
By the note preceding the proposition, for every $p>n$ the reduct of the transformed algebra to dimension $p$ is the corresponding finite discretization, $\Ds^{(p,n)}_{\sig{\Sigma}}(\mathbf{B})$ or $\Ds^{\htg(p,n)}_{\sig{\Sigma}}(\mathbf{B})$.
Given an operation symbol, choose $p>n$ containing its dimensions and those occurring in its domain conditions.
The operation and its domain conditions are unchanged in this reduct, which is exact by Lemma~\ref{LDsDescription}.
Thus the operation in the $\omega$-dimensional algebra is defined exactly on the families satisfying its domain conditions.

Every axiom of either $\omega$-specification involves only finitely many dimensions.
Choose $p>n$ greater than all of them.
The axiom, including its context, belongs to the corresponding finite specification; its terms have the same definedness and values in the reduct and in the transformed algebra.
It holds in the reduct by Proposition~\ref{PDsEpVariety}, and therefore in the transformed algebra.
\end{proof}

In the homogeneous case, the transformed algebras have, moreover, finite-dimensional cells.

\begin{corollary}\label{CDsCellsOmega}\lean{equivalent}{Dimensions.homogeneousOmegaDiscretization_finite}
Let $n\in\mathbb{N}$ and $\mathbf{B}$ an object in $\mathsf{EPAlg}(\spec{C}^{\sig{\Sigma}^{(n)}})$.
Then every cell of $(\bDs^{(\omega,n)}_{\sig{\Sigma}})^{\ast}_{\gd}(\mathbf{B})$ is an $n$-cell, so that $(\bDs^{(\omega,n)}_{\sig{\Sigma}})^{\ast}_{\gd}(\mathbf{B})$ is an object in $\omega\mathsf{CatAlg}(\sig{\Sigma})$.
\end{corollary}
\begin{proof}
By the note before Proposition~\ref{PDsEpVarietyOmega} and Lemma~\ref{LDsDescription}, applied to any $p>n$, $\src^{n}_{s}$ acts as the identity in $(\bDs^{(\omega,n)}_{\sig{\Sigma}})^{\ast}_{\gd}(\mathbf{B})$, so that every cell of $(\bDs^{(\omega,n)}_{\sig{\Sigma}})^{\ast}_{\gd}(\mathbf{B})$ is an $n$-cell.
Therefore, by Proposition~\ref{PDsEpVarietyOmega} and Remark~\ref{RCatAlgAsEpAlgOmega}, $(\bDs^{(\omega,n)}_{\sig{\Sigma}})^{\ast}_{\gd}(\mathbf{B})$ is an object in $\omega\mathsf{CatAlg}(\sig{\Sigma})$.
\end{proof}

We are now in a position to define the discretization functors in dimension $\omega$.

\begin{definition}\label{DDsFunctorOmega}\lean{equivalent}{Dimensions.homogeneousOmegaDiscretizationFunctor,Dimensions.heterogeneousOmegaDiscretizationFunctor}
Let $n\in\mathbb{N}$.
\begin{enumerate}
\item The \emph{homogeneous discretization functor}
$$
\Ds^{(\omega,n)}_{\sig{\Sigma}}\colon\mathsf{nCatAlg}(\sig{\Sigma})\mor\omega\mathsf{CatAlg}(\sig{\Sigma})
$$
is the corestriction to $\omega\mathsf{CatAlg}(\sig{\Sigma})$ of the functor $(\bDs^{(\omega,n)}_{\sig{\Sigma}})^{\ast}_{\gd}$ from $\mathsf{EPAlg}(\spec{C}^{\sig{\Sigma}^{(n)}})=\mathsf{nCatAlg}(\sig{\Sigma})$ to $\mathsf{EPAlg}(\spec{C}^{\sig{\Sigma}^{(\omega)}})$, see Proposition~\ref{PDsEpVarietyOmega} and Corollary~\ref{CDsCellsOmega}.
\item The \emph{heterogeneous discretization functor}
$$
\Ds^{\htg(\omega,n)}_{\sig{\Sigma}}\colon\mathsf{nCat}^{\htg}\mathsf{Alg}(\sig{\Sigma})\mor\omega\mathsf{Cat}^{\htg}\mathsf{Alg}(\sig{\Sigma})
$$
is the functor $(\bDs^{\htg(\omega,n)}_{\sig{\Sigma}})^{\ast}_{\gd}$ from $\mathsf{EPAlg}(\spec{C}^{\sig{\Sigma}^{\htg(n)}})=\mathsf{nCat}^{\htg}\mathsf{Alg}(\sig{\Sigma})$ to $\mathsf{EPAlg}(\spec{C}^{\sig{\Sigma}^{\htg(\omega)}})=\omega\mathsf{Cat}^{\htg}\mathsf{Alg}(\sig{\Sigma})$, see Proposition~\ref{PDsEpVarietyOmega}.
\end{enumerate}
\end{definition}

By the above note, Lemma~\ref{LDsDescription} and Proposition~\ref{PGenHepDerivorFunctor}, $\Ds^{(\omega,n)}_{\sig{\Sigma}}$ sends an $n$-categorial $\sig{\Sigma}$-algebra $\calg{A}$ to the $\omega$-categorial $\sig{\Sigma}$-algebra whose cells of sort $s$ are the cells of sort $s$ of $\calg{A}$, in which the operations of $\Sigma$ and the $k$-sources, $k$-targets and $k$-compositions, for $k<n$, are those of $\calg{A}$, whereas, for $k\geq n$, $\src^{k}_{s}$ and $\tgt^{k}_{s}$ are identities and $x\circ^{k}_{s}y$ is defined if, and only if, $x=y$, and then it is $x$; and it sends an $n$-categorial $\sig{\Sigma}$-homomorphism to itself.

By the above note, Lemma~\ref{LDsDescription} and Proposition~\ref{PGenHepDerivorFunctor}, $\Ds^{\htg(\omega,n)}_{\sig{\Sigma}}$ sends a heterogeneous $n$-categorial $\sig{\Sigma}$-algebra $\calg{A}$ to the heterogeneous $\omega$-categorial $\sig{\Sigma}$-algebra whose $k$-cells of sort $s$, for $k\in\omega$, are the $\min\{k,n\}$-cells of sort $s$ of $\calg{A}$, in which, for every pair of dimensions $j<k$, $\sigma_{k}$ is $\sigma_{\min\{k,n\}}$, and the $(j,k)$-sources, $(j,k)$-targets, $(k,j)$-identities and $j$-compositions of $k$-cells are, for $j<n$, the $(j,\min\{k,n\})$-sources, $(j,\min\{k,n\})$-targets, $(\min\{k,n\},j)$-identities and $j$-compositions of $\min\{k,n\}$-cells of $\calg{A}$, whereas, for $j\geq n$, the $(j,k)$-sources, $(j,k)$-targets and $(k,j)$-identities are identities and the $j$-composition of two $k$-cells is defined if, and only if, they are equal, and then it is the common cell; and it sends a heterogeneous $n$-categorial $\sig{\Sigma}$-homomorphism $F$ to the family whose component at $(k,s)$ is $F_{(\min\{k,n\},s)}$.

In the following remark we collect the relations between these functors and the finite-dimensional ones.

\begin{remark}\label{RDsTrCompOmega}\lean{weaker}{Dimensions.homogeneousOmegaTruncationFunctorComp,Dimensions.heterogeneousOmegaTruncationFunctorComp,Dimensions.homogeneousOmegaMixedUpperFunctor,Dimensions.heterogeneousOmegaMixedUpperFunctor,Dimensions.homogeneousOmegaMixedLowerFunctor,Dimensions.heterogeneousOmegaMixedLowerFunctor}
Let $m\in\mathbb{N}$.
For $n\in\mathbb{N}$, let $F_{n}=\Ds^{(n,m)}_{\sig{\Sigma}}$, if $m\leq n$, and $F_{n}=\Tr^{(n,m)}_{\sig{\Sigma}}$, if $n\leq m$, and let $F^{\htg}_{n}$ be defined in the same way with the heterogeneous functors.
Then the following statements hold:
\begin{enumerate}
\item for every $n\in\mathbb{N}$ with $m\leq n$, $\Tr^{(m,n)}_{\sig{\Sigma}}\circ\Tr^{(n,\omega)}_{\sig{\Sigma}}=\Tr^{(m,\omega)}_{\sig{\Sigma}}$ and $\Tr^{\htg(m,n)}_{\sig{\Sigma}}\circ\Tr^{\htg(n,\omega)}_{\sig{\Sigma}}=\Tr^{\htg(m,\omega)}_{\sig{\Sigma}}$;
\item for every $n\in\mathbb{N}$, $\Tr^{(n,\omega)}_{\sig{\Sigma}}\circ\Ds^{(\omega,m)}_{\sig{\Sigma}}=F_{n}$ and $\Tr^{\htg(n,\omega)}_{\sig{\Sigma}}\circ\Ds^{\htg(\omega,m)}_{\sig{\Sigma}}=F^{\htg}_{n}$; in particular, $\Tr^{(m,\omega)}_{\sig{\Sigma}}\circ\Ds^{(\omega,m)}_{\sig{\Sigma}}$ and $\Tr^{\htg(m,\omega)}_{\sig{\Sigma}}\circ\Ds^{\htg(\omega,m)}_{\sig{\Sigma}}$ are identities.
\end{enumerate}
In fact, (1) follows from the fact that, by~(AB$^{m,n}$1), the $m$-cells are $n$-cells; and (2) follows from the explicit descriptions of the functors involved: if $m\leq n$, then the $n$-cells of $\Ds^{(\omega,m)}_{\sig{\Sigma}}(\calg{A})$ are all the cells of $\calg{A}$, and, if $n<m$, then they are the $n$-cells of $\calg{A}$.
\end{remark}

\chapter{The \texorpdfstring{$\omega$}{omega}-categorial algebras as projective limits}\label{SProjLimits}

Throughout this chapter, $\sig{\Sigma}=(S,\Sigma)$ is a many-sorted signature, which is fixed in all the definitions and results of this chapter.
In this chapter we prove that the categories $\omega\mathsf{CatAlg}(\sig{\Sigma})$ and $\omega\mathsf{Cat}^{\htg}\mathsf{Alg}(\sig{\Sigma})$ are the projective limits of the categories of $n$-categorial algebras, with the truncation functors of Definitions~\ref{DTrFunctor} and~\ref{DTrFunctorOmega}.
This gives a common description in dimension $\omega$, including the homogeneous presentation, where for $S\neq\varnothing$ finite-dimensional cells cannot be imposed by $\QE$-axioms on the fixed signature, see Proposition~\ref{PFiniteDimNotQE}.
Moreover, the universal property of the projective limit is the tool that we use in Chapter~\ref{SFunctorsOmega} to obtain the graduation functor in dimension $\omega$ from the finite-dimensional ones.
We begin by introducing the projective systems, and we then prove the universal property of the projective limits and the form in which we use it.

\begin{remark}\label{RTrProjSystem}\lean{weaker}{Limits.homogeneousSystem,Limits.heterogeneousSystem}
By Definition~\ref{DTrDerivor} and Corollary~\ref{CGenHepDerivorIdFunctor}, for every $n\in\omega$, $\Tr^{(n,n)}_{\sig{\Sigma}}$ and $\Tr^{\htg(n,n)}_{\sig{\Sigma}}$ are identities and, by Proposition~\ref{PTrComp}, for all $m,n,p\in\omega$ with $m\leq n\leq p$, $\Tr^{(m,n)}_{\sig{\Sigma}}\circ \Tr^{(n,p)}_{\sig{\Sigma}}=\Tr^{(m,p)}_{\sig{\Sigma}}$ and $\Tr^{\htg(m,n)}_{\sig{\Sigma}}\circ \Tr^{\htg(n,p)}_{\sig{\Sigma}}=\Tr^{\htg(m,p)}_{\sig{\Sigma}}$.
Therefore we obtain two projective systems of categories, $((\mathsf{nCatAlg}(\sig{\Sigma}))_{n\in\omega},(\Tr^{(m,n)}_{\sig{\Sigma}})_{m\leq n})$ and $((\mathsf{nCat}^{\htg}\mathsf{Alg}(\sig{\Sigma}))_{n\in\omega},(\Tr^{\htg(m,n)}_{\sig{\Sigma}})_{m\leq n})$.
By Remark~\ref{RDsTrCompOmega}, the truncation functors from dimension $\omega$ form cones over them.
\end{remark}

It remains to prove that these cones satisfy the universal property.

\begin{proposition}\label{PProjLimit}\lean{weaker}{Limits.homogeneousProjectiveLimit,Limits.heterogeneousProjectiveLimit}
Let $\mathsf{L}$ be a category.
Then the following statements hold:
\begin{enumerate}
\item for every family $(F_{n}\colon\mathsf{L}\mor\mathsf{nCatAlg}(\sig{\Sigma}))_{n\in\omega}$ of functors such that $F_{m}=\Tr^{(m,n)}_{\sig{\Sigma}}\circ F_{n}$ for all $m\leq n$ in $\omega$, there exists a unique functor $\langle F_{k}\rangle_{k\in\omega}\colon\mathsf{L}\mor\omega\mathsf{CatAlg}(\sig{\Sigma})$ such that $F_{n}=\Tr^{(n,\omega)}_{\sig{\Sigma}}\circ\langle F_{k}\rangle_{k\in\omega}$ for every $n\in\omega$; and
\item for every family $(F_{n}\colon\mathsf{L}\mor\mathsf{nCat}^{\htg}\mathsf{Alg}(\sig{\Sigma}))_{n\in\omega}$ of functors such that $F_{m}=\Tr^{\htg(m,n)}_{\sig{\Sigma}}\circ F_{n}$ for all $m\leq n$ in $\omega$, there exists a unique functor $\langle F_{k}\rangle_{k\in\omega}\colon\mathsf{L}\mor\omega\mathsf{Cat}^{\htg}\mathsf{Alg}(\sig{\Sigma})$ such that $F_{n}=\Tr^{\htg(n,\omega)}_{\sig{\Sigma}}\circ\langle F_{k}\rangle_{k\in\omega}$ for every $n\in\omega$.
\end{enumerate}
\end{proposition}
\begin{proof}

(1) We first fix some terminology.
We say that a family $(\calg{A}_{n})_{n\in\omega}$, with $\calg{A}_{n}$ an $n$-categorial $\sig{\Sigma}$-algebra, is \emph{compatible} if $\Tr^{(m,n)}_{\sig{\Sigma}}(\calg{A}_{n})=\calg{A}_{m}$ for all $m\leq n$ in $\omega$, i.e., if, for all $m\leq n$ in $\omega$, the cells of $\calg{A}_{m}$ are the $m$-cells of $\calg{A}_{n}$ and its operations are the restrictions of those of $\calg{A}_{n}$.
Similarly, a family of homomorphisms $(H_{n}\colon\calg{A}_{n}\mor\calg{A}'_{n})_{n\in\omega}$ between compatible families is \emph{compatible} if $\Tr^{(m,n)}_{\sig{\Sigma}}(H_{n})=H_{m}$ for all $m\leq n$ in $\omega$, i.e., if $H_{m}$ is the restriction of $H_{n}$ to the cells of $\calg{A}_{m}$.

In order to define the functor $\langle F_{k}\rangle_{k\in\omega}$, we prove the following two claims.

\begin{claim}\label{ClProjLimitA}
For every compatible family $(\calg{A}_{n})_{n\in\omega}$ there exists a unique $\omega$-categorial $\sig{\Sigma}$-algebra $\varprojlim_{n}\calg{A}_{n}$ such that, for every $n\in\omega$, $\Tr^{(n,\omega)}_{\sig{\Sigma}}(\varprojlim_{m}\calg{A}_{m})=\calg{A}_{n}$.
\end{claim}

Let $(\calg{A}_{n})_{n\in\omega}$ be a compatible family.
For every $s\in S$, by compatibility, the cells of sort $s$ of the $\calg{A}_{n}$ form an increasing chain of sets, since, for $m\leq n$, the cells of $\calg{A}_{m}$ are the $m$-cells of $\calg{A}_{n}$.
We let $\varprojlim_{n}\calg{A}_{n}$ be the partial $\Sigma^{(\omega)}$-algebra whose carrier of sort $s$ is the union of the carriers of sort $s$ of the $\calg{A}_{n}$, for $n\in\omega$, in which, for every $(\mathbf{s},s)\in S^{\star}\times S$ and every $\sigma\in\Sigma_{\mathbf{s},s}$, $\sigma$ acts on a family of cells as $\sigma^{\calg{A}_{n}}$, for any $n\in\omega$ such that $\calg{A}_{n}$ contains them, and in which, for every $k\in\omega$ and every $s\in S$, $\src^{k}_{s}$, $\tgt^{k}_{s}$ and $\circ^{k}_{s}$ act on cells as $\src^{k\,\calg{A}_{n}}_{s}$, $\tgt^{k\,\calg{A}_{n}}_{s}$ and $\circ^{k\,\calg{A}_{n}}_{s}$, for any $n\in\omega$ with $n>k$ such that $\calg{A}_{n}$ contains them.

We first prove that this definition does not depend on $n$.
If $n\leq n'$, then, by compatibility, the operations of $\Sigma$ and the operations of dimension $k<n$ of $\calg{A}_{n}$ are the restrictions of those of $\calg{A}_{n'}$ to the cells of $\calg{A}_{n}$, with the restricted domains, see the proof of Proposition~\ref{PTrEpVariety}.

Next, we prove that $\varprojlim_{n}\calg{A}_{n}$ is a $\sig{\Sigma}^{(\omega)}$-algebra.
Every operation symbol of $\Sigma^{(\omega)}$ has finitely many arguments and its domain conditions involve at most one dimension.
Hence, for every family of cells of $\varprojlim_{n}\calg{A}_{n}$, there exists $n\in\omega$ such that $\calg{A}_{n}$ contains the cells and the dimension, and the operation is defined at the family in $\varprojlim_{n}\calg{A}_{n}$ if, and only if, it is defined in $\calg{A}_{n}$, i.e., since $\calg{A}_{n}$ is a $\sig{\Sigma}^{(n)}$-algebra, if, and only if, the family satisfies the domain conditions in $\calg{A}_{n}$, or, equivalently, in $\varprojlim_{n}\calg{A}_{n}$.

Now we prove that the axioms hold.
Every axiom of $\spec{C}^{\sig{\Sigma}^{(\omega)}}$ involves finitely many variables and dimensions.
Given a valuation of its variables in $\varprojlim_{n}\calg{A}_{n}$, choose $n\in\omega$ greater than these dimensions such that $\calg{A}_{n}$ contains the values of the variables.
Since the operations of $\calg{A}_{n}$ take their values in $\calg{A}_{n}$, the terms of the axiom have the same definedness and values in $\calg{A}_{n}$ and in $\varprojlim_{n}\calg{A}_{n}$, so that the axiom holds in $\varprojlim_{n}\calg{A}_{n}$, because it holds in $\calg{A}_{n}$.

Next, we prove that every cell has a finite dimension.
Let $x$ be a cell of sort $s$ of $\calg{A}_{n}$.
Since $\calg{A}_{n}=\Tr^{(n,n+1)}_{\sig{\Sigma}}(\calg{A}_{n+1})$, $x$ is an $n$-cell of $\calg{A}_{n+1}$, i.e., $\src^{n\,\calg{A}_{n+1}}_{s}(x)=x$, so that $x$ is an $n$-cell of $\varprojlim_{m}\calg{A}_{m}$.
Therefore, by Remark~\ref{RCatAlgAsEpAlgOmega}, $\varprojlim_{n}\calg{A}_{n}$ is an $\omega$-categorial $\sig{\Sigma}$-algebra.

Now we compute the truncations.
Let $n\in\omega$.
By the above, every cell of $\calg{A}_{n}$ is an $n$-cell of the limit; conversely, an $n$-cell of the limit belongs to some $\calg{A}_{p}$, with $p>n$, and it is an $n$-cell of $\calg{A}_{p}$, i.e., a cell of $\Tr^{(n,p)}_{\sig{\Sigma}}(\calg{A}_{p})=\calg{A}_{n}$.
Moreover, by the description of $\Tr^{(n,\omega)}_{\sig{\Sigma}}$ after Definition~\ref{DTrFunctorOmega}, the operations of $\Tr^{(n,\omega)}_{\sig{\Sigma}}(\varprojlim_{m}\calg{A}_{m})$ are the restrictions of those of the limit, i.e., those of $\calg{A}_{n}$.
Hence $\Tr^{(n,\omega)}_{\sig{\Sigma}}(\varprojlim_{m}\calg{A}_{m})=\calg{A}_{n}$.

Finally, we prove the uniqueness.
Every $\omega$-categorial $\sig{\Sigma}$-algebra $\calg{B}$ is the limit of its truncations, since every cell of $\calg{B}$ is an $n$-cell for some $n\in\omega$ and the operations of $\Tr^{(n,\omega)}_{\sig{\Sigma}}(\calg{B})$ are the restrictions of those of $\calg{B}$.
Hence, if $\Tr^{(n,\omega)}_{\sig{\Sigma}}(\calg{B})=\calg{A}_{n}$ for every $n\in\omega$, then $\calg{B}=\varprojlim_{n}\calg{A}_{n}$.
Claim~\ref{ClProjLimitA} follows.

\begin{claim}\label{ClProjLimitB}
For every compatible family $(H_{n}\colon\calg{A}_{n}\mor\calg{A}'_{n})_{n\in\omega}$ of homomorphisms there exists a unique $\omega$-categorial $\sig{\Sigma}$-homomorphism $\varprojlim_{n}H_{n}$ from $\varprojlim_{n}\calg{A}_{n}$ to $\varprojlim_{n}\calg{A}'_{n}$ such that, for every $n\in\omega$, $\Tr^{(n,\omega)}_{\sig{\Sigma}}(\varprojlim_{m}H_{m})=H_{n}$.
\end{claim}

We let $\varprojlim_{n}H_{n}$ be the $S$-sorted mapping that acts on the cells of $\calg{A}_{n}$ as $H_{n}$, for every $n\in\omega$, which is well defined by compatibility.
By Remark~\ref{RCatAlgAsEpAlgOmega}, it suffices to prove that it is a homomorphism of partial $\Sigma^{(\omega)}$-algebras.
Every operation of $\varprojlim_{n}\calg{A}_{n}$ at a family of cells is computed in some $\calg{A}_{n}$, and the corresponding operation of $\varprojlim_{n}\calg{A}'_{n}$ at the image family is computed in $\calg{A}'_{n}$, where $H_{n}$ is a homomorphism; hence $\varprojlim_{n}H_{n}$ is a homomorphism.
By construction, $\Tr^{(n,\omega)}_{\sig{\Sigma}}(\varprojlim_{m}H_{m})=H_{n}$, for every $n\in\omega$, and $\varprojlim_{n}H_{n}$ is the unique mapping with this property, since every cell of the limit belongs to some $\calg{A}_{n}$.
Claim~\ref{ClProjLimitB} follows.

Since $F_{m}=\Tr^{(m,n)}_{\sig{\Sigma}}\circ F_{n}$ for all $m\leq n$ in $\omega$, for every object $X$ of $\mathsf{L}$, the family $(F_{n}(X))_{n\in\omega}$ is compatible, and, for every morphism $f$ of $\mathsf{L}$, so is $(F_{n}(f))_{n\in\omega}$.
We let
\begin{enumerate}
	\item
	$\langle F_{k}\rangle_{k\in\omega}(X)=\varprojlim_{n}F_{n}(X)$; and

	\item
	$\langle F_{k}\rangle_{k\in\omega}(f)=\varprojlim_{n}F_{n}(f)$.
\end{enumerate}

We first prove that $\langle F_{k}\rangle_{k\in\omega}$ is a functor.
By the uniqueness in Claim~\ref{ClProjLimitB}, $\langle F_{k}\rangle_{k\in\omega}$ preserves identities and composites, since every $F_{n}$ does.

Now we prove that it factors the functors $F_{n}$.
By Claims~\ref{ClProjLimitA} and~\ref{ClProjLimitB}, $\Tr^{(n,\omega)}_{\sig{\Sigma}}\circ\langle F_{k}\rangle_{k\in\omega}=F_{n}$, for every $n\in\omega$.

Finally, we prove its uniqueness.
If $G\colon\mathsf{L}\mor\omega\mathsf{CatAlg}(\sig{\Sigma})$ also satisfies $\Tr^{(n,\omega)}_{\sig{\Sigma}}\circ G=F_{n}$ for every $n\in\omega$, then, by the uniqueness in Claims~\ref{ClProjLimitA} and~\ref{ClProjLimitB}, $G(X)=\varprojlim_{n}F_{n}(X)$ and $G(f)=\varprojlim_{n}F_{n}(f)$, so that $G=\langle F_{k}\rangle_{k\in\omega}$.

(2) We say that a family $(\calg{A}_{n})_{n\in\omega}$, with $\calg{A}_{n}$ a heterogeneous $n$-categorial $\sig{\Sigma}$-algebra, is \emph{compatible} if $\Tr^{\htg(m,n)}_{\sig{\Sigma}}(\calg{A}_{n})=\calg{A}_{m}$ for all $m\leq n$ in $\omega$, i.e., if, for all $m\leq n$ in $\omega$, $\calg{A}_{m}$ is the reduct of $\calg{A}_{n}$ to the sorts of dimension at most $m$ and to the operation symbols of $\Sigma^{\htg(m)}$, and similarly for families of homomorphisms.
We prove the following two claims.

\begin{claim}\label{ClProjLimitHtgA}
For every compatible family $(\calg{A}_{n})_{n\in\omega}$ there exists a unique heterogeneous $\omega$-categorial $\sig{\Sigma}$-algebra $\varprojlim_{n}\calg{A}_{n}$ such that, for every $n\in\omega$, $\Tr^{\htg(n,\omega)}_{\sig{\Sigma}}(\varprojlim_{m}\calg{A}_{m})=\calg{A}_{n}$.
\end{claim}
Let $(\calg{A}_{n})_{n\in\omega}$ be a compatible family.
We let $\varprojlim_{n}\calg{A}_{n}$ be the partial $\Sigma^{\htg(\omega)}$-algebra whose carrier of sort $(k,s)$, for every $k\in\omega$ and every $s\in S$, is the carrier of sort $(k,s)$ of $\calg{A}_{k}$, which, by compatibility, is the carrier of sort $(k,s)$ of $\calg{A}_{n}$, for every $n\geq k$, and in which every operation symbol of $\Sigma^{\htg(\omega)}$ whose dimensions are at most $n$ acts as in $\calg{A}_{n}$.
By compatibility, this definition does not depend on $n$.
For instance, for every $j<k$ in $\omega$ and every $s\in S$, the $(j,k)$-source of a $k$-cell $x$ of sort $s$ is $\src^{(j,k)\,\calg{A}_{k}}_{s}(x)$, which is an element of the carrier of sort $(j,s)$ of $\calg{A}_{k}$, i.e., of $\calg{A}_{j}$; and, for every $(\mathbf{s},s)\in S^{\star}\times S$ and every $\sigma\in\Sigma_{\mathbf{s},s}$, $\sigma_{k}$ acts as $\sigma^{\calg{A}_{k}}_{k}$.
Every operation symbol of $\Sigma^{\htg(\omega)}$ and every axiom of $\spec{C}^{\sig{\Sigma}^{\htg(\omega)}}$ involve only finitely many dimensions, so that, as in~(1), every operation is defined exactly at the families that satisfy its domain conditions and every axiom holds, since this is the case in $\calg{A}_{n}$, for $n$ large enough.
Hence $\varprojlim_{n}\calg{A}_{n}$ is a heterogeneous $\omega$-categorial $\sig{\Sigma}$-algebra.
By construction and the description of $\Tr^{\htg(n,\omega)}_{\sig{\Sigma}}$ after Definition~\ref{DTrFunctorOmega}, $\Tr^{\htg(n,\omega)}_{\sig{\Sigma}}(\varprojlim_{m}\calg{A}_{m})=\calg{A}_{n}$, for every $n\in\omega$.
Finally, every heterogeneous $\omega$-categorial $\sig{\Sigma}$-algebra $\calg{B}$ is the limit of its truncations, i.e., $\calg{B}=\varprojlim_{n}\Tr^{\htg(n,\omega)}_{\sig{\Sigma}}(\calg{B})$, which gives the uniqueness.
Claim~\ref{ClProjLimitHtgA} follows.

\begin{claim}\label{ClProjLimitHtgB}
For every compatible family $(H_{n}\colon\calg{A}_{n}\mor\calg{A}'_{n})_{n\in\omega}$ of homomorphisms there exists a unique heterogeneous $\omega$-categorial $\sig{\Sigma}$-homomorphism $\varprojlim_{n}H_{n}$ such that, for every $n\in\omega$, $\Tr^{\htg(n,\omega)}_{\sig{\Sigma}}(\varprojlim_{m}H_{m})=H_{n}$.
\end{claim}
We let $\varprojlim_{n}H_{n}$ be the family whose component at $(k,s)$, for every $k\in\omega$ and every $s\in S$, is the component of $H_{k}$ at $(k,s)$, which, by compatibility, is the component of $H_{n}$ at $(k,s)$, for every $n\geq k$.
It is a homomorphism of partial $\Sigma^{\htg(\omega)}$-algebras, since every operation symbol of $\Sigma^{\htg(\omega)}$ involves finitely many dimensions, at most $n$, say, and $H_{n}$ is a homomorphism; by Remark~\ref{RCatAlgAsEpAlgHtg}, it is a heterogeneous $\omega$-categorial $\sig{\Sigma}$-homomorphism.
By construction, its truncations are the $H_{n}$, and it is the unique family with this property.
Claim~\ref{ClProjLimitHtgB} follows.

Then, as in~(1), the functor $\langle F_{k}\rangle_{k\in\omega}$ is defined by $\langle F_{k}\rangle_{k\in\omega}(X)=\varprojlim_{n}F_{n}(X)$ and $\langle F_{k}\rangle_{k\in\omega}(f)=\varprojlim_{n}F_{n}(f)$, and, by Claims~\ref{ClProjLimitHtgA} and~\ref{ClProjLimitHtgB}, it is the unique functor such that $\Tr^{\htg(n,\omega)}_{\sig{\Sigma}}\circ\langle F_{k}\rangle_{k\in\omega}=F_{n}$, for every $n\in\omega$. This completes the proof.
\end{proof}

We use the universal property in the following form.

\begin{corollary}\label{CProjLimitNat}\lean{weaker}{Limits.homogeneousHomEquiv,Limits.homogeneousNatTransEquiv}
Let $\mathsf{L}$ be a category.
Then the following statements hold:
\begin{enumerate}
\item for all $\omega$-categorial $\sig{\Sigma}$-algebras $\calg{A}$ and $\calg{B}$, the mapping that sends every $\omega$-categorial $\sig{\Sigma}$-homomorphism $H$ from $\calg{A}$ to $\calg{B}$ to the family $(\Tr^{(n,\omega)}_{\sig{\Sigma}}(H))_{n\in\omega}$ is a bijection from the set of $\omega$-categorial $\sig{\Sigma}$-homomorphisms from $\calg{A}$ to $\calg{B}$ to the set of families $(H_{n}\colon \Tr^{(n,\omega)}_{\sig{\Sigma}}(\calg{A})\mor \Tr^{(n,\omega)}_{\sig{\Sigma}}(\calg{B}))_{n\in\omega}$ such that $\Tr^{(m,n)}_{\sig{\Sigma}}(H_{n})=H_{m}$ for all $m\leq n$ in $\omega$; and
\item for all functors $F,G\colon\mathsf{L}\mor\omega\mathsf{CatAlg}(\sig{\Sigma})$ and every family $(\varphi^{n})_{n\in\omega}$, with $\varphi^{n}$ a natural transformation from $\Tr^{(n,\omega)}_{\sig{\Sigma}}\circ F$ to $\Tr^{(n,\omega)}_{\sig{\Sigma}}\circ G$, such that $\Tr^{(m,n)}_{\sig{\Sigma}}\ast\varphi^{n}=\varphi^{m}$ for all $m\leq n$ in $\omega$, there exists a unique natural transformation $\varphi$ from $F$ to $G$ such that $\Tr^{(n,\omega)}_{\sig{\Sigma}}\ast\varphi=\varphi^{n}$ for every $n\in\omega$, and $\varphi$ is a natural isomorphism if every $\varphi^{n}$ is.
\end{enumerate}
\end{corollary}
\begin{proof}
(1) By the uniqueness part of the proof of Claim~\ref{ClProjLimitA}, see Proposition~\ref{PProjLimit}, $\calg{A}$ and $\calg{B}$ are the limits of their truncations, and the families of truncations are compatible, by Remark~\ref{RDsTrCompOmega}~(1).
Hence the statement follows from Claim~\ref{ClProjLimitB}.

(2) By~(1), for every object $X$ of $\mathsf{L}$ there exists a unique homomorphism $\varphi_{X}$ from $F(X)$ to $G(X)$ such that $\Tr^{(n,\omega)}_{\sig{\Sigma}}(\varphi_{X})=\varphi^{n}_{X}$ for every $n\in\omega$.

Now we prove the naturality.
For a morphism $f\colon X\mor Y$, the homomorphisms $G(f)\circ\varphi_{X}$ and $\varphi_{Y}\circ F(f)$ have the same truncations, by the naturality of every $\varphi^{n}$, so that, by~(1), they are equal.

Finally, we prove that $\varphi$ is a natural isomorphism if every $\varphi^{n}$ is.
If every $\varphi^{n}_{X}$ is an isomorphism, then the inverses $(\varphi^{n}_{X})^{-1}$ form a family as in~(1), whose image under the inverse of the bijection of~(1) is the inverse of $\varphi_{X}$. This completes the proof.
\end{proof}

The same holds in the heterogeneous presentation.

\begin{corollary}\label{CProjLimitNatHtg}\lean{weaker}{Limits.heterogeneousHomEquiv,Limits.heterogeneousNatTransEquiv}
Let $\mathsf{L}$ be a category.
Then the following statements hold:
\begin{enumerate}
\item for all heterogeneous $\omega$-categorial $\sig{\Sigma}$-algebras $\calg{A}$ and $\calg{B}$, the mapping that sends every heterogeneous $\omega$-categorial $\sig{\Sigma}$-homomorphism $H$ from $\calg{A}$ to $\calg{B}$ to the family $(\Tr^{\htg(n,\omega)}_{\sig{\Sigma}}(H))_{n\in\omega}$ is a bijection from the set of heterogeneous $\omega$-categorial $\sig{\Sigma}$-homomorphisms from $\calg{A}$ to $\calg{B}$ to the set of families $(H_{n}\colon \Tr^{\htg(n,\omega)}_{\sig{\Sigma}}(\calg{A})\mor \Tr^{\htg(n,\omega)}_{\sig{\Sigma}}(\calg{B}))_{n\in\omega}$ such that $\Tr^{\htg(m,n)}_{\sig{\Sigma}}(H_{n})=H_{m}$ for all $m\leq n$ in $\omega$; and
\item for all functors $F,G\colon\mathsf{L}\mor\omega\mathsf{Cat}^{\htg}\mathsf{Alg}(\sig{\Sigma})$ and every family $(\varphi^{n})_{n\in\omega}$, with $\varphi^{n}$ a natural transformation from $\Tr^{\htg(n,\omega)}_{\sig{\Sigma}}\circ F$ to $\Tr^{\htg(n,\omega)}_{\sig{\Sigma}}\circ G$, such that $\Tr^{\htg(m,n)}_{\sig{\Sigma}}\ast\varphi^{n}=\varphi^{m}$ for all $m\leq n$ in $\omega$, there exists a unique natural transformation $\varphi$ from $F$ to $G$ such that $\Tr^{\htg(n,\omega)}_{\sig{\Sigma}}\ast\varphi=\varphi^{n}$ for every $n\in\omega$, and $\varphi$ is a natural isomorphism if every $\varphi^{n}$ is.
\end{enumerate}
\end{corollary}
\begin{proof}
Statement~(1) follows, as in the proof of Corollary~\ref{CProjLimitNat}~(1), from Claims~\ref{ClProjLimitHtgA} and~\ref{ClProjLimitHtgB} in the proof of Proposition~\ref{PProjLimit} and Remark~\ref{RDsTrCompOmega}~(1), and statement~(2) follows from~(1) exactly as in the proof of Corollary~\ref{CProjLimitNat}~(2).
\end{proof}

\chapter{The functors in dimension \texorpdfstring{$\omega$}{omega}}\label{SFunctorsOmega}

Throughout this chapter, $\sig{\Sigma}=(S,\Sigma)$ is a many-sorted signature, which is fixed in all the definitions and results of this chapter.
In this chapter we obtain the $\omega$-dimensional versions of the functors of Chapters~\ref{SGdUf} and~\ref{SDsTr}.
The truncation and discretization functors in dimension $\omega$ have been defined in Section~\ref{SDsTrOmega}, and, by Chapter~\ref{SProjLimits}, the truncation functors are the projections of projective limits.
We begin by obtaining the graduation functor, by the universal property of Proposition~\ref{PProjLimit}~(2), from a compatible family of the finite-dimensional ones.
Next we obtain the unification functor by means of an inductive limit of sets of cells.
Next we prove that graduation and unification are mutually quasi-inverse equivalences and that discretization is left adjoint to truncation also in dimension $\omega$.
We end this chapter by showing that these adjunctions form an algebraic square between a finite dimension and $\omega$.

\section{Graduation and unification in dimension \texorpdfstring{$\omega$}{omega}}\label{SGdUfOmega}

In the following lemma we prove that the finite-dimensional graduation functors give rise to a compatible family.

\begin{lemma}\label{LCompatGdTrOmega}\lean{weaker}{Omega.graduationFamily}
For all $m,n\in\mathbb{N}$ with $m\leq n$, the following equality holds
$$
\Tr^{\htg(m,n)}_{\sig{\Sigma}}\circ\Gd^{(n)}_{\sig{\Sigma}}\circ \Tr^{(n,\omega)}_{\sig{\Sigma}}=\Gd^{(m)}_{\sig{\Sigma}}\circ \Tr^{(m,\omega)}_{\sig{\Sigma}}.
$$
\end{lemma}
\begin{proof}
The following chain of equalities holds
\begin{align*}
\Tr^{\htg(m,n)}_{\sig{\Sigma}}\circ\Gd^{(n)}_{\sig{\Sigma}}\circ \Tr^{(n,\omega)}_{\sig{\Sigma}}
&=
\Gd^{(m)}_{\sig{\Sigma}}\circ \Tr^{(m,n)}_{\sig{\Sigma}}\circ \Tr^{(n,\omega)}_{\sig{\Sigma}}
\tag{1}
\\&=
\Gd^{(m)}_{\sig{\Sigma}}\circ \Tr^{(m,\omega)}_{\sig{\Sigma}}.
\tag{2}
\end{align*}
The equality~(1) follows from Proposition~\ref{PCompatGdTr}, and the equality~(2) follows from Remark~\ref{RDsTrCompOmega}~(1).
\end{proof}

\begin{remark}\label{RCompatGdTrOmegaDerivors}\lean{weaker}{Omega.graduation_compatibility_derivor,Omega.graduation_compatibility_ne,Omega.graduation_compatibility_sameCells}
The equality of Lemma~\ref{LCompatGdTrOmega} does not, in general, come from an equality in $\mathsf{HSig}_{\gd}$.
Its two sides are induced by the generalized hep-derivors $\bTr^{(n,\omega)}_{\sig{\Sigma}}\circ\bGd^{(n)}_{\sig{\Sigma}}\circ\bTr^{\htg(m,n)}_{\sig{\Sigma}}$ and $\bTr^{(m,\omega)}_{\sig{\Sigma}}\circ\bGd^{(m)}_{\sig{\Sigma}}$ from $\sig{\Sigma}^{\htg(m)}$ to $\sig{\Sigma}^{(\omega)}$, which have the same underlying derivor: both mappings on sorts send $(k,s)$ to $s$, and, since the truncation derivors are inclusions, both send every operation symbol of $\Sigma^{\htg(m)}$ to the term assigned to it in Definition~\ref{DGdDerivor}, which does not depend on the dimension of the graduation derivor.
However, if $m<n$, then, by Definitions~\ref{DGenHepDerivorCompId}, \ref{DGdDerivor}, \ref{DTrDerivor} and~\ref{DTrDerivorOmega}, for every $s\in S$, the domain of sort $(k,s)$ of the first one is $\{(v^s_0,\src^n_s(v^s_0)),(v^s_0,\src^k_s(v^s_0))\}$, for every $k\leq m$, whereas that of the second one is $\{(v^s_0,\src^m_s(v^s_0)),(v^s_0,\src^k_s(v^s_0))\}$, if $k<m$, and $\{(v^s_0,\src^m_s(v^s_0))\}$, if $k=m$.
Thus both generalized hep-derivors are equal if $m=n$, but not if $m<n$.
On every $\omega$-categorial $\sig{\Sigma}$-algebra both domains of sort $(k,s)$ select the $k$-cells of sort $s$, since every $k$-cell is an $m$-cell and every $m$-cell is an $n$-cell, by~(AB$^{k,m}$1) and~(AB$^{m,n}$1), see Remark~\ref{RSameCellsPattern}.
\end{remark}

We are now in a position to define the graduation functor in dimension $\omega$.

\begin{definition}\label{DGdFunctorOmega}\lean{weaker}{Omega.graduationFunctor,Omega.graduation_unique_induced}
The \emph{graduation functor}
$$
\Gd^{(\omega)}_{\sig{\Sigma}}\colon\omega\mathsf{CatAlg}(\sig{\Sigma})\mor\omega\mathsf{Cat}^{\htg}\mathsf{Alg}(\sig{\Sigma})
$$
is the unique functor such that, for every $n\in\mathbb{N}$,
$$
\Tr^{\htg(n,\omega)}_{\sig{\Sigma}}\circ\Gd^{(\omega)}_{\sig{\Sigma}}=\Gd^{(n)}_{\sig{\Sigma}}\circ \Tr^{(n,\omega)}_{\sig{\Sigma}}.
$$
\end{definition}

Let us note that the functor exists and is unique by Proposition~\ref{PProjLimit}~(2), since, by Lemma~\ref{LCompatGdTrOmega}, the family $(\Gd^{(n)}_{\sig{\Sigma}}\circ \Tr^{(n,\omega)}_{\sig{\Sigma}})_{n\in\omega}$ is compatible.
Explicitly, by Claim~\ref{ClProjLimitHtgA} in the proof of Proposition~\ref{PProjLimit}, the cells of sort $(k,s)$ of $\Gd^{(\omega)}_{\sig{\Sigma}}(\calg{A})$ are the $k$-cells of sort $s$ of $\calg{A}$, for every $j<k$ in $\omega$, the $(j,k)$-sources, $(j,k)$-targets and $j$-compositions of $k$-cells are the restrictions of $\src^{j\,\calg{A}}_{s}$, $\tgt^{j\,\calg{A}}_{s}$ and $\circ^{j\,\calg{A}}_{s}$, the $(k,j)$-identities are inclusions, and, for every $(\mathbf{s},s)\in S^{\star}\times S$ and every $\sigma\in\Sigma_{\mathbf{s},s}$, $\sigma_{k}$ acts on $k$-cells as $\sigma^{\calg{A}}$.

We next consider unification.
In contrast with graduation, it cannot be obtained in the same way.

\begin{remark}\label{RUfNotDerivorOmega}\lean{stronger}{Omega.unification_not_induced}
If $S\neq\varnothing$, the unification functor cannot be obtained as the functor induced by a generalized hep-derivor.
In fact, let $\mathbf{c}=(\varphi,c,\Dom^{(\varphi,c)})$ be a generalized hep-derivor from $\sig{\Sigma}^{(\omega)}$ to $\sig{\Sigma}^{\htg(\omega)}$; then, for every heterogeneous $\omega$-categorial $\sig{\Sigma}$-algebra $\calg{A}$ and every $s\in S$, the cells of sort $s$ of $\mathbf{c}^{\ast}_{\gd}(\calg{A})$ are elements of the single carrier $A_{\varphi(s)}$ of $\calg{A}$.
Let $s\in S$ and $\varphi(s)=(k,s')$, and let $\mathsf{M}$ be the heterogeneous $\omega$-category with exactly one cell in every dimension up to $k$, whose $(k+1)$-cells form the two-element monoid $\{e,f\}$, with $f\circ f=f$, for all the compositions, and whose cells of dimension greater than $k+1$ are the identities of the $(k+1)$-cells.
Let $\calg{A}$ be the heterogeneous $\omega$-categorial $\sig{\Sigma}$-algebra in which the heterogeneous $\omega$-category of every sort is $\mathsf{M}$ and every operation symbol of $\Sigma$ is interpreted as the constant heterogeneous $\omega$-functor whose value in every dimension is the identity cell of that dimension, i.e., the unique cell up to dimension $k$ and $e$ and its identities above it; since this cell is an identity in every dimension and is idempotent for all the compositions, this constant family is a heterogeneous $\omega$-functor.
Then $\mathbf{c}^{\ast}_{\gd}(\calg{A})$ has at most one cell of sort $s$, since $A_{(k,s')}$ has exactly one element, whereas the unification of $\calg{A}$ has two distinct cells of sort $s$, which correspond to $e$ and $f$.
If $S=\varnothing$, then all the algebras involved are empty, and the unification functor is trivially induced by the unique generalized hep-derivor from $\sig{\Sigma}^{(\omega)}$ to $\sig{\Sigma}^{\htg(\omega)}$.

Nor can we apply directly the universal property of Proposition~\ref{PProjLimit} to the family $(\Uf^{(n)}_{\sig{\Sigma}}\circ \Tr^{\htg(n,\omega)}_{\sig{\Sigma}})_{n\in\omega}$, since, by Proposition~\ref{PCompatUfTr}, it is only compatible up to the natural isomorphisms $\beta^{(m,n)}$: the $m$-cells of $\Uf^{(n)}_{\sig{\Sigma}}(\calg{B})$ are not the $m$-cells of $\calg{B}$, but their $(n,m)$-identities.
These isomorphisms, however, are the inclusions of an inductive system of sets of cells, and we construct the unification of a heterogeneous $\omega$-categorial $\sig{\Sigma}$-algebra on the inductive limit of these sets, as for $\omega$-categories.
\end{remark}

We begin by defining the $\omega$-cells of a heterogeneous $\omega$-categorial $\sig{\Sigma}$-algebra.

\begin{definition}\label{DUfCellsOmega}\lean{equivalent}{Omega.Cell,Omega.cellSetoid,Omega.cellIn_eq_iff_comparable}
Let $\calg{B}$ be a heterogeneous $\omega$-categorial $\sig{\Sigma}$-algebra and $s\in S$.
The set of \emph{$\omega$-cells} of sort $s$ of $\calg{B}$ is $\coprod_{k\in\omega}B_{(k,s)}$, whose elements are the pairs $(x,k)$ with $k\in\omega$ and $x\in B_{(k,s)}$.
We let $\Theta_{s}$ be the binary relation on $\coprod_{k\in\omega}B_{(k,s)}$ defined by $((x,j),(y,k))\in\Theta_{s}$ if
\begin{enumerate}
\item $j<k$ and $y=\mathrm{i}^{(k,j)\,\calg{B}}_{s}(x)$; or
\item $j=k$ and $x=y$; or
\item $j>k$ and $x=\mathrm{i}^{(j,k)\,\calg{B}}_{s}(y)$.
\end{enumerate}
By~(EA$^{h,j,k}$2) and Lemma~\ref{LUfIdInj}, $\Theta_{s}$ is an equivalence relation; the injectivity of the identities is used for the transitivity.
We let $[B]_{s}$ stand for the quotient set $\coprod_{k\in\omega}B_{(k,s)}/\Theta_{s}$, $[x]_{k}$ for the class of $(x,k)$, and $\mathrm{in}^{[k]}_{s}$ for the mapping from $B_{(k,s)}$ to $[B]_{s}$ that sends $x$ to $[x]_{k}$.
\end{definition}

In the following lemma we collect the properties of the $\omega$-cells that we need.

\begin{lemma}\label{LUfCellsOmega}\lean{equivalent}{Omega.cell_inductive_limit,Omega.cellIn_injective,Omega.cell_cover}
Let $\calg{B}$ be a heterogeneous $\omega$-categorial $\sig{\Sigma}$-algebra and $s\in S$.
Then the following statements hold:
\begin{enumerate}
\item $[B]_{s}$, together with the mappings $\mathrm{in}^{[k]}_{s}$, is the inductive limit of the inductive system of sets $(B_{(k,s)},\mathrm{i}^{(k+1,k)\,\calg{B}}_{s})_{k\in\omega}$, i.e., for every set $L$ and every family of mappings $f_{k}\colon B_{(k,s)}\mor L$ such that $f_{j}=f_{k}\circ\mathrm{i}^{(k,j)\,\calg{B}}_{s}$ for all $j<k$ in $\omega$, there exists a unique mapping $[f_{k}]_{k\in\omega}\colon[B]_{s}\mor L$ such that $[f_{k}]_{k\in\omega}\circ\mathrm{in}^{[j]}_{s}=f_{j}$ for every $j\in\omega$;
\item every $\mathrm{in}^{[k]}_{s}$ is injective; and
\item $[B]_{s}$ is the increasing union of the images of the $\mathrm{in}^{[k]}_{s}$; in particular, every finite family of $\omega$-cell classes has representatives in a common $B_{(p,s)}$.
\end{enumerate}
\end{lemma}
\begin{proof}
(1) It is the usual construction of the inductive limit of an inductive system of sets, since, by~(EA$^{h,j,k}$2), $\mathrm{i}^{(k,j)\,\calg{B}}_{s}$ is the composite of the mappings $\mathrm{i}^{(l+1,l)\,\calg{B}}_{s}$, for $j\leq l<k$, and $\Theta_{s}$ identifies $(x,j)$ with $(\mathrm{i}^{(k,j)\,\calg{B}}_{s}(x),k)$, for $j<k$.

(2) By Definition~\ref{DUfCellsOmega}~(2), two elements $x$, $y$ of $B_{(k,s)}$ satisfy $((x,k),(y,k))\in\Theta_{s}$ if, and only if, $x=y$, so that $[x]_{k}=[y]_{k}$ if, and only if, $x=y$.

(3) By Definition~\ref{DUfCellsOmega}~(1), $[x]_{j}=[\mathrm{i}^{(k,j)\,\calg{B}}_{s}(x)]_{k}$ for all $j<k$ in $\omega$ and every $x\in B_{(j,s)}$, so that the image of $\mathrm{in}^{[j]}_{s}$ is contained in that of $\mathrm{in}^{[k]}_{s}$; and every $\omega$-cell class is $[x]_{k}$ for some $k\in\omega$.
\end{proof}

Let $\calg{B}$ be a heterogeneous $\omega$-categorial $\sig{\Sigma}$-algebra.
In the formulas below, $\mathrm{i}^{(k,k)\,\calg{B}}_s$ denotes the identity of $B_{(k,s)}$.
For every $n\in\omega$, $\Uf^{(n)}_{\sig{\Sigma}}(\Tr^{\htg(n,\omega)}_{\sig{\Sigma}}(\calg{B}))$ is an $n$-categorial $\sig{\Sigma}$-algebra whose cells of sort $s$ are the elements of $B_{(n,s)}$.
By Proposition~\ref{PCompatUfTr} and Remark~\ref{RDsTrCompOmega}~(1), for all $n\leq n'$ in $\omega$, $\mathrm{i}^{(n',n)}=(\mathrm{i}^{(n',n)\,\calg{B}}_{s})_{s\in S}$ is an isomorphism
$$
\mathrm{i}^{(n',n)}\colon\Uf^{(n)}_{\sig{\Sigma}}\left(\Tr^{\htg(n,\omega)}_{\sig{\Sigma}}(\calg{B})\right)\mor\Tr^{(n,n')}_{\sig{\Sigma}}\left(\Uf^{(n')}_{\sig{\Sigma}}\left(\Tr^{\htg(n',\omega)}_{\sig{\Sigma}}(\calg{B})\right)\right);
$$
that is, $\Uf^{(n)}_{\sig{\Sigma}}(\Tr^{\htg(n,\omega)}_{\sig{\Sigma}}(\calg{B}))$ is, up to $\mathrm{i}^{(n',n)}$, the algebra of the $n$-cells of $\Uf^{(n')}_{\sig{\Sigma}}(\Tr^{\htg(n',\omega)}_{\sig{\Sigma}}(\calg{B}))$.
We next define the unification of $\calg{B}$.

\begin{definition}\label{DUfFunctorOmega}\lean{equivalent}{Omega.unifiedAlgebra}
Let $\calg{B}$ be a heterogeneous $\omega$-categorial $\sig{\Sigma}$-algebra.
The \emph{unification} of $\calg{B}$ is the partial $\Sigma^{(\omega)}$-algebra $\Uf^{(\omega)}_{\sig{\Sigma}}(\calg{B})$ with underlying $S$-sorted set $[B]=([B]_{s})_{s\in S}$ in which, for every operation symbol $\tau$ of $\Sigma^{(\omega)}$, of arity $\mathbf{s}$ and coarity $s$, and every family $([x_{j}]_{p})_{j\in\bb{\mathbf{s}}}$ of $\omega$-cell classes represented in a common level $p$ such that $\tau$ is an operation symbol of $\Sigma^{(p)}$, $([x_{j}]_{p})_{j\in\bb{\mathbf{s}}}\in\Dom(\tau^{\Uf^{(\omega)}_{\sig{\Sigma}}(\calg{B})})$ if, and only if, $(x_{j})_{j\in\bb{\mathbf{s}}}\in\Dom(\tau^{\Uf^{(p)}_{\sig{\Sigma}}(\Tr^{\htg(p,\omega)}_{\sig{\Sigma}}(\calg{B}))})$, and then
$$
\tau^{\Uf^{(\omega)}_{\sig{\Sigma}}(\calg{B})}(([x_{j}]_{p})_{j\in\bb{\mathbf{s}}})=[\tau^{\Uf^{(p)}_{\sig{\Sigma}}(\Tr^{\htg(p,\omega)}_{\sig{\Sigma}}(\calg{B}))}((x_{j})_{j\in\bb{\mathbf{s}}})]_{p}.
$$
For a heterogeneous $\omega$-categorial $\sig{\Sigma}$-homomorphism $F\colon\calg{B}\mor\calg{B}'$, $\Uf^{(\omega)}_{\sig{\Sigma}}(F)$ is the $S$-sorted mapping from $[B]$ to $[B']$ that sends $[x]_{k}$ to $[F_{(k,s)}(x)]_{k}$.
\end{definition}

In the following lemma we prove that the above definition is correct.

\begin{lemma}\label{LUfWellDefOmega}\lean{equivalent}{Omega.unifiedHom}
Let $\calg{B}$ be a heterogeneous $\omega$-categorial $\sig{\Sigma}$-algebra and $F$ a heterogeneous $\omega$-categorial $\sig{\Sigma}$-homomorphism.
Then $\Uf^{(\omega)}_{\sig{\Sigma}}(\calg{B})$ and $\Uf^{(\omega)}_{\sig{\Sigma}}(F)$ are well defined.
\end{lemma}
\begin{proof}
We first prove that the definition does not depend on the level.
Let $p\leq p'$.
Then, for every $s\in S$, $\mathrm{in}^{[p]}_{s}=\mathrm{in}^{[p']}_{s}\circ\mathrm{i}^{(p',p)\,\calg{B}}_{s}$, and, by Proposition~\ref{PCompatUfTr}, $\mathrm{i}^{(p',p)}$ is an isomorphism from $\Uf^{(p)}_{\sig{\Sigma}}(\Tr^{\htg(p,\omega)}_{\sig{\Sigma}}(\calg{B}))$ onto the algebra of the $p$-cells of $\Uf^{(p')}_{\sig{\Sigma}}(\Tr^{\htg(p',\omega)}_{\sig{\Sigma}}(\calg{B}))$, whose operations are the restrictions of those of $\Uf^{(p')}_{\sig{\Sigma}}(\Tr^{\htg(p',\omega)}_{\sig{\Sigma}}(\calg{B}))$.
Hence computing at the level $p$ or at the level $p'$ gives the same domain and the same values.

Now we prove that $\Uf^{(\omega)}_{\sig{\Sigma}}(F)$ is well defined.
Let $F$ be a heterogeneous $\omega$-categorial $\sig{\Sigma}$-homomorphism from $\calg{B}$ to $\calg{B}'$ and $s\in S$.
Since $F$ commutes with the identities, $F_{(k,s)}\circ\mathrm{i}^{(k,j)\,\calg{B}}_{s}=\mathrm{i}^{(k,j)\,\calg{B}'}_{s}\circ F_{(j,s)}$ for all $j<k$ in $\omega$, so that the family of mappings $(\mathrm{in}^{[k]}_{s}\circ F_{(k,s)})_{k\in\omega}$ satisfies the compatibility condition of Lemma~\ref{LUfCellsOmega}~(1), and $\Uf^{(\omega)}_{\sig{\Sigma}}(F)$ is well defined. This completes the proof.
\end{proof}

In other words, the underlying $S$-sorted set of $\Uf^{(\omega)}_{\sig{\Sigma}}(\calg{B})$ is the inductive limit of the inductive system of $S$-sorted sets underlying $(\Uf^{(n)}_{\sig{\Sigma}}(\Tr^{\htg(n,\omega)}_{\sig{\Sigma}}(\calg{B})),\mathrm{i}^{(n+1,n)})_{n\in\omega}$, see Lemma~\ref{LUfCellsOmega}, and its operations are computed at a common finite level; by Claim~\ref{ClUfFunctorOmegaC} below, for every $n\in\omega$, its algebra of $n$-cells is isomorphic to $\Uf^{(n)}_{\sig{\Sigma}}(\Tr^{\htg(n,\omega)}_{\sig{\Sigma}}(\calg{B}))$.
Since these algebras live in categories of different dimensions, we do not claim a universal property of $\Uf^{(\omega)}_{\sig{\Sigma}}(\calg{B})$ in a category of categorial algebras; the term inductive limit refers only to the sets of cells.
In the following proposition we prove that it is an $\omega$-categorial $\sig{\Sigma}$-algebra.

\begin{proposition}\label{PUfFunctorOmega}\lean{equivalent}{Omega.unified_isModel,Omega.unified_isCell_iff,Omega.unificationFunctor}
Let $\calg{B}$ be a heterogeneous $\omega$-categorial $\sig{\Sigma}$-algebra.
Then $\Uf^{(\omega)}_{\sig{\Sigma}}(\calg{B})$ is an $\omega$-categorial $\sig{\Sigma}$-algebra whose $k$-cells of sort $s$ are the classes $[x]_{k}$, with $x\in B_{(k,s)}$, and, for every $n\in\omega$, $\mathrm{in}^{[n]}=(\mathrm{in}^{[n]}_{s})_{s\in S}$ is an isomorphism, natural in $\calg{B}$, from $\Uf^{(n)}_{\sig{\Sigma}}(\Tr^{\htg(n,\omega)}_{\sig{\Sigma}}(\calg{B}))$ onto $\Tr^{(n,\omega)}_{\sig{\Sigma}}(\Uf^{(\omega)}_{\sig{\Sigma}}(\calg{B}))$.
Moreover, $\Uf^{(\omega)}_{\sig{\Sigma}}$ is a functor from $\omega\mathsf{Cat}^{\htg}\mathsf{Alg}(\sig{\Sigma})$ to $\omega\mathsf{CatAlg}(\sig{\Sigma})$, the \emph{unification functor}.
\end{proposition}
\begin{proof}
We prove the following claims.

\begin{claim}\label{ClUfFunctorOmegaA}
$\Uf^{(\omega)}_{\sig{\Sigma}}(\calg{B})$ is a $\sig{\Sigma}^{(\omega)}$-algebra in which every axiom of $\spec{C}^{\sig{\Sigma}^{(\omega)}}$ holds.
\end{claim}
Every operation symbol $\tau$ of $\Sigma^{(\omega)}$ has finitely many arguments and its domain conditions involve at most one dimension.
Hence, by Lemma~\ref{LUfCellsOmega}~(3), every family of $\omega$-cell classes is represented at a level $p$ such that $\tau$ is an operation symbol of $\Sigma^{(p)}$, and, by Definition~\ref{DUfFunctorOmega} and Lemma~\ref{LUfWellDefOmega}, $\tau^{\Uf^{(\omega)}_{\sig{\Sigma}}(\calg{B})}$ is defined at the family if, and only if, $\tau$ is defined at the representatives in the $\sig{\Sigma}^{(p)}$-algebra $\Uf^{(p)}_{\sig{\Sigma}}(\Tr^{\htg(p,\omega)}_{\sig{\Sigma}}(\calg{B}))$, i.e., if, and only if, the representatives satisfy the domain conditions of $\tau$ there, or, equivalently, the family satisfies them in $\Uf^{(\omega)}_{\sig{\Sigma}}(\calg{B})$.
Thus $\Uf^{(\omega)}_{\sig{\Sigma}}(\calg{B})$ is a $\sig{\Sigma}^{(\omega)}$-algebra.
Every axiom of $\spec{C}^{\sig{\Sigma}^{(\omega)}}$ involves finitely many dimensions and finitely many variables.
Given a valuation of its variables in $[B]$, by Lemma~\ref{LUfCellsOmega}~(3), we choose a level $p$ greater than all these dimensions at which all the values of the variables are represented.
Since $\mathrm{in}^{[p]}$ is injective and, by Definition~\ref{DUfFunctorOmega}, it preserves and reflects the definedness and the values of the operations of $\Sigma^{(p)}$, the axiom holds in $\Uf^{(\omega)}_{\sig{\Sigma}}(\calg{B})$ under the given valuation, because it is an axiom of $\spec{C}^{\sig{\Sigma}^{(p)}}$, which holds in the $p$-categorial $\sig{\Sigma}$-algebra $\Uf^{(p)}_{\sig{\Sigma}}(\Tr^{\htg(p,\omega)}_{\sig{\Sigma}}(\calg{B}))$.

Claim~\ref{ClUfFunctorOmegaA} follows.

\begin{claim}\label{ClUfFunctorOmegaB}
The $k$-cells of sort $s$ of $\Uf^{(\omega)}_{\sig{\Sigma}}(\calg{B})$ are the classes $[x]_{k}$, with $x\in B_{(k,s)}$.
In particular, $\Uf^{(\omega)}_{\sig{\Sigma}}(\calg{B})$ is an $\omega$-categorial $\sig{\Sigma}$-algebra.
\end{claim}

Let $s\in S$, $x\in B_{(k,s)}$ and $p>k$.
By Lemma~\ref{LUfDescription}, the $k$-source of sort $s$ of $\Uf^{(p)}_{\sig{\Sigma}}(\Tr^{\htg(p,\omega)}_{\sig{\Sigma}}(\calg{B}))$ acts as $\mathrm{i}^{(p,k)\,\calg{B}}_{s}\circ\src^{(k,p)\,\calg{B}}_{s}$, and, by Lemma~\ref{LUfIdInj}, the following equality holds
$$
\mathrm{i}^{(p,k)\,\calg{B}}_{s}\left(\src^{(k,p)\,\calg{B}}_{s}\left(\mathrm{i}^{(p,k)\,\calg{B}}_{s}(x)\right)\right)
=
\mathrm{i}^{(p,k)\,\calg{B}}_{s}(x).
$$
Thus, by Definition~\ref{DUfFunctorOmega}, $[x]_{k}=[\mathrm{i}^{(p,k)\,\calg{B}}_{s}(x)]_{p}$ is a $k$-cell.
Conversely, if $[y]_{p}$, with $y\in B_{(p,s)}$ and $p>k$, is a $k$-cell, then, by Definition~\ref{DUfFunctorOmega}, Lemma~\ref{LUfCellsOmega}~(2) and Lemma~\ref{LUfDescription}, $y=\mathrm{i}^{(p,k)\,\calg{B}}_{s}(\src^{(k,p)\,\calg{B}}_{s}(y))$, so that $[y]_{p}=[\src^{(k,p)\,\calg{B}}_{s}(y)]_{k}$; and, if $p\leq k$, then $[y]_{p}=[\mathrm{i}^{(k,p)\,\calg{B}}_{s}(y)]_{k}$.
Hence every $\omega$-cell class is a $k$-cell for some $k$, and, by Claim~\ref{ClUfFunctorOmegaA} and Remark~\ref{RCatAlgAsEpAlgOmega}, $\Uf^{(\omega)}_{\sig{\Sigma}}(\calg{B})$ is an $\omega$-categorial $\sig{\Sigma}$-algebra.
Claim~\ref{ClUfFunctorOmegaB} follows.

\begin{claim}\label{ClUfFunctorOmegaC}
For every $n\in\omega$, $\mathrm{in}^{[n]}$ is an isomorphism from $\Uf^{(n)}_{\sig{\Sigma}}(\Tr^{\htg(n,\omega)}_{\sig{\Sigma}}(\calg{B}))$ onto $\Tr^{(n,\omega)}_{\sig{\Sigma}}(\Uf^{(\omega)}_{\sig{\Sigma}}(\calg{B}))$, natural in $\calg{B}$.
\end{claim}

By Claim~\ref{ClUfFunctorOmegaB} and Lemma~\ref{LUfCellsOmega}~(2), $\mathrm{in}^{[n]}$ is a bijection from the cells of $\Uf^{(n)}_{\sig{\Sigma}}(\Tr^{\htg(n,\omega)}_{\sig{\Sigma}}(\calg{B}))$ onto the $n$-cells of $\Uf^{(\omega)}_{\sig{\Sigma}}(\calg{B})$, i.e., onto the cells of $\Tr^{(n,\omega)}_{\sig{\Sigma}}(\Uf^{(\omega)}_{\sig{\Sigma}}(\calg{B}))$.
By Definition~\ref{DUfFunctorOmega}, it preserves and reflects the definedness and the values of the operations of $\Sigma^{(n)}$, so that it is an isomorphism.
It is natural in $\calg{B}$, since, by Definition~\ref{DUfFunctorOmega}, for every heterogeneous $\omega$-categorial $\sig{\Sigma}$-homomorphism $F$ from $\calg{B}$ to $\calg{B}'$, $\Uf^{(\omega)}_{\sig{\Sigma}}(F)$ sends $[x]_{n}$ to $[F_{(n,s)}(x)]_{n}$, for every $s\in S$ and every $x\in B_{(n,s)}$.
Claim~\ref{ClUfFunctorOmegaC} follows.

\begin{claim}\label{ClUfFunctorOmegaD}
$\Uf^{(\omega)}_{\sig{\Sigma}}$ is a functor.
\end{claim}

Let $F\colon\calg{B}\mor\calg{B}'$ be a heterogeneous $\omega$-categorial $\sig{\Sigma}$-homomorphism.
For every level $p$, the following equality holds
$$
\Uf^{(\omega)}_{\sig{\Sigma}}(F)\circ\mathrm{in}^{[p]}=\mathrm{in}^{[p]}\circ\Uf^{(p)}_{\sig{\Sigma}}\left(\Tr^{\htg(p,\omega)}_{\sig{\Sigma}}(F)\right),
$$
so that, since $\Uf^{(p)}_{\sig{\Sigma}}(\Tr^{\htg(p,\omega)}_{\sig{\Sigma}}(F))$ is a homomorphism, $\Uf^{(\omega)}_{\sig{\Sigma}}(F)$ commutes with the operations of $\Sigma^{(p)}$.
Thus it is a homomorphism.
Moreover, $\Uf^{(\omega)}_{\sig{\Sigma}}$ preserves identities and composites, since it acts componentwise on representatives.
Claim~\ref{ClUfFunctorOmegaD} follows. This completes the proof.
\end{proof}

Explicitly, for every $s\in S$, every $p\in\omega$ and every $x,y\in B_{(p,s)}$, if $k<p$, then the $k$-source and the $k$-target of $[y]_{p}$ are $[\src^{(k,p)\,\calg{B}}_{s}(y)]_{k}$ and $[\tgt^{(k,p)\,\calg{B}}_{s}(y)]_{k}$, and the $k$-composition of $[x]_{p}$ and $[y]_{p}$ is $[x\circ^{(k,p)\,\calg{B}}_{s}y]_{p}$, defined if, and only if, $\src^{(k,p)\,\calg{B}}_{s}(x)=\tgt^{(k,p)\,\calg{B}}_{s}(y)$; if $k\geq p$, then $[y]_{p}$ is a $k$-cell, so that its $k$-source and $k$-target are $[y]_{p}$ itself; and, for every $(\mathbf{s},s)\in S^{\star}\times S$ and every $\sigma\in\Sigma_{\mathbf{s},s}$, $\sigma$ acts on classes represented at the level $p$ as $\sigma^{\calg{B}}_{p}$.

\section{Equivalence and adjunctions in dimension \texorpdfstring{$\omega$}{omega}}\label{SGdUfDsTrRelOmega}

In this section we prove that the $\omega$-dimensional functors are related as the finite-dimensional ones.
We begin with graduation and unification.

\begin{proposition}\label{PGdUfIsoOmega}\lean{equivalent}{Omega.zeta,Omega.theta}
For every $\omega$-categorial $\sig{\Sigma}$-algebra $\calg{A}$, let $\zeta_{\calg{A}}$ be the family whose component at the sort $s\in S$ sends a class $[x]_{k}$ of $\Uf^{(\omega)}_{\sig{\Sigma}}(\Gd^{(\omega)}_{\sig{\Sigma}}(\calg{A}))$ to $x$; and, for every heterogeneous $\omega$-categorial $\sig{\Sigma}$-algebra $\calg{B}$, let $\theta_{\calg{B}}$ be the family $(\mathrm{in}^{[k]}_{s})_{(k,s)\in\omega\times S}$, which sends $x\in B_{(k,s)}$ to $[x]_{k}$, see Definition~\ref{DUfCellsOmega}.
Then $\zeta=(\zeta_{\calg{A}})_{\calg{A}}$ is a natural isomorphism from $\Uf^{(\omega)}_{\sig{\Sigma}}\circ\Gd^{(\omega)}_{\sig{\Sigma}}$ to $\Id_{\omega\mathsf{CatAlg}(\sig{\Sigma})}$ and $\theta=(\theta_{\calg{B}})_{\calg{B}}$ is a natural isomorphism from $\Id_{\omega\mathsf{Cat}^{\htg}\mathsf{Alg}(\sig{\Sigma})}$ to $\Gd^{(\omega)}_{\sig{\Sigma}}\circ\Uf^{(\omega)}_{\sig{\Sigma}}$.
In particular, the following natural isomorphisms hold
\begin{enumerate}
	\item
	$\Uf^{(\omega)}_{\sig{\Sigma}}\circ\Gd^{(\omega)}_{\sig{\Sigma}}\cong\Id_{\omega\mathsf{CatAlg}(\sig{\Sigma})}$; and

	\item
	$\Gd^{(\omega)}_{\sig{\Sigma}}\circ\Uf^{(\omega)}_{\sig{\Sigma}}\cong\Id_{\omega\mathsf{Cat}^{\htg}\mathsf{Alg}(\sig{\Sigma})}$.
\end{enumerate}
\end{proposition}
\begin{proof}
We prove the following claims.

\begin{claim}\label{ClGdUfIsoOmegaA}
$\Uf^{(\omega)}_{\sig{\Sigma}}\circ\Gd^{(\omega)}_{\sig{\Sigma}}\cong\Id_{\omega\mathsf{CatAlg}(\sig{\Sigma})}$.
\end{claim}

Let $\calg{A}$ be an $\omega$-categorial $\sig{\Sigma}$-algebra.
Since the identities of $\Gd^{(\omega)}_{\sig{\Sigma}}(\calg{A})$ are inclusions, two $\omega$-cells $(x,j)$ and $(y,k)$ of $\Gd^{(\omega)}_{\sig{\Sigma}}(\calg{A})$ are $\Theta_{s}$-related if, and only if, $x=y$.
Hence the mapping $\zeta_{\calg{A},s}$ that sends $[x]_{k}$ to $x$ is a bijection from the set of the cells of sort $s$ of $\Uf^{(\omega)}_{\sig{\Sigma}}(\Gd^{(\omega)}_{\sig{\Sigma}}(\calg{A}))$ onto the set of the cells of sort $s$ of $\calg{A}$, all of which are $k$-cells for some $k$.
For every level $p$, the following chain of equalities holds
\begin{align*}
\Uf^{(p)}_{\sig{\Sigma}}\left(\Tr^{\htg(p,\omega)}_{\sig{\Sigma}}\left(\Gd^{(\omega)}_{\sig{\Sigma}}(\calg{A})\right)\right)
&=
\Uf^{(p)}_{\sig{\Sigma}}\left(\Gd^{(p)}_{\sig{\Sigma}}\left(\Tr^{(p,\omega)}_{\sig{\Sigma}}(\calg{A})\right)\right)
\tag{1}
\\&=
\Tr^{(p,\omega)}_{\sig{\Sigma}}(\calg{A}).
\tag{2}
\end{align*}
The equality~(1) follows from Definition~\ref{DGdFunctorOmega}, and the equality~(2) follows from Corollary~\ref{CUfGdId}.
Therefore, by Definition~\ref{DUfFunctorOmega}, $\zeta_{\calg{A}}=(\zeta_{\calg{A},s})_{s\in S}$ preserves and reflects the definedness and the values of the operations of $\Sigma^{(p)}$.
Hence it is an isomorphism from $\Uf^{(\omega)}_{\sig{\Sigma}}(\Gd^{(\omega)}_{\sig{\Sigma}}(\calg{A}))$ onto $\calg{A}$, and $\zeta=(\zeta_{\calg{A}})_{\calg{A}\in\mathrm{Ob}(\omega\mathsf{CatAlg}(\sig{\Sigma}))}$ is a natural transformation from $\Uf^{(\omega)}_{\sig{\Sigma}}\circ\Gd^{(\omega)}_{\sig{\Sigma}}$ to $\Id_{\omega\mathsf{CatAlg}(\sig{\Sigma})}$, since $\Gd^{(\omega)}_{\sig{\Sigma}}$ and $\Uf^{(\omega)}_{\sig{\Sigma}}$ act on homomorphisms by restriction and on representatives, respectively.
Claim~\ref{ClGdUfIsoOmegaA} follows.

\begin{claim}\label{ClGdUfIsoOmegaB}
$\Gd^{(\omega)}_{\sig{\Sigma}}\circ\Uf^{(\omega)}_{\sig{\Sigma}}\cong\Id_{\omega\mathsf{Cat}^{\htg}\mathsf{Alg}(\sig{\Sigma})}$.
\end{claim}

Let $\calg{B}$ be a heterogeneous $\omega$-categorial $\sig{\Sigma}$-algebra.
By Proposition~\ref{PUfFunctorOmega}, the cells of sort $(k,s)$ of $\Gd^{(\omega)}_{\sig{\Sigma}}(\Uf^{(\omega)}_{\sig{\Sigma}}(\calg{B}))$ are the classes $[x]_{k}$, with $x\in B_{(k,s)}$, so that the family $\theta_{\calg{B}}=(\mathrm{in}^{[k]}_{s})_{(k,s)\in\omega\times S}$ is a bijection from the cells of $\calg{B}$ onto those of $\Gd^{(\omega)}_{\sig{\Sigma}}(\Uf^{(\omega)}_{\sig{\Sigma}}(\calg{B}))$.
It is a homomorphism that reflects definedness.
In fact, for every $s\in S$, all $j<k$ in $\omega$ and all $x,y\in B_{(k,s)}$,
\begin{enumerate}
\item the $(j,k)$-source of $[x]_{k}$ is its $j$-source, which, by the explicit description after Proposition~\ref{PUfFunctorOmega}, is $[\src^{(j,k)\,\calg{B}}_{s}(x)]_{j}$;
\item the $(j,k)$-target of $[x]_{k}$ is its $j$-target, which, by the same description, is $[\tgt^{(j,k)\,\calg{B}}_{s}(x)]_{j}$;
\item the $(k,j)$-identity of $[u]_{j}$, with $u\in B_{(j,s)}$, is $[u]_{j}=[\mathrm{i}^{(k,j)\,\calg{B}}_{s}(u)]_{k}$;
\item the $j$-composition of $[x]_{k}$ and $[y]_{k}$ is $[x\circ^{(j,k)\,\calg{B}}_{s}y]_{k}$, defined if, and only if, $x\circ^{(j,k)\,\calg{B}}_{s}y$ is defined; and
\item for every $(\mathbf{s},s)\in S^{\star}\times S$ and every $\sigma\in\Sigma_{\mathbf{s},s}$, $\sigma_{k}$ acts on $k$-cells as $\sigma$, i.e., as $\sigma^{\calg{B}}_{k}$ on representatives.
\end{enumerate}
Hence $\theta_{\calg{B}}$ is an isomorphism from $\calg{B}$ onto $\Gd^{(\omega)}_{\sig{\Sigma}}(\Uf^{(\omega)}_{\sig{\Sigma}}(\calg{B}))$ and, by Definition~\ref{DUfFunctorOmega}, $\theta=(\theta_{\calg{B}})_{\calg{B}\in\mathrm{Ob}(\omega\mathsf{Cat}^{\htg}\mathsf{Alg}(\sig{\Sigma}))}$ is a natural isomorphism from $\Id_{\omega\mathsf{Cat}^{\htg}\mathsf{Alg}(\sig{\Sigma})}$ to $\Gd^{(\omega)}_{\sig{\Sigma}}\circ\Uf^{(\omega)}_{\sig{\Sigma}}$.
Claim~\ref{ClGdUfIsoOmegaB} follows.

All in all, the statement follows from Claims~\ref{ClGdUfIsoOmegaA} and~\ref{ClGdUfIsoOmegaB}.
\end{proof}

From the above proposition we obtain that graduation and unification are equivalences also in dimension $\omega$.

\begin{corollary}\label{CGdUfEquivOmega}\lean{equivalent}{Omega.presentationEquivalence}
The functors $\Gd^{(\omega)}_{\sig{\Sigma}}$ and $\Uf^{(\omega)}_{\sig{\Sigma}}$ are mutually quasi-inverse equivalences of categories.
More precisely, let $\zeta$ and $\theta$ be the natural isomorphisms of Proposition~\ref{PGdUfIsoOmega}.
Then $(\Uf^{(\omega)}_{\sig{\Sigma}},\Gd^{(\omega)}_{\sig{\Sigma}},\zeta^{-1},\theta^{-1})$ is an adjoint equivalence, i.e., $\Gd^{(\omega)}_{\sig{\Sigma}}$ is left adjoint to $\Uf^{(\omega)}_{\sig{\Sigma}}$ with unit $\zeta^{-1}$ and counit $\theta^{-1}$.
\end{corollary}
\begin{proof}
The first statement follows from Proposition~\ref{PGdUfIsoOmega}.
For the second one, since $\zeta^{-1}$ and $\theta^{-1}$ are natural isomorphisms, it suffices to prove the triangular identities of Definition~\ref{DAdj}.
Let us recall that $\zeta_{\calg{A}}$ sends $[x]_{k}$ to $x$, so that $\zeta^{-1}_{\calg{A}}$ sends a cell $x$ of $\calg{A}$ to $[x]_{k}$, for any $k$ such that $x$ is a $k$-cell, and that $\theta_{\calg{B}}$ is $\mathrm{in}^{[k]}_{s}$ at the sort $(k,s)$, see Proposition~\ref{PGdUfIsoOmega}, so that $\theta^{-1}_{\calg{B}}$ sends $[x]_{k}$, with $x\in B_{(k,s)}$, regarded as a cell of sort $(k,s)$, to $x$.

(TI1) Let $\calg{B}$ be a heterogeneous $\omega$-categorial $\sig{\Sigma}$-algebra and $[x]_{k}$ a cell of sort $s$ of $\Uf^{(\omega)}_{\sig{\Sigma}}(\calg{B})$, with $x\in B_{(k,s)}$.
By Proposition~\ref{PUfFunctorOmega}, $[x]_{k}$ is a $k$-cell, so that $\zeta^{-1}_{\Uf^{(\omega)}_{\sig{\Sigma}}(\calg{B})}([x]_{k})$ is the class, at level $k$, of $[x]_{k}$ regarded as a cell of sort $(k,s)$ of $\Gd^{(\omega)}_{\sig{\Sigma}}(\Uf^{(\omega)}_{\sig{\Sigma}}(\calg{B}))$.
By Definition~\ref{DUfFunctorOmega}, $\Uf^{(\omega)}_{\sig{\Sigma}}(\theta^{-1}_{\calg{B}})$ sends this class to $[\theta^{-1}_{\calg{B}}([x]_{k})]_{k}=[x]_{k}$.
Hence the following equality holds
$$
(\Uf^{(\omega)}_{\sig{\Sigma}}\ast\theta^{-1})\circ(\zeta^{-1}\ast\Uf^{(\omega)}_{\sig{\Sigma}})=\id_{\Uf^{(\omega)}_{\sig{\Sigma}}}.
$$

(TI2) Let $\calg{A}$ be an $\omega$-categorial $\sig{\Sigma}$-algebra and $x$ a cell of sort $(k,s)$ of $\Gd^{(\omega)}_{\sig{\Sigma}}(\calg{A})$, i.e., a $k$-cell of sort $s$ of $\calg{A}$.
Since $\Gd^{(\omega)}_{\sig{\Sigma}}$ acts on homomorphisms by restriction, $\Gd^{(\omega)}_{\sig{\Sigma}}(\zeta^{-1}_{\calg{A}})$ sends $x$ to $[x]_{k}$, which is a $k$-cell of $\Uf^{(\omega)}_{\sig{\Sigma}}(\Gd^{(\omega)}_{\sig{\Sigma}}(\calg{A}))$, i.e., a cell of sort $(k,s)$ of $\Gd^{(\omega)}_{\sig{\Sigma}}(\Uf^{(\omega)}_{\sig{\Sigma}}(\Gd^{(\omega)}_{\sig{\Sigma}}(\calg{A})))$, and $\theta^{-1}_{\Gd^{(\omega)}_{\sig{\Sigma}}(\calg{A})}$ sends it back to $x$.
Hence the following equality holds
$$
(\theta^{-1}\ast\Gd^{(\omega)}_{\sig{\Sigma}})\circ(\Gd^{(\omega)}_{\sig{\Sigma}}\ast\zeta^{-1})=\id_{\Gd^{(\omega)}_{\sig{\Sigma}}}.
$$
This completes the proof.
\end{proof}

We next consider discretization and truncation.

\begin{proposition}\label{PDsTrAdjOmega}\lean{weaker}{Omega.homogeneousRoundTripFunctor,Omega.homogeneousAdjunction,Omega.heterogeneousAdjunction}
Let $m\in\mathbb{N}$.
Then $\Tr^{(m,\omega)}_{\sig{\Sigma}}\circ\Ds^{(\omega,m)}_{\sig{\Sigma}}=\Id_{\mathsf{mCatAlg}(\sig{\Sigma})}$ and $\Ds^{(\omega,m)}_{\sig{\Sigma}}$ is left adjoint to $\Tr^{(m,\omega)}_{\sig{\Sigma}}$, with the identity as unit and, as counit, the natural transformation $\varepsilon$ whose component at an $\omega$-categorial $\sig{\Sigma}$-algebra $\calg{B}$ is the inclusion of the $m$-cells.
The same holds for $\Ds^{\htg(\omega,m)}_{\sig{\Sigma}}$ and $\Tr^{\htg(m,\omega)}_{\sig{\Sigma}}$, the component of the counit at a heterogeneous $\omega$-categorial $\sig{\Sigma}$-algebra $\calg{B}$ and at the sort $(k,s)$ being the identity of $B_{(k,s)}$, if $k\leq m$, and $\mathrm{i}^{(k,m)\,\calg{B}}_{s}$, if $k>m$.
\end{proposition}
\begin{proof}
We first prove the equation of the unit.
By Remark~\ref{RDsTrCompOmega}~(2), for $n=m$, $\Tr^{(m,\omega)}_{\sig{\Sigma}}\circ\Ds^{(\omega,m)}_{\sig{\Sigma}}$ and $\Tr^{\htg(m,\omega)}_{\sig{\Sigma}}\circ\Ds^{\htg(\omega,m)}_{\sig{\Sigma}}$ are identities.

Now we define the counit.
Let $\calg{B}$ be an $\omega$-categorial $\sig{\Sigma}$-algebra and $n\in\mathbb{N}$.
By Remark~\ref{RDsTrCompOmega}, $\Tr^{(n,\omega)}_{\sig{\Sigma}}(\Ds^{(\omega,m)}_{\sig{\Sigma}}(\Tr^{(m,\omega)}_{\sig{\Sigma}}(\calg{B})))$ is $\Ds^{(n,m)}_{\sig{\Sigma}}(\Tr^{(m,n)}_{\sig{\Sigma}}(\Tr^{(n,\omega)}_{\sig{\Sigma}}(\calg{B})))$, if $m\leq n$, and $\Tr^{(n,\omega)}_{\sig{\Sigma}}(\calg{B})$, if $n<m$.
Let $\varepsilon^{n}_{\calg{B}}$ be the component of the counit of Proposition~\ref{PDsTrAdj} at $\Tr^{(n,\omega)}_{\sig{\Sigma}}(\calg{B})$, if $m\leq n$, and the identity, if $n<m$.
All of them are inclusions of $m$-cells, or identities, and, for $n\leq n'$, the restriction of $\varepsilon^{n'}_{\calg{B}}$ to the $n$-cells is $\varepsilon^{n}_{\calg{B}}$.
Hence they satisfy the compatibility condition of Corollary~\ref{CProjLimitNat}~(2).
Therefore they induce a natural transformation $\varepsilon$ from $\Ds^{(\omega,m)}_{\sig{\Sigma}}\circ\Tr^{(m,\omega)}_{\sig{\Sigma}}$ to the identity, whose component at $\calg{B}$ is the inclusion of the $m$-cells of $\calg{B}$.

In the heterogeneous case, the components of the counits of Proposition~\ref{PDsTrAdj} at the sorts of dimension at most $n$ are identities and $(k,m)$-identities, see Lemma~\ref{LDsTrCounit}.
Hence they are compatible, and, by Corollary~\ref{CProjLimitNatHtg}~(2), they induce the natural transformation $\varepsilon^{\htg}$ of the statement.

Finally, we prove the triangular identities.
Since the unit is the identity, the triangular identities state that $\Tr^{(m,\omega)}_{\sig{\Sigma}}(\varepsilon_{\calg{B}})$ and $\varepsilon_{\Ds^{(\omega,m)}_{\sig{\Sigma}}(\calg{A})}$ are identities, for every $\omega$-categorial $\sig{\Sigma}$-algebra $\calg{B}$ and every $m$-categorial $\sig{\Sigma}$-algebra $\calg{A}$, and similarly in the heterogeneous case.
They hold as in Proposition~\ref{PDsTrAdj}.
In fact, $\Tr^{(m,\omega)}_{\sig{\Sigma}}(\varepsilon_{\calg{B}})$ is the restriction of $\varepsilon_{\calg{B}}$ to the $m$-cells, which is an identity.
By Corollary~\ref{CDsCellsOmega}, every cell of $\Ds^{(\omega,m)}_{\sig{\Sigma}}(\calg{A})$ is an $m$-cell, so that $\varepsilon_{\Ds^{(\omega,m)}_{\sig{\Sigma}}(\calg{A})}$ is an identity.
In the heterogeneous case, the components of $\varepsilon^{\htg}_{\calg{B}}$ at the dimensions $k\leq m$ are identities.
And the components of $\varepsilon^{\htg}_{\Ds^{\htg(\omega,m)}_{\sig{\Sigma}}(\calg{A})}$ at the dimensions $k>m$ are the $(k,m)$-identities of $\Ds^{\htg(\omega,m)}_{\sig{\Sigma}}(\calg{A})$, which are identities, by the description after Definition~\ref{DDsFunctorOmega}~(2).
\end{proof}

\section{An algebraic square in dimension \texorpdfstring{$\omega$}{omega}}\label{SAlgSquareOmega}

In this section we show that the algebraic square of Proposition~\ref{PAlgSquareFinite} also exists between a finite dimension and $\omega$.

\begin{proposition}\label{PAlgSquareOmega}\lean{weaker}{Omega.dimensionAlgebraicSquare}
Let $m\in\mathbb{N}$, and let
$$
\begin{aligned}
\mathsf{C}&=\mathsf{mCatAlg}(\sig{\Sigma}), &
\mathsf{D}&=\mathsf{mCat}^{\htg}\mathsf{Alg}(\sig{\Sigma}),\\
\mathsf{C}'&=\omega\mathsf{CatAlg}(\sig{\Sigma}), &
\mathsf{D}'&=\omega\mathsf{Cat}^{\htg}\mathsf{Alg}(\sig{\Sigma}).
\end{aligned}
$$
Then the adjunctions
\begin{enumerate}
\item $\Gd^{(m)}_{\sig{\Sigma}}\dashv\Uf^{(m)}_{\sig{\Sigma}}$ of Proposition~\ref{PGdUfAdj} and the adjoint equivalence $\Gd^{(\omega)}_{\sig{\Sigma}}\dashv\Uf^{(\omega)}_{\sig{\Sigma}}$ of Corollary~\ref{CGdUfEquivOmega}, with unit $\zeta^{-1}$ and counit $\theta^{-1}$; and
\item $\Ds^{(\omega,m)}_{\sig{\Sigma}}\dashv \Tr^{(m,\omega)}_{\sig{\Sigma}}$ and $\Ds^{\htg(\omega,m)}_{\sig{\Sigma}}\dashv \Tr^{\htg(m,\omega)}_{\sig{\Sigma}}$ of Proposition~\ref{PDsTrAdjOmega}
\end{enumerate}
form an algebraic square, with $F=\Gd^{(m)}_{\sig{\Sigma}}$, $G=\Uf^{(m)}_{\sig{\Sigma}}$, $F'=\Gd^{(\omega)}_{\sig{\Sigma}}$, $G'=\Uf^{(\omega)}_{\sig{\Sigma}}$, $J=\Ds^{(\omega,m)}_{\sig{\Sigma}}$, $K=\Tr^{(m,\omega)}_{\sig{\Sigma}}$, $H=\Ds^{\htg(\omega,m)}_{\sig{\Sigma}}$ and $I=\Tr^{\htg(m,\omega)}_{\sig{\Sigma}}$, in which $\alpha$ is the identity of $\Gd^{(\omega)}_{\sig{\Sigma}}\circ\Ds^{(\omega,m)}_{\sig{\Sigma}}=\Ds^{\htg(\omega,m)}_{\sig{\Sigma}}\circ\Gd^{(m)}_{\sig{\Sigma}}$ and $\beta$ is the natural isomorphism $\mathrm{in}^{[m]}$ of Proposition~\ref{PUfFunctorOmega}, from $\Uf^{(m)}_{\sig{\Sigma}}\circ \Tr^{\htg(m,\omega)}_{\sig{\Sigma}}$ to $\Tr^{(m,\omega)}_{\sig{\Sigma}}\circ\Uf^{(\omega)}_{\sig{\Sigma}}$, which sends every $m$-cell $x$ to its class $[x]_{m}$.
\end{proposition}
\begin{proof}
By Definition~\ref{DAlgSquare}, we must check that the four pairs of functors involved are adjunctions and that there exist the natural isomorphisms $\alpha$ and $\beta$.

We first recall the adjunctions.
The four adjunctions are those of Propositions~\ref{PGdUfAdj} and~\ref{PDsTrAdjOmega} and Corollary~\ref{CGdUfEquivOmega}.

Now we prove that $\alpha$ is a natural isomorphism.
For every $p\in\mathbb{N}$, by Definition~\ref{DGdFunctorOmega} and Remark~\ref{RDsTrCompOmega}~(2), the following equalities hold
\begin{gather*}
\Tr^{\htg(p,\omega)}_{\sig{\Sigma}}\circ\Gd^{(\omega)}_{\sig{\Sigma}}\circ\Ds^{(\omega,m)}_{\sig{\Sigma}}=\Gd^{(p)}_{\sig{\Sigma}}\circ F_{p},\\
\Tr^{\htg(p,\omega)}_{\sig{\Sigma}}\circ\Ds^{\htg(\omega,m)}_{\sig{\Sigma}}\circ\Gd^{(m)}_{\sig{\Sigma}}=F^{\htg}_{p}\circ\Gd^{(m)}_{\sig{\Sigma}},
\end{gather*}
with $F_{p}$ and $F^{\htg}_{p}$ as in Remark~\ref{RDsTrCompOmega}.
Moreover, if $m\leq p$, $\Gd^{(p)}_{\sig{\Sigma}}\circ\Ds^{(p,m)}_{\sig{\Sigma}}=\Ds^{\htg(p,m)}_{\sig{\Sigma}}\circ\Gd^{(m)}_{\sig{\Sigma}}$, and, if $p\leq m$, $\Gd^{(p)}_{\sig{\Sigma}}\circ \Tr^{(p,m)}_{\sig{\Sigma}}=\Tr^{\htg(p,m)}_{\sig{\Sigma}}\circ\Gd^{(m)}_{\sig{\Sigma}}$, by Propositions~\ref{PCompatGdDs} and~\ref{PCompatGdTr}.
Hence both functors have the same truncations and, by the uniqueness in Proposition~\ref{PProjLimit}~(2), they are equal, so that we can take $\alpha$ to be the identity.

Finally, we prove that $\beta$ is a natural isomorphism.
By Proposition~\ref{PUfFunctorOmega}, $\mathrm{in}^{[m]}$ is a natural isomorphism from $\Uf^{(m)}_{\sig{\Sigma}}\circ \Tr^{\htg(m,\omega)}_{\sig{\Sigma}}$ to $\Tr^{(m,\omega)}_{\sig{\Sigma}}\circ\Uf^{(\omega)}_{\sig{\Sigma}}$. This completes the proof.
\end{proof}

\begin{remark}\label{RAlgSquareOmegaMate}\lean{equivalent}{Omega.truncationUnificationFunctor_eq_conjugate}
Also in Proposition~\ref{PAlgSquareOmega}, $\beta=\mathrm{in}^{[m]}$ is the natural isomorphism $\beta^{\alpha}$ of Remark~\ref{RAlgSquareMate}.
In fact, $\alpha$ is an identity and the unit of $\Ds^{(\omega,m)}_{\sig{\Sigma}}\dashv\Tr^{(m,\omega)}_{\sig{\Sigma}}$ is the identity, so that the unit of the composite adjunction is $\eta^{\circ}=\Tr^{(m,\omega)}_{\sig{\Sigma}}\ast\zeta^{-1}\ast\Ds^{(\omega,m)}_{\sig{\Sigma}}$, and
$$
\beta^{\alpha}=\left(\Tr^{(m,\omega)}_{\sig{\Sigma}}\circ\Uf^{(\omega)}_{\sig{\Sigma}}\ast\varepsilon^{\circ}\right)\circ\left(\eta^{\circ}\ast\Uf^{(m)}_{\sig{\Sigma}}\circ\Tr^{\htg(m,\omega)}_{\sig{\Sigma}}\right),
$$
where $\varepsilon^{\circ}=\varepsilon^{\htg}\circ(\Ds^{\htg(\omega,m)}_{\sig{\Sigma}}\ast\theta^{-1}\ast\Tr^{\htg(m,\omega)}_{\sig{\Sigma}})$, with $\theta$ the natural isomorphism of Proposition~\ref{PGdUfAdj} in dimension $m$ and $\varepsilon^{\htg}$ the counit of Proposition~\ref{PDsTrAdjOmega}.
Let $\calg{B}$ be a heterogeneous $\omega$-categorial $\sig{\Sigma}$-algebra, $s\in S$ and $x$ a cell of sort $s$ of $\Uf^{(m)}_{\sig{\Sigma}}(\Tr^{\htg(m,\omega)}_{\sig{\Sigma}}(\calg{B}))$, i.e., an element of $B_{(m,s)}$.
By Corollary~\ref{CDsCellsOmega}, $x$ is an $m$-cell of $\Ds^{(\omega,m)}_{\sig{\Sigma}}(\Uf^{(m)}_{\sig{\Sigma}}(\Tr^{\htg(m,\omega)}_{\sig{\Sigma}}(\calg{B})))$, so that, by the proof of Corollary~\ref{CGdUfEquivOmega}, the component of $\eta^{\circ}$ at $\Uf^{(m)}_{\sig{\Sigma}}(\Tr^{\htg(m,\omega)}_{\sig{\Sigma}}(\calg{B}))$ sends $x$ to $[x]_{m}$, the class of $x$ regarded as a cell of sort $(m,s)$ of
\begin{flushleft}
$
\Gd^{(\omega)}_{\sig{\Sigma}}\left(\Ds^{(\omega,m)}_{\sig{\Sigma}}\left(\Uf^{(m)}_{\sig{\Sigma}}\left(\Tr^{\htg(m,\omega)}_{\sig{\Sigma}}(\calg{B})\right)\right)\right)
$
\begin{align*}
&=
\Ds^{\htg(\omega,m)}_{\sig{\Sigma}}\left(\Gd^{(m)}_{\sig{\Sigma}}\left(\Uf^{(m)}_{\sig{\Sigma}}\left(\Tr^{\htg(m,\omega)}_{\sig{\Sigma}}(\calg{B})\right)\right)\right).
\end{align*}
\end{flushleft}
The component of $\varepsilon^{\circ}_{\calg{B}}$ at the sort $(m,s)$ is the identity, since so are the components of $\varepsilon^{\htg}_{\calg{B}}$, by Proposition~\ref{PDsTrAdjOmega}, and of $\Ds^{\htg(\omega,m)}_{\sig{\Sigma}}(\theta^{-1})$ at this sort, the latter being that of $\theta^{-1}$ at the top sort $(m,s)$, see Claim~\ref{ClGdUfIsoC}.
Hence, by Definition~\ref{DUfFunctorOmega}, $\beta^{\alpha}_{\calg{B}}(x)=[x]_{m}=\mathrm{in}^{[m]}_{s}(x)$.
\end{remark}

\chapter{Change of signature}\label{SChangeSig}

Throughout this chapter, unless otherwise stated, $\sig{\Sigma}=(S,\Sigma)$ and $\sig{\Lambda}=(T,\Lambda)$ are many-sorted signatures and $\mathbf{c}=(\varphi,c)$ is a derivor from $\sig{\Sigma}$ to $\sig{\Lambda}$, see Definition~\ref{DDerivor}; they are fixed in the definitions and results of this chapter that do not quantify over the many-sorted signatures.
In this chapter we show that the constructions of this work are compatible with the change of the many-sorted signature.
For every $n\in\mathbb{N}$, the derivor $\mathbf{c}$ determines a generalized hep-derivor $\mathbf{c}^{(n)}$ from $\sig{\Sigma}^{(n)}$ to $\sig{\Lambda}^{(n)}$ and a generalized hep-derivor $\mathbf{c}^{\htg(n)}$ from $\sig{\Sigma}^{\htg(n)}$ to $\sig{\Lambda}^{\htg(n)}$, which translate the operation symbols of $\Sigma$ by means of $c$ and leave the categorial structure unchanged.
We prove that the induced functors send $n$-categorial $\sig{\Lambda}$-algebras to $n$-categorial $\sig{\Sigma}$-algebras, that the assignments $\mathbf{c}\mapsto\mathbf{c}^{(n)}$ and $\mathbf{c}\mapsto\mathbf{c}^{\htg(n)}$ preserve compositions and identities, and that the induced functors commute with graduation, unification, discretization and truncation.
We first consider the finite dimensions and then the dimension $\omega$.

\section{The generalized hep-derivors \texorpdfstring{$\mathbf{c}^{(n)}$ and $\mathbf{c}^{\htg(n)}$}{c(n) and chtg(n)}}\label{SChangeSigDerivors}

We first fix the notation for the terms of $\Lambda$ inside the categorial signatures determined by $\Lambda$.
Since, for every $n\in\mathbb{N}$, $\Lambda\subseteq\Lambda^{(n)}$ and every operation symbol $\tau$ of $\Lambda$ has a copy $\tau_{k}$ in $\Lambda^{\htg(n)}$, for every $k\in n+1$, the terms of $\Lambda$ can be regarded as terms of these signatures by means of the following derivors.

\begin{definition}\label{DInclusionDerivors}\lean{equivalent}{ChangeOfSignature.homogeneousInclusion,ChangeOfSignature.heterogeneousInclusion}
Let $n\in\mathbb{N}$.
\begin{enumerate}
\item The \emph{homogeneous inclusion derivor} is the derivor $\boldsymbol{\iota}^{(n)}=(\id_{T},\iota^{(n)})$ from $(T,\Lambda)$ to $(T,\Lambda^{(n)})$ that sends every $\tau\in\Lambda_{\mathbf{t},t}$ to $\tau((v^{\mathbf{t}}_{p})_{p\in\bb{\mathbf{t}}})$.
\item For every $k\in n+1$, the \emph{heterogeneous inclusion derivor} at dimension $k$ is the derivor $\boldsymbol{\iota}^{\htg(n)}_{k}$ from $(T,\Lambda)$ to $((n+1)\times T,\Lambda^{\htg(n)})$ whose mapping on sorts sends $t\in T$ to $(k,t)$ and which sends every $\tau\in\Lambda_{\mathbf{t},t}$ to $\tau_{k}((v^{(k,\mathbf{t})}_{p})_{p\in\bb{\mathbf{t}}})$, see Definition~\ref{DCatSigHtg}.
For every $\mathbf{t}\in T^{\star}$, every $t\in T$ and every $P\in\mathrm{T}_{\Lambda}(\vs{\mathbf{t}})_{t}$, we let $P_{k}$ stand for the term $(\iota^{\htg(n)}_{k})^{\sharp}(P)\in\mathrm{T}_{\Lambda^{\htg(n)}}(\vs{(k,\mathbf{t})})_{(k,t)}$.
\end{enumerate}
\end{definition}

Let us note that, by~(Der1) and~(Der2), $(\iota^{(n)})^{\sharp}$ sends every term $P$ of $\Lambda$ to $P$, regarded as a term of $\Lambda^{(n)}$, and that $P_{k}$ is obtained from $P$ by replacing every operation symbol $\tau$ by $\tau_{k}$ and every variable $v^{\mathbf{t}}_{p}$ by $v^{(k,\mathbf{t})}_{p}$, i.e.,
\begin{align*}
\left(v^{\mathbf{t}}_{p}\right)_{k}&=v^{(k,\mathbf{t})}_{p},
&
\left(\tau\left((Q_{i})_{i\in\bb{\mathbf{r}}}\right)\right)_{k}&=\tau_{k}\left(\left((Q_{i})_{k}\right)_{i\in\bb{\mathbf{r}}}\right),
\end{align*}
for every $p\in\bb{\mathbf{t}}$, every $\tau\in\Lambda_{\mathbf{r},t}$ and every $(Q_{i})_{i\in\bb{\mathbf{r}}}\in\mathrm{T}_{\Lambda}(\vs{\mathbf{t}})_{\mathbf{r}}$.
Since $\Lambda^{\htg(n)}$ embeds canonically into $\Lambda^{\htg(m)}$, for $n\leq m$, and into $\Lambda^{\htg(\omega)}$, see the note after Definition~\ref{DCatSigHtg}, the term $P_{k}$ does not depend, up to these embeddings, on the dimension $n$, and we omit it from the notation.

We next define the generalized hep-derivors determined by a derivor.

\begin{definition}\label{DChangeSigDerivor}\lean{equivalent}{ChangeOfSignature.homogeneous,ChangeOfSignature.heterogeneous}
Let $n\in\mathbb{N}$.
\begin{enumerate}
\item The \emph{homogeneous $n$-categorial derivor determined by $\mathbf{c}$} is the generalized hep-derivor
$$
\mathbf{c}^{(n)}=\left(\varphi,c^{(n)},\varnothing\right)
$$
from $\sig{\Sigma}^{(n)}$ to $\sig{\Lambda}^{(n)}$, with empty domains, whose mapping on sorts is $\varphi$ and whose assignment of terms is given below: the formula for $c^{(n)}(\sigma)$ applies to every $\sigma\in\Sigma_{\mathbf{s},s}$, and the structural formulas apply to every $k\in n$ and $s\in S$:
\begin{align*}
c^{(n)}(\sigma)&=c(\sigma), &
c^{(n)}(\circ^{k}_{s})&=v^{\varphi(s)}_{0}\circ^{k}_{\varphi(s)}v^{\varphi(s)}_{1},\\
c^{(n)}(\src^{k}_{s})&=\src^{k}_{\varphi(s)}\left(v^{\varphi(s)}_{0}\right), &
c^{(n)}(\tgt^{k}_{s})&=\tgt^{k}_{\varphi(s)}\left(v^{\varphi(s)}_{0}\right).
\end{align*}
\item The \emph{heterogeneous $n$-categorial derivor determined by $\mathbf{c}$} is the generalized hep-derivor
$$
\mathbf{c}^{\htg(n)}=\left(\varphi^{\htg(n)},c^{\htg(n)},\varnothing\right)
$$
from $\sig{\Sigma}^{\htg(n)}$ to $\sig{\Lambda}^{\htg(n)}$, with empty domains, whose mapping on sorts sends $(k,s)$ to $\varphi^{\htg(n)}(k,s)=(k,\varphi(s))$ and whose assignment of terms is given below: the formula for $c^{\htg(n)}(\sigma_k)$ applies to every $\sigma\in\Sigma_{\mathbf{s},s}$ and $k\in n+1$, and the structural formulas apply to every pair $j<k$ in $n+1$ and $s\in S$:
\begin{align*}
c^{\htg(n)}(\sigma_{k})&=c(\sigma)_{k}, &
c^{\htg(n)}(\mathrm{i}^{(k,j)}_{s})&=\mathrm{i}^{(k,j)}_{\varphi(s)}\left(v^{(j,\varphi(s))}_{0}\right),\\
c^{\htg(n)}(\src^{(j,k)}_{s})&=\src^{(j,k)}_{\varphi(s)}\left(v^{(k,\varphi(s))}_{0}\right), &
c^{\htg(n)}(\tgt^{(j,k)}_{s})&=\tgt^{(j,k)}_{\varphi(s)}\left(v^{(k,\varphi(s))}_{0}\right),\\
c^{\htg(n)}(\circ^{(j,k)}_{s})&=v^{(k,\varphi(s))}_{0}\circ^{(j,k)}_{\varphi(s)}v^{(k,\varphi(s))}_{1}.
\end{align*}
\end{enumerate}
\end{definition}

Let us note that $c(\sigma)$ is a term of $\Lambda$ with variables in $\vs{\varphi^{\star}(\mathbf{s})}$, hence of $\Lambda^{(n)}$, and that, since $(\varphi^{\htg(n)})^{\star}((k,\mathbf{s}))=(k,\varphi^{\star}(\mathbf{s}))$, $c(\sigma)_{k}$ is a term of $\Lambda^{\htg(n)}$ of sort $(k,\varphi(s))$ with variables in $\vs{(k,\varphi^{\star}(\mathbf{s}))}$, so that the above definition is correct.
To simplify the notation, we write $\mathbf{c}^{(n)\ast}_{\gd}$ for the functor $(\mathbf{c}^{(n)})^{\ast}_{\gd}$ induced by $\mathbf{c}^{(n)}$, see Proposition~\ref{PGenHepDerivorFunctor}, and $\mathbf{c}^{\htg(n)\ast}_{\gd}$ for the functor $(\mathbf{c}^{\htg(n)})^{\ast}_{\gd}$ induced by $\mathbf{c}^{\htg(n)}$, and similarly for $\omega$ in place of $n$.
For $n=0$, since $\Sigma^{(0)}=\Sigma$ and $\Lambda^{(0)}=\Lambda$, $\mathbf{c}^{(0)}$ is the derivor $\mathbf{c}$, regarded as the generalized hep-derivor $(\varphi,c,\varnothing)$, see Remark~\ref{RHepAlgCat}.

In the following lemma we prove that, in a higher-order categorial algebra, the term operations of the terms of the base signature $\Lambda$, which are total, are higher-order functors.
This does not extend to arbitrary terms of $\Lambda^{(n)}$, which may contain partial compositions.

\begin{lemma}\label{LTermOpCatAlg}\lean{equivalent}{Homogeneous.CategorialAlgebra.termFunctor,Heterogeneous.CategorialAlgebra.termFunctor}
Let $n\in\mathbb{N}$, $\mathbf{t}\in T^{\star}$, $t\in T$ and $P\in\mathrm{T}_{\Lambda}(\vs{\mathbf{t}})_{t}$.
Then the following statements hold:
\begin{enumerate}
\item for every object $\mathbf{B}$ in $\mathsf{EPAlg}(\spec{C}^{\sig{\Lambda}^{(n)}})$, $P^{\mathbf{B}}$ is total and, for every $k\in n$ and every pair of families $(x_{l})_{l\in\bb{\mathbf{t}}}$ and $(y_{l})_{l\in\bb{\mathbf{t}}}$ in $B_{\mathbf{t}}$, the following equalities hold
\begin{align*}
\src^{k\,\mathbf{B}}_{t}\left(P^{\mathbf{B}}\left((x_{l})_{l\in\bb{\mathbf{t}}}\right)\right)&=P^{\mathbf{B}}\left(\left(\src^{k\,\mathbf{B}}_{t_{l}}(x_{l})\right)_{l\in\bb{\mathbf{t}}}\right),
\\
\tgt^{k\,\mathbf{B}}_{t}\left(P^{\mathbf{B}}\left((x_{l})_{l\in\bb{\mathbf{t}}}\right)\right)&=P^{\mathbf{B}}\left(\left(\tgt^{k\,\mathbf{B}}_{t_{l}}(x_{l})\right)_{l\in\bb{\mathbf{t}}}\right),
\end{align*}
and, if $x_{l}\circ^{k\,\mathbf{B}}_{t_{l}}y_{l}$ is defined for every $l\in\bb{\mathbf{t}}$, then
\begin{enumerate}
\item $P^{\mathbf{B}}((x_{l})_{l\in\bb{\mathbf{t}}})\circ^{k\,\mathbf{B}}_{t}P^{\mathbf{B}}((y_{l})_{l\in\bb{\mathbf{t}}})$ is defined; and
\item the following equality holds
$$
P^{\mathbf{B}}\left((x_{l})_{l\in\bb{\mathbf{t}}}\right)\circ^{k\,\mathbf{B}}_{t}P^{\mathbf{B}}\left((y_{l})_{l\in\bb{\mathbf{t}}}\right)=P^{\mathbf{B}}\left(\left(x_{l}\circ^{k\,\mathbf{B}}_{t_{l}}y_{l}\right)_{l\in\bb{\mathbf{t}}}\right);
$$
\end{enumerate}
\item for every object $\mathbf{B}$ in $\mathsf{EPAlg}(\spec{C}^{\sig{\Lambda}^{\htg(n)}})$ and every $k\in n+1$, $P_{k}^{\mathbf{B}}$ is total and, for every $j<k$, every pair of families $(x_{l})_{l\in\bb{\mathbf{t}}}$ and $(y_{l})_{l\in\bb{\mathbf{t}}}$ in $B_{(k,\mathbf{t})}$ and every family $(u_{l})_{l\in\bb{\mathbf{t}}}$ in $B_{(j,\mathbf{t})}$, the following equalities hold
\begin{align*}
\src^{(j,k)\,\mathbf{B}}_{t}\left(P_{k}^{\mathbf{B}}\left((x_{l})_{l\in\bb{\mathbf{t}}}\right)\right)&=P_{j}^{\mathbf{B}}\left(\left(\src^{(j,k)\,\mathbf{B}}_{t_{l}}(x_{l})\right)_{l\in\bb{\mathbf{t}}}\right),
\\
\tgt^{(j,k)\,\mathbf{B}}_{t}\left(P_{k}^{\mathbf{B}}\left((x_{l})_{l\in\bb{\mathbf{t}}}\right)\right)&=P_{j}^{\mathbf{B}}\left(\left(\tgt^{(j,k)\,\mathbf{B}}_{t_{l}}(x_{l})\right)_{l\in\bb{\mathbf{t}}}\right),
\\
\mathrm{i}^{(k,j)\,\mathbf{B}}_{t}\left(P_{j}^{\mathbf{B}}\left((u_{l})_{l\in\bb{\mathbf{t}}}\right)\right)&=P_{k}^{\mathbf{B}}\left(\left(\mathrm{i}^{(k,j)\,\mathbf{B}}_{t_{l}}(u_{l})\right)_{l\in\bb{\mathbf{t}}}\right),
\end{align*}
and, if $x_{l}\circ^{(j,k)\,\mathbf{B}}_{t_{l}}y_{l}$ is defined for every $l\in\bb{\mathbf{t}}$, then
\begin{enumerate}
\item $P_{k}^{\mathbf{B}}((x_{l})_{l\in\bb{\mathbf{t}}})\circ^{(j,k)\,\mathbf{B}}_{t}P_{k}^{\mathbf{B}}((y_{l})_{l\in\bb{\mathbf{t}}})$ is defined; and
\item the following equality holds
$$
P_{k}^{\mathbf{B}}\left((x_{l})_{l\in\bb{\mathbf{t}}}\right)\circ^{(j,k)\,\mathbf{B}}_{t}P_{k}^{\mathbf{B}}\left((y_{l})_{l\in\bb{\mathbf{t}}}\right)=P_{k}^{\mathbf{B}}\left(\left(x_{l}\circ^{(j,k)\,\mathbf{B}}_{t_{l}}y_{l}\right)_{l\in\bb{\mathbf{t}}}\right).
$$
\end{enumerate}
\end{enumerate}
\end{lemma}
\begin{proof}
(1) Let $\mathbf{B}$ be an object in $\mathsf{EPAlg}(\spec{C}^{\sig{\Lambda}^{(n)}})$, let $k\in n$ and let $(x_{l})_{l\in\bb{\mathbf{t}}}$ and $(y_{l})_{l\in\bb{\mathbf{t}}}$ be families in $B_{\mathbf{t}}$.
We proceed by induction on the term $P$, see Proposition~\ref{PTermChar} and Definition~\ref{DTermOp}.

If $P=v^{\mathbf{t}}_{q}$, for $q\in\bb{\mathbf{t}}$, then $P^{\mathbf{B}}$ is total and sends $(x_{l})_{l\in\bb{\mathbf{t}}}$ to $x_{q}$, and the statements are trivial.

Let us assume that $P=\tau((Q_{i})_{i\in\bb{\mathbf{r}}})$, with $\tau\in\Lambda_{\mathbf{r},t}$, and that the statements hold for every $Q_{i}$.
By~(A0), $\tau^{\mathbf{B}}$ is total, so that, by the induction hypothesis, $P^{\mathbf{B}}$ is total.
Moreover, the following chain of equalities holds
\begin{align*}
\src^{k\,\mathbf{B}}_{t}\left(P^{\mathbf{B}}\left((x_{l})_{l\in\bb{\mathbf{t}}}\right)\right)
&=
\src^{k\,\mathbf{B}}_{t}\left(\tau^{\mathbf{B}}\left(\left(Q_{i}^{\mathbf{B}}\left((x_{l})_{l\in\bb{\mathbf{t}}}\right)\right)_{i\in\bb{\mathbf{r}}}\right)\right)
\tag{1}
\\&=
\tau^{\mathbf{B}}\left(\left(\src^{k\,\mathbf{B}}_{r_{i}}\left(Q_{i}^{\mathbf{B}}\left((x_{l})_{l\in\bb{\mathbf{t}}}\right)\right)\right)_{i\in\bb{\mathbf{r}}}\right)
\tag{2}
\\&=
\tau^{\mathbf{B}}\left(\left(Q_{i}^{\mathbf{B}}\left(\left(\src^{k\,\mathbf{B}}_{t_{l}}(x_{l})\right)_{l\in\bb{\mathbf{t}}}\right)\right)_{i\in\bb{\mathbf{r}}}\right)
\tag{3}
\\&=
P^{\mathbf{B}}\left(\left(\src^{k\,\mathbf{B}}_{t_{l}}(x_{l})\right)_{l\in\bb{\mathbf{t}}}\right).
\tag{4}
\end{align*}
The equalities~(1) and~(4) unravel Definition~\ref{DTermOp}, the equality~(2) follows from~(A$^{k}$7), and the equality~(3) follows from the induction hypothesis.
The equality for the $k$-targets is proved in the same way, with $\tgt$ in place of $\src$, since~(A$^{k}$7) also states that $\tgt^{k\,\mathbf{B}}_{t}(\tau^{\mathbf{B}}((z_{i})_{i\in\bb{\mathbf{r}}}))=\tau^{\mathbf{B}}((\tgt^{k\,\mathbf{B}}_{r_{i}}(z_{i}))_{i\in\bb{\mathbf{r}}})$, for every family $(z_{i})_{i\in\bb{\mathbf{r}}}$ in $B_{\mathbf{r}}$.

Let us assume that $x_{l}\circ^{k\,\mathbf{B}}_{t_{l}}y_{l}$ is defined for every $l\in\bb{\mathbf{t}}$.
By the induction hypothesis, for every $i\in\bb{\mathbf{r}}$, $Q_{i}^{\mathbf{B}}((x_{l})_{l\in\bb{\mathbf{t}}})\circ^{k\,\mathbf{B}}_{r_{i}}Q_{i}^{\mathbf{B}}((y_{l})_{l\in\bb{\mathbf{t}}})$ is defined and equal to $Q_{i}^{\mathbf{B}}((x_{l}\circ^{k\,\mathbf{B}}_{t_{l}}y_{l})_{l\in\bb{\mathbf{t}}})$.
Hence, by~(A$^{k}$8), $P^{\mathbf{B}}((x_{l})_{l\in\bb{\mathbf{t}}})\circ^{k\,\mathbf{B}}_{t}P^{\mathbf{B}}((y_{l})_{l\in\bb{\mathbf{t}}})$ is defined and the following chain of equalities holds
\begin{flushleft}
$
P^{\mathbf{B}}\left((x_{l})_{l\in\bb{\mathbf{t}}}\right)\circ^{k\,\mathbf{B}}_{t}P^{\mathbf{B}}\left((y_{l})_{l\in\bb{\mathbf{t}}}\right)
$
\begin{align*}
&=
\tau^{\mathbf{B}}\left(\left(Q_{i}^{\mathbf{B}}\left((x_{l})_{l\in\bb{\mathbf{t}}}\right)\circ^{k\,\mathbf{B}}_{r_{i}}Q_{i}^{\mathbf{B}}\left((y_{l})_{l\in\bb{\mathbf{t}}}\right)\right)_{i\in\bb{\mathbf{r}}}\right)
\tag{1}
\\&=
\tau^{\mathbf{B}}\left(\left(Q_{i}^{\mathbf{B}}\left(\left(x_{l}\circ^{k\,\mathbf{B}}_{t_{l}}y_{l}\right)_{l\in\bb{\mathbf{t}}}\right)\right)_{i\in\bb{\mathbf{r}}}\right)
\tag{2}
\\&=
P^{\mathbf{B}}\left(\left(x_{l}\circ^{k\,\mathbf{B}}_{t_{l}}y_{l}\right)_{l\in\bb{\mathbf{t}}}\right).
\tag{3}
\end{align*}
\end{flushleft}
The equality~(1) follows from Definition~\ref{DTermOp} and~(A$^{k}$8), the equality~(2) follows from the induction hypothesis, and the equality~(3) unravels Definition~\ref{DTermOp}.

(2) Let $\mathbf{B}$ be an object in $\mathsf{EPAlg}(\spec{C}^{\sig{\Lambda}^{\htg(n)}})$, let $j<k$ in $n+1$ and let $(x_{l})_{l\in\bb{\mathbf{t}}}$ and $(y_{l})_{l\in\bb{\mathbf{t}}}$ be families in $B_{(k,\mathbf{t})}$.
By the note after Definition~\ref{DInclusionDerivors}, $(v^{\mathbf{t}}_{q})_{k}=v^{(k,\mathbf{t})}_{q}$ and $(\tau((Q_{i})_{i\in\bb{\mathbf{r}}}))_{k}=\tau_{k}(((Q_{i})_{k})_{i\in\bb{\mathbf{r}}})$, so that
$$
\left(\tau\left((Q_{i})_{i\in\bb{\mathbf{r}}}\right)\right)_{k}^{\mathbf{B}}\left((x_{l})_{l\in\bb{\mathbf{t}}}\right)=\tau^{\mathbf{B}}_{k}\left(\left((Q_{i})_{k}^{\mathbf{B}}\left((x_{l})_{l\in\bb{\mathbf{t}}}\right)\right)_{i\in\bb{\mathbf{r}}}\right),
$$
and similarly at the dimension $j$.
Hence the proof is the same as that of~(1), by induction on $P$, with the operations $\tau^{\mathbf{B}}_{k}$, $\src^{(j,k)\,\mathbf{B}}$, $\tgt^{(j,k)\,\mathbf{B}}$, $\mathrm{i}^{(k,j)\,\mathbf{B}}$ and $\circ^{(j,k)\,\mathbf{B}}$ of $\mathbf{B}$ and the axioms~(EA0), (EA$^{j,k}$7) and~(EA$^{j,k}$8) in place of $\tau^{\mathbf{B}}$, $\src^{k\,\mathbf{B}}$, $\tgt^{k\,\mathbf{B}}$ and $\circ^{k\,\mathbf{B}}$ and the axioms~(A0), (A$^{k}$7) and~(A$^{k}$8); the third equation of~(EA$^{j,k}$7) gives the statement for the identities. This completes the proof.
\end{proof}

We next prove that the generalized hep-derivors of Definition~\ref{DChangeSigDerivor} send categorial algebras to categorial algebras.

\begin{proposition}\label{PChangeSigEpVariety}\lean{weaker}{ChangeOfSignature.homogeneous_isModel,ChangeOfSignature.heterogeneous_isModel,ChangeOfSignature.homogeneousRealization,ChangeOfSignature.heterogeneousRealization}
Let $n\in\mathbb{N}$.
Then the following statements hold:
\begin{enumerate}
\item the restriction of $\mathbf{c}^{(n)\ast}_{\gd}$ to $\mathsf{EPAlg}(\spec{C}^{\sig{\Lambda}^{(n)}})$ corestricts to $\mathsf{EPAlg}(\spec{C}^{\sig{\Sigma}^{(n)}})$; moreover, for every object $\mathbf{B}$ in $\mathsf{EPAlg}(\spec{C}^{\sig{\Lambda}^{(n)}})$, for every $s\in S$, the cells of sort $s$ of $\mathbf{c}^{(n)\ast}_{\gd}(\mathbf{B})$ are the elements of $B_{\varphi(s)}$, for every $(\mathbf{s},s)\in S^{\star}\times S$, every $\sigma\in\Sigma_{\mathbf{s},s}$ acts as $c(\sigma)^{\mathbf{B}}$ and, for every $k\in n$ and every $s\in S$, $\src^{k}_{s}$, $\tgt^{k}_{s}$ and $\circ^{k}_{s}$ act as $\src^{k\,\mathbf{B}}_{\varphi(s)}$, $\tgt^{k\,\mathbf{B}}_{\varphi(s)}$ and $\circ^{k\,\mathbf{B}}_{\varphi(s)}$, respectively; and
\item the restriction of $\mathbf{c}^{\htg(n)\ast}_{\gd}$ to $\mathsf{EPAlg}(\spec{C}^{\sig{\Lambda}^{\htg(n)}})$ corestricts to $\mathsf{EPAlg}(\spec{C}^{\sig{\Sigma}^{\htg(n)}})$; moreover, for every object $\mathbf{B}$ in $\mathsf{EPAlg}(\spec{C}^{\sig{\Lambda}^{\htg(n)}})$, for every $k\in n+1$ and every $s\in S$, the cells of sort $(k,s)$ of $\mathbf{c}^{\htg(n)\ast}_{\gd}(\mathbf{B})$ are the elements of $B_{(k,\varphi(s))}$, for every $(\mathbf{s},s)\in S^{\star}\times S$, every $\sigma\in\Sigma_{\mathbf{s},s}$ and every $k\in n+1$, $\sigma_{k}$ acts as $c(\sigma)_{k}^{\mathbf{B}}$ and, for every pair of dimensions $j<k$ in $n+1$ and every $s\in S$, $\src^{(j,k)}_{s}$, $\tgt^{(j,k)}_{s}$, $\mathrm{i}^{(k,j)}_{s}$ and $\circ^{(j,k)}_{s}$ act as $\src^{(j,k)\,\mathbf{B}}_{\varphi(s)}$, $\tgt^{(j,k)\,\mathbf{B}}_{\varphi(s)}$, $\mathrm{i}^{(k,j)\,\mathbf{B}}_{\varphi(s)}$ and $\circ^{(j,k)\,\mathbf{B}}_{\varphi(s)}$, respectively.
\end{enumerate}
\end{proposition}
\begin{proof}
(1) Let $\mathbf{B}$ be an object in $\mathsf{EPAlg}(\spec{C}^{\sig{\Lambda}^{(n)}})$.
Since the domains of $\mathbf{c}^{(n)}$ are empty, by Remark~\ref{RDomainSetEmpty}, the cells of sort $s$ of $\mathbf{c}^{(n)\ast}_{\gd}(\mathbf{B})$ are the elements of $B_{\varphi(s)}$.
By Definitions~\ref{DCorestriction} and~\ref{DEpVariety}, we have to prove that $\mathbf{c}^{(n)\ast}_{\gd}(\mathbf{B})$ is a $\sig{\Sigma}^{(n)}$-algebra in which every axiom of $\mathcal{C}^{\sig{\Sigma}^{(n)}}$ holds.

We first prove that $\mathbf{c}^{(n)\ast}_{\gd}(\mathbf{B})$ is a $\sig{\Sigma}^{(n)}$-algebra.
By Proposition~\ref{PTransformedEpAlg}, we check that, for every operation symbol of $\Sigma^{(n)}$, the corresponding partial term operation is defined on the families that satisfy its domain conditions in $\mathbf{c}^{(n)\ast}_{\gd}(\mathbf{B})$; the values are cells, since every element of $B_{\varphi(s)}$ is a cell.
For every $(\mathbf{s},s)\in S^{\star}\times S$ and every $\sigma\in\Sigma_{\mathbf{s},s}$, the set of domain conditions of $\sigma$ is empty and the partial term operation $c(\sigma)^{\mathbf{B}}$ is total, by Lemma~\ref{LTermOpCatAlg}~(1), also when $n=0$.

For the remaining operation symbols, fix $k\in n$ and $s\in S$.
The sets of domain conditions of $\src^{k}_{s}$ and $\tgt^{k}_{s}$ are empty, and the partial term operations $\src^k_{\varphi(s)}(v^{\varphi(s)}_0)^{\mathbf{B}}$ and $\tgt^k_{\varphi(s)}(v^{\varphi(s)}_0)^{\mathbf{B}}$, which are $\src^{k\,\mathbf{B}}_{\varphi(s)}$ and $\tgt^{k\,\mathbf{B}}_{\varphi(s)}$, are total, by~(A$^k$1).

Next, we consider the compositions.
By the above, two cells $x$, $y$ of sort $s$ satisfy the domain condition of $\circ^{k}_{s}$ in $\mathbf{c}^{(n)\ast}_{\gd}(\mathbf{B})$ if, and only if, $\src^{k\,\mathbf{B}}_{\varphi(s)}(x)=\tgt^{k\,\mathbf{B}}_{\varphi(s)}(y)$, in which case, by~(A$^{k}$3), $x\circ^{k\,\mathbf{B}}_{\varphi(s)}y$ is defined.

Therefore, by Proposition~\ref{PTransformedEpAlg}, $\mathbf{c}^{(n)\ast}_{\gd}(\mathbf{B})$ is a $\sig{\Sigma}^{(n)}$-algebra and, by Proposition~\ref{PTransformedAlg}, its operations are as stated.

Now we prove that the axioms hold.
By the above, the axioms of $\spec{C}^{\sig{\Sigma}^{(n)}}$ in $\mathbf{c}^{(n)\ast}_{\gd}(\mathbf{B})$ become statements on $\mathbf{B}$.
The axiom (A0) holds by Lemma~\ref{LTermOpCatAlg}~(1).
For every $s\in S$, the axioms (A$^{k}$1)--(A$^{k}$6), for $k\in n$, and (AB$^{l,k}$1)--(AB$^{l,k}$3), for $l<k$ in $n$, only involve the sources, the targets and the compositions of sort $s$, which are those of sort $\varphi(s)$ of $\mathbf{B}$; hence they become the same axioms for the sort $\varphi(s)$, which hold in $\mathbf{B}$.
Finally, for every $\sigma\in\Sigma_{\mathbf{s},s}$, the axioms (A$^{k}$7) and (A$^{k}$8) for $\sigma$ become the statements of Lemma~\ref{LTermOpCatAlg}~(1) for the term $c(\sigma)$, which hold in $\mathbf{B}$.

All in all, $\mathbf{c}^{(n)\ast}_{\gd}(\mathbf{B})$ is an object in $\mathsf{EPAlg}(\spec{C}^{\sig{\Sigma}^{(n)}})$.

(2) The proof is the same as that of~(1).
The domains of $\mathbf{c}^{\htg(n)}$ are empty, so that the cells of sort $(k,s)$ of $\mathbf{c}^{\htg(n)\ast}_{\gd}(\mathbf{B})$ are the elements of $B_{(k,\varphi(s))}$.
The operation symbols $\sigma_{k}$, $\src^{(j,k)}_{s}$, $\tgt^{(j,k)}_{s}$ and $\mathrm{i}^{(k,j)}_{s}$ have empty sets of domain conditions and, by Lemma~\ref{LTermOpCatAlg}~(2) and~(EA$^{j,k}$1), they are sent to terms that are defined everywhere in $\mathbf{B}$; and, by~(EA$^{j,k}$3), $x\circ^{(j,k)\,\mathbf{B}}_{\varphi(s)}y$ is defined for every pair of cells $x$, $y$ of sort $(k,s)$ that satisfies the domain condition of $\circ^{(j,k)}_{s}$.
Hence, by Propositions~\ref{PTransformedEpAlg} and~\ref{PTransformedAlg}, $\mathbf{c}^{\htg(n)\ast}_{\gd}(\mathbf{B})$ is a $\sig{\Sigma}^{\htg(n)}$-algebra whose operations are as stated.
As for the axioms, (EA0) holds by Lemma~\ref{LTermOpCatAlg}~(2); the axioms (EA$^{j,k}$1)--(EA$^{j,k}$6) and (EA$^{h,j,k}$1)--(EA$^{h,j,k}$5) for the sort $s$ become the same axioms for the sort $\varphi(s)$; and, for every operation symbol $\sigma$ of $\Sigma$, the axioms (EA$^{j,k}$7) and (EA$^{j,k}$8) for $\sigma$ become the statements of Lemma~\ref{LTermOpCatAlg}~(2) for the term $c(\sigma)$. This completes the proof.
\end{proof}

From Proposition~\ref{PChangeSigEpVariety} we obtain the functors induced by a derivor between the categories of higher-order categorial algebras.

\begin{corollary}\label{CChangeSigFunctor}\lean{equivalent}{ChangeOfSignature.homogeneousFunctor,ChangeOfSignature.heterogeneousFunctor}
Let $n\in\mathbb{N}$.
Then the restriction of $\mathbf{c}^{(n)\ast}_{\gd}$ to $\mathsf{nCatAlg}(\sig{\Lambda})$ corestricts to $\mathsf{nCatAlg}(\sig{\Sigma})$, and the restriction of $\mathbf{c}^{\htg(n)\ast}_{\gd}$ to $\mathsf{nCat}^{\htg}\mathsf{Alg}(\sig{\Lambda})$ corestricts to $\mathsf{nCat}^{\htg}\mathsf{Alg}(\sig{\Sigma})$.
\end{corollary}
\begin{proof}
It follows from Proposition~\ref{PChangeSigEpVariety} and Remarks~\ref{RCatAlgAsEpAlg} and~\ref{RCatAlgAsEpAlgHtg}.
\end{proof}

By abuse of notation, we denote the corestrictions of Corollary~\ref{CChangeSigFunctor} by
\begin{enumerate}
	\item
	$ \mathbf{c}^{(n)\ast}_{\gd}\colon\mathsf{nCatAlg}(\sig{\Lambda})\mor\mathsf{nCatAlg}(\sig{\Sigma}) $; and

	\item
	$ \mathbf{c}^{\htg(n)\ast}_{\gd}\colon\mathsf{nCat}^{\htg}\mathsf{Alg}(\sig{\Lambda})\mor\mathsf{nCat}^{\htg}\mathsf{Alg}(\sig{\Sigma}) $.
\end{enumerate}
Explicitly, $\mathbf{c}^{(n)\ast}_{\gd}$ sends an $n$-categorial $\sig{\Lambda}$-algebra $\calg{B}$ to the $n$-categorial $\sig{\Sigma}$-algebra whose homogeneous $n$-category of sort $s$ is the homogeneous $n$-category of sort $\varphi(s)$ of $\calg{B}$ and in which every operation symbol $\sigma$ of $\Sigma$ acts as the $n$-functor $c(\sigma)^{\calg{B}}$, see Lemma~\ref{LTermOpCatAlg}; and, by Proposition~\ref{PGenHepDerivorFunctor}, it sends an $n$-categorial $\sig{\Lambda}$-homomorphism $F$ to $F_{\varphi}=(F_{\varphi(s)})_{s\in S}$.
Similarly, $\mathbf{c}^{\htg(n)\ast}_{\gd}$ sends a heterogeneous $n$-categorial $\sig{\Lambda}$-algebra $\calg{B}$ to the heterogeneous $n$-categorial $\sig{\Sigma}$-algebra whose heterogeneous $n$-category of sort $s$ is that of sort $\varphi(s)$ of $\calg{B}$ and in which every $\sigma_{k}$, for $\sigma$ an operation symbol of $\Sigma$ and $k\in n+1$, acts as $c(\sigma)_{k}^{\calg{B}}$; and it sends $F$ to $F_{\varphi^{\htg(n)}}=(F_{(k,\varphi(s))})_{(k,s)\in(n+1)\times S}$.
For $n=0$, $\mathbf{c}^{(0)\ast}_{\gd}$ is the usual functor from $\mathsf{Alg}(\Lambda)$ to $\mathsf{Alg}(\Sigma)$ induced by $\mathbf{c}$, see Remark~\ref{RHepAlgCat}.

In the following lemma we collect the translations of the terms $P$ and $P_{k}$ along the generalized hep-derivors of this work that we need.
By Corollary~\ref{CDerivorCatLaws}, for every pair of derivors $\mathbf{d}=(\psi,d)$ from $(T,\Lambda)$ to $(U,\Omega)$ and $\mathbf{e}=(\chi,e)$ from $(U,\Omega)$ to $(W,\Xi)$, $(e^{\sharp}\circ d)^{\sharp}=e^{\sharp}\circ d^{\sharp}$, so that $e^{\sharp}\circ d^{\sharp}$ only depends on the composite derivor $(\chi\circ\psi,e^{\sharp}\circ d)$, see Definition~\ref{DDerivorCompId}.
Hence, in order to prove an equality of the form $e^{\sharp}\circ d^{\sharp}=e'^{\sharp}\circ d'^{\sharp}$, it suffices to prove that the composite derivors have the same mapping on sorts and that $e^{\sharp}(d(\tau))=e'^{\sharp}(d'(\tau))$, for every operation symbol $\tau$ of $\Lambda$.

\begin{lemma}\label{LChangeSigTerms}\lean{equivalent}{ChangeOfSignature.homogeneousInclusion_comp,ChangeOfSignature.heterogeneousInclusion_comp,ChangeOfSignature.heterogeneousInclusion_comp_graduation,ChangeOfSignature.homogeneousInclusion_comp_unification,ChangeOfSignature.homogeneousInclusion_comp_discretization,ChangeOfSignature.homogeneousInclusion_comp_truncation,ChangeOfSignature.heterogeneousInclusion_comp_dimension}
Let $\sig{\Omega}=(U,\Omega)$ be a many-sorted signature, $\mathbf{d}=(\psi,d)$ a derivor from $\sig{\Lambda}$ to $\sig{\Omega}$, $m,n\in\mathbb{N}$ with $m\leq n$, $\mathbf{t}\in T^{\star}$, $t\in T$ and $P\in\mathrm{T}_{\Lambda}(\vs{\mathbf{t}})_{t}$.
Then the following statements hold:
\begin{enumerate}
\item $(d^{(n)})^{\sharp}(P)=d^{\sharp}(P)$ and, for every $k\in n+1$, $(d^{\htg(n)})^{\sharp}(P_{k})=(d^{\sharp}(P))_{k}$;
\item for every $k\in n+1$, $(\Gd^{(n)}_{\sig{\Lambda}})^{\sharp}(P_{k})=P$, and $(\Uf^{(n)}_{\sig{\Lambda}})^{\sharp}(P)=P_{n}$;
\item $(\Ds^{(n,m)}_{\sig{\Lambda}})^{\sharp}(P)=P$ and $(\Tr^{(m,n)}_{\sig{\Lambda}})^{\sharp}(P)=P$; and
\item for every $k\in n+1$, $(\Ds^{\htg(n,m)}_{\sig{\Lambda}})^{\sharp}(P_{k})=P_{\min\{k,m\}}$, and, for every $k\in m+1$, $(\Tr^{\htg(m,n)}_{\sig{\Lambda}})^{\sharp}(P_{k})=P_{k}$.
\end{enumerate}
\end{lemma}
\begin{proof}
By the note before this lemma, in each case it suffices to compare the composite derivors.
Let $\tau\in\Lambda_{\mathbf{r},t'}$.

(1) The equality $(d^{(n)})^{\sharp}\circ(\iota^{(n)})^{\sharp}=(\iota^{(n)})^{\sharp}\circ d^{\sharp}$ follows from the fact that both composite derivors have $\psi$ as mapping on sorts and that, by Definition~\ref{DChangeSigDerivor},
$$
(d^{(n)})^{\sharp}\left(\tau\left((v^{\mathbf{r}}_{p})_{p\in\bb{\mathbf{r}}}\right)\right)=d^{(n)}(\tau)=d(\tau)=(\iota^{(n)})^{\sharp}(d(\tau)).
$$
Similarly, the equality $(d^{\htg(n)})^{\sharp}\circ(\iota^{\htg(n)}_{k})^{\sharp}=(\iota^{\htg(n)}_{k})^{\sharp}\circ d^{\sharp}$ follows from the fact that both composite derivors send $t''\in T$ to $(k,\psi(t''))$ and that
$$
(d^{\htg(n)})^{\sharp}\left(\tau_{k}\left((v^{(k,\mathbf{r})}_{p})_{p\in\bb{\mathbf{r}}}\right)\right)=d^{\htg(n)}(\tau_{k})=d(\tau)_{k}=(\iota^{\htg(n)}_{k})^{\sharp}(d(\tau)).
$$

(2) By Definition~\ref{DGdDerivor}, $\Gd^{(n)}_{\sig{\Lambda}}(\tau_{k})=\tau((v^{\mathbf{r}}_{p})_{p\in\bb{\mathbf{r}}})$ and the mapping on sorts of $\bGd^{(n)}_{\sig{\Lambda}}$ sends $(k,t'')$ to $t''$.
Hence the composite of $\boldsymbol{\iota}^{\htg(n)}_{k}$ and $\bGd^{(n)}_{\sig{\Lambda}}$ is $\boldsymbol{\iota}^{(n)}$.
Similarly, by Definition~\ref{DUfDerivor}, the composite of $\boldsymbol{\iota}^{(n)}$ and $\bUf^{(n)}_{\sig{\Lambda}}$ is $\boldsymbol{\iota}^{\htg(n)}_{n}$, since $\Uf^{(n)}_{\sig{\Lambda}}(\tau)=\tau_{n}((v^{(n,\mathbf{r})}_{p})_{p\in\bb{\mathbf{r}}})$ and the mapping on sorts of $\bUf^{(n)}_{\sig{\Lambda}}$ sends $t''$ to $(n,t'')$.

(3) By Definitions~\ref{DDsDerivor} and~\ref{DTrDerivor}, both $\bDs^{(n,m)}_{\sig{\Lambda}}$ and $\bTr^{(m,n)}_{\sig{\Lambda}}$ have $\id_{T}$ as mapping on sorts and send $\tau$ to $\tau((v^{\mathbf{r}}_{p})_{p\in\bb{\mathbf{r}}})$.
Hence the composites of $\boldsymbol{\iota}^{(n)}$ and $\bDs^{(n,m)}_{\sig{\Lambda}}$, and of $\boldsymbol{\iota}^{(m)}$ and $\bTr^{(m,n)}_{\sig{\Lambda}}$, are $\boldsymbol{\iota}^{(m)}$ and $\boldsymbol{\iota}^{(n)}$, respectively.

(4) By Definition~\ref{DDsDerivor}~(2), $\Ds^{\htg(n,m)}_{\sig{\Lambda}}(\tau_{k})=\tau_{\min\{k,m\}}((v^{(\min\{k,m\},\mathbf{r})}_{p})_{p\in\bb{\mathbf{r}}})$ and the mapping on sorts of $\bDs^{\htg(n,m)}_{\sig{\Lambda}}$ sends $(k,t'')$ to $(\min\{k,m\},t'')$.
Hence the composite of $\boldsymbol{\iota}^{\htg(n)}_{k}$ and $\bDs^{\htg(n,m)}_{\sig{\Lambda}}$ is $\boldsymbol{\iota}^{\htg(m)}_{\min\{k,m\}}$.
Similarly, by Definition~\ref{DTrDerivor}, the composite of $\boldsymbol{\iota}^{\htg(m)}_{k}$ and $\bTr^{\htg(m,n)}_{\sig{\Lambda}}$ is $\boldsymbol{\iota}^{\htg(n)}_{k}$.

For $m=n$, the statements~(3) and~(4) are trivial, since the derivors involved are identities. This completes the proof.
\end{proof}

In the following proposition we prove that the assignments $\mathbf{c}\mapsto\mathbf{c}^{(n)}$ and $\mathbf{c}\mapsto\mathbf{c}^{\htg(n)}$ preserve compositions and identities.

\begin{proposition}\label{PChangeSigDerivorComp}\lean{equivalent}{ChangeOfSignature.homogeneous_comp,ChangeOfSignature.homogeneous_identity,ChangeOfSignature.heterogeneous_comp,ChangeOfSignature.heterogeneous_identity}
Let $\sig{\Omega}=(U,\Omega)$ be a many-sorted signature, $\mathbf{d}=(\psi,d)$ a derivor from $\sig{\Lambda}$ to $\sig{\Omega}$ and $n\in\mathbb{N}$.
Then the following equalities hold in $\mathsf{HSig}_{\gd}$
\begin{align*}
(\mathbf{d}\circ\mathbf{c})^{(n)}&=\mathbf{d}^{(n)}\circ\mathbf{c}^{(n)},
&
(\id_{\sig{\Sigma}})^{(n)}&=\id_{\sig{\Sigma}^{(n)}},
\\
(\mathbf{d}\circ\mathbf{c})^{\htg(n)}&=\mathbf{d}^{\htg(n)}\circ\mathbf{c}^{\htg(n)},
&
(\id_{\sig{\Sigma}})^{\htg(n)}&=\id_{\sig{\Sigma}^{\htg(n)}}.
\end{align*}
\end{proposition}
\begin{proof}
We first consider the composition in the homogeneous case.
By Definitions~\ref{DDerivorCompId}, \ref{DGenHepDerivorCompId} and~\ref{DChangeSigDerivor}, both generalized hep-derivors have $\psi\circ\varphi$ as mapping on sorts and empty domains, since $\Dom^{\mathbf{d}^{(n)}\circ\mathbf{c}^{(n)}}_{s}=\varnothing\cup(d^{(n)})^{\sharp}[\varnothing]=\varnothing$, for every $s\in S$.
We compare their assignments of terms.
Let $\sigma\in\Sigma_{\mathbf{s},s}$.
Then the following chain of equalities holds
\begin{align*}
(d^{(n)})^{\sharp}\left(c^{(n)}(\sigma)\right)
&=
(d^{(n)})^{\sharp}\left(c(\sigma)\right)
\tag{1}
\\&=
d^{\sharp}\left(c(\sigma)\right)
\tag{2}
\\&=
(d\circ c)^{(n)}(\sigma).
\tag{3}
\end{align*}
The equality~(1) unravels Definition~\ref{DChangeSigDerivor}, the equality~(2) follows from Lemma~\ref{LChangeSigTerms}~(1), and the equality~(3) follows from Definitions~\ref{DDerivorCompId} and~\ref{DChangeSigDerivor}, where $d\circ c$ stands for $d^{\sharp}\circ c$.
Moreover, for every $k\in n$ and every $s'\in S$, by~(Der1), (Der2) and Definition~\ref{DChangeSigDerivor}, the following equalities hold
\begin{align*}
(d^{(n)})^{\sharp}\left(c^{(n)}(\src^{k}_{s'})\right)&=\src^{k}_{\psi(\varphi(s'))}\left(v^{\psi(\varphi(s'))}_{0}\right)=(d\circ c)^{(n)}(\src^{k}_{s'}),\\
(d^{(n)})^{\sharp}\left(c^{(n)}(\tgt^{k}_{s'})\right)&=\tgt^{k}_{\psi(\varphi(s'))}\left(v^{\psi(\varphi(s'))}_{0}\right)=(d\circ c)^{(n)}(\tgt^{k}_{s'}),\\
(d^{(n)})^{\sharp}\left(c^{(n)}(\circ^{k}_{s'})\right)&=v^{\psi(\varphi(s'))}_{0}\circ^{k}_{\psi(\varphi(s'))}v^{\psi(\varphi(s'))}_{1}=(d\circ c)^{(n)}(\circ^{k}_{s'}).
\end{align*}

Now we consider the identities.
The identity derivor at $\sig{\Sigma}$ is $(\id_{S},\eta^{\Sigma})$, see Definition~\ref{DDerivorCompId}.
Hence $(\id_{\sig{\Sigma}})^{(n)}$ has $\id_{S}$ as mapping on sorts, empty domains, and sends every operation symbol of $\Sigma^{(n)}$ to itself, applied to variables, i.e., it is $(\id_{S},\eta^{\Sigma^{(n)}},\varnothing)=\id_{\sig{\Sigma}^{(n)}}$, see Definition~\ref{DGenHepDerivorCompId}.

Next, we consider the heterogeneous case.
The proof is the same, with Lemma~\ref{LChangeSigTerms}~(1) applied to the terms $c(\sigma)_{k}$: for every operation symbol $\sigma$ of $\Sigma$ and every $k\in n+1$, $(d^{\htg(n)})^{\sharp}(c(\sigma)_{k})=(d^{\sharp}(c(\sigma)))_{k}=(d\circ c)^{\htg(n)}(\sigma_{k})$; and $\eta^{\Sigma}(\sigma)_{k}=\sigma_{k}((v^{(k,\mathbf{s})}_{p})_{p\in\bb{\mathbf{s}}})$; the structural operation symbols are handled as in the homogeneous case.
\end{proof}

From the above proposition we obtain that the assignments $\mathbf{c}\mapsto\mathbf{c}^{(n)\ast}_{\gd}$ and $\mathbf{c}\mapsto\mathbf{c}^{\htg(n)\ast}_{\gd}$ preserve compositions and identities.

\begin{corollary}\label{CChangeSigFunctorComp}\lean{weaker}{ChangeOfSignature.homogeneousFunctorComp,ChangeOfSignature.homogeneousFunctorIdentity,ChangeOfSignature.heterogeneousFunctorComp,ChangeOfSignature.heterogeneousFunctorIdentity}
Let $\sig{\Omega}=(U,\Omega)$ be a many-sorted signature, $\mathbf{d}=(\psi,d)$ a derivor from $\sig{\Lambda}$ to $\sig{\Omega}$ and $n\in\mathbb{N}$.
Then the following equalities hold
\begin{align*}
(\mathbf{d}\circ\mathbf{c})^{(n)\ast}_{\gd}&=\mathbf{c}^{(n)\ast}_{\gd}\circ\mathbf{d}^{(n)\ast}_{\gd},
&
(\id_{\sig{\Sigma}})^{(n)\ast}_{\gd}&=\Id_{\mathsf{nCatAlg}(\sig{\Sigma})},
\\
(\mathbf{d}\circ\mathbf{c})^{\htg(n)\ast}_{\gd}&=\mathbf{c}^{\htg(n)\ast}_{\gd}\circ\mathbf{d}^{\htg(n)\ast}_{\gd},
&
(\id_{\sig{\Sigma}})^{\htg(n)\ast}_{\gd}&=\Id_{\mathsf{nCat}^{\htg}\mathsf{Alg}(\sig{\Sigma})}.
\end{align*}
\end{corollary}
\begin{proof}
Let $\calg{C}$ be an $n$-categorial $\sig{\Omega}$-algebra.
Since all the domains involved are empty, by Remark~\ref{RDomainSetEmpty}, both $\mathbf{c}^{(n)\ast}_{\gd}(\mathbf{d}^{(n)\ast}_{\gd}(\calg{C}))$ and $(\mathbf{d}^{(n)}\circ\mathbf{c}^{(n)})^{\ast}_{\gd}(\calg{C})$ have $C_{\psi\circ\varphi}$ as underlying $S$-sorted set and, by Proposition~\ref{PChangeSigEpVariety}, applied twice, the former is a $\sig{\Sigma}^{(n)}$-algebra.
Therefore, by Proposition~\ref{PChangeSigDerivorComp} and Corollary~\ref{CGenHepDerivorCompEqual}, $(\mathbf{d}\circ\mathbf{c})^{(n)\ast}_{\gd}(\calg{C})=\mathbf{c}^{(n)\ast}_{\gd}(\mathbf{d}^{(n)\ast}_{\gd}(\calg{C}))$, and the same holds for the $n$-categorial $\sig{\Omega}$-homomorphisms.
Finally, by Proposition~\ref{PChangeSigDerivorComp} and Corollary~\ref{CGenHepDerivorIdFunctor}, $(\id_{\sig{\Sigma}})^{(n)\ast}_{\gd}$ is the identity.
The heterogeneous case is proved in the same way, with $\mathbf{c}^{\htg(n)}$, $\mathbf{d}^{\htg(n)}$ and $C_{\psi^{\htg(n)}\circ\varphi^{\htg(n)}}$ in place of $\mathbf{c}^{(n)}$, $\mathbf{d}^{(n)}$ and $C_{\psi\circ\varphi}$. This completes the proof.
\end{proof}

From Corollaries~\ref{CSigdCat} and~\ref{CChangeSigFunctorComp} we obtain the following split indexed categories.

\begin{corollary}\label{CChangeSigIndexed}\lean{weaker}{ChangeOfSignature.homogeneousIndexed,ChangeOfSignature.heterogeneousIndexed,ChangeOfSignature.homogeneousFunctorReductIso,ChangeOfSignature.heterogeneousFunctorReductIso}
Let $n\in\mathbb{N}$.
Then the assignment that sends a many-sorted signature $\sig{\Sigma}$ to $\mathsf{nCatAlg}(\sig{\Sigma})$ and a derivor $\mathbf{c}$ to $\mathbf{c}^{(n)\ast}_{\gd}$ is a split $\mathsf{Sig}_{\mathfrak{d}}$-indexed category, see Definition~\ref{DSplitIndexed}, and so is the assignment that sends $\sig{\Sigma}$ to $\mathsf{nCat}^{\htg}\mathsf{Alg}(\sig{\Sigma})$ and $\mathbf{c}$ to $\mathbf{c}^{\htg(n)\ast}_{\gd}$.
We denote these split $\mathsf{Sig}_{\mathfrak{d}}$-indexed categories by $\mathsf{nCatAlg}$ and $\mathsf{nCat}^{\htg}\mathsf{Alg}$, respectively.
\end{corollary}
\begin{proof}
By Corollary~\ref{CSigdCat}, $\mathsf{Sig}_{\mathfrak{d}}$ is a category and, by Corollary~\ref{CChangeSigFunctor}, both assignments send a derivor to a functor in the opposite direction.
By Corollary~\ref{CChangeSigFunctorComp}, both assignments preserve compositions and identities, i.e., they are contravariant functors from $\mathsf{Sig}_{\mathfrak{d}}$ to $\mathsf{Cat}$.
\end{proof}

\begin{remark}\label{RChangeSigIndexedSize}
When $S\neq\varnothing$, the categories $\mathsf{nCatAlg}(\sig{\Sigma})$ and $\mathsf{nCat}^{\htg}\mathsf{Alg}(\sig{\Sigma})$ are not $\boldsymbol{\mathcal{U}}$-small: even the one-element carriers, with arbitrary labels in $\boldsymbol{\mathcal{U}}$, give a class of objects not belonging to $\boldsymbol{\mathcal{U}}$. If $S=\varnothing$, both categories are terminal and hence $\boldsymbol{\mathcal{U}}$-small.
However, they are objects of $\mathsf{Cat}$, which, see Chapter~\ref{SIntroduction}, is the category of all the $\boldsymbol{\mathcal{U}}$-categories, i.e., of the categories $\mathsf{C}$ such that $\mathrm{Ob}(\mathsf{C})\subseteq\boldsymbol{\mathcal{U}}$ and, for every $x,y\in\mathrm{Ob}(\mathsf{C})$, $\mathrm{Hom}(x,y)\in\boldsymbol{\mathcal{U}}$, so that no larger universe is needed in Corollary~\ref{CChangeSigIndexed}.
Here we use the coding convention of Chapter~\ref{SIntroduction}: a family whose codomain is $\boldsymbol{\mathcal{U}}$, such as an $S$-sorted set $A\colon S\mor\boldsymbol{\mathcal{U}}$ or an $S$-sorted signature, is coded, when it is regarded as an element of $\boldsymbol{\mathcal{U}}$, by its graph, since the triple $(S,F,\boldsymbol{\mathcal{U}})$ itself is not an element of $\boldsymbol{\mathcal{U}}$.
By Definitions~\ref{DSortedSet} and~\ref{DSig}, $S\in\boldsymbol{\mathcal{U}}$ and the values of $\Sigma$ are in $\boldsymbol{\mathcal{U}}$, so that the graphs of $S$-sorted sets and of $\Sigma$ are $\boldsymbol{\mathcal{U}}$-small.
Hence every $n$-categorial $\sig{\Sigma}$-algebra, which, by Remark~\ref{RCatAlgAsEpAlg}, is a partial $\Sigma^{(n)}$-algebra $(A,F)$ with $A$ an $S$-sorted set, is coded by an element of $\boldsymbol{\mathcal{U}}$, and, for every pair of $n$-categorial $\sig{\Sigma}$-algebras, the set of the homomorphisms between them, coded in the same way, is a subset of an element of $\boldsymbol{\mathcal{U}}$, hence an element of $\boldsymbol{\mathcal{U}}$.
The same holds for $\mathsf{nCat}^{\htg}\mathsf{Alg}(\sig{\Sigma})$.
\end{remark}

For $n=0$, the first split indexed category of Corollary~\ref{CChangeSigIndexed} is the split indexed category of~\cite{ClimentSoliveres2010b} that sends $\sig{\Sigma}$ to $\mathsf{Alg}(\Sigma)$.

\section[Compatibility with the presentation and dimension functors]{Compatibility with graduation, unification, discretization and truncation}\label{SChangeSigCompat}

In this section we prove that the functors induced by the generalized hep-derivors $\mathbf{c}^{(n)}$ and $\mathbf{c}^{\htg(n)}$ commute with graduation, unification, discretization and truncation.
As in Chapters~\ref{SGdUf} and~\ref{SDsTr}, we compare the composites of the functors involved by means of the compositions of the corresponding generalized hep-derivors, which, in this section, are equal in $\mathsf{HSig}_{\gd}$, and of Lemma~\ref{LGenHepDerivorCommSquare}.

We begin with graduation and unification.

\begin{proposition}\label{PChangeSigGdUf}\lean{equivalent}{ChangeOfSignature.graduation_comp,ChangeOfSignature.unification_comp}
Let $n\in\mathbb{N}$.
Then the following equalities hold in $\mathsf{HSig}_{\gd}$
\begin{enumerate}
	\item
	$\mathbf{c}^{(n)}\circ\bGd^{(n)}_{\sig{\Sigma}}=\bGd^{(n)}_{\sig{\Lambda}}\circ\mathbf{c}^{\htg(n)}$; and

	\item
	$\mathbf{c}^{\htg(n)}\circ\bUf^{(n)}_{\sig{\Sigma}}=\bUf^{(n)}_{\sig{\Lambda}}\circ\mathbf{c}^{(n)}$.
\end{enumerate}
\end{proposition}
\begin{proof}
(1) By Definitions~\ref{DGdDerivor} and~\ref{DChangeSigDerivor}, both mappings on sorts send $(k,s)$ to $\varphi(s)$.
We compare the assignments of terms.
Let $\sigma\in\Sigma_{\mathbf{s},s}$ and $k\in n+1$.
By~(Der1) and~(Der2), $(c^{(n)})^{\sharp}(\Gd^{(n)}_{\sig{\Sigma}}(\sigma_{k}))=c^{(n)}(\sigma)=c(\sigma)$ and, by Lemma~\ref{LChangeSigTerms}~(2), $(\Gd^{(n)}_{\sig{\Lambda}})^{\sharp}(c^{\htg(n)}(\sigma_{k}))=(\Gd^{(n)}_{\sig{\Lambda}})^{\sharp}(c(\sigma)_{k})=c(\sigma)$.
For every pair $j<k$ in $n+1$ and every $s'\in S$, the following equalities hold
\begin{align*}
(c^{(n)})^{\sharp}\left(\Gd^{(n)}_{\sig{\Sigma}}(\src^{(j,k)}_{s'})\right)
&=
\src^{j}_{\varphi(s')}\left(v^{\varphi(s')}_{0}\right)
=
(\Gd^{(n)}_{\sig{\Lambda}})^{\sharp}\left(c^{\htg(n)}(\src^{(j,k)}_{s'})\right),
\\
(c^{(n)})^{\sharp}\left(\Gd^{(n)}_{\sig{\Sigma}}(\tgt^{(j,k)}_{s'})\right)
&=
\tgt^{j}_{\varphi(s')}\left(v^{\varphi(s')}_{0}\right)
=
(\Gd^{(n)}_{\sig{\Lambda}})^{\sharp}\left(c^{\htg(n)}(\tgt^{(j,k)}_{s'})\right),
\\
(c^{(n)})^{\sharp}\left(\Gd^{(n)}_{\sig{\Sigma}}(\mathrm{i}^{(k,j)}_{s'})\right)
&=
v^{\varphi(s')}_{0}
=
(\Gd^{(n)}_{\sig{\Lambda}})^{\sharp}\left(c^{\htg(n)}(\mathrm{i}^{(k,j)}_{s'})\right),
\\
(c^{(n)})^{\sharp}\left(\Gd^{(n)}_{\sig{\Sigma}}(\circ^{(j,k)}_{s'})\right)
&=
v^{\varphi(s')}_{0}\circ^{j}_{\varphi(s')}v^{\varphi(s')}_{1}
=
(\Gd^{(n)}_{\sig{\Lambda}})^{\sharp}\left(c^{\htg(n)}(\circ^{(j,k)}_{s'})\right),
\end{align*}
by Definitions~\ref{DGdDerivor} and~\ref{DChangeSigDerivor}, (Der1) and~(Der2).
Finally, by Definition~\ref{DGenHepDerivorCompId}, for every $k\in n$ and every $s'\in S$,
\begin{align*}
\Dom^{\mathbf{c}^{(n)}\circ\bGd^{(n)}_{\sig{\Sigma}}}_{(k,s')}
&=
\varnothing\cup(c^{(n)})^{\sharp}\left[\left\lbrace\left(v^{s'}_{0},\src^{k}_{s'}(v^{s'}_{0})\right)\right\rbrace\right]
=
\left\lbrace\left(v^{\varphi(s')}_{0},\src^{k}_{\varphi(s')}(v^{\varphi(s')}_{0})\right)\right\rbrace,
\\
\Dom^{\bGd^{(n)}_{\sig{\Lambda}}\circ\mathbf{c}^{\htg(n)}}_{(k,s')}
&=
\Dom^{\Gd^{(n)}_{\sig{\Lambda}}}_{(k,\varphi(s'))}\cup(\Gd^{(n)}_{\sig{\Lambda}})^{\sharp}[\varnothing]
=
\left\lbrace\left(v^{\varphi(s')}_{0},\src^{k}_{\varphi(s')}(v^{\varphi(s')}_{0})\right)\right\rbrace.
\end{align*}
For every $s'\in S$, both domains of the sort $(n,s')$ are empty, including when $n=0$.

(2) By Definitions~\ref{DUfDerivor} and~\ref{DChangeSigDerivor}, both mappings on sorts send $s$ to $(n,\varphi(s))$ and all the domains are empty.
Let $\sigma\in\Sigma_{\mathbf{s},s}$.
By~(Der1) and~(Der2), $(c^{\htg(n)})^{\sharp}(\Uf^{(n)}_{\sig{\Sigma}}(\sigma))=c^{\htg(n)}(\sigma_{n})=c(\sigma)_{n}$ and, by Lemma~\ref{LChangeSigTerms}~(2), $(\Uf^{(n)}_{\sig{\Lambda}})^{\sharp}(c^{(n)}(\sigma))=(\Uf^{(n)}_{\sig{\Lambda}})^{\sharp}(c(\sigma))=c(\sigma)_{n}$.
For every $k\in n$ and every $s'\in S$, the following equalities hold
\begin{align*}
(c^{\htg(n)})^{\sharp}\left(\Uf^{(n)}_{\sig{\Sigma}}(\src^{k}_{s'})\right)
&=
\mathrm{i}^{(n,k)}_{\varphi(s')}\left(\src^{(k,n)}_{\varphi(s')}\left(v^{(n,\varphi(s'))}_{0}\right)\right)
=
(\Uf^{(n)}_{\sig{\Lambda}})^{\sharp}\left(c^{(n)}(\src^{k}_{s'})\right),
\\
(c^{\htg(n)})^{\sharp}\left(\Uf^{(n)}_{\sig{\Sigma}}(\tgt^{k}_{s'})\right)
&=
\mathrm{i}^{(n,k)}_{\varphi(s')}\left(\tgt^{(k,n)}_{\varphi(s')}\left(v^{(n,\varphi(s'))}_{0}\right)\right)
=
(\Uf^{(n)}_{\sig{\Lambda}})^{\sharp}\left(c^{(n)}(\tgt^{k}_{s'})\right),
\\
(c^{\htg(n)})^{\sharp}\left(\Uf^{(n)}_{\sig{\Sigma}}(\circ^{k}_{s'})\right)
&=
v^{(n,\varphi(s'))}_{0}\circ^{(k,n)}_{\varphi(s')}v^{(n,\varphi(s'))}_{1}
=
(\Uf^{(n)}_{\sig{\Lambda}})^{\sharp}\left(c^{(n)}(\circ^{k}_{s'})\right),
\end{align*}
by Definitions~\ref{DUfDerivor} and~\ref{DChangeSigDerivor}, (Der1) and~(Der2). This completes the proof.
\end{proof}

From the above proposition we obtain that the functors induced by $\mathbf{c}^{(n)}$ and $\mathbf{c}^{\htg(n)}$ commute with graduation and unification.

\begin{corollary}\label{CChangeSigGdUf}\lean{weaker}{ChangeOfSignature.graduationFunctorChange,ChangeOfSignature.unificationFunctorChange}
Let $n\in\mathbb{N}$.
Then the following equalities hold
\begin{enumerate}
	\item
	$\Gd^{(n)}_{\sig{\Sigma}}\circ\mathbf{c}^{(n)\ast}_{\gd}=\mathbf{c}^{\htg(n)\ast}_{\gd}\circ\Gd^{(n)}_{\sig{\Lambda}}$; and

	\item
	$\Uf^{(n)}_{\sig{\Sigma}}\circ\mathbf{c}^{\htg(n)\ast}_{\gd}=\mathbf{c}^{(n)\ast}_{\gd}\circ\Uf^{(n)}_{\sig{\Lambda}}$.
\end{enumerate}
\end{corollary}
\begin{proof}
(1) Let $\calg{B}$ be an $n$-categorial $\sig{\Lambda}$-algebra.
We apply Lemma~\ref{LGenHepDerivorCommSquare} to the compositions $\mathbf{c}^{(n)}\circ\bGd^{(n)}_{\sig{\Sigma}}$ and $\bGd^{(n)}_{\sig{\Lambda}}\circ\mathbf{c}^{\htg(n)}$, which, by Proposition~\ref{PChangeSigGdUf}~(1), are equal.
By Proposition~\ref{PChangeSigEpVariety} and Lemma~\ref{LGdDescription}, the cells of sort $(k,s)$ of both $\Gd^{(n)}_{\sig{\Sigma}}(\mathbf{c}^{(n)\ast}_{\gd}(\calg{B}))$ and $\mathbf{c}^{\htg(n)\ast}_{\gd}(\Gd^{(n)}_{\sig{\Lambda}}(\calg{B}))$ are the $k$-cells of sort $\varphi(s)$ of $\calg{B}$, for $k\in n$, and all the cells of sort $\varphi(s)$ of $\calg{B}$, for $k=n$; by Proposition~\ref{PChangeSigGdUf}~(1), they are also the cells of $(\mathbf{c}^{(n)}\circ\bGd^{(n)}_{\sig{\Sigma}})^{\ast}_{\gd}(\calg{B})$.
Moreover, by Propositions~\ref{PGdEpVariety} and~\ref{PChangeSigEpVariety}, both partial algebras are $\sig{\Sigma}^{\htg(n)}$-algebras.
Hence, by Lemma~\ref{LGenHepDerivorCommSquare}, they are equal, and the same holds for the $n$-categorial $\sig{\Lambda}$-homomorphisms.

(2) Let $\calg{B}$ be a heterogeneous $n$-categorial $\sig{\Lambda}$-algebra.
By Proposition~\ref{PChangeSigGdUf}~(2), Lemma~\ref{LUfDescription} and Proposition~\ref{PChangeSigEpVariety}, the cells of sort $s$ of $\Uf^{(n)}_{\sig{\Sigma}}(\mathbf{c}^{\htg(n)\ast}_{\gd}(\calg{B}))$, of $\mathbf{c}^{(n)\ast}_{\gd}(\Uf^{(n)}_{\sig{\Lambda}}(\calg{B}))$ and of $(\mathbf{c}^{\htg(n)}\circ\bUf^{(n)}_{\sig{\Sigma}})^{\ast}_{\gd}(\calg{B})$ are the $n$-cells of sort $\varphi(s)$ of $\calg{B}$, and, by Propositions~\ref{PUfEpVariety} and~\ref{PChangeSigEpVariety}, the first two are $\sig{\Sigma}^{(n)}$-algebras.
Hence, by Lemma~\ref{LGenHepDerivorCommSquare}, they are equal, and the same holds for the heterogeneous $n$-categorial $\sig{\Lambda}$-homomorphisms. This completes the proof.
\end{proof}

From Corollaries~\ref{CChangeSigIndexed} and~\ref{CChangeSigGdUf} we obtain that graduation and unification are natural transformations between split indexed categories.

\begin{corollary}\label{CChangeSigGdUfNat}\lean{weaker}{ChangeOfSignature.graduationNaturality,ChangeOfSignature.unificationNaturality}
Let $n\in\mathbb{N}$.
Then the families
\begin{enumerate}
	\item 
	$\Gd^{(n)}=(\Gd^{(n)}_{\sig{\Sigma}})_{\sig{\Sigma} \in \mathrm{Ob}(\mathsf{Sig}_{\mathfrak{d}})}\colon\mathsf{nCatAlg}\Longrightarrow\mathsf{nCat}^{\htg}\mathsf{Alg}$ and
	
	\item
	$\Uf^{(n)}=(\Uf^{(n)}_{\sig{\Sigma}})_{\sig{\Sigma} \in \mathrm{Ob}(\mathsf{Sig}_{\mathfrak{d}})}\colon\mathsf{nCat}^{\htg}\mathsf{Alg}\Longrightarrow\mathsf{nCatAlg}$
\end{enumerate}
are natural transformations.
\end{corollary}
\begin{proof}
By Corollary~\ref{CChangeSigIndexed}, $\mathsf{nCatAlg}$ and $\mathsf{nCat}^{\htg}\mathsf{Alg}$ are contravariant functors from $\mathsf{Sig}_{\mathfrak{d}}$ to $\mathsf{Cat}$.
By Corollaries~\ref{CGdFunctor} and~\ref{CUfFunctor}, for every many-sorted signature $\sig{\Sigma}$, $\Gd^{(n)}_{\sig{\Sigma}}$ is a functor from $\mathsf{nCatAlg}(\sig{\Sigma})$ to $\mathsf{nCat}^{\htg}\mathsf{Alg}(\sig{\Sigma})$ and $\Uf^{(n)}_{\sig{\Sigma}}$ is a functor from $\mathsf{nCat}^{\htg}\mathsf{Alg}(\sig{\Sigma})$ to $\mathsf{nCatAlg}(\sig{\Sigma})$.
Finally, for all many-sorted signatures $\sig{\Sigma}$ and $\sig{\Lambda}$ and every derivor $\mathbf{c}$ from $\sig{\Sigma}$ to $\sig{\Lambda}$, the naturality conditions at $\mathbf{c}$, i.e., the commutativity of the diagrams
$$
\xymatrix@C=40pt@R=35pt{
\mathsf{nCatAlg}(\sig{\Lambda})\ar[r]^-{\Gd^{(n)}_{\sig{\Lambda}}}\ar[d]_-{\mathbf{c}^{(n)\ast}_{\gd}} &
\mathsf{nCat}^{\htg}\mathsf{Alg}(\sig{\Lambda})\ar[d]^-{\mathbf{c}^{\htg(n)\ast}_{\gd}}\\
\mathsf{nCatAlg}(\sig{\Sigma})\ar[r]_-{\Gd^{(n)}_{\sig{\Sigma}}} &
\mathsf{nCat}^{\htg}\mathsf{Alg}(\sig{\Sigma})
}
\qquad
\xymatrix@C=40pt@R=35pt{
\mathsf{nCat}^{\htg}\mathsf{Alg}(\sig{\Lambda})\ar[r]^-{\Uf^{(n)}_{\sig{\Lambda}}}\ar[d]_-{\mathbf{c}^{\htg(n)\ast}_{\gd}} &
\mathsf{nCatAlg}(\sig{\Lambda})\ar[d]^-{\mathbf{c}^{(n)\ast}_{\gd}}\\
\mathsf{nCat}^{\htg}\mathsf{Alg}(\sig{\Sigma})\ar[r]_-{\Uf^{(n)}_{\sig{\Sigma}}} &
\mathsf{nCatAlg}(\sig{\Sigma})
}
$$
are the equalities of Corollary~\ref{CChangeSigGdUf}.
\end{proof}

We next consider discretization and truncation.

\begin{proposition}\label{PChangeSigDsTr}\lean{equivalent}{ChangeOfSignature.homogeneousDiscretization_comp,ChangeOfSignature.homogeneousTruncation_comp,ChangeOfSignature.heterogeneousDimension_comp}
Let $m,n\in\mathbb{N}$ with $m\leq n$.
Then the following equalities hold in $\mathsf{HSig}_{\gd}$
\begin{align*}
\mathbf{c}^{(m)}\circ\bDs^{(n,m)}_{\sig{\Sigma}}&=\bDs^{(n,m)}_{\sig{\Lambda}}\circ\mathbf{c}^{(n)},
\\
\mathbf{c}^{(n)}\circ\bTr^{(m,n)}_{\sig{\Sigma}}&=\bTr^{(m,n)}_{\sig{\Lambda}}\circ\mathbf{c}^{(m)},
\\
\mathbf{c}^{\htg(m)}\circ\bDs^{\htg(n,m)}_{\sig{\Sigma}}&=\bDs^{\htg(n,m)}_{\sig{\Lambda}}\circ\mathbf{c}^{\htg(n)},
\\
\mathbf{c}^{\htg(n)}\circ\bTr^{\htg(m,n)}_{\sig{\Sigma}}&=\bTr^{\htg(m,n)}_{\sig{\Lambda}}\circ\mathbf{c}^{\htg(m)}.
\end{align*}
\end{proposition}
\begin{proof}
If $m=n$, all the discretization and truncation derivors involved are identities, and the statements follow from Proposition~\ref{PHepSigCat}.
Assume that $m<n$.
In the four cases both mappings on sorts are equal, see Definitions~\ref{DDsDerivor}, \ref{DTrDerivor} and~\ref{DChangeSigDerivor}: in the two homogeneous cases, both are $\varphi$; in the heterogeneous discretization, both send $(k,s)$ to $(\min\{k,m\},\varphi(s))$; and, in the heterogeneous truncation, both send $(k,s)$ to $(k,\varphi(s))$.
We compare the assignments of terms and the domains.

We first consider the homogeneous discretization.
Let $\sigma$ be an operation symbol of $\Sigma$.
The following chain of equalities holds
\begin{align*}
(c^{(m)})^{\sharp}\left(\Ds^{(n,m)}_{\sig{\Sigma}}(\sigma)\right)
&=
c(\sigma)
\tag{1}
\\&=
(\Ds^{(n,m)}_{\sig{\Lambda}})^{\sharp}\left(c(\sigma)\right)
\tag{2}
\\&=
(\Ds^{(n,m)}_{\sig{\Lambda}})^{\sharp}\left(c^{(n)}(\sigma)\right).
\tag{3}
\end{align*}
The first equality follows from~(Der1), (Der2) and Definitions~\ref{DDsDerivor}~(1) and~\ref{DChangeSigDerivor}; the second equality follows from Lemma~\ref{LChangeSigTerms}~(3); finally, the last equality follows from Definition~\ref{DChangeSigDerivor}.
For the structural symbols, fix $k\in n$ and $s\in S$.
If $k<m$, both generalized hep-derivors send $\src^{k}_{s}$, $\tgt^{k}_{s}$ and $\circ^{k}_{s}$ to $\src^{k}_{\varphi(s)}(v^{\varphi(s)}_{0})$, $\tgt^{k}_{\varphi(s)}(v^{\varphi(s)}_{0})$ and $v^{\varphi(s)}_{0}\circ^{k}_{\varphi(s)}v^{\varphi(s)}_{1}$, respectively; and, if $k\geq m$, both send them to $v^{\varphi(s)}_{0}$.
All the domains are empty.

Now we consider the homogeneous truncation.
Let $\tau$ be an operation symbol of $\Sigma^{(m)}$.
The following chain of equalities holds
\begin{align*}
(c^{(n)})^{\sharp}\left(\Tr^{(m,n)}_{\sig{\Sigma}}(\tau)\right)
&=
c^{(n)}(\tau)
\tag{4}
\\&=
c^{(m)}(\tau)
\tag{5}
\\&=
(\Tr^{(m,n)}_{\sig{\Lambda}})^{\sharp}\left(c^{(m)}(\tau)\right),
\tag{6}
\end{align*}
where $c^{(m)}(\tau)$ is regarded as a term of $\Lambda^{(n)}$.
The first equality follows from~(Der1), (Der2) and Definition~\ref{DTrDerivor}; the second equality follows from Definition~\ref{DChangeSigDerivor}, since $c^{(n)}$ and $c^{(m)}$ agree on the operation symbols of $\Sigma^{(m)}$; finally, the last equality follows from~(Der1), (Der2) and Definition~\ref{DTrDerivor}, since the assignment of terms of $\bTr^{(m,n)}_{\sig{\Lambda}}$ is the inclusion, so that $(\Tr^{(m,n)}_{\sig{\Lambda}})^{\sharp}$ is the inclusion of the terms of $\Lambda^{(m)}$ into those of $\Lambda^{(n)}$.
Moreover, by Definition~\ref{DGenHepDerivorCompId}, for every $s\in S$,
\begin{align*}
\Dom^{\mathbf{c}^{(n)}\circ\bTr^{(m,n)}_{\sig{\Sigma}}}_{s}
&=
\varnothing\cup(c^{(n)})^{\sharp}\left[\left\lbrace\left(v^{s}_{0},\src^{m}_{s}(v^{s}_{0})\right)\right\rbrace\right]
=
\left\lbrace\left(v^{\varphi(s)}_{0},\src^{m}_{\varphi(s)}(v^{\varphi(s)}_{0})\right)\right\rbrace,
\\
\Dom^{\bTr^{(m,n)}_{\sig{\Lambda}}\circ\mathbf{c}^{(m)}}_{s}
&=
\Dom^{\Tr^{(m,n)}_{\sig{\Lambda}}}_{\varphi(s)}\cup(\Tr^{(m,n)}_{\sig{\Lambda}})^{\sharp}[\varnothing]
=
\left\lbrace\left(v^{\varphi(s)}_{0},\src^{m}_{\varphi(s)}(v^{\varphi(s)}_{0})\right)\right\rbrace.
\end{align*}

Next, we consider the heterogeneous discretization.
Let $\sigma$ be an operation symbol of $\Sigma$ and $k\in n+1$.
The following chain of equalities holds
\begin{align*}
(c^{\htg(m)})^{\sharp}\left(\Ds^{\htg(n,m)}_{\sig{\Sigma}}(\sigma_{k})\right)
&=
c(\sigma)_{\min\{k,m\}}
\tag{7}
\\&=
(\Ds^{\htg(n,m)}_{\sig{\Lambda}})^{\sharp}\left(c(\sigma)_{k}\right)
\tag{8}
\\&=
(\Ds^{\htg(n,m)}_{\sig{\Lambda}})^{\sharp}\left(c^{\htg(n)}(\sigma_{k})\right).
\tag{9}
\end{align*}
The first equality follows from~(Der1), (Der2) and Definitions~\ref{DDsDerivor}~(2) and~\ref{DChangeSigDerivor}; the second equality follows from Lemma~\ref{LChangeSigTerms}~(4); finally, the last equality follows from Definition~\ref{DChangeSigDerivor}.
For the structural symbols, fix $j<k$ in $n+1$ and $s\in S$.
If $j<m$, both generalized hep-derivors send $\src^{(j,k)}_{s}$, $\tgt^{(j,k)}_{s}$, $\mathrm{i}^{(k,j)}_{s}$ and $\circ^{(j,k)}_{s}$ to the corresponding operation symbols of sort $\varphi(s)$ and dimensions $j$ and $\min\{k,m\}$, applied to variables; and, if $j\geq m$, both send them to $v^{(m,\varphi(s))}_{0}$.
All the domains are empty.

Finally, we consider the heterogeneous truncation.
Let $\tau$ be an operation symbol of $\Sigma^{\htg(m)}$.
As in the homogeneous case, by~(Der1), (Der2) and Definition~\ref{DTrDerivor}, $(\Tr^{\htg(m,n)}_{\sig{\Lambda}})^{\sharp}$ is the inclusion of the terms of $\Lambda^{\htg(m)}$ into those of $\Lambda^{\htg(n)}$ determined by the canonical embedding, and, by Definition~\ref{DChangeSigDerivor}, $c^{\htg(n)}$ and $c^{\htg(m)}$ agree on the operation symbols of $\Sigma^{\htg(m)}$.
Hence both generalized hep-derivors send $\tau$ to $c^{\htg(m)}(\tau)$, regarded as a term of $\Lambda^{\htg(n)}$.
All the domains are empty. This completes the proof.
\end{proof}

From the above proposition we obtain that the functors induced by $\mathbf{c}^{(n)}$ and $\mathbf{c}^{\htg(n)}$ commute with discretization and truncation.

\begin{corollary}\label{CChangeSigDsTr}\lean{weaker}{ChangeOfSignature.homogeneousDiscretizationFunctorChange,ChangeOfSignature.homogeneousTruncationFunctorChange,ChangeOfSignature.heterogeneousDiscretizationFunctorChange,ChangeOfSignature.heterogeneousTruncationFunctorChange}
Let $m,n\in\mathbb{N}$ with $m\leq n$.
The following equalities hold
\begin{enumerate}
	\item
	$\Ds^{(n,m)}_{\sig{\Sigma}}\circ\mathbf{c}^{(m)\ast}_{\gd}=\mathbf{c}^{(n)\ast}_{\gd}\circ\Ds^{(n,m)}_{\sig{\Lambda}}$;

	\item
	$\Tr^{(m,n)}_{\sig{\Sigma}}\circ\mathbf{c}^{(n)\ast}_{\gd}=\mathbf{c}^{(m)\ast}_{\gd}\circ\Tr^{(m,n)}_{\sig{\Lambda}}$;

	\item
	$\Ds^{\htg(n,m)}_{\sig{\Sigma}}\circ\mathbf{c}^{\htg(m)\ast}_{\gd}=\mathbf{c}^{\htg(n)\ast}_{\gd}\circ\Ds^{\htg(n,m)}_{\sig{\Lambda}}$; and

	\item
	$\Tr^{\htg(m,n)}_{\sig{\Sigma}}\circ\mathbf{c}^{\htg(n)\ast}_{\gd}=\mathbf{c}^{\htg(m)\ast}_{\gd}\circ\Tr^{\htg(m,n)}_{\sig{\Lambda}}$.
\end{enumerate}
\end{corollary}
\begin{proof}
If $m=n$, all the discretization and truncation functors involved are identities, by Corollary~\ref{CGenHepDerivorIdFunctor}.
Assume that $m<n$.
In the four cases we apply Lemma~\ref{LGenHepDerivorCommSquare} by means of Proposition~\ref{PChangeSigDsTr}.
In each case, let $\calg{B}$ be an algebra in the common domain of the two composites, i.e., an $m$-categorial, $n$-categorial, heterogeneous $m$-categorial or heterogeneous $n$-categorial $\sig{\Lambda}$-algebra, respectively.

We first consider the homogeneous discretization.
By Lemma~\ref{LDsDescription} and Proposition~\ref{PChangeSigEpVariety}, the cells of sort $s$ of $\Ds^{(n,m)}_{\sig{\Sigma}}(\mathbf{c}^{(m)\ast}_{\gd}(\calg{B}))$, of $\mathbf{c}^{(n)\ast}_{\gd}(\Ds^{(n,m)}_{\sig{\Lambda}}(\calg{B}))$ and of the transformed algebra along the composite are all the cells of sort $\varphi(s)$ of $\calg{B}$; and, by Propositions~\ref{PDsEpVariety} and~\ref{PChangeSigEpVariety}, the first two are $\sig{\Sigma}^{(n)}$-algebras.

Now we consider the homogeneous truncation.
By the note after Definition~\ref{DTrDerivor} and Proposition~\ref{PChangeSigEpVariety}, the cells of sort $s$ of $\Tr^{(m,n)}_{\sig{\Sigma}}(\mathbf{c}^{(n)\ast}_{\gd}(\calg{B}))$, of $\mathbf{c}^{(m)\ast}_{\gd}(\Tr^{(m,n)}_{\sig{\Lambda}}(\calg{B}))$ and, by the computation of the domains in the proof of Proposition~\ref{PChangeSigDsTr}, of the transformed algebra along the composite are the $m$-cells of sort $\varphi(s)$ of $\calg{B}$; and, by Propositions~\ref{PTrEpVariety} and~\ref{PChangeSigEpVariety}, the first two are $\sig{\Sigma}^{(m)}$-algebras.

Next, we consider the heterogeneous discretization.
By Lemma~\ref{LDsDescription} and Proposition~\ref{PChangeSigEpVariety}, the cells of sort $(k,s)$ of the three algebras are the elements of $B_{(\min\{k,m\},\varphi(s))}$; and, by Propositions~\ref{PDsEpVariety} and~\ref{PChangeSigEpVariety}, the first two are $\sig{\Sigma}^{\htg(n)}$-algebras.

Finally, we consider the heterogeneous truncation.
By the note after Definition~\ref{DTrDerivor} and Proposition~\ref{PChangeSigEpVariety}, the cells of sort $(k,s)$ of the three algebras, for $k\in m+1$, are the elements of $B_{(k,\varphi(s))}$; and, by Propositions~\ref{PTrEpVariety} and~\ref{PChangeSigEpVariety}, the first two are $\sig{\Sigma}^{\htg(m)}$-algebras.

In the four cases, by Lemma~\ref{LGenHepDerivorCommSquare}, both algebras are equal, and the same holds for the homomorphisms. This completes the proof.
\end{proof}

From Corollaries~\ref{CChangeSigIndexed} and~\ref{CChangeSigDsTr} we obtain that discretization and truncation are natural transformations between split indexed categories.

\begin{corollary}\label{CChangeSigDsTrNat}\lean{weaker}{ChangeOfSignature.homogeneousDiscretizationNaturality,ChangeOfSignature.homogeneousTruncationNaturality,ChangeOfSignature.heterogeneousDiscretizationNaturality,ChangeOfSignature.heterogeneousTruncationNaturality}
Let $m,n\in\mathbb{N}$ with $m\leq n$.
Then the families
\begin{enumerate}
	\item 
	$\Ds^{(n,m)}=(\Ds^{(n,m)}_{\sig{\Sigma}})_{\sig{\Sigma} \in \mathrm{Ob}(\mathsf{Sig}_{\mathfrak{d}})}\colon\mathsf{mCatAlg}\Longrightarrow\mathsf{nCatAlg}$;
	
	\item
	$\Tr^{(m,n)}=(\Tr^{(m,n)}_{\sig{\Sigma}})_{\sig{\Sigma} \in \mathrm{Ob}(\mathsf{Sig}_{\mathfrak{d}})}\colon\mathsf{nCatAlg}\Longrightarrow\mathsf{mCatAlg}$;
	
	\item 
	$\Ds^{\htg(n,m)}=(\Ds^{\htg(n,m)}_{\sig{\Sigma}})_{\sig{\Sigma} \in \mathrm{Ob}(\mathsf{Sig}_{\mathfrak{d}})}\colon\mathsf{mCat}^{\htg}\mathsf{Alg}\Longrightarrow\mathsf{nCat}^{\htg}\mathsf{Alg}$; and
	
	\item 
	$\Tr^{\htg(m,n)}=(\Tr^{\htg(m,n)}_{\sig{\Sigma}})_{\sig{\Sigma} \in \mathrm{Ob}(\mathsf{Sig}_{\mathfrak{d}})}\colon\mathsf{nCat}^{\htg}\mathsf{Alg}\Longrightarrow\mathsf{mCat}^{\htg}\mathsf{Alg}$
\end{enumerate}
are natural transformations.
\end{corollary}
\begin{proof}
By Corollary~\ref{CChangeSigIndexed}, the four assignments $\mathsf{mCatAlg}$, $\mathsf{nCatAlg}$, $\mathsf{mCat}^{\htg}\mathsf{Alg}$ and $\mathsf{nCat}^{\htg}\mathsf{Alg}$ are contravariant functors from $\mathsf{Sig}_{\mathfrak{d}}$ to $\mathsf{Cat}$.
By Corollaries~\ref{CDsFunctor} and~\ref{CTrFunctor}, for every many-sorted signature $\sig{\Sigma}$, the components at $\sig{\Sigma}$ of the four families are functors between the corresponding categories.
Finally, for all many-sorted signatures $\sig{\Sigma}$ and $\sig{\Lambda}$ and every derivor $\mathbf{c}$ from $\sig{\Sigma}$ to $\sig{\Lambda}$, the naturality conditions at $\mathbf{c}$ are the equalities of Corollary~\ref{CChangeSigDsTr}.
\end{proof}

In the following remark we show that, as natural transformations between split indexed categories, the presentation and dimension functors also act on algebras of different signatures, by means of the Grothendieck construction, see Definition~\ref{DGrothendieck} and Remark~\ref{RGrothendieckNat}.

\begin{remark}\label{RChangeSigGrothendieck}
Let $m,n\in\mathbb{N}$ with $m\leq n$.
By Corollaries~\ref{CChangeSigGdUfNat} and~\ref{CChangeSigDsTrNat}, the families $\Gd^{(n)}$, $\Uf^{(n)}$, $\Ds^{(n,m)}$, $\Tr^{(m,n)}$, $\Ds^{\htg(n,m)}$ and $\Tr^{\htg(m,n)}$ are natural transformations between split $\mathsf{Sig}_{\mathfrak{d}}$-indexed categories.
Hence, by Remark~\ref{RGrothendieckNat}, they extend to functors between the Grothendieck constructions of the split indexed categories of Corollary~\ref{CChangeSigIndexed} that commute with the canonical split fibrations.
For instance, graduation extends to the functor from $\int^{\mathsf{Sig}_{\mathfrak{d}}}\mathsf{nCatAlg}$ to $\int^{\mathsf{Sig}_{\mathfrak{d}}}\mathsf{nCat}^{\htg}\mathsf{Alg}$ that sends an object $(\sig{\Sigma},\calg{A})$ to $(\sig{\Sigma},\Gd^{(n)}_{\sig{\Sigma}}(\calg{A}))$ and a morphism $(\mathbf{c},F)$ from $(\sig{\Sigma},\calg{A})$ to $(\sig{\Lambda},\calg{B})$ to $(\mathbf{c},\Gd^{(n)}_{\sig{\Sigma}}(F))$.
Here $F\colon\calg{A}\mor\mathbf{c}^{(n)\ast}_{\gd}(\calg{B})$.
By Corollary~\ref{CChangeSigGdUf}, $\Gd^{(n)}_{\sig{\Sigma}}(F)$ is a homomorphism from $\Gd^{(n)}_{\sig{\Sigma}}(\calg{A})$ to $\Gd^{(n)}_{\sig{\Sigma}}(\mathbf{c}^{(n)\ast}_{\gd}(\calg{B}))=\mathbf{c}^{\htg(n)\ast}_{\gd}(\Gd^{(n)}_{\sig{\Lambda}}(\calg{B}))$, so that $(\mathbf{c},\Gd^{(n)}_{\sig{\Sigma}}(F))$ is a morphism from $(\sig{\Sigma},\Gd^{(n)}_{\sig{\Sigma}}(\calg{A}))$ to $(\sig{\Lambda},\Gd^{(n)}_{\sig{\Lambda}}(\calg{B}))$.
\end{remark}

The natural transformations of Chapters~\ref{SGdUf} and~\ref{SDsTr} are also compatible with the change of signature.

\begin{remark}\label{RChangeSigNat}
Let $m,n\in\mathbb{N}$ with $m\leq n$.
The natural isomorphisms $\theta$ of Proposition~\ref{PGdUfAdj} and $\beta^{(m,n)}$ of Proposition~\ref{PCompatUfTr}, and the counits $\varepsilon$ and $\varepsilon^{\htg}$ of Lemma~\ref{LDsTrCounit}, are given, at every algebra, by families of identities, inclusions and $(k,j)$-identities indexed by the sorts.
By Proposition~\ref{PChangeSigEpVariety}, the functors $\mathbf{c}^{(n)\ast}_{\gd}$ and $\mathbf{c}^{\htg(n)\ast}_{\gd}$ reindex the sorts along $\varphi$ and keep the cells, the sources, the targets, the identities and the compositions.
Therefore, by Corollaries~\ref{CChangeSigGdUf} and~\ref{CChangeSigDsTr}, for every heterogeneous $n$-categorial $\sig{\Lambda}$-algebra $\calg{B}$,
\begin{align*}
\theta^{\sig{\Sigma}}_{\mathbf{c}^{\htg(n)\ast}_{\gd}(\calg{B})}&=\mathbf{c}^{\htg(n)\ast}_{\gd}\left(\theta^{\sig{\Lambda}}_{\calg{B}}\right),
&
\beta^{(m,n),\sig{\Sigma}}_{\mathbf{c}^{\htg(n)\ast}_{\gd}(\calg{B})}&=\mathbf{c}^{(m)\ast}_{\gd}\left(\beta^{(m,n),\sig{\Lambda}}_{\calg{B}}\right),
\end{align*}
where the superscripts $\sig{\Sigma}$ and $\sig{\Lambda}$ indicate the signature, and, similarly, $\varepsilon^{\sig{\Sigma}}_{\mathbf{c}^{(n)\ast}_{\gd}(\calg{B})}=\mathbf{c}^{(n)\ast}_{\gd}(\varepsilon^{\sig{\Lambda}}_{\calg{B}})$, for every $n$-categorial $\sig{\Lambda}$-algebra $\calg{B}$, and $\varepsilon^{\htg,\sig{\Sigma}}_{\mathbf{c}^{\htg(n)\ast}_{\gd}(\calg{B})}=\mathbf{c}^{\htg(n)\ast}_{\gd}(\varepsilon^{\htg,\sig{\Lambda}}_{\calg{B}})$, for every heterogeneous $n$-categorial $\sig{\Lambda}$-algebra $\calg{B}$.
For instance, $\beta^{(m,n),\sig{\Sigma}}_{\mathbf{c}^{\htg(n)\ast}_{\gd}(\calg{B})}$ is the family $(\mathrm{i}^{(n,m)\,\calg{B}}_{\varphi(s)})_{s\in S}$ of $(n,m)$-identities of $\calg{B}$, which is $(\beta^{(m,n),\sig{\Lambda}}_{\calg{B}})_{\varphi}$; and both homomorphisms have the same domain and codomain, since, by Corollaries~\ref{CChangeSigGdUf} and~\ref{CChangeSigDsTr},
\begin{enumerate}
	\item
	$ \Uf^{(m)}_{\sig{\Sigma}}\circ\Tr^{\htg(m,n)}_{\sig{\Sigma}}\circ\mathbf{c}^{\htg(n)\ast}_{\gd}=\mathbf{c}^{(m)\ast}_{\gd}\circ\Uf^{(m)}_{\sig{\Lambda}}\circ\Tr^{\htg(m,n)}_{\sig{\Lambda}} $; and

	\item
	$ \Tr^{(m,n)}_{\sig{\Sigma}}\circ\Uf^{(n)}_{\sig{\Sigma}}\circ\mathbf{c}^{\htg(n)\ast}_{\gd}=\mathbf{c}^{(m)\ast}_{\gd}\circ\Tr^{(m,n)}_{\sig{\Lambda}}\circ\Uf^{(n)}_{\sig{\Lambda}} $.
\end{enumerate}
Thus the adjunctions of Propositions~\ref{PGdUfAdj} and~\ref{PDsTrAdj}, and the algebraic square of Proposition~\ref{PAlgSquareFinite}, are compatible with the change of signature.
\end{remark}

Summing up, in every finite dimension the functors induced by $\mathbf{c}^{(n)}$ and $\mathbf{c}^{\htg(n)}$ commute with graduation, unification, discretization and truncation, and the natural transformations, adjunctions and algebraic square that relate these functors are compatible with the change of signature.
In the next section we show that the same holds in dimension $\omega$.

\section{Change of signature in dimension \texorpdfstring{$\omega$}{omega}}\label{SChangeSigOmega}

In this section we extend the results of Sections~\ref{SChangeSigDerivors} and~\ref{SChangeSigCompat} to the dimension $\omega$.
Every operation symbol and every axiom of the $\omega$-dimensional specifications involves only finitely many dimensions, so that, as in Section~\ref{SDsTrOmega}, the proofs reduce to the finite dimensions.
Graduation and unification in dimension $\omega$ are treated by means of the universal property of Proposition~\ref{PProjLimit} and of the explicit description of Definition~\ref{DUfFunctorOmega}, respectively.

We begin by extending the generalized hep-derivors of Definition~\ref{DChangeSigDerivor} to the dimension $\omega$: they are given by the same formulas, now for all the dimensions in $\omega$, and they translate the operation symbols of $\Sigma$ by means of $c$ and leave the categorial structure unchanged.

\begin{definition}\label{DChangeSigDerivorOmega}\lean{equivalent}{ChangeOfSignature.homogeneous,ChangeOfSignature.heterogeneous}
\begin{enumerate}
\item The \emph{homogeneous $\omega$-categorial derivor determined by $\mathbf{c}$} is the generalized hep-derivor
$$
\mathbf{c}^{(\omega)}=\left(\varphi,c^{(\omega)},\varnothing\right)
$$
from $\sig{\Sigma}^{(\omega)}$ to $\sig{\Lambda}^{(\omega)}$, with empty domains, whose mapping on sorts is $\varphi$ and whose assignment of terms is given by the formulas of Definition~\ref{DChangeSigDerivor}~(1), for every $k\in\omega$ and every $s\in S$.
\item The \emph{heterogeneous $\omega$-categorial derivor determined by $\mathbf{c}$} is the generalized hep-derivor
$$
\mathbf{c}^{\htg(\omega)}=\left(\varphi^{\htg(\omega)},c^{\htg(\omega)},\varnothing\right)
$$
from $\sig{\Sigma}^{\htg(\omega)}$ to $\sig{\Lambda}^{\htg(\omega)}$, with empty domains, whose mapping on sorts sends $(k,s)$ to $\varphi^{\htg(\omega)}(k,s)=(k,\varphi(s))$ and whose assignment of terms is given by the formulas of Definition~\ref{DChangeSigDerivor}~(2), for every $k\in\omega$, every pair $j<k$ in $\omega$ and every $s\in S$.
\end{enumerate}
\end{definition}

Let us note that, for every $p\in\mathbb{N}$, the restriction of the assignment of terms of $\mathbf{c}^{(\omega)}$ to $\Sigma^{(p)}\subseteq\Sigma^{(\omega)}$ is that of $\mathbf{c}^{(p)}$, and the restriction of the assignment of terms of $\mathbf{c}^{\htg(\omega)}$ to the image of $\Sigma^{\htg(p)}$ in $\Sigma^{\htg(\omega)}$ is that of $\mathbf{c}^{\htg(p)}$, since the defining formulas coincide; moreover, all the domains are empty.
Therefore, by Proposition~\ref{PTransformedAlg}, for every partial $\Lambda^{(\omega)}$-algebra $\mathbf{B}$, the reduct of $\mathbf{c}^{(\omega)\ast}_{\gd}(\mathbf{B})$ to $\Sigma^{(p)}$ is $\mathbf{c}^{(p)\ast}_{\gd}$ applied to the reduct of $\mathbf{B}$ to $\Lambda^{(p)}$; and, for every partial $\Lambda^{\htg(\omega)}$-algebra $\mathbf{B}$, the reduct of $\mathbf{c}^{\htg(\omega)\ast}_{\gd}(\mathbf{B})$ to the sorts of dimension at most $p$ and to the operation symbols of $\Sigma^{\htg(p)}$ is $\mathbf{c}^{\htg(p)\ast}_{\gd}$ applied to the corresponding reduct of $\mathbf{B}$.
In fact, in both algebras the cells are the same and every operation symbol is interpreted by means of the same term and the same domain conditions.

In the following proposition we prove that these generalized hep-derivors send categorial algebras to categorial algebras.

\begin{proposition}\label{PChangeSigEpVarietyOmega}\lean{weaker}{ChangeOfSignature.homogeneous_finiteCells,ChangeOfSignature.homogeneousOmegaFunctor,ChangeOfSignature.heterogeneous_isModel,ChangeOfSignature.homogeneousRealization,ChangeOfSignature.heterogeneousRealization}
The following statements hold:
\begin{enumerate}
\item the restriction of $\mathbf{c}^{(\omega)\ast}_{\gd}$ to $\omega\mathsf{CatAlg}(\sig{\Lambda})$ corestricts to $\omega\mathsf{CatAlg}(\sig{\Sigma})$; moreover, for every $\omega$-categorial $\sig{\Lambda}$-algebra $\calg{B}$, for every $s\in S$, the cells of sort $s$ of $\mathbf{c}^{(\omega)\ast}_{\gd}(\calg{B})$ are the elements of $B_{\varphi(s)}$, for every $(\mathbf{s},s)\in S^{\star}\times S$, every $\sigma\in\Sigma_{\mathbf{s},s}$ acts as $c(\sigma)^{\calg{B}}$ and, for every $k\in\omega$ and every $s\in S$, $\src^{k}_{s}$, $\tgt^{k}_{s}$ and $\circ^{k}_{s}$ act as $\src^{k\,\calg{B}}_{\varphi(s)}$, $\tgt^{k\,\calg{B}}_{\varphi(s)}$ and $\circ^{k\,\calg{B}}_{\varphi(s)}$, respectively; and
\item the restriction of $\mathbf{c}^{\htg(\omega)\ast}_{\gd}$ to $\omega\mathsf{Cat}^{\htg}\mathsf{Alg}(\sig{\Lambda})$ corestricts to $\omega\mathsf{Cat}^{\htg}\mathsf{Alg}(\sig{\Sigma})$; moreover, for every heterogeneous $\omega$-categorial $\sig{\Lambda}$-algebra $\calg{B}$, for every $k\in\omega$ and every $s\in S$, the cells of sort $(k,s)$ of $\mathbf{c}^{\htg(\omega)\ast}_{\gd}(\calg{B})$ are the elements of $B_{(k,\varphi(s))}$, for every $(\mathbf{s},s)\in S^{\star}\times S$, every $\sigma\in\Sigma_{\mathbf{s},s}$ and every $k\in\omega$, $\sigma_{k}$ acts as $c(\sigma)_{k}^{\calg{B}}$ and, for every pair of dimensions $j<k$ in $\omega$ and every $s\in S$, $\src^{(j,k)}_{s}$, $\tgt^{(j,k)}_{s}$, $\mathrm{i}^{(k,j)}_{s}$ and $\circ^{(j,k)}_{s}$ act as $\src^{(j,k)\,\calg{B}}_{\varphi(s)}$, $\tgt^{(j,k)\,\calg{B}}_{\varphi(s)}$, $\mathrm{i}^{(k,j)\,\calg{B}}_{\varphi(s)}$ and $\circ^{(j,k)\,\calg{B}}_{\varphi(s)}$, respectively.
\end{enumerate}
\end{proposition}
\begin{proof}
(1) Let $\calg{B}$ be an $\omega$-categorial $\sig{\Lambda}$-algebra and let $\mathbf{A}$ stand for $\mathbf{c}^{(\omega)\ast}_{\gd}(\calg{B})$.
Since the domains of $\mathbf{c}^{(\omega)}$ are empty, by Remark~\ref{RDomainSetEmpty}, the cells of sort $s$ of $\mathbf{A}$ are the elements of $B_{\varphi(s)}$.
Let $p\in\mathbb{N}$.
By Remark~\ref{RCatAlgAsEpAlgOmega}, $\calg{B}$ is an object in $\mathsf{EPAlg}(\spec{C}^{\sig{\Lambda}^{(\omega)}})$; hence, as in the proof of Proposition~\ref{PTrEpVarietyOmega}~(1), its reduct to $\Lambda^{(p)}$ is an object in $\mathsf{EPAlg}(\spec{C}^{\sig{\Lambda}^{(p)}})$, since its operations and domain conditions are those of $\calg{B}$ and the axioms of $\spec{C}^{\sig{\Lambda}^{(p)}}$ are among those of $\spec{C}^{\sig{\Lambda}^{(\omega)}}$.
Hence, by the above note and Proposition~\ref{PChangeSigEpVariety}~(1), the reduct of $\mathbf{A}$ to $\Sigma^{(p)}$ is an object in $\mathsf{EPAlg}(\spec{C}^{\sig{\Sigma}^{(p)}})$, in which the operation symbols of $\Sigma^{(p)}$ act as stated.

Now we prove that $\mathbf{A}$ is an object in $\mathsf{EPAlg}(\spec{C}^{\sig{\Sigma}^{(\omega)}})$.
Given an operation symbol of $\Sigma^{(\omega)}$, choose $p\in\mathbb{N}$ such that it is an operation symbol of $\Sigma^{(p)}$.
The operation and its domain conditions are the same in $\mathbf{A}$ and in its reduct to $\Sigma^{(p)}$, so that it acts as stated and it is defined exactly on the families that satisfy its domain conditions.
Every axiom of $\spec{C}^{\sig{\Sigma}^{(\omega)}}$ involves only finitely many dimensions; choose $p\in\mathbb{N}$ greater than all of them.
The axiom, including its context, belongs to $\spec{C}^{\sig{\Sigma}^{(p)}}$, and its terms have the same definedness and values in $\mathbf{A}$ and in its reduct to $\Sigma^{(p)}$, in which it holds; hence it holds in $\mathbf{A}$.

Finally, we prove that every cell of $\mathbf{A}$ has a finite dimension.
Let $x$ be a cell of sort $s$ of $\mathbf{A}$.
Then $x$ is a cell of sort $\varphi(s)$ of $\calg{B}$, hence a $k$-cell of $\calg{B}$ for some $k\in\omega$, i.e., $\src^{k\,\calg{B}}_{\varphi(s)}(x)=x$.
Since $\src^{k}_{s}$ acts in $\mathbf{A}$ as $\src^{k\,\calg{B}}_{\varphi(s)}$, $x$ is a $k$-cell of $\mathbf{A}$.
Therefore, by Remark~\ref{RCatAlgAsEpAlgOmega}, $\mathbf{A}$ is an $\omega$-categorial $\sig{\Sigma}$-algebra, and, since $\omega\mathsf{CatAlg}(\sig{\Sigma})$ is a full subcategory of $\mathsf{EPAlg}(\spec{C}^{\sig{\Sigma}^{(\omega)}})$, the restriction corestricts as stated.

(2) Let $\calg{B}$ be a heterogeneous $\omega$-categorial $\sig{\Lambda}$-algebra and let $\mathbf{A}$ stand for $\mathbf{c}^{\htg(\omega)\ast}_{\gd}(\calg{B})$.
Since the domains of $\mathbf{c}^{\htg(\omega)}$ are empty, the cells of sort $(k,s)$ of $\mathbf{A}$ are the elements of $B_{(k,\varphi(s))}$.
Let $p\in\mathbb{N}$.
By the proof of Proposition~\ref{PTrEpVarietyOmega}~(2), the reduct of $\calg{B}$ to the sorts of dimension at most $p$ and to the operation symbols of $\Lambda^{\htg(p)}$ is $\Tr^{\htg(p,\omega)}_{\sig{\Lambda}}(\calg{B})$, which is a heterogeneous $p$-categorial $\sig{\Lambda}$-algebra.
Hence, by the above note and Proposition~\ref{PChangeSigEpVariety}~(2), the corresponding reduct of $\mathbf{A}$ is $\mathbf{c}^{\htg(p)\ast}_{\gd}(\Tr^{\htg(p,\omega)}_{\sig{\Lambda}}(\calg{B}))$, an object in $\mathsf{EPAlg}(\spec{C}^{\sig{\Sigma}^{\htg(p)}})$ in which the operation symbols of $\Sigma^{\htg(p)}$ act as stated.
As in~(1), since every operation symbol and every axiom of $\spec{C}^{\sig{\Sigma}^{\htg(\omega)}}$ involves only finitely many dimensions, $\mathbf{A}$ is an object in $\mathsf{EPAlg}(\spec{C}^{\sig{\Sigma}^{\htg(\omega)}})=\omega\mathsf{Cat}^{\htg}\mathsf{Alg}(\sig{\Sigma})$, in which the operation symbols act as stated. This completes the proof.
\end{proof}

By abuse of notation, we denote the corestrictions of Proposition~\ref{PChangeSigEpVarietyOmega} by
\begin{enumerate}
	\item
	$ \mathbf{c}^{(\omega)\ast}_{\gd}\colon\omega\mathsf{CatAlg}(\sig{\Lambda})\mor\omega\mathsf{CatAlg}(\sig{\Sigma}) $; and

	\item
	$ \mathbf{c}^{\htg(\omega)\ast}_{\gd}\colon\omega\mathsf{Cat}^{\htg}\mathsf{Alg}(\sig{\Lambda})\mor\omega\mathsf{Cat}^{\htg}\mathsf{Alg}(\sig{\Sigma}) $.
\end{enumerate}
By Proposition~\ref{PGenHepDerivorFunctor}, they send a homomorphism $F$ to $F_{\varphi}=(F_{\varphi(s)})_{s\in S}$ and to $F_{\varphi^{\htg(\omega)}}=(F_{(k,\varphi(s))})_{(k,s)\in\omega\times S}$, respectively.

The assignments $\mathbf{c}\mapsto\mathbf{c}^{(\omega)}$ and $\mathbf{c}\mapsto\mathbf{c}^{\htg(\omega)}$ preserve compositions and identities.

\begin{proposition}\label{PChangeSigDerivorCompOmega}\lean{equivalent}{ChangeOfSignature.homogeneous_comp,ChangeOfSignature.homogeneous_identity,ChangeOfSignature.heterogeneous_comp,ChangeOfSignature.heterogeneous_identity}
Let $\sig{\Omega}=(U,\Omega)$ be a many-sorted signature and $\mathbf{d}=(\psi,d)$ a derivor from $\sig{\Lambda}$ to $\sig{\Omega}$.
Then the following equalities hold in $\mathsf{HSig}_{\gd}$
\begin{align*}
(\mathbf{d}\circ\mathbf{c})^{(\omega)}&=\mathbf{d}^{(\omega)}\circ\mathbf{c}^{(\omega)},
&
(\id_{\sig{\Sigma}})^{(\omega)}&=\id_{\sig{\Sigma}^{(\omega)}},
\\
(\mathbf{d}\circ\mathbf{c})^{\htg(\omega)}&=\mathbf{d}^{\htg(\omega)}\circ\mathbf{c}^{\htg(\omega)},
&
(\id_{\sig{\Sigma}})^{\htg(\omega)}&=\id_{\sig{\Sigma}^{\htg(\omega)}}.
\end{align*}
\end{proposition}
\begin{proof}
The proof of Proposition~\ref{PChangeSigDerivorComp} applies with $\omega$ in place of $n$.
In fact, it compares the mappings on sorts, the domains, which are empty, and the assignments of terms symbol by symbol, by means of~(Der1), (Der2) and the formulas of Definition~\ref{DChangeSigDerivor}, which, by Definition~\ref{DChangeSigDerivorOmega}, also define $\mathbf{c}^{(\omega)}$ and $\mathbf{c}^{\htg(\omega)}$; and of Lemma~\ref{LChangeSigTerms}~(1), whose proof, by induction on terms, only uses that $d^{(n)}$ and $d^{\htg(n)}$ send every operation symbol $\lambda$ of $\Lambda$ to $d(\lambda)$ and to $d(\lambda)_{k}$, respectively, and therefore also applies to $\mathbf{d}^{(\omega)}$ and $\mathbf{d}^{\htg(\omega)}$.
\end{proof}

From the above proposition we obtain the corresponding split indexed categories.

\begin{corollary}\label{CChangeSigIndexedOmega}\lean{weaker}{ChangeOfSignature.homogeneousOmegaFunctorComp,ChangeOfSignature.homogeneousOmegaFunctorIdentity,ChangeOfSignature.homogeneousOmegaIndexed,ChangeOfSignature.heterogeneousIndexed}
Let $\sig{\Omega}=(U,\Omega)$ be a many-sorted signature and $\mathbf{d}=(\psi,d)$ a derivor from $\sig{\Lambda}$ to $\sig{\Omega}$.
Then the following equalities hold
\begin{align*}
(\mathbf{d}\circ\mathbf{c})^{(\omega)\ast}_{\gd}&=\mathbf{c}^{(\omega)\ast}_{\gd}\circ\mathbf{d}^{(\omega)\ast}_{\gd},
&
(\id_{\sig{\Sigma}})^{(\omega)\ast}_{\gd}&=\Id_{\omega\mathsf{CatAlg}(\sig{\Sigma})},
\\
(\mathbf{d}\circ\mathbf{c})^{\htg(\omega)\ast}_{\gd}&=\mathbf{c}^{\htg(\omega)\ast}_{\gd}\circ\mathbf{d}^{\htg(\omega)\ast}_{\gd},
&
(\id_{\sig{\Sigma}})^{\htg(\omega)\ast}_{\gd}&=\Id_{\omega\mathsf{Cat}^{\htg}\mathsf{Alg}(\sig{\Sigma})}.
\end{align*}
Consequently, the assignment that sends a many-sorted signature $\sig{\Sigma}$ to $\omega\mathsf{CatAlg}(\sig{\Sigma})$ and a derivor $\mathbf{c}$ to $\mathbf{c}^{(\omega)\ast}_{\gd}$ is a split $\mathsf{Sig}_{\mathfrak{d}}$-indexed category, and so is the assignment that sends $\sig{\Sigma}$ to $\omega\mathsf{Cat}^{\htg}\mathsf{Alg}(\sig{\Sigma})$ and $\mathbf{c}$ to $\mathbf{c}^{\htg(\omega)\ast}_{\gd}$.
We denote these split $\mathsf{Sig}_{\mathfrak{d}}$-indexed categories by $\omega\mathsf{CatAlg}$ and $\omega\mathsf{Cat}^{\htg}\mathsf{Alg}$, respectively.
\end{corollary}
\begin{proof}
Let $\calg{C}$ be an $\omega$-categorial $\sig{\Omega}$-algebra.
Since all the domains involved are empty, by Remark~\ref{RDomainSetEmpty}, both $\mathbf{c}^{(\omega)\ast}_{\gd}(\mathbf{d}^{(\omega)\ast}_{\gd}(\calg{C}))$ and $(\mathbf{d}^{(\omega)}\circ\mathbf{c}^{(\omega)})^{\ast}_{\gd}(\calg{C})$ have $C_{\psi\circ\varphi}$ as underlying $S$-sorted set and, by Proposition~\ref{PChangeSigEpVarietyOmega}~(1), applied twice, the former is a $\sig{\Sigma}^{(\omega)}$-algebra.
Therefore, by Proposition~\ref{PChangeSigDerivorCompOmega} and Corollary~\ref{CGenHepDerivorCompEqual}, $(\mathbf{d}\circ\mathbf{c})^{(\omega)\ast}_{\gd}(\calg{C})=\mathbf{c}^{(\omega)\ast}_{\gd}(\mathbf{d}^{(\omega)\ast}_{\gd}(\calg{C}))$, and the same holds for the $\omega$-categorial $\sig{\Omega}$-homomorphisms.
By Proposition~\ref{PChangeSigDerivorCompOmega} and Corollary~\ref{CGenHepDerivorIdFunctor}, $(\id_{\sig{\Sigma}})^{(\omega)\ast}_{\gd}$ is the identity.
The heterogeneous case is proved in the same way, with $\mathbf{c}^{\htg(\omega)}$, $\mathbf{d}^{\htg(\omega)}$ and $C_{\psi^{\htg(\omega)}\circ\varphi^{\htg(\omega)}}$ in place of $\mathbf{c}^{(\omega)}$, $\mathbf{d}^{(\omega)}$ and $C_{\psi\circ\varphi}$.
Finally, as in the proof of Corollary~\ref{CChangeSigIndexed}, by Corollary~\ref{CSigdCat}, Proposition~\ref{PChangeSigEpVarietyOmega} and the above equalities, both assignments are contravariant functors from $\mathsf{Sig}_{\mathfrak{d}}$ to $\mathsf{Cat}$. This completes the proof.
\end{proof}

We next prove that the change of signature commutes with truncation and discretization in dimension $\omega$.

\begin{proposition}\label{PChangeSigDsTrOmega}\lean{equivalent}{ChangeOfSignature.homogeneousOmegaTruncation_comp,ChangeOfSignature.homogeneousOmegaDiscretization_comp,ChangeOfSignature.heterogeneousDimension_comp}
Let $n\in\mathbb{N}$.
Then the following equalities hold in $\mathsf{HSig}_{\gd}$
\begin{align*}
\mathbf{c}^{(n)}\circ\bDs^{(\omega,n)}_{\sig{\Sigma}}&=\bDs^{(\omega,n)}_{\sig{\Lambda}}\circ\mathbf{c}^{(\omega)},
\\
\mathbf{c}^{(\omega)}\circ\bTr^{(n,\omega)}_{\sig{\Sigma}}&=\bTr^{(n,\omega)}_{\sig{\Lambda}}\circ\mathbf{c}^{(n)},
\\
\mathbf{c}^{\htg(n)}\circ\bDs^{\htg(\omega,n)}_{\sig{\Sigma}}&=\bDs^{\htg(\omega,n)}_{\sig{\Lambda}}\circ\mathbf{c}^{\htg(\omega)},
\\
\mathbf{c}^{\htg(\omega)}\circ\bTr^{\htg(n,\omega)}_{\sig{\Sigma}}&=\bTr^{\htg(n,\omega)}_{\sig{\Lambda}}\circ\mathbf{c}^{\htg(n)}.
\end{align*}
\end{proposition}
\begin{proof}
By Definitions~\ref{DTrDerivorOmega}, \ref{DDsDerivorOmega} and~\ref{DChangeSigDerivorOmega}, the generalized hep-derivors involved are defined by the clauses of Definitions~\ref{DDsDerivor}, \ref{DTrDerivor} and~\ref{DChangeSigDerivor}, with $\omega$ in place of the greatest dimension and $n$ in place of the smallest one.
Moreover, Lemma~\ref{LChangeSigTerms}~(3) and~(4) also hold, with the same proof, for $\bDs^{(\omega,n)}_{\sig{\Lambda}}$, $\bTr^{(n,\omega)}_{\sig{\Lambda}}$, $\bDs^{\htg(\omega,n)}_{\sig{\Lambda}}$ and $\bTr^{\htg(n,\omega)}_{\sig{\Lambda}}$.
Hence the proof of Proposition~\ref{PChangeSigDsTr} applies verbatim.
In the four cases both mappings on sorts are equal: in the two homogeneous cases, both are $\varphi$; in the heterogeneous discretization, both send $(k,s)$ to $(\min\{k,n\},\varphi(s))$; and, in the heterogeneous truncation, both send $(k,s)$ to $(k,\varphi(s))$.
The assignments of terms are compared by the computations of the proof of Proposition~\ref{PChangeSigDsTr}, for every $k\in\omega$ and every pair $j<k$ in $\omega$.
Finally, all the domains are empty, except in the homogeneous truncation, where, by Definition~\ref{DGenHepDerivorCompId}, for every $s\in S$, the following equalities hold
\begin{align*}
\Dom^{\mathbf{c}^{(\omega)}\circ\bTr^{(n,\omega)}_{\sig{\Sigma}}}_{s}
&=
\varnothing\cup(c^{(\omega)})^{\sharp}\left[\left\lbrace\left(v^{s}_{0},\src^{n}_{s}(v^{s}_{0})\right)\right\rbrace\right]
=
\left\lbrace\left(v^{\varphi(s)}_{0},\src^{n}_{\varphi(s)}(v^{\varphi(s)}_{0})\right)\right\rbrace,
\\
\Dom^{\bTr^{(n,\omega)}_{\sig{\Lambda}}\circ\mathbf{c}^{(n)}}_{s}
&=
\Dom^{\Tr^{(n,\omega)}_{\sig{\Lambda}}}_{\varphi(s)}\cup(\Tr^{(n,\omega)}_{\sig{\Lambda}})^{\sharp}[\varnothing]
=
\left\lbrace\left(v^{\varphi(s)}_{0},\src^{n}_{\varphi(s)}(v^{\varphi(s)}_{0})\right)\right\rbrace.
\end{align*}
This completes the proof.
\end{proof}

From the above proposition we obtain that the functors induced by $\mathbf{c}^{(\omega)}$ and $\mathbf{c}^{\htg(\omega)}$ commute with truncation and discretization.

\begin{corollary}\label{CChangeSigDsTrOmega}\lean{weaker}{ChangeOfSignature.homogeneousOmegaTruncationFunctorChange,ChangeOfSignature.homogeneousOmegaDiscretizationFunctorChange,ChangeOfSignature.heterogeneousOmegaTruncationFunctorChange,ChangeOfSignature.heterogeneousOmegaDiscretizationFunctorChange}
Let $n\in\mathbb{N}$.
The following equalities hold
\begin{enumerate}
	\item
	$\Tr^{(n,\omega)}_{\sig{\Sigma}}\circ\mathbf{c}^{(\omega)\ast}_{\gd}=\mathbf{c}^{(n)\ast}_{\gd}\circ\Tr^{(n,\omega)}_{\sig{\Lambda}}$;

	\item
	$\Ds^{(\omega,n)}_{\sig{\Sigma}}\circ\mathbf{c}^{(n)\ast}_{\gd}=\mathbf{c}^{(\omega)\ast}_{\gd}\circ\Ds^{(\omega,n)}_{\sig{\Lambda}}$;

	\item
	$\Tr^{\htg(n,\omega)}_{\sig{\Sigma}}\circ\mathbf{c}^{\htg(\omega)\ast}_{\gd}=\mathbf{c}^{\htg(n)\ast}_{\gd}\circ\Tr^{\htg(n,\omega)}_{\sig{\Lambda}}$; and

	\item
	$\Ds^{\htg(\omega,n)}_{\sig{\Sigma}}\circ\mathbf{c}^{\htg(n)\ast}_{\gd}=\mathbf{c}^{\htg(\omega)\ast}_{\gd}\circ\Ds^{\htg(\omega,n)}_{\sig{\Lambda}}$.
\end{enumerate}
\end{corollary}
\begin{proof}
In the four cases we apply Lemma~\ref{LGenHepDerivorCommSquare} by means of Proposition~\ref{PChangeSigDsTrOmega}.
In each case, let $\calg{B}$ be an algebra in the common domain of the two composites, i.e., an $\omega$-categorial, $n$-categorial, heterogeneous $\omega$-categorial or heterogeneous $n$-categorial $\sig{\Lambda}$-algebra, respectively.

(1) By the description of $\Tr^{(n,\omega)}$ after Definition~\ref{DTrFunctorOmega} and Propositions~\ref{PChangeSigEpVariety} and~\ref{PChangeSigEpVarietyOmega}, the cells of sort $s$ of $\Tr^{(n,\omega)}_{\sig{\Sigma}}(\mathbf{c}^{(\omega)\ast}_{\gd}(\calg{B}))$, of $\mathbf{c}^{(n)\ast}_{\gd}(\Tr^{(n,\omega)}_{\sig{\Lambda}}(\calg{B}))$ and, by the computation of the domains in the proof of Proposition~\ref{PChangeSigDsTrOmega}, of the transformed algebra along the composite are the $n$-cells of sort $\varphi(s)$ of $\calg{B}$; and, by Propositions~\ref{PTrEpVarietyOmega}, \ref{PChangeSigEpVariety} and~\ref{PChangeSigEpVarietyOmega}, the first two are $\sig{\Sigma}^{(n)}$-algebras.

(2) By the description of $\Ds^{(\omega,n)}$ after Definition~\ref{DDsFunctorOmega} and Propositions~\ref{PChangeSigEpVariety} and~\ref{PChangeSigEpVarietyOmega}, the cells of sort $s$ of the three algebras are all the cells of sort $\varphi(s)$ of $\calg{B}$; and, by Propositions~\ref{PDsEpVarietyOmega}, \ref{PChangeSigEpVariety} and~\ref{PChangeSigEpVarietyOmega} and Corollary~\ref{CDsCellsOmega}, the first two are $\sig{\Sigma}^{(\omega)}$-algebras.

(3) By the description of $\Tr^{\htg(n,\omega)}$ after Definition~\ref{DTrFunctorOmega} and Propositions~\ref{PChangeSigEpVariety} and~\ref{PChangeSigEpVarietyOmega}, the cells of sort $(k,s)$ of the three algebras, for $k\in n+1$, are the elements of $B_{(k,\varphi(s))}$; and, by Propositions~\ref{PTrEpVarietyOmega}, \ref{PChangeSigEpVariety} and~\ref{PChangeSigEpVarietyOmega}, the first two are $\sig{\Sigma}^{\htg(n)}$-algebras.

(4) By the description of $\Ds^{\htg(\omega,n)}$ after Definition~\ref{DDsFunctorOmega} and Propositions~\ref{PChangeSigEpVariety} and~\ref{PChangeSigEpVarietyOmega}, the cells of sort $(k,s)$ of the three algebras, for $k\in\omega$, are the elements of $B_{(\min\{k,n\},\varphi(s))}$; and, by Propositions~\ref{PDsEpVarietyOmega}, \ref{PChangeSigEpVariety} and~\ref{PChangeSigEpVarietyOmega}, the first two are $\sig{\Sigma}^{\htg(\omega)}$-algebras.

In the four cases, by Lemma~\ref{LGenHepDerivorCommSquare}, both algebras are equal, and the same holds for the homomorphisms. This completes the proof.
\end{proof}

Graduation and unification in dimension $\omega$ are not induced by generalized hep-derivors, see Remark~\ref{RUfNotDerivorOmega}; we prove their compatibility with the change of signature by means of their definitions.

\begin{corollary}\label{CChangeSigGdUfOmega}\lean{weaker}{ChangeOfSignature.omegaGraduationFunctorChange,ChangeOfSignature.omegaUnificationFunctorChange}
The following equalities hold
\begin{enumerate}
	\item
	$\Gd^{(\omega)}_{\sig{\Sigma}}\circ\mathbf{c}^{(\omega)\ast}_{\gd}=\mathbf{c}^{\htg(\omega)\ast}_{\gd}\circ\Gd^{(\omega)}_{\sig{\Lambda}}$; and

	\item
	$\Uf^{(\omega)}_{\sig{\Sigma}}\circ\mathbf{c}^{\htg(\omega)\ast}_{\gd}=\mathbf{c}^{(\omega)\ast}_{\gd}\circ\Uf^{(\omega)}_{\sig{\Lambda}}$.
\end{enumerate}
\end{corollary}
\begin{proof}
(1) Both sides are functors from $\omega\mathsf{CatAlg}(\sig{\Lambda})$ to $\omega\mathsf{Cat}^{\htg}\mathsf{Alg}(\sig{\Sigma})$.
Let $n\in\mathbb{N}$.
The following chain of equalities holds
\begin{align*}
\Tr^{\htg(n,\omega)}_{\sig{\Sigma}}\circ\Gd^{(\omega)}_{\sig{\Sigma}}\circ\mathbf{c}^{(\omega)\ast}_{\gd}
&=
\Gd^{(n)}_{\sig{\Sigma}}\circ\Tr^{(n,\omega)}_{\sig{\Sigma}}\circ\mathbf{c}^{(\omega)\ast}_{\gd}
\tag{1}
\\&=
\Gd^{(n)}_{\sig{\Sigma}}\circ\mathbf{c}^{(n)\ast}_{\gd}\circ\Tr^{(n,\omega)}_{\sig{\Lambda}}
\tag{2}
\\&=
\mathbf{c}^{\htg(n)\ast}_{\gd}\circ\Gd^{(n)}_{\sig{\Lambda}}\circ\Tr^{(n,\omega)}_{\sig{\Lambda}}
\tag{3}
\\&=
\mathbf{c}^{\htg(n)\ast}_{\gd}\circ\Tr^{\htg(n,\omega)}_{\sig{\Lambda}}\circ\Gd^{(\omega)}_{\sig{\Lambda}}
\tag{4}
\\&=
\Tr^{\htg(n,\omega)}_{\sig{\Sigma}}\circ\mathbf{c}^{\htg(\omega)\ast}_{\gd}\circ\Gd^{(\omega)}_{\sig{\Lambda}}.
\tag{5}
\end{align*}
The equalities~(1) and~(4) follow from Definition~\ref{DGdFunctorOmega}; the equalities~(2) and~(5) follow from Corollary~\ref{CChangeSigDsTrOmega}~(1) and~(3), respectively; finally, the equality~(3) follows from Corollary~\ref{CChangeSigGdUf}~(1).
Hence both functors have the same truncations and, by the uniqueness in Proposition~\ref{PProjLimit}~(2), they are equal.

(2) Let $\calg{B}$ be a heterogeneous $\omega$-categorial $\sig{\Lambda}$-algebra and let $\calg{C}$ stand for $\mathbf{c}^{\htg(\omega)\ast}_{\gd}(\calg{B})$.
By Proposition~\ref{PChangeSigEpVarietyOmega}~(2), for every $k\in\omega$ and every $s\in S$, $C_{(k,s)}=B_{(k,\varphi(s))}$ and, for every $j<k$ in $\omega$, $\mathrm{i}^{(k,j)\,\calg{C}}_{s}=\mathrm{i}^{(k,j)\,\calg{B}}_{\varphi(s)}$.
Hence, by Definition~\ref{DUfCellsOmega}, the relation $\Theta_{s}$ of $\calg{C}$ is the relation $\Theta_{\varphi(s)}$ of $\calg{B}$, and $[C]_{s}=[B]_{\varphi(s)}$.
Therefore, by Definition~\ref{DUfFunctorOmega} and Proposition~\ref{PChangeSigEpVarietyOmega}~(1), both $\Uf^{(\omega)}_{\sig{\Sigma}}(\calg{C})$ and $\mathbf{c}^{(\omega)\ast}_{\gd}(\Uf^{(\omega)}_{\sig{\Lambda}}(\calg{B}))$ have $([B]_{\varphi(s)})_{s\in S}$ as underlying $S$-sorted set.

We next compare the operations.
Let $\tau$ be an operation symbol of $\Sigma^{(\omega)}$, of arity $\mathbf{s}$, and $([x_{j}]_{p})_{j\in\bb{\mathbf{s}}}$ a family of $\omega$-cell classes represented in a common level $p$ such that $\tau$ is an operation symbol of $\Sigma^{(p)}$, see Lemma~\ref{LUfCellsOmega}~(3).
Let $\calg{U}_{p}$ stand for $\Uf^{(p)}_{\sig{\Lambda}}(\Tr^{\htg(p,\omega)}_{\sig{\Lambda}}(\calg{B}))$.
By Corollaries~\ref{CChangeSigDsTrOmega}~(3) and~\ref{CChangeSigGdUf}~(2), the following equalities hold
$$
\Uf^{(p)}_{\sig{\Sigma}}\left(\Tr^{\htg(p,\omega)}_{\sig{\Sigma}}(\calg{C})\right)
=
\Uf^{(p)}_{\sig{\Sigma}}\left(\mathbf{c}^{\htg(p)\ast}_{\gd}\left(\Tr^{\htg(p,\omega)}_{\sig{\Lambda}}(\calg{B})\right)\right)
=
\mathbf{c}^{(p)\ast}_{\gd}(\calg{U}_{p}),
$$
in which, by Proposition~\ref{PChangeSigEpVariety}~(1), $\tau$ acts as the term operation $c^{(p)}(\tau)^{\calg{U}_{p}}$.
Hence, by Definition~\ref{DUfFunctorOmega}, $\tau$ is defined in $\Uf^{(\omega)}_{\sig{\Sigma}}(\calg{C})$ at $([x_{j}]_{p})_{j\in\bb{\mathbf{s}}}$ if, and only if, $c^{(p)}(\tau)^{\calg{U}_{p}}$ is defined at $(x_{j})_{j\in\bb{\mathbf{s}}}$, and then its value is $[c^{(p)}(\tau)^{\calg{U}_{p}}((x_{j})_{j\in\bb{\mathbf{s}}})]_{p}$.
On the other hand, by Proposition~\ref{PChangeSigEpVarietyOmega}~(1), $\tau$ acts in $\mathbf{c}^{(\omega)\ast}_{\gd}(\Uf^{(\omega)}_{\sig{\Lambda}}(\calg{B}))$ as the term operation of $c^{(\omega)}(\tau)=c^{(p)}(\tau)$ in $\Uf^{(\omega)}_{\sig{\Lambda}}(\calg{B})$.
By Proposition~\ref{PUfFunctorOmega} and the description of $\Tr^{(p,\omega)}_{\sig{\Lambda}}$ after Definition~\ref{DTrFunctorOmega}, $\mathrm{in}^{[p]}$ is a homomorphism of partial $\Lambda^{(p)}$-algebras from $\calg{U}_{p}$ to the reduct of $\Uf^{(\omega)}_{\sig{\Lambda}}(\calg{B})$ to $\Lambda^{(p)}$ that reflects definedness, since its image is the closed relative subalgebra of the $p$-cells.
Hence $c^{(p)}(\tau)$ is defined in $\Uf^{(\omega)}_{\sig{\Lambda}}(\calg{B})$ at $([x_{j}]_{p})_{j\in\bb{\mathbf{s}}}$ if, and only if, it is defined in $\calg{U}_{p}$ at $(x_{j})_{j\in\bb{\mathbf{s}}}$, and then its value is the class at level $p$ of the value in $\calg{U}_{p}$.
Therefore both algebras are equal.

Finally, let $F$ be a heterogeneous $\omega$-categorial $\sig{\Lambda}$-homomorphism.
By Definition~\ref{DUfFunctorOmega} and Proposition~\ref{PGenHepDerivorFunctor}, both $\Uf^{(\omega)}_{\sig{\Sigma}}(\mathbf{c}^{\htg(\omega)\ast}_{\gd}(F))$ and $\mathbf{c}^{(\omega)\ast}_{\gd}(\Uf^{(\omega)}_{\sig{\Lambda}}(F))$ send, at every sort $s\in S$, a class $[x]_{k}$, with $x\in B_{(k,\varphi(s))}$, to $[F_{(k,\varphi(s))}(x)]_{k}$. This completes the proof.
\end{proof}

From Corollaries~\ref{CChangeSigIndexed}, \ref{CChangeSigIndexedOmega}, \ref{CChangeSigDsTrOmega} and~\ref{CChangeSigGdUfOmega} we obtain that the functors in dimension $\omega$ are natural transformations between split indexed categories.

\begin{corollary}\label{CChangeSigNatOmega}\lean{weaker}{ChangeOfSignature.omegaGraduationIndexed,ChangeOfSignature.omegaUnificationIndexed,ChangeOfSignature.homogeneousOmegaTruncationNaturality,ChangeOfSignature.homogeneousOmegaDiscretizationNaturality,ChangeOfSignature.heterogeneousOmegaTruncationNaturality,ChangeOfSignature.heterogeneousOmegaDiscretizationNaturality}
Let $n\in\mathbb{N}$.
Then the families
\begin{enumerate}
	\item
	$\Gd^{(\omega)}=(\Gd^{(\omega)}_{\sig{\Sigma}})_{\sig{\Sigma}\in\mathrm{Ob}(\mathsf{Sig}_{\mathfrak{d}})}\colon\omega\mathsf{CatAlg}\Longrightarrow\omega\mathsf{Cat}^{\htg}\mathsf{Alg}$;

	\item
	$\Uf^{(\omega)}=(\Uf^{(\omega)}_{\sig{\Sigma}})_{\sig{\Sigma}\in\mathrm{Ob}(\mathsf{Sig}_{\mathfrak{d}})}\colon\omega\mathsf{Cat}^{\htg}\mathsf{Alg}\Longrightarrow\omega\mathsf{CatAlg}$;

	\item
	$\Tr^{(n,\omega)}=(\Tr^{(n,\omega)}_{\sig{\Sigma}})_{\sig{\Sigma}\in\mathrm{Ob}(\mathsf{Sig}_{\mathfrak{d}})}\colon\omega\mathsf{CatAlg}\Longrightarrow\mathsf{nCatAlg}$;

	\item
	$\Ds^{(\omega,n)}=(\Ds^{(\omega,n)}_{\sig{\Sigma}})_{\sig{\Sigma}\in\mathrm{Ob}(\mathsf{Sig}_{\mathfrak{d}})}\colon\mathsf{nCatAlg}\Longrightarrow\omega\mathsf{CatAlg}$;

	\item
	$\Tr^{\htg(n,\omega)}=(\Tr^{\htg(n,\omega)}_{\sig{\Sigma}})_{\sig{\Sigma}\in\mathrm{Ob}(\mathsf{Sig}_{\mathfrak{d}})}\colon\omega\mathsf{Cat}^{\htg}\mathsf{Alg}\Longrightarrow\mathsf{nCat}^{\htg}\mathsf{Alg}$; and

	\item
	$\Ds^{\htg(\omega,n)}=(\Ds^{\htg(\omega,n)}_{\sig{\Sigma}})_{\sig{\Sigma}\in\mathrm{Ob}(\mathsf{Sig}_{\mathfrak{d}})}\colon\mathsf{nCat}^{\htg}\mathsf{Alg}\Longrightarrow\omega\mathsf{Cat}^{\htg}\mathsf{Alg}$
\end{enumerate}
are natural transformations.
\end{corollary}
\begin{proof}
By Corollaries~\ref{CChangeSigIndexed} and~\ref{CChangeSigIndexedOmega}, the assignments $\mathsf{nCatAlg}$, $\mathsf{nCat}^{\htg}\mathsf{Alg}$, $\omega\mathsf{CatAlg}$ and $\omega\mathsf{Cat}^{\htg}\mathsf{Alg}$ are contravariant functors from $\mathsf{Sig}_{\mathfrak{d}}$ to $\mathsf{Cat}$.
By Definitions~\ref{DGdFunctorOmega}, \ref{DTrFunctorOmega} and~\ref{DDsFunctorOmega} and Proposition~\ref{PUfFunctorOmega}, for every many-sorted signature $\sig{\Sigma}$, the components at $\sig{\Sigma}$ of the six families are functors between the corresponding categories.
Finally, for all many-sorted signatures $\sig{\Sigma}$ and $\sig{\Lambda}$ and every derivor $\mathbf{c}$ from $\sig{\Sigma}$ to $\sig{\Lambda}$, the naturality conditions at $\mathbf{c}$ are the equalities of Corollaries~\ref{CChangeSigGdUfOmega} and~\ref{CChangeSigDsTrOmega}.
\end{proof}

As in Section~\ref{SChangeSigCompat}, the natural transformations of Chapter~\ref{SFunctorsOmega} are also compatible with the change of signature.

\begin{remark}\label{RChangeSigNatOmega}
Let $m\in\mathbb{N}$.
The natural isomorphisms $\zeta$ and $\theta$ of Corollary~\ref{CGdUfEquivOmega} and $\mathrm{in}^{[m]}$ of Proposition~\ref{PUfFunctorOmega}, and the counits $\varepsilon$ and $\varepsilon^{\htg}$ of Proposition~\ref{PDsTrAdjOmega}, are given, at every algebra, by families of identities, inclusions, $(k,m)$-identities and the mappings that send a cell $x$ of dimension $k$ to its class $[x]_{k}$, or such a class to $x$, indexed by the sorts.
By Proposition~\ref{PChangeSigEpVarietyOmega} and Corollary~\ref{CChangeSigGdUfOmega}~(2), the functors $\mathbf{c}^{(\omega)\ast}_{\gd}$ and $\mathbf{c}^{\htg(\omega)\ast}_{\gd}$ reindex the sorts along $\varphi$ and keep the cells, the classes $[x]_{k}$, the sources, the targets, the identities and the compositions.
Therefore, by Corollaries~\ref{CChangeSigDsTrOmega} and~\ref{CChangeSigGdUfOmega}, with the notation of Remark~\ref{RChangeSigNat}, the following equalities hold
\begin{align*}
\zeta^{\sig{\Sigma}}_{\mathbf{c}^{(\omega)\ast}_{\gd}(\calg{A})}&=\mathbf{c}^{(\omega)\ast}_{\gd}\left(\zeta^{\sig{\Lambda}}_{\calg{A}}\right),
&
\theta^{\sig{\Sigma}}_{\mathbf{c}^{\htg(\omega)\ast}_{\gd}(\calg{B})}&=\mathbf{c}^{\htg(\omega)\ast}_{\gd}\left(\theta^{\sig{\Lambda}}_{\calg{B}}\right),
\\
\varepsilon^{\sig{\Sigma}}_{\mathbf{c}^{(\omega)\ast}_{\gd}(\calg{A})}&=\mathbf{c}^{(\omega)\ast}_{\gd}\left(\varepsilon^{\sig{\Lambda}}_{\calg{A}}\right),
&
\varepsilon^{\htg,\sig{\Sigma}}_{\mathbf{c}^{\htg(\omega)\ast}_{\gd}(\calg{B})}&=\mathbf{c}^{\htg(\omega)\ast}_{\gd}\left(\varepsilon^{\htg,\sig{\Lambda}}_{\calg{B}}\right),
\end{align*}
and $\mathrm{in}^{[m],\sig{\Sigma}}_{\mathbf{c}^{\htg(\omega)\ast}_{\gd}(\calg{B})}=\mathbf{c}^{(m)\ast}_{\gd}(\mathrm{in}^{[m],\sig{\Lambda}}_{\calg{B}})$, for every $\omega$-categorial $\sig{\Lambda}$-algebra $\calg{A}$ and every heterogeneous $\omega$-categorial $\sig{\Lambda}$-algebra $\calg{B}$.
For instance, $\theta^{\sig{\Sigma}}_{\mathbf{c}^{\htg(\omega)\ast}_{\gd}(\calg{B})}$ is, at the sort $(k,s)$, the mapping $\mathrm{in}^{[k]}_{\varphi(s)}$ of $\calg{B}$, which sends $x\in B_{(k,\varphi(s))}$ to $[x]_{k}$, and this is the component at $(k,s)$ of $(\theta^{\sig{\Lambda}}_{\calg{B}})_{\varphi^{\htg(\omega)}}$.
Thus the adjoint equivalence of Corollary~\ref{CGdUfEquivOmega}, the adjunctions of Proposition~\ref{PDsTrAdjOmega} and the algebraic square of Proposition~\ref{PAlgSquareOmega} are compatible with the change of signature.
\end{remark}

Moreover, by Remark~\ref{RGrothendieckNat}, as in Remark~\ref{RChangeSigGrothendieck}, the natural transformations of Corollary~\ref{CChangeSigNatOmega} extend to functors between the Grothendieck constructions of the split indexed categories of Corollaries~\ref{CChangeSigIndexed} and~\ref{CChangeSigIndexedOmega} that commute with the canonical split fibrations.

Summing up, every derivor $\mathbf{c}$ determines, for every $k\in\mathbb{N}\cup\{\omega\}$, generalized hep-derivors $\mathbf{c}^{(k)}$ and $\mathbf{c}^{\htg(k)}$ whose induced functors send categorial algebras to categorial algebras, preserve compositions and identities, and commute with graduation, unification, discretization and truncation; moreover, the natural transformations, adjunctions and algebraic squares that relate these functors are compatible with the change of signature.
Together with Chapters~\ref{SGdUf}--\ref{SFunctorsOmega}, this completes the three directions of Figure~\ref{FFunctorsSummary}, which we discuss in the concluding remarks.

\chapter*{Concluding remarks}

\begin{figure}[tp]
\centering
\begin{tikzpicture}
    [ACliment/.style={-{To [angle'=45, length=5.75pt, width=4pt, round]}},
    back/.style={ACliment, gray!70},
    depth/.style={ACliment, densely dashed},
    lbl/.style={font=\scriptsize, fill=white, inner sep=1pt},
    cat/.style={fill=white, inner sep=2pt},
    adj/.style={font=\footnotesize, inner sep=0pt},
    scale=1]
    \begin{scope}[shift={(3.9,2.3)}]
        \node[cat] (Lm) at (0,7.4) {$\mathsf{mCatAlg}(\sig{\Lambda})$};
        \node[cat] (Ln) at (0,3.7) {$\mathsf{nCatAlg}(\sig{\Lambda})$};
        \node[cat] (Lw) at (0,0) {$\omega\mathsf{CatAlg}(\sig{\Lambda})$};
        \node[cat] (Lhm) at (5.8,7.4) {$\mathsf{mCat}^{\htg}\mathsf{Alg}(\sig{\Lambda})$};
        \node[cat] (Lhn) at (5.8,3.7) {$\mathsf{nCat}^{\htg}\mathsf{Alg}(\sig{\Lambda})$};
        \node[cat] (Lhw) at (5.8,0) {$\omega\mathsf{Cat}^{\htg}\mathsf{Alg}(\sig{\Lambda})$};
    \end{scope}
    \foreach \a/\b in {Lm/Lhm, Ln/Lhn, Lw/Lhw}{
        \draw[back] ([yshift=6pt]\a.east) to ([yshift=6pt]\b.west);
        \draw[back] ([yshift=-6pt]\b.west) to ([yshift=-6pt]\a.east);}
    \foreach \a/\b in {Lm/Ln, Ln/Lw, Lhm/Lhn, Lhn/Lhw}{
        \draw[back] ([xshift=6pt]\a.south) to ([xshift=6pt]\b.north);
        \draw[back] ([xshift=-6pt]\b.north) to ([xshift=-6pt]\a.south);}
    \node[cat] (Sm) at (0,7.4) {$\mathsf{mCatAlg}(\sig{\Sigma})$};
    \node[cat] (Sn) at (0,3.7) {$\mathsf{nCatAlg}(\sig{\Sigma})$};
    \node[cat] (Sw) at (0,0) {$\omega\mathsf{CatAlg}(\sig{\Sigma})$};
    \node[cat] (Shm) at (5.8,7.4) {$\mathsf{mCat}^{\htg}\mathsf{Alg}(\sig{\Sigma})$};
    \node[cat] (Shn) at (5.8,3.7) {$\mathsf{nCat}^{\htg}\mathsf{Alg}(\sig{\Sigma})$};
    \node[cat] (Shw) at (5.8,0) {$\omega\mathsf{Cat}^{\htg}\mathsf{Alg}(\sig{\Sigma})$};
    \foreach \a/\b/\k in {Lm/Sm/(m), Ln/Sn/(n), Lw/Sw/(\omega), Lhm/Shm/\htg(m), Lhn/Shn/\htg(n), Lhw/Shw/\htg(\omega)}{
        \draw[depth] (\a) to node[lbl, pos=0.45] {$\mathbf{c}^{\k\ast}_{\gd}$} (\b);}
    \foreach \a/\b/\k in {Sm/Shm/m, Sn/Shn/n, Sw/Shw/\omega}{
        \draw[ACliment] ([yshift=6pt]\a.east) to node[lbl, above] {$\Gd^{(\k)}_{\sig{\Sigma}}$} ([yshift=6pt]\b.west);
        \draw[ACliment] ([yshift=-6pt]\b.west) to node[lbl, below] {$\Uf^{(\k)}_{\sig{\Sigma}}$} ([yshift=-6pt]\a.east);}
    \foreach \a/\b/\j/\k in {Sm/Sn/m/n, Sn/Sw/n/\omega}{
        \draw[ACliment] ([xshift=6pt]\a.south) to node[lbl, right, pos=0.6] {$\Ds^{(\k,\j)}_{\sig{\Sigma}}$} ([xshift=6pt]\b.north);
        \draw[ACliment] ([xshift=-6pt]\b.north) to node[lbl, left, pos=0.4] {$\Tr^{(\j,\k)}_{\sig{\Sigma}}$} ([xshift=-6pt]\a.south);}
    \foreach \a/\b/\j/\k in {Shm/Shn/m/n, Shn/Shw/n/\omega}{
        \draw[ACliment] ([xshift=6pt]\a.south) to node[lbl, right, pos=0.6] {$\Ds^{\htg(\k,\j)}_{\sig{\Sigma}}$} ([xshift=6pt]\b.north);
        \draw[ACliment] ([xshift=-6pt]\b.north) to node[lbl, left, pos=0.4] {$\Tr^{\htg(\j,\k)}_{\sig{\Sigma}}$} ([xshift=-6pt]\a.south);}
    \foreach \a/\b in {Sm/Shm, Sn/Shn, Sw/Shw}{
        \node[adj] at ($(\a.east)!0.5!(\b.west)$) {\rotatebox{90}{$\vdash$}};}
    \foreach \a/\b in {Sm/Sn, Sn/Sw, Shm/Shn, Shn/Shw}{
        \node[adj] at ($(\a.south)!0.5!(\b.north)$) {$\vdash$};}
\end{tikzpicture}
\caption{The categories of higher-order categorial algebras and the functors between them, for $m\leq n$ in $\mathbb{N}$ and a derivor $\mathbf{c}$ from $\sig{\Sigma}$ to $\sig{\Lambda}$: the horizontal arrows change the presentation, the vertical arrows change the dimension and the dashed arrows change the signature. On the front face, the symbols $\perp$ and $\vdash$ indicate the adjunctions $\Gd\dashv\Uf$, which are adjoint equivalences, and $\Ds\dashv\Tr$, respectively. In dimension $\omega$, the heterogeneous categorial algebras form the category canonically associated with an ep-variety, whereas, if the set of sorts is nonempty, the homogeneous ones do not, see Proposition~\ref{PFiniteDimNotQE}.}\label{FFunctorsSummary}
\end{figure}
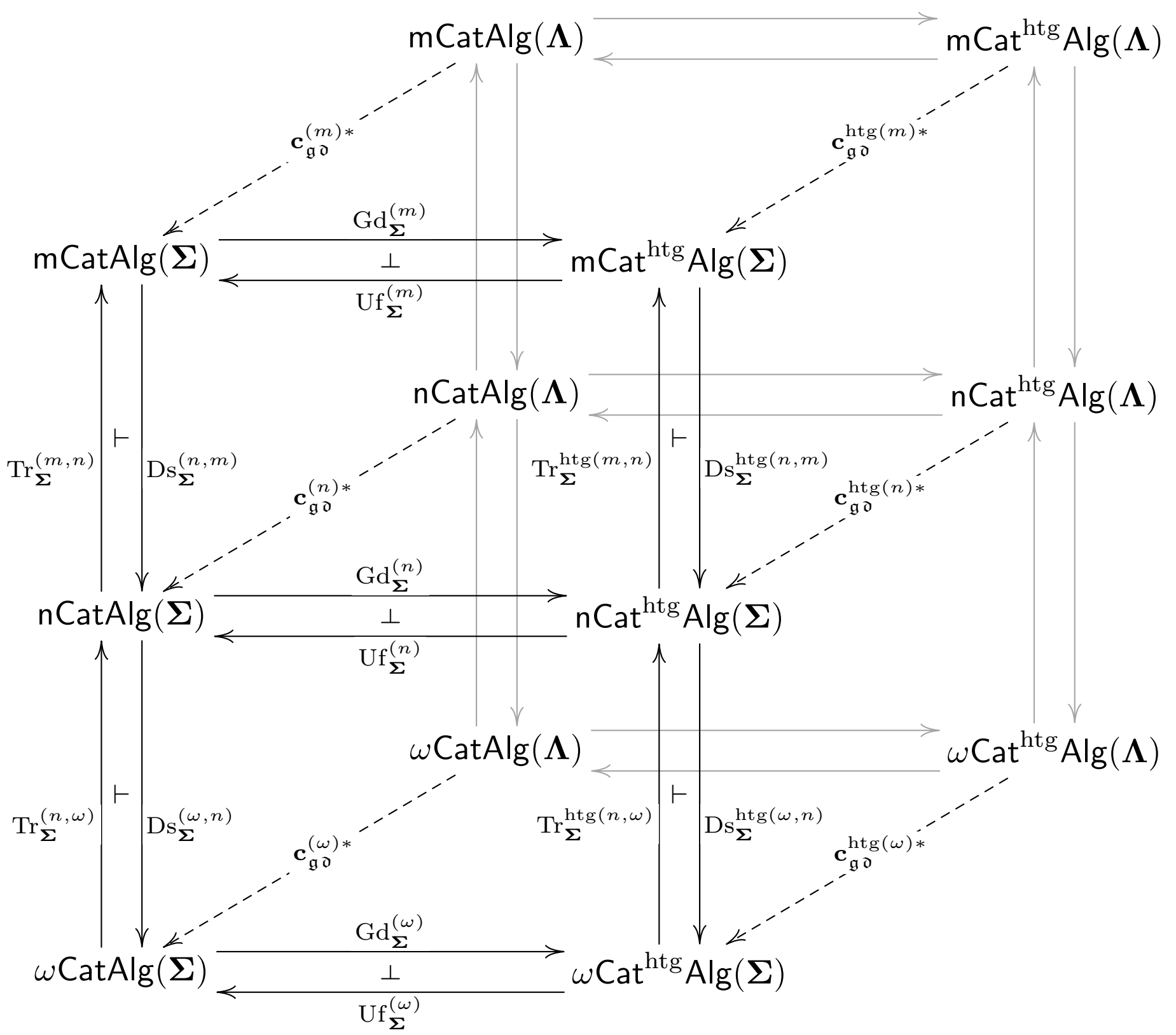

Figure~\ref{FFunctorsSummary} gathers the categories of higher-order categorial algebras of this work and the functors between them, for $m\leq n$ in $\mathbb{N}$ and a derivor $\mathbf{c}=(\varphi,c)$ from $\sig{\Sigma}$ to $\sig{\Lambda}$.
It can be read in three independent directions.
The horizontal arrows change the presentation: graduation and unification are mutually quasi-inverse equivalences, see Corollaries~\ref{CGdUfEquiv} and~\ref{CGdUfEquivOmega}.
The vertical arrows change the dimension: discretization is left adjoint to truncation, see Propositions~\ref{PDsTrAdj} and~\ref{PDsTrAdjOmega}, and the categories of dimension $\omega$ are the projective limits of the categories of finite dimension and the truncation functors, see Chapter~\ref{SProjLimits}; together with graduation and unification, these adjunctions form the algebraic squares of Propositions~\ref{PAlgSquareFinite} and~\ref{PAlgSquareOmega}.
The arrows in depth change the signature: they are induced by the generalized hep-derivors $\mathbf{c}^{(k)}$ and $\mathbf{c}^{\htg(k)}$, for $k\in\{m,n,\omega\}$, and go in the direction opposite to that of $\mathbf{c}$, see Chapter~\ref{SChangeSig}.
The back face, for $\sig{\Lambda}$, contains the same functors as the front one, with $\sig{\Lambda}$ in place of $\sig{\Sigma}$.
Every square formed by an arrow in depth and a horizontal or vertical arrow commutes, see Corollaries~\ref{CChangeSigGdUf}, \ref{CChangeSigDsTr}, \ref{CChangeSigDsTrOmega} and~\ref{CChangeSigGdUfOmega}; this means that graduation, unification, discretization and truncation, in finite dimensions and in dimension $\omega$, are natural transformations between split $\mathsf{Sig}_{\mathfrak{d}}$-indexed categories, see Corollaries~\ref{CChangeSigGdUfNat}, \ref{CChangeSigDsTrNat} and~\ref{CChangeSigNatOmega}.

The generalized hep-derivors on which these functors rest separate three ingredients: the interpretation of the operations as term operations on the selected carriers, the domains of these operations, which are evaluated recursively in the transformed algebra, and the satisfaction of the source axioms.
There is, in general, only a natural inclusion of the composite of the induced functors into the functor induced by a composite of generalized hep-derivors, see Proposition~\ref{PGenHepDerivorComp}: the composite can keep cells that the successive transformations remove, and the inclusion can be strict, see Example~\ref{EStrictComparison}.
This is why the induced functors form a normal lax indexed category with values in the categories of weak ep-algebras, see Proposition~\ref{PHepLaxIndexed}.
Example~\ref{EDiscretizationDiagonal} shows, in turn, why the recursively evaluated domains are needed: in a discretization functor, the restriction of a composition to the diagonal is imposed by a domain condition alone, not by a term operation.
In the categorial applications, the comparison of the carriers and the exactness of the domains, see Corollary~\ref{CGenHepDerivorCompEqual}, turn these inclusions into the equalities or the natural isomorphisms used in the equivalences, adjunctions and algebraic squares.

The difference between finite and infinite dimensions appears in unification.
In finite dimensions, unification keeps the cells of top dimension and is induced by a generalized hep-derivor, see Definition~\ref{DUfDerivor}.
In dimension $\omega$, it assembles all the finite levels by identifying every cell with its iterated identities and, when the set of sorts is nonempty, it is not induced by any generalized hep-derivor, see Remark~\ref{RUfNotDerivorOmega}; it is obtained instead on an inductive limit of sets of cells, see Definition~\ref{DUfFunctorOmega}.
In contrast, truncation and discretization in dimension $\omega$ are induced by generalized hep-derivors, see Definitions~\ref{DTrFunctorOmega} and~\ref{DDsFunctorOmega}, and graduation in dimension $\omega$ is obtained from the universal property of the projective limit, see Definition~\ref{DGdFunctorOmega}.
The dimension $\omega$ also breaks the symmetry between the two presentations.
The heterogeneous $\omega$-categorial algebras form the category canonically associated with an ep-variety, see Remark~\ref{RCatAlgAsEpAlgHtg}, as in finite dimensions.
In contrast, if the set of sorts is nonempty, the homogeneous $\omega$-categorial algebras form a proper full subcategory of the corresponding ep-variety, since the condition that every cell has a finite dimension is not preserved by products and therefore cannot be imposed by $\QE$-equations, see Proposition~\ref{PFiniteDimNotQE} and Corollary~\ref{CProperSubCatAlgOmega}.
This is why, in the homogeneous presentation, the projective limits of Chapter~\ref{SProjLimits} give the description of dimension $\omega$ that the ep-varieties do not provide.

Several questions remain open.
First, for strict higher-order categories, truncation also has a right adjoint, which adds a unique cell between every pair of parallel cells of the top dimension.
We do not know whether $\Tr^{(m,n)}_{\sig{\Sigma}}$ and $\Tr^{\htg(m,n)}_{\sig{\Sigma}}$ have right adjoints in the presence of the operations of $\Sigma$ and, if so, whether they are induced by generalized hep-derivors; as with unification in dimension $\omega$, see Remark~\ref{RUfNotDerivorOmega}, their cells would not be elements of the carriers of the given algebra.
Second, for total algebras, the functor induced by a derivor has a left adjoint.
It remains to determine under which conditions the functors $\mathbf{c}^{(n)\ast}_{\gd}$ and $\mathbf{c}^{\htg(n)\ast}_{\gd}$ have left adjoints, which is related to the last question below.
Third, it would be interesting to characterize the functors between categories of partial algebras that are induced by generalized hep-derivors; Remark~\ref{RUfNotDerivorOmega} shows that unification in dimension $\omega$ is not one of them.
The results concern strict higher-order categories whose base operations are finitary, although the sets of domain conditions of the operations may be infinite; weak higher-order categories would require additional coherence data.
Finally, finding hypotheses under which the ep-varieties have finitary essentially algebraic presentations with useful presentability properties remains an open question.

\bibliographystyle{amsplain-initials}
\bibliography{references}

\end{document}